\documentclass[reqno]{article}
\usepackage[T1]{fontenc}
\usepackage{amsmath}
\usepackage{mathtools}
\usepackage{hyperref, bbm, amssymb, url, enumerate,caption}
\usepackage{amsthm}
\usepackage{enumitem}

\usepackage[a4paper,top=3cm,bottom=3cm,left=3cm,right=3cm, includeheadfoot]{geometry}
\newlist{romanenum}{enumerate}{1}
\setlist[romanenum]{label=(\roman*),ref=(\roman*)}
\usepackage[scr]{rsfso}
\usepackage{comment}
\usepackage{bm}
\usepackage{tikz-cd}
\usetikzlibrary{matrix,arrows,decorations.pathmorphing, cd}
\usepackage[style=alphabetic, backend=biber]{biblatex}
\usepackage{pgfplots}
\pgfplotsset{compat=1.18}
\usepackage{todonotes}
\usepackage{accents}
\usepackage{trimclip}
\usepackage{thmtools}

\DeclareMathAlphabet{\mathpzc}{OT1}{pzc}{m}{it}

\newcommand{\cshv}[1]{c\left(#1\right)}

\usepackage[capitalise]{cleveref}

\theoremstyle{plain}
\newtheorem{thm}{Theorem}[subsection]
\newtheorem{cor}[thm]{Corollary}
\newtheorem{lemma}[thm]{Lemma}
\newtheorem{proposition}[thm]{Proposition}
\newtheorem{thmintro}{Theorem}

\newtheorem{lemmaintro}[thmintro]{Lemma}

\newtheorem*{propositionintro}{Proposition}
\theoremstyle{definition}
\newtheorem{definition}[thm]{Definition}
\newtheorem{construction}[thm]{Construction}
\newtheorem{recollection}[thm]{Recollection}

\theoremstyle{remark}
\newtheorem{rem}[thm]{Remark}
\newtheorem{notation}[thm]{Notation}
\newtheorem{rems}[thm]{Remarks}
\newtheorem*{remintro}{Remark}
\newtheorem{ex}[thm]{Example}

\crefname{thm}{Theorem}{Theorems}
\crefname{cor}{Corollary}{Corollaries}
\crefname{lemma}{Lemma}{Lemmas}
\crefname{proposition}{Proposition}{Propositions}
\crefname{thmintro}{Theorem}{Theorems}
\crefname{lemmaintro}{Lemma}{Lemmas}
\crefname{cor}{Corollary}{Corollaries}
\crefname{definition}{Definition}{Definitions}
\crefname{construction}{Construction}{Constructions}
\crefname{recollection}{Recollection}{Recollections}
\crefname{rem}{Remark}{Remarks}
\crefname{notation}{Notation}{Notations}
\crefname{rems}{Remarks}{Remarks}
\crefname{ex}{Example}{Examples}
\crefname{exs}{Examples}{Examples}

\newcommand{\Cov}{\operatorname{Cov}}
\newcommand{\ladg}[1]{{#1}^*}
\newcommand{\radg}[1]{{#1}_*}
\newcommand{\stradg}[1]{{#1}_{*, \Sp}}

\newcommand{\stladg}[1]{{#1}^*_{\Sp}}

\newcommand{\xradg}[2]{{#1}^{#2}_{*}}
\newcommand{\xladg}[2]{{#1}^*_{#2}}

\newcommand{\Gradg}[1]{{#1}^G_{*}}
\newcommand{\Gladg}[1]{{#1}^{G,*}}
\newcommand{\Gstradg}[1]{{#1}^G_{*, \Sp}}
\newcommand{\xstradg}[2]{{#1}^{#2}_{*, \Sp}}

\newcommand{\xstladg}[2]{{#1}^{*}_{#2,\Sp}}
\newcommand{\cof}{\operatorname{cof}}
\newcommand{\Cone}{\operatorname{Cone}}

\newcommand{\an}{\operatorname{An}}
\newcommand{\Pro}{\operatorname{Pro}}
\newcommand{\Ind}{\operatorname{Ind}}
\newcommand{\Fin}{\operatorname{\operatorname{Fin}}}
\newcommand{\Set}{\operatorname{Set}}
\newcommand{\Mon}{\operatorname{Mon}}

\newcommand{\Grp}{\operatorname{Grp}}
\newcommand{\CGrp}{\operatorname{CGrp}}
\newcommand{\Ab}{\operatorname{Ab}}
\newcommand{\CAlg}{\operatorname{CAlg}}
\newcommand{\aCAlg}{\operatorname{(C)Alg}} 
\newcommand{\Comm}{\operatorname{Comm}}
\newcommand{\Alg}{\operatorname{Alg}}
\newcommand{\Assoc}{\operatorname{Assoc}}
\newcommand{\LM}{\mathcal{L}\mathcal M}
\newcommand{\RM}{\mathcal{R}\mathcal M}

\newcommand{\BM}{\mathcal{B}\mathcal M}
\newcommand{\LMon}{\operatorname{LMon}}
\newcommand{\colim}[1]{{\underset{#1}{\operatorname{colim}}}}

\newcommand{\clim}[1]{{\underset{#1}{\operatorname{lim}}}}
\newcommand{\St}{\operatorname{\operatorname{st}}}
\newcommand{\Hom}{\operatorname{Hom}}
\newcommand{\Map}{\operatorname{Map}}
\newcommand{\Fun}{\operatorname{Fun}}

\newcommand{\map}{\operatorname{map}}
\newcommand{\imap}{\underline{\operatorname{map}}}
\newcommand{\iMap}{\underline{\operatorname{Map}}}
\newcommand{\iHom}{\underline{\operatorname{Hom}}}
\newcommand{\Cat}{\operatorname{Cat}_{\infty}}
\newcommand{\vlCat}{\widehat{\operatorname{Cat}}_{\infty}}
\newcommand{\cat}[1]{\mathcal{#1}}

\newcommand{\Coker}{\operatorname{Coker}}
\newcommand{\Ext}{\operatorname{Ext}}
\newcommand{\Fib}{\operatorname{Fib}}
\newcommand{\Eq}{\operatorname{Eq}}
\newcommand{\im}{\operatorname{Im}}
\newcommand{\Cofib}{\operatorname{Cofib}}
\newcommand{\essim}{\operatorname{ess. im}}
\newcommand{\oc}[2]{{{#1}_{/#2}}}
\newcommand{\stab}[1]{\operatorname{Sp}(#1)}
\newcommand{\Sp}{\operatorname{Sp}}

\newcommand{\heart}{\heartsuit}

\newcommand{\gr}{\operatorname{gr}}
\newcommand{\Coeq}{\operatorname{Coeq}}
\newcommand{\grAb}{\operatorname{grAb}}
\newcommand{\Ch}{\operatorname{Ch}}
\newcommand{\Tkcont}{\operatorname{T1 Top}_k}

\newcommand{\mt}[1]{\mathcal{#1}}
\newcommand{\Op}{\operatorname{Op}}
\newcommand{\BMod}[2]{\operatorname{LMod}_{#1}(\operatorname{Sp}(\mt{B}_{#2}))}

\newcommand{\LMod}[2]{\operatorname{LMod}_{#1}({#2})}
\newcommand{\RMod}[2]{\operatorname{RMod}_{#1}({#2})}
\newcommand{\BiMod}{\operatorname{BiMod}}
\newcommand{\Mod}[2]{\operatorname{Mod}_{#1}({#2})}
\newcommand{\topo}[1]{\mathcal{#1}}
\newcommand{\Shv}{\operatorname{Shv}}
\newcommand{\Psh}{\mathcal P}
\newcommand{\CMon}{\operatorname{CMon}}
\newcommand{\triv}{\operatorname{triv}}
\newcommand{\Act}[2]{{#2}^{#1}}
\newcommand{\lAct}[3]{{#2}^{#1}_{#3}}
\newcommand{\ind}[2]{\operatorname{ind}^{#1}_{#2}}
\newcommand{\res}[2]{\operatorname{res}^{#1}_{#2}}
\newcommand{\CMod}[2]{\operatorname{LMod}_{#1}(\operatorname{Cond}_{#2}(\Sp))}

\newcommand{\cond}[1]{\operatorname{#1-cond}}
\newcommand{\condo}{\operatorname{cond}}
\newcommand{\grostop}{\mathcal T^{LS}_{\lambda}}

\newcommand{\light}{\operatorname{light}}
\newcommand{\solid}{{\square}}
\newcommand{\Cond}[1]{\operatorname{Cond}_{#1}}
\newcommand{\Sol}[1]{\operatorname{Solid}_{#1}}

\newcommand{\cH}[1]{\mathbb H_{#1}}

\newcommand{\icH}[1]{\underline{\mathbb H}_{#1}}
\newcommand{\ccH}{\mathbb H_{\operatorname{cond}}}
\newcommand{\ckH}{\mathbb H_{\operatorname{\kappa-cond}}}

\newcommand{\cckH}[1]{\mathbb H_{\operatorname{#1cond}}}

\newcommand{\icckH}[1]{\underline{\mathbb H}_{\operatorname{#1cond}}}
\newcommand{\ickH}{\underline{\mathbb H}_{\kappa\operatorname{-cond}}}
\newcommand{\iccH}{\underline{\mathbb H}_{\operatorname{cond}}}
\newcommand{\intgrpcoh}[1]{\underline{H}_{#1}}

\newcommand{\sheaf}{\operatorname{sheaf}}
\newcommand{\Zar}{\operatorname{Zar}}
\newcommand{\sing}{\operatorname{\operatorname{sing}}}

\newcommand{\condgrpcoh}[1]{H_{\operatorname{#1cond}}}
\newcommand{\icondgrpcoh}[1]{\underline{H}_{\operatorname{#1cond}}}

\newcommand{\solgrpcoh}[1]{H_{\operatorname{#1sol}}}
\newcommand{\grpcoh}[1]{H_{#1}}

\newcommand{\contgrpcoh}[1]{H_{\operatorname{#1cont}}}
\newcommand{\Top}{\operatorname{Top}}
\newcommand{\cont}{\mathcal C}
\newcommand{\CH}{\operatorname{CH}}

\newcommand{\edCH}[1]{\operatorname{edCH}_{#1}}
\newcommand{\kappacont}{\operatorname{Top}^{\kappa\text{-}c}}
\newcommand{\wt}{\operatorname{wt}}

\renewcommand{\Pr}{\operatorname{Pr}}
\newcommand{\id}{\operatorname{id}}
\newcommand{\algebra}[1]{{#1}}
\newcommand{\sphere}[1]{\mathbb T^{#1}}

\newcommand{\Cont}[1]{\operatorname{Cont}^{#1}}
\newcommand{\To}{\operatorname{T1}}
\newcommand{\TCont}[1]{\operatorname{T1Cont}^{#1}}
\newcommand{\Bun}[3]{\ensuremath{\operatorname{Bun}^{#1}_{#3}(#2)}}
\tikzset{%
    symbol/.style={%
        draw=none,
        every to/.append style={%
            edge node={node [sloped, allow upside down, auto=false]{$#1$}}}
    }
}
\newcommand{\mywidehatshv}{\operatorname{(hyp)Shv}}
\newcommand{\hypershv}{\operatorname{hypShv}}
\newcommand{\hypershvacc}{\operatorname{hypShv}^{\operatorname{acc}}}
\newcommand{\shvacc}{\operatorname{Shv}^{\operatorname{acc}}}
\newcommand{\mywidehatshvacc}[1]{\operatorname{(hyp)Shv}^{\operatorname{acc}}_{#1}}
\newcommand{\mywidehatshvkacc}[2]{\operatorname{(hyp)Shv}^{{#2}\text{-}\operatorname{colim}}_{#1}}
\title{Condensed Group Cohomology}
\author{Emma Brink\footnote{\href{mailto:brink@math.uni-bonn.de}{brink@math.uni-bonn.de}}}
\begin{document}
\maketitle
\begin{abstract} Condensed mathematics as developed by Clausen and Scholze yields a version of derived functors over the category of continuous $G$-modules for a Hausdorff topological group $G$. In this article, we study the resulting notion of group cohomology and its relation to continuous group cohomology and the condensed/sheaf/singular cohomology of classifying spaces. 
While condensed group cohomology is generally a more refined invariant than continuous group cohomology, we show that for a broad class of topological groups, continuous group cohomology with solid coefficients, such as locally profinite continuous $G$-modules, can be realized as a derived functor in the condensed setting.
We also revisit cornerstones of condensed mathematics from \cite{Scholzecondensed}, paying special attention to set-theoretic size issues. 
To this end, we review a framework for working with accessible (hyper)\-sheaves on large sites satisfying suitable accessibility conditions and show that the associated categories retain many topos-like properties. Moreover, we generalize identifications of condensed with sheaf cohomology from \cite{Scholzecondensed}. 
\end{abstract}
\tableofcontents
\section*{Introduction}

\numberwithin{equation}{section}
\addcontentsline{toc}{section}{Introduction}
For a topological group $G$, one would like to implement group cohomology as a derived fixed-point functor over the category of continuous $G$-modules, but since the category of continuous $G$-modules is not abelian, there is no obvious candidate for such a functor. 
One strategy to address this is to embed topological spaces into a topos. This is effective since module categories in a topos are abelian and hence accessible to homological algebra. 
Condensed mathematics as introduced in \cite{Scholzecondensed} works along these lines. 
Up to set-theoretic size issues, condensed sets are a 1-topos and there is a finite products-preserving functor \begin{align*}\underline{(-)}\colon \To\Top\to\Cond{}(\Set)\end{align*} from $\To$ topological spaces to condensed sets which is fully faithful on compactly generated topological spaces. 
For a Hausdorff\footnote{As $\To$ topological groups are Hausdorff, we need to assume that $G$ is Hausdorff to evaluate $\underline{(-)}$ at $G$.} topological group $G$, $\underline{(-)}$ induces a functor from the category of Hausdorff continuous $G$-modules to the category of condensed $\mathbb Z[\underline{G}]$-modules which is fully faithful on continuous $G$-modules whose underlying topological space is compactly generated. The category $\Cond{}(\mathbb Z[\underline{G}])$ of condensed $\mathbb Z[\underline{G}]$-modules is abelian, and we define condensed group cohomology as \[\condgrpcoh{}^*(\underline{G},-)\coloneqq \Ext_{\Cond{}(\mathbb Z[\underline{G}])}(\mathbb Z,-).\footnote{There also exists a derived fixed-point functor $\mathcal D(\Cond{}(\mathbb Z[\underline{G}]))\to\mathcal D(\Cond{}(\Ab))$ which refines condensed group cohomology. We briefly discuss this in \cref{section:condensedgroupcohomology}.}\]
In this article, we study its relation to continuous group cohomology (defined in terms of continuous cochains) and the condensed/sheaf/singular cohomology of classifying spaces.  
\subsection*{Main results}
\addcontentsline{toc}{subsection}{Main results}
For a Hausdorff topological group $G$, the \v{C}ech-to-cohomology spectral sequence (\cref{Bousfieldkanspectralsequencehomology}) for the cover $\underline{G}\to *$ yields a natural transformation \[\contgrpcoh{}^*(G,-)\to \condgrpcoh{}^*(\underline{G},-)\circ \underline{(-)}\] relating continuous group cohomology to condensed group cohomology. This is an isomorphism in two important cases: For locally profinite groups\footnote{We call a group locally profinite if it is locally compact Hausdorff and totally disconnected.} and solid (\cref{definitionsolid}) coefficients, for example locally profinite continuous $G$-modules (\cref{condensedequalscontinuousprofinite}), as well as for locally compact groups and finite-dimensional, continuous real $G$-representations as coefficients (\cref{continuousequalscondensedonvectorspacecoefficients}). 

In general, however, condensed group cohomology is a significantly more refined invariant than continuous group cohomology, as the following result illustrates: 
\begin{thmintro}[\cref{condensedcohomologyiscohomologyofclassifyingspaces}, \cref{condensedgroupcohomologyissheafcohomologyofbg}]\label{groupcohomologycohomologyclassifyingspacesintro1}
    Suppose $G$ is a Hausdorff topological group satisfying one of the following conditions:
    \begin{romanenum}
        \item $G$ is homotopy equivalent to a locally compact Hausdorff space.
        \item $G$ is homotopy equivalent to a locally contractible topological space. 
    \end{romanenum}
    For a discrete abelian group $M$, viewed as continuous $G$-module with trivial $G$-action, condensed group cohomology of $\underline{G}$ with coefficients in $\underline{M}$ is isomorphic to the singular/sheaf cohomology of a classifying space of numerable principal $G$-bundles $BG$,
    \[ H^*_{\operatorname{sing}}(BG,M)\cong \cH{\sheaf}^*(BG,M)\cong \condgrpcoh{}^*(\underline{G}, \underline{M}).\]
\end{thmintro}
By contrast, in the situation of the above theorem, continuous group cohomology only depends on the group of connected components $\pi_0G$ of $G$.

The idea to define group cohomology by embedding topological spaces into a topos goes back to \cite[Expos{\'e} IV, section 2.5]{SGA4}, who worked with the gros topos. Gros topos group cohomology was studied by Flach (\cite{flach}), who established analogues of the above results.   
A main advantage of condensed group cohomology over gros topos group cohomology is that more short exact sequences of continuous $G$-modules induce long exact sequences in group cohomology. There is a natural comparison map from gros topos group cohomology to condensed group cohomology, which is an isomorphism in many cases (\cref{comparisongrostoposgroupcohomologyl}). 

A central feature of the condensed formalism is that there is an extremely well-behaved notion of (non-archimedean) completeness for condensed abelian groups, called \textit{solidity}. We review it in \cref{section:solidmodules}. Importantly, if $M$ is a locally profinite abelian group, then $\underline{M}$, the image of $M$ under the functor $\underline{(-)}$, is solid. 
For a condensed ring $\algebra{R}$, the category of solid $\algebra{R}$-modules $\Sol{}(\algebra{R})$ (\cref{definitionsolid}) is abelian, has enough projectives, and the forget functor $\Sol{}(\algebra{R})\subseteq \Cond{}(\algebra{R})$ has a left adjoint, called \emph{solidification} (\cref{solidenoughprojectiveskappa,solidenoughprojectiveswithoutkappa,underivedsolidification,underivedsolidificationwithoutkappa}).
We show that, despite being very different from condensed group cohomology, continuous group cohomology with solid coefficients can be realized as a derived functor in the condensed setup for a large class of topological groups: 
\begin{thmintro}[\cref{solidequalscontinuousgoodgroups}]\label{solidgroupcohomologyintro} Suppose that $G$ is a Hausdorff topological group which is a finite product of groups of the following types: 
\begin{romanenum}
    \item Groups homotopy equivalent to a coproduct of compact Hausdorff spaces, e.g.\ locally compact abelian groups and locally profinite groups. 
    \item Groups homotopy equivalent to a locally contractible space. 
    \item Groups which are locally connected and locally compact.
\end{romanenum}
Then for all solid continuous $G$-modules $M$, 
\[\contgrpcoh{}^*(G,M)\cong \Ext^*_{\Sol{}(\mathbb Z[\underline{G}])}(\mathbb Z, \underline{M}).\]
\end{thmintro}

It is natural to define the cohomology of a condensed set $X$ as \[\ccH^*(X,-)\coloneqq \Ext^*_{\Cond{}(\Ab)}(\mathbb Z[X],-).\] Building on results of \cite{Scholzecondensed}, we show that condensed cohomology evaluated at the condensed set $\underline{X}$ associated to a $\To$ topological space $X$ via the functor $\underline{(-)}$ recovers sheaf cohomology of $X$ in many cases: 
\begin{thmintro}[\cref{condensedandsheafcohomology1}]\label{cohomologycondensedintro1}
    Suppose that $X$ is a locally compact Hausdorff space and $M$ is a product of a discrete abelian group and a finite-dimensional normed $\mathbb R$-vector space. 
    Then condensed and sheaf cohomology of $X$ with coefficients in $M$ are naturally isomorphic: \[ \cH{\sheaf}^*(X,M)\cong \ccH^*(\underline{X}, \underline{M}).\]
\end{thmintro}
\begin{thmintro}[\cref{condensedandsheafcohomologylocallycontractible}]\label{cohomologycondensedintro2}
    Suppose that $M$ is a Hausdorff topological group such that $\underline{M}$ is \emph{solid} (\cref{definitionsolid}), e.g.\ $M$ is discrete or more generally locally profinite. Denote by $M^{\delta}$ the underlying discrete abelian group.
    Suppose $X$ is a $\To$ topological space which is homotopy equivalent to a locally contractible space, e.g.\ a CW-complex. Then there is a natural isomorphism 
    \[\cH{\sheaf}^*(X,M^{\delta})\cong \ccH^*(\underline{X}, \underline{M}).\]
\end{thmintro}
 These identifications of sheaf with condensed cohomology are the essential ingredient for our comparison of condensed cohomology with continuous cohomology and the cohomology of classifying spaces.

 Condensed sets are sheaves on a large site and therefore not a topos. To obtain a reasonably well-behaved category, \cite{Scholzecondensed} chose to work with the colimit, in the very large category of large categories, over all categories of sheaves on small subsites. In the condensed setting, this is equivalent to restricting to sheaves which are accessible. We discuss this in a more general context in the first chapter and, following \cite{Waterhouse-fpqc-sheafification} and \cite{barwick2019pyknoticobjectsibasic}, formulate accessibility conditions on a large site under which this identification is valid. 
While the resulting categories of accessible (hyper)sheaves are typically not presentable, they retain many topos-like features and for instance satisfy all of Girauds's axioms except accessibility. 
We also show that their categories of spectrum objects behave very similar to the stabilization of a presentable $\infty$-category/$\infty$-topos, which leads to a well-behaved notion of spectrum-valued (group) cohomology. 
In the second chapter, we review cornerstones of condensed mathematics from \cite[Lectures 1-6]{Scholzecondensed} and prove \cref{cohomologycondensedintro1,cohomologycondensedintro2}.
The final chapter studies condensed group cohomology. We first discuss a notion of group cohomology in the general context of \emph{big topoi} (\cref{definitionbigpresentable}) introduced in the first chapter, and then specialize to the condensed setting and prove \cref{groupcohomologycohomologyclassifyingspacesintro1,solidgroupcohomologyintro}.

This article is written in the language of $\infty$-categories as as set down in \cite{higheralgebra}, \cite{highertopostheory}, and so by a category, we will from now on always mean an $(\infty,1)$-category. Analogously, topos will mean $\infty$-topos. We denote by $\an$ the category of animae/spaces/homotopy types. 
\subsection*{Outline}
\addcontentsline{toc}{subsection}{Outline}
  We now give a more detailed overview of the contents of this article and sketch the proofs of our main results. 
\subsubsection*{I Accessible sheaves and Higher Algebra in big presentable categories}
The category underlying a site is usually required to be small to guarantee the existence of a (hyper)sheafification functor. 
However, there are many contexts where one wishes to study sheaves on large categories, and where no general sheafification functor exists, for example the fpqc topology on affine schemes (\cite{Waterhouse-fpqc-sheafification}). 
Following \cite{Waterhouse-fpqc-sheafification}, \cite{barwick2019pyknoticobjectsibasic}, we describe accessibility conditions on an \textit{explicit covering site} (\cref{hyperaccessiblecoveringsitesdefinition}) under which accessible presheaves can be (hyper)sheafified. 
It is natural to work with accessible presheaves over a coaccessible category $\mathcal C$, since these are precisely the small colimits of representables. 

The accessibility conditions ensure that the resulting categories of (hyper)accessible sheaves are, despite typically not presentable, reasonably well-behaved and retain many topos-like properties. For instance, they satisfy all of Giraud's axioms except accessibility. 
Adapting a result of \cite{LucasMannthesis}, we show that for a presentably symmetric monoidal category $\mathcal C$, categories of $\mathcal C$-valued accessible (hyper)sheaves on (hyper)accessible explicit covering sites inherit a closed symmetric monoidal structure (\cref{closednessmonoidalstructureaccessiblesheaves}). 
This all follows from a description of the category of accessible (hyper)sheaves on (hyper)accessible explicit covering sites as a large filtered colimit of topoi along fully faithful, left-exact left adjoints (in the very large category of large categories). 
We call (large) categories arising this way \textit{big topoi} (\cref{definitionbigpresentable}). 
Topoi and categories of accessible (hyper)sheaves on a (hyper)accessible explicit covering site are big topoi. Moreover, if $G$ is a group object in a big topos, the category of $G$-objects is itself a big topos (\cref{Gobjectsisbigtopos}). 
More generally, we call a (possibly large) category \emph{big presentable} if it arises as a (possibly
large) filtered colimit, taken in the category of large categories, of presentable categories along
fully faithful, left-exact left adjoints. This notion is tailored to guarantee that the category of
spectrum objects in a big presentable category behaves like the stabilization of a presentable category in many aspects:
\begin{thmintro} Suppose $\mt{B}_{\infty}$ is a big presentable category. 
\begin{romanenum}
    \item The category $\stab{\mt{B}_{\infty}}$ of spectrum objects in $\mt{B}_{\infty}$ is stable
(\cref{spectrumobjectsinbigcatsstable}), and big presentable (\cref{stabilizationbigpresentablecategoriesbigpresentable}). 
    \item The infinite loop space functor $\Omega^{\infty}\colon \stab{\mt{B}_{\infty}}\to \mt{B_{\infty}}$
admits a left adjoint $\Sigma^{\infty}_{+}\colon \mt{B}_{\infty}\to\stab{\mt{B}_{\infty}},$ which factors into left adjoints
$\mt{B}_{\infty}\to \CMon(\mt{B}_{\infty})\to\CGrp(\mt{B}_{\infty})\to \stab{\mt{B}_{\infty}}$
(\cref{existencesuspension}).
\item There is a $t$-structure on $\stab{\mt{B}_{\infty}}$ whose connective part is generated by the essential image of $\Sigma^{\infty}_{+}$ under small colimits and extensions (\cref{tstructurespectrumobjects}).
\item If $\mt{B}_{\infty}$ is a big topos, then $\Sigma^{\infty}_{+}$ factors over an equivalence $\CGrp(\mt{B}_{\infty})\cong \stab{\mt{B}_{\infty}}_{\geq 0}$. \item 

For a stable, big presentable category $\mt{C}_{\infty}$, pullback along $\Sigma^{\infty}_{+}\colon \mt{B}_{\infty}\to\stab{\mt{B}_{\infty}}$ defines a fully faithful functor 
\[(\Sigma^{\infty}_{+})^*\colon \Fun^{\operatorname{colim}}(\stab{\mt{B}_{\infty}}, \mt{C}_\infty)\hookrightarrow \Fun^{\operatorname{colim}}(\mt{B}_{\infty}, \mt{C}_{\infty}),\] where $\Fun^{\operatorname{colim}}(-,-)$ refers to small colimits preserving functors. If $\mt{B}_{\infty}$ and $\mt{C}_{\infty}$ admit small coproducts, then $(\Sigma^{\infty}_{+})^*$ is an equivalence (\cref{universalpropertyspectrumobjectsbigpresentable}). 
\end{romanenum} 
\end{thmintro}
The stabilization functor $\Sigma^{\infty}_{+}\colon \mt{B}_{\infty}\to \stab{\mt{B}_{\infty}}$ moreover interacts well with (symmetric) monoidal structures. We define big presentably (symmetric) monoidal categories (\cref{presentablymonoidalcat}) and show that the (symmetric) monoidal structure on a big presentably (symmetric) monoidal category $\mt{B}_{\infty}^{\otimes}$ induces a (symmetric) monoidal structure on $\stab{\mt{B}_{\infty}}^{\otimes}$ such that $\Sigma^{\infty}_{+}$ enhances to a (symmetric) monoidal functor $\Sigma^{\infty, \otimes}_{+}$ (\cref{symmetricmonoidalstructureonspectrumobjects}). This in particular applies to big topoi with the cartesian monoidal structure. For categories of accessible (hyper)sheaves, it recovers the closed symmetric monoidal structure induced by the smash product on spectra (\cref{symmetricmonoidalstructureaccessiblesheavesofspectraclosed}). 
We describe a universal property of the (symmetric) monoidal structure on the stabilization (\cref{symmetricmonoidalstructureonspectrumobjects}), and show the following: 
\begin{lemmaintro}[{\cref{universalpropertyspectrumobjectsbigpresentablemonoidal}}]
For a stable, big presentably (symmetric) monoidal category $\mt{C}_{\infty}^{\otimes}$, pullback along $\Sigma^{\infty, \otimes}_{+}\colon \mt{B}_{\infty}^{\otimes}\to\stab{\mt{B}_{\infty}}^{\otimes}$ induces a fully faithful functor 
\[ (\Sigma^{\infty, \otimes}_{+})^*\colon \Fun^{\operatorname{colim}, \otimes}(\stab{\mt{B}_{\infty}}^{\otimes}, \mt{C}_{\infty}^{\otimes})\hookrightarrow \Fun^{\operatorname{colim}, \otimes}(\mt{B}_{\infty}^{\otimes}, \mt{C}_{\infty}^{\otimes}), \] where $\Fun^{\operatorname{colim}, \otimes}(-,-)$ refers to the category of small colimits preserving, (symmetric) monoidal functors.
Its essential image consists of those colimit-preserving, symmetric monoidal functors $\mt{B}_{\infty}^{\otimes}\to \mt{C}_{\infty}^{\otimes}$ which factor over $\mt{B}_{\infty}\xrightarrow{\Sigma^{\infty}_{+}}\stab{\mt{B}_{\infty}}\to\mt{C}_{\infty}$.  
If $\stab{\mt{B}_{\infty}}$ and $\mt{C}_{\infty}$ admit small coproducts, then $(\Sigma^{\infty, \otimes}_{+})^*$ is an equivalence.
\end{lemmaintro}

In \cref{Spmodulestructure}, we describe a natural spectral enrichment of the category of spectrum objects $\stab{\mt{X}}$ in a big presentable category $\topo{X}$. We use the spectrum-valued mapping functors $\map_{\stab{\mt{X}}}(-,-)$ to define cohomology as \[\cH{\topo{\mt{X}}}(-,-)\coloneqq \map_{\stab{\mathcal X}}(\Sigma^{\infty}_{+}-,-)\colon \topo{X}^{\operatorname{op}}\times \stab{\topo{X}}\to \Sp.\]
This notion enjoys the expected properties, including a \v{C}ech-to-cohomology spectral sequence (\cref{Bousfieldkanspectralsequencehomology}) for the cohomology groups $\cH{\topo{X}}^*(-,-)\coloneqq \pi_{-*}\cH{\topo{X}}(-,-)$. Since every adjunction between stable, big presentable categories is spectrally enriched (\cref{adjunctionsspectrallyenriched}), every geometric morphism $L\colon \topo{X}\leftrightarrows \topo{Y}\colon R$ (\cref{definitiongeometricmorphism}) between big topoi induces an equivalence \begin{align}\label{cohomologyinvariantundergeometricmorphismintro}\cH{\topo{X}}(-,R_{\Sp}-)\cong \cH{\topo{Y}}(L-,-), \end{align} where $R_{\Sp}$ denotes the \textit{stabilization} (\cref{definitionstabilization}) of $R$. In particular, cohomology in a big topos $\topo{X}$ can always be computed in a subcategory $\topo{X}_{\lambda}\subseteq \topo{X}$ which is a topos. 

In \cref{section:modulecategories} we record structural properties of module categories in big presentably monoidal categories which we will use in our discussion of group cohomology in big topoi. 
We show that for an algebra $\algebra{R}$ in a big presentably monoidal category $\mathcal C^{\otimes}$, the category of left $\algebra{R}$-modules $\LMod{\algebra{R}}{\mathcal{C}}$ is itself big presentable (\cref{filteredcolimitsmodules}). 
Moreover, if $\mathcal C^{\otimes}$ is big presentably symmetric monoidal, for a commutative algebra $\algebra{R}\in\CAlg(\mathcal C)$, $\LMod{\algebra{R}}{\mathcal C}$ inherits the structure of a big presentably symmetric monoidal category (\cref{symmetricmonoidalstructurebigtoposmodules}).
Under mild conditions on $\mathcal C^{\otimes}$, this monoidal structure is closed:
\begin{lemmaintro}[\cref{internalhommodules}]
    Suppose that $\mathcal C^{\otimes}$ is a potentially large, closed symmetric monoidal category which has all $\Delta^{\operatorname{op}}$-indexed colimits and all $\Delta$-indexed limits. 
    For $\algebra{R}\in \CAlg(\mathcal C)$, the induced symmetric monoidal structure (\cref{symmetricmonoidalstructure}) on $\LMod{\algebra{R}}{\mathcal C}$ is closed. \end{lemmaintro}

For a big presentably monoidal category $\topo{X}$ and an algebra $R\in\Alg(\stab{\mathcal X}_{\geq 0})$, the category of left $R$-modules in $\stab{\mathcal X}$ inherits a $t$-structure (\cref{tstructureonmodules}). 
If $\topo{X}$ is a big topos, the heart of this $t$-structure is equivalent to the 1-category $\LMod{\pi_0\algebra{R}}{\Ab(\tau_{\leq 0}\topo{X})}$ of underived $\pi_0\algebra{R}$-modules in $\Ab(\tau_{\leq 0}\topo{X})$.
(Un)derived module categories in big topoi behave different than in ordinary topoi. For instance, the category of discrete abelian group objects $\Ab(\tau_{\leq 0}\topo{X})$ in a topos $\topo{X}$ is Grothendieck abelian (\cite[Expos{\'e} 2, 6.7]{SGA4}) and in particular has functorial injective embeddings (\cite[\href{https://stacks.math.columbia.edu/tag/079H}{Tag 079H}]{stacks-project}), whereas the category of condensed abelian groups has no non-zero injectives (\cite{noinjectivesincondensed}). 
Adapting \cite[Theorem 2.1.2.2]{SAG}, we describe conditions on a big topos $\topo{X}$ under which \[\mathcal D(\LMod{\pi_0\algebra{R}}{\Ab(\tau_{\leq 0}\topo{X})})\cong \LMod{\pi_0\algebra{R}}{\stab{\topo{X}}}, \] see \cref{derivedcategorymodules}. This applies in particular to $\Cond{(\kappa)}(\an)$: \begin{lemmaintro}[{\cref{condensedmodulesarederivedcats}}]
For  $\algebra{R}\in\Alg(\Cond{(\kappa)}(\Ab))$, 
\[\mathcal D(\LMod{\algebra{R}}{\Cond{(\kappa)}(\Ab)})\cong \LMod{\algebra{R}}{\Cond{(\kappa)}(\Sp)}.\] 
\end{lemmaintro}
Such an equivalence is instructive as we have better control about the categorical properties of the right-hand side. It also allows to identify our definition of cohomology with more classical definitions in terms of Ext-groups: 
If $\topo{X}$ is a big topos satisfying \cref{derivedcategorymodules}, for $X\in\tau_{\leq 0}\topo{X}$ and $A\in\Ab(\tau_{\leq 0}\topo{X})\cong \stab{\topo{X}}^{\heart}$, \[ \cH{\topo{X}}^*(X,A)=\Ext_{\Ab(\tau_{\leq 0}\topo{X})}^*(\mathbb Z[X], \mathbb Z).\] 
\subsubsection*{II Condensed mathematics}
In the second chapter, we review cornerstones of condensed mathematics from \cite[Lectures 1-6]{Scholzecondensed}, paying special attention to set-theoretic size issues and using the language of $\infty$-categories. 
As preparation, we recall the notion of \textit{weight} of topological spaces and related categorical properties of the categories of compact Hausdorff spaces and profinite sets.  
A condensed anima is roughly a hypersheaf on the site of profinite sets with finite, jointly surjective covers. 
There are different ways to deal with the fact that profinite sets are a large category:
\begin{itemize}
\item \cite{barwick2019pyknoticobjectsibasic} fix universes $\mathcal U\in\mathcal V$ and work with hypersheaves of $\mathcal V$-small animae on the category $\Pro(\Fin)_{\mathcal U}$ of $\mathcal U$-small profinite sets. This yields the category of pyknotic animae $\operatorname{Pyk}(\an)_{\mathcal U}^{\mathcal V}$ (with respect to $\mathcal U\in\mathcal V$).
\item More generally, one can consider sheaves on the small category of $\kappa$-light profinite sets (\cref{definitionweight}) for a fixed small, uncountable cardinal $\kappa$, this yields the category of $\kappa$-condensed animae $\Cond{\kappa}(\an)$. \cite{Analyticstacks} work with $\Cond{\aleph_1}(\an)$ for $\aleph_1$ the smallest uncountable regular cardinal, and call this the category of light condensed animae. 
\item \cite{Scholzecondensed} work with accessible, condensed hypersheaves $\Pro(\Fin)^{\operatorname{op}}\to\an$ (\cref{sheafconditionexplicitcoveringsite}, \cref{definitioncondensed}). This yields the category of condensed animae $\Cond{}(\an)$. 
\end{itemize} 
\begin{remintro}Suppose we work internal to a universe $\mathcal U_0$. 
If $\kappa\in\mathcal U_0$ is a small regular strong limit cardinal and $\mathcal U_{\kappa}\subseteq \mathcal U_0$ is the universe associated to ${\kappa}$, i.e. the set of $\kappa$-small sets, then \[\Cond{\kappa}(\an)=\operatorname{Pyk}(\an)_{\mathcal U_0}^{\mathcal U_1}\] is the category of pyknotic animae corresponding to $\mathcal U_{\kappa}\in\mathcal U_0$. If $\mathcal U_0\in\mathcal U_1\in\mathcal U_2$ are further universes, there are fully faithful left adjoints \[\Cond{}(\an)_{\mathcal U_0}\hookrightarrow \operatorname{Pyk}(\an)_{\mathcal U_0}^{\mathcal U_1}\hookrightarrow \Cond{}(\an)_{\mathcal U_1},\] where $\Cond{}(\an)_{\mathcal U_i}$ refers to the categories of condensed animae computed in $\mathcal U_i$, see \cref{pyknoticversuscondensed}. 
\end{remintro}
The category of condensed animae is not presentable and in particular not a topos; however our discussion from the first chapter implies that it retains many topos-like properties, making it a reasonable category to work in. 
More precisely, the category of condensed animae is a big topos and can be expressed as large filtered colimit of the categories $\Cond{\kappa}(\an)$ along fully faithful, left-exact left adjoints, where $\kappa$ runs over all small uncountable regular cardinals, or alternatively all small strong limit cardinals (\cref{condensedascolimit}).
The latter was the original definition of $\Cond{}(\an)$ given in \cite{Scholzecondensed} and is instructive as for strong limit cardinals $\kappa$, $\Cond{\kappa}(\an)\cong \mathcal P_{\Sigma}(\edCH{\kappa})$ enjoys particularly favorable properties (\cref{condensedoncompactextremallydisconnected}). For instance, $\Cond{\kappa}(\Ab)$ has enough projectives for all strong limit cardinals $\kappa$, which implies that the same holds for $\Cond{}(\Ab)$ (\cref{enoughprojectives}). It follows from these colimit-descriptions of $\Cond{}(\an)$ that many computations such as small limits, colimits, and cohomology can be carried out within a topos $\Cond{\kappa}(\an)\subseteq \Cond{}(\an)$ for a sufficiently large regular/strong limit cardinal $\kappa$.
We recall that ($\kappa$-)condensed sets are a good approximation to topological spaces: 
\begin{propositionintro}[{\cite{Scholzecondensed}, \cref{kappacontinuousfullyfaithfullyintocondensed}, \cref{section:condensedsetsandtopologicalspaces}}]
\leavevmode\par
    \begin{romanenum}
    \item For every uncountable cardinal $\kappa$, there is a right adjoint \[\underline{(-)}_{\kappa}\colon \Top\to \Cond{\kappa}(\Set).\] This functor is fully faithful on $\kappa$-compactly generated (\cref{definitionkcontinuous}) spaces. Its left adjoint sends a $\kappa$-condensed set $X$ to a $\kappa$-compactly generated topological space with set of points $X(*)$. 
\item There is a functor \[ \underline{(-)}\colon \To\Top\to\Cond{}(\Set)\] which is fully faithful on compactly generated topological spaces and admits a partially defined left adjoint. 
    \end{romanenum}
\end{propositionintro}
In \cref{section:othermodelscondensedcat}, we apply the functor $\underline{(-)}_{\kappa}$ to describe $\kappa$-condensed animae as hypersheaves on the category of all $\lambda$-small topological spaces (for $\lambda\geq \kappa$), which we apply to generalize comparison results between condensed and sheaf cohomology from \cite{Scholzecondensed}.
\paragraph*{Comparison of condensed with sheaf cohomology and proof of \cref{cohomologycondensedintro1,cohomologycondensedintro2}}
\cite{Scholzecondensed} constructed a comparison map from sheaf to condensed cohomology for compact Hausdorff spaces. Our description of $\kappa$-condensed animae as hyper\-sheaves on the category of all $\lambda$-small topological spaces allows to extend this comparison map to arbitrary topological spaces.
More concretely, for $\lambda\geq \kappa$ we construct a left-exact left adjoint \[ \ladg{j}\colon \grostop\to\Cond{\kappa}(\an)\] from the gros topos of $\lambda$-small topological spaces (\cref{grostoposdefinition}), which sends a $\lambda$-small topological space $X$ to the associated $\kappa$-condensed set $\underline{X}_{\kappa}$. 
For a $\lambda$-small topological space $X$, gros topos cohomology recovers sheaf cohomology (see \cref{grostoposrecoversheafcohomology} for a precise statement). 
Using the geometric morphism $\ladg{j}\dashv \radg{j}$ and (\ref{cohomologyinvariantundergeometricmorphismintro}), we obtain for a $\kappa$-condensed abelian group $A$ a comparison map 
\[\cH{\sheaf}(X,A_{\kappa})\cong \cH{\grostop}(X,\pi_0j_{*,\Sp}A)\to \cH{\grostop}(X,j_{*,\Sp}A)\cong \cckH{}(\underline{X}_{\kappa},A),\] where $A_{\kappa}$ denotes the sheaf \[\Op(X)^{\operatorname{op}}\ni U\mapsto \Hom_{\Cond{\kappa}(\Ab)}(\underline{U}_{\kappa},A).\] This comparison map is by construction natural in the $\lambda$-small topological space $X$. If $X$ is $\To$, for regular cardinals $\kappa>|X|$, $\cckH{}(\underline{X}_{\kappa},A)\cong \ccH{}(\underline{X},A)$ (\cref{condensedcohomologycanbecomputedonfinitestage}). 
We then apply descent of gros topos and condensed cohomology along \textit{local section covers} (\cref{localsectioncover}) and results from \cite{Scholzecondensed} on the condensed cohomology of compact Hausdorff spaces and the solidification of free condensed abelian groups on topological spaces to prove \cref{cohomologycondensedintro1,cohomologycondensedintro2}.
\paragraph*{Solid modules}
In the second half of the second chapter, we review the notion of solidity, the condensed analogue of (non-archimedean) completeness for abelian groups. In preparation for that, we recall basic properties of derived and underived condensed module categories (\cref{section:condensedmodulecategories}) and some compactness properties of condensed abelian groups (\cref{section:compactcondensedabeliangroups}).
For a condensed ring $\algebra{R}$, we define the category of solid $\algebra{R}$-modules $\Sol{(\kappa)}(\algebra{R})\subseteq \Cond{(\kappa)}(\algebra{R})$ as the full subcategory of condensed $\algebra{R}$-modules whose underlying condensed abelian group is solid (\cref{definitionsolid}). 
This enjoys excellent formal properties: 
 \begin{romanenum}
\item 
$\Sol{(\kappa)}(\algebra{R})\subseteq \Cond{(\kappa)}(\algebra{R})$ is an abelian subcategory closed under small limits, colimits and extensions (\cref{solidclosedunderlimitscolimitskappa,underivedsolidificationwithoutkappa}), and has enough projectives (\cref{solidenoughprojectiveskappa,solidenoughprojectiveswithoutkappa}). 
\item The forget functor $\Sol{(\kappa)}(\algebra{R})\to \Cond{(\kappa)}(\algebra{R})$ admits a left adjoint $(-)^{\solid\algebra{R}}$, called solidification. (\cref{underivedsolidification,underivedsolidificationwithoutkappa}). For a commutative condensed ring $\algebra{R}$, $(-)^{\solid\algebra{R}}$ is a symmetric monoidal localisation, i.e. there exists a symmetric monoidal structure on $\Sol{(\kappa)}(\algebra{R})$ such that $(-)^{\solid\algebra{R}}$ enhances to a symmetric monoidal functor (\cref{solidmodulesismodulesinsolid}). 
\item In many cases, the forget functor $\mathcal D(\Sol{(\kappa)}(\algebra{R}))\to\mathcal D(\Cond{(\kappa)}(\algebra{R}))$ admits a left adjoint $(-)^{L\solid\algebra{R}}$, for instance for $\kappa=\aleph_1$, for strong limit cardinals $\kappa$, for $\Cond{}(\algebra{R})$, as well as for $\algebra{R}=\mathbb Z$ and more generally for $\algebra{R}=\mathbb Z[\underline{G}_{(\kappa)}]$ for a Hausdorff topological group $G$ (\cref{derivedsolidificationexistssflatrings,lightringssflat,examplessflatrings,derivedsolificationexistsstronglimitcardinal,derivedsolidificationwithoutkappa,derivedsolidificationexistsabeliangroupskappa}). We call the left-adjoint $(-)^{L\solid\algebra{R}}$ derived solidification. It agrees with the derived functor of solidification at points where the latter is defined (\cref{derivedsolidificationisderivedfunctorofsolidification}).
\item $\mathcal D(\Sol{(\kappa)}(\mathbb Z))\subseteq \mathcal D(\Cond{(\kappa)}(\Ab))$ is a local subcategory and the left adjoint $(-)^{L\solid\mathbb Z}$ is a symmetric monoidal localisation. The induced symmetric monoidal structure on $\mathcal D(\Sol{(\kappa)})$ is closed (\cref{symmetricmonoidalstructureonderivedsolidabeliangroups}).
 \end{romanenum}
To prove these statements, we slightly adapt the arguments from \cite{Scholzecondensed}, who proved the above results for $\Cond{}(\Ab)$ and $\Cond{\kappa}(\Ab)$ for strong limit cardinals $\kappa$. 
We review computations of the derived and underived solidification of free condensed abelian groups on compact Hausdorff spaces and CW-complexes (\cref{solidificationcompacthausdorffspace,solidcw}). These computations imply that condensed cohomology with solid coefficients is homotopy invariant (\cref{solidcohomologyhomotopyinvariant}) and serve as the essential ingredient for the proofs of \cref{cohomologycondensedintro2} and \cref{solidgroupcohomologyintro}. 

We moreover characterize ($\kappa$-)condensed rings for which the forget functor \[\mathcal D(\Sol{(\kappa)}(\algebra{R}))\to\mathcal D(\Sol{(\kappa)}(\mathbb Z))\] factors over an equivalence \[\mathcal D(\Sol{(\kappa)}(\algebra{R}))\cong \LMod{\algebra{R^{\solid\mathbb Z}}}{\mathcal D(\Sol{(\kappa)}(\mathbb Z))}\] (\cref{derivedsolidrmodulesismodulecategory}) and call them ($\kappa$-)$s$-flat. Condensed rings represented by ($\kappa$-light) profinite topological rings (\cref{profiniteringssflat}), as well as group rings $\mathbb Z[\underline{G}_{(\kappa)}]$ for Hausdorff topological groups $G$  are ($\kappa$-)$s$-flat (\cref{examplessflatrings}). All light condensed rings are $\aleph_1$-$s$-flat (\cref{lightringssflat}). 

($\kappa$-)$s$-flatness implies that the forget functor $\mathcal D(\Sol{(\kappa)}(\algebra{R}))\to\mathcal D(\Cond{(\kappa)}(\algebra{R}))$ admits a left adjoint $(-)^{L\solid\algebra{R}}$ (\cref{derivedsolidificationexistssflatrings}). For commutative ($\kappa$-)$s$-flat condensed rings, there is a closed symmetric monoidal structure on $\mathcal D(\Sol{(\kappa)}(\algebra{R}))$ such that $(-)^{L\solid\algebra{R}}$ enhances to a symmetric monoidal functor (\cref{symmetricmonoidalstructuresolidrmodules}). 
\subsubsection*{III Group Cohomology}
In the third chapter, we finally discuss group cohomology. 
We first review the category of $G$-objects in a big topos. This category is itself a big topos (\cref{Gobjectsisbigtopos}), and we define group cohomology as cohomology of its terminal object.\footnote{If $\topo{X}$ is a big topos as in \cref{derivedcategorymodules} and $\stab{\topo{X}}$ has $\Delta^{\operatorname{op}}$-indexed colimits, for a discrete group $G\in\Grp(\tau_{\leq 0}\topo{X})$ and $A\in\LMod{\mathbb Z[G]}{\Ab(\tau_{\leq 0}\topo{X})}$, \[\grpcoh{\topo{X}}^*(G,A)=\Ext^*_{\LMod{\mathbb Z[G]}{\Ab(\tau_{\leq 0}\topo{X})}}(\mathbb Z,A),\] see \cref{groupcohomologyisext}. This in particular applies to $\Cond{(\kappa)}(\an)$.} 
Under mild assumptions on a big topos $\topo{X}$, there exists a \textit{derived fixed-point functor} \[(-)^G\colon \LMod{\mathbb S[G]}{\stab{\topo{X}}}\to \stab{\topo{X}}\] refining group cohomology (\cref{internalgroupcohomologytotalisaton}). We discuss this very briefly before specializing to ($\kappa$-)condensed group cohomology. For the purpose of readability, we focus on condensed group cohomology in this introduction, but we will also prove $\kappa$-condensed analogues of all results below. 
\paragraph*{Comparison to the cohomology of classifying spaces and proof of \cref{groupcohomologycohomologyclassifyingspacesintro1}} 
By \cite{Principalinftybundles}, for a group object $G$ in a big topos $\topo{X}$, the homotopy orbit functor (\cref{deftrivfixedpoints}) 
\[ -//G\colon \Act{G}{\topo{X}}\to \topo{X}\] exhibits the category of $G$-objects in $\topo{X}$ as classifying topos for principal $\infty$-$G$-bundles and factors over an equivalence $\Act{G}{\topo{X}}\cong \oc{\topo{X}}{\mathbb BG}$, where $\mathbb BG\coloneqq *//G$ denotes the classifying space of principal $\infty$-$G$-bundles in $\topo{X}$ (\cref{characterizationbundles}). 
In particular, group cohomology with coefficients in a $G$-module $M$ with trivial $G$-action is isomorphic to the cohomology of $\mathbb BG$ with coefficients in $M$ (\cref{grpcohomologyiscohomologyonclassifyingspaces}). 

For a Hausdorff topological group $G$, every principal $G$-bundle $E\to B$ (with $B$ T1) defines a principal $\infty$-$\underline{G}$-bundle $\underline{E}\to\underline{B}$ in condensed animae which is classified by a map $\underline{B}\to \mathbb BG$ (\cref{bundlesarecondensedbundles}). Pullback along this map yields a comparison map \[ \condgrpcoh{}^*(\underline{G}_{}, \triv-)\cong \ccH^*(\mathbb B\underline{G},-)\to \ccH^{*}(\underline{B},-)\] from condensed group cohomology of $G$ with coefficients in modules with trivial $G$-action to condensed cohomology of $B$. 
If $0\neq G$, the map $\underline{BG}\to\mathbb  BG$ classifying a universal numerable principal $G$-bundle $EG\to BG$ is not an equivalence as $\underline{BG}\in\tau_{\leq 0}\Cond{}(\an)$ whereas $\pi_1\mathbb BG=\underline{G}$. However, homotopy invariance of condensed cohomology with coefficients in solid abelian groups implies the following: 
\begin{thmintro}[\cref{condensedcohomologyiscohomologyofclassifyingspaces}]\label{groupcohomologycohomologyclassifyingspacesintro}
Suppose that $G$ is a Hausdorff topological group and $BG$ is a classifying space for numerable principal $G$-bundles which is $\To$.\footnote{This exists by \cref{classifyingspacet1}.} 
For a solid abelian group $M$, \[\condgrpcoh{}^*(\underline{G},\triv M)\cong \ccH^*(\underline{BG},M).\] 
\end{thmintro} 
For many topological groups, the condensed cohomology of a classifying space $BG$ of numerable principal $G$-bundles with coefficients in a discrete abelian group can be identified with its singular and sheaf cohomology, see \cref{section:Cohomologyclassifyingspaces}. This then implies \cref{groupcohomologycohomologyclassifyingspacesintro1}. 
\paragraph*{Solid group cohomology and proof of \cref{solidgroupcohomologyintro}}
When restricted to solid coefficients, $\Ext^*_{\Sol{}(\mathbb Z[G])}(\mathbb Z,-)$ is a natural alternative to condensed group cohomology. We refer to this as solid group cohomology. As we show in \cref{section:solidmodules}, solid $\mathbb Z[G]$-modules admit plenty projectives, which makes solid group cohomology more computationally accessible than condensed group cohomology. 
Projective solid abelian groups are stable under solid tensor products (\cref{solidprojectivestensorproduct}). This implies the following: 
\begin{lemmaintro}[{\cref{solidequalscontinuousifzgprojective}}]
    Suppose that $G$ is a Hausdorff topological group such that $\mathbb Z[\underline{G}]^{\solid}$ is projective in $\Sol{}$. 
    Then $k$-continuous (\cref{definitioncontinuousgroupcohomology}) group cohomology is isomorphic to solid group cohomology, \[ \contgrpcoh{k\text{-}}^*(G,-)\cong \Ext^*_{\Sol{(\mathbb Z[\underline{G}])}}(\mathbb Z, \underline{(-)}).\]
\end{lemmaintro}
Our computations from \cref{section:solidmodules} moreover show that $\mathbb Z[\underline{G}_{}]^{\solid}$ is projective for large classes of topological groups (\cref{goodgroups}), which then implies \cref{solidgroupcohomologyintro}.

\paragraph*{Notation and conventions}This article is written in the language of $\infty$-categories as set down in \cite{higheralgebra}, \cite{highertopostheory}, and so by a category, we always mean an $(\infty,1)$-category unless stated explicitly otherwise. Analogously, topos means $\infty$-topos. 
We need to assume the existence of Grothendieck universes $\mathcal U_0\in\mathcal U_1\in\mathcal U_2$ to prove structural results on condensed animae and more generally big topoi/big presentable categories. We fix these universes once and for all and refer to categories in $\mathcal U_0$ as small, categories in $\mathcal U_1$ as large and categories in $\mathcal U_2$ as very large. We denote by $\vlCat$ the very large category of large categories. 

We moreover use the following notation: 
\begin{itemize}
\item We write $\an$ for the category of small animae/$\infty$-groupoids/spaces/(weak) homotopy types. 
\item  $\Pr^L$ denotes the category of presentable categories (in $\mathcal U_0$) and left adjoint functors. 
\item Unless stated otherwise, all cardinals are assumed to be small. If $\kappa$ and $\lambda$ are regular cardinals, we write $\lambda\ll\kappa$ if for all cardinals $\lambda_0<\lambda$ and $\kappa_0<\kappa$, $\kappa_0^{\lambda_0}<\kappa$, see the discussion around \cite[Definition 5.4.2.8]{highertopostheory} for examples.
\item For a symmetric monoidal category $\mathcal C^{\otimes}\to N(\Fin_*)$, we denote by $\CAlg(\mathcal C)\coloneqq \Alg_{\mathbb E_{\infty}}(\mathcal C)$ its category of commutative algebra objects. 
\item If $\mathcal C^{\otimes}\to \Assoc^{\otimes}$ is a monoidal category, we denote by $\Alg(\mathcal C)\coloneqq \Alg_{\mathbb A_{\infty}}(\mathcal C)$ the category of associative algebras in $\mathcal C^{\otimes}$. 
\item For a topos $\topo{X}$, we denote by $\operatorname{hyp}\topo{X}\subseteq \topo{X}$ the subcategory on hypercomplete objects (\cite[page 666]{highertopostheory}). 
\item For a presentable category $\mathcal C$ and an uncountable cardinal $\kappa$, we denote by $\Cond{\kappa}(\mathcal C)$ the category of $\kappa$-condensed objects in $\mathcal C$ (\cref{kappacondensedtopology}) and by $\Cond{}(\mathcal C)$ the category of condensed objects (\cref{definitioncondensed}). We use $\Cond{(\kappa)}(\mathcal C)$ as placeholder for any of the two categories. 
\item We denote by \[\underline{(-)}_{\kappa}\colon \Top\to \Cond{\kappa}(\Set)\] and \[\underline{(-)}\colon \To\Top\to\Cond{}(\Set)\] the functors defined in \cref{underlinefunctorontop} and \cref{underlinefunctor}, respectively. We will always view $(\To)$ topological spaces as ($\kappa$-)condensed sets via those functors, and drop the $\underline{(-)}_{(\kappa)}$ occasionally. 
\item We denote by \[\cH{\topo{X}}(-,-)\colon \topo{X}^{\operatorname{op}}\times\stab{\topo{X}}\to\Sp\] the cohomology functor (\cref{definitioncohomologyinatopos}) of a big topos $\topo{X}$, and by $\cH{\topo{X}}^*(-,-)\coloneqq \pi_{-*}\cH{\topo{X}}(-,-)$ the underlying cohomology groups. 
If it exists, internal cohomology in $\topo{X}$ (\cref{definitioncohomologyinatopos}) is denoted \[\icH{\topo{X}}(-,-)\colon \topo{X}^{\operatorname{op}}\times\stab{\topo{X}}\to\stab{\topo{X}}.\]
We denote by \[\cckH{(\kappa\text{-})}(-,-)\coloneqq \cH{\Cond{(\kappa)}(\an)}(-,-)\text{ and }\icckH{(\kappa\text{-})}(-,-)\coloneqq \icckH{\Cond{(\kappa)}(\an)}(-,-)\] the cohomology and internal cohomology in $\Cond{(\kappa)}(\an)$, respectively.    
\item  For a group object $G$ in a big topos $\topo{X}$, $\grpcoh{\topo{X}}(G,-)$ denotes the group cohomology functor of $G$ (\cref{definitiongroupcohomologybigtopos}). If it exists, we denote the internal group cohomology by $\intgrpcoh{\topo{X}}(G,-)$ (\cref{definitiongroupcohomologybigtopos}). For a ($\kappa$-)condensed group $G\in\Grp(\Cond{(\kappa)}(\an)$, we denote by \[\condgrpcoh{(\kappa\text{-})}(G,-)\coloneqq \cH{\Cond{(\kappa)}(\an)}(G,-)\text{ and }\icondgrpcoh{(\kappa\text{-})}(G,-)\coloneqq \icH{\Cond{(\kappa)}(\an)}(G,-)\] its group cohomology and internal group cohomology, respectively. 
\item We call totally disconnected compact Hausdorff spaces profinite. This is justified by \cref{profinitetdch}. Locally profinite means locally compact Hausdorff and totally disconnected. 
\end{itemize}
Our main motivation to work with $\infty$-categories is that it makes the comparison of condensed group cohomology with cohomology of classifying spaces conceptually transparent. 
\paragraph*{Relation to other work}
Our discussion of accessible (hyper)\-sheaves is an adaptation of \cite{barwick2019pyknoticobjectsibasic}.
The second chapter is to a large extent a reformulation of \cite[Lectures 1-6]{Scholzecondensed}, we generalized their results to condensed categories for arbitrary uncountable cutoff cardinals. 
\cite{haine2022descentsheavescompacthausdorff} and \cite{CatrinMairPhDthesis} obtained identifications of condensed with sheaf cohomology similar to \cref{cohomologycondensedintro1,cohomologycondensedintro2} using shape theory. 
For locally profinite groups and solid coefficients, continuous, solid and condensed group cohomology were identified in \cite{Anschuetzsolidhomology}. 
\cite{flach} studied group cohomology in the gros topos and its relation to continuous group cohomology and the sheaf cohomology of classifying spaces. Our proofs are analogous to his, except that we work with spectrum-valued cohomology. 

\paragraph*{Acknowledgements}This article is based on my Master's thesis. I owe special thanks to my Master's thesis advisor Markus Land for introducing me to condensed mathematics, for suggesting the topic of this thesis, for his guidance throughout, and for feedback on an earlier draft. I am very grateful to Ko Aoki, William Balderrama, Jack Davies, Dominik Kirstein, Christian Kremer, Tim Kuppel, Klaus Mattis, Lucas Piessevaux and Stefan Schwede for enlightening discussions about this project, and to Ko Aoki, Jack Davies, Linus Kramer, Klaus Mattis, and Stefan Schwede for valuable feedback on a draft. 
I am supported by the Hausdorff Center of Mathematics at the University of Bonn (DFG GZ 2047/1, project ID 390685813).
\newpage

\numberwithin{equation}{thm}
\section{Accessible sheaves and Higher Algebra in big presentable categories}\label{section:Higheralgebra}
In the definition of sites, the underlying category is usually required to be small to guarantee the existence of a (hyper)sheafification functor and that the associated categories of (hyper)sheaves form a topos. However, there are many situations where this smallness assumption fails and no general sheafification functor exists, for example the fpqc topology on affine schemes, cf. \cite{Waterhouse-fpqc-sheafification}. 
Following \cite{Waterhouse-fpqc-sheafification} and \cite[Section 1.4]{barwick2019pyknoticobjectsibasic}, we formulate accessibility conditions on an \textit{explicit covering site} (\cref{hyperaccessiblecoveringsitesdefinition}) under which accessible presheaves can be sheafified.  
The accessibility conditions ensure that the resulting categories of (hyper)accessible sheaves are, despite typically not presentable, reasonably well-behaved and retain many topos-like properties. For instance, they satisfy all of Giraud's axioms except accessibility. Moreover, their categories of spectrum objects as well as module categories therein behave very similar to those of presentable categories (\cref{section:spectrumobjectsofbigpresentablecategories}, \cref{section:modulecategories}). 
To prove this, it is convenient to work with a class of so-called \emph{big presentable categories} (\cref{definitionbigpresentable}), which encompasses categories of accessible (hyper)sheaves, presentable categories, as well as module categories therein. 
The notion of big presentable category is arranged so that their categories of spectrum objects closely resemble those of presentable categories.  
This chapter is structured as follows: We first discuss accessibility conditions of explicit covering sites and basic properties of the associated categories of (hyper)complete sheaves. 
Then we show that spectrum objects in big presentable categories behave similar as for presentable categories (\cref{section:spectrumobjectsofbigpresentablecategories}). 
In \cref{section:modulecategories}, we record structural properties of module categories, in particular the existence of (closed) symmetric monoidal structures. We also briefly review derived categories. 
There is a natural notion of spectrum-valued cohomology in a \emph{big topos} which we discuss in \cref{section:cohomologyinabigtopos}. 

\subsection{Sheaves on (hyper)accessible explicit covering sites}\label{explicitcoveringssitessection}
We now recall from \cite[section A.3.2]{SAG} a construction of Grothendieck topologies which allows a very explicit characterization of (hyper)complete sheaves. This will lead us to conditions on a site under which the category of accessible (hyper)\-sheaves retains many properties of a topos and accessible presheaves can be (hyper)sheafified. 
\begin{definition}\label{hyperaccessiblecoveringsitesdefinition}{\cite[Definition A.3.12]{LucasMannthesis}, \cite[section A.3.2]{SAG}}
An \emph{explicit covering site} $(\mathcal C,S)$ consists of a (possibly large) category $\mathcal C$ and a wide subcategory  $S\subseteq \mathcal C$ satisfying the following conditions: 
\begin{romanenum}
\item $\mathcal C$ has finite coproducts and pullbacks.
\item Finite coproducts in $\mathcal C$ are universal, i.e.\ given a diagram $\coprod_{i=1}^n c_i\to d\leftarrow e$, the canonical map $\coprod_{i=1}^n (c_i\times_d e)\to (\coprod_{i=1}^n c_i)\times_d e$ is an equivalence.
\item Coproducts in $\mathcal C$ are disjoint, i.e.\ for $c,d\in\mathcal C, c\times_{c\sqcup d} d$ is an initial object of $\mathcal C$. 
\item $S$ contains all equivalences and $f\circ g\in S\Rightarrow f\in S$.
\item $S$ is stable under pullbacks, i.e.\ for $s\colon A\to B\in S$ and $Y\to B\in\mathcal C$, $A\times_B Y\to Y\in S$.
\item $S$ is stable under finite coproducts, i.e.\ for $\{s_i\colon x_i\to y_i\}_{i=1}^n\subseteq S$, \[\sqcup_{i=1}^n s_i\colon \sqcup_{i=1}^n x_i\to \sqcup_{i=1}^n y_i \in  S.\] 
\end{romanenum}
\end{definition} 
If $(\mathcal C,S)$ is a small explicit covering site (i.e.\ $\mathcal C$ is a small category), the sieves $\{ c_i\to c\}_{i\in I}\subseteq \oc{\mathcal C}{c}$ for which there exists a finite subset $F\subseteq I$ with $\sqcup_{i\in F}c_i\to c\in S$ constitute a Grothendieck topology on $\mathcal C$ by \cite[Proposition A.3.2.1]{SAG}. 

\begin{ex}$(\text{Affine schemes}, \text{faithfully flat morphisms})$ is an explicit covering site. The associated site is the fpqc topology on affine schemes. 
\end{ex}(Hyper)\-sheaves in the Grothendieck topology associated to an explicit covering site admit a very explicit description which we recall now.
\begin{definition}[Hypercovers, {\cite[Definition A.5.7.1]{SAG}}]
    \label{Definitionhypercover}
    Suppose that $(\mathcal C,S)$ is an explicit covering site.  
    Denote by $\Delta_{s,+}$ the category with objects the linearly ordered sets \[[n]\coloneqq \{ 0\leq \ldots\leq n\}, n\in \mathbb Z_{\geq -1}\] and morphisms the strictly increasing maps. 
    A functor $U\colon \Delta_{s,+}^{\operatorname{op}}\to\mathcal C$ is called an \emph{augmented semi-simplicial object} in $\mathcal C$. 
    For $n\in\mathbb N_0$, denote by $\Delta_{s,+,< n}\subseteq \Delta_{s,+}$ the full subcategory on linearly ordered sets $[k]$, $k<n$ and by
    \[M_n\colon \Fun(\Delta_{s,+}^{\operatorname{op}}, \mathcal C)\xrightarrow{i^{*}} \Fun(\Delta_{s,+,<n}^{\operatorname{op}}, \mathcal C)\xrightarrow{i_*} \Fun(\Delta_{s,+}^{\operatorname{op}}, \mathcal C)\] the composition of restriction and right Kan extension along $i\colon \Delta_{s,+,< n}\subseteq \Delta_{s,+}$. 
    The right Kan extension exists since $\mathcal C$ has finite limits and ${\Delta_{s,+,<n}}_{/k}$ is finite for all $k\in\Delta_{s,+}$.

    An augmented semi-simplicial object $U_*\colon \Delta_{s,+}^{\operatorname{op}}\to\mathcal C$ is an $S$-\emph{hypercover} if for all $n\in\mathbb N_0$, the unit $U\to M_nU$ evaluates to a morphism $U([n])\to M_nU([n])$ in $S$. 
    Denote by $\Delta_{s}\subseteq \Delta_{s,+}$ the full subcategory on objects $[n], n\in\mathbb N_0$. 
\end{definition}
\begin{definition}[Sheaves on explicit covering sites]\label{sheavesonexplicitcoveringsitesdefinition}
Suppose $(\mathcal C,S)$ is an explicit covering site and $\mathcal D$ is a category. 
\begin{romanenum}
\item The category of $\mathcal D$-valued $S$-sheaves \[ {\Shv}_{S}(\mathcal C, \mathcal D)\subseteq \Fun(\mathcal C^{\operatorname{op}}, \mathcal D)\] is the full subcategory of finite products-preserving functors $\mathcal C^{\operatorname{op}}\to\mathcal D$ such that 
for all $s\colon x\to y\in S$, 
    \[ \Delta_{+}\xrightarrow{\check{C}(s)}\mathcal C^{\operatorname{op}}\xrightarrow{F}\mathcal D\] is a limit diagram, where $\check{C}(s)\colon \Delta_{+}\to\mathcal C^{\operatorname{op}}$ denotes the (opposite of) the \v{C}ech nerve of $s$.  
\item The category of $\mathcal D$-valued hypercomplete $S$-sheaves 
\[ \hypershv_{S}(\mathcal C, \mathcal D)\subseteq \Fun(\mathcal C^{\operatorname{op}}, \mathcal D)\] is the full subcategory of finite products-preserving functors $\mathcal C^{\operatorname{op}}\to\mathcal D$ such that for all $S$-hypercovers $U_*\colon \Delta_{s,+}^{\operatorname{op}}\to\mathcal C$, 
    \[ \Delta_{s,+}\xrightarrow{U_*}\mathcal C^{\operatorname{op}}\xrightarrow{F}\mathcal D\] is a limit diagram.  
\end{romanenum}
\end{definition}
This terminology is compatible with \cite[Notation 6.3.5.16]{highertopostheory}:
\begin{proposition}[{\cite[Proposition A.3.3.1, Proposition A.5.7.2]{SAG}}]\label{sheafconditionexplicitcoveringsite}
    Suppose that $(\mathcal C,S)$ is a small explicit covering site (i.e.\ $\mathcal C$ is small), denote by $\tau_S$ the associated Grothendieck topology on $\mathcal C$ (\cite[Proposition A.3.2.1]{SAG}), and by $\hypershv_{\tau_S}(\mathcal C)\subseteq \Shv_{\tau_{S}}(\mathcal C)$ the topoi of (hypercomplete, \cite[p. 666]{highertopostheory}) $\tau_S$-sheaves. 
    Suppose that $\mathcal D$ is a category with small limits.  
    \begin{romanenum}
    \item Restriction along the sheafified Yoneda embedding $\mathcal C^{\operatorname{op}}\hookrightarrow \mathcal P(\mathcal C)^{\operatorname{op}}\to \Shv_{\tau_S}(\mathcal C)^{\operatorname{op}}$ induces an equivalence 
\[\Fun^{\operatorname{lim}}(\Shv_{\tau_S}(\mathcal C)^{\operatorname{op}}, \mathcal D)\cong \Shv_{S}(\mathcal C, \mathcal D).\] 
\item Restriction along the hypersheafified Yoneda embedding $\mathcal C^{\operatorname{op}}\hookrightarrow\mathcal P(\mathcal C)^{\operatorname{op}}\to \hypershv_{\tau_S}(\mathcal C)^{\operatorname{op}}$ is an equivalence \[\Fun^{\operatorname{lim}}(\hypershv_{\tau_S}(\mathcal C)^{\operatorname{op}}, \mathcal D)\cong \hypershv_{S}(\mathcal C, \mathcal D).\] 
    \end{romanenum}
\end{proposition}

\begin{proof}
    By \cite[Proposition A.3.3.1]{SAG}, a functor $F\colon \mathcal C^{\operatorname{op}}\to \mathcal D$ is in $\Shv_S(\mathcal C, \mathcal D)$ if and only if for all covering sieves $S\subseteq \oc{\mathcal C}{c}$, $F$ exhibits $F(c)$ as limit of the diagram $F|_S$.
    Since $\mathcal D$ has all small colimits and $\mathcal C$ is small, $\Fun(\mathcal C^{\operatorname{op}}, \mathcal D)\cong \Fun^{\operatorname{lim}}(\mathcal P(\mathcal C)^{\operatorname{op}}, \mathcal D)$. 
    As $\Shv_S(\mathcal C)$ is the localization of $\mathcal P(\mathcal C)$ at the maps $S\to h_c$ for covering sieves $S\subseteq \oc{\mathcal C}{c}$, this implies that the above equivalence restricts to an equivalence \[\Fun^{\operatorname{lim}}(\Shv_{\tau_S}(\mathcal C)^{\operatorname{op}}, \mathcal D)\cong \Shv_{S}(\mathcal C, \mathcal D).\] 

    By \cite[Proposition A.5.7.2]{SAG}, a functor $F\colon \mathcal C^{\operatorname{op}}\to\mathcal D$ is in $\hypershv_S(\mathcal C^{\operatorname{op}}, \mathcal D)$ if and only if for all $d\in\mathcal D$, 
    \[\Map_{\mathcal D}(d,F-)\colon\mathcal C^{\operatorname{op}}\to\an\] is a hypercomplete $\tau_S$-sheaf, i.e.\ local with respect to $\infty$-connective morphisms in $\Shv_{\tau_S}(\mathcal C, \an)$. 
    As the functors $\Map_{\mathcal D}(d,-),d\in\mathcal D$ are jointly conservative and $\hypershv_{\tau_S}(\mathcal C)\subseteq \Shv_{\tau_S}(\mathcal C)$ is the localization at the $\infty$-connective morphisms, it follows that the above equivalence restricts to an equivalence \[\Fun^{\operatorname{lim}}(\hypershv_{\tau_S}(\mathcal C)^{\operatorname{op}}, \mathcal D)\cong \hypershv_{S}(\mathcal C, \mathcal D).\qedhere\] 
\end{proof}
\begin{rem}\label{identifysheaveswithtensorproductinprl} 
    If in the situation of the above lemma, $\mathcal D$ is presentable, then \[\Fun^{\operatorname{lim}}(\mywidehatshv_{\tau_S}(\mathcal C, \mathcal D)\cong \mywidehatshv_S(\mathcal C, \an)\otimes_{\Pr^L}\mathcal D\] is the tensor product of presentable categories (\cite[Proposition 4.8.1.15]{higheralgebra}):
    By \cite[Theorem 6.1.0.6]{highertopostheory}, the topos $\mywidehatshv_S(\mathcal C, \an)$ is presentable and by \cite[Proposition 4.8.1.17]{higheralgebra} and \cite[Proposition 5.2.6.2, Remark 5.5.2.10]{highertopostheory}, for presentable categories $\mathcal X, \mathcal D$, \[\mathcal X\otimes_{\Pr^L}\mathcal D\cong \Fun^R(\mathcal D^{\operatorname{op}}, \mathcal X)\cong \Fun^L(\mathcal X, \mathcal D^{\operatorname{op}})^{\operatorname{op}}\cong \Fun^{\operatorname{colim}}(\mathcal X, \mathcal D^{\operatorname{op}})^{\operatorname{op}}\cong \Fun^{\operatorname{lim}}(\mathcal X^{\operatorname{op}}, \mathcal D).\] 
\end{rem}
The above description of (hyper)\-sheaves on explicit covering sites leads to conditions on an explicit covering site $(\mathcal C,S)$ under which its category of accessible (hyper)\-sheaves retains many properties of a topos and accessible presheaves can be (hyper)sheafified:
\begin{definition}\label{accessiblecoveringsitedefinition}Suppose $\kappa$ is a regular cardinal.\footnote{A cardinal $\kappa$ is regular if $\kappa$-small colimits of $\kappa$-small sets are $\kappa$-small. A set $S$ is $\kappa$-small if it has cardinality less than $\kappa$.} 
    \begin{romanenum}
\item An explicit covering site $(\mathcal C,S)$ $\kappa$-\emph{accessible} if $\mathcal C$ is $\kappa$-coaccessible and for all $s\colon X\to Y\in S$, there exists a $\kappa$-cofiltered diagram \[s_*\colon I\to \Fun(\Delta^1, \mathcal C_{\kappa})\cap S\] such that $s=\clim{i\in I}s_i$, where $\mathcal C_{\kappa}\subseteq \mathcal C$ denotes the full subcategory of $\kappa$-cocompact objects.\footnote{As $\mathcal C$ is $\kappa$-coaccessible, it has $\kappa$-filtered limits. In particular, for every small category $K$, $\Fun(K, \mathcal C)$ has $\kappa$-filtered limits and they can be computed pointwise.}

An explicit covering site is \emph{accessible} if for all small regular cardinals $\kappa$, there exists a small regular cardinal $\lambda\geq \kappa$ such that $(\mathcal C,S)$ is $\lambda$-accessible. 
If $(\mathcal C,S)$ is an accessible explicit covering site, denote by 
\[ \Shv_{S}^{\operatorname{acc}}(\mathcal C, \mathcal D)\coloneqq \Shv_S(\mathcal C, \mathcal D)\cap\Fun^{\operatorname{acc}}(\mathcal C^{\operatorname{op}}, \mathcal D)\subseteq \Fun(\mathcal C^{\operatorname{op}}, \mathcal D)\] the full subcategory on accessible $S$-sheaves $F\colon \mathcal C^{\operatorname{op}}\to\mathcal D$. 

\item An explicit covering site $(\mathcal C,S)$ is $\kappa$-\emph{hyperaccessible} if $\mathcal C$ is $\kappa$-coaccessible and for all $S$-hy\-per\-co\-ve\-rings $X_*\colon \Delta^{\operatorname{op}}_{+}\to \mathcal C$, there exists a $\kappa$-cofiltered diagram of $S$-hy\-per\-co\-ve\-rings \[s_*\colon I\to \Fun(\Delta_{+}^{\operatorname{op}}, \mathcal C_{\kappa})\] with $\clim{i\in I}(s_i\colon X_*^i\to X^i)=(X_*\to X)$. Again, $\mathcal C_{\kappa}\subseteq \mathcal C$ denotes the category of $\kappa$-cocompact objects. 

An explicit covering site is \emph{hyperaccessible} if for all small regular cardinals $\kappa$ there exists a small regular cardinal $\lambda\geq \kappa$ such that $(\mathcal C,S)$ is $\lambda$-hyperaccessible. 
If $(\mathcal C,S)$ is a hyperaccessible explicit covering site, denote by \[ \hypershv_{S}^{\operatorname{acc}}(\mathcal C, \mathcal D)\coloneqq \hypershv_S(\mathcal C, \mathcal D)\cap\Fun^{\operatorname{acc}}(\mathcal C^{\operatorname{op}}, \mathcal D)\subseteq \Fun(\mathcal C^{\operatorname{op}}, \mathcal D)\] the full subcategory on accessible, hypercomplete $S$-sheaves.
    \end{romanenum} 
For a (hyper)accessible explicit covering site $(\mathcal C,S)$, we denote by \[\mywidehatshvacc{S}(\mathcal C)\coloneqq \mywidehatshvacc{S}(\mathcal C, \an)\] the category of $\an$-valued accessible (hyper)sheaves. 
\end{definition}
\begin{rems}
\begin{enumerate}
\item It is natural to consider accessible presheaves on a coaccessible, locally small category since this is its free cocompletion under small colimits, see \ \cite[Proposition A.2]{HesselholtPstragowski2024}.  
More generally, if  $\mathcal C$ is small and $i\colon \mathcal C\to \mathcal C^{\operatorname{idem}}$ is an idempotent completion (\cite[Proposition 5.4.2.18]{highertopostheory}), then $\mathcal C^{\operatorname{idem}}$ is small (\cite[Lemma 5.4.2.4]{highertopostheory}), and left Kan extension along $i$ is an equivalence \[\mathcal P(\mathcal C)\cong \mathcal P^{\operatorname{acc}}(\mathcal C^{\operatorname{idem}}).\] In particular, every presheaf category arises as category of accessible presheaves. If $\mathcal C$ is coaccessible and small, then every functor $\mathcal C^{\operatorname{op}}\to \an$ is accessible by \cite[Corollary 5.4.3.6]{highertopostheory}.
\item (Hyper)accessibility of an explicit covering site $(\mathcal C,S)$ implies that 
\begin{align}\label{generalcolimitdescriptionaccessiblesheaves}\colim{\mathcal C^{'}\subseteq C}\mywidehatshv_{S}(\mathcal C^{'})\cong\mywidehatshvacc{S}(\mathcal C,\mathcal D),\end{align} where the colimit on the left is over the large poset of all small subcategories of $\mathcal C$ which are closed under finite limits, and the transition maps are the left adjoints of restriction, see \cref{accessiblesheavesisbigtoposrem} below. The colimit is computed in $\vlCat$.  
In particular, a (hyper)sheaf on a (hyper)accessible explicit covering site is accessible if and only if it is left Kan extended from a sheaf on a small subcategory. 
\item (Hyper)accessibility of an explicit covering site $(\mathcal C,S)$ ensures that there exists a (hyper)\-sheafi\-fi\-ca\-tion functor for accessible presheaves, i.e.\ 
\[\mywidehatshvacc{S}(\mathcal C)\subseteq \Fun^{\operatorname{acc}}(\mathcal C^{\operatorname{op}}, \an)\] has a left-exact left adjoint, see \cref{accessiblesheafification} below. This follows from the colimit-description \ref{generalcolimitdescriptionaccessiblesheaves}.  
\item Note that in general, (hyper)accessibility is a stronger condition than being $\kappa$-(hyper)ac\-ces\-si\-ble for some cardinal $\kappa$. However, if $(\mathcal C,S)$ is an explicit covering site such that the collection of $S$-(hyper)covers is closed under $\kappa$-filtered limits, then $\kappa$-(hyper)accessibility implies (hyper)accessibility, see \cref{kappaccessibilityimpliesaccessibilityifSclosedunderfilteredlimits} below. 

\item \cite{barwick2019pyknoticobjectsibasic} call $\kappa$-accessible explicit covering sites $\kappa$-accessible presites. They define accessible explicit covering sites as explicit covering sites which are accessible for some regular cardinal $\kappa$, but our stronger condition seems to be necessary to obtain the colimit-description of accessible sheaves (\cref{accessiblesheavesisbigtopos}), which ensures that this retains many properties of a topos and accessible presheaves can be sheafified.  
\end{enumerate}
\end{rems}
\begin{ex}
The explicit covering site $(\text{Affine schemes}, \text{faithfully flat morphisms})$ is $\kappa$-accessible for every uncountable cardinal by \cite[Theorem 3.3]{Waterhouse-fpqc-sheafification}.  
\end{ex}

\begin{notation}
If $(\mathcal C,S)$ is a $\kappa$-(hyper)ac\-ces\-si\-ble explicit covering site, denote by $\mathcal C_{\kappa}\subseteq \mathcal C$ the full subcategory on $\kappa$-cocompact objects. Then $(\mathcal C_{\kappa}, S\cap \Fun(\Delta^1, \mathcal C_{\kappa}))$ is an explicit covering site. 
For a category $\mathcal D$ denote by $\mywidehatshv_S(\mathcal C_{\kappa}, \mathcal D)$ the associated category of $\mathcal D$-valued (hypercomplete) sheaves on $\mathcal C_{\kappa}$. \end{notation}
The definition of $\kappa$-(hyper)accessibility is tailored so that the following holds:
\begin{thm}[{\cite[Proposition 1.4.3]{barwick2019pyknoticobjectsibasic}, \cite[Lemma A.2.8]{LucasMannthesis}}]\label{kappacondensedkappaacc}
    Suppose that $(\mathcal C,S)$ is a $\kappa$-accessible explicit covering site and $\mathcal D$ is a category with small limits and colimits and $\Delta$-indexed limits commute with $\kappa$-filtered colimits in $\mathcal D$. 
    Denote by $\mathcal C_{\kappa}\subseteq \mathcal C$ the full subcategory of $\kappa$-cocompact objects and by $\Fun^{\kappa\text{-}\colim{}}(\mathcal C^{\operatorname{op}}, \mathcal D)\subseteq \Fun(\mathcal C^{\operatorname{op}}, \mathcal D)$ the full subcategory on functors which preserve $\kappa$-filtered colimits. 
\begin{romanenum}
    \item The restriction $r^{\kappa}\colon \Fun^{\operatorname{acc}}(\mathcal C^{\operatorname{op}}, \mathcal D)\to \Fun(\mathcal C_{\kappa}^{\operatorname{op}}, \mathcal D)$ has a fully faithful left adjoint $i_{\kappa}$ which factors over an equivalence 
    \[ \Fun(\mathcal C_{\kappa}^{\operatorname{op}}, \mathcal D)\cong \Fun^{\kappa\text{-}\operatorname{colim}}(\mathcal C^{\operatorname{op}}, \mathcal D)\subseteq \Fun^{\operatorname{acc}}(\mathcal C^{\operatorname{op}}, \mathcal D).\] 
    \item This restricts to an adjunction \[ r^{\kappa}\colon \Shv_S^{\operatorname{acc}}(\mathcal C, \mathcal D)\leftrightarrows \Shv_S(\mathcal C_{\kappa}, \mathcal D)\colon i_{\kappa}, \] and $i_{\kappa}$ factors over an equivalence \[ \Shv_{S}(\mathcal C_{\kappa},\mathcal D)\cong \Fun^{\kappa\text{-}\colim{}}(\mathcal C^{\operatorname{op}}, \mathcal D)\cap \Shv_{S}(\mathcal C,\mathcal D).\]

    \item If $(\mathcal C,S)$ is $\kappa$-hyperaccessible, then the adjunction $r^{\kappa}\vdash i_{\kappa}$ restricts to an adjoint pair \[ r^{\kappa}\colon \hypershv_S^{\operatorname{acc}}(\mathcal C, \mathcal D)\leftrightarrows \hypershv_S(\mathcal C_{\kappa}, \mathcal D)\colon i_{\kappa},\] and $i_{\kappa}$ factors over an equivalence \[ \hypershv_{S}(\mathcal C_{\kappa},\mathcal D)\cong \Fun^{\kappa\text{-}\colim{}}(\mathcal C^{\operatorname{op}}, \mathcal D)\cap \hypershv_{S}(\mathcal C,\mathcal D).\]
\end{romanenum}
\end{thm}
\begin{rem}
If $\mathcal D$ is a $\mu$-accessible category, then $\mu$-small limits commute with $\mu$-filtered colimits in $\mathcal D$, see e.g.\ \cite[Lemma A.2.8]{LucasMannthesis}. 
In particular, for uncountable cardinals $\mu\leq \kappa$, $\mu$-accessible categories satisfy the condition of the above theorem.    
\end{rem}
\begin{proof}
    This is a very straightforward generalization of \cite[Proposition 1.4.3]{barwick2019pyknoticobjectsibasic} or \cite[Proposition 2.1.8]{LucasMannthesis}.
    We first explain that restriction has a left adjoint. 
    As $\mathcal C$ is $\kappa$-coaccessible, for all $c\in\mathcal C$, ${c\backslash}\mathcal C_{\kappa}$ is essentially small (\cite[Proposition 5.4.2.2]{highertopostheory}). 
    Since $\mathcal D$ has all small colimits, this implies that for all $F\colon \mathcal C_{\kappa}^{\operatorname{op}}\to\mathcal D$, the left Kan extension along $\mathcal C_{\kappa}\subseteq \mathcal C$ exists, i.e.\ the restriction \[ \Fun(\mathcal C^{\operatorname{op}}, \mathcal D)\to \Fun(\mathcal C^{\operatorname{op}}_{\kappa}, \mathcal D)\] has a left adjoint $i_{\kappa}$. This is fully faithful. 
    As $\mathcal C$ is $\kappa$-coaccessible, $\mathcal C^{\operatorname{op}}\cong \Ind_{\kappa}(\mathcal C_{\kappa}^{\operatorname{op}})$ (\cite[Proposition 5.4.2.2]{highertopostheory}), whence the essential image of the left Kan extension \[i_{\kappa}\colon \Fun(\mathcal C_{\kappa}^{\operatorname{op}}, \mathcal D)\to \Fun(\mathcal C^{\operatorname{op}}, \mathcal D)\] consists precisely of the $\kappa$-filtered colimits preserving functors $\mathcal C^{\operatorname{op}}\to\mathcal D$.
    
    The restriction $r^{\kappa}$ obviously restricts to a functor $\Shv_S(\mathcal C, \mathcal D)\to \Shv_S(\mathcal C_{\kappa}, \mathcal D)$. We claim that its left adjoint $i_{\kappa}$ also preserves $S$-sheaves. Suppose that $F\in\Shv_S(\mathcal C_{\kappa}, \mathcal D)$ and $s\colon X\to Y$ is a morphism in $S$. By $\kappa$-accessibility of $(\mathcal C,S)$ there exists a $\kappa$-filtered diagram $s_*\colon I\to\Fun(\Delta^1, \mathcal C_{\kappa})$ with $s=\clim{i\in I}(s_i\colon X_i\to Y_i)$. 
    As $I$ is filtered, $\check{C}(s)\cong \clim{i\in I}\check{C}(s_i)$, and since finite limits of $\kappa$-cocompact objects are $\kappa$-cocompact, for all $i\in I$ and $[n]\in\Delta^{\operatorname{op}}$, $\check{C}(s_i)([n])\in \mathcal C_{\kappa}$. 
    As $i_{\kappa}F\colon\mathcal C^{\operatorname{op}}\to\mathcal D$ commutes with $\kappa$-filtered colimits, this implies that \[(i_{\kappa}F)(\check{C}(s))\cong \colim{i\in I}\left(i_{\kappa}F(\check{C}(s_i))\right)\cong \colim{i\in I}F(\check{C}(s_i))\] and $i_{\kappa}F(y)\cong \colim{i\in I}F(y_i)$. Since $F$ is an $S$-sheaf, $\clim{\Delta}F(\check{C}(s_*))\cong F(y_*)\in \Fun(I, \mathcal D)$. 
    As $I$ is $\kappa$-filtered and $\Delta$-indexed limits commute with $\kappa$-filtered colimits in $\mathcal C$, it now follows that \[F(y)\cong \colim{i\in I}F(y_i)\cong \colim{i\in I}\clim{\Delta}F(\check{C}(s_i))\cong \clim{\Delta}\colim{i\in I}F(\check{C}(s_i))\cong \clim{\Delta}F(\check{C}(s)), \] which shows that $i_{\kappa}F$ is an $S$-sheaf. 
    This proves that $r^{\kappa}\vdash i_{\kappa}$ restricts to an adjoint pair \[ r^{\kappa}\colon \Shv_S^{\operatorname{acc}}(\mathcal C, \mathcal D)\leftrightarrows \Shv_S(\mathcal C_{\kappa}, \mathcal D)\colon i_{\kappa}.\] 
    By the above and fully faithfulness of $i_{\kappa}$, a functor $F\colon\mathcal C_{\kappa}^{\operatorname{op}}\to\mathcal D$ is an $S$-sheaf if and only if its left Kan extension $i_{\kappa}F\colon \mathcal C^{\operatorname{op}}\to\mathcal D$ is. It follows that $i_{\kappa}$ factors over an equivalence \[ \Shv_{S}(\mathcal C_{\kappa},\mathcal D)\cong \Fun^{\kappa\text{-}\colim{}}(\mathcal C^{\operatorname{op}}, \mathcal D)\cap \Shv_{S}(\mathcal C,\mathcal D).\]

    The same arguments show that if $(\mathcal C,S)$ is $\kappa$-hyperaccessible, restriction and left Kan extension restrict to an  adjunction \[ r^{\kappa}\colon \hypershv_S^{\operatorname{acc}}(\mathcal C, \mathcal D)\leftrightarrows \hypershv_S(\mathcal C_{\kappa}, \mathcal D)\colon i_{\kappa}\] and $i_{\kappa}$ factors over an equivalence \[ \hypershv_{S}(\mathcal C_{\kappa},\mathcal D)\cong \Fun^{\kappa\text{-}\colim{}}(\mathcal C^{\operatorname{op}}, \mathcal D)\cap \hypershv_{S}(\mathcal C,\mathcal D).\qedhere\]
\end{proof}
\begin{cor}\label{fullyfaithfulnessleftkanextension}
    Suppose that $\lambda\geq \kappa$ are regular cardinals and $(\mathcal C,S)$ is an explicit covering site which is $\kappa$- and $\lambda$-(hyper)ac\-ces\-si\-ble. 
    Suppose that $\mathcal D$ is a category with small limits and colimits, and $\Delta$-indexed limits commute with $\kappa$-filtered colimits in $\mathcal D$. Then the restriction \[\mywidehatshv_S(\mathcal C_{\lambda}, \mathcal D)\to \mywidehatshv_S(\mathcal C_{\kappa}, \mathcal D)\] has a fully faithful left adjoint.  
\end{cor}
\begin{proof}
    By \cref{kappacondensedkappaacc}, for $\tau=\kappa, \lambda$, \[\mywidehatshv_S(\mathcal C_{\tau}, \mathcal D)\cong \mywidehatshvkacc{S}{\tau}(\mathcal C, \mathcal D)\] via the left adjoint $i_{\tau}$ of the restriction. 
    This implies that $i_{\kappa}$ factors over a fully faithful functor 
    $\mywidehatshv_S(\mathcal C_{\kappa}, \mathcal D)\to \mywidehatshv_S(\mathcal C_{\lambda}, \mathcal D)$ which is left adjoint to the restriction. 
\end{proof}
The proof of \cref{kappacondensedkappaacc} shows that $\kappa$-and $\lambda$-(hyper)accessibility implies that left Kan extension along $\mathcal C_{\kappa}\subseteq\mathcal C_{\lambda}$ restricts to a functor 
    $\mywidehatshv_S(\mathcal C_{\kappa}, \mathcal D)\to \mywidehatshv_S(\mathcal C_{\lambda}, \mathcal D)$.

We now apply \cref{fullyfaithfulnessleftkanextension} to prove basic categorical properties of categories of accessible (hyper)\-sheaves on (hyper)accessible explicit covering sites. 
If $(\mathcal C,S)$ is a (hyper)ac\-ces\-si\-ble explicit covering site and $\mathcal C$ is not small, then $\mywidehatshv_S(\mathcal C)$ is typically not accessible: 
\begin{cor}\label{accessiblesheavesnottopos}
Suppose that $(\mathcal C,S)$ is a (hyper)ac\-ces\-si\-ble 
covering site and $S$ is subcanonical, i.e.\ $\Map_{\mathcal C}(-,c)$ is an $S$-(hyper)sheaf for all $c\in\mathcal C$. 
Then $\mywidehatshvacc{S}(\mathcal C)$ is accessible if and only if $\mathcal C$ is essentially small.  
\end{cor}
\begin{proof}
If $\mathcal C$ is small and $\tilde{\mathcal C}\subseteq\mathcal C$ is a small subcategory such that the inclusion is an equivalence, then $(\tilde C,\tilde C\cap S)$ is an explicit covering site and \[ \mywidehatshvacc{S}(\mathcal C)\cong \mywidehatshvacc{S\cap \tilde{\mathcal C}}(\tilde{\mathcal C})\] via restriction. By \cite[Proposition A.2]{HesselholtPstragowski2024} and \cref{sheafconditionexplicitcoveringsite}, \[\mywidehatshvacc{S\cap \tilde{\mathcal C}}(\tilde{\mathcal C})\cong \mywidehatshv_{S\cap \tilde{\mathcal C}}(\tilde{\mathcal C})\cong \mywidehatshv_{\tau_{S\cap \tilde{\mathcal C}}}(\tilde{\mathcal C})\] is a topos and in particular accessible (\cite[Proposition 6.2.2.7]{highertopostheory}).

Conversely, if $\mywidehatshv_S(\mathcal C)$ is accessible, there exists a small subcategory $G\subseteq \mywidehatshv_{S}(\mathcal C)$ which generates $\Shv_S(\mathcal C)$ under small colimits.  
Choose a regular cardinal $\kappa$ such that $\mathcal C$ is $\kappa$-coaccessible and all elements in $G$ preserve $\kappa$-filtered colimits. 
\cref{kappacondensedkappaacc} implies that the restriction 
\[\mywidehatshvacc{S}(\mathcal C)\to \mywidehatshv_S(\mathcal C_{\kappa})\] is conservative. 
As the restriction has a fully faithful left adjoint (\cref{kappacondensedkappaacc}), this implies that the restriction is an equivalence. (It follows from the triangle identities that the unit and counit are equivalences). 
Hence by \cref{kappacondensedkappaacc}, every element in $\mywidehatshvacc{S}(\mathcal C)$ is $\kappa$-accessible.
In particular, since $S$ is subcanonical, $\Map_{\mathcal C}(-,c)$ is $\kappa$-accessible for all $c\in\mathcal C$, i.e.\ $\mathcal C=\mathcal C_{\kappa}$. The category $\mathcal C_{\kappa}$ is essentially small by \cite[Proposition 5.4.2.2]{highertopostheory}. 
\end{proof}

\cref{kappacondensedkappaacc} moreover implies: 
\begin{cor}\label{limitscolimitscanbecomputedfinitestagesheaves}
Suppose that $\mathcal D$ is a $\mu$-accessible category with small limits and colimits and $\kappa\geq\mu$ is a regular cardinal. 
\begin{romanenum}
\item If $(\mathcal C,S)$ is a $\kappa$-(hyper)ac\-ces\-si\-ble explicit covering site, then \[\mywidehatshvkacc{S}{\kappa}(\mathcal C, \mathcal D)\subseteq\mywidehatshvacc{S}(\mathcal C, \mathcal D)\] is closed under $\mu$-small limits and small colimits. 
\item If $(\mathcal C,S)$ is (hyper)ac\-ces\-si\-ble, then $\mywidehatshvacc{S}(\mathcal C, \mathcal D)$ has all small limits and colimits. 
\end{romanenum}
\end{cor}
\begin{proof}
    As $\mathcal D$ has all small limits, small limits in $\Fun(\mathcal C^{\operatorname{op}}, \mathcal D)$ are computed pointwise (\cite[Corollary 5.1.2.3]{highertopostheory}) and hence $\mywidehatshv_S(\mathcal C, \mathcal D)\subseteq \Fun(\mathcal C^{\operatorname{op}}, \mathcal D)$ is closed under small limits. Since $\mathcal D$ is $\mu$-accessible, $\mu$-filtered colimits commute with $\mu$-small limits in $\mathcal C$ (see e.g.\ \cite[Lemma A.2.8]{LucasMannthesis}). This implies that for $\kappa\geq \mu$, \[i_{\kappa}\colon \mywidehatshvkacc{S}{\kappa}(\mathcal C, \mathcal D)\subseteq\mywidehatshvacc{S}(\mathcal C, \mathcal D)\] preserves $\mu$-small limits.
    Being a left adjoint, $i_{\kappa}$ also preserves small colimits. 

    Suppose now that $(\mathcal C,S)$ is (hyper)ac\-ces\-si\-ble and $F\colon I\to \mywidehatshv_S(\mathcal C)$ is a small diagram. Choose an uncountable regular cardinal $\mu_F\geq \mu$ such that $(\mathcal C,S)$ is $\mu_F$-(hyper)accessible and \[F(i)\in \mywidehatshvkacc{S}{\mu_F}(\mathcal C, \mathcal D)\] for all $i\in I$. Choose a regular cardinal $\tilde{\mu}\gg\mu$ with $\tilde{\mu}> \max\{ |I|, \mu_F\}$ and a regular cardinal $\tilde{\kappa}\geq \tilde{\mu}$ such that $(\mathcal C,S)$ is $\tilde{\kappa}$-accessible. By \cite[Lemma 5.4.2.10]{highertopostheory}, $\mathcal D$ is $\tilde{\mu}$-accessible, whence \[\mywidehatshvkacc{S}{\tilde{\kappa}}(\mathcal C, \mathcal D)\subseteq\mywidehatshvacc{S}(\mathcal C, \mathcal D)\] is closed under $\tilde{\mu}$-small limits and small colimits by the above. In particular, $F$ admits a limit and colimit and both can be computed in $\mywidehatshvkacc{S}{\tilde{\kappa}}(\mathcal C, \mathcal D)$. 
\end{proof}

Suppose that $(\mathcal C,S)$ is a (hyper)ac\-ces\-si\-ble explicit covering site. For all $\lambda\in\Lambda_{\mathcal C}$, $\mathcal C_{\lambda}$ is essentially small since $\mathcal C$ is $\lambda$-coaccessible (\cite[Proposition 5.4.2.2]{highertopostheory}). \cref{kappacondensedkappaacc} therefore implies:  
\begin{cor}\label{kappaccessibleistopos}
If $(\mathcal C,S)$ is a $\kappa$-(hyper)ac\-ces\-si\-ble explicit covering site, then $\mywidehatshvkacc{S}{\kappa}(\mathcal C)$ is a topos. 
\end{cor}
\begin{proof}
By \cref{kappacondensedkappaacc}, \[\mywidehatshvkacc{S}{\kappa}(\mathcal C)\cong\mywidehatshv_S(\mathcal C_{\kappa}).\] 
As $\mathcal C$ is $\kappa$-coaccessible, $\mathcal C_{\kappa}$ is essentially small (\cite[Proposition 5.4.2.2]{highertopostheory}). Choose a small category $\tilde{\mathcal C}_{\kappa}$ with an equivalence $f\colon \tilde{\mathcal C}_{\kappa}\cong \mathcal C_{\kappa}$ and denote by $f^*S\subseteq \Fun(\Delta^1, \tilde{\mathcal C}_{\kappa})$ the full subcategory on morphisms $s$ with $f(s)\in S$.
It is straightforward to check that $(\tilde{\mathcal C}_{\kappa}, f^*S)$ is an explicit covering site and that the pullback $f^*\colon\Fun(\mathcal C_{\kappa}^{\operatorname{op}}, \an)\cong \Fun(\tilde{\mathcal C}_{\kappa}^{\operatorname{op}}, \an)$ restricts to an equivalence \[\mywidehatshv_S(\mathcal C_{\kappa})\cong \mywidehatshv_{f^*S}(\tilde{\mathcal C}_{\kappa}).\]
By \cref{sheafconditionexplicitcoveringsite} and \cite[Propositions 6.2.2.7, 6.5.1.16, 6.2.1.1]{highertopostheory}), $\mywidehatshv_{f^*S}(\tilde{\mathcal C}_{\kappa})$ is a topos.
\end{proof}
Together with \cref{limitscolimitscanbecomputedfinitestagesheaves}, this implies that small limits and colimits in $\mywidehatshvacc{S}(\mathcal C)$ behave like in a topos:
\begin{cor}[Giraud's axioms]\label{Giraudsaxiomsaccessiblesheaves}
    Suppose $(\mathcal C,S)$ is a (hyper)ac\-ces\-si\-ble explicit covering site.  
    The category of accessible, (hypercomplete) $S$-sheaves $\mywidehatshvacc{S}(\mathcal C)$ is a locally small category and satisfies all of Giraud's axioms ({\cite[Proposition 6.1.3.19]{highertopostheory}}) except accessibility, i.e.\ 
\begin{romanenum}
    \item \label{limitscolimitsaccessiblesheaveslist} $\mywidehatshvacc{S}(\mathcal C)$ has all small limits and colimits. 
    \item \label{universalitycolimitsaccessiblesheaveslist} Small colimits in $\mywidehatshvacc{S}(\mathcal C)$ are universal.
    \item \label{coproductsdisjointaccessiblesheaveslist} Small coproducts in $\mywidehatshvacc{S}(\mathcal C)$ are disjoint.
    \item \label{effectivitygroupoidsaccessiblesheaveslist} Every groupoid object in $\mywidehatshvacc{S}(\mathcal C)$ is effective. 
\end{romanenum}
\end{cor}
    
\begin{proof}The category $\Fun^{\operatorname{acc}}(\mathcal C^{\operatorname{op}}, \an)$ is locally small as every accessible functor is a small colimit of representables, cf.\ \cite[Proposition A.9]{HesselholtPstragowski2024}.
By \cref{limitscolimitscanbecomputedfinitestagesheaves}, $\mywidehatshvacc{S}(\mathcal C)$ has all small limits and colimits. 
By \cref{kappaccessibleistopos} and \cite[Theorem 6.1.0.6]{highertopostheory}, for all uncountable regular cardinals $\kappa$ such that $(\mathcal C,S)$ is $\kappa$-(hyper)ac\-ces\-si\-ble, \ref{universalitycolimitsaccessiblesheaveslist}, \ref{coproductsdisjointaccessiblesheaveslist}, and \ref{effectivitygroupoidsaccessiblesheaveslist} hold in $\mywidehatshvkacc{S}{\kappa}(\mathcal C)$. 
For every small diagram $F\colon S\to\mywidehatshvacc{S}(\mathcal C)$, there exists a regular cardinal $\kappa$ such that $F$ factors over \[ \mywidehatshvkacc{S}{\kappa}(\mathcal C)\subseteq\mywidehatshvacc{S}(\mathcal C)\] and $(\mathcal C,S)$ is $\kappa$-accessible. As $\mywidehatshvkacc{S}{\kappa}(\mathcal C)\subseteq\mywidehatshvacc{S}(\mathcal C)$ is closed under small colimits and finite limits (\cref{limitscolimitscanbecomputedfinitestagesheaves}), this implies \ref{universalitycolimitsaccessiblesheaveslist}, \ref{coproductsdisjointaccessiblesheaveslist}, and \ref{effectivitygroupoidsaccessiblesheaveslist}.  
\end{proof}

By \cref{accessiblesheavesnottopos}, categories of accessible (hyper)\-sheaves are typically not presentable and in particular not topoi (\cite[Theorem 6.1.0.6]{highertopostheory}). We now explain that they can be \textit{exhausted} by topoi, from which we will deduce that accessible presheaves can be sheafified. 
\begin{notation}
If $(\mathcal C,S)$ is a (hyper)ac\-ces\-si\-ble explicit covering site, denote by $\Lambda_{\mathcal C}$ the large poset of small regular cardinals $\kappa$ such that $(\mathcal C,S)$ is $\kappa$-(hyper)ac\-ces\-si\-ble, ordered by $\leq$.
\end{notation}
The poset $\Lambda_{\mathcal C}$ depends on $S$ and on whether we consider $(\mathcal C,S)$ as a hyperaccessible or as an accessible explicit covering site.
\begin{construction}\label{chainpresentablecategoriesaccessiblesheaves}
    Suppose that $\mathcal D$ is a category with small colimits and $\mu$ is a small regular cardinal such that $\mu$-filtered colimits commute with $\mu$-small limits in $\mathcal D$. 
    For a (hyper)ac\-ces\-si\-ble explicit covering site $(\mathcal C,S)$ denote by $\Lambda_{\mathcal C, \geq \mu}\subseteq \Lambda_{\mathcal C}$ the full subcategory on cardinals $\geq \mu$ and by $\Lambda_{\mathcal C, \geq \mu}^{\triangleright}\coloneqq \Lambda_{\mathcal C, \geq \mu}\cup\{\infty\}$ its cone. 
    The restrictions \[\mywidehatshvacc{S}(\mathcal C, \mathcal D)\to \mywidehatshv_S(\mathcal C_{\lambda}, \mathcal D)\to \mywidehatshv_S(\mathcal C_{\kappa}, \mathcal D), \,  \lambda\geq \kappa\in \Lambda_{\mathcal C, \geq \mu} \] assemble into a functor 
    \begin{align*} \mywidehatshv_S(\mathcal C_*, \mathcal D)\colon (\Lambda_{\mathcal C, \geq \mu}^{\triangleright})^{\operatorname{op}}& \to \vlCat, \\  \lambda& \mapsto \begin{cases} \mywidehatshv_S(\mathcal C_{\lambda}^{\operatorname{op}}, \mathcal D)& \lambda\neq \infty\\ 
    \mywidehatshvacc{S}(\mathcal C^{\operatorname{op}}, \mathcal D)& \lambda=\infty.\end{cases}\end{align*} By \cref{fullyfaithfulnessleftkanextension}, this enhances to a functor 
    \[  (\Lambda_{\mathcal C, \geq \mu}^{\triangleright})^{\operatorname{op}}\to \vlCat^R\] to the very large category of large categories and right adjoint functors. 
    Taking opposites yields a functor \[\mywidehatshv_S(\mathcal  C_*, \mathcal D)\colon \Lambda_{\mathcal C, \geq \mu}^{\triangleright}\to {\vlCat^L}\] to the very large category of large categories and left adjoint functors.
    We used that $\vlCat^L\cong (\vlCat^R)^{\operatorname{op}}$, see e.g.\ \cite[Theorem B]{HaugsengHebestreitLinskensNuiten} in the large universe $\mathcal U_1$.  
\end{construction}
\begin{rem}\label{identificationsheavesfunctorial}
If $\mathcal D$ is presentable, then $\mywidehatshv_S(\mathcal  C_*, \mathcal D)$ factors over $\Pr^L\subseteq \vlCat^L$ and the equivalence from \cref{identifysheaveswithtensorproductinprl} yields an equivalence \[\mywidehatshv_S(\mathcal C_*, \mathcal D)\cong \mywidehatshv_S(\mathcal C_*, \an)\otimes_{\Pr^L}\mathcal D\in \Fun(\Lambda_{\mathcal C, \geq \mu}, \Pr^L).\]  
\end{rem}
\begin{cor}\label{accessiblesheavesisbigtopos}

    Suppose that $\mathcal D$ is a category with small colimits and there exists a small regular cardinal $\mu$ such that $\mu$-filtered colimits commute with $\mu$-small limits in $\mathcal D$. 
    Then  
     \[\mywidehatshv_S(\mathcal  C_*, \mathcal D)\colon \Lambda_{\mathcal C, \geq \mu}^{\triangleright}\to {\vlCat^L}\to \vlCat\] 
 is a colimit diagram.
\end{cor}

\begin{proof}[Proof of \cref{accessiblesheavesisbigtopos}.]
By \cite[Lemma 4.8.4.2]{higheralgebra} (in the large universe $\mathcal U_1$), the colimit \[\colim{\lambda\in\Lambda_{\mathcal C}}\mywidehatshv_S(\mathcal C_{\lambda}, \mathcal D)\] exists in $\vlCat$. 
Since $\Lambda_{\mathcal C}$ is filtered, \cref{filteredcolimitsofcategoriesappendix}, \cref{fullyfaithfulnessleftkanextension}, and \cref{kappacondensedkappaacc} imply that the functor 
\[\colim{\lambda\in \Lambda_{\mathcal C, \geq \mu}}\mywidehatshv_S(\mathcal C_{\lambda}, \mathcal D)\cong \colim{\lambda\in \Lambda_{\mathcal C}}\mywidehatshv_S(\mathcal C_{\lambda}, \mathcal D)\to\mywidehatshvacc{S}(\mathcal C, \mathcal D)\] is fully faithful. 
By \cref{kappacondensedkappaacc}, it is essentially surjective. 
\end{proof}

\begin{rem}\label{accessiblesheavesisbigtoposrem}
\cref{accessiblesheavesisbigtopos} implies that if $(\mathcal C,S)$ is a (hyper)accessible explicit covering site, the category of accessible (hyper)sheaves is the colimit (in the category of large categories) of all categories of (hyper)sheaves on an essentially small subsite: 

If $\mathcal C^{'}\subseteq \mathcal C$ is a subcategory closed under finite limits and colimits, then $(\mathcal C^{'},S\cap\mathcal C^{'})$ is an explicit covering site, and if $\mathcal C^{'}\subseteq \mathcal C^{''}$ are two such subcategories which are essentially small, then the inclusion $i\colon \mathcal C^{'}\subseteq \mathcal C^{''}$ defines a morphism of sites, whence restriction $\mathcal P(\mathcal C^{''})\to\mathcal P(\mathcal C^{'})$ restricts to a functor $i^*\colon \mywidehatshv_S(\mathcal C^{''})\to\mywidehatshv_S(\mathcal C^{'})$.
This admits a left adjoint given by (hyper)sheafification after left Kan extension along $i$, which is left-exact by \cite[Proposition 6.2.3.20]{highertopostheory}. 
Denote by $\mathcal P_{\mathcal C}$ the poset of essentially small subcategories of $\mathcal C$ which are closed under finite limits and colimits. 
The left adjoints of restriction assemble into a diagram \[\mathcal P_{\mathcal C}\to \Pr^L,\,  \mathcal C^{'}\mapsto \mywidehatshv_S(\mathcal C^{'}).\] 
As $\mathcal C$ is $\kappa$-accessible, for all regular cardinals $\lambda\geq \kappa$, $\mathcal C_{\lambda}$ is essentially small (\cite[Proposition 5.4.2.2]{highertopostheory}), and every small subcategory is contained in $\mathcal C_{\lambda}$ for some small regular cardinal $\lambda$ (\cref{everyobjectcompactinanaccessiblecategory}). 
It therefore follows from the above corollary that for a (hyper)accessible explicit covering site $(\mathcal C,S)$ and a presentable category $\mathcal D$, \[\colim{\mathcal C^{'}\in \mathcal P_{\mathcal C}}\mywidehatshv_S(\mathcal C^{'},\mathcal D)\cong \mywidehatshvacc{S}(\mathcal C,\mathcal D).\] 

More generally, for an arbitrary, not necessarily accessible explicit covering site $(\mathcal C,S)$, one has a functor \[\mywidehatshv_S(-)\colon \mathcal P_{\mathcal C}\to \Pr^L, \mathcal C^{'}\mapsto \mywidehatshv_S(\mathcal C^{'}),\] and one can show that its colimit in $\vlCat$ satisfies Giraud's axioms except accessibility, see \cite[section 2.1]{Affinenessreconstructionincpgeometry}. If $(\mathcal C,S)$ is not (hyper)accessible, then the colimit $\colim{\mathcal C^{'}\subseteq \mathcal C}\mywidehatshv_S(\mathcal C^{'})$ (in $\vlCat$) seems to be better behaved than the category $\mywidehatshv_S^{\operatorname{acc}}(\mathcal C)$ of (hyper)accessible sheaves. However, this has the drawback that in general there need not exist a forget functor to the category of presheaves $\colim{\mathcal C^{'}\subseteq \mathcal C}\mathcal P(\mathcal C^{'})$.
\end{rem}
The above corollary (\ref{accessiblesheavesisbigtopos}) implies that accessible presheaves can be (hyper)-sheafified: 
\begin{cor}\label{accessiblesheafification}
 If $(\mathcal C,S)$ is a (hyper)ac\-ces\-si\-ble explicit covering site and $\mathcal D$ is a presentable category, the $\tau_S$-(hyper)\-sheafi\-fi\-ca\-tions $\Fun(\mathcal C_{\lambda}^{\operatorname{op}}, \mathcal D)\to \Shv_{\tau_S}(\mathcal C_{\lambda}, \mathcal D)$ assemble into a left adjoint \[\Fun^{\operatorname{acc}}(\mathcal C^{\operatorname{op}}, \mathcal D)\to \mywidehatshv_S^{\operatorname{acc}}(\mathcal C, \mathcal D)\] of the forget functor.  
\end{cor}
\begin{rem}
This was observed for fpqc-sheaves on affine schemes by \cite{Waterhouse-fpqc-sheafification} and stated for sheaves in \cite{barwick2019pyknoticobjectsibasic}. 
\end{rem}
\begin{proof}
We will deduce this from \cref{adjunctionfiberwise}.
We first recall that for all $\lambda\in\Lambda_{\mathcal C}$, the forget functor 
\[ \mywidehatshv_S(\mathcal C_{\lambda}, \mathcal D)\subseteq \Fun(\mathcal C_{\lambda}^{\operatorname{op}}, \mathcal D)\] admits a left adjoint. 
As $\mathcal C$ is $\lambda$-coaccessible, $\mathcal C_{\lambda}$ is essentially small. Choose a small category $\tilde{\mathcal C_{\lambda}}$ with an equivalence $f\colon \tilde{\mathcal C_{\lambda}}\cong \mathcal C_{\lambda}$ and denote by $f^*S\subseteq \Fun(\Delta^1, \tilde{\mathcal C_{\lambda}})$ the full subcategory on morphisms $s$ with $f(s)\in S$.
It is straightforward to check that $(\tilde{\mathcal C}_{\lambda}, f^*S)$ is an explicit covering site and that the pullback $f^*\colon\Fun(\mathcal C_{\lambda}^{\operatorname{op}}, \mathcal D)\cong \Fun(\tilde{\mathcal C}_{\lambda}^{\operatorname{op}}, \mathcal D)$ restricts to an equivalence \[\mywidehatshv_S(\mathcal C_{\lambda}, \mathcal D)\cong \mywidehatshv_{f^*S}(\tilde{\mathcal C}_{\lambda}, \mathcal D).\]  
The forget functor $\mywidehatshv_{f^*S}(\tilde{\mathcal C}_{\lambda}, \mathcal D)\subseteq \Fun(\tilde{\mathcal C}_{\lambda}^{\operatorname{op}}, \mathcal D)$ admits a left adjoint by e.g.\ \cite[Proposition A.3.7, Proposition A.3.14, Proposition A.3.16]{LucasMannthesis}. 
This implies that \[ \mywidehatshv_S(\mathcal C_{\lambda}, \mathcal D)\subseteq \Fun(\mathcal C_{\lambda}^{\operatorname{op}}, \mathcal D)\] admits a left adjoint $L_{\lambda}$. 

For $\kappa\leq\lambda\in\Lambda_{\mathcal C}$ denote by \[ i^{\lambda}_{\kappa}\colon \mywidehatshv_S(\mathcal C_{\kappa}, \mathcal D)\rightleftarrows \mywidehatshv_S(\mathcal C_{\lambda}, \mathcal D)\colon r^{\lambda}_{\kappa}\] the adjunction induced by restriction and left Kan extension (\cref{fullyfaithfulnessleftkanextension}). 
    Under the identification of \cref{accessiblesheavesisbigtopos}, the colimit of \[ \mywidehatshv_S(\mathcal C_{*}, \mathcal D)\xhookrightarrow{f_*}\Fun(\mathcal C_{*}^{\operatorname{op}}, \mathcal D)\in \Fun(\Lambda_{\mathcal C}^{\operatorname{op}}, \vlCat)\] over $\Lambda_{\mathcal C}$ is the forget functor $\mywidehatshvacc{S}(\mathcal C, \mathcal D)\subseteq \Fun^{\operatorname{acc}}(\mathcal C^{\operatorname{op}}, \mathcal D)$. 
    The Beck-Chevalley transformation $f_{\kappa}r^{\kappa}\Rightarrow r^{\kappa}f_{\lambda}$ associated to the commutative diagram \begin{center}
    \begin{tikzcd}
        \mywidehatshv_S(\mathcal C_{\kappa}, \mathcal D)\arrow[r,"i^{\lambda}_{\kappa}"]\arrow[d,"f_{\kappa}"]&  \mywidehatshv_S(\mathcal C_{\lambda}, \mathcal D)\arrow[d,"f_{\lambda}"]\\ 
         \Fun(\mathcal C_{\kappa}^{\operatorname{op}}, \mathcal D)\arrow[r,"i^{\lambda}_{\kappa}"] &\Fun(\mathcal C_{\lambda}^{\operatorname{op}}, \mathcal D)
    \end{tikzcd}
    \end{center} is obviously an equivalence, whence so is its opposite $L_{\lambda}i_{\kappa}\Rightarrow i_{\kappa}L_{\kappa}$ by \cite[Remark 4.7.4.14]{higheralgebra}. 
    It now follows from \cref{adjunctionsbigpresentable} (applied to the opposite categories $\Fun^{\operatorname{acc}}(\mathcal C_{*}, \mathcal D)^{\operatorname{op}}$ and \linebreak[4]$\mywidehatshv_S(\mathcal C_*, \mathcal D)^{\operatorname{op}}$), that there exists a left adjoint \[L\colon \Fun^{\operatorname{acc}}(\mathcal C^{\operatorname{op}}, \mathcal D)\to\mywidehatshvacc{S}(\mathcal C, \mathcal D)\] such that for all $\lambda\in\Lambda_{\mathcal C}$, $L|_{\Fun(\mathcal C_{\lambda}^{\operatorname{op}}, \mathcal D)}$ factors as \[\Fun(\mathcal C_{\lambda}^{\operatorname{op}}, \mathcal D)\xrightarrow{L_{\lambda}}\mywidehatshv_S(\mathcal C_{\lambda}, \mathcal D)\xhookrightarrow{i_{\lambda}}\mywidehatshvacc{S}(\mathcal C, \mathcal D).\qedhere\] 
\end{proof}
\begin{cor}
Suppose that $(\mathcal C,S)$ is a (hyper)ac\-ces\-si\-ble explicit covering site. 
The left adjoint of the forget functor $\mywidehatshvacc{S}(\mathcal C)\subseteq \Fun^{\operatorname{acc}}(\mathcal C^{\operatorname{op}}, \an)$ is left-exact.  
\end{cor}
\begin{proof}
   Since finite limits commute with filtered colimits in $\an$, for all uncountable cardinals $\kappa$, $\Fun^{\kappa\text{-}\operatorname{colim}}(\mathcal C^{\operatorname{op}}, \an)\subseteq \Fun^{\operatorname{acc}}(\mathcal C^{\operatorname{op}}, \an)$ is closed under finite limits which are computed pointwise, cf.\ \cite[Proposition A.9]{HesselholtPstragowski2024}. This implies that \[\mywidehatshvkacc{S}{\kappa}(\mathcal C)\subseteq\mywidehatshvacc{S}(\mathcal C)\] is closed under finite limits as well.
    Suppose that $F\colon I\to \Fun^{\operatorname{acc}}(\mathcal C^{\operatorname{op}}, \mathcal D)$ is a finite diagram, and choose a regular cardinal $\kappa$ such that $(\mathcal C,S)$ is $\kappa$-(hyper)ac\-ces\-si\-ble and $F(i)\in \Fun^{\kappa\text{-}\operatorname{colim}}(\mathcal C^{\operatorname{op}}, \an)$ for all $i\in I$.
    By construction, the left adjoint $L$ of the forget functor restricts to a functor \[\Fun^{\kappa\text{-}\operatorname{colim}}(\mathcal C^{\operatorname{op}}, \an)\cong \Fun(\mathcal C_{\kappa}^{\operatorname{op}}, \an)\to \mywidehatshv_S(\mathcal C_{\kappa}, \mathcal D)\cong  \mywidehatshvkacc{S}{\kappa}(\mathcal C),\] which is left adjoint to the forget functor. We recall below that this functor is left-exact, then it follows from the above that $L(\clim{i\in I}F)\cong \clim{i\in I}L(F_i)$. 

    As $\mathcal C_{\kappa}$ is essentially small, we can argue as in the proof of \cref{accessiblesheafification} to reduce to the case where $\mathcal C_{\kappa}$ is small. The identification of \cref{sheafconditionexplicitcoveringsite} and \cite[Proposition 6.2.2.7]{highertopostheory} imply that the left adjoint \[\Fun(\mathcal C_{\kappa}^{\operatorname{op}}, \an)\to \Shv_S(\mathcal C_{\kappa}, \an)\] is left-exact. Since $\infty$-connective morphisms are stable under pullbacks (\cite[Proposition 6.5.1.16]{highertopostheory}), the hypersheafification \[\Shv_S(\mathcal C_{\kappa}, \an)\to \hypershv_S(\mathcal C_{\kappa}, \an)\] is left-exact as well by \cite[Proposition 6.2.1.1]{highertopostheory}. 
\end{proof}
The definition of (hyper)accessible explicit covering sites is rigged so that \cref{accessiblesheavesisbigtopos} holds. Alternatively, one could require $\kappa$-accessibility for some regular $\kappa$ and that the collection of (hyper)covers is closed under $\kappa$-filtered limits: 
\begin{lemma}\label{kappaccessibilityimpliesaccessibilityifSclosedunderfilteredlimits}\begin{romanenum}Suppose that $(\mathcal C,S)$ is an explicit covering site. 
\item If $(\mathcal C,S)$ is $\kappa$-accessible and $S\subseteq \Fun(\Delta^1, \mathcal C)$ is closed under $\kappa$-filtered limits, then $(\mathcal C,S)$ is $\lambda$-accessible for all regular cardinals $\lambda\gg \kappa$.

\item Denote by $S_{\Delta}\subseteq \Fun(\Delta_{+}^{\operatorname{op}}, \mathcal C)$ the full subcategory on $S$-hypercovers. 
If $(\mathcal C,S)$ is $\kappa$-hy\-per\-ac\-ces\-sib\-le and $S_{\Delta}\subseteq \Fun(\Delta_{s,+}^{\operatorname{op}}, \mathcal C)$ is closed under $\kappa$-filtered limits, then $(\mathcal C,S)$ is $\lambda$-hy\-per\-ac\-ces\-sib\-le for all regular cardinals $\lambda\gg \kappa$.
\end{romanenum}
\end{lemma}
\begin{proof}We only show the statement about accessibility, the hyperaccessibility statement can be proven completely analogously.     
Fix a regular cardinal $\lambda\gg \kappa$. 
By \cite[Proposition 5.4.2.11]{highertopostheory}, $\mathcal C$ is $\lambda$-coaccessible. The proof of \cite[Proposition 5.4.4.3]{highertopostheory} shows that for every small category $I$, 
\[\Fun(I^{\operatorname{op}}, \mathcal C)^{\operatorname{op}}\cong \Fun(I, \mathcal C^{\operatorname{op}})\] is $\kappa$-accessible, and the left Kan extension \[\mathcal P(\Fun(I^{\operatorname{op}}, \mathcal C_{\kappa})^{\operatorname{op}})\xhookrightarrow{i_{\kappa}} \mathcal P(\Fun(I^{\operatorname{op}}, \mathcal C)^{\operatorname{op}})\] along $\Fun(I^{\operatorname{op}}, \mathcal C_{\kappa})^{\operatorname{op}}\subseteq \Fun(I^{\operatorname{op}}, \mathcal C)^{\operatorname{op}}$ restricts to an equivalence \[\Ind_{\kappa}(\Fun(I^{\operatorname{op}}, \mathcal C_{\kappa})^{\operatorname{op}})\cong  \Fun(I^{\operatorname{op}}, \mathcal C)^{\operatorname{op}}, \] where we identified the right-hand side with its Yoneda image. 
It follows that $(\mathcal C,S)$ is $\lambda$-accessible if and only if (the Yoneda image of) $S^{\operatorname{op}}\subseteq \Fun(\Delta^{1}, \mathcal C)^{\operatorname{op}}$ is contained in the image of the composite \[\Ind_{\lambda}\left((\Fun(\Delta^1, \mathcal C_{\lambda})\cap S)^{\operatorname{op}}\right)\to \Ind_{\lambda}\left(\Fun(\Delta^1, \mathcal C_{\lambda})^{\operatorname{op}}\right)\subseteq  \mathcal P(\Fun(\Delta^1, \mathcal C_{\lambda})^{\operatorname{op}})\xrightarrow{i_{\lambda}} \mathcal P(\Fun(\Delta^1, \mathcal C)^{\operatorname{op}}).\] 
The left Kan extension $\mathcal P(\Fun(\Delta^1, \mathcal C_{\kappa})^{\operatorname{op}})\xhookrightarrow{i_{\kappa}} \mathcal P(\Fun(\Delta^1, \mathcal C)^{\operatorname{op}})$ factors into left adjoints \[\mathcal P(\Fun(\Delta^1, \mathcal C_{\kappa})^{\operatorname{op}})\xhookrightarrow{i_{\kappa}^{\lambda}} \mathcal P(\Fun(\Delta^1, \mathcal C_{\lambda})^{\operatorname{op}})\xhookrightarrow{i_{\lambda}} \mathcal P(\Fun(\Delta^1, \mathcal C)^{\operatorname{op}}),\] and by \cite[Lemma 5.4.2.10]{highertopostheory}, $i_{\kappa}^{\lambda}$ restricts to an equivalence \[\Ind_{\kappa}(\Fun(\Delta^1, \mathcal C_{\kappa})^{\operatorname{op}})\cong \Ind_{\lambda}(\Fun(\Delta^1, \mathcal C_{\lambda})^{\operatorname{op}}).\] 
Since $S\subseteq  \Fun(\Delta^1, \mathcal C)$ is closed under $\kappa$-filtered limits, \cite[Lemma 5.4.2.10]{highertopostheory} implies that \[ \Ind_{\kappa}\left((\Fun(\Delta^1, \mathcal C_{\kappa})\cap S)^{\operatorname{op}}\right)\to \Ind_{\kappa}(\Fun(\Delta^1, \mathcal C_{\kappa})^{\operatorname{op}})\cong \Ind_{\lambda}\left((\Fun(\Delta^1, \mathcal C_{\lambda}))^{\operatorname{op}}\right)\] factors over $\Ind_{\lambda}\left((\Fun(\Delta^1, \mathcal C_{\lambda})\cap S)^{\operatorname{op}}\right)$, which shows that $(\mathcal C,S)$ is $\lambda$-accessible.
\end{proof}
We now deduce from \cref{accessiblesheavesisbigtopos} that for a presentably symmetric monoidal category $\mathcal D^{\otimes}$, $\mywidehatshvacc{S}(\mathcal C, \mathcal D)$ inherits a symmetric monoidal structure:
\begin{construction}\label{constructionmonoidalstructureaccessiblesheaves}
Endow $\Pr^L$ with the Lurie tensor product (\cite[Proposition 4.8.1.15]{higheralgebra}). 
For all cardinals $\kappa\in\Lambda_{\mathcal C}$, the cartesian monoidal structure on the topos $\mywidehatshv_S(\mathcal C_{\kappa})$ is cocontinuous (\cite[Theorem 6.1.0.6.3]{highertopostheory}) and for $\kappa\leq\lambda\in\Lambda_{\mathcal C}$, the left adjoint \[\mywidehatshv_S(\mathcal C_{\kappa})\to \mywidehatshv_S(\mathcal C_{\lambda})\] is left-exact by \cref{limitscolimitscanbecomputedfinitestagesheaves} and \cref{kappacondensedkappaacc}.  
\cite[Corollary 2.4.1.9]{higheralgebra} now implies that $\mywidehatshv_S(\mathcal C_*)$ enhances to a functor \[\mywidehatshv_S(\mathcal C_{*})^{\times}\colon \Lambda_{\mathcal C}\to \CAlg(\Pr^L).\] 

By \cite[Example 3.2.4.4]{higheralgebra}, $\CAlg(\Pr^L)$ admits a symmetric monoidal structure and the forget functor $\CAlg(\Pr^L)\to \Pr^L$ enhances to a symmetric monoidal functor. In particular, if $\mathcal C^{\otimes}, \mathcal D^{\otimes}$ are two presentably symmetric monoidal categories, then $\mathcal C^{\otimes}\otimes_{\CAlg(\Pr^L)} \mathcal D^{\otimes}\in \CAlg(\Pr^L)$ is a presentably symmetric monoidal enhancement of $\mathcal C\otimes_{\Pr^L} \mathcal D$.
Let 
\begin{align*}\mywidehatshv_S(\mathcal C_{*},-)^{\otimes}\colon \Lambda_{\mathcal C}\times \CAlg(\Pr^L) &\to \CAlg(\Pr^L), \\  (\lambda, \mathcal D^{\otimes}) &\mapsto \mywidehatshv_S(\mathcal C_{\lambda})^{\times}\otimes_{\CAlg(\Pr^L)}\mathcal D^{\otimes}.\end{align*}
By \cite[Lemma 4.8.4.2]{higheralgebra} in the universe $\mathcal U_1$, the category $\vlCat$ has all $\mathcal U_1$-small colimits, and the tensor product on $\vlCat$ is compatible with $\mathcal U_1$-small colimits (\cite[Definition 3.1.1.18]{higheralgebra}). 
Hence by \cite[Proposition 3.2.3.1]{higheralgebra} (in $\mathcal U_1$), $\CAlg(\vlCat)$ has all large sifted colimits and \[\CAlg(\vlCat)\to \vlCat\] preserves large sifted colimits. 
Together with \cref{accessiblesheavesisbigtopos} and \cref{identificationsheavesfunctorial}, this implies that by taking the colimit of $\mywidehatshv_S(\mathcal C_{*},-)^{\otimes}$ over $\Lambda_{\mathcal C}$ in $\CAlg(\vlCat)$, we obtain a functor 
\[ \CAlg(\Pr^L)\to \CAlg(\vlCat), \, \, \mathcal D^{\otimes}\mapsto \colim{\lambda\in \Lambda}\mywidehatshv_S(\mathcal C_{\lambda}, \mathcal D)^{\otimes}\] which sends a presentably symmetric monoidal category $\mathcal D$ to a symmetric monoidal enhancement of $\mywidehatshvacc{S}(\mathcal C, \mathcal D)$. 
\end{construction}
\begin{thm}\label{closednessmonoidalstructureaccessiblesheaves}
    Suppose that $(\mathcal C,S)$ is a (hyper)accessible explicit covering site.    
    For a presentably symmetric monoidal category $\mathcal D^{\otimes}$, the symmetric monoidal structure on $\mywidehatshvacc{S}(\mathcal C, \mathcal D)$ described above is closed. 
    \end{thm}
    This is a straightforward adaptation of \cite[Proposition 2.1.11]{LucasMannthesis} and relies on the following observation.  
    \begin{lemma}\label{generatorsmonoidalstructure}
    Suppose that $(\mathcal C,S)$ is a (hyper)ac\-ces\-si\-ble explicit covering site and $\mathcal D^{\otimes}$ is a presentably symmetric monoidal category. 
    Suppose that $\kappa$ is a regular cardinal such that $(\mathcal C,S)$ is $\kappa$-(hyper)ac\-ces\-si\-ble. 
    For $c\in\mathcal C_{\kappa}$, the functor \[\operatorname{ev}_c\colon\mywidehatshvacc{S}(\mathcal C, \mathcal D)\to \mathcal D, \,  F\mapsto F(c)\] admits a left adjoint $(-)[c]$.
    This has the following properties: 
    \begin{romanenum}
    \item For $d\in\mathcal D$, $d[c]\in \mywidehatshvkacc{S}{\kappa}(\mathcal C, \mathcal D)$. 
    \item Suppose that $G\subseteq \mathcal D$ is a small generating set for $\mathcal D$. 
    Then $\{d[c],c\in\mathcal C_{\kappa}, d\in G\}$ generates $\mywidehatshv^{\kappa\text{-}\operatorname{colim}}_S(\mathcal C, \mathcal D)$ under small colimits. 
    \item For $c, \tilde c\in\mathcal C_{\kappa}$ and $d, \tilde d\in \mathcal D$, 
    \[ d[c]\otimes \tilde d[\tilde c]\cong (d\otimes_{\mathcal D} \tilde d)[c\times \tilde c]\] naturally.  
    \end{romanenum}
    \end{lemma}
    \begin{proof}
    We first explain that $\operatorname{ev}_c$ admits a left adjoint for all $c\in\mathcal C_{\kappa}$. 
    As $\operatorname{ev}_c$ factors as \[ \mywidehatshvacc{S}(\mathcal C,\mathcal D)\to \mywidehatshv_S(\mathcal C_{\kappa},\mathcal D)\xrightarrow{\operatorname{ev}_c}\mathcal D\] and the restriction $\mywidehatshvacc{S}(\mathcal C,\mathcal D)\to \mywidehatshv_S(\mathcal C_{\kappa},\mathcal D)$ has a left adjoint $i_{\kappa}$ (\cref{kappacondensedkappaacc}), it suffices to show that \[\mywidehatshv_S(\mathcal C_{\kappa},\mathcal D)\xrightarrow{\operatorname{ev}_c}\mathcal D\] admits a left adjoint.
    We showed in the proof of \cref{accessiblesheavesnottopos} that the forget functor \[\mywidehatshv_S(\mathcal C_{\kappa}, \mathcal D)\subseteq \Fun(\mathcal C_{\kappa}^{\operatorname{op}}, \mathcal D)\] has a left adjoint $L_{\kappa}$ for all $\kappa\in\Lambda_{\mathcal C}$. 
    The category $\Fun(\mathcal C_{\kappa}^{\operatorname{op}}, \mathcal D)$ is presentable by \cite[Proposition 5.4.4.3]{highertopostheory}. 
    As $\mathcal D$ has small limits and colimits (\cite[Theorem 5.5.1.1, Corollary 5.5.2.4]{highertopostheory}), \[ \operatorname{ev}_c\colon \Fun(\mathcal C_{\kappa}^{\operatorname{op}}, \mathcal D)\to\mathcal D\] preserves small limits and colimits by \cite[Corollary 5.1.2.3]{highertopostheory} and in particular admits a left adjoint $p_c$ by the adjoint functor theorem (\cite[Corollary 5.5.2.9]{highertopostheory}). 
    This shows that \[(-)[c]\coloneqq i_{\kappa}\circ L_{\kappa}\circ p_c\] is left adjoint to $\operatorname{ev}_c\colon \mywidehatshvacc{S}(\mathcal C, \mathcal D)\to \mathcal D$, and in particular, $d[c]\in \mywidehatshv^{\kappa\text{-}\operatorname{colim}}_S(\mathcal C, \mathcal D)$ for all $d\in\mathcal D$.  

    Suppose now that $G$ is a small generating set for $\mathcal D$. This exists since $\mathcal D$ is presentable. Then the functors \[\Map_{\mywidehatshv_S(\mathcal C_{\kappa}, \mathcal D)}(c[g],-)\cong \Map_{\mywidehatshv_S(\mathcal C_{\kappa}, \mathcal D)}(g, \operatorname{ev}_c-), \, c\in\mathcal C_{\kappa},g\in G\] are jointly conservative.
    As $\mathcal C_{\kappa}$ is essentially small (\cite[Proposition 5.4.2.2]{highertopostheory}) and the category $\mywidehatshv_S(\mathcal C_{\kappa}, \mathcal D)$ is accessible (\cref{kappaccessibleistopos}), by \cref{everyobjectcompactinanaccessiblecategory} there exists a regular cardinal $\mu$ such that $c[g]$ is $\mu$-compact for all $c\in\mathcal C_{\kappa}, g\in G$. It now follows from \cite[Lemma A.2.1]{LucasMannthesis} that $c[g],c\in\mathcal C, g\in G$ generate \[\mywidehatshv_S(\mathcal C_{\kappa}, \mathcal D)\cong \mywidehatshv^{\kappa\text{-}\operatorname{colim}}_S(\mathcal C, \mathcal D)\] under colimits. 

    It remains to prove the relation \[ d[c]\otimes \tilde d[\tilde c]\cong (d\otimes\tilde d)[c\otimes \tilde c].\]
    Denote by $l\colon \Fun(\mathcal C_{\kappa}^{\operatorname{op}}, \an)\to \mywidehatshv_S(\mathcal C_{\kappa}, \an)$ the $S$-(hyper)\-sheafi\-fi\-ca\-tion. The proof of \cite[Proposition 4.8.1.15]{highertopostheory} and \cref{sheafconditionexplicitcoveringsite} imply that \[ l\otimes_{\Pr^L}\id_{\mathcal D}\colon \Fun(\mathcal C_{\kappa}^{\operatorname{op}}, \mathcal D)\to \mywidehatshv_S(\mathcal C_{\kappa}, \mathcal D)\] is the left adjoint of the inclusion, i.e. equivalent to $L_{\kappa}$.  
    Since $l$ is cartesian, $l$ enhances to a symmetric monoidal functor \[\Fun(\mathcal C_{\kappa}^{\operatorname{op}}, \an)^{\times}\to  \mywidehatshv_S(\mathcal C_{\kappa}, \an)^{\times}\in\CAlg(\Pr^L), \] cf.\ \cite[Corollary 2.4.1.9]{higheralgebra}. Hence by construction of the symmetric monoidal structure on $\mywidehatshv_S(\mathcal C_{\kappa}, \mathcal D)$, $L_{\kappa}$ enhances to a symmetric monoidal functor.
    By construction of the symmetric monoidal structure on $\mywidehatshvacc{S}(\mathcal C,\mathcal D)$, $i_{\kappa}$ is symmetric monoidal as well. 
    It therefore suffices to show that for $c, \tilde c\in\mathcal C_{\kappa},d, \tilde d\in\mathcal D$, \[p_c(d)\otimes p_{\tilde c}(\tilde d)\cong p_{c\otimes \tilde c}(d\otimes \tilde d).\] 
    For $t\in\mathcal C, x\in\mathcal D$, $\Map_{\Fun(\mathcal C^{\operatorname{op}}, \mathcal D)}(p_t(x),F)=\Map_{\mathcal D}(x,F(t))$, whence $p_t(x)\colon\mathcal C^{\operatorname{op}}\to\mathcal D$ is the left Kan extension of $\{x\}\to\mathcal D$ along $\{t\}\to \mathcal C^{\operatorname{op}}$.  
    As the symmetric monoidal structure on $\mathcal C^{\operatorname{op}}$ is cocartesian, this implies that  $p_{c\otimes \tilde c}(d\otimes \tilde d)$ is the left Kan extension of 
    \begin{align*}p_c(d)\otimes^{ptw}p_d(\tilde d)\colon \mathcal C^{\operatorname{op}}\times\mathcal C^{\operatorname{op}}&\to \mathcal D,\\ (x,y)& \mapsto p_c(d)(x)\otimes_{\mathcal D} p_{\tilde c}(\tilde d)(y)\end{align*} along $\mathcal C^{\operatorname{op}}\times\mathcal C^{\operatorname{op}}\xrightarrow{-\otimes-}\mathcal C^{\operatorname{op}}$, which shows that $p_{c\otimes \tilde c}(d\otimes \tilde d)\cong p_c(d)\otimes p_{\tilde c}(\tilde d),$ cf.\ \cref{formulasymmetricmonoidalstructureontensorproductofpresentablecatsasdayconcvolution}. 
    \end{proof}
    \begin{proof}[Proof of \cref{closednessmonoidalstructureaccessiblesheaves}]
    Since small colimits can always be computed on some stage \[\mywidehatshvkacc{S}{\kappa}(\mathcal C, \mathcal D)\subseteq \mywidehatshvacc{S}(\mathcal C, \mathcal D)\] (\cref{limitscolimitscanbecomputedfinitestagesheaves}), the tensor product on $\mywidehatshvacc{S}(\mathcal C)$ preserves small colimits in both variables. As $\mywidehatshvacc{S}(\mathcal C, \mathcal D)$ has small limits and colimits, it follows from \cref{leftadjointsstableundercolimits} that \[\{ M\in\mywidehatshvacc{S}(\mathcal C, \mathcal D)\, |\, M\otimes -\text{ admits right adjoint}\}\subseteq \mywidehatshvacc{S}(\mathcal C, \mathcal D)\] is closed under small colimits.

    Fix a regular cardinal $\mu$ such that $\mathcal D$ is $\mu$-accessible. 
    For $\lambda\geq \kappa\in\Lambda_{\mathcal C, \geq \mu}$, denote by \[\mywidehatshv_S(\mathcal C_{\kappa}, \mathcal D)\xrightarrow{i^{\lambda}_{\kappa}}\mywidehatshv_S(\mathcal C_{\lambda}, \mathcal D)\xrightarrow{i^{\lambda}}\mywidehatshvacc{S}(\mathcal C, \mathcal D)\] the left adjoints of the restriction, cf.\ \cref{fullyfaithfulnessleftkanextension}, 
    \cref{kappacondensedkappaacc}. \cref{generatorsmonoidalstructure,kappacondensedkappaacc} imply that 
    \[ i_{\kappa}d_{\kappa}[c], c\in\mathcal C_{\kappa},d\in\mathcal D_{\mu}, \kappa\in\Lambda_{\mathcal C, \geq \mu}\]  generates $\mywidehatshvacc{S}(\mathcal C, \mathcal D)$ under small colimits, where $(-)_{\kappa}[c]$ denotes the left adjoint of \[\operatorname{ev}_c\colon \mywidehatshv_S(\mathcal C_{\kappa}, \mathcal D)\to \mathcal D.\] It therefore suffices to show that $ i_{\kappa}d_{\kappa}[c]\otimes-$ is a left adjoint for all $\kappa\in\Lambda_{\mathcal C, \geq \mu}, c\in\mathcal C_{\kappa},d\in\mathcal D_{\mu}$. 
    Fix $d\in\mathcal D_{\mu}$ and choose a regular cardinal $\eta\geq \mu$ such that $\imap_{\mathcal D}(d,-)\colon\mathcal D\to\mathcal D$ preserves $\eta$-filtered colimits. This exists by the adjoint functor theorem (\cite[Corollary 5.5.2.9]{highertopostheory}). 
    Choose $\kappa\in\Lambda_{\mathcal C, \geq\eta}$, and $c\in\mathcal C_{\kappa}$. We first show that for all $\lambda\geq\kappa\in\Lambda_{\mathcal C}$ and 
    $N\in\mywidehatshv_S(\mathcal C_{\kappa}, \mathcal D)$, \[\imap_{\mywidehatshv_S(\mathcal C_{\lambda}, \mathcal D)}(i_{\kappa}^{\lambda}d_{\kappa}[c],i_{\kappa}^{\lambda}N)\colon C_{\lambda}^{\operatorname{op}}\to\mathcal D\] preserves $\kappa$-filtered colimits. For $e\in\mathcal D$, 
\begin{align*}&\Map_{\mathcal D}(e, \imap_{\mywidehatshv_S(\mathcal C_{\lambda}, \mathcal D)}(i_{\kappa}^{\lambda}d_{\kappa}[c],i_{\kappa}^{\lambda}N)(-))\\& \cong \Map_{\mywidehatshv_S(\mathcal C_{\lambda}, \mathcal D)}(e_{\lambda}[-], \imap_{\mywidehatshv_S(\mathcal C_{\lambda}, \mathcal D)}(i^{\lambda}_{\kappa}d_{\kappa}[c],i^{\lambda}_{\kappa}N))\\ & \cong \Map_{\mywidehatshv_S(\mathcal C_{\lambda}, \mathcal D)}(e_{\lambda}[-]\otimes i^{\lambda}_{\kappa}d_{\kappa}[c],i^{\lambda}_{\kappa}N), \end{align*} and since $i^{\lambda}_{\kappa}d_{\kappa}[c]\cong d_{\lambda}[c]$, 
    \[e_{\lambda}[-]\otimes i^{\lambda}_{\kappa}d_{\kappa}[c]\cong (e\otimes d)_{\lambda}[-\times c]\] by \cref{generatorsmonoidalstructure}. 
    This implies that \begin{align*} \Map_{\mathcal D}(e, \imap_{\mywidehatshv_S(\mathcal C_{\lambda}, \mathcal D)}(i_{\kappa}^{\lambda}d_{\kappa}[c],i_{\kappa}^{\lambda}N)(-))& \cong \Map_{\mathcal D}((e\otimes d)_{\lambda}[-\times c],i^{\lambda}_{\kappa}N)\\ &\cong \Map_{\mathcal D}(e\otimes d,i^{\lambda}_{\kappa}N(-\times c))\\ & \cong \Map_{\mathcal D}(e,\imap_{\mathcal D}(d,i^{\lambda}_{\kappa}N(-\times c))),\end{align*}
    which shows that \[\imap_{\mywidehatshv_S(\mathcal C_{\lambda},\mathcal D)}(i_{\kappa}^{\lambda}d_{\kappa}[c],i_{\kappa}^{\lambda}N)=(\mathcal C_{\lambda}^{\operatorname{op}}\ni x \mapsto \imap_{\mathcal D}(d,i^{\lambda}_{\kappa}N(x\times c)).\] 

    Since $i^{\lambda}_{\kappa}N\colon\mathcal C_{\lambda}^{\operatorname{op}}\to\mathcal D$ preserves $\kappa$-filtered colimits and the symmetric monoidal structure on $\mathcal C_{\lambda}$ is cartesian, this functor commutes with $\kappa$-filtered colimits by choice of $\eta$. 
    The symmetric monoidal enhancement of $i^{\lambda}_{\kappa}$ from \cref{constructionmonoidalstructureaccessiblesheaves} provides a commutative diagram \begin{center}\label{commtativediagramsymmetricmonoidalstructure}
        \begin{tikzcd}
        \mywidehatshv_S(\mathcal C_{\kappa}, \mathcal D)\arrow[rr,"d_{\kappa}\lbrack c\rbrack \otimes-"]\arrow[d,"i_{\kappa}^{\lambda}"] &&\arrow[d] \mywidehatshv_S(\mathcal C_{\kappa}, \mathcal D)\arrow[d,"i_{\kappa}^{\lambda}"]\\ 
        \mywidehatshv_S(\mathcal C_{\lambda}, \mathcal D)\arrow[rr,"i_{\kappa}^{\lambda}d_{\kappa}\lbrack c\rbrack\otimes-"] && \mywidehatshv_S(\mathcal C_{\lambda}, \mathcal D).
        \end{tikzcd}
    \end{center}
    We claim that the associated Beck-Chevalley transformation 
    \[ \beta\colon  i^{\lambda}_{\kappa}\imap_{\mywidehatshv_{S}(\mathcal C_{\kappa}, \mathcal D)}(d_{\kappa}[c],-)\Rightarrow\imap_{\mywidehatshv_{S}(\mathcal C_{\lambda}, \mathcal D)}(i^{\lambda}_{\kappa}d_{\kappa}[c],i^{\lambda}_{\kappa}-)\] is an equivalence, then it follows from \cref{adjunctionsbigpresentable} that $i_{\kappa}d_{\kappa}[c]\otimes_{\mywidehatshvacc{S}(\mathcal C, \mathcal D)}-$ admits a right adjoint. 
    Since for all $N\in \mywidehatshv_S(\mathcal C_{\kappa},\mathcal D)$, \[\imap_{\mywidehatshv_{S}(\mathcal C_{\lambda}, \mathcal D)}(i^{\lambda}_{\kappa}d_{\kappa}[c],i^{\lambda}_{\kappa}N)\colon \mathcal C_{\lambda}^{\operatorname{op}}\to\mathcal D\] preserves $\kappa$-filtered colimits (see above), by \cref{kappacondensedkappaacc} it suffices to show that $r^{\lambda}_{\kappa}(\beta)$ is an equivalence.
    As right adjoint of a symmetric monoidal functor, $r^{\lambda}_{\kappa}$ enhances to a lax symmetric monoidal functor (\cref{adjunctionfiberwise}).
    The lax monoidal structure yields a $2$-cell  
    \begin{center}
   \begin{tikzcd}
        \mywidehatshv_S(\mathcal C_{\lambda}, \mathcal D)\arrow[rr,"i^{\lambda}_{\kappa}d_{\kappa}\lbrack c\rbrack\otimes-"]\arrow[d,"r^{\lambda}_{\kappa}"] &&\arrow[d] \mywidehatshv_S(\mathcal C_{\lambda}, \mathcal D)\arrow[d,"r^{\lambda}_{\kappa}"]\\ 
        \mywidehatshv_S(\mathcal C_{\kappa}, \mathcal D)\arrow[rr,"r^{\lambda}_{\kappa}i^{\lambda}_{\kappa}d_{\kappa}\lbrack c\rbrack\otimes-"] \arrow[urr, Rightarrow]&& \mywidehatshv_S(\mathcal C_{\kappa}, \mathcal D).
        \end{tikzcd}
    \end{center}
    Denote by \[\gamma\colon r^{\lambda}_{\kappa}\imap_{\mywidehatshv_S(\mathcal C_{\lambda}, \mathcal D)}(i^{\lambda}_{\kappa}d_{\kappa}[c],-)\Rightarrow \imap_{\mywidehatshv_S(\mathcal C_{\kappa}, \mathcal D)}(r^{\lambda}_{\kappa}i^{\lambda}_{\kappa}d_{\kappa}[c],-)\circ r^{\lambda}_{\kappa}\] the associated Beck-Chevalley transformation.
    By fully faithfulness of $i^{\lambda}_{\kappa}$ (\cref{kappacondensedkappaacc}), the unit for $i^{\lambda}_{\kappa}\dashv r^{\lambda}_{\kappa}$ induces an equivalence \[\eta^*\colon \imap_{\mywidehatshv_S(\mathcal C_{\kappa}, \mathcal D)}(r^{\lambda}_{\kappa}i^{\lambda}_{\kappa}d_{\kappa}[c],-)\circ r^{\lambda}_{\kappa}\cong \imap_{\mywidehatshv_S(\mathcal C_{\kappa}, \mathcal D)}(d_{\kappa}[c],-)\circ r^{\lambda}_{\kappa}.\]
    As $\eta^*\circ \gamma$ is the total mate (\cite[Definition F.7]{BastiaanCnossenPhDthesis}) of the equivalence $i^{\lambda}_{\kappa}d_{\kappa}[c]\otimes -\cong i^{\lambda}_{\kappa}(d_{\kappa}[c]\otimes-)$, it is an equivalence. 
    This shows that $\gamma$ is an equivalence. 
    By the pasting law for mates (\cite[Lemma 2.2.4]{carmeli2022ambidexterity}), $\gamma\circ \beta$ is the Beck-Chevalley transformation induced by  
      \begin{center}
   \begin{tikzcd}
        \mywidehatshv_S(\mathcal C_{\kappa}, \mathcal D)\arrow[rr,"d_{\kappa}\lbrack c\rbrack\otimes-"]\arrow[d,"r^{\lambda}_{\kappa}i_{\kappa}"] &&\mywidehatshv_S(\mathcal C_{\kappa}, \mathcal D)\arrow[d,"r^{\lambda}_{\kappa}i_{\kappa}"]\\ 
        \mywidehatshv_S(\mathcal C_{\kappa}, \mathcal D)\arrow[rr,"r^{\lambda}_{\kappa}i^{\lambda}_{\kappa}d_{\kappa}\lbrack c\rbrack \otimes-"]\arrow[urr,Rightarrow] && \mywidehatshv_S(\mathcal C_{\kappa}, \mathcal D),
        \end{tikzcd}
    \end{center} where the 2-cell comes from the lax monoidal structure on $r^{\lambda}_{\kappa}i_{\kappa}^{\lambda}$. By fully faithfulness of $i^{\lambda}_{\kappa}$ (\cref{fullyfaithfulnessleftkanextension}), the 2-cell and the vertical arrows are equivalences, which implies that $\gamma\circ \beta$ is an equivalence. 
    This shows that $r^{\lambda}_{\kappa}\beta$ is an equivalence. 
    \end{proof}
    \begin{rem}\label{forkappalargeenoughclosedmonoidal}
        The above proof also shows that for $M,N\in\mywidehatshvacc{S}(\mathcal C, \mathcal D)$, there exists a regular cardinal $\kappa$ such that for all $\kappa\leq \lambda\in\Lambda_{\mathcal C}$, 
        \[i_{\lambda}\imap_{\mywidehatshv_{S}(\mathcal C_{\lambda}, \mathcal D)}(M,N)\cong \imap_{\mywidehatshvacc{S}(\mathcal C, \mathcal D)}(M,N).\] 
        Indeed, for $M,N\in\mywidehatshvacc{S}(\mathcal C, \mathcal D)$, choose a regular cardinal $\kappa_0$ such that $(\mathcal C,S)$ is $\kappa_0$-(hyper)ac\-ces\-si\-ble, $M$ and $N$ preserve $\kappa_0$-filtered colimits, and there exists an uncountable regular cardinal $\mu\leq \kappa_0$ such that $\mathcal D$ is $\mu$-accessible.
        By \cref{kappacondensedkappaacc} and \cref{generatorsmonoidalstructure}, $M$ is a small colimit of elements of the form $i_{\kappa_0}d_{\kappa_0}[c], c\in\mathcal C_{\kappa_0}, d\in\mathcal D$. 
        Choose a small regular cardinal $\kappa_1\geq\kappa_0$ such that $M=\colim{j\in J}i_{\kappa_1}d^j_{\kappa_1}[c^j]$ is a $\kappa_1$-small colimit of such elements and $\mathcal D$ is $\kappa_1$-accessible. 
        Choose a regular cardinal $\kappa\geq \kappa_1$ such that $\imap_{\mathcal D}(d^j,-)$ preserves $\kappa$-filtered colimits for all $j\in J$. 
        We showed in the proof of \cref{closednessmonoidalstructureaccessiblesheaves} that for $\lambda\geq \kappa\in\Lambda_{\mathcal C}$, \[i_{\lambda}\imap_{\mywidehatshv_{S}(\mathcal C_{\lambda}, \mathcal D)}(i_{\kappa}^{\lambda}d^j_{\kappa}[c^j],N)\cong \imap_{\mywidehatshvacc{S}(\mathcal C, \mathcal D)}(i_{\kappa}d^j_{\kappa}[c^j],N).\]
        Since $\mywidehatshv_{S}(\mathcal C_{\lambda}, \mathcal D)\subseteq\mywidehatshvacc{S}(\mathcal C, \mathcal D)$ is closed under $\kappa$-small limits (\cref{limitscolimitscanbecomputedfinitestagesheaves}), it follows that \[i_{\lambda}\imap_{\mywidehatshv_{S}(\mathcal C_{\lambda}, \mathcal D)}(M,N)\cong \imap_{\mywidehatshvacc{S}(\mathcal C, \mathcal D)}(M,N)\] for $\kappa\leq \lambda\in\Lambda_{\mathcal C}$.  
    \end{rem}
\subsection{Big presentable categories}
Categories of accessible (hyper)\-sheaves on a (hyper)ac\-ces\-si\-ble explicit covering site are in general not accessible (\cref{accessiblesheavesnottopos}) and therefore not topoi, cf.\ \cite[Theorem 6.1.0.6]{highertopostheory}. 
Our next goal is to show that however, their category of spectrum objects behaves like the stabilization of a presentably symmetric monoidal category in many aspects. 
We chose to work with the following more general setup to also include the category of $G$-objects in a category of accessible (hyper)sheaves for a group object $G$ in our discussion. 
\begin{definition}\label{definitionbigpresentable}
\begin{romanenum}
    \item A \emph{chain of presentable} \emph{categories} is a functor 
    $\mt{B}_*\colon \Lambda\to \Pr^L$ from a possibly large filtered set $\Lambda$ such that for all $\kappa\to \lambda\in \Lambda$, $\topo{B_{\kappa}}\to \topo{B_{\lambda}}$ is a fully faithful, left-exact left adjoint. 
    A possibly large category $\mt{B_{\infty}}$ is \emph{big presentable} if there exists a chain of presentable categories $\mathcal B\colon \Lambda\to \Pr^L$ with $\mt{B_{\infty}}=\colim{\Lambda}\mathcal B_*$, where the colimit is computed in the very large category $\vlCat$ of large categories. 
    We call such a chain an exhaustion of $\mt{B_{\infty}}$ (by presentable categories). 

    \item A \emph{chain of} \emph{topoi} is a chain of presentable categories $\mt{B}_*\colon \Lambda \to \Pr^L$ such that for all $\lambda\in \Lambda$, $\mt{B_{\lambda}}$ is a topos. A \emph{big topos} is a possibly large category $\mt{B_{\infty}}$ such that there exists a chain of topoi $\mt{B}\colon \Lambda\to \Pr^L$ with $\mt{B_{\infty}}=\colim{\Lambda}\mt{B}_*$, where the colimit is computed in the very large category of large categories. 
    We call such a chain an exhaustion of $\mt{B_{\infty}}$ (by topoi).
\end{romanenum}
\end{definition}
\begin{rem}The definitions of big topoi/big presentable categories are rigged so that their categories of spectrum objects behave much like for presentable categories. 
    The left-exactness of the transition functors may seem unnatural, however it is crucial in several arguments, for example to ensure that the stabilization of a big presentable category is stable and big presentable. We also used it to \textit{glue} left adjoints which exist on the level of presentable categories, for instance in the proofs of \cref{truncationbigtopos,existencesuspension}. 
    The fully faithfulness (and left-exactness) of the transition functors seems necessary to obtain a corresponding colimit-description of the category of spectrum objects (\cref{stabilizationbigpresentablecategoriesbigpresentable}).  

    Our discussion of monoidal structures in big presentable categories and module categories therein could be streamlined by restricting to those big presentable categories which admit small colimits (or at least  countable colimits). While this condition is satisfied in all examples relevant to this article, simply adding it to the definition would be inconvenient as it is generally not preserved by passing to categories of spectrum objects. A more suitable alternative would be to work with chains of presentable categories indexed on diagrams which are $\lambda$-filtered for every small regular cardinal $\lambda$ (or at least $\omega_1$-filtered), rather than merely filtered. By \cref{fullyfaithfulnessandpreservationoflimits}, this implies that the resulting big presentable categories admit all small (respectively $\omega_1$-small) colimits, while still encompassing all examples we care about. However, we did not adopt this convention to keep the general framework as flexible as possible. 
\end{rem}
\begin{lemma}\label{hyperaccessilesheavesarebigtopos}    If $(\mathcal C,S)$ is a (hyper)ac\-ces\-si\-ble explicit covering site, the category of accessible $S$-(hyper)\-sheaves on $\mathcal C$ is a big topos. 
\end{lemma}
\begin{proof} 
    Denote by $\Lambda_{\mathcal C}$ the poset of regular cardinals $\lambda$ such that $(\mathcal C,S)$ is $\lambda$-(hyper)ac\-ces\-si\-ble. 
    In \cref{chainpresentablecategoriesaccessiblesheaves}, we described a functor \[\mywidehatshv_S(\mathcal C_{*}, \an)\colon \Lambda_{\mathcal C}\to \Pr^L, \, \lambda\mapsto \mywidehatshv_S(\mathcal C_{\lambda}, \an).\] For all $\lambda\in \Lambda_{\mathcal C}$, $\mywidehatshv_S(\mathcal C_{\lambda}, \an)$ is a topos by \cref{kappaccessibleistopos}, and by \cref{fullyfaithfulnessleftkanextension}, for all $\kappa\leq\lambda\in\Lambda_{\mathcal C}$, 
    $\mywidehatshv_S(\mathcal C_{\kappa}, \an)\to \mywidehatshv_S(\mathcal C_{\lambda}, \an)$ is a fully faithful left adjoint. 
    \cref{kappacondensedkappaacc} and \cref{limitscolimitscanbecomputedfinitestagesheaves} imply that it is also left-exact. 
    It now follows from \cref{accessiblesheavesisbigtopos} that $\mywidehatshv_S(\mathcal C_{*}, \an)$ is an exhaustion for $\mywidehatshvacc{S}(\mathcal C, \an)$ by topoi. 
\end{proof}
More generally, \textit{macrotopoi} \cite[Definition 1.4.10]{barwick2019pyknoticobjectsibasic}, i.e. left-exact accessible localizations of categories of accessible presheaves are big topoi by \cite[Example 1.4.11]{barwick2019pyknoticobjectsibasic}. 
Macrotopoi are arguably a more natural generalization of topoi than big topoi, but the definition of big topoi/big presentable categories captures the essential structure needed to generalize results from presentable categories. 

Before turning to spectrum objects in big presentable categories, we record some basic categorical properties of big presentable categories.  
\begin{lemma}\label{fullyfaithfulnessandpreservationoflimits}
Suppose $\mt{B}_*\colon \Lambda\to \Pr^L$ is an exhaustion for a big presentable category $\mt{B}_{\infty}$. 
\begin{romanenum}
\item For all $\lambda\in\Lambda$, the canonical functor $\topo{B_{\lambda}}\to \mt{B_{\infty}}$ is a fully faithful, left-exact left adjoint.
In particular, $\mt{B}_{\infty}$ is locally small.

\item If $I$ is a category such that for all $\kappa\to\lambda\in\Lambda$, $\mt{B}_{\kappa}\to \mt{B_{\lambda}}$ preserves $I$-indexed limits, then $\topo{B_{\kappa}}\to \mt{B_{\infty}}$ preserves $I$-indexed limits. In particular, $\mt{B_{\infty}}$ has all finite limits. 

\item If $I$ is a category such that for all diagrams $\phi\colon I\to \mt{B_{\infty}}$, there exists $\lambda\in\Lambda$ such that $\phi$ factors over $I\to\mt{B}_{\lambda}\subseteq \mt{B}_{\infty}$, then $\mt{B_{\infty}}$ has all $I$-indexed colimits. 
If in addition, for all $\kappa\to \lambda\in \Lambda$, $\mt{B}_{\kappa}\to\mt{B}_{\lambda}$ preserves $I$-indexed limits, then $\mt{B}_{\infty}$ also has all $I$-indexed limits. 
\end{romanenum}
\end{lemma}
\begin{proof}This is an immediate consequence of \cref{filteredcolimitsofcategoriesappendix}. 
\end{proof}

\begin{lemma}\label{truncationbigtopos}
    Suppose that $\mt{B}_{\infty}$ is a big presentable category. 
    \begin{romanenum}
    \item The full subcategory $\tau_{\leq n}\mt{B}_{\infty}$ of $n$-truncated objects (\cite[Definition 5.5.6.1]{highertopostheory}) is a big presentable category.     The inclusion $\tau_{\leq n}\mt{B}_{\infty}\hookrightarrow \mt{B}_{\infty}$ admits a left adjoint $\tau_{\leq n}$. 

    \item More concretely, if $\mt{B}_*\colon\Lambda\to\Pr^L$ is an exhaustion for $\mt{B}_{\infty}$ by presentable categories, then $\tau_{\leq n}\mt{B}_*\subseteq \mt{B}_{*}$ defines a subfunctor which is an exhaustion for $\tau_{\leq n}\mt{B}_{\infty}$ by presentable categories and for $\lambda\in\Lambda$, \[\tau_{\leq n}\mt{B}_{\lambda}\cong \tau_{\leq n}\mt{B}_{\infty}\times_{\mt{B}_{\infty}}\mt{B}_{\lambda}.\]  
    For all $\lambda\in\Lambda$, the functor $\tau_{\leq n}\colon \mt{B}_{\infty}\colon \mt{B}_{\infty}\to \tau_{\leq n}\mt{B}_{\infty}$ restricts to the truncation functor $\tau_{\leq n}^{\mt{B}_{\lambda}}\colon\mt{B}_{\lambda}\to \tau_{\leq n}\mt{B}_{\lambda}$ of $\mt{B}_{\lambda}$. 
    \end{romanenum}
\end{lemma}
\begin{proof}
     Choose an exhaustion $\mt{B}_{*}\colon\Lambda\to \Pr^L$ of $\mt{B}_{\infty}$ by presentable categories. 
    By \cite[Proposition 5.5.6.16]{highertopostheory}, for all $\kappa\to\lambda\in\Lambda$, the left-exact functor $\mt{B}_{\kappa}\to\mt{B}_{\lambda}$ restricts to a functor $\tau_{\leq n}\mt{B}_{\kappa}\to\tau_{\leq n}\mt{B}_{\lambda}$. 
    We therefore obtain a functor \[\Lambda\to\vlCat, \lambda\mapsto \tau_{\leq n}\mt{B}_{\lambda}.\] 
    By \cite[Proposition 5.5.6.18]{highertopostheory}, for all $\lambda\in\Lambda$, $\tau_{\leq n}\mt{B}_{\lambda}\subseteq\mt{B}_{\lambda}$ admits a left adjoint $\tau_{\leq n}^{\mt{B}_{\lambda}}$, and by \cite[Proposition 5.5.6.28]{highertopostheory}, the Beck-Chevalley transformation 
    \[\tau_{\leq n}^{\mt{B}_{\lambda}}\circ \mt{B}(\kappa\to\lambda)\to (\tau_{\leq n}\mt{B})(\kappa\to\lambda)\circ \tau_{\leq n}^{\mt{B}_{\kappa}}\] coming from the commutative diagram \begin{center}
        \begin{tikzcd}
        \tau_{\leq n}\mt{B}_{\kappa}\arrow[r, hookrightarrow ]\arrow[d] &\mt{B}_{\kappa}\arrow[d]\\ 
        \tau_{\leq n}\mt{B}_{\lambda}\arrow[r, hookrightarrow] &\mt{B}_{\lambda}
        \end{tikzcd}
    \end{center} is an equivalence for all $\kappa\to\lambda\in\Lambda$. 
    It now follows from \cref{adjunctionsbigpresentable} that $\colim{\lambda\in\Lambda}\tau_{\leq n}\mt{B}_{\lambda}\to \mt{B}_{\infty}$ admits a left adjoint $\tau_{\leq n}^{\mt{B}_{\infty}}=\colim{\lambda\in\Lambda}\tau_{\leq n}^{\mt{B}_{\lambda}}$.
    By \cref{filteredcolimitsofcategoriesappendix}, the functor $c\colon \colim{\lambda\in\Lambda}\tau_{\leq n}\mt{B}_{\lambda}\to \mt{B}_{\infty}$ is fully faithful. 
    As for all $\lambda\in \Lambda$, $\mt{B}_{\lambda}\subseteq \mt{B}_{\infty}$ is left-exact (\cref{fullyfaithfulnessandpreservationoflimits}), 
    $\tau_{\leq n}\mt{B}_{\lambda}\subseteq \tau_{\leq n}\mt{B}_{\infty}$ by \cite[Proposition 5.5.6.15]{highertopostheory}. 
    As $\mt{B}_{\lambda}\subseteq \mt{B}_{\infty}$ is a fully faithful left adjoint (\cref{fullyfaithfulnessandpreservationoflimits}), \[(\tau_{\leq n}\mt{B}_{\infty})\cap \mt{B}_{\lambda}\subseteq \tau_{\leq n}\mt{B}_{\lambda}, \] which shows that \[\tau_{\leq n}\mt{B}_{\lambda}\cong \tau_{\leq n}\mt{B}_{\infty}\times_{\mt{B}_{\infty}}\mt{B}_{\lambda}.\] 
    It follows that $c$ factors over an equivalence 
     \[\colim{\lambda\in\Lambda}\tau_{\leq n}\mt{B}_{\lambda}\cong  \tau_{\leq n}\mt{B}_{\infty}.\] 

    It remains to show that $\tau_{\leq n}\mt{B}_{*}$ is a chain of presentable categories. 
    As localization of the presentable category $\mt{B}_{\lambda}$, $\tau_{\leq n}\mt{B}_{\lambda}$ is presentable for all $\lambda\in\Lambda$. 
    For all $\kappa\leq\lambda\in\Lambda$, $\tau_{\leq n}\mt{B}_{\kappa}\to \tau_{\leq n}\mt{B}_{\lambda}$ is fully faithful as restriction of a fully faithful functor. As for all $\mu\in\Lambda$, colimits in $\tau_{\leq n}\mt{B}_{\mu}$ can be computed by applying $\tau_{\leq n}^{\mt{B}_{\mu}}$ to the corresponding colimit in $\mt{B}_{\mu}$ and for $\kappa\leq\lambda\in\Lambda$, \[\tau_{\leq n}^{\mt{B}_{\lambda}}\circ \mt{B}(\kappa\to\lambda)\cong \mt{B}(\kappa\to\lambda)\tau_{\leq n}^{\mt{B}_{\kappa}}\] (see above), cocontinuity of $\mt{B}_{\kappa}\to\mt{B}_{\lambda}$ implies that $\tau_{\leq n}\mt{B}_{\kappa}\to\tau_{\leq n}\mt{B}_{\lambda}$ preserve small colimits for all $\kappa\to\lambda\in\Lambda$. 
    As for all $\mu\in\Lambda$, $\tau_{\leq n}\mt{B}_{\mu}\subseteq \mt{B}_{\mu}$ is closed under small limits, and for all $\kappa\to\lambda\in\Lambda$, $\mt{B}_{\kappa}\to\mt{B}_{\lambda}$ is left-exact, the functors $\tau_{\leq n}\mt{B}_{\kappa}\to\tau_{\leq n}\mt{B}_{\lambda}$ are left-exact for all $\kappa\leq\lambda$. 
\end{proof}

\subsection{Spectrum objects in big presentable categories}\label{section:spectrumobjectsofbigpresentablecategories}
\begin{definition}[Spectrum objects]{\cite[Definition 1.4.2.8]{higheralgebra}}\label{definitionspectrumobjects}
    For a category $\mathcal C$, the category of \emph{spectrum objects} in $\mathcal C$, \[\stab{\cat{C}}\subseteq \Fun(\an^{\operatorname{fin}}_*, \cat{C})\] is the full subcategory on reduced, excisive functors $\an^{\operatorname{fin}}_*\to\mathcal C$, i.e.\ functors $F\colon \an^{\operatorname{fin}}_*\to\mathcal C$ which send pushouts to pullbacks and satisfy $F(*)=*$.

    Evaluation at $S^0\in \an^{\operatorname{fin}}_*$, \begin{align*}\Omega^{\infty}\colon \stab{\cat{C}}&\to\cat{C}\\ F & \mapsto F(S^0)\end{align*} is called \emph{infinite loop space functor}.
\end{definition}
In this section, we show that the category of spectrum objects of a big presentable category behaves in many aspects similar to the stabilization of a presentable category. This leads to well-behaved notions of (group) cohomology in big topoi which we will discuss in \cref{section:cohomologyinabigtopos,section:groupcohomologybigtopos,section:groupcohomologybigtopos}.  
This section is structured as follows: 
We first show that the category of spectrum objects in a big presentable category is stable
(\cref{spectrumobjectsinbigcatsstable}) and big presentable (\cref{stabilizationbigpresentablecategoriesbigpresentable}).
Then we show that the infinite loop space functor admits a left adjoint which factors over the categories of (grouplike) commutative monoids (\cref{existencesuspension}). We use the left adjoint to define a $t$-structure on the category of spectrum objects (\cref{tstructurespectrumobjects}). 

For a presentable category $\mathcal C$, the left adjoint $\mathcal C\to \Sp(\mathcal C)$ of the infinite loop space functor exhibits $\Sp(\mathcal C)$ as the initial stable category in $\Pr^L$ over $\mathcal C$. We establish a similar universal property for big presentable categories (\cref{universalpropertyspectrumobjectsbigpresentable}, \cref{actualuniversalpropertystabilization}). 
In \cref{section:monoidalstructureonspectrumobjects}, we describe conditions on a (symmetric) monoidal structure on a big presentable category $\mathcal C$ under which the left adjoint $\mathcal C\to \stab{\mathcal C}$ of $\Omega^{\infty}$ inherits a (symmetric) monoidal structure. We establish a universal property of this (symmetric) monoidal structure and show that it is compatible with the $t$-structure from \cref{tstructurespectrumobjects}. 
In \cref{section:spectralenrichment}, we show that the category of spectrum objects in a big presentable category is naturally spectrally enriched (\cref{Spmodulestructure}), and that every adjunction between stable, big presentable categories enhances to an adjunction of spectrally enriched categories (\cref{adjunctionsspectrallyenriched}). 
We conclude with basic properties of the \textit{constant sheaf functor} into a big presentable category. 

\begin{lemma}\label{spectrumobjectsinbigcatsstable}
If $\mathcal B$ is a big presentable category, its category of spectrum objects $\stab{\mathcal B}$ is stable. 
\end{lemma}
\begin{proof}
    By \cref{fullyfaithfulnessandpreservationoflimits}, $\mt{B}$ has all finite limits and hence $\stab{\mt{B}}$ is stable by \cite[Corollary 1.4.2.17]{higheralgebra}.
\end{proof}
\begin{definition}\label{definitionstabilization}
If $f\colon \mathcal C\to\mathcal D$ is a finite limits preserving functor, then pushforward along 
\[f_*\colon \Fun(\an^{\operatorname{fin}}_*, \cat{C})\to\Fun(\an^{\operatorname{fin}}_*, \mathcal D)\] restricts to a functor $\Sp(f)\colon \stab{\cat{C}}\to\stab{\cat{D}}$. 
We call $\Sp(f)$ the \emph{stabilization} of $f$.
\end{definition}
If $f$ is fully faithful, so is $f_*$ and hence $\Sp(f)$. 

\begin{ex}\label{geometricmorphismstabilization}
Suppose $L\colon \topo{X}\rightleftarrows \topo{Y}\colon R$ is an adjunction and $L$ preserves finite limits.  
Then pushforward along $L$ and $R$ restricts to functors \[L_{\Sp}\colon\stab{\mathcal X}\to\stab{\mathcal Y}\text{ and } R_{\Sp}\colon \stab{\mathcal Y}\to\stab{\mathcal X}.\] 
An adjunction datum $(\eta\colon \id\to RL, \epsilon\colon LR\to \id)$ yields natural transformations \[ \id\to R_{\Sp}L_{\Sp}, \, \, L_{\Sp}R_{\Sp}\to \id\] which exhibit $L_{\Sp}$ as left adjoint to $R_{\Sp}$.
By construction, 
\begin{equation}\label{diagrambeckchevalleystabilisedadjunction}
\begin{tikzcd}
     \mathcal X\arrow[r,"L"]& \mathcal Y\\ 
      \stab{\topo{Y}}\arrow[u,"\Omega^{\infty}"]\arrow[r,"L_{\Sp}"]& \stab{\topo{X}}\arrow[u,"\Omega^{\infty}"']
\end{tikzcd}
\end{equation} commutes, and by construction of the adjunction datum for $L_{\Sp}\dashv R_{\Sp}$, its mate 
\begin{equation}\label{diagrambeckchevalleystabilisedadjunctionmate1}
\begin{tikzcd}
           \mathcal X & \mathcal Y\arrow[l,"R"']\\ 
          \stab{\topo{X}}\arrow[u,"\Omega^{\infty}"]& \stab{\topo{Y}}\arrow[l,"R_{\Sp}"']\arrow[u,"\Omega^{\infty}"']
        \end{tikzcd}\end{equation} commutes as well.
    \end{ex}
 
\begin{lemma}\label{stabilizationbigpresentablecategoriesbigpresentable}
    Suppose that $\mt{B}_*\colon \Lambda\to\Pr^L$ is an exhaustion for a big presentable category $\mt{B_{\infty}}$. 
    Then $\Fun(\an^{\operatorname{fin}}_*, \mt{B}_{*})\colon \Lambda\to \Cat$ has \begin{align*}\Sp(\mt{B}_*)\colon\Lambda&\to \Cat, \\\lambda&\mapsto \stab{\mt{B}_{\lambda}}\end{align*} as a subfunctor. This functor enhances to a functor $\Sp(\mt{B}_*)\colon\Lambda\to \Pr^L$ which is an exhaustion for $\stab{\mt{B}_{\infty}}$.  
    In particular, $\colim{\lambda\in\Lambda}\stab{\mt{B}_{\lambda}}\cong \stab{\mt{B}_{\infty}}$ is a big presentable category.
\end{lemma}
\begin{proof}
    Since for all $\kappa\to\lambda\in\Lambda$, $\mt{B}_{\kappa}\to\mt{B}_{\lambda}$ preserves finite limits, $\stab{\mt{B}_*}\subseteq \Fun(\an^{\operatorname{fin}}_*, \mt{B}_{*})$ is a subfunctor. 
    We explained above that for all $\kappa\to\lambda\in\Lambda$, $\stab{\mt{B}_{\kappa}}\to\stab{\mt{B}_{\lambda}}$ is a fully faithful left adjoint. Since its source and target are stable, it is also left-exact. 
    By \cite[Example 4.8.1.23]{higheralgebra}, for every presentable category $\topo{X}$, $\Sp(\topo{X})\cong \topo{X}\otimes_{\Pr^L}\Sp$ is presentable, whence $\stab{\mt{B}_*}$ lifts to a functor $\Lambda\to \Pr^L$. 
    
    As for all $\kappa\to\lambda\in\Lambda$, the functors $\stab{\mt{B_{\kappa}}}\to \stab{\mt{B_{\lambda}}}$ and $\stab{\mt{B_{\kappa}}}\to\stab{\mt{B_{\infty}}}$ are fully faithful, \[\colim{\kappa\in\Lambda}\stab{\mt{B_{\kappa}}}\to \stab{\mt{B_{\infty}}}\] is fully faithful by \cref{filteredcolimitsofcategoriesappendix}. Its essential image consists of the reduced, excisive functors $\an^{\operatorname{fin}}_*\to\mt{B_{\infty}}$ which factor over $\mt{B_{\kappa}}\subseteq\mt{B_{\infty}}$ for some $\kappa\in\Lambda$. Suppose $E\colon \an^{\operatorname{fin}}_*\to\mt{B_{\infty}}$ is a reduced and excisive functor and choose $\kappa\in \Lambda$ with $E(S^0)\in \mt{B_{\kappa}}\subseteq\mt{B_{\infty}}$. Since $E$ is excisive and $\mt{B_{\kappa}}\subseteq\mt{B_{\infty}}$ is closed under finite limits, \[\an^{\operatorname{fin}}_*(E, \kappa)\coloneqq \{ X\in \an^{\operatorname{fin}}_*\, | \, E(X)\in \mt{B_{\kappa}}\}\] is closed under pushouts. As $E$ is reduced, $*\in \an^{\operatorname{fin}}_{*}(E, \kappa)$. As $\an^{\operatorname{fin}}_*$ is generated under pushouts by $S^0$ and $*$, this implies that \[\an^{\operatorname{fin}}_*(E, \kappa)=\an^{\operatorname{fin}}_*, \] i.e.\ $E\in  \stab{\mt{B_{\kappa}}}\subseteq \stab{\mt{B_{\infty}}}$.
\end{proof}
For a presentable category $\mt{B}$, the infinite loop space functor has a left adjoint $\Sigma^{\infty}_{+}\colon\mt{B}\to \stab{\mt{B}}$ which exhibits $\stab{\mt{B}}$ as the initial stable category over $\mt{B}$ in $\Pr^L$. This implies the following: 
\begin{lemma}\label{identifystabilizationwithtensorproductwithspectra}If $\mt{B}_*\colon\Lambda\to \Pr^L$ is a chain of presentable categories, there is an equivalence  
    $\Sp\otimes_{\Pr^L}\mt{B}_*\cong \stab{\mt{B}_*}\in\Fun(\Lambda, \Pr^L)$ such that for all $\lambda\in\Lambda$, the functor \[\mt{B}_{\lambda}\cong \an\otimes_{\Pr^L}\mt{B}_{\lambda}\to \Sp\otimes_{\Pr^L}\mt{B}_{\lambda}\cong \stab{\mt{B}_{\lambda}}\] induced by $\Sigma^{\infty}_{+}\colon \an\to \Sp\in\Pr^L$ is left adjoint to \[ \Omega^{\infty}\colon\stab{\mt{B}_{\lambda}}\to\mt{B}_{\lambda}, \, F\mapsto F(S^0).\]  
\end{lemma}
\begin{proof}
Evaluation at $S^0$ induces a natural transformation 
$\Omega^{\infty}_*\colon \stab{\mt{B}_*}\to \mt{B}_*\in\Fun(\Lambda, \vlCat)$, cf.\ \cref{geometricmorphismstabilization}. 
For all $\lambda\in\Lambda$, $\Omega^{\infty}_{\lambda}$ has a left adjoint $\Sigma^{\infty}_{+, \lambda}$ by \cite[Proposition 1.4.4.4]{higheralgebra} and for all $\kappa\to \lambda\in\Lambda$, the mate \begin{center}
\begin{tikzcd}
    \mt{B}_{\kappa}\arrow[rr,"\mt{B}(\kappa\to\lambda)"]\arrow[d,"\Sigma^{\infty}_{+, \kappa}"'] && \mt{B}_{\lambda}\arrow[d,"\Sigma^{\infty}_{+, \lambda}"]\arrow[dll,Rightarrow]\\
    \stab{\mt{B}_{\lambda}}\arrow[rr,"\stab{\mt{B}}(\kappa\to\lambda)"'] && \stab{\mt{B}_{\lambda}} 
\end{tikzcd}
\end{center} of \ref{diagrambeckchevalleystabilisedadjunction} commutes as its opposite \ref{diagrambeckchevalleystabilisedadjunctionmate1} is commutative, cf.\ \cite[Remark 4.7.4.14]{higheralgebra}. 
Whence by \cref{adjunctionsbigpresentable}, the left adjoints assemble into a natural transformation \[\mt{B}_*\to \stab{\mt{B}_*}\in\Fun(\Lambda, \Pr^L).\] 
Since $\stab{\mt{B}_{\lambda}}$ is presentable and stable for all $\lambda\in\Lambda$, by \cite[Proposition 4.8.2.18]{higheralgebra}, this factors uniquely over a natural transformation $\mt{B}_*\to \Sp\otimes_{\Pr^L}\mt{B}_*\to \stab{\mt{B}_*}$ and by \cite[Corollary 1.4.4.5]{higheralgebra}, $\Sp\otimes_{\Pr^L}\mt{B}_{\lambda}\to\stab{\mt{B}_{\lambda}}$ is an equivalence for all $\lambda\in\Lambda$. 
\end{proof}
\begin{cor}\label{stabilizationisspectrumvaluedsheaves}
For a (hyper)ac\-ces\-si\-ble explicit covering site $(\mathcal C,S)$,  \[\stab{\mywidehatshvacc{S}(\mathcal C, \an)}\cong\mywidehatshvacc{S}(\mathcal C, \Sp).\] 
\end{cor}
\begin{proof} Denote by $\Lambda_{\mathcal C}$ the poset of regular cardinals $\kappa$ such that $(\mathcal C,S)$ is $\kappa$-(hyper)accessible and by \[ \mywidehatshv_S(\mathcal C_*,\an),\mywidehatshv_S(\mathcal C_*,\Sp)\colon\Lambda_{\mathcal C}\to \Pr^L\] the functors described in \cref{chainpresentablecategoriesaccessiblesheaves}. 
By \cref{identificationsheavesfunctorial},
\[\Sp\otimes_{\Pr^L}{\mywidehatshv_S(\mathcal C_{*}, \an)}\cong \mywidehatshv_S(\mathcal C_{*}, \Sp)\in \Fun(\Lambda_{\mathcal C}, \vlCat), \] and by \cref{identifystabilizationwithtensorproductwithspectra}, \[\Sp\otimes_{\Pr^L}{\mywidehatshv_S(\mathcal C_{*}, \an)}\cong \stab{\mywidehatshv_S(\mathcal C_{*}, \an)}\in \Fun(\Lambda_{\mathcal C}, \vlCat).\] The statement now follows from \cref{accessiblesheavesisbigtopos} and \cref{stabilizationbigpresentablecategoriesbigpresentable}. 
\end{proof}
\begin{cor}\label{infinitelooppreservesfilteredcolimits}
    \begin{romanenum}
    \item If $\topo{X}$ is a topos, then 
\[\Omega^{\infty-i}\colon \stab{\topo{X}}\to\topo{X}\] preserves filtered colimits for all $i\in\mathbb Z$. 
\item Suppose that $(\mathcal C,S)$ is a (hyper)ac\-ces\-si\-ble explicit covering site. 
For all $i\in\mathbb Z$, 
\[ \Omega^{\infty-i}\colon \stab{\mywidehatshvacc{S}(\mathcal C, \an)}\to\mywidehatshvacc{S}(\mathcal C, \an)\] preserves filtered colimits. 
    \end{romanenum}
\end{cor}
\begin{proof}
    By \cite[Example 7.3.4.7]{highertopostheory}, for a topos $\topo{X}$, \[\Omega\colon\topo{X}_*\to\topo{X}_*\] preserves filtered colimits. By \cite[Remark 1.4.2.25]{higheralgebra}, this implies that $\Omega^{\infty-i}\colon \stab{\topo{X}}\to \topo{X}$ preserves filtered colimits for all $i\in\mathbb Z$. 
    It now follows from \cref{stabilizationisspectrumvaluedsheaves}, \cref{kappacondensedkappaacc}, \cref{kappaccessibleistopos} and \cref{limitscolimitscanbecomputedfinitestagesheaves} that for a (hyper)accessible explicit covering site $(\mathcal C,S)$, \[\Omega^{\infty-i} \colon\stab{\mywidehatshvacc{S}(\mathcal C, \an)}\to\mywidehatshvacc{S}(\mathcal C, \an)\] preserves filtered colimits for all $i\in\mathbb Z$. 
\end{proof}
We now show that for every big presentable category $\mathcal C$, the infinite loop space functor has a left adjoint which factors over the categories of (grouplike) commutative monoids in $\mathcal C$. We will apply this below to describe a $t$-structure on spectrum objects and to prove a recognition principle for big topoi. 
\begin{definition}[{\cite[Definition 1.2]{GepnerGrothNikolaus}}]
    Suppose $\mathcal C$ is a cartesian monoidal category. 
    A commutative monoid $A\in\CMon(\mathcal C)$ is grouplike if the underlying commutative monoid in the homotopy category $\operatorname{ho}(\mathcal C)$ is a group, i.e.\ the shear map $A\times A\to A\times A$ is an isomorphism. 
    Denote by $\CGrp(\mathcal C)$ the full subcategory of $\CMon(\mathcal C)$ on grouplike commutative monoids.
\end{definition}
See \cite[Proposition 1.1]{GepnerGrothNikolaus} for other characterizations of grouplike commutative monoids.  
\begin{lemma}\label{existencesuspension}
    Suppose that $\mt{B}_{\infty}$ is a big presentable category. 
    The functor $\Omega^{\infty}\colon \stab{\mt{B_{\infty}}}\to\mt{B_{\infty}}$ admits a left adjoint $\Sigma^{\infty}_{+}\colon \mt{B_{\infty}}\to \stab{\mt{B}_{\infty}}$ which factors into left adjoints 
    \[ \mt{B_{\infty}}\to \CMon(\mt{B_{\infty}})\to \CGrp(\mt{B_{\infty}})\to \stab{\mt{B_{\infty}}}, \] where   
     \[\mt{B_{\infty}}\to \CMon(\mt{B_{\infty}}), \text{ and } \CMon(\mt{B_{\infty}})\to \CGrp(\mt{B_{\infty}})\] are left adjoint to the forget functors. 
\end{lemma}
\begin{proof}
    This is a straightforward consequence of \cite[Corollary 4.10]{GepnerGrothNikolaus} and \cref{adjunctionsbigpresentable}. 
    Recall from \cite[Corollary 4.7]{GepnerGrothNikolaus} that for a presentable category $\mathcal C$, \[\mathcal C\cong \an\otimes_{\Pr^L}\mathcal C, \, \CMon(\mathcal C)\cong \CMon(\an)\otimes_{\Pr^L}\mathcal C, \text{ and } \CGrp(\mathcal C)\cong \CGrp(\an)\otimes_{\Pr^L}\mathcal C.\] 
    Under these identifications, the functors \[(\an\to \CMon(\an))\otimes_{\Pr^L} \id_{\mathcal C}\colon \mathcal C\to\CMon(\an)\] and \[(\CMon(\an)\to \CGrp(\an))\otimes_{\Pr^L} \operatorname{id}_{\mathcal C}\colon \CMon(\mathcal C)\to\CGrp(\mathcal C)\] are left adjoint to the forget functors $\CGrp(\mathcal C)\to \CMon(\mathcal C)\to \mathcal C$. This follows from \cite[Corollary 4.10, Corollary 4.8]{GepnerGrothNikolaus}.  
    By \cite[Corollary 4.10]{GepnerGrothNikolaus}, $\Sp\otimes_{\Pr^L}-\colon \Pr^L\to \Pr^L$ factors as \[ \Pr^L\xrightarrow{\CMon(\an)\otimes_{\Pr^L}-}\Pr^L\xrightarrow{\CGrp(\an)\otimes_{\Pr^L}-}\Pr^L\xrightarrow{\Sp\otimes_{\Pr^L}-}\Pr^L.\]
    By \cref{identifystabilizationwithtensorproductwithspectra}, we now obtain natural transformations \[\mt{B}_*\to \CMon(\mt{B}_*)\to \CGrp(\mt{B}_*)\to \stab{\mt{B}_*}\in\Fun(\Lambda, \Pr^L)\] such that for all $\lambda\in\Lambda$, \[\mt{B_{\lambda}}\to \CMon(\mt{B}_{\lambda})\to \CGrp(\mt{B}_{\lambda})\to\stab{\mt{B_{\lambda}}}\] is left adjoint to $\Omega^{\infty}$ and $\mt{B_{\lambda}}\to \CMon(\mt{B}_{\lambda}), \CMon(\mt{B}_{\lambda})\to \CGrp(\mt{B}_{\lambda})$ are left adjoint to the forget functors. 
    For $\lambda\to\mu\in\Lambda$, this yields commutative diagrams 
    \begin{center}
        \begin{tikzcd}
                 \stab{\mt{B_{\lambda}}}\arrow[d]&\arrow[l] \CGrp(\mt{B}_{\lambda})\arrow[d] & \CMon(\mt{B}_{\lambda})\arrow[d]\arrow[l] &\mt{B}_{\lambda}\arrow[d]\arrow[l]\\ 
                \stab{\mt{B}_{\mu}} & \CGrp(\mt{B}_{\mu})\arrow[l] & \CMon(\mt{B}_{\mu})\arrow[l]& \mt{B_{\mu}}\arrow[l]
            \end{tikzcd}
        \end{center}
    We claim that their mates 
    \begin{center}
    \begin{tikzcd}
             \stab{\mt{B_{\lambda}}}\arrow[r]\arrow[d] & \arrow[dl, Rightarrow,"\alpha"] \CGrp(\mt{B}_{\lambda})\arrow[d]\arrow[r]&\arrow[dl, Rightarrow,"\beta"] \CMon(\mt{B}_{\lambda})\arrow[d] \arrow[r]& \mt{B}_{\lambda}\arrow[d]\arrow[dl, Rightarrow,"\gamma"]\\ 
            \stab{\mt{B}_{\mu}}\arrow[r] & \CGrp(\mt{B}_{\mu})\arrow[r] & \CMon(\mt{B}_{\mu})\arrow[r] &  \mt{B_{\mu}}
        \end{tikzcd}
    \end{center} commute as well, where the two right horizontal arrows are the forget functors.  
    Since $\mt{B}_{\lambda}\to\mt{B}_{\mu}$ preserves finite limits, 
    \cite[Lemma 6.1]{GepnerGrothNikolaus} implies that $\beta$ and $\gamma$ are equivalences, so it remains to show that $\alpha$ is an equivalence. 
    As the forget functor $f_{\mu}\colon\CGrp(\mt{B}_{\mu})\to\mt{B}_{\mu}$ is conservative, it suffices to show that $f_{\mu}\alpha$ is an equivalence.  
    By the pasting law for mates (\cite[Lemma 2.2.4]{carmeli2022ambidexterity}) and since $\beta, \gamma$ are equivalences, $f_{\mu}\alpha$ is equivalent to the Beck-Chevalley transformation $\alpha\beta\gamma$
    \begin{center}
        \begin{tikzcd}
            \stab{\mt{B}_{\lambda}}\arrow[r] \arrow[d]& {\mt{B}}_{\lambda}\arrow[d]\arrow[dl, Rightarrow]\\ 
            \stab{\mt{B}_{\mu}}\arrow[r] &\mt{B}_{\mu},
        \end{tikzcd}
    \end{center} which is an equivalence by \ref{diagrambeckchevalleystabilisedadjunction} since $\mt{B}_{\lambda}\to\mt{B}_{\mu}$ preserves finite limits. 
    It now follows from \cref{adjunctionsbigpresentable} that after taking the colimit over $\Lambda$ in $\vlCat$, we obtain left adjoints \[\mt{B_{\infty}}\to \colim{\lambda\in \Lambda}\CMon(\mt{B_{\lambda}})\to \colim{\lambda}\CGrp(\mt{B}_{\lambda})  \to \colim{\lambda\in\Lambda}\Sp(\mt{B_{\lambda}}).\] 
    
    By \cref{stabilizationbigpresentablecategoriesbigpresentable}, $\colim{\lambda\in\Lambda}\stab{\mt{B}_{\lambda}}\cong \stab{\mt{B}_{\infty}}$.     
    Since for all $\lambda\in\Lambda$, $\mt{B}_{\lambda}^{\times}\subseteq \mt{B}_{\infty}^{\times}$ is a symmetric monoidal subcategory and $\Lambda$ is filtered, $\colim{\lambda\in \Lambda}\CMon(\mt{B_{\lambda}})\cong \CMon(\mt{B}_{\infty})$.
    As for $\lambda\in\Lambda$, $\operatorname{ho}(\mt{B}_{\lambda})\subseteq \operatorname{ho}(\mt{B}_{\infty})$ is a full cartesian subcategory, this restricts to an equivalence $\colim{\lambda\in \Lambda}\CGrp(\mt{B}_{\lambda})\cong \CGrp(\mt{B}_{\infty}).$ 
    We have therefore constructed left adjoints 
    \[\mt{B}_{\infty}\to \CMon(\mt{B}_{\infty})\to \CGrp(\mt{B}_{\infty})\to \stab{\mt{B}_{\infty}}.\] 
    Their composition is left adjoint to $\colim{\lambda\in\Lambda}(\Omega^{\infty}_{\lambda}\colon \stab{\mt{B}_{\lambda}}\to\mt{B}_{\lambda})=\Omega^{\infty}_{B_{\infty}}$. 
    The functors \[\mt{B}_{\infty}\to \CMon(\mt{B}_{\infty}), \CMon(\mt{B}_{\infty})\to \CGrp(\mt{B}_{\infty})\] are left adjoint to the $\Lambda$-indexed colimits of the forget functors \[\CGrp(\mt{B}_{\lambda})\to \CMon(\mt{B}_{\lambda}), \CMon(\mt{B}_{\lambda})\to\mt{B}_{\lambda}, \] which equal the forget functors. 
\end{proof}

\begin{rem}
    The above proof and \cref{accessiblesheavesisbigtopos} imply that if $(\mathcal C,S)$ is a (hyper)ac\-ces\-si\-ble explicit covering site, then \[\CMon(\mywidehatshvacc{S}(\mathcal C, \an))\cong \mywidehatshvacc{S}(\mathcal C, \CMon(\an))\] and 
\[\CGrp(\mywidehatshvacc{S}(\mathcal C))\cong \mywidehatshvacc{S}(\mathcal C, \CGrp(\an)).\] 
\end{rem}

\begin{lemma}\label{tstructurespectrumobjects}
Suppose that $\mathcal\mt{B}_{\infty}$ is a big presentable category.
Then \[\stab{\mt{B_{\infty}}}_{\leq -1}\coloneqq \{ X\in \stab{\mt{B_{\infty}}}\,|\, \Omega^{\infty}X=0\}\]determines a $t$-structure on $\stab{\mt{B_{\infty}}}$.
Its connective part is generated by the essential image of $\Sigma^{\infty}_{+}$ under small colimits and extensions, and its coconnective part equals \[\stab{\mt{B}_{\infty}}_{\leq 0}=\{ X\in\stab{\mt{B}_{\infty}}\, |\, \Omega^{\infty}(X)\in\tau_{\leq 0}\mt{B}_{\infty}\}.\] 
If $\mt{B}_*\colon \Lambda\to \Pr^L$ is an exhaustion of $\mt{B}_{\infty}$ by presentable categories, for all $\lambda\in\Lambda$, 
\[ \stab{\mt{B}_{\lambda}}_{\geq 0}=\stab{\mt{B}_{\lambda}}\cap \stab{\mt{B}_{\infty}}_{\geq 0}\text{ and } \stab{\mt{B}_{\lambda}}_{\leq 0}=\stab{\mt{B}_{\lambda}}\cap \stab{\mt{B}_{\infty}}_{\leq 0}.\] 
\end{lemma}
\begin{proof}
    Denote by $\stab{\mt{B}_{\infty}}_{\geq 0}\subseteq \stab{\mt{B}_{\infty}}$ the subcategory generated by the essential image of $\Sigma^{\infty}_{+}$ under small colimits and extensions.
    As for $X\in \stab{\topo{B_{\infty}}}$, \[L_X\coloneqq \{ T\in \stab{\mt{B}_{\infty}}\, | \, \Map_{\stab{\mt{B_{\infty}}}}(T,X)=0\}\subseteq \stab{\topo{B}_{\infty}}\] is closed under small colimits and extensions, for $X\in\stab{\mt{B}_{\infty}}_{\leq -1}$ and $T\in\stab{\mt{B}_{\infty}}_{\geq 0}$,
    \begin{align}\label{orthogonalitytstructure}\Map_{\stab{\mt{B}_{\infty}}}(T,X)=0.\end{align}  
    Choose an exhaustion $\mt{B}_*\colon\Lambda\to \Pr^L$ of $\mt{B}_{\infty}$ by presentable categories.
    As for all $\lambda\in\Lambda$, $\Omega^{\infty}$ restricts to $\Omega^{\infty}_{\mt{B}_{\lambda}}\colon\stab{\mt{B}_{\lambda}}\to\mt{B}_{\lambda}$, \[\stab{\mt{B}_{\lambda}}_{\leq -1}\coloneqq\stab{\mt{B}_{\infty}}_{\leq -1}\cap \stab{\mt{B}_{\lambda}}\] is the $(-1)$-coconnective part of a $t$-structure on $\stab{\mt{B}_{\lambda}}$ by \cite[Proposition 1.4.3.4]{higheralgebra}. 
    We claim that its connective part $\stab{\mt{B}_{\lambda}}_{\geq 0}$ equals $\stab{\mt{B}_{\infty}}_{\geq 0}\cap \stab{\mt{B}_{\lambda}}$. As $\stab{\mt{B}_{\lambda}}\subseteq \stab{\mt{B}_{\infty}}$ is a full subcategory, \ref{orthogonalitytstructure} implies that 
    $\stab{\mt{B}_{\lambda}}\cap \stab{\mt{B}_{\infty}}_{\geq 0}\subseteq \stab{\mt{B}_{\lambda}}_{\geq 0}$. The proof of \cite[Proposition 1.4.3.4]{higheralgebra} shows that $\stab{\mt{B}_{\lambda}}_{\geq 0}$ is generated under small colimits and extensions by the essential image of $\Sigma^{\infty}_{+, \mt{B}_{\lambda}}$. 
    By construction, $\stab{\mt{B}_{\lambda}}\cap \stab{\mt{B}_{\infty}}_{\geq 0}$ contains the essential image of $\Sigma^{\infty}_{+, \mt{B}_{\lambda}}$. Since $\stab{\mt{B}_{\lambda}}\subseteq \stab{\mt{B}_{\infty}}$ is a stable subcategory closed under small colimits, it is also closed under extensions, which implies that $\stab{\mt{B}_{\lambda}}_{\geq 0}=\stab{\mt{B}_{\lambda}}\cap \stab{\mt{B}_{\infty}}_{\geq 0}$. 

    We now deduce from this the statement of the lemma. 
    As $\stab{\mt{B}_{\infty}}_{\geq 0}\subseteq \stab{\mt{B}_{\infty}}$ is closed under small colimits and extensions, the category
    \[ \{X\in \stab{\mt{B}_{\infty}}_{\geq 0}\, | \, \Sigma X\in \stab{\mt{B}_{\infty}}_{\geq 0}\}\] is closed under small colimits and extensions. It contains the essential image of $\Sigma^{\infty}_{+}$, which implies that \[\Sigma(\stab{\mt{B}_{\infty}}_{\geq 0})\subseteq \stab{\mt{B}_{\infty}}_{\geq 0}.\] Clearly, $\Omega(\stab{\mt{B}_{\infty}}_{\leq 0})\subseteq \stab{\mt{B}_{\infty}}_{\leq 0}$. 
    For $x\in\mt{B}_{\infty}$ choose $\lambda\in\Lambda$ with $x\in\mt{B}_{\lambda}$ and a fiber sequence $x_{+}\to x\to x_{-}\in\stab{\mt{B}_{\lambda}}$ with $x_{+}\in \stab{\mt{B}_{\lambda}}_{\geq 0}, x_{-}\in\stab{\mt{B}_{\lambda}}_{\leq -1}$. 
    Since $\stab{\mt{B}_{\lambda}}\subseteq \stab{\mt{B}_{\infty}}$ is exact, $x_{+}\to x\to x_{-}$ is a fiber sequence in $\stab{\mt{B}_{\infty}}$, and by the above, $x_{+}\in \stab{\mt{B}_{\infty}}_{\geq 0}, x_{-}\in\stab{\mt{B}_{\infty}}_{\leq -1}$.  This shows that $(\stab{\mt{B}_{\infty}}_{\geq 0}, \stab{\mt{B}_{\infty}}_{\leq -1})$ constitutes a $t$-structure. 
    We have shown above that for all $\lambda\in\Lambda$, 
    \[ \stab{\mt{B}_{\lambda}}_{\geq 0}=\stab{\mt{B}_{\lambda}}\cap \stab{\mt{B}_{\infty}}_{\geq 0}\text{ and } \stab{\mt{B}_{\lambda}}_{\leq -1}=\stab{\mt{B}_{\lambda}}\cap \stab{\mt{B}_{\infty}}_{\leq -1}.\] 
    As \[\stab{\mt{B}_{\lambda/\infty}}_{\leq 0}=\{ X\in\stab{\mt{B}_{\lambda/\infty}}\,|\, \Omega X\in\stab{\mt{B}_{\lambda/\infty}}_{\leq -1}\}\] and $\stab{\mt{B}_{\lambda}}\subseteq \stab{\mt{B}_{\infty}}$ is a stable subcategory, this implies that \[\stab{\mt{B}_{\lambda}}_{\leq 0}=\stab{\mt{B}_{\lambda}}\cap \stab{\mt{B}_{\infty}}_{\leq 0}.\] 
    For $X\in\stab{\mt{B}_{\infty}}$, 
    \[\mathcal L_{X}^{+}\coloneqq \{ B\in\stab{\mt{B}_{\infty}}\, |\, \Map_{\stab{\mt{B}_{\infty}}}(B,X)\in\tau_{\leq 0}\an\}\] is closed under small colimits and extensions, whence $X\in\stab{\mt{B}_{\infty}}_{\leq 0}$ if and only if $\Sigma^{\infty}_{+}(\mt{B}_{\infty})\subseteq \mathcal L_{X}^{+}$.  
    As $\Sigma^{\infty}_{+}\dashv \Omega^{\infty}$, it follows that \[\stab{\mt{B}_{\infty}}_{\leq 0}=\{ X\in\stab{\mt{B}_{\infty}}\, |\, \Omega^{\infty}(X)\in\tau_{\leq 0}\mt{B}_{\infty}\}.\qedhere\] 
\end{proof}
\begin{cor}\label{recognitionprinciplebigtopoi}
If $\topo{B}_{\infty}$ is a big topos, the left adjoint $\CGrp(\mt{B}_{\infty})\to \stab{\mt{B}_{\infty}}$ from \cref{existencesuspension} factors over an equivalence 
\[ \CGrp(\mt{B}_{\infty})\cong \stab{\mt{B}_{\infty}}_{\geq 0}.\]
This restricts to an equivalence $\Ab(\tau_{\leq 0}\mt{B_{\infty}})\cong \stab{\mt{B_{\infty}}}^{\heart}$ between the 1-category of abelian group objects in $\tau_{\leq 0}\mt{B_{\infty}}$ and the heart of the $t$-structure.
\end{cor}
\begin{proof}
    Choose an exhaustion $\mt{B}_{*}\colon\Lambda\to \Pr^L$ of $\mt{B}_{\infty}$ by topoi. 
    By \cite[Proposition 3.8]{haine2025standardtstructures}, for all $\lambda\in\Lambda$ the left adjoint $\CGrp(\mt{B}_{\lambda})\to \stab{\mt{B}_{\lambda}}$ is fully faithful with essential image $\stab{\mt{B}_{\lambda}}_{\geq 0}$. We thank Lucas Piessevaux for pointing out this reference to us. 
    It now follows from \cref{filteredcolimitsofcategoriesappendix} and the construction of the left adjoint $\CGrp(\mt{B}_{\infty})\to \stab{\mt{B}_{\infty}}$ that it is fully faithful. 
    Since for all $\lambda\in\Lambda$, $\stab{\mt{B}_{\lambda}}_{\geq 0}\cong \stab{\mt{B}_{\lambda}}\cap \stab{\mt{B}_{\infty}}_{\geq 0}$, its essential image equals $\stab{\mt{B}_{\infty}}_{\geq 0}$. 

    As \[\stab{\mt{B}_{\infty}}^{\heart}=\{ X\in\stab{\mt{B}_{\infty}}_{\geq 0}\, |\, \Omega^{\infty}(X)\in \tau_{\leq 0}\mt{B}_{\infty}\}, \] the equivalence $\CGrp(\mt{B}_{\infty})\cong \stab{\mt{B}_{\infty}}_{\geq 0}$ restricts to an equivalence 
    \[\CGrp(\mt{B}_{\infty})\times_{\mt{B}_{\infty}}\tau_{\leq 0}\mt{B}_{\infty}\cong \stab{\mt{B}_{\infty}}^{\heart}, \] and since $\tau_{\leq 0}\mt{B}_{\infty}\subseteq \mt{B}_{\infty}$ is closed under limits, 
    \[\CGrp(\mt{B}_{\infty})\times_{\mt{B}_{\infty}}\tau_{\leq 0}\mt{B}_{\infty}\cong \CGrp(\tau_{\leq 0}\mt{B}_{\infty}).\qedhere\]
\end{proof}
\begin{definition}\label{defeilenbergmaclanefunctorhomotopygroups}
    For a big topos $\mt{B_{\infty}}$, the functor
    \[H\colon \Ab(\tau_{\leq 0}\mt{B_{\infty}})\cong \stab{\mt{B_{\infty}}}^{\heart}\hookrightarrow \stab{\mt{B_{\infty}}}\] is called \emph{Eilenberg-Mac Lane functor}. 
\end{definition}

In addition to \cref{geometricmorphismstabilization}, stabilization also commutes with left-exact left adjoints: 
\begin{cor}\label{stabilizationcommuteswithsuspension}
Suppose that $L\colon \topo{X}\rightleftarrows \topo{Y}\colon R$, $L\dashv R$ is an adjunction between big presentable categories such that $L$ preserves finite limits. 
Then $\Sigma^{\infty}_{+}L\cong L_{\Sp}\Sigma^{\infty}_{+}$ via the Beck-Chevalley transformation of the commutative diagram 
\begin{center}
    \begin{tikzcd}
    \topo{X}\arrow[r,"L"]& \topo{Y}\\ 
    \stab{\topo{X}}\arrow[r,"L_{\Sp}"]\arrow[u,"\Omega^{\infty}"] & \stab{\topo{Y}}\arrow[u,"\Omega^{\infty}"']
    \end{tikzcd}. 
\end{center}
\end{cor}
\begin{proof}
By \cite[Remark 4.7.4.14]{higheralgebra}, the Beck-Chevalley transformation is an equivalence if and only if the opposite Beck-Chevalley transformation 
$\Omega^{\infty}R_{\Sp}\to R\Omega^{\infty}$ is, which holds by construction of the adjunction $L_{\Sp}\dashv R_{\Sp}$.  
\end{proof}
\begin{lemma}\label{texactnessgeometricmorphism}Suppose $\topo{X}, \topo{Y}$ are big presentable categories. 
\begin{romanenum}
\item If $F\colon \topo{X}\to\topo{Y}$ is a finite limits preserving functor, its stabilization $F_{\Sp}\colon\stab{\topo{X}}\to \stab{\topo{Y}}$ is left t-exact, i.e. $F_{\Sp}$ is exact and $F_{\Sp}(\stab{\topo{X}}_{\leq 0})\subseteq \stab{\topo{Y}}_{\leq 0}$. 
\item If $L\colon\topo{X}\to\topo{Y}$ is a left-exact left adjoint, its stabilization $L_{\Sp}\colon \stab{\topo{X}}\to\stab{\topo{Y}}$ is t-exact. 
\end{romanenum}
\end{lemma}
\begin{proof}By \cite[Proposition 5.5.6.16]{highertopostheory}, a finite limits preserving functor $F\colon\topo{X}\to\topo{Y}$ preserves $0$-truncated objects. 
By definition of the stabilization $F_{\Sp}$, $F\circ \Omega^{\infty}=\Omega^{\infty}F_{\Sp}$. The characterization of the coconnective part of the $t$-structure from \cref{tstructurespectrumobjects} now implies that \[F_{\Sp}(\stab{\topo{X}}_{\leq 0})\subseteq \stab{\topo{Y}}_{\leq 0}.\]
Since finite limits in \[\stab{\mathcal X}\subseteq \Fun(\an^{\operatorname{fin}}_*,\mathcal X) \text{ and } \stab{\mathcal Y}\subseteq \Fun(\an^{\operatorname{fin}}_*,\mathcal Y)\] are computed pointwise, $F_{\Sp}$ preserves finite limits. As $\stab{\mathcal X}$ and $\stab{\mathcal Y}$ are stable (\cref{spectrumobjectsinbigcatsstable}), it follows from \cite[Proposition 1.1.4.1]{higheralgebra} that $F_{\Sp}$ is exact. 

Suppose now that $L\colon\topo{X}\to\topo{Y}$ is a left-exact left adjoint. 
Since $L_{\Sp}$ is a left adjoint (\cref{geometricmorphismstabilization}) and $\stab{\topo{Y}}_{\geq 0}\subseteq \stab{\topo{Y}}$ is closed under small colimits and extensions, 
    \[\mathcal L_{\geq 0}\coloneqq \{ X\in \stab{\topo{X}}\, | \, L_{\Sp}(X)\in\stab{\topo{Y}}_{\geq 0}\}\subseteq \stab{\topo{X}}\] is closed under small colimits and extensions. As $L_{\Sp}\circ \Sigma^{\infty}_{+}\cong \Sigma^{\infty}_{+}L$ (\cref{stabilizationcommuteswithsuspension}), $\mathcal L_{\geq 0}$ contains the essential image of $\Sigma^{\infty}_{+}\colon \topo{X}\to \stab{\topo{X}}$, which implies that $\stab{\topo{X}}_{\geq 0}\subseteq \mathcal L_{\geq 0}$.
    It now follows from the above that $L_{\Sp}$ is $t$-exact. 
\end{proof}

For a presentable category $\mt{C}$, the stabilization functor $\mt{C}\to\stab{\mt{C}}$ is the initial cocontinuous functor into a presentable stable category. 
Our next goal is to show that for a big presentable category $\mt{B}_{\infty}$, the \textit{stabilization} \[\Sigma^{\infty}_{+}\colon\mt{B}_{\infty}\to \stab{\mt{B}_{\infty}}\] has a similar universal property, see \cref{universalpropertyspectrumobjectsbigpresentable}, \cref{actualuniversalpropertystabilization}. 

An object $b$ in a stable category $\mt{B}$ is a generator if $\pi_0\Map_{\mt{B}}(b,T)=0$ implies that $T=0$. By \cite[Corollary 1.4.4.2]{higheralgebra}, every stable, presentable category admits a generator.
\begin{lemma}\label{presentablecategoriessinglegenerator}
    Suppose that $\mt{B}$ is a presentable, stable category and $b\in\mt{B}$ is a generator. 
    If $\mathcal C\subseteq \mt{B}$ is a stable, full subcategory closed under small colimits with $b\in \mathcal C$, then the inclusion $\mathcal C\subseteq \mt{B}$ is an equivalence.    
\end{lemma}
\begin{proof}
    Suppose that $\mathcal C\subseteq \mt{B}$ is a stable, full subcategory closed under small colimits and $b\in\mathcal C$. 
    Then $\mathcal C$ is presentable by \cite[Corollary 1.4.4.2]{higheralgebra}, whence by the adjoint functor theorem \cite[Corollary 5.5.2.9]{highertopostheory}, $l\colon \mathcal C\to\mathcal D$ admits a right adjoint $r$. 
    Since $l\colon \mathcal C\to\mathcal D$ is fully faithful, the unit is an equivalence. 
    For $x\in \mathcal D$, 
    \[\Map_{\mathcal D}(lb,lrx)\cong \Map_{\mathcal C}(b,rx)\cong \Map_{\mathcal D}(lb,x), \] whence 
    $\Map_{\mathcal D}(lb, \Fib(lrx\to x))=0$. 
    As $lb$ is a generator, this implies that $\Fib(lrx\to x)=0$, which shows that the counit $lr\to \id$ is an equivalence. 
\end{proof}
\begin{cor}\label{presentablestablesubcategorycontainedinfinitestage}
Suppose that $T$ is a stable, presentable category, $\mt{B}_{\infty}$ is a big presentable category, and $F\colon T\to \mt{B}_{\infty}$ is a cocontinuous, left-exact functor. 
If $\mt{B}_*\colon\Lambda\to \Pr^L$ is an exhaustion of $\mt{B}_{\infty}$ by presentable categories, there exists $\lambda\in\Lambda$ such that $F$ factors over $T\to{\mt{B}_{\lambda}}\subseteq {\mt{B}_{\infty}}$.  
\end{cor}
\begin{proof}For $\lambda\in\Lambda$ denote by $\overline{{\mt{B}_{\lambda}}}$ the essential image of $\mt{B}_{\lambda}\to \mt{B}_{\infty}$. Choose a generator $t\in T$ and $\lambda\in\Lambda$ with $F(t)\in \overline{{\mt{B}_{\lambda}}}$.
Since $\overline{{\mt{B}_{\lambda}}}\subseteq {\mt{B}_{\infty}}$ is closed under small colimits and finite limits and $F$ preserves finite limits and small colimits, 
    \[ \{ x\in T\, |\, F(x)\in \overline{\mt{B}_{\lambda}}\}\subseteq T\] is a stable subcategory closed under small colimits. 
    \cref{presentablecategoriessinglegenerator} now implies that $F$ factors over $\mt{B}_{\lambda}\subseteq \mt{B}_{\infty}$. 
\end{proof}
\cref{presentablecategoriessinglegenerator} also implies that the category of spectrum objects $\stab{\mt{B}_{\infty}}$ of a big presentable category $\mt{B}_{\infty}$ is as a stable category generated by the image of $\Sigma^{\infty}_{+}\colon \mt{B}_{\infty}\to\stab{\mt{B}_{\infty}}$ under small colimits:
\begin{cor}\label{suspensionspectragenerate}
Suppose that $\mt{B}_{\infty}$ is a big presentable category.  
If $\mathcal C\subseteq \stab{\mt{B}_{\infty}}$ is a stable subcategory closed under small colimits and equivalences such that $\essim(\Sigma^{\infty}_{+})\subseteq \mathcal C$, then \[\mathcal C=\stab{\mt{B}_{\infty}}.\] 
\end{cor}
\begin{proof}
    Choose an exhaustion $\mt{B}_{*}\colon\Lambda\to \Pr^L$ of $\mt{B}_{\infty}$, and for $\lambda\in\Lambda$ denote by $i_{\lambda}\colon\stab{\mt{B}_{\lambda}}\to \stab{\mt{B}_{\infty}}$ the canonical functor. 
    As $i_{\lambda}$ is a left adjoint and $\mathcal C\subseteq \stab{\mt{B}_{\infty}}$ is a stable subcategory closed under small colimits, 
    \[\tilde{\mt{B}_{\lambda}}\coloneqq \stab{\mt{B}_{\lambda}}\times_{\stab{\mt{B}_{\infty}}}\mathcal C\cong  \{ b\in \stab{\mt{B}_{\lambda}}\, |\, i_{\lambda}(b)\in\mathcal C\}\subseteq \stab{\mt{B}_{\lambda}}\] is a stable subcategory closed under small colimits. 
    Since $\essim(\Sigma^{\infty}_{+})\subseteq \mathcal C$, the proof of \cite[Corollary 1.4.4.2]{higheralgebra} implies that $\tilde{\mt{B}_{\lambda}}\subseteq \mt{B}_{\lambda}$ contains a generator of $\stab{\mt{B}_{\lambda}}$. Hence by \cref{presentablecategoriessinglegenerator}, $\stab{\mt{B}_{\lambda}}\times_{\stab{\mt{B}_{\infty}}}\mathcal C\to \stab{\mt{B}_{\lambda}}$ is an equivalence for all $\lambda\in\Lambda$, which by \cref{stabilizationbigpresentablecategoriesbigpresentable} implies that the inclusion $\mathcal C\to \stab{\mt{B}_{\infty}}$ is an equivalence. 
    As $\mathcal C$ is closed under equivalences, it follows that $\mathcal C=\stab{\mt{B}_{\infty}}$. 
\end{proof}
For potentially large categories $\mathcal B, \mathcal C$ denote by $\Fun^{\operatorname{colim}}(\mathcal B, \mathcal C)\subseteq \Fun(\mathcal B, \mathcal C)$ the full subcategory of small colimits preserving functors. 
The above implies that the \textit{stabilization} 
\[ \Sigma^{\infty}_{+}\colon \mt{B}_{\infty}\to\stab{\mt{B}_{\infty}}\] of a big presentable category $\mt{B}_{\infty}$ almost has a universal property: 
\begin{cor}\label{universalpropertyspectrumobjectsbigpresentable}
Suppose $\mt{B}_{\infty}, \mt{C}_{\infty}$ are two big presentable categories and $\mt{C}_{\infty}$ is stable. 

Pullback along $\Sigma^{\infty}_{+}\colon \mt{B}_{\infty}\to \stab{\mt{B}_{\infty}}$ defines a fully faithful functor 
\[ \Fun^{\operatorname{colim}}(\stab{\mt{B}_{\infty}}, \mt{C}_{\infty})\hookrightarrow \Fun^{\operatorname{colim}}(\mt{B}_{\infty}, \mt{C}_{\infty}).\] 
If $\mt{B}_{\infty}$ and $\mt{C}_{\infty}$ admit small coproducts, then this is an equivalence. 
\end{cor}

\begin{proof}
    Choose an exhaustion $\mt{B}_*\colon\Lambda\to \Pr^L$ of $\mt{B}_{\infty}$ by presentable categories. 
    Since $\mt{B}_{\lambda}\subseteq \mt{B}_{\infty}$ and $\stab{\mt{B}_{\lambda}}\subseteq \stab{\mt{B}_{\infty}}$ are closed under small colimits for all $\lambda\in\Lambda$ (\cref{stabilizationbigpresentablecategoriesbigpresentable} and \cref{fullyfaithfulnessandpreservationoflimits}), 
    \cref{stabilizationbigpresentablecategoriesbigpresentable} and the construction of $\Sigma^{\infty}_{+}$ \cref{existencesuspension} imply that 
\[(\Sigma^{\infty}_{+})^*\colon \Fun^{\operatorname{colim}}(\stab{\mt{B}_{\infty}}, \mt{C})\to \Fun(\mt{B_{\infty}}, \mt{C}_{\infty})\] factors as 
\[ \Fun^{\operatorname{colim}}(\stab{\mt{B}_{\infty}}, \mt{C}_{\infty})\hookrightarrow \clim{\lambda\in\Lambda}\Fun^{\operatorname{colim}}(\stab{\mt{B}_{\lambda}}, \mt{C}_{\infty})\to \clim{\lambda\in\Lambda}\Fun^{\operatorname{colim}}(\mt{B}_{\lambda}, \mt{C}_{\infty})\hookrightarrow \Fun(\mt{B}_{\infty}, \mt{C}_{\infty})\] where the left functor is induced by the restrictions along $\stab{\mt{B}_{\lambda}}\subseteq \stab{\mt{B}_{\infty}}$, the middle functor is induced by pullback along the functors $\Sigma^{\infty}_{+}\colon\mt{B}_{\lambda}\to \stab{\mt{B}_{\lambda}}$, and the right functor is induced by the inclusions $\Fun^{\operatorname{colim}}(\mt{B}_{\lambda}, \mt{C}_{\infty})\subseteq \Fun(\mt{B}_{\lambda}, \mt{C}_{\infty})$. 
The left and right functor are fully faithful: $\clim{\lambda\in\Lambda}\Fun^{\operatorname{colim}}(\stab{\mt{B}_{\lambda}}, \mt{C}_{\infty})$ is equivalent to the full subcategory of functors $F\colon \stab{\mt{B}_{\infty}}\to\mt{C}_{\infty}$ such that for all $\lambda\in\Lambda$, ${\stab{\mt{B}_{\lambda}}}\to\stab{\mt{B}_{\infty}}\to \mathcal C_{\infty}$ is cocontinuous, and $\clim{\lambda\in\Lambda}\Fun^{\operatorname{colim}}(\mt{B}_{\lambda}, \mt{C}_{\infty})$ is equivalent to the full subcategory of functors $F\colon \mt{B}_{\infty}\to\mt{C}_{\infty}$ such that for all $\lambda\in\Lambda$, $F|_{{\mt{B}_{\lambda}}}$ is cocontinuous, cf.\ \cref{limitsincat}. 
We claim that for all $\lambda\in\Lambda$, 
\[(\Sigma^{\infty}_{+})^*\colon \Fun^{\operatorname{colim}}(\stab{\mt{B}_{\lambda}}, \mt{C}_{\infty})\to \Fun^{\operatorname{colim}}(\mt{B}_{\lambda}, \mt{C}_{\infty})\] is fully faithful.   
Choose an exhaustion $\mt{C}_*\colon F\to \Pr^L$ of $\mt{C}_{\infty}$. 
Since for all $f\in F$, $\mt{C}_{f}\subseteq \mt{C}_{\infty}$ is a presentable category closed under finite limits and colimits, it is a stable presentable category. \cite[Corollary 1.4.4.5]{higheralgebra} implies that for all $\lambda\in\Lambda$, pullback along $\Sigma^{\infty}_{+}\colon \mt{B}_{\lambda}\to \stab{\mt{B}_{\lambda}}$ induces an equivalence
\[ \colim{f\in F}\Fun^{\operatorname{colim}}(\stab{\mt{B}_{\lambda}}, \mt{C}_f)\cong \colim{f\in F}\Fun^{\operatorname{colim}}(\mt{B}_{\lambda}, \mt{C}_f).\]
Since for all $f\to g\in F$, $\mt{C}_f\to \mt{C}_g$ is fully faithful and cocontinuous, \[\colim{f\in F}\Fun^{\operatorname{colim}}(\mt{B}_{\lambda}, \mt{C}_f)\to \Fun^{\operatorname{colim}}(\mt{B}_{\lambda}, \mt{C}_{\infty})\] and  
\[\colim{f\in F}\Fun^{\operatorname{colim}}(\stab{\mt{B}_{\lambda}}, \mt{C}_f)\to \Fun^{\operatorname{colim}}(\stab{\mt{B}_{\lambda}}, \mt{C}_{\infty})\] are fully faithful by \cref{filteredcolimitsofcategoriesappendix}. 
By \cref{presentablestablesubcategorycontainedinfinitestage}, \[\colim{f\in F}\Fun^{\operatorname{colim}}(\stab{\mt{B}_{\lambda}}, \mt{C}_f)\to \Fun^{\operatorname{colim}}(\stab{\mt{B}_{\lambda}}, \mt{C}_{\infty})\] is essentially surjective and hence an equivalence, which shows that 
\[(\Sigma^{\infty}_{+})^*\colon \Fun^{\operatorname{colim}}(\stab{\mt{B}_{\lambda}}, \mt{C}_{\infty})\to \Fun^{\operatorname{colim}}(\mt{B}_{\lambda}, \mt{C}_{\infty})\] is fully faithful for all $\lambda\in\Lambda$.  
This proves that \[(\Sigma^{\infty}_{+})^*\colon \Fun^{\operatorname{colim}}(\stab{\mt{B}_{\infty}}, \mt{C})\to \Fun(\mt{B_{\infty}}, \mt{C}_{\infty})\] is fully faithful and factors as 
\begin{align*}
\Fun^{\operatorname{colim}}(\stab{\mt{B}_{\infty}}, \mt{C}_{\infty})\hookrightarrow \clim{\lambda\in\Lambda}\Fun^{\operatorname{colim}}(\stab{\mt{B}_{\lambda}}, \mt{C}_{\infty})& \cong \clim{\lambda\in\Lambda}\colim{f\in F}\Fun^{\operatorname{colim}}(\stab{\mt{B}_{\lambda}}, \mt{C}_f)\\ & \cong \clim{\lambda\in\Lambda}\colim{f\in F}\Fun^{\operatorname{colim}}({\mt{B}_{\lambda}}, \mt{C}_f)\\ 
& \hookrightarrow \clim{\lambda\in\Lambda}\Fun^{\operatorname{colim}}(\mt{B}_{\lambda}, \mt{C}_{\infty})\\ 
&\hookrightarrow \Fun(\mt{B}_{\infty}, \mt{C}_{\infty}).
\end{align*}

Suppose now that $\stab{\mt{B}_{\infty}}$ and $\mt{C}_{\infty}$ admit small coproducts.
We first want to show that \begin{align}\label{cocontinuitycanbecheckedlevelwise} \Fun^{\operatorname{colim}}(\stab{\mt{B}_{\infty}}, \mt{C}_{\infty})\cong \clim{\lambda\in\Lambda}\Fun^{\operatorname{colim}}(\stab{\mt{B}_{\lambda}}, \mt{C}_{\infty}), \end{align} i.e.\ that a functor $G\colon \stab{\mt{B}_{\infty}}\to\mt{C}_{\infty}$ is cocontinuous if for all $\lambda\in\Lambda$, $G|_{\stab{\mt{B}_{\lambda}}}$ is cocontinuous. 
 Suppose that $T\colon I\to \stab{\mt{B}_{\infty}}$ is a small diagram and choose $\lambda\in\Lambda$ with \[\sqcup_{i\in I}T(i)\in \essim(\stab{\mt{B}_{\lambda}}\to \stab{\mt{B}_{\infty}}).\] 
 Since $\essim(\stab{\mt{B}_{\lambda}}\to \stab{\mt{B}_{\infty}})\subseteq \stab{\mt{B}_{\infty}}$ is a stable subcategory closed under small colimits, for all $i\in I$, \[T(i)=\Cofib(\sqcup_{\substack{j\in I\\ j\neq i}}T(i)\to \sqcup_{{j\in I}}T(j))\in  \essim(\stab{\mt{B}_{\lambda}}\to \stab{\mt{B}_{\infty}}), \] i.e.\ $T$ factors over $\stab{\mt{B}_{\lambda}}\subseteq \stab{\mt{B}_{\infty}}$ and we can compute $\colim{}^{\stab{\mt{B_{\infty}}}}T$ in $\stab{\mt{B}_{\lambda}}\subseteq \stab{\mt{B}_{\infty}}$, which shows \ref{cocontinuitycanbecheckedlevelwise}. 
 As a consequence, 
 \[\clim{\lambda\in\Lambda}\colim{f\in F}\Fun^{\operatorname{colim}}(\mt{B}_{\lambda}, \mt{C}_{f})\to \Fun(\mt{B}_{\infty}, \mt{C_{\infty}})\] factors over a fully faithful functor $\clim{\lambda\in\Lambda}\colim{f\in F}\Fun^{\operatorname{colim}}(\mt{B}_{\lambda}, \mt{C}_{f})\to \Fun^{\operatorname{colim}}(\mt{B}_{\infty}, \mt{C_{\infty}}).$ 
 We claim that this functor is essentially surjective, i.e.\ for every small colimits preserving functor $G\colon \mt{B}_{\infty}\to \mt{C}_{\infty}$ and $\lambda\in\Lambda$, there exists $f(\lambda)\in F$ such that $F|_{\mt{B}_{\lambda}}$ factors over $\mt{C}_f\subseteq \mt{C}_{\infty}$. 
 Since $\mt{B}_{\lambda}$ is presentable, there exists a regular cardinal $\kappa$ and a small set of $\kappa$-compact objects $(b_i)_{i\in I}$ which generates $\mt{B}_{\lambda}$ under small colimits (\cite[Theorem 5.5.1.1]{highertopostheory}). 
 Choose $f\in F$ with $\sqcup_{i\in I}G(b_i)\in \essim(\mt{C}_{f}\to \mt{C}_{\infty})$. Since $\essim(\mt{C}_{f}\to \mt{C}_{\infty})\subseteq \mt{C}_{\infty}$ is a stable subcategory closed under colimits, by the same argument as above, if follows that \[\{ b\in\mt{B}_{\lambda}\, | \, G(b)\in \essim(\mt{C}_{f}\to \mt{C}_{\infty})\}\] is closed under small colimits and contains $b_i$ for all $i\in I$, i.e.\ $G$ factors over $\mt{B}_{\lambda}\to \mt{C}_f\to \mt{C}_{\infty}$. 
\end{proof}
\begin{rem}\label{actualuniversalpropertystabilization}
The above proof also shows a universal property of $\Sigma^{\infty}_{+}\colon \mt{B}_{\infty}\to \stab{\mt{B}_{\infty}}$:
For big presentable categories $\mt{B}_{\infty}, \mt{C}_{\infty}$ denote by \[\Fun^{p\operatorname{-colim}}(\mt{B}_{\infty}, \mt{C}_{\infty})\subseteq \Fun(\mt{B}_{\infty}, \mt{C}_{\infty})\] the full subcategory of functors $F\colon \mt{B}_{\infty}\to \mt{C}_{\infty}$ such that for all presentable categories $P$ and all small colimits preserving, left-exact functors $p\colon P\to \mt{B}_{\infty}$, $F\circ p$ preserves small colimits. 
The above proof and \cref{presentablestablesubcategorycontainedinfinitestage} show that if $\mt{B}_{\infty}, \mt{C}_{\infty}$ are big presentable categories and $\mt{C}_{\infty}$ is stable, pullback along $\Sigma^{\infty}_{+}$ induces a fully faithful functor 
\[(\Sigma^{\infty}_+)^*\colon\Fun^{p\operatorname{-colim}}(\stab{\mt{B}_{\infty}}, \mt{C}_{\infty})\hookrightarrow \Fun(\mt{B}_{\infty}, \mt{C}_{\infty}).\]
Its essential image consists of those functors $F\colon \mt{B}_{\infty}\to\mt{C}_{\infty}$ such that there exists an exhaustion $\mt{B}_*^F\colon\Lambda\to \Pr^L$ of $\mt{B}_{\infty}$ by topoi such that $\mt{B}_{\lambda}^F\subseteq \mt{B}_{\infty}\xrightarrow{F} \mt{C}_{\infty}$ preserves small colimits for all $\lambda\in\Lambda$. 
\end{rem}

\subsubsection{Symmetric monoidal structure}\label{section:monoidalstructureonspectrumobjects}
We now formulate conditions under which a (symmetric) monoidal structure on a big presentable category induces a (symmetric) monoidal category on its category of spectrum objects. 
\begin{definition}\label{presentablymonoidalcat}
    A \emph{chain} of presentably (symmetric) monoidal categories is a functor 
    \[\mt{B}^{\otimes}\colon \Lambda\to \aCAlg(\Pr^L)\] from a large filtered set $\Lambda$ such that for all $\kappa\to \lambda\in \Lambda$, the underlying functor of categories $\topo{B_{\kappa}}\to \topo{B_{\lambda}}$ is a fully faithful, left-exact left adjoint. 
    
    The forget functor $\Pr^L\to \vlCat$ is lax symmetric monoidal (\cite[Propositions 4.8.1.15, 4.8.1.4]{higheralgebra}), and therefore lifts to a functor \[ \aCAlg(\Pr^L)\to \aCAlg(\vlCat).\]
    A \emph{big presentably (symmetric) monoidal category} is a (symmetric) monoidal category $\mt{B}_{\infty}^{\otimes}$ such that there exists a chain of presentably (symmetric) monoidal categories $\mt{B}^{\otimes}_*\colon \Lambda\to \aCAlg(\Pr^L)$ with colimit $\mt{B}_{\infty}^{\otimes}$ in $\aCAlg(\vlCat)$. 
\end{definition}
 
\begin{lemma}\label{colimitmonoidalcategories}
    Suppose that $\mt{B}^{\otimes}_*\colon \Lambda\to\aCAlg(\Pr^L)$ is a chain of presentably (symmetric) monoidal categories. The colimit $\mt{B}_{\infty}^{\otimes}\coloneqq \colim{\lambda\in\Lambda}\mt{B}_{\lambda}^{\otimes}$ in $\aCAlg(\vlCat)$ exists and is a (symmetric) monoidal enhancement of $\colim{\lambda\in\Lambda}{\mt{B}_{\lambda}}$. 
\end{lemma}
\begin{proof}
    The category $\vlCat$ has all $\mathcal U_1$-small colimits, and the tensor product on $\vlCat$ is compatible with $\mathcal U_1$-small colimits (\cite[Definition 3.1.1.18]{higheralgebra}) by \cite[Lemma 4.8.4.2]{higheralgebra} in the universe $\mathcal U_1$.
    Since $\Lambda$ is $\mathcal U_1$-small and filtered, \cite[Proposition 3.2.3.1]{higheralgebra} (in $\mathcal U_1)$ now implies that the colimit of $\mt{B}_*^{\otimes}$ in $\aCAlg(\vlCat)$ exists and is a (symmetric) monoidal enhancement of $\mt{B_{\infty}}$. 
    \end{proof}   
    \begin{definition}\label{bigpresentablymonoidalcatdef}
    Suppose $\mathcal D^{\otimes}$ is a possibly large (symmetric) monoidal category. 
    A subcategory $\mathcal C\subseteq \mathcal D$ is a (symmetric) monoidal subcategory if it contains the unit and for $b,c\in\mathcal C$, $b\otimes_{\mathcal D}c\in\mathcal C$. 
    By \cite[Proposition 2.2.1.1]{higheralgebra}, this implies that there exists an operad $\mathcal C^{\otimes}\subseteq\mathcal D^{\otimes}$ which exhibits $\mathcal C$ as a (symmetric) monoidal category. 
    
    If $\mathcal D^{\otimes}$ is a presentably (symmetric) monoidal category, a \emph{presentably (symmetric) monoidal subcategory} $\mathcal C\subseteq\mathcal D$ is a (symmetric) monoidal subcategory which is presentable, closed under small colimits and finite limits, and such that the induced tensor product on $\mathcal C$ preserves small colimits in both variables. 
    \end{definition}

By \cref{colimitmonoidalcategories}, a (symmetric) monoidal category $\mt{B}_{\infty}^{\otimes}$ is a big presentably (symmetric) monoidal category if and only if its underlying category is a big presentable category and admits an exhaustion by presentably (symmetric) monoidal subcategories. 

\begin{ex}
    If $\mt{X}$ is a big topos, the cartesian monoidal structure on $\mt{X}$ exhibits it a big presentably symmetric monoidal category: 
    If $\topo{X}_*\colon\Lambda\to\Pr^L$ is an exhaustion for $\topo{X}$ by topoi, for all $\lambda\in\Lambda$, $\topo{X}_{\lambda}\subseteq \topo{X}$ is closed under finite limits (\cref{fullyfaithfulnessandpreservationoflimits}), whence $\topo{X}_{\lambda}^{\times}\subseteq \topo{X}^{\times}$ is a symmetric monoidal subcategory. By universality of colimits in the topos $\topo{X}_{\lambda}$, $\topo{X}_{\lambda}^{\times}$ is presentably symmetric monoidal. 
\end{ex}
The main result of this section is the following: 
\begin{proposition}\label{symmetricmonoidalstructureonspectrumobjects}
    Suppose that $\mt{B}_{\infty}^{\otimes}$ is a big presentably (symmetric) monoidal category.
    There exists a (symmetric) monoidal structure on $\stab{\mt{B_{\infty}}}$ such that \[ \Sigma^{\infty}_{+}\colon \mt{B_{\infty}}\to \stab{\mt{B_{\infty}}}\] promotes to a (symmetric) monoidal functor and for all presentably (symmetric) monoidal subcategories $\mathcal C\subseteq \mt{B}_{\infty}$, 
    $\stab{\mathcal C}\hookrightarrow \stab{\mt{B}_{\infty}}$ is a presentably (symmetric) monoidal subcategory.
    This is unique in the following sense: If $\mt{B}_{\infty}^{\otimes}\to\stab{\mt{B}_{\infty}}^{\otimes}, \mt{B}_{\infty}^{\otimes}\to\stab{\mt{B}_{\infty}}^{\tilde\otimes}$ are two (symmetric) monoidal enhancements as above, there is an essentially unique (symmetric) monoidal equivalence $\stab{\mt{B}_{\infty}}^{\otimes}\cong \stab{\mt{B}_{\infty}}^{\tilde \otimes}$ of (symmetric) monoidal categories over $\mt{B}_{\infty}^{\otimes}$, and this is a (symmetric) monoidal enhancement of the identity $\id_{\stab{\mt{B}_{\infty}}}$. 
\end{proposition}
We will prove this below. If $\stab{\mt{B}_{\infty}}$ has small colimits, the above (symmetric) monoidal structure has a handier characterization: 
\begin{cor}\label{monoidalstructurespectrumobjectscocontinuous}Suppose that $\mt{B}_{\infty}^{\otimes}$ is a big presentably (symmetric) monoidal category and denote by $\mt{B}_{\infty}$ its underlying category. Suppose that $\stab{\mt{B}_{\infty}}$ has small coproducts. 

The (symmetric) monoidal enhancement $\mt{B}_{\infty}^{\otimes}\to \stab{\mt{B}_{\infty}}^{\otimes}$ from \cref{symmetricmonoidalstructureonspectrumobjects} is the unique one so that the tensor product on $\stab{\mt{B}_{\infty}}$ preserves small colimits in both variables.   
\end{cor} 

\begin{proof}[Proof of \cref{symmetricmonoidalstructureonspectrumobjects}]
   We only prove the statement for symmetric monoidal structures here, the statement for monoidal structures can be shown completely analogously. We first prove existence of the symmetric monoidal structure. 
   The symmetric monoidal localization \[-\otimes_{\Pr^L}\Sp\colon \Pr^L\to \Pr^L\] (\cite[Corollary 4.8.2.18]{higheralgebra}) enhances to a localization \[\Sp(-)^{\otimes}\colon \CAlg(\Pr^L)\to \CAlg(\Pr^L), \] cf.\ \cref{adjuncttioninducesadjunctiononalgebraobjects}.
   Choose an exhaustion $\mt{B}^{\otimes}_*\colon\Lambda\to \CAlg(\Pr^L)$ of $\mt{B}_{\infty}^{\otimes}$ by presentably symmetric monoidal categories.
   The adjunction unit for $\Sp(-)^{\otimes}$ yields an enhancement \[\mt{B}^{\otimes}\to \stab{\mt{B}^{\otimes}}^{\otimes}\in \Fun(\Lambda, \CAlg(\Pr^L))\] of $\mt{B}_*\to \Sp\otimes_{\Pr^L}\mt{B}_*$, cf.\ \cref{adjuncttioninducesadjunctiononalgebraobjects} or \cite[Theorem 5.1/Corollary 5.5.(ii)]{GepnerGrothNikolaus}. 
   Composing with the forget functor $\CAlg(\Pr^L)\to \CAlg(\vlCat)$ yields a natural transformation \[\mt{B}^{\otimes}_*\to \stab{\mt{B}_*^{\otimes}}^{\otimes}\in \Fun(\Lambda, \CAlg(\vlCat)).\] 
By \cite[Lemma 4.8.4.2]{higheralgebra} in the large universe $\mathcal U_1$, $\vlCat$ has all large colimits and the symmetric monoidal structure is compatible with large colimits. Hence, by \cite[Corollary 3.2.3.2]{higheralgebra} (in $\mathcal U_1$), the colimit \[\colim{\Lambda}\mt{B}^{\otimes}_*\to \colim{\Lambda}\stab{\mt{B}^{\otimes}_*}^{\otimes}\in \CAlg(\vlCat)\] exists and is a symmetric monoidal enhancement of $\Sigma^{\infty}_{+}$. By construction, $ \stab{\mt{B_{\lambda}}}\subseteq \stab{\mt{B_{\infty}}}$ is a symmetric monoidal subcategory for all $\lambda\in\Lambda$. 
Suppose that $\mathcal C\hookrightarrow \mt{B_{\infty}}$ is a presentably symmetric monoidal subcategory.
Since $\mathcal C\subseteq \mt{B}_{\infty}$ is closed under small colimits, $\mathcal C$ is presentable and $\mt{B}_{\infty}$ is locally small (\cref{fullyfaithfulnessandpreservationoflimits}), $\mt{C}\subseteq \mt{B}_{\infty}$ is a left adjoint by \cite[Remark 5.5.2.10]{highertopostheory}.
As $\mt{C}\subseteq {\mt{B}_{\infty}}$ is closed under finite limits, the induced functor $\stab{\mt{C}}\to\stab{\mt{B}_{\infty}}$ is a fully faithful left adjoint and \[\mathcal C\subseteq \mt{B}_{\infty}\xrightarrow{\Sigma^{\infty}_{+}}\stab{\mt{B}_{\infty}}\] factors over \[\mathcal C\xrightarrow{\Sigma^{\infty}_{+}}\stab{\mathcal C}\subseteq \stab{\mt{B}_{\infty}}\] by \cref{stabilizationcommuteswithsuspension}. 
As $\mathcal C\subseteq \mt{B}_{\infty}\xrightarrow{\Sigma^{\infty}_{+}}\stab{\mt{B}_{\infty}}$ is symmetric monoidal, this implies that $\stab{\mathcal C}$ contains the unit of $\stab{\mt{B}_{\infty}}$. 
By \cref{presentablestablesubcategorycontainedinfinitestage}, there exists $\lambda\in\Lambda$ such that $\stab{\mathcal C}\subseteq \stab{\mt{B}_{\lambda}}$. 
For $x\in \stab{\mt{\mathcal C}}$ denote by $\mathcal G_x\subseteq \stab{\mt{C}}$ the full subcategory on objects $g\in \stab{\mathcal C}$ such that $x\otimes g, g\otimes x\in\stab{\mt{C}}$.  
Since \[x\otimes_{\stab{\mt{B}_{\lambda}}} -, -\otimes_{\stab{\mt{B}}_{\lambda}} x\colon \stab{\mt{B}_{\lambda}}\to \stab{\mt{B}_{\lambda}}\] are exact, cocontinuous functors and $\stab{\mt{C}}\subseteq\stab{\mt{B}_{\lambda}}$ is a stable subcategory closed under small colimits, $\mathcal \mathcal G_x\subseteq \stab{\mt{C}}$ is a stable subcategory closed under small colimits. 
It therefore suffices to show that for all $c\in \mt{C}$, $\Sigma^{\infty}_{+}c\in \mathcal G_x$, then it follows from \cref{suspensionspectragenerate} that $\mathcal G_x=\stab{\mt{C}}$. If $x=\Sigma^{\infty}_{+}d, d\in\mathcal C$, this holds by monoidality of $\Sigma^{\infty}_{+}\circ (\mathcal C\subseteq \mt{B}_{\infty})$. 
By the same argument as above, 
\[\mathcal G\coloneqq \{ x\in \stab{\mt{C}}\, |\, \mathcal G_x=\stab{\mt{C}}\}\subseteq\stab{\mt{C}}\] is a stable subcategory closed under small colimits. Since it contains the image of $\Sigma^{\infty}_{+}$, $\mathcal G=\stab{\mt{\mathcal C}}$ by \cref{suspensionspectragenerate}, which proves that $\stab{\mt{C}}\subseteq \stab{\mt{B_{\lambda}}}$ is a symmetric monoidal subcategory. 
As $\stab{\mt{C}}\subseteq \stab{\mt{B}_{\lambda}}$ preserves small colimits, the induced tensor product on $\stab{\mt{C}}$ is cocontinuous in both variables.
This proves the existence of a symmetric monoidal structure with the claimed properties. 

It remains to show the uniqueness statement. 
Suppose that $\stab{\mt{B_{\infty}}}^{\tilde\otimes}$ is another symmetric monoidal structure and $\mt{B}_{\infty}^{\otimes}\to\stab{\mt{B}_{\infty}}^{\tilde\otimes}$ is a symmetric monoidal enhancement of $\Sigma^{\infty}_{+}$ such that for all $\lambda\in\Lambda$, $ \stab{\mt{B_{\lambda}}}\subseteq\stab{\mt{B_{\infty}}}^{\tilde \otimes}$ is a presentably symmetric monoidal subcategory. 
Denote by \[\Map^{c, \Lambda}_{\vlCat}(\stab{\mt{B}_{\infty}}, \stab{\mt{B}_{\infty}})\subseteq \Map_{\vlCat}(\stab{\mt{B}_{\infty}}, \stab{\mt{B}_{\infty}})\] the full subspace on functors $F\colon \stab{\mt{B}_{\infty}}\to\stab{\mt{B}_{\infty}}$ such that for all $\lambda\in\Lambda$, 
$F|_{\stab{\mt{B}_{\lambda}}}$ preserves small colimits, and by  $\vlCat^c\subseteq \vlCat$ the wide subcategory on small colimits preserving functors. 
Then \begin{align*}&\Map_{\CAlg_{\mt{B_{\infty}^{\otimes}}}(\vlCat)}(\stab{\mt{B_{\infty}}}^{\otimes}, \stab{\mt{B_{\infty}}}^{\tilde\otimes})\times_{\Map_{\vlCat}(\stab{\mt{B}_{\infty}}, \stab{\mt{B}_{\infty}})}\Map_{\vlCat}^{c, \Lambda}(\stab{\mt{B}_{\infty}}, \stab{\mt{B}_{\infty}}) \\ &\cong \clim{\kappa\in\Lambda}\left(\Map_{{\CAlg_{\mt{B_{\kappa}^{\otimes}}}(\vlCat)}}(\stab{\mt{B_{\kappa}}}^{\otimes}, \stab{\mt{B_{\infty}}}^{\tilde\otimes})\times_{\Map_{\vlCat}(\stab{\mt{B}_{\kappa}}, \stab{\mt{B}_{\infty}})}\Map_{\vlCat^c}(\stab{\mt{B}_{\kappa}}, \stab{\mt{B}_{\infty}})\right)\\ &\cong  \clim{\kappa\in\Lambda}\left(\Map_{{\CAlg_{\mt{B_{\kappa}^{\otimes}}}(\vlCat)}}(\stab{\mt{B_{\kappa}}}^{\otimes}, \stab{\mt{B_{\kappa}}}^{\tilde\otimes})\times_{\Map_{\vlCat}(\stab{\mt{B}_{\kappa}}, \stab{\mt{B}_{\kappa}})}\Map_{\vlCat^c}(\stab{\mt{B}_{\kappa}}, \stab{\mt{B}_{\kappa}})\right), \end{align*} where we used in the last step that every cocontinuous functor $\mt{B}_{\kappa}\to\stab{\mt{B}_{\infty}}$ factors over $\stab{\mt{B}_{\mu}}\subseteq \stab{\mt{B}_{\infty}}$ for some $\mu\geq\kappa$ (\cref{presentablestablesubcategorycontainedinfinitestage}), and that $\Lambda$ is filtered. 
By assumption, $ \mt{B}_{\kappa}\to \stab{\mt{B_{\kappa}}}^{\otimes}, \mt{B}_{\kappa}\to \stab{\mt{B_{\kappa}}}^{\tilde\otimes|_{ \stab{\mt{B_{\kappa}}}}}\in \CAlg(\Pr^L)$, whence the right-hand side is equivalent to \begin{align*}\clim{\kappa\in\Lambda} \Map_{{\CAlg_{\mt{B_{\kappa}^{\otimes}}}(\Pr^L)}}(\stab{\mt{B_{\kappa}}}^{\otimes}, \stab{\mt{B_{\kappa}}}^{\tilde\otimes}).
    \end{align*}
As $\Sp(-)^{\otimes}\colon \CAlg(\Pr^L)\to \CAlg(\Pr^L)$ is a localization, $ \Map_{{\CAlg_{\mt{B_{\kappa}^{\otimes}}}(\Pr^L)}}(\stab{\mt{B_{\kappa}}}^{\otimes}, \stab{\mt{B_{\kappa}}}^{\tilde\otimes|_{ \stab{\mt{B_{\kappa}}}}})$ is contractible for all $\kappa\in\Lambda$. 
This shows that \[ \Map_{\CAlg_{\mt{B_{\infty}^{\otimes}}}}(\stab{\mt{B_{\infty}}}^{\otimes}, \stab{\mt{B_{\infty}}}^{\tilde\otimes})\times_{\Map_{\vlCat}(\stab{\mt{B}_{\infty}}, \stab{\mt{B}_{\infty}})}\Map_{\vlCat}^{c, \Lambda}(\stab{\mt{B}_{\infty}}, \stab{\mt{B}_{\infty}})\cong \clim{\Lambda}*\] is contractible, i.e.\ there exists an essentially unique symmetric monoidal functor \[F^{\otimes}\colon \stab{\mt{B}_{\infty}}^{\otimes}\to \stab{\mt{B}_{\infty}}^{\tilde \otimes}\] over $\mt{B}_{\infty}^{\otimes}$ such that $F|_{\stab{\mt{B}_{\lambda}}}$ is cocontinuous for all $\lambda\in\Lambda$.
As \[\Sp(-)^{\otimes}\colon \CAlg(\Pr^L)\to \CAlg(\Pr^L)\] is a lift of the localization $\Sp(-)\colon \Pr^L\to\Pr^L$, for all $\kappa\in\Lambda$, the essentially unique cocontinuous symmetric monoidal functor \[\stab{\mt{B}_{\kappa}}^{\otimes}\to \stab{\mt{B}_{\kappa}}^{\tilde{\otimes}}\in {\CAlg_{\mt{B_{\kappa}^{\otimes}}}(\Pr^L)}\] is a symmetric monoidal enhancement of the identity for all $\kappa\in\Lambda$. 
As $\CAlg(\vlCat)\to \vlCat$ preserves large sifted colimits (\cite[Lemma 4.8.4.2, Corollary 3.2.3.2]{higheralgebra}), this implies that \[F\colon \stab{\mt{B}_{\infty}}^{\otimes}\to \stab{\mt{B}_{\infty}}^{\tilde \otimes}\] is a symmetric monoidal enhancement of the identity and in particular an equivalence. 
As every equivalence has to preserve small colimits and $\stab{\mt{B}_{\lambda}}\subseteq \stab{\mt{B}_{\infty}}$ is closed under small colimits, the space of symmetric monoidal equivalences \[\stab{\mt{B}_{\infty}}^{\otimes}\to \stab{\mt{B_{\infty}}}^{\tilde\otimes}\] in $\CAlg_{\mt{B_{\infty}^{\otimes}}}$ is a subspace of \[\Map_{\CAlg_{\mt{B_{\infty}^{\otimes}}}}(\stab{\mt{B_{\infty}}}^{\otimes}, \stab{\mt{B_{\infty}}}^{\tilde\otimes})\times_{\Map_{\vlCat}(\stab{\mt{B}_{\infty}}, \stab{\mt{B}_{\infty}})}\Map_{\vlCat}^{c, \Lambda}(\stab{\mt{B}_{\infty}}, \stab{\mt{B}_{\infty}})\] and hence contractible. (It is non-empty by the above). 
\end{proof}
\begin{proof}[Proof of \cref{monoidalstructurespectrumobjectscocontinuous}]
    We only prove the statement for symmetric monoidal structures here, the statement for monoidal structures can be shown completely analogously. 
    We first show that the tensor product on $\stab{\mt{B}_{\infty}}$ is cocontinuous in both variables. 
    Suppose that $S$ is a small category and $F\colon S\to \stab{\mt{B}_{\infty}}$. 
    Choose an essentially surjective functor $i\colon T\to S$ from a set $T$ and $\lambda\in\Lambda$ with  $\oplus_{t\in T}F(it)\in \essim(\stab{\mt{B}_{\lambda}}\subseteq \stab{\mt{B}_{\infty}})=:\overline{\stab{\mt{B}_{\lambda}}}$. 
    Since $\overline{\stab{\mt{B}_{\lambda}}}\subseteq\stab{\mt{B}_{\infty}}$ is closed under colimits, for all $t\in T$, 
    \[F(it)=\Cofib(\oplus_{s\in T}F(is)\to \oplus_{\substack{s\in T\\ s\neq t}}F(it)\hookrightarrow \oplus_{s\in T}F(is))\in \overline{\stab{\mt{B}_{\lambda}}}.\] 
    This implies that $F$ factors over $S\to \stab{\mt{B}_{\lambda}}\subseteq \stab{\mt{B}_{\infty}}$. In particular, $\colim{S}F$ exists in $\stab{\mt{B}_{\infty}}$ and can be computed in $\stab{\mt{B}_{\lambda}}\subseteq \stab{\mt{B}_{\infty}}$. 
    For $b\in\stab{\mt{B}_{\infty}}$ choose $\lambda\to\mu\in\Lambda$ such that $b\in\stab{\mt{B}_{\mu}}\subseteq \stab{\mt{B}_{\infty}}$. 
    As $\stab{\mt{B}_{\mu}}\subseteq \stab{\mt{B}_{\infty}}$ is a presentably symmetric monoidal subcategory, it follows that \[\colim{S}(F\otimes_{\stab{\mt{B}_{\infty}}} b)\cong \colim{S}(F\otimes_{\stab{\mt{B}_{\mu}}} b)\cong (\colim{S}F)\otimes_{\stab{\mt{B}_{\mu}}} b\cong (\colim{S}F)\otimes_{\stab{\mt{B}_{\infty}}} b, \] which shows that the tensor product on $\stab{\mt{B}_{\infty}}$ is cocontinuous in both variables. 

    Suppose now that $\mt{B}^{\otimes}_{\infty}\to \stab{\mt{B}_{\infty}}^{\tilde \otimes}$ is a symmetric monoidal enhancement of $\Sigma^{\infty}_{+}$ such that the tensor product on $\stab{\mt{B}_{\infty}}$ is cocontinuous in both variables, and $\mathcal C\subseteq \mt{B}_{\infty}$ is a presentably symmetric monoidal subcategory. 
    Since the tensor product on $\stab{\mt{B}_{\infty}}$ is cocontinuous, one can argue as in the proof of \cref{symmetricmonoidalstructureonspectrumobjects} that $\stab{\mathcal C}\subseteq \stab{\mt{B}_{\infty}}$ is a symmetric monoidal subcategory. As $\stab{\mathcal C}\subseteq \stab{\mt{B}_{\infty}}$ is closed under finite limits and small colimits, it is a presentably symmetric monoidal subcategory.
\end{proof}
For categories of accessible (hyper)sheaves on (hyper)accessible explicit covering sites, the above recovers the symmetric monoidal structure constructed in \ref{constructionmonoidalstructureaccessiblesheaves}:
\begin{lemma}\label{symmetricmonoidalstructureaccessiblesheavesofspectraclosed}
    Suppose that $(\mathcal C,S)$ is a (hyper)ac\-ces\-si\-ble explicit covering site. 
    The equivalence from \cref{stabilizationisspectrumvaluedsheaves} enhances to  a symmetric monoidal equivalence \[\stab{\mywidehatshvacc{S}(\mathcal C, \an)^{\times}}^{\otimes}\cong\mywidehatshvacc{S}(\mathcal C, \Sp)^{\otimes},\] where the left-hand side is endowed with the symmetric monoidal structure from \cref{symmetricmonoidalstructureonspectrumobjects} and the right-hand side with the symmetric monoidal structure induced from $\Sp^{\otimes}$ via \cref{constructionmonoidalstructureaccessiblesheaves}. 
    In particular, the symmetric monoidal structure on $\stab{\mywidehatshvacc{S}(\mathcal C, \an)}$ is closed.
\end{lemma}
\begin{proof}
    By \cref{constructionmonoidalstructureaccessiblesheaves}, $\Sigma^{\infty}_{+}\colon \an^{\times}\to\Sp^{\otimes}$ induces a symmetric monoidal functor 
    \[ \Sigma^{\infty, \otimes}_{+}\colon\mywidehatshvacc{S}(\mathcal C, \an)^{\times}\to\mywidehatshvacc{S}(\mathcal C, \Sp)^{\otimes}.\] 
    Under the identification $\mywidehatshvacc{S}(\mathcal C, \Sp)\cong \stab{\mywidehatshvacc{S}(\mathcal C)}$ of \cref{stabilizationisspectrumvaluedsheaves}, the functor underlying $\Sigma^{\infty, \otimes}_{+}$ is left adjoint to $\Omega^{\infty}$.
    As the symmetric monoidal structure on $\mywidehatshvacc{S}(\mathcal C, \Sp)$ is closed (\cref{closednessmonoidalstructureaccessiblesheaves}) and $\mywidehatshvacc{S}(\mathcal C, \Sp)$ has all small colimits (\cref{limitscolimitscanbecomputedfinitestagesheaves}), the statement now follows from \cref{monoidalstructurespectrumobjectscocontinuous}. 
\end{proof}
\begin{rem}\label{propertiesinternalhomspectra}
Suppose that $\topo{B_{\infty}}$ is a big presentable symmetric monoidal category such that the induced symmetric monoidal structure on $\stab{\mt{B}_{\infty}}$ is closed. 
The internal Hom \[\imap_{\stab{\mt{B}_{\infty}}}(-,-)\colon \stab{\mt{B}_{\infty}}^{\operatorname{op}}\times \stab{\mt{B}_{\infty}}\to \stab{\mt{B}_{\infty}}\] has the following properties: 
\begin{romanenum}
\item If $I$ is a small category such that $\stab{\mt{B}_{\infty}}$ has $I$-indexed colimits and $I^{\operatorname{op}}$-indexed limits, then \linebreak[1]$\imap_{\stab{\mt{B}_{\infty}}}(-,-)$ preserves $I$-indexed limits in the first variable by \cref{leftadjointsstableundercolimits} and \cite[Proposition 5.2.6.2]{highertopostheory}. 
\item Since $\imap_{\stab{\mt{B}_{\infty}}}(-,-)$ preserves limits in the second variable, stability implies that for $i\in\mathbb Z$, $\Sigma^i\imap_{\stab{\mt{B}_{\infty}}}(-,-)\cong \imap_{\stab{\mt{B}_{\infty}}}(-, \Sigma^i-)$. 

\item As $\Sigma^{\infty, \otimes}_{+}\colon \topo{B}_{\infty}^{\otimes}\to \stab{\topo{B}_{\infty}}^{\otimes}$ is symmetric monoidal, for $b\in \mt{B}_{\infty}$, 
\[ \Omega^{\infty}\imap_{\stab{\mt{B}_{\infty}}}(\Sigma^{\infty}_{+}b,-)\] is right adjoint to $\Sigma^{\infty}_{+}(b\otimes -)$. 
In particular, if the symmetric monoidal structure on $\mt{B}_{\infty}$ is closed with internal mapping functors $\iMap_{\mt{B}_{\infty}}(-,-)$, then 
\[ \Omega^{\infty}\imap_{\stab{\mt{B}_{\infty}}}(\Sigma^{\infty}_{+}-,-)\cong \iMap_{\mt{B}_{\infty}}(-, \Omega^{\infty}-).\] 
\end{romanenum}
\end{rem}

The universal property of the stabilization functor (\cref{universalpropertyspectrumobjectsbigpresentable}) has a (symmetric) monoidal analogue: 
For potentially large (symmetric) monoidal categories $\mathcal B^{\otimes}, \mathcal C^{\otimes}$ denote by $\Fun^{\otimes}(\mathcal B^{\otimes}, \mathcal C^{\otimes})$ the category of (symmetric) monoidal functors and by \[\Fun^{\operatorname{colim}, \otimes}(\mathcal B^{\otimes}, \mathcal C^{\otimes})\subseteq \Fun^{\otimes}(\mathcal B^{\otimes}, \mathcal C^{\otimes})\] the full subcategory on (symmetric) monoidal functors whose underlying functor preserves small colimits.
\begin{cor}\label{universalpropertyspectrumobjectsbigpresentablemonoidal}
Suppose $\mt{B}_{\infty}^{\otimes}, \mt{C}_{\infty}^{\otimes}$ are two big presentably (symmetric) monoidal categories and $\mt{C}_{\infty}$ is stable. 
Then pullback along $\Sigma^{\infty, \otimes}_{+}\colon \mt{B}_{\infty}^{\otimes}\to \stab{\mt{B}_{\infty}}^{\otimes}$ defines a fully faithful functor 
\[ \Fun^{\operatorname{colim}, \otimes}(\stab{\mt{B}_{\infty}}^{\otimes}, \mt{C}_{\infty}^{\otimes})\to \Fun^{\operatorname{colim}, \otimes}(\mt{B}_{\infty}^{\otimes}, \mt{C}_{\infty}^{\otimes}
)\] with essential image  
\[ \Fun^{\operatorname{colim}}(\stab{\mt{B}_{\infty}}, \mt{C}_{\infty})\times_{\Fun^{\operatorname{colim}}({\mt{B}_{\infty}}, \mt{C}_{\infty})}\Fun^{\operatorname{colim}, \otimes}(\mt{B}_{\infty}^{\otimes}, \mt{C}_{\infty}^{\otimes}).\]  
In particular, if $\stab{\mt{B}_{\infty}}$ and $\mt{C}_{\infty}$ admit small coproducts, then \[ \Fun^{\operatorname{colim}, \otimes}(\stab{\mt{B}_{\infty}}^{\otimes}, \mt{C}_{\infty}^{\otimes})\cong \Fun^{\operatorname{colim}, \otimes}(\mt{B}_{\infty}^{\otimes}, \mt{C}_{\infty}^{\otimes}).\] 
\end{cor}
\begin{proof}We first prove fully faithfulness.
    Choose exhaustions \[\mt{B}_{*}^{\otimes}\colon \Lambda\to \aCAlg(\Pr^L), \mt{C}_{*}^{\otimes}\colon F\to \aCAlg(\Pr^L)\] of $\mt{B}_{\infty}^{\otimes}, \mt{C}_{\infty}^{\otimes}$ by presentably (symmetric) monoidal categories. 
    Since for all $\lambda\in\Lambda$, $\mt{B}_{\lambda}\subseteq \mt{B}_{\infty}$ is closed under colimits and $\stab{\mt{B}_{\infty}}^{\otimes}\cong \colim{\lambda\in\Lambda}\stab{\mt{B}_{\lambda}}^{\otimes}$, 
    \[ \Fun^{\operatorname{colim}, \otimes}(\stab{\mt{B}_{\infty}}^{\otimes}, \mt{C}_{\infty}^{\otimes})\hookrightarrow \clim{\lambda\in\Lambda}\Fun^{\operatorname{colim}, \otimes}(\stab{\mt{B}_{\lambda}}^{\otimes}, \mt{C}_{\infty}^{\otimes})\] is fully faithful.  
    Since $\mt{C}_{\infty}$ is stable, \cref{presentablestablesubcategorycontainedinfinitestage,filteredcolimitsofcategoriesappendix} imply that for all $\lambda\in\Lambda$, 
\[ \colim{f\in F}\Fun^{\operatorname{colim}, \otimes}(\stab{\mt{B}_{\lambda}}^{\otimes}, \mt{C}_{f}^{\otimes})\cong \Fun^{\operatorname{colim}, \otimes}(\stab{\mt{B}_{\lambda}}^{\otimes}, \mt{C}_{\infty}^{\otimes}).\] 
For all $f\in F$, $\mathcal C_f\subseteq \mathcal C$ is closed under finite limits and colimits and hence stable. 
As $\Sp\in \Pr^L$ is idempotent, this implies that for all $\lambda\in\Lambda$ and $f\in F$, 
\[\Fun^{\operatorname{colim}, \otimes}(\stab{\mt{B}_{\lambda}}^{\otimes}, \mt{C}_{f}^{\otimes})\cong \Fun^{\operatorname{colim}, \otimes}(\mt{B}_{\lambda}^{\otimes}, \mt{C}_f^{\otimes})\] via pullback along $\Sigma^{\infty,\otimes}_{+,\mt{B}_{\lambda}}\colon\mt{B}_{\lambda}^{\otimes}\to\stab{\mt{B}_{\lambda}}^{\otimes}$, cf.\ \cite[Proposition 5.4]{GepnerGrothNikolaus}. 
By \cref{filteredcolimitsofcategoriesappendix}, $\colim{f\in F}\Fun^{\operatorname{colim}, \otimes}(\mt{B}_{\lambda}^{\otimes}, \mt{C}_f^{\otimes})\hookrightarrow \Fun^{\operatorname{colim}, \otimes}(\mt{B}_{\lambda}^{\otimes}, \mt{C}_{\infty}^{\otimes})$ is fully faithful. 
By construction (see the proof of \cref{symmetricmonoidalstructureonspectrumobjects}), \begin{align}\label{colimitdescriptionmonoidalstructuresuspension}\Sigma^{\infty, \otimes}_{+}\colon \stab{\mt{B}_{\infty}}^{\otimes}\cong \colim{\lambda\in\Lambda}(\Sigma^{\infty, \otimes}_{+, \lambda}\colon \mt{B}_{\lambda}^{\otimes}\to \stab{\mt{B}_{\lambda}^{\otimes}})\in\Fun(\Delta^1, \aCAlg(\vlCat)).\end{align}
This implies that pullback along $\Sigma^{\infty, \otimes}_{+}$ factors as \begin{align*}\Fun^{\operatorname{colim}, \otimes}(\stab{\mt{B}_{\infty}}^{\otimes}, \mt{C}_{\infty}^{\otimes}) \hookrightarrow \clim{\lambda\in\Lambda}\Fun^{\operatorname{colim}, \otimes}(\stab{\mt{B_{\lambda}}}^{\otimes}, \mt{C}_{\infty}^{\otimes})& \cong \clim{\lambda\in\Lambda}\colim{f\in F} \Fun^{\operatorname{colim}, \otimes}(\mt{B}_{\lambda}^{\otimes}, \mt{C}_{f}^{\otimes})\\ &\hookrightarrow \clim{\lambda\in\Lambda}\Fun^{\operatorname{colim}, \otimes}(\mt{B}_{\lambda}^{\otimes}, \mt{C}_{\infty}^{\otimes})
    \\ &\hookrightarrow \Fun^{\otimes}(\mt{B}_{\infty}^{\otimes}, \mt{C}_{\infty}^{\otimes}), \end{align*} and is in particular fully faithful. 
Pullback along $\Sigma^{\infty, \otimes}_{+}$ also factors over
\begin{align*}\Fun^{\operatorname{colim}, \otimes}(\stab{\mt{B}_{\infty}}^{\otimes}, \mt{C}_{\infty}^{\otimes})\xrightarrow{i}\Fun^{\operatorname{colim}}(\stab{\mt{B}_{\infty}}, \mt{C}_{\infty})&\times_{\Fun^{\operatorname{colim}}(\mt{B}_{\infty}, \mt{C}_{\infty})}\Fun^{\operatorname{colim}, \otimes}(\mt{B}_{\infty}^{\otimes}, \mt{C}_{\infty}^{\otimes})\\ &\to \Fun^{\operatorname{colim}, \otimes}(\mt{B}_{\infty}^{\otimes}, \mt{C}_{\infty}^{\otimes}).\end{align*} 
By \cref{universalpropertyspectrumobjectsbigpresentable}, the right functor is fully faithful, hence so is $i$ by the above.
We claim that $i$ is an equivalence. 
Suppose that $G\colon \stab{\mt{B}_{\infty}}\to \mt{C}_{\infty}$ is a cocontinuous functor and there exists a cocontinuous (symmetric) monoidal functor $H^{\otimes}\colon \mt{B}_{\infty}^{\otimes}\to \mt{C}_{\infty}^{\otimes}$ with $H\cong G\circ \Sigma^{\infty}_{+}$. 

\ref{colimitdescriptionmonoidalstructuresuspension} implies that 
\[ \Fun^{\otimes}(\stab{\mt{B}_{\infty}^{\otimes}}, \mt{C}_{\infty}^{\otimes})\times_{\Fun^{\otimes}(\mt{B}_{\infty}^{\otimes}, \mt{C}_{\infty}^{\otimes})}\{ H^{\otimes}\}\cong \clim{\lambda\in\Lambda}\bigl(\Fun^{\otimes}(\stab{\mt{B}_{\lambda}^{\otimes}}, \mt{C}_{\infty}^{\otimes})\times_{\Fun^{\otimes}(\mt{B_{\lambda}^{\otimes}}, \mt{C}_{\infty}^{\otimes})}\{ H^{\otimes}|_{\mt{B}_{\lambda}^{\otimes}}\}\bigr).\] 
For $\lambda\in\Lambda$, choose $f\in F$ such that $G(\stab{\mt{B}_{\lambda}})\subseteq \mathcal C_f$. This is possible by \cref{presentablestablesubcategorycontainedinfinitestage}. Then $H^{\otimes}(\mt{B}_{\lambda})\subseteq \mathcal C_f^{\otimes}$, and hence \begin{align*}\Fun^{\operatorname{colim}, \otimes}(\stab{\mt{B}_{\lambda}}^{\otimes}, \mt{C}_{\infty}^{\otimes})\times_{\Fun^{\operatorname{colim}, \otimes}(\mt{B_{\lambda}^{\otimes}}, \mt{C}_{\infty}^{\otimes})}\{ H^{\otimes}|_{\mt{B}_{\lambda}^{\otimes}}\}\\ \cong \Fun^{\operatorname{colim}, \otimes}(\stab{\mt{B}_{\lambda}}^{\otimes}, \mt{C}_{f}^{\otimes})\times_{\Fun^{\otimes}(\mt{B}_{\lambda}^{\otimes}, \mt{C}_f^{\otimes})}\{H^{\otimes}|_{\mt{B}_{\lambda}^{\otimes}}\}.\end{align*}
The right-hand side is contractible since $\Sp$ is an idempotent in $\Pr^L$,  cf.\ \cite[Proposition 5.4]{GepnerGrothNikolaus}.
This implies that \[*\cong \clim{\lambda\in\Lambda}*\cong \clim{\lambda\in\Lambda} \Fun^{\operatorname{colim}, \otimes}(\stab{\mt{B}_{\lambda}}^{\otimes}, \mt{C}_{\infty}^{\otimes})\times_{\Fun^{\operatorname{colim}, \otimes}(\mt{B_{\lambda}^{\otimes}}, \mt{C}_{\infty}^{\otimes})}\{ H^{\otimes}|_{\mt{B}_{\lambda}^{\otimes}}\}, \] i.e.\ there exists an essentially unique (symmetric) monoidal functor $\tilde G^{\otimes}\colon \stab{\mt{B}_{\infty}}\to \stab{\mt{C}_{\infty}}$ with $\tilde G^{\otimes}\circ \Sigma^{\infty, \otimes}_{+}=H^{\otimes}$ such that $\tilde{G}^{\otimes}|_{\stab{\mt{B}_{\lambda}}}$ is cocontinuous for all $\lambda\in\Lambda$. 
We have shown in the proof of \cref{universalpropertyspectrumobjectsbigpresentable} that \[(\Sigma^{\infty}_{+})^*\colon \clim{\lambda\in\Lambda}\Fun^{\operatorname{colim}}(\stab{\mt{B}_{\lambda}}, \mathcal C_{\infty})\to \clim{\lambda\in\Lambda}\Fun^{\operatorname{colim}}(\mt{B}_{\lambda}, \mathcal C_{\infty})\] is fully faithful.
As $\tilde G\circ \Sigma^{\infty}_{+}\cong H\cong G\circ \Sigma^{\infty}_{+}$, this shows that $\tilde G^{\otimes}$ is a symmetric monoidal enhancement of $G$, which proves that 
\[ \Fun^{\operatorname{colim}, \otimes}(\stab{\mt{B}_{\infty}}^{\otimes}, \mt{C}_{\infty}^{\otimes})\xrightarrow{i}\Fun^{\operatorname{colim}}(\stab{\mt{B}_{\infty}}, \mt{C}_{\infty})\] is an equivalence. 
\cref{universalpropertyspectrumobjectsbigpresentable} now implies that if $\stab{\mt{B}_{\infty}}$ and $\mt{C}_{\infty}$ admit small coproducts, then \[ \Fun^{\operatorname{colim}, \otimes}(\stab{\mt{B}_{\infty}}^{\otimes}, \mt{C}_{\infty}^{\otimes})\cong\Fun^{\operatorname{colim}, \otimes}(\mt{B}_{\infty}^{\otimes}, \mt{C}_{\infty}^{\otimes}).\qedhere\] 
\end{proof}
\begin{cor}\label{monoidalleftadjointsstabilisetomonoidalleftadjoints}
Suppose that $F^{\otimes}\colon \mt{B}_{\infty}^{\otimes}\to \mt{C}_{\infty}^{\otimes}$ is a (symmetric) monoidal functor between big presentably (symmetric) monoidal categories such that the underlying functor $F\colon \mt{B}_{\infty}\to \mt{C}_{\infty}$ is a left-exact left adjoint. 
There exists an essentially unique filler \[F_{\Sp}^{\otimes}\colon \stab{\mt{B}_{\infty}}^{\otimes}\to \stab{\mt{C}_{\infty}}^{\otimes}\] for the diagram 
    \begin{center}
        \begin{tikzcd}
       \mt{B}_{\infty}^{\otimes}\arrow[d,"F^{\otimes}"'] \arrow[r,"\Sigma^{\infty, \otimes}_{+}"] & \stab{\mt{B}_{\infty}}^{\otimes}\arrow[d,dashed, "F_{\Sp}^{\otimes}"', "\exists!"]\\ 
        {C}_{\infty}^{\otimes}\arrow[r,"\Sigma^{\infty, \otimes}_{+}"] & \stab{\mt{C}_{\infty}}^{\otimes}.
        \end{tikzcd}
        \end{center} The functor $F_{\Sp}^{\otimes}$ is a (symmetric) monoidal enhancement of $F_{\Sp}$. 
\end{cor}
\begin{proof}This follows from \cref{universalpropertyspectrumobjectsbigpresentablemonoidal} since $F_{\Sp}\circ \Sigma^{\infty}_{+}\cong \Sigma^{\infty}_{+}\circ F$ by \cref{stabilizationcommuteswithsuspension}. 
\end{proof}

\begin{lemma}\label{monoidalstructurecompatibletstructure}
    Suppose $\mathcal B^{\otimes}\colon \Lambda\to\CAlg(\Pr^L)$ is a chain of presentably (symmetric) monoidal categories. 
    The induced (symmetric) monoidal structure on  $\stab{\mt{B}_{\infty}}$ is \textit{compatible} with the $t$-structure from \cref{tstructurespectrumobjects}, i.e.\ 
    \[ \stab{\mt{B_{\infty}}}_{\geq 0}\subseteq \stab{\mt{B_{\infty}}}\] is a (symmetric) monoidal subcategory (\cref{bigpresentablymonoidalcatdef}).
\end{lemma}
\begin{proof}
    Since $\stab{\mt{B}_{\infty}}_{\geq 0}$ contains the unit, it suffices to show that for all $x,y\in\stab{\mt{B}_{\infty}}_{\geq 0}$, \[x\otimes y,y\otimes x\in\stab{\mt{B}_{\infty}}_{\geq 0}.\] 
    Choose an exhaustion $\mt{B}_*\colon\Lambda\to \CAlg(\Pr^L)$ for $\mt{B}_{\infty}^{\otimes}$ by presentably (symmetric) monoidal subcategories. 
    For all $\lambda\in\Lambda$, 
    $\stab{\mt{B}_{\lambda}}^{\otimes}\subseteq \stab{\mt{B}_{\infty}}^{\otimes}$ is a (symmetric) monoidal subcategory by \cref{symmetricmonoidalstructureonspectrumobjects}, and by \cref{tstructurespectrumobjects}, for all $\lambda\in\Lambda$, \[\stab{\mt{B}_{\lambda}}_{\geq 0}=\stab{\mt{B}_{\lambda}}\cap \stab{\mt{B}_{\infty}}_{\geq 0}.\] 
    As $\Lambda$ is filtered, it therefore suffices to show that $\stab{\mt{B}_{\lambda}}_{\geq 0}$ is closed under $-\otimes_{\stab{\mt{B}_{\lambda}}}-$ for all $\lambda\in\Lambda$. 
    Fix $\lambda\in\Lambda$ and for $x\in \stab{\mt{B}_{\lambda}}$ let  
    \[ \mathcal G_x\coloneqq \{ y\in \stab{\mt{B}_{\lambda}}_{\geq 0}\, |\, x\otimes y\in\stab{\mt{B}_{\lambda}}_{\geq 0}\}.\] 
    As the (symmetric) monoidal structure on $\stab{\mt{B}_{\lambda}}_{\geq 0}$ is cocontinuous and $\stab{\mt{B}_{\lambda}}_{\geq 0}\subseteq\stab{\mt{B}_{\lambda}}$ is closed under small colimits and extensions, $\mathcal G_x$ is closed under small colimits and extensions. 
    
    Since $\Sigma^{\infty}_{+}$ enhances to a (symmetric) monoidal functor, this implies that for $x=\Sigma^{\infty}_{+}x_0$, $\mathcal G_x$ contains the essential image of $\Sigma^{\infty}_{+}$ and hence $\mathcal G_x=\stab{\mt{B}_{\lambda}}_{\geq 0}$ by \cref{suspensionspectragenerate}. 
    By the same reasoning as above, \[ \mathcal G\coloneqq \{ x\in \stab{\mt{B}_{\lambda}}_{\geq 0}\, |\, \mathcal G_x=\stab{\mt{B}_{\lambda}}_{\geq 0}\}\] is closed under small colimits and extensions, whence $\mathcal G=\stab{\mt{B}_{\lambda}}_{\geq 0}$ by \cref{suspensionspectragenerate}. 
\end{proof}

In particular, if $\mt{B}_{\infty}$ is a big presentably symmetric monoidal category, then $\stab{\mt{B}_{\infty}}^{\heart}$ inherits a symmetric monoidal structure such that $\tau_{\leq 0}\colon \stab{\mt{B}_{\infty}}_{\geq 0}\to \stab{\mt{B}_{\infty}}^{\heart}$ enhances to a (symmetric) monoidal functor $\tau_{\leq 0}^{\otimes}$ by \cite[Example 2.2.1.10]{higheralgebra}.
Being right adjoint to the (symmetric) monoidal functor $\tau_{\leq 0}^{\otimes}$, the inclusion $\stab{\mt{B}_{\infty}}^{\heart}\subseteq \stab{\mt{B}_{\infty}}_{\geq 0}$ enhances to fully faithful right adjoints
\[ \aCAlg(\Ab(\tau_{\leq 0}\mt{B}_{\infty}))\to \aCAlg(\stab{\mt{B}_{\infty}}_{\geq 0}), \] cf.\ \cref{adjuncttioninducesadjunctiononalgebraobjects}.  
\begin{notation}\label{Eilenbergmaclanealgebra}
If $\mt{B}_{\infty}$ is a big topos, endow \[\Ab(\tau_{\leq 0}\mt{B}_{\infty})\cong \stab{\mt{B}_{\infty}}^{\heart}\text{ (\cref{recognitionprinciplebigtopoi})}\] with the induced symmetric monoidal structure. The Eilenberg-Mac Lane functor (\cref{defeilenbergmaclanefunctorhomotopygroups}) induces fully faithful functors 
\begin{align*}& H\colon \Alg(\Ab(\tau_{\leq 0}\mt{B}_{\infty}))\hookrightarrow \Alg(\stab{\mt{B}_{\infty}}_{\geq 0})\subseteq \Alg(\stab{\mt{B}_{\infty}}) \text{ and }\\ & H\colon \CAlg(\Ab(\tau_{\leq 0}\mt{B}_{\infty}))\hookrightarrow \CAlg(\stab{\mt{B}_{\infty}}_{\geq 0})\subseteq \CAlg(\stab{\mt{B}_{\infty}}).\end{align*} 
We will in the sequel freely identify (commutative) algebras in $\Ab(\tau_{\leq 0}\mt{B}_{\infty})$ with their images under these functors. 
\end{notation}
\begin{rem}\label{symmetricmonoidalstructureonspheartisinducedbyab}
If $\mt{B}_*\colon\Lambda\to \aCAlg(\Pr^L)$ is an chain of presentably (symmetric) monoidal categories exhausting $\mt{B}_{\infty}^{\otimes}$, then $\tau_{\leq 0}^{\otimes}\colon \stab{\mt{B}_{\infty}}_{\geq 0}^{\otimes}\to (\stab{\mt{B}_{\infty}}^{\heart})^{\otimes}$ is the colimit of the functor  \[ \id_{\mt{B}_*}\otimes_{\aCAlg(\Pr^L)} \tau_{\leq 0}\colon \mt{B}_*^{\times}\otimes_{\aCAlg(\Pr^L)}\Sp_{\geq 0}\to \mt{B}_*^{\times}\otimes_{\aCAlg(\Pr^L)}\Ab^{\otimes}\] in $\CAlg(\vlCat)$. 
In particular, for a (hyper)accessible explicit covering site $(\mathcal C,S)$, there is a symmetric monoidal equivalence \[ \mywidehatshvacc{S}(\mathcal C,\Sp)^{\heart}\cong \mywidehatshvacc{S}(\mathcal C,\Ab),\] where the right-hand side is endowed with the symmetric monoidal structure from \cref{constructionmonoidalstructureaccessiblesheaves}.  
\end{rem}
    \subsubsection{Spectral enrichment}\label{section:spectralenrichment}
        In this section, we explain that every big presentable, stable category is naturally spectrally enriched and record some basic properties of this spectral enrichment. 
     The main result is that every adjunction between big presentable, stable categories enhances to an adjunction of spectrally enriched categories.
     We will use the spectral enrichment described below to define cohomology in big topoi. 
\begin{construction}\label{Spmodulestructure}By \cite[Corollary 4.8.2.18]{higheralgebra}, the forget functor $\LMod{\Sp}{\Pr^L}\to\Pr^L$ is fully faithful with essential image the category $\Pr^L_{\operatorname{st}}$ of stable, presentable categories. A chain of stable, presentable categories $\mathcal B_*\colon\Lambda\to \Pr^L$ therefore enhances essentially uniquely to a functor $\mathcal B^{\Sp}_*\colon \Lambda\to \LMod{\Sp}{\Pr^L}$. 
By \cite[Proposition 4.8.1.15, Remark 4.8.1.8]{higheralgebra}, the forget functor $\Pr^L\to \vlCat$ enhances to a lax symmetric monoidal functor (where the left-hand side is endowed with the Lurie tensor product \cite[Proposition 4.8.1.15]{higheralgebra} and the right-hand side with the cartesian monoidal structure). This induces a functor $\LMod{\Sp}{\Pr^L}\to \LMod{\Sp}{\vlCat}$. 
Since $\vlCat$ has all large colimits and the symmetric monoidal structure on $\vlCat$ is compatible with large colimits (\cite[Lemma 4.8.4.2]{higheralgebra} in the large universe), $\LMod{\Sp}{\vlCat}$ has all large colimits and the forget functor $\LMod{\Sp}{\vlCat}\to \vlCat$ reflects large colimits by \cite[Corollary 4.2.3.5]{higheralgebra} in the very large universe. In particular, $\colim{\Lambda}\mathcal B_*^{\Sp}$ exhibits the colimit $\mathcal B_{\infty}$ as a left $\Sp$-module.

If $\mt{B}_{\infty}$ is a big presentable category which is stable and $\mt{B}_*\colon\Lambda\to \Pr^L$ is an exhaustion for $\mt{B}_{\infty}$ by presentable categories, then $\mt{B}_{\lambda}\subseteq \mt{B}_{\infty}$ is closed under finite limits and colimits for all $\lambda\in\Lambda$ and in particular stable. 
The above construction therefore yields a $\Sp$-module structure on $\mt{B}_{\infty}$.  
\end{construction}
    \begin{lemma}
        Suppose $\mt{B}_{\infty}$ is a stable, big presentable category. The $\Sp$-module structure of $\mt{B}_{\infty}$ constructed above is independent of the chosen exhaustion of $\mathcal B_{\infty}$ by presentable categories.\end{lemma}
    \begin{proof}Suppose $\mt{X}_*\colon\Lambda\to \Pr^L_{\operatorname{st}}, \mt{Y}_*\colon M\to \Pr^L_{\operatorname{st}}$ are two exhaustion for $\mt{B}_{\infty}$ by presentable stable categories. 
    By \cref{presentablestablesubcategorycontainedinfinitestage}, for $m\in M$, there exists $\kappa(m)\in \Lambda$ such that $\mt{Y}_m\to \mt{B}_{\infty}$ factors over $\mt{X}_{\kappa}\subseteq \mt{B}_{\infty}$. 
    For $(\lambda,m)\in\Lambda\times M$ choose $\kappa\in\Lambda$ with maps  $\lambda\to \kappa, \kappa(m)\to \kappa$. Then \[\mt{X}_m\times_{\mt{B}_{\infty}}\mt{Y}_{\lambda}\cong \mt{X}_m\times_{\mt{Y}_{\kappa}}\mt{Y}_{\lambda}\] is presentable. 
    As for all $\lambda\in\Lambda$, $\mt{X}_{\lambda}\subseteq \mt{B}_{\infty}$ is closed under small colimits and finite limits (\cref{fullyfaithfulnessandpreservationoflimits}), $\mt{X}_{\lambda}\times_{\mt{B}_{\infty}}\mt{Y}_{m}\subseteq \mt{Y}_m$ is closed under small colimits and finite limits, and analogously for $\mt{X}_{\lambda}\times_{\mt{B}_{\infty}}\mt{Y}_m\subseteq \mt{X}_{\lambda}$. 
    It follows that the functor 
    \[ M\times \Lambda \to \vlCat, (\lambda,m)\mapsto \mt{X}_{\lambda}\times_{\mt{B}_{\infty}}\mt{Y}_{m}\] enhances to a functor $\mt{X}_*\times_{\mt{B}_{\infty}}\mt{Y}_*\colon M\times \Lambda\to \Pr^L_{\operatorname{st}}$ which is an exhaustion of $\mt{B}_{\infty}$ by presentable categories. 
    Let \[\pi_{\Lambda}^*\mt{B}_*\colon M\times \Lambda\to \Lambda\to \Pr^L_{\operatorname{st}}\] and \[\pi_{M}^*\mt{C}_*\colon M\times \Lambda\to M\to \Pr^L_{\operatorname{st}}.\] 
    As $\LMod{\Sp}{\Pr^L}\cong \Pr^L_{\operatorname{st}}$ (\cite[Corollary 4.8.2.18]{highertopostheory}), we obtain natural transformations \[\pi_M^*\mt{C}_* \leftarrow \mt{B}_*\times_{\mt{B}}\mt{C}_*\to \pi^*\mt{B}_* \in \Fun(\Lambda\times M,\LMod{\Sp}{\Pr^L}).\] 
    Taking the colimit in $\LMod{\Sp}{\vlCat^{\times}}$ yields functors 
    \[\colim{M}\mt{C}^{\Sp}_*\leftarrow \colim{\Lambda\times M}(\mt{B}_*\times_{\mt{B}}\mt{C}_*)\to \colim{\lambda\in\Lambda}\mt{B}_{\lambda}^{\Sp}.\] By \cite[Lemma 4.8.4.2, Corollary 4.2.3.5]{higheralgebra}, the forget functor $\LMod{\Sp}{\vlCat}\to \vlCat$ is conservative and preserves large colimits. This implies that \[\colim{M}\mathcal Y_*\cong \colim{\Lambda}\mt{X}_{*}\in\LMod{\Sp}{\vlCat}.\qedhere\] 
    \end{proof}
    \begin{lemma}\label{bigpresentablespectrallyenriched}
        Suppose $\mt{B}_{\infty}$ is a stable, big presentable category. 
    The left $\Sp$-module structure on $\mt{B}_{\infty}$ described above exhibits $\mt{B}_{\infty}$ as enriched (\cite[Definition 4.2.1.28]{higheralgebra}) over $\Sp$. 
    Denote by \[\map_{\mathcal B_{\infty}}(-,-)\colon\mathcal B_{\infty}^{\operatorname{op}}\times\mathcal B_{\infty} \to \Sp\] the associated mapping spectrum functor (\cite[Remark 4.2.1.31]{higheralgebra}). 
    Then \[\Omega^{\infty}\map_{\mathcal B_{\infty}}(-,-)\cong \Map_{\mathcal B_{\infty}}(-,-).\] 
    \end{lemma}
    \begin{proof}
        For $M\in \mt{B_{\infty}}$ choose $\lambda\in\Lambda$ with $M\in \mt{B_{\lambda}}$. 
        Since $\mt{B_{\lambda}}\subseteq \mt{B_{\infty}}$ is a $\Sp$-module map, \[-\otimes M\colon \Sp\to \mathcal B_{\infty}\] factors as $-\otimes M\colon \Sp\to \mt{B}_{\lambda}\subseteq \mt{B}_{\infty}$, where $-\otimes M\colon \Sp\to \mt{B}_{\lambda}$ is induced by the left $\Sp$-tensoring of $\mt{B}_{\lambda}$. As $\mt{B}_{\lambda}^{\Sp}$ is in the image of the forget functor $\LMod{\Sp}{\Pr^L}\to \LMod{\Sp}{\vlCat}$, \[-\otimes M\colon \Sp\to \mt{B}_{\lambda}\] is a left adjoint functor. By \cref{fullyfaithfulnessandpreservationoflimits}, $\mt{B}_{\lambda}\subseteq \mt{B}_{\infty}$ is a left adjoint, which implies that \[-\otimes M\colon \Sp\to \mathcal B_{\infty}\] is a left adjoint. 
        This shows that $\mt{B}_{\infty}$ is enriched in $\Sp$. 
        By definition of the mapping spectrum functor (\cite[Remark 4.2.1.31]{higheralgebra}), \[\Map_{\Sp}(\mathbb S, \map_{\mt{B_{\infty}}}(-_1,-_2))\cong \Map_{\mt{B_{\infty}}}(\mathbb S\otimes -_1,-_2)\cong \Map_{\mt{B_{\infty}}}(-_1,-_2)\in \Fun(\mathcal B_{\infty}^{\operatorname{op}}\times\mathcal B_{\infty}, \an).\] The left-hand side is equivalent to $\Omega^{\infty}\map_{\mt{B}_{\infty}}(-_1,-_2)$.  
    \end{proof}
\begin{cor}\label{spectralenrichmentcocontinuous}
Suppose that $\mt{B_{\infty}}$ is a stable, big presentable category. 
The mapping spectrum functor 
\[ \map_{\mt{B_{\infty}}}(-,-)\colon \mt{B}_{\infty}^{\operatorname{op}}\times \mt{B}_{\infty}\to \Sp\] preserves small limits in both variables.  
\end{cor}
\begin{proof}
    For $b\in\mt{B_{\infty}}$, $\map_{\mt{B_{\infty}}}(b,-)\colon \mt{B}_{\infty}\to \Sp$ is a right adjoint and in particular preserves small limits. 
    Since $\mt{B}_{\infty}$ and $\Sp$ are stable, this implies that \[\Sigma^n\map_{\mt{B}_{\infty}}(-,-)\cong \map_{\mt{B}_{\infty}}(-, \Sigma^n-)\] for all $n\in\mathbb Z$.  
    Suppose $F\colon I\to B_{\infty}$ is some small diagram and denote by $\colim{I}F$ its colimit in $\mt{B_{\infty}}$. 
    For $b\in \mt{B}_{\infty}$, the canonical map 
    \[ \map_{\mt{B_{\infty}}}(\colim{I}F,b)\to \clim{i\in I}\map_{\mt{B_{\infty}}}(F,b)\] is an equivalence since for all $n\in\mathbb Z$, 
    \begin{align*}\Omega^{\infty-n}\map_{\mt{B_{\infty}}}(\colim{I}F,b)&\cong \Map_{\mt{B_{\infty}}}(\colim{I}F, \Sigma^nb)\\&\cong \clim{I}\Map_{\mt{B_{\infty}}}(F, \Sigma^nb)\\&\cong \clim{I}\Omega^{\infty-n}\map_{\mt{B_{\infty}}}(F,b) \\& \cong \Omega^{\infty-n}\clim{I}\map_{\mt{B_{\infty}}}(F,b)\end{align*} and the functors $\Omega^{\infty-n}\colon\Sp\to \an, n\in\mathbb N_0$ are jointly conservative. 
\end{proof}

Next, we show that every adjunction between stable, big presentable categories is spectrally enriched. 
    For a category $\mathcal T$ with small limits denote by $\Fun^{\operatorname{lim}}(\mathcal B_{\infty}^{\operatorname{op}}, \mathcal T)\subseteq \Fun(\mathcal B_{\infty}^{\operatorname{op}}, \mathcal T)$ the full subcategory on small limits-preserving functors.   
\begin{lemma}\label{spectraclosetoidempotentbigpresentable}
    Suppose $\mathcal B_{\infty}$ is a stable, big presentable category.    
    Pushforward along $\Omega^{\infty}$ defines a fully faithful functor
    \[ \Fun^{\operatorname{lim}}(\mathcal B_{\infty}^{\operatorname{op}}, \Sp)\hookrightarrow \Fun^{\operatorname{lim}}(\mathcal B_{\infty}^{\operatorname{op}}, \an).\]  
    \end{lemma}
    \begin{proof}Choose an exhaustion $\mt{B}_*\colon\Lambda\to \Pr^L$ of $\mt{B}_{\infty}$ by presentable categories. 
    Since $\mt{B}_{\infty}^{\operatorname{op}}=\colim{\lambda\in\Lambda}\mt{B}_{\lambda}^{\operatorname{op}}$ and for all $\lambda\in\Lambda$, $\mt{B}_{\lambda}\subseteq \mt{B}_{\infty}$ is closed under small colimits, for a category $\mathcal T$, we have a fully faithful functor \[\Fun^{\operatorname{lim}}(\mathcal B_{\infty}^{\operatorname{op}}, \mathcal T)\hookrightarrow\clim{\lambda\in \Lambda}\Fun^{\operatorname{lim}}(\mt{B}_{\lambda}^{\operatorname{op}}, \mathcal T).\] 
    For all $\lambda\in\Lambda$, $\mt{B}_{\lambda}\subseteq \mt{B}_{\infty}$ is closed under finite limits and colimits and in particular stable. Hence $\Omega^{\infty}$ induces an equivalence \[ \Fun^{\operatorname{lim}}(\mt{B}_{\lambda}^{\operatorname{op}}, \Sp)\cong \Fun^{\operatorname{lim}}(\mt{B}_{\lambda}^{\operatorname{op}}, \an), \] by \cite[Proposition 1.4.2.21]{higheralgebra}.
    This shows that pushforward along $\Omega^{\infty}$ defines a fully faithful functor \begin{align*}\Fun^{\operatorname{lim}}(\mathcal B_{\infty}^{\operatorname{op}}, \Sp)\hookrightarrow \clim{\Lambda}\Fun^{\operatorname{lim}}(\mt{B}_{\lambda}^{\operatorname{op}}, \Sp)& \cong \clim{\Lambda}\Fun^{\operatorname{lim}}(\mt{B}_{\lambda}^{\operatorname{op}}, \an)\\ & \hookrightarrow \clim{\Lambda}\Fun(\mt{B}_{\lambda}^{\operatorname{op}}, \an)\cong \Fun(\mt{B}_{\infty}^{\operatorname{op}}, \an).\end{align*}
        Since $\Omega^{\infty}$ preserves small limits, this factors over $\Fun^{\operatorname{lim}}(\mt{B}_{\infty}^{\operatorname{op}}, \an)\subseteq \Fun(\mt{B}_{\infty}^{\operatorname{op}}, \an)$. 
    \end{proof}
    \begin{cor}\label{adjunctionsspectrallyenriched}
        Every adjunction $L\colon \mathcal C\to \mathcal D\colon R$ between big presentable, stable categories is spectrally enriched, i.e.\ every adjunction datum \[ \Map_{\mathcal C}(L-,-)\cong \Map_{\mathcal C}(-,R-)\colon \mathcal C^{\operatorname{op}}\times\mathcal D\to \an\] lifts essentially uniquely to an equivalence of functors
        \[ \map_{\mathcal C}(L-,-)\cong \map_{\mathcal C}(-,R-)\colon \mathcal C^{\operatorname{op}}\times\mathcal D\to \Sp.\]
    \end{cor}
    \begin{proof}
        For $\mathcal T=\an, \Sp$ denote by $\Fun^{\operatorname{lim},v}(\mathcal C^{\operatorname{op}}\times\mathcal D, \mathcal T)\subseteq \Fun(\mathcal C^{\operatorname{op}}\times\mathcal D, \mathcal T)$ the category of functors which preserve small limits in both variables. 
        By \cref{spectraclosetoidempotentbigpresentable}, pushforward along $\Omega^{\infty}\colon\Sp\to \an$ induces a fully faithful functor 
        \begin{align*}\Fun^{\operatorname{lim},v}(\mathcal C^{\operatorname{op}}\times \mathcal D, \Sp)\cong \Fun^{\operatorname{lim}}(\mathcal D, \Fun^{\operatorname{lim}}(\mathcal C^{\operatorname{op}}, \Sp))\hookrightarrow & \Fun^{\operatorname{lim}}(\mathcal D, \Fun^{\operatorname{lim}}(\mathcal C^{\operatorname{op}}, \an))\\ &\cong \Fun^{\operatorname{lim},v}(\mathcal C^{\operatorname{op}}\times \mathcal D, \an).\end{align*}
        By \cref{spectralenrichmentcocontinuous} and \cref{bigpresentablespectrallyenriched}, this implies that every equivalence \[\Map_{\mathcal C}(L-,-)\cong \Map_{\mathcal D}(-,R-)\] lifts essentially uniquely to an equivalence
        $\map_{\mathcal C}(L-,-)\cong \map_{\mathcal D}(-,R-)$. 
    \end{proof}
\subsubsection{Constant sheaves}\label{section:constantsheaffunctor}

\begin{lemma}\label{existenceconstantsheaf}
    For a big presentable category $\topo{B}$, there exists a unique small colimits preserving functor $c\colon \an\to \mt{B}$ with $c(*)=*$. 
    This is a left adjoint. 
    \end{lemma} 
    \begin{proof}
        Choose an exhaustion $\topo{B}_*\colon\Lambda\to\Pr^L$ of $\topo{B}$ by presentable categories and $\lambda\in\Lambda$. 
        Since $\topo{B}_{\lambda}\subseteq \topo{B}$ is closed under small colimits and $*\in \topo{B}_{\lambda}$, every small colimits preserving functor $c\colon\an\to \topo{B}$ with $c(*)=*$ factors over a colimits preserving functor $\an\to \topo{B_{\lambda}}$. Since $\topo{B_{\lambda}}$ is cocomplete, there exists a unique such functor. This proves uniqueness and existence. 
        Since $\topo{B_{\lambda}}$ is presentable, the cocontinuous functor $\an\to \topo{B_{\lambda}}$ admits a right adjoint $\Gamma_{\lambda}$. 
        By \cref{fullyfaithfulnessandpreservationoflimits}, the inclusion $\topo{B_{\lambda}}\to \topo{B}$ admits a right adjoint $r^{\lambda}$.
        The composite $\Gamma_{\lambda}\circ r^{\lambda}$ is right adjoint to $c$. 
    \end{proof}
    \begin{definition}\label{constantsheaffunctordefinition} Suppose $\topo{B}$ is a big presentable category. The unique small colimits preserving functor $c\colon\an\to \topo{B}$ with $c(*)=*$ is called \emph{constant sheaf functor}. Its right adjoint $\Gamma$ is called \emph{global sections functor}.   
    \end{definition}

\begin{cor}
    If $\topo{X}$ is a big topos, the constant sheaf functor $c\colon \an\to \topo{X}$ is left-exact. 
\end{cor}
\begin{proof}Choose an exhaustion $\topo{X}_*\colon\Lambda\to \Pr^L$ of $\topo{X}$ by topoi and fix $\lambda\in\Lambda$. 
    Choose a small category $\mathcal C$ with a left-exact localisation $L\colon\mathcal P(\mathcal C)\to\topo{X}_{\lambda}$. 
    As $\topo{X}_{\lambda}\subseteq \topo{X}$ is closed under finite limits and small colimits, \cref{existenceconstantsheaf} implies that the constant sheaf functor factors as \[\an\xrightarrow{c_{\mathcal P(\mathcal C)}}\mathcal P(\mathcal C)\xrightarrow{L}\topo{X}_{\lambda}\to \topo{X},\] where $c_{\mathcal P(\mathcal C)}$ denotes the constant sheaf functor of $\mathcal P(\mathcal C)$.
    The functor $c_{\mathcal P(\mathcal C)}\colon \an\to \mathcal P(\mathcal C)$ is given by $X\mapsto (c\mapsto X)$ and in particular left-exact as limits in $\mathcal P(\mathcal C)$ are computed pointwise.  
\end{proof}
In particular, for a big topos $\topo{B}_{\infty}$, the constant sheaf functor enhances essentially uniquely to a symmetric monoidal functor $c^{\times}\colon \an^{\times}\to \topo{B}_{\infty}$ (\cite[Corollary 2.4.1.9]{higheralgebra}). \cref{monoidalleftadjointsstabilisetomonoidalleftadjoints} implies: 
    \begin{cor}\label{constantsheafsymmetricmonoidal}
        Suppose that $\topo{B}_{\infty}^{\otimes}$ is a big topos. 

        There is an essentially unique filler $c_{\Sp}^{\otimes}\colon \Sp^{\otimes}\to  \stab{\mt{B}_{\infty}}$ for the diagram 
        \begin{center}
            \begin{tikzcd}\an \arrow[d,"\Sigma^{\infty,\otimes}_{+}"]\arrow[r,"c^{\times}"] & \topo{B_{\infty}}^{\otimes}\arrow[d,"\Sigma^{\infty,\otimes}_{+}"']\\ 
            \Sp\arrow[r,dashed,"\exists!"',"c_{\Sp}^{\otimes}"] & \stab{\topo{B}_{\infty}}^{\otimes}
            \end{tikzcd}
        \end{center} of (symmetric) monoidal categories such that the underlying functor $c_{\Sp}\colon\Sp\to \stab{\mt{B}_{\infty}}$ is cocontinuous. The functor $c_{\Sp}^{\otimes}$ is a (symmetric) monoidal enhancement of the stabilization of the constant sheaf functor.
    \end{cor}
    For many big topoi, the constant sheaf functor commutes with totalizations of bounded above cosimplicial objects: 
    \begin{lemma}\label{constantsheafandtotalization}Suppose $\topo{X}$ is a big topos such that $\stab{\topo{X}}$ admits countable products and denote by $c_{\Sp}\colon \Sp\to\stab{\topo{X}}$ the stabilization of the constant sheaf functor. 
        If $A\colon \Delta\to \Sp_{\leq N}$ is a bounded above cosimplicial spectrum, 
        the canonical map 
        \[ c_{\Sp}(\clim{\Delta}A)\to \clim{\Delta}c_{\Sp}(A)\] is an equivalence. 
        \end{lemma}
        \begin{proof}Since $\stab{\topo{X}}$ is stable and admits countable products, it admits $\Delta$-indexed limits by \cite[Proposition 4.4.2.6]{highertopostheory}.
            As $c_{\Sp}$ is exact, for all $p\in\mathbb N_0$, 
            $c_{\Sp}\clim{\Delta_{\leq p}}A\cong \clim{\Delta_{\leq p}}c_{\Sp}A$, whence it suffices to show that the canonical map $c_{\Sp}(\clim{\Delta}A)\to \clim{p\in\mathbb N_0}c_{\Sp}(\clim{\Delta_{\leq p}}(A))$ is an equivalence. 
            By \cite[Remark 1.2.4.3]{higheralgebra}, for all $p\in\mathbb N_1$, $\Fib(\clim{\Delta_{\leq p}}A\to \clim{\Delta_{\leq p-1}}A)=\Omega^pA(p-1)$.
            Choose $N\in\mathbb N_0$ with $A\in\Fun(\Delta, \Sp_{\leq N})$. 
            Then \[\Omega^pA(p-1)\in \Sp_{\leq N-p}, \] whence 
            \[\tau_{>N-p}\clim{\Delta_{\leq p}}A\cong \tau_{>N-p}\clim{\Delta_{\leq p-1}}A.\] 
            As $c_{\Sp}$ is $t$-exact, this implies that \[\tau_{>N-p}\clim{\Delta_{\leq p}}c_{\Sp}A\cong \tau_{>N-p}\clim{\Delta_{\leq p-1}}c_{\Sp}A.\] Using $t$-exactness of $c_{\Sp}$ again, it follows that for all $m\in\mathbb Z$, \[\tau_{\geq m}c_{\Sp}\clim{\Delta}A\cong c_{\Sp}\tau_{\geq m}\clim{\Delta}A\cong c_{\Sp}(\clim{k}^{\stab{\topo{X}}_{\geq m}}\tau_{\geq m}\clim{\Delta_{\leq k}}(A))\cong c_{\Sp}\tau_{\geq m}\clim{\Delta_{\leq N-m}}A\] and \[ \tau_{\geq m}\clim{\Delta}(c_{\Sp}A)\cong \clim{k}^{\stab{\topo{X}}_{\geq m}}\tau_{\geq m}\clim{\Delta_{\leq k}}c_{\Sp}A\cong \tau_{\geq m}\clim{\Delta_{\leq N-m}}c_{\Sp}A, \] which shows that \[ \tau_{\geq m}\Fib(c_{\Sp}(\clim{\Delta}A)\to \clim{p\in\mathbb N_0}c_{\Sp}(\clim{\Delta_{\leq p}}A))=0\] for all $m\in\mathbb Z$.         
            Choose an exhaustion $\topo{X}_{*}\colon\Lambda\to \Pr^L$ of $\topo{X}$ by topoi and $\lambda\in\Lambda$ with $F\in\stab{\topo{X}_{\lambda}}$. 
            As for all $m\in\mathbb N_0$, $\stab{\topo{X}_{\lambda}}\cap \stab{\topo{X}}_{\leq m}=\stab{\topo{X}_{\lambda}}_{\leq m}$, the above implies that \[F\coloneqq \Fib(c_{\Sp}(\clim{\Delta}A)\to \clim{p\in\mathbb N_0}c_{\Sp}(\clim{\Delta_{\leq p}}A))\in \bigcap_{m\in\mathbb Z}\tau_{<m}\stab{\topo{X}_{\lambda}}.\] 
            Since the $t$-structure on $\stab{\topo{X}_{\lambda}}$ is right-separated (\cite[Proposition 1.3.2.7]{SAG}), it follows that $F=0$, i.e.\ 
             \[ c_{\Sp}(\clim{\Delta}A)\cong \clim{p\in\mathbb N_0}c_{\Sp}(\clim{\Delta_{\leq p}}A).\qedhere\] 
        \end{proof}

    \begin{lemma}\label{spectralenrichmentcomesfromconstantsheaf}Suppose that $\topo{X}$ is a big topos. The monoidal functor $c_{\Sp}^{\otimes}\colon\Sp^{\otimes}\to \stab{\topo{X}}^{\otimes}$ provided by \cref{constantsheafsymmetricmonoidal} exhibits $\stab{\topo{X}}$ as left $\Sp$-module in $\vlCat$. 
    This $\Sp$-module structure on $\stab{\topo{X}}$ agrees with the $\Sp$-module structure described in \cref{Spmodulestructure}. 
    \end{lemma}
    \begin{proof}
        If $\mt{X}_*\colon\Lambda\to\Pr^L$ is an exhaustion of $\mt{X}_{\infty}$ by topoi, $c_{\Sp}^{\otimes}$ factors over a symmetric monoidal functor $\Sp^{\otimes}\to \stab{\mt{X}_{\lambda}}^{\otimes}$ for all $\lambda\in\Lambda$. 
        This implies that $\stab{\mt{X}_*}^{\otimes}$ enhances to a functor $\stab{\mt{X}_*}^{\otimes}\colon \Lambda\to \phantom{}_{\Sp\backslash}\CAlg(\Pr^L)$. By \cite[Corollary 3.4.1.5, Corollary 4.5.1.6]{higheralgebra}, \[\phantom{}_{\Sp\backslash}\CAlg(\Pr^L)\cong \CAlg(\LMod{\Sp}{\Pr^L})\text{ and }\phantom{}_{\Sp\backslash}\CAlg(\vlCat)\cong \CAlg(\LMod{\Sp}{\vlCat}).\] 
        Denote by \[u\colon \CAlg(\LMod{\Sp}{\Pr^L})\to \LMod{\Sp}{\Pr^L}\to \LMod{\Sp}{\vlCat}\] the forget functor. 
        As $u$ factors as \[ \CAlg(\LMod{\Sp}{\Pr^L})\xrightarrow{v} \phantom{}_{\Sp\backslash}\CAlg(\vlCat)\cong \CAlg(\LMod{\Sp}{\vlCat})\xrightarrow{w} \LMod{\Sp}{\vlCat}\] (where $v$, $w$ are the forget functors), and $w$ preserves large filtered colimits (\cite[Corollary 3.2.3.1, Lemma 4.8.4.2]{higheralgebra}), 
        $\colim{\lambda\in\Lambda}u\stab{\mt{X}_{\lambda}}^{\otimes}$ is the $\Sp$-module structure on $\stab{\mt{X}_{\infty}}$ provided by the symmetric monoidal functor $c_{\Sp}^{\otimes}$. 
        As every presentable stable category admits an essentially unique $\Sp^{\otimes}$-module structure (\cite[Corollary 4.8.2.18]{higheralgebra}), 
        $\colim{\lambda\in\Lambda}u\stab{\mt{X}_{*}}^{\otimes}$ is the $\Sp$-module constructed in \cref{Spmodulestructure}.
    \end{proof}
    Denote by $\Gamma_{\Sp}\colon\stab{\topo{\mt{B}_{\infty}}}\to \Sp$ the stabilization of the global sections functor, which is right adjoint to $c_{\Sp}$. 
\begin{cor}\label{internalhomrecoversspectralenrichment}
Suppose that $\topo{B}_{\infty}^{\otimes}$ is a big topos such that the symmetric monoidal structure on $\stab{\mt{B_{\infty}}}$ induced by the cartesian monoidal structure on $\mt{B}_{\infty}$ is closed and denote by \[\imap_{\stab{\mt{B_{\infty}}}}(-,-)\colon \stab{\mt{B_{\infty}}}^{\operatorname{op}}\times \stab{\mt{B_{\infty}}}\to \Sp\] the internal Hom. 
Then $\Gamma_{\Sp}\circ \imap_{\stab{\mt{B_{\infty}}}}(-,-)\cong \map_{\stab{\mt{B_{\infty}}}}(-,-)$ recovers the $\Sp$-enrichment of $\stab{\mt{B}_{\infty}}$ from \cref{bigpresentablespectrallyenriched}.  
\end{cor}
\subsection{Module categories}\label{section:modulecategories}

In this section, we record fundamental structural results on module categories (in big presentable categories). 
Throughout, we freely use the notation of \cite[section 4.2.1]{higheralgebra}.  
We begin by recalling conditions on a symmetric monoidal category $\mathcal C^{\otimes}$ under which the category of left $A$-modules $\LMod{A}{\mathcal C}$ inherits a (closed) symmetric monoidal structure for all commutative algebras $A\in\CAlg(\mathcal C)$.
Then we establish analogues of the following two results for big presentable categories:
\begin{romanenum}\label{propertiesmodulespresentable}
    \item 
If $A$ is an algebra in a presentably monoidal category $\mathcal C$, the category of left $A$-modules is again presentable. (\cite[Corollary 4.2.3.7]{higheralgebra})
\item If $\mathcal C^{\otimes}$ is presentably symmetric monoidal, for a commutative algebra $A$ in $\mathcal C$, the category of left $A$-modules $\LMod{A}{\mathcal C}$ inherits the structure of a presentably symmetric monoidal category. (\cite[Theorem 3.4.4.2]{higheralgebra}) 
\end{romanenum} 
Next, we show that if $\mathcal C^{\otimes}$ is a big presentably monoidal category, for $A\in\Alg(\stab{\mathcal C}_{\geq 0})$, the category of left $A$-modules $\LMod{A}{\stab{\mathcal C}}$ inherits a $t$-structure such that the forget functor $\LMod{A}{\stab{\mathcal C}}\to\stab{\mathcal C}$ is $t$-exact. 
Here $\stab{\mathcal C}$ is endowed with the $t$-structure and the monoidal structure described in \cref{symmetricmonoidalstructureonspectrumobjects,tstructurespectrumobjects}, respectively. 
In \cref{section:Derivedcategories}, we recall basic facts on derived categories and, adapting \cite[Proposition 2.1.2.2]{SAG}, describe conditions on a big topos $\topo{X}$ under which \[\mathcal D(\LMod{\algebra{R}}{\Ab(\tau_{\leq 0}\topo{X})})\cong \LMod{\algebra{R}}{\stab{\topo{X}}}\] for all $\algebra{R}\in \Alg(\Ab(\tau_{\leq 0}\topo{X}))$. 
Such an equivalence is instructive since the right-hand side has more tractable categorical properties.  
Finally, we record some basic properties of module categories over constant rings in big topoi. 

\begin{proposition}[{\cite[Chapter 4]{higheralgebra}}]\label{forgetfreeadjunctionmodules}
    Suppose that $\mathcal C^{\otimes}\in \Alg(\vlCat)$ is a monoidal category and $A\in\Alg(\mathcal C^{\otimes})$ is an associative algebra. 
    \begin{romanenum}
    \item The forget functor 
    $f\colon \LMod{A}{\mathcal C}\to \mathcal C$ is conservative and reflects limits.
    \item If $K$ is a category such that $\mathcal C$ has $K$-indexed colimits and the tensor product $-\otimes_{\mathcal C}-$ preserves $K$-indexed colimits in both variables, then the forget functor preserves all $K$-indexed colimits.
    \item The forget functor has a left adjoint $A[-]$ with $f\circ A[-]=A\otimes_{\mathcal C} -$.
    \end{romanenum}
\end{proposition} 

\begin{proof}
   The forget functor is conservative by \cite[Proposition 4.8.5.8]{higheralgebra} (in the universe $\mathcal U_1$).
   By \cite[Corollary 4.2.3.3]{higheralgebra} (in $\mathcal U_1$), it preserves limits. 
   The preservation of $K$-indexed colimits is \cite[Corollary 4.2.3.5]{higheralgebra} (in $\mathcal U_1$). 
   By \cite[Corollary 4.2.4.8]{higheralgebra} (in the universe $\mathcal U_1$), the forget functor has a left adjoint $A[-]$ with $f\circ A[-]=A\otimes_{\mathcal C}-$. 
\end{proof}
\begin{recollection}
\label{symmetricmonoidalstructure}
Suppose that $\mathcal C^{\otimes}\to N(\Fin_*)$ is a symmetric monoidal category such that $\mathcal C$ has $\Delta^{\operatorname{op}}$-indexed colimits and the tensor product $-\otimes_{\mathcal C}-$ preserves $\Delta^{\operatorname{op}}$-indexed colimits in both variables. 
For a commutative algebra $A\in\CAlg(\mathcal C^{\otimes})$ denote by $\Mod{A}{\mathcal C}^{\otimes}\to N(\Fin_*)$ the operad defined in \cite[Definition 4.5.1.1]{higheralgebra}. By \cite[Theorem 4.5.2.1]{higheralgebra}, this is a symmetric monoidal category and by \cite[Corollary 4.5.1.6]{higheralgebra}, the category underlying $\Mod{A}{\mathcal C}^{\otimes}$ is equivalent to $\LMod{A}{\mathcal C}$. We endow $\LMod{A}{\mathcal C}$ with the induced symmetric monoidal structure. 

By \cite[Theorem 4.5.3.1/Remark 4.5.3.2]{higheralgebra}, a morphism $A\to B\in\CAlg(A)$ induces a strong symmetric monoidal functor $\Mod{A}{\mathcal C}^{\otimes}\to \Mod{B}{\mathcal C}^{\otimes}$. 
The symmetric monoidal structure on $\Mod{-}{-}$ is moreover natural with respect to symmetric monoidal functors: 
Suppose $\mathcal D^{\otimes}$ is another symmetric monoidal category with $\Delta^{\operatorname{op}}$-indexed colimits and the tensor product $-\otimes_{\mathcal D}-$ preserves $\Delta^{\operatorname{op}}$-indexed colimits in both variables. 
A lax symmetric monoidal functor $\phi\colon \mathcal C^{\otimes}\to\mathcal D^{\otimes}$ induces a functor $\phi^{\CAlg}\colon \CAlg(\mathcal C^{\otimes})\to \CAlg(\mathcal D^{\otimes})$. 
For $A\in\CAlg(\mathcal C)$, $\phi$ induces a functor $\phi^{A, \otimes}_{\LM}\colon\Mod{A}{\mathcal C}^{\otimes}\to \Mod{\phi^{\CAlg}A}{\mathcal D}^{\otimes}$ of categories over $N(\Fin_*)$. This is immediate from the definition of $\Mod{-}{-}^{\otimes}$.
It follows from \cite[Proposition 3.3.3.10]{higheralgebra}, that \[\phi^{A, \otimes}_{\LM}\colon\Mod{A}{\mathcal C}^{\otimes}\to \Mod{\phi^{\CAlg}A}{\mathcal D}^{\otimes}\] is lax symmetric monoidal (i.e.\ preserves inert morphisms).
By construction of $\phi^{A, \otimes}_{\LM}$, the underlying functor $\LMod{A}{\mathcal C}\to\LMod{\phi^{\CAlg}A}{\mathcal D}$ is pushforward along the lax $\LM$-monoidal functor \[\phi\times \id_{\LM}\colon\mathcal C^{\otimes}\times_{N(\Fin_*)}\LM^{\otimes}\to \mathcal D^{\otimes}\times_{N(\Fin_*)}\LM^{\otimes}.\]
\end{recollection}
\begin{lemma}\label{symmetricmonoidalstructuremodulesnatural} 
    Suppose $\mathcal C^{\otimes}, \mathcal D^{\otimes}$ are symmetric monoidal categories whose underlying categories admit $\Delta^{\operatorname{op}}$-indexed colimits and the tensor products $-\otimes_{\mathcal C}-$ and $-\otimes_{\mathcal D}-$ preserve $\Delta^{\operatorname{op}}$-indexed colimits in both variables. 
    If $\phi\colon\mathcal C^{\otimes}\to\mathcal D^{\otimes}$ is a symmetric monoidal functor such that the underlying functor $\phi\colon\mathcal C\to\mathcal D$ preserves $\Delta^{\operatorname{op}}$-indexed colimits, then the induced functor \[\phi^{A, \otimes}_{\LM}\colon\Mod{A}{\mathcal C}^{\otimes}\to \Mod{\phi^{\CAlg}A}{\mathcal D}^{\otimes}\] is symmetric monoidal. 
\end{lemma}
\begin{proof}
Since $\phi^{A, \otimes}_{\LM}\colon\Mod{A}{\mathcal C}^{\otimes}\to \Mod{\phi^{\CAlg}A}{\mathcal D}^{\otimes}$ is lax symmetric monoidal, it suffices to show that for $M,N\in\Mod{A}{\mathcal C}$, the map \[ \phi_{N,M}\colon \phi^A_{\LM}(M\otimes_A N)\to \phi^A_{\LM}(M)\otimes_{\phi^{\CAlg}A} \phi^A_{\LM}(N)\] provided by the lax symmetric monoidal structure is an equivalence.
By construction of $\phi^{A}_{\LM}$ and symmetric monoidality of $\phi$, this holds if $M$ and $N$ are free modules. 

As the forget functors $f_A\colon \LMod{A}{\mathcal C}\to \mathcal C, f_{\phi^{\CAlg}A}\colon \LMod{\phi^{\CAlg}A}{\mathcal D}\to \mathcal D$ reflect colimits, (\cref{forgetfreeadjunctionmodules}) and 
\[ \phi^A_{\LM} \circ f_{A}\cong f_{\phi^{\CAlg}A}\phi^A_{\LM}, \] $\phi^{A}_{\LM}\colon \LMod{A}{\mathcal C}\to \LMod{\phi^{\Alg}A}{\mathcal D}$ preserves $\Delta^{\operatorname{op}}$-indexed colimits.
By \cref{monoidalstructuremodulescocontinuous} below, the tensor products on $\LMod{A}{\mathcal C}$ and $\LMod{\phi^{\CAlg}A}{\mathcal D}$ preserve $\Delta^{\operatorname{op}}$-indexed colimits in both variables. 
Since every module is a $\Delta^{\operatorname{op}}$-indexed colimit of free $A$-modules (\cite[Proposition 4.7.3.14]{higheralgebra}), it now follows that \[ \phi^A_{\LM}(M\otimes_A N)\to \phi^A_{\LM}(M)\otimes_{\phi^{\CAlg}A} \phi^A_{\LM}(N)\] is an equivalence for all $M,N\in\LMod{A}{\mathcal C}$. 
\end{proof}
\begin{lemma}\label{monoidalstructuremodulescocontinuous}
Suppose $\mathcal C^{\otimes}$ is a symmetric monoidal category such that $\mathcal C$ has $\Delta^{\operatorname{op}}$-indexed colimits and $-\otimes_{\mathcal C}-$
preserves $\Delta^{\operatorname{op}}$-indexed colimits in both variables. 
If $K$ is a category such that $\mathcal C$ admits $K$-indexed colimits and the tensor product $-\otimes_{\mathcal C}-$ preserves $K$-indexed colimits in both variables, then $\LMod{A}{\mathcal C}$ has $K$-indexed colimits and $-\otimes_A-$ preserves $K$-indexed colimits in both variables. 
\end{lemma}
\begin{proof}The category $\LMod{A}{\mathcal C}$ has $K$-indexed colimits by \cite[Corollary 4.2.3.5]{higheralgebra}. 
By \cite[Corollary 4.4.2.15, Theorem 4.5.2.1.(ii)]{higheralgebra}, the tensor product preserves $K$-indexed colimits in both variables. 
\end{proof}
\begin{lemma}\label{internalhommodules}
    Suppose that $\mathcal C^{\otimes}$ is a potentially large, closed symmetric monoidal category which has all $\Delta^{\operatorname{op}}$-indexed colimits and all $\Delta$-indexed limits.

    For every commutative algebra $\algebra{R}\in \CAlg(\mathcal C)$, the induced symmetric monoidal structure (\cref{symmetricmonoidalstructure}) on $\LMod{\algebra{R}}{\mathcal C}$ is closed. 
    Denote by $\imap_{\LMod{\algebra{R}}{\mathcal C}}(-,-)$ the internal Hom and by $f\colon \LMod{\algebra{R}}{\mathcal C}\to \mathcal C$ the forget functor. 
    Then \[f\circ \imap_{\LMod{\algebra{R}}{\mathcal C}}(\algebra{R}[-],-)\cong \imap_{\mathcal C}(-,f-).\]   
\end{lemma} 
\begin{proof}We first show that the symmetric monoidal structure on $\LMod{\algebra{R}}{\mathcal C}$ is closed. 
    Choose a universe $\mathcal U$ in which $\mathcal C$ is small and denote by $\an_{+}$ the category of $\mathcal U$-small animae. Denote by $\mathcal P^{\Delta}_{+}(\mathcal C)\subseteq \Fun(\mathcal C^{\operatorname{op}}, \an_{+})$ the full subcategory of functors which preserve $\Delta$-indexed limits and by $y\colon \mathcal C\to \Fun(\mathcal C^{\operatorname{op}}, \an_{+})$ the Yoneda embedding. 
    As \[S\coloneqq \{ \colim{\Delta^{\operatorname{op}}}\,yF\to y(\colim{\Delta^{\operatorname{op}}}F)\, |\, F\in\Fun(\Delta^{\operatorname{op}}, \mathcal C)\}\] is $\mathcal U$-small and $\mathcal P^{\Delta}_{+}(\mathcal C)\subseteq \Fun(\mathcal C^{\operatorname{op}}, \an_{+})$ is the category of $S$-local objects, it is presentable (in $\mathcal U$) by \cite[Proposition 5.5.4.15]{highertopostheory} in the universe $\mathcal U$.
    By \cite[Proposition 4.8.1.10]{higheralgebra} (in the universe $\mathcal U$), there exists a symmetric monoidal structure on $\mathcal P^{\Delta}_{+}(\mathcal C)$ which is compatible with $\mathcal U$-small colimits, such that the Yoneda embedding enhances to a strong symmetric monoidal functor 
    \[ y^{\otimes}\colon \mathcal C^{\otimes}\hookrightarrow\mathcal P^{\Delta}_{+}(\mathcal C)^{\otimes}.\] 
    This induces a fully faithful functor $y_{LM}^{\algebra{R}}\colon \LMod{\algebra{R}}{\mathcal C}\hookrightarrow \LMod{y\algebra{R}}{\mathcal P^{\Delta}_{+}(\mathcal C)}$.

    By naturality of the symmetric monoidal structure on module categories (\cref{symmetricmonoidalstructuremodulesnatural}), $y_{LM}^{\algebra{R}}$ enhances to a symmetric monoidal functor $\LMod{\algebra{R}}{\mathcal C}^{\otimes}\hookrightarrow \LMod{y\algebra{R}}{\mathcal P^{\Delta}_{+}(\mathcal C)}^{\otimes}$. 
    In particular, for $M\in\LMod{\algebra{R}}{\mathcal C}$, we obtain a commutative diagram 
    \begin{center}
    \begin{tikzcd}
    \LMod{\algebra{R}}{\mathcal C}\arrow[r,"y_{LM}^{\algebra{R}}"]\arrow[d,"M\otimes_{\algebra{R}}-"] & \LMod{y\algebra{R}}{\mathcal P^{\Delta}_{+}(\mathcal C)}\arrow[d,"y_{LM}^{\algebra{R}}M\otimes_{y\algebra{R}}-"]\\ 
    \LMod{\algebra{R}}{\mathcal C}\arrow[r,"y_{LM}^{\algebra{R}}"] & \LMod{y\algebra{R}}{\mathcal P^{\Delta}_{+}(\mathcal C)}.
    \end{tikzcd}
    \end{center}
   Since $\mathcal P^{\Delta}_{+}(\mathcal U)$ is a presentable category in $\mathcal U$ and its symmetric monoidal structure is compatible with all $\mathcal U$-small colimits, the same holds for $\LMod{y\algebra{R}}{\mathcal P^{\Delta}_{+}(\mathcal C)}$ by \cite[Corollary 4.2.3.7]{higheralgebra} in the universe $\mathcal U$, i.e.\ $\LMod{y\algebra{R}}{\mathcal P^{\Delta}_{+}(\mathcal C)}$ is presentably symmetric monoidal in $\mathcal U$. 
    In particular, the right vertical functor has a right adjoint $\imap_{\LMod{y\algebra{R}}{\mathcal P^{\Delta}_{+}(\mathcal C)}}(yM,-)$. 
    We claim that $\imap_{\LMod{y\algebra{R}}{\mathcal P^{\Delta}_{+}(\mathcal C)}}(yM,-)\circ y$ factors over $\LMod{\algebra{R}}{\mathcal C}\hookrightarrow \LMod{y\algebra{R}}{\mathcal P^{\Delta}_{+}(\mathcal C)}$, then the induced functor $\LMod{\algebra{R}}{\mathcal C}\to \LMod{\algebra{R}}{\mathcal C}$ is right adjoint to $M\otimes_{\algebra{R}}-$.
    Denote by \[\mathcal L\coloneqq \{ M\in\LMod{\algebra{R}}{\mathcal C}\, |\, \imap_{\LMod{y\algebra{R}}{\mathcal P_{+}^{\Delta}(\mathcal C)}}(yM,-)\circ y\in \Fun(\LMod{\algebra{R}}{\mathcal C}, \LMod{\algebra{R}}{\mathcal C})\} \] the full subcategory on $\algebra{R}$-modules for which this holds. 
    Since $\LMod{y\algebra{R}}{\mathcal P^{\Delta}_{+}(\mathcal C)}$ has all $\mathcal U$-small limits and colimits, \[\Fun^{L}(\LMod{y\algebra{R}}{\mathcal P^{\Delta}_{+}(\mathcal C)}, \LMod{y\algebra{R}}{\mathcal P^{\Delta}_{+}(\mathcal C)})\subseteq\Fun(\LMod{y\algebra{R}}{\mathcal P^{\Delta}_{+}(\mathcal C)}, \LMod{y\algebra{R}}{\mathcal P^{\Delta}_{+}(\mathcal C)})\] and \[ \Fun^{R}(\LMod{y\algebra{R}}{\mathcal P^{\Delta}_{+}(\mathcal C)}, \LMod{y\algebra{R}}{\mathcal P^{\Delta}_{+}(\mathcal C)})\subseteq\Fun(\LMod{y\algebra{R}}{\mathcal P^{\Delta}_{+}(\mathcal C)}, \LMod{y\algebra{R}}{\mathcal P^{\Delta}_{+}(\mathcal C)})\] are closed under $\mathcal U$-small colimits and they can be computed pointwise by \cref{leftadjointsstableundercolimits} in $\mathcal U$.

    Denote by \[f_{y\algebra{R}}\colon\LMod{y\algebra{R}}{\mathcal P^{\Delta}_{+}(\mathcal C)}\to\mathcal P^{\Delta}_{+}(\mathcal C)\] the forget functor. 
    As $\LMod{y\algebra{R}}{\mathcal C}\subseteq\LMod{y\algebra{R}}{\mathcal P^{\Delta}_{+}(\mathcal C)}$ is the full subcategory on objects $M\in\LMod{y\algebra{R}}{\mathcal P^{\Delta}_{+}(\mathcal C)}$ with $f_{y\algebra{R}}M\in\mathcal C\subseteq\mathcal P^{\Delta}_{+}(\mathcal C)$, by \cref{forgetfreeadjunctionmodules}, $\LMod{\algebra{R}}{\mathcal C}\subseteq \LMod{y\algebra{R}}{\mathcal P^{\Delta}_{+}(\mathcal C)}$ is closed under small limits and $\Delta^{\operatorname{op}}$-indexed colimits. This implies that  $\mathcal L\subseteq \LMod{\algebra{R}}{\mathcal C}$ is closed under $\Delta^{\operatorname{op}}$-indexed colimits. 
    By \cite[Proposition 4.7.3.14]{higheralgebra}, every $\algebra{R}$-module $M\in\LMod{\algebra{R}}{\mathcal C}$ is a $\Delta^{\operatorname{op}}$-indexed colimit of free modules, so it suffices to show that for $c\in\mathcal C$, $\algebra{R}[c]\in\mathcal L$, then it follows that $\mathcal L=\LMod{\algebra{R}}{\mathcal C}$. Note that \[\mathcal L=\{M\in\LMod{\algebra{R}}{\mathcal C}\, | \, f_{y\algebra{R}}\circ \imap_{\LMod{y\algebra{R}}{\mathcal P^{\Delta}_{+}(\mathcal C)}}(M,-)\circ y\in\Fun(\LMod{\algebra{R}}{\mathcal C}, \mathcal C)\}.\] 

    For $c\in\mathcal C$, $y_{LM}^{\algebra{R}}\algebra{R}[c]=y\algebra{R}[yc]$ is the free $y\algebra{R}$-module on $yc$, whence 
    \[f_{y\algebra{R}}\circ \imap_{\LMod{y\algebra{R}}{\mathcal P^{\Delta}_{+}(\mathcal C)}}(y_{LM}^{\algebra{R}}\algebra{R}[c],-) \vdash (y\algebra{R}[yc]\otimes_{y\algebra{R}}-)\circ y\algebra{R}[-].\] 
    The symmetric monoidal enhancement of the free $y\algebra{R}$-module functor \[y\algebra{R}[-]\colon\mathcal P^{\Delta}_{+}(\mathcal C)\to\LMod{y\algebra{R}}{\mathcal P^{\Delta}_{+}(\mathcal C)}\] yields an equivalence \[(y\algebra{R}[yc]\otimes_{y\algebra{R}}-)\circ y\algebra{R}[-]\cong y\algebra{R}[-]\circ (yc\otimes_{\mathcal P^{\Delta}_{+}(\mathcal C)}-)\dashv \imap_{\mathcal P^{\Delta}_{+}(\mathcal C)}(c,-)\circ f_{y\algebra{R}}, \] whence \[ \imap_{\mathcal P^{\Delta}_{+}(\mathcal C)}(c,-)\circ f_{y\algebra{R}}\cong f_{y\algebra{R}}\circ \imap_{\LMod{y\algebra{R}}{\mathcal P^{\Delta}_{+}(\mathcal C)}}(y_{LM}^{\algebra{R}}\algebra{R}[c],-).\]
      Since $f_{y\algebra{R}}\circ y=y\circ f_{\algebra{R}}$, we are reduced to showing that for all $c\in\mathcal C$, \[\imap_{\mathcal P^{\Delta}_{+}(\mathcal C)}(yc,-)\circ y\in \Fun(\mathcal C, \mathcal C)\subseteq \Fun(\mathcal C, \mathcal P^{\Delta}_{+}(\mathcal C)).\]
      By fully faithfulness and symmetric monoidality of the Yoneda embedding $y$, 
      \begin{align*} \Map_{\mathcal P^{\Delta}_{+}(\mathcal C)}(-,y\imap_{\mathcal C}(c,-))\circ y& \cong \Map_{\mathcal C}(-, \imap_{\mathcal C}(c,-))\\& \cong \Map_{\mathcal C}(-\otimes c,-) \\&\cong \Map_{\mathcal P^{\Delta}_{+}(\mathcal C)}(y(-)\otimes yc,y-)\\& \cong \Map_{\mathcal P^{\Delta}_{+}(\mathcal C)}(y-, \imap_{\mathcal P^{\Delta}_{+}(\mathcal C)}(yc,-)\circ y), \end{align*} which implies that 
      \[ \imap_{\mathcal P^{\Delta}_{+}(\mathcal C)}(yc,-)\circ y\cong y\circ \imap_{\mathcal C}(c,-)\in \Fun(\mathcal C, \mathcal C)\subseteq \Fun(\mathcal C, \mathcal P^{\Delta}_{+}(\mathcal C)).\]
      This shows that the symmetric monoidal structure on $\LMod{\algebra{R}}{\mathcal C}$ is closed. 
      
    The symmetric monoidal structure of the free $\algebra{R}$-module functor $\algebra{R}[-]$ provides an equivalence \[\algebra{R}[-]\otimes_{\LMod{\algebra{R}}{\mathcal C}}\algebra{R}[-]\cong \algebra{R}[-\otimes_{\mathcal C}-], \] which defines an equivalence \[f\circ \imap_{\LMod{\algebra{R}}{\mathcal C}}(\algebra{R}[-],-)\cong \imap_{\mathcal C}(-,f-).\qedhere\] 
\end{proof}
If $(\mathcal C,S)$ is a (hyper)ac\-ces\-si\-ble explicit covering site and $\mathcal D^{\otimes}$ is a presentably symmetric monoidal category, then $\mywidehatshvacc{S}(\mathcal C, \mathcal D)$ has all small limits and colimits by \cref{limitscolimitscanbecomputedfinitestagesheaves}. 
By \cref{closednessmonoidalstructureaccessiblesheaves}, $\mywidehatshv_{S}(\mathcal C, \mathcal D)$ inherits a closed symmetric monoidal structure. 
\cref{internalhommodules} therefore implies:  
\begin{cor}\label{modulesaccessiblesheavesareclosed}
    Suppose that $(\mathcal C,S)$ is a (hyper)ac\-ces\-si\-ble explicit covering site and $\mathcal D^{\otimes}$ is a presentably symmetric monoidal category. For $\algebra{R}\in\CAlg(\mywidehatshvacc{S}(\mathcal C, \mathcal D))$, the symmetric monoidal structure on \[\LMod{\algebra{R}}{\mywidehatshvacc{S}(\mathcal C, \mathcal D)}\] is closed. 
\end{cor}
We now recall that under mild assumptions on a closed symmetric monoidal category $\mathcal C$, for an algebra $\algebra{R}\in\Alg(\mathcal C)$, the category of left $\algebra{R}$-modules is naturally enriched in $\mathcal C$. 
\begin{recollection}\label{modulecategorieslefttensored}
If $\mathcal C^{\otimes}$ is a monoidal category, for $\algebra{R}\in \Alg(\mathcal C^{\otimes})$, the category of left $\algebra{R}$-modules $\LMod{\algebra{R}}{\mathcal C}$ is right-tensored over $\mathcal C$ by \cite[Remark 4.3.3.7, Example 4.3.1.15]{higheralgebra}. 
By \cite[Proposition 4.6.3.15, Corollary 4.3.2.8]{higheralgebra},
\[\LMod{\mathcal C_{\operatorname{rev}}}{\Cat}\cong \RMod{\mathcal C}{\Cat}, \] so we can consider $\LMod{\algebra{R}}{\mathcal C}$ as left tensored over the reverse $\mathcal C^{\otimes}_{\operatorname{rev}}$ (\cite[Remark 4.1.1.7]{higheralgebra}). 
\end{recollection}
\begin{lemma}\label{modulecategoryisagainenriched}Suppose that $\mathcal C^{\otimes}$ is a closed symmetric monoidal category which has all $\Delta^{\operatorname{op}}$-indexed colimits and all $\Delta$-indexed limits. For $\algebra{R}\in \Alg(\mathcal C)$, the left-tensoring of $\LMod{\algebra{R}}{\mathcal C}$ over $\mathcal C_{\operatorname{rev}}^{\otimes}\cong\mathcal C^{\otimes}$ exhibits $\LMod{\algebra{R}}{\mathcal C}$ as enriched (\cite[Definition 4.2.1.28]{higheralgebra}) in $\mathcal C$. 
\end{lemma}
\begin{proof}
    By \cref{forgetfreeadjunctionmodules}, $\mathcal C$ and $\LMod{\algebra{R}}{\mathcal C}$ admit all $\Delta^{\operatorname{op}}$-indexed colimits and all $\Delta$-indexed limits. 
    In particular, $\Fun^{L}(\mathcal C, \LMod{\algebra{R}}{\mathcal C})\subseteq \Fun(\mathcal C, \LMod{\algebra{R}}{\mathcal C})$ is closed under $\Delta^{\operatorname{op}}$-indexed colimits and they can be computed pointwise by \cref{leftadjointsstableundercolimits}. 
    By \cite[Proposition 4.8.5.8]{higheralgebra}, the action map \[-\otimes-\colon \mathcal C\times \LMod{\algebra{R}}{\mathcal C}\to \LMod{\algebra{R}}{\mathcal C}\] preserves $\Delta^{\operatorname{op}}$-indexed colimits in both variables.  
    As every left $\algebra{R}$-module is a $\Delta^{\operatorname{op}}$-indexed colimit of free $\algebra{R}$-modules (\cite[Proposition 4.7.3.14]{higheralgebra}), it therefore suffices to show that for all free $\algebra{R}$-modules $\algebra{R}[c], c\in\mathcal C$, 
    \[-\otimes\algebra{R}[c]\colon \mathcal C\to \LMod{\algebra{R}}{\mathcal C}\] is a left adjoint.  
    By \cite[Construction 4.8.3.24]{higheralgebra}, the algebra map $1_{\mathcal C}\to\algebra{R}$ determines a $\mathcal C_{\operatorname{rev}}$-linear enhancement of the free $\algebra{R}$-module functor $\algebra{R}[-]$. 
    This yields an equivalence of $-\otimes\algebra{R}[c]$ with  
    $\mathcal C\xrightarrow{-\otimes_{\mathcal C} c}\mathcal C\xrightarrow{\algebra{R}[-]}\LMod{\algebra{R}}{\mathcal C}$. This functor is a left adjoint as the symmetric monoidal structure on $\mathcal C$ is closed. 
\end{proof}

We now show that module categories in big presentably monoidal categories behave similar as for presentable categories. More concretely, we prove analogues of the statements listed on page \pageref{propertiesmodulespresentable}. 
Suppose $\mathcal B^{\otimes}\colon\Lambda\to \aCAlg(\Pr^L)$ is a chain of presentably (symmetric) monoidal categories. By \cref{colimitmonoidalcategories}, this has a colimit in $\aCAlg(\vlCat)$. 
The (symmetric) monoidal functors $\mt{B}_{\kappa}^{\otimes}\to \mt{B}_{\infty}^{\otimes}$ induce an equivalence
 \[\colim{\kappa}\aCAlg(\mt{B_{\kappa}^{\otimes}})\cong \aCAlg(\mt{B_{\infty}^{\otimes}}).\]   
\begin{lemma}\label{filteredcolimitsmodules}
     Suppose $\mt{B_{\infty}}$ is a big presentably monoidal category. 
     \begin{romanenum}\item For $\algebra{R}\in\Alg(\mt{B_{\infty}})$, $\LMod{\algebra{R}}{\mt{B}_{\infty}}$ is a big presentable category. 
    \item If $\mt{B}_*^{\otimes}\colon\Lambda\to \Alg(\Pr^L)$ is an exhaustion by presentably monoidal categories and $\mu\in\Lambda$ is such that $\algebra{R}\in\Alg(\mt{B}_{\mu})$, then \[\LMod{\algebra{R}}{{\mt{B}_*}}\colon \phantom{}_{\mu\backslash }\Lambda\to \Pr^L\] is an exhaustion for $\LMod{\algebra{R}}{\mt{B}_{\infty}}$ by presentable categories and for $\lambda\in\phantom{}_{\mu\backslash}\Lambda$, \[\LMod{\algebra{R}}{\mt{B}_{\lambda}}\cong \LMod{\algebra{R}}{\mt{B}_{\infty}}\times_{\mt{B}_{\infty}}\mt{B}_{\lambda}.\] \end{romanenum}
\end{lemma} 
\begin{proof}
    Choose an exhaustion $\mt{B}_{*}\colon\Lambda\to \Alg(\Pr^L)$ by presentably monoidal categories and $\mu\in\Lambda$ with $\algebra{R}\in\Alg(\mt{B}_{\mu})\subseteq \Alg(\mt{B}_{\infty})$. 
    As for all $\lambda\in \phantom{}_{\mu\backslash}\Lambda$, $\mt{B}_{\lambda}\subseteq \mt{B}_{\infty}$ is a monoidal subcategory, we obtain a functor \begin{align*} \phantom{}_{\mu\backslash }\Lambda & \to \widehat{\operatorname{Cat}}_{\infty/\LMod{\algebra{R}}{\mt{B}_{\infty}}}\\ \mu \to \kappa& \mapsto \LMod{\algebra{R}}{\mt{B}_{\kappa}}\subseteq \LMod{\algebra{R}}{\mt{B}_{\infty}},\end{align*} cf. \cref{localisationinduceslocalisationonmodulecategories}. By \cref{localisationinduceslocalisationonmodulecategories,naturalitylocalizationmodulecategories} and \cref{forgetfreeadjunctionmodules}, for all $\kappa\to \lambda\in \phantom{}_{\mu\backslash}\Lambda$, 
\[\LMod{\algebra{R}}{\mt{B}_{\kappa}}\to\LMod{\algebra{R}}{\mt{B}_{\lambda}}\to \LMod{\algebra{R}}{\mt{B}_{\infty}}\] are fully faithful left adjoints and left-exact. 
\cref{naturalitylocalizationmodulecategories} moreover implies that for all $\lambda\in\Lambda_{\mu\backslash}$, the essential image of $\LMod{\algebra{R}}{\mt{B}_{\lambda}}\to \LMod{\algebra{R}}{\mt{B}_{\infty}}$ is equivalent to $\LMod{\algebra{R}}{\mt{B}_{\infty}}\times_{\mt{B}_{\infty}}\mt{B}_{\lambda}$. 

It now follows from 
\cref{filteredcolimitsofcategoriesappendix} that \[\colim{\lambda\in\phantom{}_{\mu\backslash}\Lambda}\LMod{\algebra{R}}{\mt{B}_{\lambda}}\to\LMod{\algebra{R}}{\mt{B}_{\infty}}\] is an equivalence. 
By \cite[Corollary 4.2.3.7]{higheralgebra}, $\LMod{\algebra{R}}{\mt{B_{\lambda}}}$ is presentable for all $\lambda\in\Lambda$.  
\end{proof}

\begin{recollection}Suppose $A$ is an algebra in a presentably symmetric monoidal category $\mt{B}_{\infty}^{\otimes}$ and \[\mt{B}_*^{\otimes}\colon\Lambda\to \CAlg(\Pr^L)\] is an exhaustion of $\mt{B}_{\infty}$ by presentably symmetric monoidal categories. Choose $\lambda\in\Lambda$ with $A\in\Alg(\mt{B}_{\lambda})$. 
As $\mt{B}_{\lambda}\in\CAlg(\Pr^L)$, $\LMod{\mt{B}_{\lambda}}{\Pr^L}$ inherits a symmetric monoidal structure by \cite[Remark 4.8.1.24]{higheralgebra} and \cref{symmetricmonoidalstructure}.
Recall from \cite[Corollary 3.4.1.5, Corollary 4.5.1.6]{higheralgebra} that \[\CAlg(\LMod{\mt{B}_{\lambda}}{\Pr^L})\cong \CAlg(\Pr^L)_{\mt{B}_{\lambda}\backslash}.\]
This implies that $\mt{B}_*^{\otimes}$ enhances to a diagram 
    \[\Lambda_{\lambda\backslash}\to \CAlg(\LMod{\mt{B}_{\lambda}}{\Pr^L}).\]
By \cite[section 4]{higheralgebra}/\cref{modulesarelefttensored} $\LMod{A}{\mt{B}_{\lambda}}$ is left-tensored over $\mt{B}_{\lambda}$, i.e. enhances to an element in $\LMod{\mt{B}_{\lambda}}{\Pr^L}$. 
    We therefore obtain a diagram 
     \[\mt{B_{*}}\otimes_{\LMod{\mt{B}_{\lambda}}{\Pr^L}}\LMod{A}{\mt{B}_{\lambda}}\in\Fun(\Lambda_{\lambda\backslash},\Pr^L).\] 
\end{recollection}
\begin{cor}\label{modulecategoriestensoredup}In the above situation,
\[\colim{\Lambda} (\mt{B_{*}}\otimes_{\LMod{\mt{B}_{\lambda}}{\Pr^L}}\LMod{A}{\mt{B}_{\lambda}})\cong \BMod{A}{\infty}, \] where the colimit on the left-hand side is computed in $\vlCat$.   
\end{cor}
\begin{proof}
\cite[Theorem 4.8.4.6]{higheralgebra} yields an equivalence \begin{align*}\mt{B}_*\otimes_{\LMod{\mt{B}_{\lambda}}{\Pr^L}}\LMod{A}{\mt{B}_{\lambda}}\cong \LMod{A}{\mt{B}_*}\end{align*} in $\Fun(\Lambda_{\lambda\backslash }, \Pr^L)$, note that the construction of the equivalence (\cite[Construction 4.8.4.4]{higheralgebra}) is natural in the presentably symmetric monoidal category $\mt{B}_{\kappa}$. 
The statement now follows from \cref{filteredcolimitsmodules}.
\end{proof}

\cite[Proposition 4.6.2.17]{higheralgebra} implies: 
\begin{cor}\label{existenceleftadjointpullback}
    Suppose $\mt{B}_{\infty}^{\otimes}$ is a big presentably monoidal category. 
    For an algebra morphism $\phi\colon R\to S\in\Alg(\mt{B}_{\infty})$, restriction of scalars $\phi^*\colon \LMod{S}{\mt{B}_{\infty}}\to\LMod{R}{\mt{B}_{\infty}}$ admits a left adjoint.  
\end{cor}
In case $\mt{B}_{\infty}$ admits $\Delta^{\operatorname{op}}$-indexed colimits and the tensor product $-\otimes_{\mt{B}_{\infty}}-$ preserves $\Delta^{\operatorname{op}}$-indexed colimits in both variables, this is  \cite[Proposition 4.6.2.17]{higheralgebra}. 
\begin{proof}
    Choose an exhaustion $\mt{B}_*\colon\Lambda\to\Alg(\Pr^L)$ of $\mt{B}_{\infty}$ by presentably monoidal categories such that $R,S\in\Alg(\mt{B}_{\lambda})$ for all $\lambda\in\Lambda$. (This can be achieved by passing to a suitable slice $_{\mu\backslash} \Lambda$).
    For $\kappa\to\lambda\in\Lambda$, the symmetric monoidal functor $i\colon\mt{B}_{\kappa}\to\mt{B}_{\lambda}$ yields a commutative diagram
    \begin{equation}\label{diagramphistar}
        \begin{tikzcd}
        \LMod{S}{\mt{B}_{\kappa}}\arrow[r, "i_{LM}^S"]\arrow[d,"\phi^*"]& \LMod{S}{\mt{B}_{\lambda}}\arrow[d,"\phi^*"]\\
        \LMod{R}{\mt{B}_{\kappa}}\arrow[r,"i_{LM}^R"]& \LMod{R}{\mt{B}_{\lambda}},
        \end{tikzcd}
    \end{equation} where $i_{LM}^{S,R}$ denotes the functors induced by the symmetric monoidal functor $\mt{B}_{\kappa}^{\otimes}\subseteq \mt{B}_{\lambda}^{\otimes}$. 
    Since $\mt{B}_{\lambda}$ is presentably symmetric monoidal, \cref{forgetfreeadjunctionmodules} implies that \[\phi^*_{\lambda}\colon \LMod{R}{\mt{B}_{\lambda}}\to\LMod{S}{\mt{B}_{\lambda}}\] preserves small colimits and limits for all $\lambda\in\Lambda$. 
    As $\LMod{R}{\mt{B}_{\lambda}}, \LMod{S}{\mt{B}_{\lambda}}$ are presentable (\cite[Corollary 4.2.3.7]{higheralgebra}), $\phi^*_{\lambda}$ admits a left adjoint $L_{\lambda}$ by the adjoint functor theorem.
    We claim that for all $\kappa\to\lambda\in\Lambda$, the mate \begin{equation}\label{diagramphistar1}
        \begin{tikzcd}
        \LMod{S}{B_{\kappa}}\arrow[r, "i_{LM}^S"]& \LMod{S}{\mt{B}_{\lambda}}\\
        \LMod{R}{B_{\kappa}}\arrow[r,"i_{LM}^R"]\arrow[u,"L_{\kappa}"]& \LMod{R}{\mt{B}_{\lambda}}\arrow[u,"L_{\lambda}"]\arrow[ul, Rightarrow]\end{tikzcd}
    \end{equation} of the above diagram commutes, then it follows from \cref{adjunctionsbigpresentable} that $\colim{\lambda}\, L_{\lambda}$ is left adjoint to $\phi^*=\colim{\Lambda}\, \phi^*_{\lambda}$.
    The functors $i_{LM}^{R,S}$ are left adjoints, cf.\ \cref{localisationinduceslocalisationonmodulecategories}.
    Denote by $r_{LM}^{R,S}$ their right adjoints.   
    By \cite[Remark 4.7.4.14]{higheralgebra}, it is enough to show that the opposite mate of \ref{diagramphistar1}, that is 
    \begin{center}
        \begin{tikzcd}
        \LMod{S}{B_{\kappa}}\arrow[d,"\phi^*"]\arrow[dr, Rightarrow, "\beta_{\phi}"]& \arrow[l,"r_{LM}^S"']\LMod{S}{\mt{B}_{\lambda}}\arrow[d,"\phi^*"]\\
        \LMod{R}{B_{\kappa}}& \arrow[l,"r_{LM}^R"]\LMod{R}{\mt{B}_{\lambda}}. 
        \end{tikzcd}
    \end{center} commutes. 
    Suppose first that $R=1$. By construction of the adjunction $i_{LM}^S\dashv r_{LM}^S$, $r\phi^*\cong \phi^*r_{LM}^S$ and the composition of this equivalence with the Beck-Chevalley transformation $\beta_{\phi}$ is \[\phi^*_{\kappa}(r_{LM}^S\xrightarrow{\eta r_{LM}^S} r_{LM}^Si_{LM}^Sr_{LM}^S\xrightarrow{r_{LM}^S\epsilon^S} r_{LM}^S)=\phi^*_{\kappa}(\id_{r_{LM}^S}),\] which shows that $\beta_{\phi}$ is an equivalence if $R=1$. As the forget functor $f_R^{\lambda}\colon \LMod{R}{\mt{B}_{\lambda}}\to\mt{B}_{\lambda}$ is conservative, $\beta_{\phi}$ is an equivalence if and only if ${f_{R}^{\lambda}}_*(\beta_{\phi})$ is an equivalence, and since $\beta_{1\to R}$ is an equivalence, ${f_{R}^{\lambda}}_*(\beta_{\phi})\cong \beta_{1\to R}\circ \beta_{\phi}$. 
    By the pasting law for mates (\cite[Lemma 2.2.4]{carmeli2022ambidexterity}), $\beta_{1\to R}\circ \beta_{\phi}=\beta_{1\to R\xrightarrow{\phi} S}$ which is an equivalence by the above. 
\end{proof}

We now deduce from \cite[Definition 4.5.1.1, Theorem 4.5.2.1]{higheralgebra}/\cref{symmetricmonoidalstructure} that for a commutative algebra $\algebra{R}$ in a big presentably symmetric monoidal category $\mt{B_{\infty}}$, $\LMod{\algebra{R}}{\mt{B}_{\infty}}$ enhances to a big presentably symmetric monoidal category. We can not apply \cref{symmetricmonoidalstructure} respectively \cite[Definition 4.5.1.1, Theorem 4.5.2.1]{higheralgebra} directly since big presentably monoidal categories need not admit $\Delta^{\operatorname{op}}$-indexed colimits and their tensor product need not preserve them. 
\begin{construction}\label{symmetricmonoidalstructurebigtoposmodules}
Choose an exhaustion $\mt{B}_*\colon\Lambda\to \CAlg(\Pr^L)$ by presentably symmetric monoidal categories such that $\algebra{R}\in\CAlg(\mt{B}_{\lambda})$ for all $\lambda\in\Lambda$. 
By \cite[Theorem 4.8.5.16]{higheralgebra} (in the large universe), we obtain a functor \[\LMod{\algebra{R}}{\mt{B}_{*}}^{\otimes}\colon\Lambda\to \CAlg(\vlCat^L)\to \CAlg(\vlCat)\] such that for all $\lambda\in\Lambda$, $\LMod{\algebra{R}}{\mt{B}_{\lambda}}^{\otimes}$ is the symmetric monoidal enhancement described in \cref{symmetricmonoidalstructure}. This is a symmetric monoidal enhancement of $\LMod{A}{\mt{B}_{\lambda}}$ and a presentably symmetric monoidal category by \cite[Theorem 3.4.4.2]{higheralgebra}.
For all $\kappa\to\lambda\in\Lambda$, \[\LMod{\algebra{R}}{\mt{B}_{\kappa}}\to\LMod{\algebra{R}}{\mt{B}_{\lambda}}\] is a fully faithful, left-exact left adjoint by \cref{filteredcolimitsmodules}.
This shows that the colimit \[\colim{\lambda\in\Lambda}\LMod{\algebra{R}}{\mt{B}_{\lambda}}^{\otimes}\] in $\CAlg(\vlCat)$ is a big presentably symmetric monoidal category. Since $\CAlg(\vlCat)\to\vlCat$ preserves large sifted colimits (\cite[Lemma 4.8.4.2, Corollary 3.2.3.2]{higheralgebra}), $\colim{\lambda\in\Lambda}\LMod{\algebra{R}}{\mt{B}_{\lambda}}^{\otimes}$ is a symmetric monoidal enhancement of $\LMod{\algebra{R}}{\mt{B}_{\infty}}$ by \cref{filteredcolimitsmodules} and \cite[Corollary 4.5.1.6]{higheralgebra}.\end{construction}
\begin{lemma}\label{symmetricmonoidalstructuremodulecategories}
The symmetric monoidal enhancement of $\LMod{\algebra{R}}{\mt{B}_{\infty}}$ constructed above is independent of the chosen exhaustion by presentably symmetric monoidal categories. 
\end{lemma}
\begin{proof}
    Suppose $\mt{A}_*\colon \Lambda \to \CAlg(\Pr^L), \mt{B}_*\colon J\to \CAlg(\Pr^L)$ are two exhaustions of $\mt{B_{\infty}^{\otimes}}$ by presentably symmetric monoidal subcategories such that for all $\lambda\in\Lambda$, $\algebra{R}\in\CAlg(\mt{A}_{\lambda})\subseteq \CAlg(\mt{B}_{\infty})$, and for all $j\in J$, $\algebra{R}\in\CAlg(\mt{B}_{j})\subseteq \CAlg(\mt{B}_{\infty})$. 
    For $(\lambda,j)\in\Lambda\times J$ let $C_{\lambda,j}\coloneqq \mt{A}_{\lambda}\times_{\mt{B}_{\infty}}\mt{B}_j$. 
    Since $\mt{A}_{\lambda}\to \mt{B}_{\infty}, \mt{B}_j\to \mt{B}_{\infty}$ are fully faithful, symmetric monoidal functors which preserve small colimits, $C_{\lambda,j}\hookrightarrow \mt{A}_{\lambda}$ is a symmetric monoidal subcategory closed under small colimits. As $\mt{A}_{\lambda}$ is presentably symmetric monoidal, this implies that $C_{\lambda,j}$ has all small colimits and the symmetric monoidal structure on $C_{\lambda,j}$ is cocontinuous in both variables. 
    Since $C_{\lambda,j}\subseteq \mt{B}_{\infty}^{\otimes}$ is a symmetric monoidal subcategory, \[\CAlg(C_{\lambda,j})\cong \CAlg(\mt{B}_{\infty}^{\otimes})\times_{\mt{B_{\infty}^{\otimes}}} C_{\lambda,j}, \] so in particular $\algebra{R}\in \CAlg(C_{\lambda,j})\subseteq \CAlg(\mt{B}_{\infty})$ for all $(\lambda,j)\in\Lambda\times J$.  
    Equip $\LMod{\algebra{R}}{C_{\lambda,j}}$ with the symmetric monoidal structure described in \cref{symmetricmonoidalstructure}. 
    By naturality of the symmetric monoidal structure on module categories (\cref{symmetricmonoidalstructuremodulesnatural}), $\Mod{\algebra{R}}{C_{\lambda,j}}^{\otimes}\subseteq \Mod{\algebra{R}}{\mt{A}_{\lambda}}^{\otimes}$ and $\Mod{\algebra{R}}{C_{\lambda,j}}^{\otimes}\subseteq \Mod{\algebra{R}}{\mt{B}_j}^{\otimes}$ are symmetric monoidal subcategories. 
    This implies that for all $\lambda\in\Lambda$, \[\Mod{\algebra{R}}{\mt{A}_{\lambda}}^{\otimes}\cong \colim{j\in J}\Mod{\algebra{R}}{C_{\lambda,j}}^{\otimes}\in\CAlg(\vlCat)\] and analogously for $\Mod{\algebra{R}}{{\mt{B}}_j}^{\otimes}$, which shows that \[\colim{\lambda\in\Lambda}\Mod{\algebra{R}}{\mt{A}_{\lambda}}^{\otimes}\cong \colim{(\lambda,j)\in\Lambda\times J}\Mod{\algebra{R}}{C_{\lambda,j}}^{\otimes}\cong \colim{j\in J}\Mod{\algebra{R}}{\mt{B}_j}^{\otimes}.\qedhere\] 
\end{proof}
\begin{rem}
If $\mt{B}_{\infty}$ is a big presentably symmetric monoidal category which admits $\Delta^{\operatorname{op}}$-indexed colimits and the tensor product on $\mt{B}_{\infty}$ preserves $\Delta^{\operatorname{op}}$-indexed colimits in both variables, then the symmetric monoidal structure on $\LMod{\algebra{R}}{\mt{B}_{\infty}}$ constructed above agrees with the symmetric monoidal structure described in \cref{symmetricmonoidalstructure}. 
Indeed: Choose an exhaustion $\mt{B}_{*}\colon\Lambda\to \CAlg(\Pr^L)$ of $\mt{B}_{\infty}$ by presentably symmetric monoidal categories such that $\algebra{R}\in\CAlg(\mt{B}_{\lambda})$ for all $\lambda\in\Lambda$. The symmetric monoidal functor $\mt{B}_{\lambda}\to\mt{B}_{\infty}$ is cocontinuous for all $\lambda\in\Lambda$, and hence induces a symmetric monoidal functor $\Mod{\algebra{R}}{\mt{B}_{\lambda}}^{\otimes}\to\Mod{\algebra{R}}{\mt{B}_{\infty}}^{\otimes}$ by \cref{symmetricmonoidalstructuremodulesnatural}. These functors induce a symmetric monoidal functor $\colim{\lambda\in\Lambda}\Mod{\algebra{R}}{\mt{B}_{\lambda}}^{\otimes}\to\Mod{\algebra{R}}{\mt{B}_{\infty}}$ which is an equivalence by \cref{filteredcolimitsmodules}. 
\end{rem}
\begin{cor}\label{symmetricmonoidalfunctorsfrompresentableenhancetosymmetricmonoidalfunctorsonmodulecategories}
    Suppose $\mt{B}_{\infty}^{\otimes}$ is a big presentably symmetric monoidal category, $T^{\otimes}$ is a stable, presentably symmetric monoidal category, and $F^{\otimes}\colon T^{\otimes }\to \mt{B}_{\infty}^{\otimes}$ is a cocontinuous, symmetric monoidal functor. 
    For $A\in\CAlg(T)$, 
    $F^{\otimes}$ induces a symmetric monoidal functor \[\LMod{A}{T}^{\otimes}\to \LMod{F^{\CAlg}A}{\mt{B}_{\infty}}^{\otimes}.\]
\end{cor}
\begin{proof}Choose an exhaustion $\mt{B}_*\colon\Lambda\to \CAlg(\Pr^L)$ of $\mt{B}_{\infty}$ by presentably symmetric monoidal categories.
    By \cref{presentablestablesubcategorycontainedinfinitestage}, there exists $\lambda\in\Lambda$ such that $F\colon T\to \mt{B}_{\infty}$ factors over $\mt{B}_{\lambda}\subseteq \mt{B}_{\infty}$. 
    Since $\mt{B}_{\lambda}\subseteq \mt{B}_{\infty}^{\otimes}$ is a symmetric monoidal subcategory, $F^{\otimes}$ defines a symmetric monoidal functor $T^{\otimes}\to \mt{B}_{\lambda}^{\otimes}$.
    As $\mt{B}_{\lambda}\subseteq \mt{B}_{\infty}$ is closed under small colimits and has all small colimits and $F$ is cocontinuous, \cref{symmetricmonoidalstructuremodulesnatural} implies that $F^{\otimes}$ induces a symmetric monoidal functor \[\LMod{A}{T}^{\otimes}\to\LMod{F^{\CAlg}A}{\mt{B}_{\lambda}}^{\otimes}.\] 
    By \cref{symmetricmonoidalstructuremodulecategories}, 
    $\LMod{F^{\CAlg}A}{\mt{B}_{\lambda}}^{\otimes}\subseteq \LMod{F^{\CAlg}A}{\mt{B}_{\infty}}^{\otimes}$ is a symmetric monoidal subcategory. 
    The composite \[ \LMod{A}{T}^{\otimes}\to\LMod{F^{\CAlg}A}{\mt{B}_{\lambda}}^{\otimes}\subseteq \LMod{F^{\CAlg}A}{\mt{B}_{\infty}}^{\otimes}\] is a symmetric monoidal enhancement of $F$. 
    As $\Lambda$ is filtered, this enhancement does not depend on the choice of $\lambda\in\Lambda$. 
    \cref{presentablestablesubcategorycontainedinfinitestage} and the argument from the proof of \cref{symmetricmonoidalstructuremodulecategories} imply that it also does not depend on the chosen exhaustion of $\mt{B}_{\infty}^{\otimes}$ by presentably symmetric monoidal subcategories. 
\end{proof}

    \begin{cor}\label{symmetricmonoidalstructuretensoredupmodules}
    Suppose $\mt{B_{\infty}}$ is a big presentably symmetric monoidal category and \[\mt{B}_*^{\otimes}\colon\Lambda\to \CAlg(\Pr^L)\] is an exhaustion of $\mt{B}_{\infty}$ by presentably symmetric monoidal categories. For $A\in\CAlg(\mt{B}_{\infty})$ choose $\lambda\in\Lambda$ with $A\in\CAlg(\mt{B}_{\lambda})$. Then the equivalence from \cref{modulecategoriestensoredup} enhances to an an equivalence of symmetric monoidal categories
\[\colim{\mu\in\Lambda_{\lambda\backslash}} (\mt{B_{\mu}}^{\otimes}\otimes_{\CAlg(\LMod{\mt{B}_{\lambda}}{\Pr^L})}\LMod{A}{\mt{B}_{\lambda}}^{\otimes})\cong \BMod{A}{\infty}^{\otimes}, \] where the colimit is computed in $\CAlg(\vlCat)$, and the right-hand side and $\LMod{A}{\mt{B}_{\lambda}}$ are endowed with the symmetric monoidal structure from \cref{symmetricmonoidalstructurebigtoposmodules}.
    \end{cor}
    \begin{proof}
        By construction of the symmetric monoidal structure on the right-hand side, \[\BMod{A}{\infty}^{\otimes}\cong \colim{\kappa\in\Lambda_{\lambda\backslash}}\BMod{A}{\kappa}^{\otimes},\] where the colimit is computed in $\CAlg(\vlCat)$.  
        By \cite[Corollary 3.4.4.2, Theorem 4.8.5.16]{higheralgebra}, the free $A$-module functors enhance to a symmetric monoidal natural transformation \[A^{\otimes}[-]_*\colon \mt{B}_{*}\to\LMod{A}{\mt{B}_*}\in \Fun(\Lambda_{\lambda\backslash},\CAlg(\LMod{\mt{B}_{\lambda}}{\Pr^L})),\] and by \cite[Theorem 4.8.5.16, Corollary 4.2.3.7]{higheralgebra}, $\LMod{A}{\mt{B}_*}$ enhances to a functor \[\LMod{A}{\mt{B}_*}\colon \Lambda_{\mu\backslash}\to \CAlg(\LMod{\mt{B}_{\mu}}{\Pr^L}).\] 
        As the tensor product on $\CAlg(\LMod{\mt{B}_{\lambda}}{\Pr^L})$ is cocartesian (\cite[Proposition 3.2.4.7]{higheralgebra}), $A^{\otimes}[-]_*$ and the canonical functor \[\operatorname{const}_{\LMod{A}{\mt{B}_{\mu}}}\to \LMod{A}{\mt{B_*}}\] determine a natural transformation 
        \[ c_{*}^{\otimes}\colon \LMod{A}{\mt{B}_{\lambda}}\otimes_{\CAlg(\LMod{\mt{B}_{\lambda}}{\Pr^L})}{\mt{B}_{*}} \to \LMod{A}{{\mt{B}_{*}}}\in\Fun(\Lambda_{\lambda\backslash},\CAlg(\LMod{\mt{B}_{\lambda}}{\Pr^L})).\] 
        For $\lambda\to \mu\in \Lambda_{\lambda\backslash}$, the functor $c_{\mu}^{\otimes}$ is a symmetric monoidal enhancement of the equivalence from \cite[Theorem 4.8.4.6]{higheralgebra}, i.e. the one used in the identification of \cref{modulecategoriestensoredup}. 
        As the forget functor $\CAlg(\vlCat)\to\vlCat$ preserves filtered colimits (\cite[Proposition 3.2.3.1]{higheralgebra}), this implies the statement. 
    \end{proof}
\begin{cor}\label{symmetricmonoidalenhancementfreefunctorpresentablymonoidalcat}
    Suppose $\mt{B}_{\infty}^{\otimes}$ is a big presentably monoidal category. 
    For a map of commutative algebras $\phi\colon R\to S\in\CAlg(\mt{B}_{\infty})$, restriction of scalars $\phi^*\colon \LMod{S}{\mt{B}_{\infty}}\to\LMod{R}{\mt{B}_{\infty}}$ admits a symmetric monoidal left adjoint. 

    In particular, for a commutative algebra $R\in\CAlg(\mt{B}_{\infty})$, the free $\algebra{R}$-module functor enhances to a symmetric monoidal functor $\algebra{R}[-]^{\otimes}\colon \mt{B}_{\infty}\to\LMod{\algebra{R}}{\mt{B}_{\infty}}$. 
\end{cor}
\begin{proof}
    Choose an exhaustion $\mt{B}_{*}\colon\Lambda\to \CAlg(\Pr^L)$ by presentably symmetric monoidal categories and $\lambda\in\Lambda$ with $\algebra{R},\algebra{S}\in \CAlg(\mt{B}_{\lambda})$. 
    By \cite[Theorem 4.5.3.1]{higheralgebra}, $\phi$ determines a symmetric monoidal functor $\phi_{!,\lambda}\colon \LMod{\algebra{R}}{\mt{B}_{\lambda}}\to\LMod{\algebra{S}}{\mt{B}_{\lambda}}$ which is left adjoint to restriction of scalars. 
    Tensoring with $\mt{B}_{*}$ in $\CAlg(\LMod{\mt{B}_{\lambda}}{\Pr^L})$ yields a natural transformation 
    \[\phi_{!,*}^{\otimes}\colon \mt{B}_*\otimes_{\CAlg(\LMod{\mt{B}_{\lambda}}{\Pr^L})}\LMod{\algebra{R}}{\mt{B}_{\lambda}}\to \mt{B}_{*}\otimes_{\CAlg(\LMod{\mt{B}_{\lambda}}{\Pr^L})}\LMod{\algebra{S}}{\mt{B}_{\lambda}}\] in $\Fun(\Lambda_{\lambda\backslash},\CAlg(\Pr^L)).$ 
    By \cref{symmetricmonoidalstructuretensoredupmodules}, taking the colimit over $\Lambda_{\lambda\backslash}$ in $\vlCat$ yields a symmetric monoidal functor \[ \phi^{\otimes}_{!}\colon \LMod{\algebra{R}}{\mt{B}_{\infty}}^{\otimes}\to\LMod{S}{\mt{B}_{\infty}}^{\otimes}.\] 
    Since $\CAlg(\Pr^L)\to\vlCat$ preserves filtered colimits (\cite[Proposition 3.2.3.1]{higheralgebra}), it follows from the proof of \cref{existenceleftadjointpullback} that $\phi^{\otimes}_{!}$ is a symmetric monoidal enhancement of the left adjoint of restriction of scalars. 

    We claim that the symmetric monoidal enhancement is independent of choices. 
    We first show that it is independent of the choice of $\lambda\in\Lambda$ with $\algebra{R},\algebra{S}\in\CAlg(\mt{B}_{\lambda})$. 
    Fix $\kappa\to \lambda\in\Lambda$ such that $\algebra{R},\algebra{S}\in \CAlg(\mt{B}_{\kappa})$ and denote by \[c_{\lambda}(\algebra{R}/\algebra{S})^{\otimes}\colon \LMod{\algebra{R/S}}{\mt{B}_{\kappa}}\otimes_{\CAlg(\LMod{\mt{B}_{\kappa}}{\Pr^L})}\mt{B}_{\lambda}\cong \LMod{\algebra{R/S}}{\mt{B}_{\lambda}}\] the equivalences chosen in the proof of \cref{symmetricmonoidalstructuretensoredupmodules}. 
    We will show that \begin{align}\label{independentoflambda}c_{\lambda}(S)^{\otimes}\circ (\phi^{\otimes}_{!,\kappa}\otimes \id_{\mt{B}_{\lambda}})=\phi^{\otimes}_{!,\lambda}\circ c_{\lambda}(\algebra{R})^{\otimes}.\end{align}
     As $\Lambda$ is filtered, it then follows that $\phi^{\otimes}_{!}$ is independent of the choice of $\lambda$ with $\algebra{R},\algebra{S}\in\CAlg(\mt{B}_{\lambda})$.  
     
     For $A\in\CAlg(\mt{B}_{\kappa})$, the composite $\mt{B}_{\lambda}\to \LMod{\algebra{A}}{\mt{B}_{\kappa}}\otimes \mt{B}_{\lambda}\xrightarrow{c_{\lambda}(A)^{\otimes}} \LMod{A}{\mt{B}_{\lambda}}$ is by construction of $c_{\lambda}(A)^{\otimes}$ the free $A$-module functor (with the symmetric monoidal structure induced from the unique map $1\to A\in\CAlg(\mt{B}_{\lambda})$ via \cite[Theorem 4.5.3.1]{higheralgebra}). 
     This determines an enhancement of $c_{\lambda}(A)^{\otimes}\circ (\mt{B}_{\lambda}\to \LMod{\algebra{A}}{\mt{B}_{\kappa}}\otimes_{\CAlg(\LMod{\mt{B}_{\kappa}}{\Pr^L})} \mt{B}_{\lambda})$ to a functor in  $_{\mt{B}_{\lambda}^{\otimes}\backslash}\CAlg(\LMod{\mt{B}_{\kappa}}{\Pr^L})$. 
     By \cite[Corollary 3.4.1.5]{higheralgebra}, 
     \[_{\mt{B}_{\lambda}^{\otimes}\backslash}\CAlg(\LMod{\mt{B}_{\kappa}}{\Pr^L})\cong \CAlg(\LMod{\mt{B}_{\lambda}}{\LMod{\mt{B}_{\kappa}}{\Pr^L}})\] and by \cite[Corollary 3.4.1.9]{higheralgebra}, $\LMod{\mt{B}_{\lambda}}{\LMod{\mt{B}_{\kappa}}{\Pr^L}}\cong \LMod{\mt{B}_{\lambda}}{{\Pr^L}}$ as symmetric monoidal categories, whence 
     \[_{\mt{B}_{\lambda}^{\otimes}\backslash}\CAlg(\LMod{\mt{B}_{\kappa}}{\Pr^L})\cong \CAlg(\LMod{\mt{B}_{\lambda}}{\Pr^L}).\] 
    In particular, \[ c_{\lambda}(S)^{\otimes}\circ (\phi^{\otimes}_{!,\kappa}\otimes \id_{\mt{B}_{\lambda}})\circ (\mt{B}_{\lambda}\to \LMod{\algebra{R}}{\mt{B}_{\kappa}}\otimes_{\CAlg(\LMod{\mt{B}_{\kappa}}{\Pr^L})}\mt{B}_{\lambda})\] and \[\phi^{\otimes}_{!,\lambda}\circ c_{\lambda}(\algebra{R})^{\otimes}\circ (\mt{B}_{\lambda}\to \LMod{\algebra{R}}{\mt{B}_{\kappa}}\otimes_{\CAlg(\LMod{\mt{B}_{\kappa}}{\Pr^L})}\mt{B}_{\lambda})\] enhance to functors in \[_{\mt{B}_{\lambda}^{\otimes}\backslash}\CAlg(\LMod{\mt{B}_{\kappa}}{\Pr^L})\cong \CAlg(\LMod{\mt{B}_{\lambda}}{\LMod{\mt{B}_{\kappa}}{\Pr^L}})\cong \CAlg(\LMod{\mt{B}_{\lambda}}{\Pr^L}).\] 
     As $\mt{B}_{\lambda}$ is the initial object of $\CAlg(\LMod{\mt{B}_{\lambda}}{\Pr^L})$, it follows that \begin{align*}  c_{\lambda}(S)^{\otimes}\circ (\phi^{\otimes}_{!,\kappa}\otimes \id_{\mt{B}_{\lambda}})\circ (\mt{B}_{\lambda}\to \LMod{\algebra{R}}{\mt{B}_{\kappa}}\otimes_{\CAlg(\LMod{\mt{B}_{\kappa}}{\Pr^L})}\mt{B}_{\lambda})\\ = \phi^{\otimes}_{!,\lambda}\circ c_{\lambda}(\algebra{R})^{\otimes}\circ (\mt{B}_{\lambda}\to \LMod{\algebra{R}}{\mt{B}_{\kappa}}\otimes_{\CAlg(\LMod{\mt{B}_{\kappa}}{\Pr^L})}\mt{B}_{\lambda}).\end{align*} 

     For $A\in\CAlg(\mt{B}_{\kappa})$, the symmetric monoidal functor $i_{\kappa}^{\lambda}\colon \mt{B}_{\kappa}\to\mt{B}_{\lambda}$ induces a symmetric monoidal left adjoint \[i_{\kappa}^{\lambda}(A)^{\otimes}\colon \LMod{\algebra{A}}{\mt{B}_{\kappa}}\to\LMod{\algebra{A}}{\mt{B}_{\lambda}}\, \, \text{(\cite[Theorem 4.8.5.16]{higheralgebra})},\] and by construction of $c_{\lambda}(\algebra{A})^{\otimes}$, the composite 
     \[\LMod{\algebra{A}}{\mt{B}_{\kappa}}\to \mt{B}_{\lambda}\otimes_{\CAlg(\LMod{\mt{B}_{\kappa}}{\Pr^L})}\LMod{\algebra{A}}{\mt{B}_{\kappa}}\xrightarrow{c_{\lambda}(A)^{\otimes}} \LMod{\algebra{A}}{\mt{B}_{\lambda}}\] is $i_{\kappa}^{\lambda}(A)^{\otimes}$. \cite[Theorem 4.8.5.16]{higheralgebra} now implies that \begin{center}
     \begin{tikzcd}
     \LMod{\algebra{R}}{\mt{B}_{\kappa}}\arrow[r]\arrow[d,"\phi^{\otimes}_{!,\kappa}"] &\mt{B}_{\lambda}\otimes_{\CAlg(\LMod{\mt{B}_{\kappa}}{\Pr^L})}\LMod{\algebra{R}}{\mt{B}_{\kappa}}\arrow[r,"c_{\kappa}(\algebra{R})^{\otimes}"] & \LMod{\algebra{R}}{\mt{B}_{\lambda}}\arrow[d,"{\phi}^{\otimes}_{!,\lambda}"]\\ 
     \LMod{\algebra{S}}{\mt{B}_{\kappa}}\arrow[r] & \mt{B}_{\lambda}\otimes_{\CAlg(\LMod{\mt{B}_{\kappa}}{\Pr^L})}\LMod{\algebra{S}}{\mt{B}_{\kappa}}\arrow[r,"c_{\kappa}(\algebra{S})^{\otimes}"] & \LMod{\algebra{S}}{\mt{B}_{\lambda}}
     \end{tikzcd}
     \end{center} commutes.
     As the symmetric monoidal structure on $\CAlg(\LMod{\mt{B}_{\kappa}}{\Pr^L})$ is cocartesian, this implies (\ref{independentoflambda}).      
    The same argument as in the proof of \cref{symmetricmonoidalstructuremodulecategories} shows that the symmetric monoidal enhancement $\phi^{\otimes}_{!}$ of $\phi_{!}$ is also independent of the chosen exhaustion of $\mt{B}_{\infty}^{\otimes}$ by presentably symmetric monoidal categories.   
\end{proof}
 
\begin{cor}\label{modulesinmodulecategoriessymmetricmonoidalbigpresentable}
Suppose $\phi\colon A\to B\in\CAlg(\topo{X})$ is a morphism of commutative algebras in a big presentably symmetric monoidal category $\topo{X}$. 
By \cite[Corollary 3.4.1.5, Corollary 4.5.1.6]{higheralgebra}, this exhibits $B\in\CAlg(\LMod{A}{\topo{X}})$.
The forget functor $\LMod{A}{\topo{X}}\to \topo{X}$ induces a symmetric monoidal equivalence 
\[ \LMod{B}{\LMod{A}{\topo{X}}}\to \LMod{B}{\topo{X}}.\]  
\end{cor}
\begin{proof}The symmetric monoidal structure on the free $A$-module functor (\cref{symmetricmonoidalenhancementfreefunctorpresentablymonoidalcat}) determines a lax symmetric monoidal structure of the forget functor \[f\colon\LMod{A}{\topo{X}}\to \topo{X},\] cf. \cref{adjunctionfiberwise}. 
By construction of the symmetric monoidal structure (\cref{symmetricmonoidalstructure}/\cite[Definition 3.3.8]{higheralgebra}) on module categories, this induces a lax symmetric monoidal functor \[f^{\otimes}\colon \LMod{B}{\LMod{A}{\stab{\topo{X}}}}\to \LMod{B}{\stab{\topo{X}}}.\] 
We claim that $f^{\otimes}$ is a symmetric monoidal equivalence. Choose an exhaustion of $\topo{X}$ by presentably symmetric monoidal subcategories $\topo{X}_*$ such that $A,B\in\CAlg(\topo{X}_{\lambda})$ for all $\lambda\in\Lambda$. For all $\lambda\in\Lambda$, $f^{\otimes}$ restricts to \[f^{\otimes}_{\lambda}\colon  \LMod{B}{\LMod{A}{\topo{X}_{\lambda}}}\to \LMod{B}{\topo{X}_{\lambda}},\] the functor induced by the lax monoidal structure of the forget functor 
\[ \LMod{A}{\topo{X}_{\lambda}}\to \topo{X}_{\lambda}.\] This is immediate from the construction (\cref{adjunctionfiberwise}).
By \cite[Theorem 3.4.1.9]{higheralgebra}, $f^{\otimes}_{\lambda}$ is strong symmetric monoidal and an equivalence.
The statement now follows from \cref{symmetricmonoidalstructurebigtoposmodules}. 
\end{proof}
Next, we show that if $\mt{B}$ is a big presentably monoidal category, for $\algebra{R}\in \Alg(\stab{\mt{B}}_{\geq 0})$, the category of left $\algebra{R}$-modules in $\stab{\mt{B}}$ inherits a $t$-structure from $\stab{\mt{B}}$ (with the $t$-structure described in \cref{tstructurespectrumobjects}). 
\begin{lemma}\label{tstructureonmodulespresentable}
    Suppose that $\mathcal C^{\otimes}$ is a stable, presentably symmetric monoidal category with an accessible $t$-structure $(\mathcal C_{\geq 0}, \mathcal C_{\leq 0})$ which is compatible with the symmetric monoidal structure, i.e.\ $\mathcal C_{\geq 0}\subseteq \mathcal C$ is a symmetric monoidal subcategory and presentable. 
    \begin{romanenum}
    \item For $\algebra{R}\in \Alg(\mathcal C_{\geq 0})$, 
    \begin{align*}\LMod{\algebra{R}}{\mathcal C}_{\geq 0}&\coloneqq \LMod{\algebra{R}}{\mathcal C}\times_{\mathcal C}\mathcal C_{\geq 0}, \\
        \LMod{\algebra{R}}{\mathcal C}_{\leq 0}&\coloneqq \LMod{\algebra{R}}{\mathcal C}\times_{\mathcal C}\mathcal C_{\leq 0}\end{align*} define a $t$-structure on $\LMod{\algebra{R}}{\mathcal C}$. 

    \item If $\algebra{R}\in\CAlg(\mathcal C_{\geq 0})$ is a commutative algebra, then this is compatible with the symmetric monoidal structure on $\LMod{\algebra{R}}{\mathcal C}$ described in \cref{symmetricmonoidalstructure}, i.e.\  
        $\LMod{\algebra{R}}{\mathcal C}_{\geq 0}\subseteq \LMod{\algebra{R}}{\mathcal C}$ is a symmetric monoidal subcategory.
    \end{romanenum} 
\end{lemma}
\begin{proof}
    The proof is completely analogous to \cite[Proposition 2.1.1.1]{SAG}. 
    Since $\mathcal C_{\geq 0}\subseteq \mathcal C$ is a symmetric monoidal subcategory, $\LMod{\algebra{R}}{\mathcal C_{\geq 0}}\cong \LMod{\algebra{R}}{\mathcal C}_{\geq 0}$. As $\mathcal C_{\geq 0}$ and $\mathcal C$ are presentable, so are $\LMod{\algebra{R}}{\mathcal C}$ and $\LMod{\algebra{R}}{\mathcal C_{\geq 0}}$ by \cite[Corollary 4.2.3.7]{higheralgebra}.  
    Since the forget functor \[f\colon \LMod{\algebra{R}}{\mathcal C}\to\mathcal C\] is exact and creates colimits (\cref{forgetfreeadjunctionmodules}), $\LMod{\algebra{R}}{\mathcal C}_{\geq 0}$ is the connective part of a $t$-structure on $\LMod{\algebra{R}}{\mathcal C}$ by \cite[Proposition 1.4.4.11]{higheralgebra}. Denote by $\LMod{\algebra{R}}{\mathcal C}_{\tilde{\leq} 0}$ the coconnective part of this $t$-structure. For $T\in \mathcal C_{\geq 0}$, $\algebra{R}[T]\in \LMod{\algebra{R}}{\mathcal C}_{\geq 0}$ since $f\algebra{R}[T]=\algebra{R}\otimes T\in \mathcal C_{\geq 0}$. This implies that for $B\in \LMod{\algebra{R}}{\mathcal C}_{\tilde{\leq 0}}$, \[\Map_{\LMod{\algebra{R}}{\mathcal C}}(\algebra{R}[T],B)\cong \Map_{\mathcal C}(T,fB)\in \tau_{\leq 0}\an, \] i.e.\ $fB\in \mathcal C_{\leq 0}$, which shows that $\LMod{\algebra{R}}{\mathcal C}_{\tilde{\leq}0}\subseteq \LMod{\algebra{R}}{\mathcal C}_{\leq 0}$. Conversely, for $B\in \LMod{\algebra{R}}{\mathcal C}_{\leq 0}$, denote by $\mathcal G_B\subseteq \LMod{\algebra{R}}{\mathcal C}_{\geq 0}$ the full subcategory on objects $G$ such that \[\Map_{\LMod{\algebra{R}}{\mathcal C}}(G,B)\in\tau_{\leq 0}\an.\] Since $\tau_{\leq 0}\an$ is closed under limits, $\mathcal G_B$ is closed under colimits. By definition, $\mathcal G_B$ contains all free modules $\algebra{R}[T],T\in \mathcal C_{\geq 0}$. As those generate $\LMod{\algebra{R}}{\mathcal C}_{\geq 0}\cong \LMod{\algebra{R}}{\mathcal C_{\geq 0}}$ under $\Delta^{\operatorname{op}}$-indexed colimits (\cite[Proposition 4.7.3.14]{higheralgebra}), it follows that $\mathcal G_B=\LMod{\algebra{R}}{\mathcal C}_{\geq 0}$, which proves that 
    $\LMod{\algebra{R}}{\mathcal C}_{\leq 0}= \LMod{\algebra{R}}{\mathcal C}_{\tilde\leq 0}$. 

    Suppose now that $\algebra{R}\in\CAlg(\mathcal C_{\geq 0})$ is a commutative algebra.
    As the symmetric monoidal functor $\mathcal C_{\geq 0}^{\otimes}\subseteq \mathcal C^{\otimes}$ preserves colimits, it induces a symmetric monoidal functor \[\LMod{\algebra{R}}{\mathcal C_{\geq 0}}^{\otimes}\to\LMod{\algebra{R}}{\mathcal C}^{\otimes}\] by \cref{symmetricmonoidalstructuremodulesnatural}.
    The underlying functor $\LMod{\algebra{R}}{\mathcal C_{\geq 0}}\to\LMod{\algebra{R}}{\mathcal C}$ is fully faithful with essential image $\LMod{\algebra{R}}{\mathcal C}_{\geq 0}$, which shows that the latter is a symmetric monoidal subcategory.
\end{proof}
\begin{cor}\label{tstructureonmodules}
    Suppose that $\mt{B}^{\otimes}_{\infty}$ is a big presentably symmetric monoidal category and endow $\stab{\mt{B_{\infty}}}$ with the induced symmetric monoidal structure (\cref{symmetricmonoidalstructureonspectrumobjects}) and the $t$-structure described in \cref{tstructurespectrumobjects}. 
    For $\algebra{R}\in \Alg(\Sp(\mt{B_{\infty}})_{\geq 0})$, the subcategories \begin{align*}\BMod{\algebra{R}}{\infty}_{\geq 0}&\coloneqq \BMod{\algebra{R}}{\infty}\times_{\stab{\mt{B_{\infty}}}}\stab{\mt{B_{\infty}}}_{\geq 0}, \\
    \BMod{\algebra{R}}{\infty}_{\leq 0}&\coloneqq \BMod{\algebra{R}}{\infty}\times_{\stab{\mt{B_{\infty}}}}\stab{\mt{B_{\infty}}}_{\leq 0}\end{align*} constitute a $t$-structure on $\BMod{\algebra{R}}{\infty}$. 
    \end{cor}
    \begin{proof}Choose an exhaustion $\mt{B}_*\colon\Lambda\to \Alg(\Pr^L)$ of $\mt{B}_{\infty}$ by presentably monoidal categories such that $\algebra{R}\in\CAlg(\stab{\mt{B}_{\lambda}})$ for all $\lambda\in\Lambda$. Since for all $\lambda\in\Lambda$, the $t$-structure on $\stab{\mt{B}_{\lambda}}$ is accessible (\cite[Proposition 1.4.3.4]{higheralgebra}) and compatible with the monoidal structure (\cref{monoidalstructurecompatibletstructure}), \[(\BMod{\algebra{R}}{\lambda}_{\geq 0}, \BMod{\algebra{R}}{\lambda}_{\leq 0})\] defines a $t$-structure on $\BMod{\algebra{R}}{\lambda}$ by \cref{tstructureonmodulespresentable}. 
    
    By \cref{tstructurespectrumobjects}, for all $\lambda\in\Lambda$, \[\stab{\mt{B}_{\lambda}}_{\geq 0}= \stab{\mt{B}_{\lambda}}\cap \stab{\mt{B}_{\infty}}_{\geq 0}\text{ and }\stab{\mt{B}_{\lambda}}_{\leq 0}= \stab{\mt{B}_{\lambda}}\cap \stab{\mt{B}_{\infty}}_{\leq 0}, \] whence \[\BMod{\algebra{R}}{\lambda}\cap \BMod{\algebra{R}}{\infty}_{\geq 0}=\BMod{\algebra{R}}{\lambda}_{\geq 0}\] and \[\BMod{\algebra{R}}{\lambda}\cap \BMod{\algebra{R}}{\infty}_{\leq 0}=\BMod{\algebra{R}}{\lambda}_{\leq 0}.\] Using \cref{filteredcolimitsmodules} and \cref{fullyfaithfulnessandpreservationoflimits}, it is now straightforward to check that \[(\BMod{\algebra{R}}{\infty}_{\geq 0}, \BMod{\algebra{R}}{\infty}_{\leq 0})\] defines a $t$-structure on $\BMod{\algebra{R}}{\infty}.$
    \end{proof}
    \begin{cor}\label{symmetricmonoidalstructurestablemodulecatcompatibletstructure}Suppose that $\mt{B}_{\infty}^{\otimes}$ is a big presentably symmetric monoidal category. 
    For a connective ring spectrum object $\algebra{R}\in\CAlg(\stab{\mt{B_{\infty}}}_{\geq 0})$, the $t$-structure from \cref{tstructureonmodules} is compatible with the symmetric monoidal structure from \cref{symmetricmonoidalstructurebigtoposmodules}, i.e.\ \[\BMod{\algebra{R}}{\infty}_{\geq 0}\subseteq \BMod{\algebra{R}}{\infty}\] is a symmetric monoidal subcategory. 
\end{cor}
\begin{proof}Choose an exhaustion $\mt{B}^{\otimes}_*\colon\Lambda\to \CAlg(\Pr^L)$ of $\mt{B}^{\otimes}_{\infty}$ by presentably symmetric monoidal categories. 
Since $\algebra{R}\in\stab{\mt{B_{\infty}}}_{\geq 0}$, $\BMod{\algebra{R}}{\infty}_{\geq 0}$ contains the unit of the symmetric monoidal structure on $\BMod{\algebra{R}}{\infty}$. 
For $M,N\in\BMod{\algebra{R}}{\infty}_{\geq 0}$ choose $\kappa\in\Lambda$ with \[\algebra{R}\in\CAlg(\stab{\mt{B_{\kappa}}})\subseteq \CAlg(\stab{\mt{B_{\infty}}})\] and $M,N\in\BMod{\algebra{R}}{\kappa}\subseteq \BMod{\algebra{R}}{\infty}$.
Then $M,N\in\BMod{\algebra{R}}{\kappa}_{\geq 0}$. Hence by \cref{tstructureonmodulespresentable}, $M\otimes_{\LMod{\algebra{R}}{\mt{B}_{\kappa}}}N\in \BMod{\algebra{R}}{\kappa}_{\geq 0}$. 
Since \[\BMod{\algebra{R}}{\kappa}\subseteq \BMod{\algebra{R}}{\infty}\] is symmetric monoidal and $t$-exact, this implies that \[M\otimes_{\algebra{R}}N\in \BMod{\algebra{R}}{\infty}_{\geq 0}.\qedhere\] 
\end{proof}

If $\mt{B}_{\infty}$ is a big topos and $\algebra{R}\in\Alg(\Ab(\tau_{\leq 0}\mt{B}_{\infty}))$, the Eilenberg-Mac Lane functor is symmetric monoidal (\cref{Eilenbergmaclanealgebra}) and therefore induces a functor \[\LMod{\algebra{R}}{\Ab(\tau_{\leq 0}\mt{B}_{\infty})}\to \LMod{\algebra{R}}{\stab{\mt{B}_{\infty}}_{\geq 0}}\subseteq \LMod{\algebra{R}}{\stab{\mt{B}_{\infty}}}.\]
\begin{lemma}\label{identificationheart}
This factors over an equivalence 
\[\LMod{\algebra{R}}{\Ab(\tau_{\leq 0}\mt{B}_{\infty})}\cong \LMod{\algebra{R}}{\stab{\mt{B}_{\infty}}}^{\heart}.\] 
\end{lemma}
\begin{proof}
    By definition of the $t$-structure on $\LMod{\algebra{R}}{\stab{\mt{B}_{\infty}}}$, the Eilenberg-Mac Lane functor factors over $\LMod{\algebra{R}}{\stab{\mt{B}_{\infty}}}^{\heart}$. 
    The symmetric monoidal functor 
    \[\pi_0^{\mt{B}_{\infty}}\colon \stab{\mt{B}_{\infty}}_{\geq 0}\to \stab{\mt{B}_{\infty}}^{\heart}\cong \Ab(\tau_{\leq 0}\mt{B}_{\infty})\] induces a functor \[ \LMod{\algebra{R}}{\stab{\mt{B}_{\infty}}}^{\heart}\hookrightarrow \LMod{\algebra{R}}{\stab{\mt{B}_{\infty}}}_{\geq 0}\cong \LMod{\algebra{R}}{\stab{\mt{B}_{\infty}}_{\geq 0}}\xrightarrow{\pi_0^{\mt{B}_{\infty}}} \LMod{\algebra{R}}{\Ab(\tau_{\leq 0}\mt{B}_{\infty})}\] which is inverse to the functor induced by the Eilenberg-Mac Lane functor. 
\end{proof}
\subsubsection{Derived categories}\label{section:Derivedcategories}
We now review derived categories of abelian categories. 
Adapting \cite[Theorem 2.1.2.2]{SAG}, we formulate conditions on a big topos $\topo{X}$ under which $\mathcal D(\LMod{\algebra{R}}{\Ab(\tau_{\leq 0}\topo{X})})\cong \LMod{\algebra{R}}{\stab{\topo{X}}}$ for all discrete rings $\algebra{R}\in \Alg(\Ab(\tau_{\leq 0}\topo{X}))$, see \cref{derivedcategorymodules}. Such an equivalence is instructive as, by our above discussion, the right-hand side enjoys good categorical properties. 
\begin{definition}Consider $\mathbb Z$ as poset via $\geq$.
    For a possibly large abelian category $\mathcal A$ denote by $\Ch(\mathcal A)\subseteq \Fun(\mathbb Z, \mathcal A)$ the (large) category of chain complexes in $\mathcal A$.
    
    A chain map $f\colon C_*\to D_*\in \Ch(\mathcal A)$ is a homology equivalence if $H_i(f)$ is an isomorphism in $\mathcal A$ for all $i\in\mathbb Z$. 
    
    Define the \emph{derived category} of $\mathcal A$ as 
    \[ \mathcal D(\mathcal A)\coloneqq \Ch(\mathcal A)[W^{-1}]\] 
    the localization of $\Ch(\mathcal A)$ at the homology equivalences in the (very large) category $\vlCat$ of large categories. 
\end{definition}
\begin{rem}If $\mathcal U_1\in\tilde{\mathcal U}_1$ are universes, the inclusion $\Cat({\mathcal U}_1)\subseteq \Cat({\tilde{\mathcal U}_1})$ of $\mathcal U_1$-small categories into $\tilde{\mathcal U}_1$-small categories preserves $\mathcal U_1$-small colimits. 
In particular, for a $\mathcal U_1$-small abelian category $A$, it sends the derived category of $A$ (computed in $\mathcal U_1$) to the derived category of $A$ (computed in $\mathcal U_2$). 
Indeed: For a universe $\mathcal V$ and a $\mathcal V$-small abelian category $A$, the derived category $\mathcal D(A)$ computed in $\mathcal V$ is the pushout of \begin{center}\begin{tikzcd}
W \arrow[r]\arrow[d] & \Ch(A)\\ 
\mid W\mid 
\end{tikzcd}\end{center} in $\Cat(\mathcal V)$, where $W\subseteq \Ch(A)$ denotes the wide subcategory on weak equivalences and $|W|$ is its \textit{geometric realization}, i.e. the image of $W$ under the left adjoint of the inclusion $\an_{\mathcal V}\to \Cat(\mathcal V)$ of $\mathcal V$-small animae into $\mathcal V$-small categories. 
Denote by $J$ the nerve of the 1-category with two objects $\{0,1\}$ and a unique isomorphism between them. 
There is a canonical map $\Delta^1\to J$, and $|W|$ is the pushout 
\begin{center}
\begin{tikzcd}
\sqcup_{w\in \Fun(\Delta^1,W)}\Delta^1\arrow[r] \arrow[d]& W\\
\sqcup_{w\in \Fun(\Delta^1,W)}J 
\end{tikzcd}\end{center} in $\Cat(\mathcal V)$, see e.g.\ \cite[Lemma 2.4.6]{landinfinity}. 
\end{rem}
The following is again folklore, but we could not find a reference. 
\begin{lemma}\label{derivedcategoryisstable}
    Suppose $\mathcal A$ is a (large) abelian category. 
    The (large) category $\mathcal D(\mathcal A)$ is stable and cofibers can be computed using the mapping cone construction on the level of chain complexes. 
\end{lemma}
\begin{proof} We thank Stefan Schwede for suggesting the following argument. 
By \cite[Definition 1.3.2.1]{higheralgebra} (in the large universe), the category $\Ch(\mathcal A)$ is a differential graded category and by \cite[Proposition 1.3.4.5]{higheralgebra}, \[\Ch(\mathcal A)[ \text{(chain homotopy equivalences)}^{-1}]\cong \operatorname{N}_{\operatorname{dg}}(\Ch(\mathcal A))\] is the \emph{differential graded nerve} (\cite[Construction 1.3.1.6]{higheralgebra}) of $\Ch(A)$. 
Since chain homotopy equivalences are homology equivalences, the homology functors descend to functors \[H_i\colon \operatorname{N}_{\operatorname{dg}}(\Ch(\mathcal A))\to \mathcal A\] and in particular, \[\mathcal D(\mathcal A)\cong \operatorname{N}_{\operatorname{dg}}(\Ch(\mathcal A))[(\text{homology equivalences})^{-1}].\] 
The category $\operatorname{N}_{\operatorname{dg}}(\Ch(\mathcal A))$ is stable by \cite[Proposition 1.3.2.10]{higheralgebra} (in the large universe), and by \cite[Remark 1.3.2.17]{higheralgebra} (in the large universe): if $a\to b\in \operatorname{N}_{\operatorname{dg}}(\Ch(\mathcal A))$ is a map represented by a map of chain complexes $f\colon A\to B\in\Ch(\mathcal A)$, then the mapping cone $\Cone(f)$ (computed in the category of chain complexes) represents the cofiber $\Cofib(f)$.
In particular, a map $f\colon A\to B$ in $\operatorname{N}_{\operatorname{dg}}(\Ch(\mathcal A))$ is a homology equivalence if and only if $H_i(\Cofib(f))=0$ for all $i\in\mathbb Z$. 
This implies that \[ \mathcal N\coloneqq \{ N\in \operatorname{N}_{\operatorname{dg}}(\Ch(A))\, |\, H_i(N)=0 \text{ for all }i\in\mathbb Z\}\subseteq \operatorname{N}_{\operatorname{dg}}(\Ch(A))\] is a stable subcategory (\cite[Definition I.3.2]{NikolausScholze}). 
Indeed, for $N\in\mathcal N$, the cofiber sequence $\Omega N\to *\to N$ implies that $\Omega N\to *$ is a homology equivalence, i.e.\ $\Omega N\in\mathcal N$. 
Suppose that $A\xrightarrow{f}B\to C\in \operatorname{N}_{\operatorname{dg}}(\mathcal A)$ is a cofiber sequence. 
If $A,B\in \mathcal N$, then $f$ is a homology equivalence and hence $C=\Cofib(f)\in \mathcal N$, and if $B,C\in \mathcal N$, the (co)fiber sequence $\Omega B\to \Omega C\to A$ implies that $A\in \mathcal N$. 
    It now follows from \cite[Theorem I.3.3]{NikolausScholze} (in the large universe) that \[\mathcal D(\mathcal A)\cong \operatorname{N}_{\operatorname{dg}}(\Ch(\mathcal A))[\text{(homology equivalences)}^{-1}]\] is a stable category and that the functor $\operatorname{N}_{\operatorname{dg}}(\Ch(\mathcal A))\to \mathcal D(\mathcal A)$ is exact. 
    In particular, cofibers can be computed using the mapping cone construction. 
\end{proof}

\begin{cor}[{\cite[section IV.4, Proposition 3]{GelfandManin}}]
    Suppose $\mathcal A$ is a (large) abelian category. 
    Then \[\mathcal D(\mathcal A)_{\geq 0}\coloneqq \{ X\in\mathcal D(\mathcal A)\, |\, H_i(X)=0\text{ for } i < 0\}, \] 
    \[\mathcal D(\mathcal A)_{\leq 0}\coloneqq\{ X\in\mathcal D(A)\, |\, H_i(X)=0 \text{ for }i >0\}\] defines a $t$-structure on $\mathcal D(A_{})$. 
\end{cor}

\begin{cor}\label{exactfunctorsinduceexactfunctorsonderivedcats}
    Suppose that $\phi\colon \mathcal A_1\to\mathcal A_2$ is an exact (i.e.\ finite limits and finite colimits preserving) functor between (large) abelian categories. 
    Then $\phi_*\colon \Ch(\mathcal A_1)\to\Ch(\mathcal A_2)$ descends to a functor $\mathcal D(\mathcal A_1)\to \mathcal D(\mathcal A_2)$ which is $t$-exact. 
\end{cor}
\begin{proof}
As $\mathcal A_1\to\mathcal A_2$ is exact, $H_i^{\mathcal A_2}\circ \phi_*=\phi_*\circ H_i^{\mathcal A_1}$ for all $i\in\mathbb Z$. 
In particular, $\phi_*$ descends to a functor 
$\mathcal D(\phi)\colon \mathcal D(\mathcal A_1)\to \mathcal D(\mathcal A_2)$ with \[H_i^{\mathcal D(\mathcal A_2)}\circ \mathcal D(\phi)=\mathcal D(\phi)\circ H_i^{\mathcal D(\mathcal A_1)}\] for all $i\in\mathbb Z$. 
Since finite coproducts in abelian categories are exact, for a finite set $I$, \[\oplus_{I}\colon \Fun(I, \Ch(\mathcal A_k))\to \Ch(\mathcal A_k)\to\mathcal D(\mathcal A_k)\] descends to a functor \[\oplus_{I}^{\mathcal D(\mathcal A_k)}\colon \Fun(I, \mathcal D(\mathcal A_k))\to \mathcal D(\mathcal A_k)\] for $k=1,2$. The unit and counit for the adjunction 
\[ \oplus_{I}\colon\Fun(I, \Ch(\mathcal A_k))\rightleftarrows \Ch(\mathcal A_k)\colon \operatorname{diag}\] induce natural transformations \[\oplus_{I}^{\mathcal D(\mathcal A_{k})}\operatorname{diag}\to \id, \, \id\to \operatorname{diag}\oplus_{I}^{\mathcal D(\mathcal A_k)}\] which exhibit $\oplus_{I}^{\mathcal D(\mathcal A_k)}$ as left adjoint to the diagonal functor \[\operatorname{diag}\colon \mathcal D(\mathcal A_k)\to \Fun(I, \mathcal D(\mathcal A_k)), \] i.e.\ finite coproducts can be computed on the level of chain complexes in $\mathcal D(\mathcal A_k)$. This implies that $\mathcal D(\phi)$ preserves finite coproducts. Since $\phi_*\colon \Ch(\mathcal A_1)\to\Ch(\mathcal A_2)$ preserves mapping cones, $\mathcal D(\phi)$ preserves cofibers by \cref{derivedcategoryisstable} and hence finite colimits. Since $\mathcal D(\mathcal A_{1})$ and $\mathcal D(\mathcal A_2)$ are stable, this shows that $f$ is exact. 
\end{proof}
\begin{notation}\label{definitionchainabeliancategories}
    A chain of abelian categories is a functor $\mt{A}_*\colon\Lambda\to \vlCat$ from a possibly large filtered set $\Lambda$ such that for all $\lambda\in\Lambda$, $\mt{A}_{\lambda}$ is an abelian category, and $\mt{A}_{\kappa}\to \mt{A}_{\lambda}$ is a fully faithful, left-exact left adjoint. 
\end{notation}
\begin{lemma}
If $\mt{A}\colon \Lambda\to \vlCat$ is a chain of abelian categories, its colimit $\colim{\Lambda}\mt{A}_{\lambda}$ in $\vlCat$ is an abelian category. 
\end{lemma}
\begin{proof}
    This holds as by \cref{filteredcolimitsofcategoriesappendix}, finite limits and finite colimits in $\mt{A}_{\infty}$ can always be computed in some stage $\mt{A}_{\lambda}$.
\end{proof}
\begin{lemma}\label{Filteredcolimitderivedcategories}
    Suppose $\mt{A}_*\colon \Lambda\to \vlCat$ is a chain of abelian categories whose colimit $\mathcal A_{\infty}$ has countable coproducts. 
    Then \[ \colim{\lambda}\mathcal D(\mathcal A_{\lambda})\cong \mathcal D(\mathcal A_{\infty})\in\vlCat.\] 
\end{lemma}

\begin{proof}As for $\lambda\geq \kappa$, $\Ch(\mathcal A_{\kappa})\to \Ch(\mathcal A_{\lambda})$ is fully faithful, 
    \[c\colon \colim{\kappa}\Ch(\mathcal A_{\kappa})\to \Ch(\mathcal A_{\infty})\] is fully faithful by \cref{filteredcolimitsofcategoriesappendix}.
    For a countable family $(C_i)_{i\in\mathbb Z}\subseteq \mathcal A_{\infty}$ and $j\in\mathbb Z$ there exists $\kappa\in \Lambda$ with $\oplus_{i\in\mathbb Z}C_i\in \mathcal A_{\kappa}\subseteq \mathcal A_{\infty}$. 
    For $j\in I$ let \[\pi_j\colon \oplus_{i\in \mathbb Z}C_i\to \oplus_{\substack{i\in I\\ i\neq j}}C_i\hookrightarrow \oplus_{i\in \mathbb Z} C_i.\] Since $\mathcal A_{\kappa}\subseteq \mathcal A_{\infty}$ is closed under colimits (\cref{filteredcolimitsofcategoriesappendix}), 
    $C_j=\Coker(\pi_j)\in \mathcal A_{\kappa}$ for all $j\in\mathbb Z$. 
    This implies that $c$ is an equivalence. 
    Denote by $W_{\infty/\kappa}\subseteq \Ch(\mathcal A_{\infty/\kappa})$ the wide subcategories on homology equivalences. As for all $\kappa\in\Lambda$,
    $\mathcal A_{\kappa}\to\mathcal A_{\infty}$ preserves finite limits and colimits and is conservative, \[W_{\infty}\cap\Ch(\mathcal A_{\kappa})=W_{\kappa}.\]
    In particular, $\colim{\kappa}\, W_{\kappa}\cong W_{\infty}$ as wide subcategories of $\Ch(\mathcal A)$.
    As left adjoint to the forget functor, $|-|\colon \vlCat\to \widehat{\an}$ commutes with large colimits. 
    This shows that $\mathcal D(\mathcal A)=\Ch(A)[W^{-1}]$ is the pushout of 
    \begin{center}
    \begin{tikzcd}
        \colim{\kappa}\, W_{\kappa}\arrow[r] \arrow[d]& \colim{\kappa}\, \mid W_{\kappa} \mid \\
        \colim{\kappa}\, \Ch(\mathcal A_{\kappa})\end{tikzcd}
    \end{center} in the very large category $\vlCat$ of large categories, and hence 
    \[\colim{\kappa}\, \mathcal D(\mathcal A_{\kappa})\cong \colim{\kappa}\, (\Ch(\mathcal A_{\kappa})\cup_{W_{\kappa}}|W_{\kappa}|)\cong \Ch(\mathcal A)\cup_{W}|W|=\mathcal D(\mathcal A).\]
\end{proof}
\begin{lemma}\label{derivedcategorymodules}
Suppose that $\mt{B}\colon\Lambda\to \Pr^L$ is a chain of hypercomplete, 1-localic topoi with colimit $\mt{B}_{\infty}\in \vlCat$ such that $\Ab(\tau_{\leq 0}\mt{B}_{\infty})$ has countable coproducts. 
For $\algebra{R}\in \Alg(\Ab(\tau_{\leq 0}\mt{B}_{\infty}))$, the \textit{inclusion}\[ \LMod{\algebra{R}}{\Ab(\tau_{\leq 0}\mt{B}_{\infty})}\cong \BMod{\algebra{R}}{\infty}^{\heart}\hookrightarrow \BMod{\algebra{R}}{\infty}\] extends to a $t$-exact equivalence
\begin{align*}\mathcal D(\LMod{\algebra{R}}{\Ab(\tau_{\leq 0}\mt{B}_{\infty})})\cong \BMod{\algebra{R}}{\infty}, \end{align*} where we view $\algebra{R}$ as an algebra in $\stab{\mt{B}_{\infty}}$ via the Eilenberg-Mac Lane functor (\cref{Eilenbergmaclanealgebra}).
\end{lemma}

\begin{proof}
Choose $\mu\in\Lambda$ with $\algebra{R}\in \Alg(\Ab(\mt{B_{\mu}}))\subseteq \Alg(\Ab({\mt{B_{\infty}}}))$. Since $\Lambda$ is filtered, $\Lambda_{\mu\backslash}\to \Lambda$ is cofinal, so we can assume that $\Lambda=\Lambda_{\mu\backslash}$ and for all $\lambda\in\Lambda$, $\algebra{R}\in\Alg(\Ab({\mt{B_{\lambda}}}))$. 
Fix $\lambda\in\Lambda$. 
By \cite[Proposition 2.1.1.1(c)]{SAG}, the $t$-structure on $\LMod{\algebra{R}}{\stab{\mt{B}_{\lambda}}}$ is right-complete (the proof there does not use commutativity of the ring $\algebra{R}$, hence also applies in our situation). 
By \cref{recognitionprinciplebigtopoi}, \[\Omega^{\infty}|_{\stab{\topo{B_{\lambda}}}_{\geq 0}}\colon \stab{\topo{B_{\lambda}}}_{\geq 0}\cong \CGrp(\topo{B_{\lambda}})\to \topo{B_{\lambda}}\] is conservative. As $\topo{B_{\lambda}}$ is hypercomplete, this implies that the $t$-structure on $\stab{\mt{B}_{\lambda}}$ is left-separated. Since the forget functor $\LMod{\algebra{R}}{\stab{\mt{B_{\lambda}}}}\to \stab{\mt{B_{\lambda}}}$ is conservative and $t$-exact, it follows that the $t$-structure on $\LMod{\algebra{R}}{\stab{\mt{B_{\lambda}}}}$ is left-separated.
Since $\LMod{\algebra{R}}{\Ab(\tau_{\leq 0}\mt{B}_{\lambda})}$ is Grothendieck abelian, the $t$-structure on $\mathcal D(\LMod{\algebra{R}}{\Ab(\tau_{\leq 0}\mt{B}_{\lambda})})$ is right-complete by \cite[Proposition 1.3.5.21]{higheralgebra}. 
It now follows from \cite[Proposition C.3.1.1, C.3.2.1, Theorem C.5.4.9]{SAG} that for all $\kappa, \lambda\in\Lambda$, restriction to hearts defines an equivalence \begin{align}\label{commutativitycanbecheckedonhearts}\Fun^{L, \operatorname{t-ex}}&(\mathcal D(\LMod{\algebra{R}}{\Ab(\tau_{\leq 0}\mt{B}_{\lambda})}), \LMod{\algebra{R}}{\stab{{B_{\kappa}}}})\\ &\cong \Fun^{L}(\LMod{\algebra{R}}{\Ab(\tau_{\leq 0}\mt{B}_{\lambda})}, \LMod{\algebra{R}}{\stab{\mt{B_{\kappa}}}}^{\heart}), \end{align} see also \cite[Proposition A.2]{cartiermodulesmattisweiss}. 
In particular, for $\kappa\to\lambda\in\Lambda$ there exists an essentially unique $t$-exact left adjoint 
\[c_{\kappa, \lambda}\colon \mathcal D(\LMod{\algebra{R}}{\Ab(\tau_{\leq 0}\mt{B}_{\kappa})})\to \LMod{\algebra{R}}{\stab{\mt{B_{\lambda}}}}\] extending the functor $\LMod{\algebra{R}}{\Ab(\tau_{\leq 0}\mt{B}_{\kappa})}\to \LMod{\algebra{R}}{\stab{\mt{B}_{\lambda}}}^{\heart}$ induced by \[\Ab(\tau_{\leq 0}\mt{B}_{\kappa})\cong \stab{\mt{B}_{\kappa}}^{\heart}\subseteq \stab{\mt{B}_{\lambda}}^{\heart}.\] 
For $\kappa\to\lambda\in\Lambda$, \ref{commutativitycanbecheckedonhearts} yields a commutative diagram 
\begin{center}
\begin{tikzcd}
    \mathcal D(\LMod{\algebra{R}}{\Ab(\tau_{\leq 0}\mt{B}_{\kappa})})\arrow[d]\arrow[r,"c_{{\kappa, \kappa}}"] & \LMod{\algebra{R}}{{\stab{\mt{B_{\kappa}}}}}\arrow[d] \\ 
    \mathcal D(\LMod{\algebra{R}}{\Ab(\tau_{\leq 0}\mt{B}_{\lambda})})\arrow[r, "c_{{\lambda, \lambda}}"] & \LMod{\algebra{R}}{\stab{\mt{B_{\lambda}}}}
\end{tikzcd}
\end{center} where the left vertical functor is induced by the fully faithful, exact left adjoint \[ \LMod{\algebra{R}}{\Ab(\tau_{\leq 0}\mt{B}_{\kappa})}\subseteq \LMod{\algebra{R}}{\Ab(\tau_{\leq 0}\mt{B}_{\lambda})}\] (cf.\ \cref{exactfunctorsinduceexactfunctorsonderivedcats}) and the right vertical morphism is induced by the fully faithful, symmetric monoidal left adjoint $\stab{\mt{B}_{\kappa}}\subseteq \stab{\mt{B}_{\lambda}}$. 
The proof of \cite[Proposition 6.4.5.7]{highertopostheory} shows that every 1-localic topos is a left-exact localisation of a presheaf category $\mathcal P(\mathcal C)$ for a 1-category $\mathcal C$. This implies that for every $\lambda\in\Lambda$ and $b\in\mt{B}_{\lambda}$, there exists $c\in \tau_{\leq 0}\mt{B}_{\lambda}$ with an effective epimorphism $c\to b$.\footnote{By \cite[Proposition 20.4.5.1, Remark 20.4.5.2]{SAG}, \cite[Proposition 6.4.5.9]{highertopostheory} this condition is satisfied for hypercomplete topoi if and only if they are 1-localic.}
Hence by \cite[Theorem 2.1.2.2]{SAG}, the maps $c_{\lambda, \lambda}$ are equivalences for all $\lambda\in\Lambda$ (the proof there also works for associative algebras), 
which shows that the induced functor 
\begin{align*}\colim{\lambda\in\Lambda}\mathcal D(\LMod{\algebra{R}}{\Ab(\tau_{\leq 0}\mt{B}_{\lambda})})\to\colim{\lambda\in\Lambda}\LMod{\algebra{R}}{\stab{\mt{B}_{\lambda}}}\end{align*} is an equivalence.
\cref{Filteredcolimitderivedcategories} and \cref{filteredcolimitsmodules} now imply that 
\[\mathcal D(\LMod{\algebra{R}}{\Ab(\tau_{\leq 0}\mt{B}_{\infty})})\cong \LMod{\algebra{R}}{\stab{\mt{B}_{\infty}}}.\] 
As $c_{\lambda, \lambda}$ is $t$-exact for all $\lambda\in\Lambda$, this equivalence is $t$-exact.  
\end{proof}
\begin{rem}
    For a topos $\topo{X}$ and an algebra $R\in\Alg(\Ab(\tau_{\leq 0}\topo{X}))$,  \cite[Proposition C.3.1.1, C.3.2.1, Theorem C.5.4.9]{SAG}/\cite[Proposition A.2]{cartiermodulesmattisweiss} yield a $t$-exact functor 
\[\mathcal D(\LMod{\algebra{R}}{\Ab(\tau_{\leq 0}\topo{X})})\to \LMod{\algebra{R}}{\stab{\widehat{\topo{X}}}}\subseteq \LMod{\algebra{R}}{\stab{\topo{X}}}, \] where $\widehat{\topo{X}}\subseteq \topo{X}$ denotes the full subcategory on hypercomplete objects. But the author is unaware of a construction of a $t$-exact functor \[\mathcal D(\LMod{\algebra{R}}{\Ab(\tau_{\leq 0}\topo{X})})\to \LMod{\algebra{R}}{\stab{\topo{X}}}\] for general big which extends the inclusion of the heart.  

If $\topo{X}_*\colon\Lambda\to \Pr^L$ is a chain of $\infty$-topoi, the transition functors $i_{\kappa}^{\lambda}\colon \topo{X}_{\lambda}\to\topo{X}_{\kappa}$ preserve $\infty$-connective morphisms. 
We therefore obtain a functor $\widehat{\topo{X}}_{*}\colon\Lambda\to \Pr^L$, where for all $\kappa\to\lambda\in\Lambda$, $\widehat{\topo{X}}_{\kappa}\to\widehat{\topo{X}_{\lambda}}$ is the composition of hypersheafification (in $\topo{X}_{\lambda})$ with $i^{\lambda}_{\kappa}$ and in particular a left-exact left adjoint. 
The argument from the above proof shows that there exists a $t$-exact functor
\[ \colim{\lambda\in\Lambda}\mathcal D(\LMod{\algebra{R}}{\Ab(\tau_{\leq 0}\topo{X}_{\lambda})})\to \colim{\lambda\in\Lambda}\LMod{\algebra{R}}{\stab{\widehat{\topo{X}}_{\lambda}}}.\] 
If $\Ab(\tau_{\leq 0}\topo{X})$ admits countable coproducts, then the left-hand side is equivalent to $\mathcal D(\Ab(\tau_{\leq 0}\topo{X}))$, but since $\widehat{\topo{X}}_*$ is not a subfunctor of $\topo{X}$, there is no obvious way to compare the right-hand side to $\LMod{\algebra{R}}{\stab{\topo{X}}}$. 

This is not an issue for categories of (hyper)accessible sheaves on hyperaccessible explicit covering sites: 
If $(\mathcal C,S)$ is a hyperaccessible explicit covering site, denote by $\Lambda_{\mathcal C}$ the poset of regular cardinals $\kappa$ such that $(\mathcal C,S)$ is $\kappa$-hyperaccessible. 
Then \[\hypershv_S(\mathcal C_{*})\subseteq \Shv_{S}(\mathcal C_*)\colon\Lambda\to \Pr^L\] is a subfunctor and $\Ab(\tau_{\leq 0}\shvacc_S(\mathcal C))=\shvacc_S(\mathcal C, \Ab)$ has small colimits. 
We therefore obtain $t$-exact functors 
\[ \mathcal D(\LMod{\algebra{R}}{\shvacc_S(\mathcal C, \Ab)})\to \LMod{\algebra{R}}{\hypershvacc_{S}(\mathcal C, \Sp)}\subseteq \LMod{\algebra{R}}{\shvacc_S(\mathcal C, \Sp)}.\] 
By the above lemma, the left functor is an equivalence if $\mathcal C$ is a 1-category. 
\end{rem}
\subsubsection{Modules over constant rings}

If $\topo{X}^{\otimes}$ is a big topos, the (symmetric) monoidal enhancement $c_{\Sp}^{\otimes}$ of the stabilization of the constant sheaf functor (\cref{constantsheafsymmetricmonoidal}) induces a functor \[c_{\Sp}\colon \aCAlg(\Sp)\to \aCAlg(\stab{\topo{X}}). \] 
We now record some basic observations about module categories over constant rings. 
\cref{modulecategoriestensoredup,symmetricmonoidalstructuretensoredupmodules} imply the following: 
\begin{lemma}\label{modulesoverconstantrings}Suppose $\mt{B_{\infty}}$ is a big topos and $\mt{B}_*\colon\Lambda\to \Pr^L$ is an exhaustion of $\mt{B}_{\infty}$ by topoi. 
    \begin{romanenum}
\item For a ring spectrum $A\in \Alg(\Sp)$, 
\[\colim{\Lambda} (\mt{B_{*}}\otimes_{\Pr^L}\LMod{A}{\Sp})\cong \BMod{c_{\Sp}A}{\infty}, \] where the colimit is computed in $\vlCat$.   
\item For a commutative ring spectrum $A\in\CAlg(\Sp)$, this  enhances to an an equivalence of symmetric monoidal categories
\[\colim{\Lambda} (\mt{B_{*}}^{\times}\otimes_{\CAlg(\Pr^L)}\LMod{A}{\Sp})\cong \BMod{c_{\Sp}A}{\infty}^{\otimes}, \] where the colimit is computed in $\CAlg(\vlCat)$, and the right-hand side is endowed with the symmetric monoidal structure from \cref{symmetricmonoidalstructureonspectrumobjects,symmetricmonoidalstructurebigtoposmodules}.    
    \end{romanenum}
\end{lemma}

\begin{proof}We first prove the first statement. 
    Choose $\lambda\in\Lambda$. As the constant sheaf functor factors over $\stab{\mt{B}_{\lambda}}\subseteq \stab{\mt{B}_{\infty}}$, by \cref{modulecategoriestensoredup}, 
    \[\BMod{c_{\Sp}A}{\infty}\cong \colim{\mu\in \Lambda_{\lambda\backslash}}\stab{\mt{B}_{\mu}}\otimes_{\LMod{\stab{\mt{B}_{\lambda}}}{\Pr^L}}\LMod{c_{\Sp}A}{\stab{\mt{B}_{\lambda}}}.\] 
    By  \cite[Theorem 4.8.4.6]{higheralgebra}, $\LMod{c_{\Sp}A}{\stab{\mt{B}_{\lambda}}}\cong \LMod{A}{\Sp}\otimes_{\LMod{\Sp}{\Pr^L}}\stab{\mt{B}_{\lambda}}$, and by \cref{identifystabilizationwithtensorproductwithspectra}, $\Sp(\mt{B}_*)\cong \Sp\otimes_{\Pr^L}\mt{B}_*$. It now follows that 
    \begin{align*}\colim{\mu\in \Lambda_{\lambda\backslash}}\stab{\mt{B}_{\mu}}\otimes_{\LMod{\stab{\mt{B}_{\lambda}}}{\Pr^L}}\LMod{c_{\Sp}A}{\stab{\mt{B}_{\lambda}}}& \cong  \colim{\mu\in \Lambda_{\lambda\backslash}}\stab{\mt{B}_{\mu}}\otimes_{\LMod{\Sp}{\Pr^L}}\LMod{A}{\Sp}\\ &\cong \colim{\mu\in \Lambda_{\lambda\backslash}}\mt{B}_{\mu}\otimes_{\Pr^L}\LMod{A}{\Sp}.\end{align*} 

    Suppose now that $A\in\CAlg(\Sp)$. 
    By \cref{symmetricmonoidalstructuretensoredupmodules}, 
        \[\BMod{c_{\Sp}A}{\infty}^{\otimes}\cong \colim{\mu\in\Lambda_{\lambda\backslash}}(\stab{\mt{B_{*}}}\otimes_{\CAlg(\LMod{\stab{\mt{B}_{\lambda}}}{\Pr^L})}\LMod{c_{\Sp}A}{\stab{\mt{B}_{\lambda}}}),\] 
        and by \cite[Theorem 4.8.4.6]{higheralgebra}, $\LMod{c_{\Sp}A}{\stab{\mt{B}_{\lambda}}}\cong \LMod{A}{\Sp}\otimes_{\CAlg(\LMod{\Sp}{\Pr^L})}\stab{\mt{B}_{\lambda}}$, whence \begin{align*}\LMod{A}{\stab{\mt{B}_{\lambda}}}&\cong \LMod{A}{\Sp}\otimes_{\CAlg(\LMod{\Sp}{\Pr^L})}\stab{\mt{B}_{\lambda}}\otimes_{\CAlg(\LMod{\stab{\mt{B}_{\lambda}}}{\Pr^L})}\stab{\mt{B}_{\mu}}\\ &\cong \LMod{A}{\Sp}\otimes_{\CAlg(\LMod{\Sp}{\Pr^L})}\stab{\mt{B}_{\mu}}\\ &\cong \LMod{A}{\Sp}\otimes_{\CAlg(\Pr^L)}\mt{B}_{\mu}\end{align*} naturally in $\mu$, where we used that $\stab{\mt{B}_{*}}^{\otimes}\cong \Sp^{\otimes}\otimes_{\CAlg(\Pr^L)}\mt{B}_{*}^{\otimes}$ (\cite[Proposition 5.4]{GepnerGrothNikolaus}).
    \end{proof}   
    \begin{cor}\label{modulesoverconstantringssymmetricmonoidalaccessibleshaves}Suppose $(\mathcal C,S)$ is a (hyper)accessible explicit covering site and denote by $c_{\Sp}\colon \CAlg(\Sp)\to \CAlg(\mywidehatshvacc{S}(\mathcal C,\Sp))$ the functor induced by the stabilization of the constant sheaf functor. 
        \begin{romanenum}
\item For $\algebra{A}\in\Alg(\Sp)$, 
        \[\mywidehatshvacc{S}(\mathcal C,\LMod{\algebra{A}}{\Sp})\cong \LMod{c_{\Sp}\algebra{A}}{\mywidehatshvacc{S}(\mathcal C,\Sp)}.\]
    \item For $\algebra{A}\in\CAlg(\Sp)$, this enhances to a symmetric monoidal equivalence 
        \[\mywidehatshvacc{S}(\mathcal C,\LMod{\algebra{A}}{\Sp}^{\otimes})\cong \LMod{c_{\Sp}\algebra{A}}{\mywidehatshvacc{S}(\mathcal C,\Sp)},\] where the left-hand side is endowed with the symmetric monoidal structure induced from $\LMod{A}{\Sp}^{\otimes}$ via \cref{constructionmonoidalstructureaccessiblesheaves}, and the right-hand side with the symmetric monoidal structure from \cref{symmetricmonoidalstructurebigtoposmodules}. 
        \end{romanenum}
    \end{cor}
    \begin{proof}The first statement holds by \cref{modulesoverconstantrings,accessiblesheavesisbigtopos,identifysheaveswithtensorproductinprl}, the second statement follows from \cref{modulesoverconstantrings,constructionmonoidalstructureaccessiblesheaves}. 
    \end{proof}
Recall from \cite[Theorem 2.1.2.2]{SAG} or \cref{derivedcategorymodules} that if $A=\pi_0A$ is a discrete ring spectrum, then $\LMod{A}{\Sp}\cong \mathcal D(A)$ is the $\infty$-derived category of $A$-modules. 

    If $\topo{X}^{\otimes}$ is a big topos, for $A\in \Sp$, the symmetric monoidal enhancement $c_{\Sp}^{\otimes}$ of the stabilization of the constant sheaf functor (\cref{constantsheafsymmetricmonoidal}) induces a left adjoint \[ c_{\Sp}\colon \LMod{A}{\Sp}\to \LMod{c_{\Sp}A}{\stab{\topo{X}}} \text{ (\cref{localisationinduceslocalisationonmodulecategories})}.\]
    We will use the following elementary observation to construct simplicial resolutions (\cref{simplicialresolutionisresolution}) and to compare internal (group) cohomology with (group) cohomology:
\begin{lemma}\label{freemoduloediscrete}
    Suppose $\topo{X}$ is a big topos and denote by $c_{\Sp}\colon\Alg(\Sp)\to\Alg(\stab{\topo{X}})$ the functor induced by the constant sheaf functor. 
    If $A\in \Alg(\Sp^{\heart})\cong \Alg(\Ab)$ is a discrete ring, for $X\in\tau_{\leq 0}\topo{X}$, 
    $c_{\Sp}A[\Sigma^{\infty}_{+}X]\in \LMod{c_{\Sp}A}{\stab{\topo{X}}}^{\heart}$. 
\end{lemma}
\begin{proof}By definition of the $t$-structure on $\LMod{c_{\Sp}A}{\stab{\topo{X}}}$, we have to show that \[c_{\Sp}A\otimes \Sigma^{\infty}_{+}X\in \stab{\topo{X}}^{\heart}.\] 
Choose an exhaustion $\topo{X}_*\colon\Lambda\to \Pr^L$ of $\topo{X}$ by topoi.  
As for all $\lambda\in\Lambda$, $\stab{\topo{X}_{\lambda}}\subseteq\stab{\topo{X}}$ is a symmetric monoidal subcategory and $\stab{\topo{X_{\lambda}}}\subseteq\stab{\topo{X}}$ is $t$-exact, it suffices to prove the statement for $\topo{X_{\lambda}}$, so we can assume that $\topo{X}$ is a topos, i.e.\ there exists a small category $\mathcal C$ with a left-exact, accessible localization $L\colon\mathcal P(\mathcal C)\to \mathcal X$. 
Since $L$ is cartesian, $c^{\times}$ factors as $\an^{\times}\xrightarrow{\operatorname{const}}\mathcal P(\mathcal C)^{\times}\xrightarrow{L} \mathcal C^{\times}$, where $\operatorname{const}$ sends an anima to the associated constant presheaf. 
This implies that \[ c_{\Sp}A[\Sigma^{\infty}_{+}X]=L_{\Sp}(c \mapsto A\otimes_{\Sp} \Sigma^{\infty}_{+}X(c)).\] 
As for $c\in\mathcal C$, $X(c)=\Map_{\mathcal X}(L y(c),X)\in \tau_{\leq 0}\an$ (where $y$ denotes the Yoneda embedding), \[\mathcal C\ni c\mapsto A\otimes_{\Sp}\Sigma^{\infty}_{+}X(c)\in \stab{\mathcal P(\mathcal C)}^{\heart}.\] 
$t$-exactness of $L_{\Sp}$ (\cref{texactnessgeometricmorphism}) now implies that $ c_{\Sp}A[\Sigma^{\infty}_{+}X]\in\stab{\topo{X}}^{\heart}$.  
\end{proof} 

We now explain that for a big topos $\topo{X}$, the category of left $A$-modules $\LMod{A}{\stab{\topo{X}}}$ over a a commutative algebra $A\in\CAlg(\stab{\topo{X}})$ is naturally enriched in $\LMod{\Gamma_{\Sp}A}{\Sp}$.  

\begin{notation}Suppose $\topo{X}$ is a big topos. 
The symmetric monoidal enhancement of the constant sheaf/global sections-adjunction (\cref{constantsheafsymmetricmonoidal,adjunctionfiberwise}) induces an adjunction \[c_{\Sp}\colon \CAlg(\Sp)\rightleftarrows\CAlg(\stab{\topo{X}})\colon \Gamma_{\Sp} \text{ (\cref{adjuncttioninducesadjunctiononalgebraobjects})}.\]
For $\algebra{A}\in\CAlg(\stab{\topo{X}})$, the counit $c_{\Sp}\Gamma_{\Sp}A\to A\in\CAlg(\stab{\topo{X}})$ determines a symmetric monoidal enhancement of the free $A$-module functor 
\[A[-]\colon\LMod{c_{\Sp}\Gamma_{\Sp}A}{\stab{\topo{X}}}\to \LMod{A}{\LMod{c_{\Sp}\Gamma_{\Sp}A}{\stab{\topo{X}}}}\] (\cref{symmetricmonoidalenhancementfreefunctorpresentablymonoidalcat}), and by \cref{modulesinmodulecategoriessymmetricmonoidalbigpresentable}, 
\[ \LMod{A}{\LMod{c_{\Sp}\Gamma_{\Sp}}{\stab{\topo{X}}}}\cong \LMod{A}{\stab{\topo{X}}}\] as symmetric monoidal categories. 

By \cref{symmetricmonoidalfunctorsfrompresentableenhancetosymmetricmonoidalfunctorsonmodulecategories}, the constant sheaf functor induces a symmetric monoidal functor 
\[ c_{\Sp}^{A}\colon \LMod{\Gamma_{\Sp}A}{\Sp}\to \LMod{c_{\Sp}\Gamma_{\Sp}A}{\stab{\topo{X}}}.\] 

Denote by \[c_A\colon \LMod{\Gamma_{\Sp}A}{\Sp}\to \LMod{A}{\stab{\topo{X}}}\] the composition
\[ \LMod{\Gamma_{\Sp}A}{\Sp}\xrightarrow{c^A_{\Sp}} \LMod{c_{\Sp}\Gamma_{\Sp}A}{\stab{\topo{X}}}\xrightarrow{A[-]} \LMod{A}{\LMod{c_{\Sp}\Gamma_{\Sp}A}{\stab{\topo{X}}}}\cong \LMod{A}{\stab{\topo{X}}}.\] 
\end{notation}
\begin{lemma}\label{modulecategoriesareenriched}Suppose that $\topo{X}$ is a big topos. For $A\in\CAlg(\stab{\topo{X}})$, the functor \[c_A\colon \LMod{\Gamma_{\Sp}A}{\Sp}\to \LMod{A}{\stab{\topo{X}}}\] exhibits $\LMod{A}{\stab{\topo{X}}}$ as enriched (\cite[Definition 4.2.1.28]{higheralgebra}) over $\LMod{\Gamma_{\Sp}A}{\Sp}^{\otimes}$. 
\end{lemma}
\begin{proof}
Choose an exhaustion \[\topo{X}_*\colon\Lambda\to \CAlg(\Pr^L)\] of ${\topo{X}}^{\otimes}$ by presentably symmetric monoidal categories with $A\in\Alg(\stab{\topo{X}_{\lambda}})$ for all $\lambda\in\Lambda$.  As $c_{\Sp}$ factors over $\mt{X}_{\lambda}$ for all $\lambda\in\Lambda$, $c_{\Sp}^{A,\otimes}$ factors as \[  \LMod{\Gamma_{\Sp}A}{\Sp}^{\otimes}\to  \LMod{c_{\Sp}\Gamma_{\Sp}A}{\stab{\topo{X}_{\lambda}}}^{\otimes}\subseteq \LMod{c_{\Sp}\Gamma_{\Sp}A}{\stab{\topo{X}}}^{\otimes}\] for all $\lambda\in\Lambda$. 
The construction of the free $A$-module functor (\cref{symmetricmonoidalenhancementfreefunctorpresentablymonoidalcat}) and the identification of \cref{modulesinmodulecategoriessymmetricmonoidalbigpresentable} imply that $c_A$ factors as 
\[ \LMod{\Gamma_{\Sp}A}{\Sp}\xrightarrow{c^{\lambda}_A} \LMod{A}{\stab{\topo{X}_{\lambda}}}\subseteq \LMod{A}{\stab{\topo{X}}},\] where \[ c^{\lambda}_A\colon\LMod{\Gamma_{\Sp}A}{\Sp}\to \LMod{A}{\stab{\topo{X}_{\lambda}}}\] is the symmetric monoidal functor induced by the constant sheaf functor for $\topo{X}_{\lambda}$. 

For $M\in\LMod{A}{\stab{\topo{X}}}$ choose $\lambda\in\Lambda$ with $M\in\LMod{A}{\stab{\topo{X}_{\lambda}}}$. By the above, 
\[-\otimes M\colon \LMod{A}{\Sp}\to \LMod{c_{\Sp}A}{\stab{\topo{X}}}\] factors as \[\LMod{A}{\Sp}\xrightarrow{-\otimes_{\lambda}M} \LMod{A}{\stab{\topo{X}_{\lambda}}}\subseteq \LMod{A}{\stab{\topo{X}}}, \] where $-\otimes_{\lambda}M$ is induced by the left-tensoring of $\LMod{A}{\stab{\topo{X}_{\lambda}}}$ over $\LMod{\Gamma_{\Sp}A}{\Sp}$ induced by the constant sheaf functor for $\topo{X}_{\lambda}$. 
By \cite[Proposition 4.3.1.33]{higheralgebra}, $-\otimes_{\lambda}M$ is a left adjoint, and since $\stab{\topo{X}_{\lambda}}\subseteq \stab{\topo{X}}$ is a symmetric monoidal left adjoint, $\LMod{A}{\stab{\topo{X}_{\lambda}}}\subseteq \LMod{A}{\stab{\topo{X}}}$ is a left adjoint by \cref{localisationinduceslocalisationonmodulecategories}. 
\end{proof}

\begin{cor}\label{continuityenrichmentorbits}Suppose $\mt{B}_{\infty}$ is a big topos such that $\stab{\mt{B}_{\infty}}$ has all small limits and colimits. For $A\in\CAlg(\Sp)$, the $\LMod{\Gamma_{\Sp}A}{\Sp}$-enriched mapping functor (\cite[Remark 4.2.1.31]{higheralgebra})
\[ \map_{\LMod{A}{\stab{\topo{B}_{\infty}}}}(-,-)\colon \LMod{A}{\mt{B}_{\infty}}^{\operatorname{op}}\times \LMod{A}{\mt{B}_{\infty}}\to \LMod{\Gamma_{\Sp}A}{\Sp}\] preserves small limits in both variables.
\end{cor}
\begin{proof}As for all $b\in \stab{\topo{B}_{\infty}}$, $\map_{\LMod{A}{\stab{\topo{B}_{\infty}}}}(b,-)$ is a right adjoint, $\map_{\LMod{A}{\stab{\topo{B}_{\infty}}}}(-,-)$ preserves small limits in the second variable. 
By \cref{monoidalstructurespectrumobjectscocontinuous} and \cref{monoidalstructuremodulescocontinuous}, the symmetric monoidal structure on $\LMod{\Gamma_{\Sp}A}{\stab{\mt{B}_{\infty}}}$ is cocontinuous. 
As the left-tensoring \[\LMod{A}{\Sp}\times \LMod{A}{\stab{\mt{B}_{\infty}}}\to \LMod{c_{\Sp}A}{\stab{\mt{B}_{\infty}}}\] is induced by the symmetric monoidal left adjoint $\LMod{\Gamma_{\Sp}A}{\Sp}\to\LMod{A}{\stab{\mt{B}_{\infty}}}$, this implies that the left-tensoring \[\LMod{\Gamma_{\Sp}A}{\Sp}\times \LMod{A}{\stab{\mt{B}_{\infty}}}\to \LMod{A}{\stab{\mt{B}_{\infty}}}\] preserves small colimits in both variables. 
It now follows from \cref{leftadjointsstableundercolimits} and \cite[Proposition 5.2.6.2]{highertopostheory} that $\map_{\LMod{A}{\stab{\topo{B}_{\infty}}}}(-,-)$ preserves small limits in the first variable. 
\end{proof}

\subsection{Cohomology in big topoi}\label{section:cohomologyinabigtopos}
Our discussion from \cref{section:spectrumobjectsofbigpresentablecategories} leads a robust notion of cohomology in big topoi. In this section, we record some  elementary observations about it. In particular, we show that cohomology is invariant under geometric morphisms (\cref{geometricmorphismcohomology}), and admits a \v{C}ech-to-cohomology spectral sequence (\cref{Bousfieldkanspectralsequencehomology}).
\begin{definition}\label{definitioncohomologyinatopos}
    Suppose $\topo{X}$ is a big topos and denote by $\Sigma^{\infty}_{+}\colon \topo{X}\to \stab{\topo{X}}$ the stabilization functor (\cref{existencesuspension}). 
\begin{romanenum}
\item Denote by $\map_{\stab{\topo{X}}}(-,-)$ the spectrally enriched mapping functor of $\stab{\topo{X}}$ (\cref{bigpresentablespectrallyenriched},  \cite[Remark 4.2.1.31]{higheralgebra}). For $X\in \topo{X}$ and $A\in \stab{\topo{X}}$, the \emph{cohomology spectrum} of $A$ at $X$ is 
\[ \cH{\topo{X}}(X,A)\coloneqq \map_{\stab{\topo{X}}}(\Sigma^{\infty}_{+}X,A).\] 
For $i\in \mathbb Z$, the ith cohomology group of $X$ with coefficients in $A$ is \[ \cH{\topo{X}}^i(X,A)\coloneqq \pi_{-i}(\cH{\topo{X}}(X,A)).\]

\item If the symmetric monoidal structure on $\stab{\topo{X}}$ provided by \cref{symmetricmonoidalstructureonspectrumobjects} is closed, denote by $\imap_{\stab{\topo{X}}}(-,-)$ the internal Hom in $\stab{\topo{X}}$ and define \emph{internal cohomology} by \[\underline{\mathbb H}_{\topo{X}}(X,A)\coloneqq \imap_{\stab{\topo{X}}}(\Sigma^{\infty}_{+}X,A)\in \stab{\topo{X}}.\] For $i\in\mathbb Z$, \[\underline{\mathbb H}^i_{\topo{X}}(X,A)\coloneqq \pi_{-i}^{\mathcal X}\imap_{\stab{\topo{X}}}(\Sigma^{\infty}_{+}X,A)\in \stab{\topo{X}}^{\heart}\cong \Ab(\tau_{\leq 0}\topo{X}).\]\end{romanenum}
\end{definition}
\begin{rems}

\begin{enumerate}
        \item Recall that if $\topo{X}$ is a topos or the category of accessible (hyper)\-sheaves on a (hyper)ac\-ces\-si\-ble explicit covering site, then the symmetric monoidal structure on $\stab{\topo{X}}$ is closed by \cref{symmetricmonoidalstructureaccessiblesheavesofspectraclosed}, whence internal cohomology in $\topo{X}$ exists. 
  
\item Denote by $\Gamma_{\Ab}\colon\Ab(\tau_{\leq 0}\topo{X})\to \Ab$, $\Gamma_{\Sp}\colon \stab{\topo{X}}\to \Sp$ the functors induced by the global sections geometric morphism $\Gamma\colon \topo{X}\to \an$ (\cref{constantsheaffunctordefinition}).
By \cref{internalhomrecoversspectralenrichment}, \[\Gamma_{\Sp}\circ \underline{\mathbb H}_{\topo{X}}(X,A)\cong \cH{\topo{X}}(X,A)\] (provided that internal cohomology in $\topo{X}$ exists).  
If $\Gamma_{\Sp}$ is $t$-exact (equivalently, if $\Gamma_{\Ab}$ is exact), this implies that \[\Gamma_{\Ab}\circ \intgrpcoh{\topo{X}}^q(X,A)=\mathbb H^q_{\topo{X}}(X,A)\text{ for all }q\in\mathbb Z.\] 
This holds for example for $\topo{X}=\Cond{(\kappa)}(\an)$ (\cref{kappacondensedtopology}, \cref{definitioncondensed}), see \cref{globalsectionstexactcocontinuous}. 
\end{enumerate}
\end{rems}
\begin{ex}\label{cohomologyglobalsections}
Suppose $\topo{X}$ is a big topos and denote by \[ c\colon \an\rightleftarrows \topo{X}\colon \Gamma\] the global sections geometric morphism. 
\cref{stabilizationcommuteswithsuspension}, and \cref{spectralenrichmentcomesfromconstantsheaf} imply that for $A\in\Sp$, 
\[ \cH{\topo{X}}(cA,-)\cong \map_{\stab{\topo{X}}}(c_{\Sp}\Sigma^{\infty}_{+}A,-)\] is right adjoint to 
\[ \Sp\xrightarrow{c_{\Sp}\Sigma^{\infty}_{+}A\otimes c_{\Sp}(-)}\stab{\topo{X}}.\] 
The symmetric monoidal structure of $c_{\Sp}$ (\cref{constantsheafsymmetricmonoidal}) yields an equivalence of this functor with $c_{\Sp}(\Sigma^{\infty}_{+}A\otimes -)$, whence \[ \cH{\topo{X}}(cA,-)\cong \imap_{\Sp}(\Sigma^{\infty}_{+}A,-)\circ \Gamma_{\Sp}.\] 
\end{ex}

\begin{rems}\label{cohomologyisext1}

It is natural to consider the following enhancements of cohomology in a big topos $\topo{X}$:
\begin{enumerate}
\item\label{cohomologyisext22} For $\algebra{R}\in \Alg(\stab{\topo{X}})$, denote by $\map_{\LMod{\algebra{R}}{\stab{\mathcal X}}}(-,-)$ the $\Sp$-enrichment of the stable, big presentable category $\LMod{R}{\stab{\topo{X}}}$ (\cref{bigpresentablespectrallyenriched}, \cite[Remark 4.2.1.31]{higheralgebra}). 

The free-forget adjunction \[\algebra{R}[-]\colon \stab{\topo{X}}\rightleftarrows \LMod{\algebra{R}}{\stab{\topo{X}}}\colon f\] is spectrally enriched by \cref{adjunctionsspectrallyenriched}, whence 
    \begin{align}\label{cohomologyisext2}\cH{\mathcal X}(-,f-)\cong \map_{\LMod{\algebra{R}}{\stab{\mathcal X}}}(\algebra{R}[\Sigma^{\infty}_{+}-],-).\end{align}
    In particular, if $\topo{X}$ is a big topos as in \cref{derivedcategorymodules}, for $\algebra{R}\in\Alg(\Ab(\tau_{\leq 0}\topo{X}))$
    \begin{align}\label{cohomologyisext5}\cH{\topo{X}}^*(-,fi-)|_{\tau_{\leq 0}\topo{X}^{\operatorname{op}}\times \mathcal D(\LMod{\algebra{R}}{\Ab(\tau_{\leq 0}\topo{X})})}\cong \Ext^*_{\LMod{\algebra{R}}{\Ab(\tau_{\leq 0}\topo{X})}}(\algebra{R}[-],-), \end{align} where \[i\colon\mathcal D(\LMod{\algebra{R}}{\Ab(\tau_{\leq 0}\topo{X})})\cong \LMod{\algebra{R}}{\stab{\topo{X}}}\] is the equivalence provided by \cref{derivedcategorymodules} and $\algebra{R}[-]\colon \tau_{\leq 0}\topo{X}\to \LMod{\algebra{R}}{\Ab(\tau_{\leq 0}\topo{X})}$ denotes the free $\algebra{R}$-module functor. 
\item Denote by $\Gamma_{\Sp}\colon \CAlg(\stab{\topo{X}})\to \CAlg(\Sp)$ the functor induced by the global sections functor. For $R\in\CAlg(\stab{\topo{X}})$, $\LMod{\algebra{R}}{\stab{\topo{X}}}$ is enriched in $\LMod{\Gamma_{\Sp}\algebra{R}}{\Sp}$ by \cref{modulecategoriesareenriched}. In particular, cohomology with coefficients in $\algebra{R}$-modules refines to a functor \[\cH{\topo{X}}(-,-)_R\colon \topo{X}^{\operatorname{op}}\times\LMod{\algebra{R}}{\stab{\topo{X}}}\to \LMod{\Gamma_{\Sp}\algebra{R}}{\Sp}.\]  
\item 
If $\stab{\topo{X}}$ is cartesian closed and has $\Delta^{\operatorname{op}}$-indexed colimits and all $\Delta$-indexed limits, then for all commutative algebras $\algebra{R}\in \CAlg(\stab{\topo{X}})$, internal cohomology with coefficients in $\algebra{R}$-modules refines to a functor \[\icH{\topo{X}}(-,-)_{\algebra{R}}\colon \topo{X}^{\operatorname{op}}\times\LMod{\algebra{R}}{\stab{\topo{X}}}\to \LMod{\algebra{R}}{\stab{\topo{X}}}\] by \cref{internalhommodules}.  
\end{enumerate}
\end{rems}
\begin{lemma}
    For $i\in\mathbb Z$, $\cH{\topo{X}}^i(-,-)\cong \pi_0\Map_{\topo{X}}(-, \Omega^{\infty}\Sigma^i-).$

    In particular, the counit $\tau_{\geq -i}\to \id$ induces an equivalence \[ \cH{\topo{X}}^i(-,-)\cong \cH{\topo{X}}^i(-, \tau_{\geq -i}-).\] 
\end{lemma}
\begin{proof}
    For $i\in\mathbb Z$, 
    \begin{align*} \cH{\topo{X}}^i(-,-)=\pi_{-i}\map_{\stab{\topo{X}}}(\Sigma^{\infty}_{+}-,-)=\pi_0\Sigma^i\map_{\stab{\topo{X}}}(\Sigma^{\infty}_{+}-,-)&=\pi_0\map_{\stab{\topo{X}}}(\Sigma^{\infty}_{+}-, \Sigma^i-)\\&=\pi_0\Map_{\stab{\topo{X}}}(-, \Omega^{\infty}\Sigma^i-).\end{align*} 
    This implies the second statement since $\Omega^{\infty}\Sigma^i=\Omega^{\infty}\tau_{\geq 0}\Sigma^i=\Omega^{\infty}\Sigma^i\tau_{\geq -i}.$
\end{proof}
\begin{definition}\label{definitiongeometricmorphism}
Suppose $\topo{X}, \topo{Y}$ are big topoi. 
A geometric morphism $L\colon \topo{X}\rightleftarrows \topo{Y}\colon R$ is an adjoint pair $L\dashv R$ such that the left adjoint $L$ preserves finite limits. 
\end{definition}
By \cref{geometricmorphismstabilization,stabilizationcommuteswithsuspension}, a geometric morphism stabilizes to an adjunction
    \[L_{\Sp}\colon\stab{\topo{X}}\rightleftarrows \stab{\topo{Y}}\colon R_{\Sp}\] with $L_{\Sp}\circ \Sigma^{\infty}_{+, \topo{X}}\cong \Sigma^{\infty}_{+, \topo{Y}}\circ L$. 
\begin{lemma}\label{geometricmorphismcohomology}
    Suppose $L\colon \topo{X}\rightleftarrows \topo{Y}\colon R$ is a geometric morphism between big topoi.
    The adjunction $L\dashv R$ induces an equivalence of functors
    \[\cH{\topo{Y}}(L-,-)\cong \cH{\topo{X}}(-,R_{\Sp}-)\colon \topo{X}^{\operatorname{op}}\times \stab{\topo{Y}}\to \Sp.\] 
\end{lemma}
\begin{proof}
    By \cref{adjunctionsspectrallyenriched}, the adjunction $L_{\Sp}\dashv R_{\Sp}$ is spectrally enriched, 
    whence \[\cH{\topo{X}}(X,R_{\Sp}A)=\map_{\stab{\topo{X}}}(\Sigma^{\infty}_{+}X,R_{\Sp}A)\cong \map_{\stab{\topo{X}}}(\underbrace{L_{\Sp}\Sigma^{\infty}_{+}X}_{\cong \Sigma^{\infty}_{+}LX},A)\cong\cH{\topo{X}}(LX,A).\qedhere\]
\end{proof}
\begin{rem}\label{existenceleftadjointmoregeneral}
    If $R\colon \topo{Y}\to \topo{X}$ is a right adjoint between presentable categories with left adjoint $L$, the stabilization $R_{\Sp}\colon \stab{\topo{X}}\to\stab{\topo{Y}}$ (\cref{definitionstabilization}) admits a left adjoint $L_{\Sp}$ by \cite[Corollary 1.4.4.4]{higheralgebra}. 
    This satisfies $L_{\Sp}\circ \Sigma^{\infty}_{+, \topo{X}}\cong \Sigma^{\infty}_{+, \topo{X}}\circ L$ since both functors are left adjoint to $\Omega^{\infty}\circ R_{\Sp}=R\circ\Omega^{\infty}$. In particular, \[\cH{\topo{X}}(L,-)\cong \cH{\topo{Y}}(-,R_{\Sp}-).\]
    However, for big topoi the functor $R_{\Sp}$ need not admit a left adjoint if its left adjoint $L$ is not left-exact. 
\end{rem}
\begin{cor}\label{cohomologycanbecomputedonfinitelelvel}
    Suppose $\mt{B}_*\colon \Lambda\to \Pr^L$ is a chain of topoi. 
    For all $\kappa\in\Lambda$, 
    \[ \cH{\mt{B}_{\infty}}(-,-)|_{\mt{B}_{\kappa}^{\operatorname{op}}\times  \stab{\mt{B_{\kappa}}}}\cong\cH{\mt{B}_{\kappa}}(-,-).\]     
\end{cor}
\begin{proof}
    Denote by $r^{\kappa}\colon\mt{B_{\infty}}\to \mt{B_{\kappa}}$ the right adjoint of the canonical functor $i_{\kappa}\colon \mt{B_{\kappa}}\hookrightarrow\mt{B_{\infty}}$. This exists by \cref{fullyfaithfulnessandpreservationoflimits}. 
    For $x\in \mt{B_{\kappa}}\subseteq\mt{B_{\infty}}$, 
    \[ \cH{\mt{B}_{\infty}}(i_{\kappa}x,-)\cong \cH{\mt{B}_{\kappa}}(x,r^{\kappa}_{\Sp}-)\] by \cref{geometricmorphismcohomology}.  
    As $\mt{B_{\kappa}}\xhookrightarrow{i_{\kappa}}\mt{B_{\infty}}\xrightarrow{r^{\kappa}}\mt{B_{\kappa}}$ is the identity (\cref{fullyfaithfulnessandpreservationoflimits}), this implies the statement. 
\end{proof}

For the convenience of the reader, we spell out a version of a \v{C}ech-to-cohomology spectral sequence in a big topos that arises as an application of \cite[section 1.2.4]{higheralgebra}. 
\begin{proposition}[{\cite[Proposition 1.2.4.5, Variant 1.2.4.9]{higheralgebra}}]\label{propositionbousfieldkanspectralsequence}
    Suppose $\mathcal C$ is a stable category with countable products and $(\mathcal C_{\geq 0}, \mathcal C_{\leq 0})$ is a $t$-structure.  
    \begin{romanenum}
    \item If $X_*\colon\Delta\to\mathcal C$ is a cosimplicial object in $\mathcal C$, 
    there exists a spectral sequence $BKSS(X_*)$
    with $E_1^{p,q}=\pi_{-q}(X_p)$ and $d_1$-differential  
    \[ d_1\colon E_1^{p,*}=\pi_{-*}(X_p)\to E_1^{p+1,*}=\pi_{-*}(X_{p+1})\]
    the alternating sum 
    \[ d_{1}^{p,*}\coloneqq \sum_{i=0}^{p}(-1)^i \pi_{-*}(s_i)\] of the degeneracies \[ s_i\colon X_p\to X_{p+1}.\]

    \item If there exists $N\in\mathbb N_0$ with $X_*\in \Fun(\Delta, \mathcal C_{\leq N})$, then the above spectral sequence converges to $\pi_{-(p+q)}(\clim{\Delta}X)$, filtered by \[\{F^p\pi_*(\clim{\Delta}X)\coloneqq \ker(\pi_*(\clim{\Delta}X)\to \pi_*(\clim{\Delta_{\leq p}}X))\}_{p\geq 0}.\] 
    This filtration is bounded, for $k\in\mathbb Z$: 
    \[0=F^{k+N-1}\pi_{-k}\clim{\Delta}X\subseteq \ldots \subseteq F^0\pi_{-k}\clim{\Delta}X\subseteq F^{-1}\pi_{-k}\clim{\Delta}X=\pi_{-k}\clim{\Delta}X.\] 
    \item This is natural in the cosimplicial object $X_*$:
    A map of cosimplicial objects \[c_*\colon X_*\to Y_*\in \Fun(\Delta, \mathcal C)\] induces
    \begin{itemize}[label={--}]
        \item A map of spectral sequences 
    \[ BKSS(c)\colon BKSS(X_*)\to BKSS(Y_*)\] which is on the $E_1$-page given by pushforward $\pi_{-q}X_p\to \pi_{-q}Y_p$.
    \item A map of filtered graded abelian groups $\pi_{*}\clim{\Delta}(c)\colon \pi_*(\clim{\Delta}X)\to \pi_*(\clim{\Delta}Y)$. \end{itemize}
    If $X_*$ and $Y_*$ are bounded above, then  
    \begin{center}
        \begin{tikzcd}
            E_{\infty}^{p,q}(X)\arrow[d, "BKSS(c)"']\arrow[r,"\cong"] & \gr^p\pi_{-(p+q)}(\clim{\Delta}X)\arrow[d,"\pi_*\clim{\Delta}(c)"]\\ 
            E_{\infty}^{p,q}(Y)\arrow[r,"\cong"] & \gr^p\pi_{-(p+q)}(\clim{\Delta}Y)
        \end{tikzcd}
    \end{center} commutes for all $p,q\in\mathbb Z$. 
\end{romanenum}
\end{proposition}
\begin{proof}
    The existence of the spectral sequence is \cite[Variant 1.2.4.9]{higheralgebra} in the category $\mathcal C^{\operatorname{op}}$ with the opposite $t$-structure. 
    The convergence statement is \cite[Proposition 1.2.4.5]{higheralgebra} in the category $\mathcal C^{\operatorname{op}}$ with opposite $t$-structure, shifted by degree $-N$. 
    Naturality follows from \cite[Proposition 1.2.4.5, Variant 1.2.4.9]{higheralgebra} in the stable category $\Fun(\Delta^1, \mathcal C)^{\operatorname{op}}$ with the pointwise $t$-structure (shifted by degree $-N$). 
\end{proof}
\begin{definition}
    A \emph{simplicial cover} of an object $X$ in a category $\mathcal C$ is an augmented simplicial object $X_*\to X\in \Fun(\Delta^{\operatorname{op}}_{+}, \mathcal C)$ which exhibits $X$ as $\colim{\Delta^{\operatorname{op}}}X_*=X$. 
    A morphism $p\colon Y\to X$ in a big topos $\topo{X}$ is an \emph{effective epimorphism} if its \v{C}ech nerve $\check{C}(p)$ is a simplicial cover of $X$. 
    \end{definition}

    \begin{rem}\label{characterizationseffectiveepi} By \cite[Proposition 7.2.1.14]{highertopostheory}, the \v{C}ech nerve of a map $p\colon Y\to X$ in a topos $\topo{X}$ is a simplicial cover if and only if $\tau_{\leq 0}f\colon \tau_{\leq 0}Y\to\tau_{\leq 0}X$ is an epimorphism in the $1$-topos $\tau_{\leq 0}\topo{X}$. 
    Since for a chain of topoi $\topo{X}_*\colon\Lambda\to \Pr^L$, $\topo{X_{\lambda}}\subseteq \topo{X}$ is closed under small colimits and finite limits (\cref{fullyfaithfulnessandpreservationoflimits}) and the left adjoint $\tau_{\leq 0}\topo{X}_{\lambda}\to\tau_{\leq 0}\topo{X}$ preserves epimorphisms, the same statement holds for big topoi.  
    \end{rem}

\begin{cor}[{\v{C}ech-to-cohomology spectral sequence}]\label{Bousfieldkanspectralsequencehomology}
    Suppose $X_*\to X$ is a simplicial cover in a big topos $\mathcal X$. 
    \begin{romanenum}
    \item Then \[\clim{\Delta}\cH{\topo{X}}(X_*,-)\cong \cH{\topo{X}}(X_*,-).\] 
    \item For $A\in \stab{\mathcal X}$, there is a spectral sequence $BKSS(X_*,A)$
    with
    \[ E^1_{p,q}\coloneqq \cH{\topo{X}}^q(X_p,A), \] and $d_1$-differential
    \[ d_{1}^{p,q}\coloneqq \sum_{i=0}^{p}(-1)^i (\partial_i)^*\colon H^q(X_p,A)\to H^q(X_{p+1},A)\] the alternating sum of the pullback along the face maps $\partial_i\colon X_{p+1}\to X_{p}.$

    \item If $A\in \Sp(\mt{X})_{\leq N}$ is a bounded above spectrum object, the spectral sequence $BKSS(X_*,A)$ converges to $\cH{\topo{X}}^{p+q}(X,A)$, filtered by \[\{F^n\cH{\topo{X}}^{p+q}(X,A)\coloneqq \ker(\cH{\topo{X}}^{p+q}(X,A)\to \cH{\topo{X}}^{p+q}(X_n,A))\}_{n\in\mathbb N_0}.\] 
    This filtration is bounded: for $k\in\mathbb Z$, 
    \[0=F^{k+N-1}\cH{\topo{X}}^k(X,A)\subseteq \ldots \subseteq F^0\cH{\topo{X}}^k(X,A)\subseteq F^{-1}\cH{\topo{X}}^k(X,A)=\cH{\topo{X}}^k(X,A).\]
    
    \item This is natural in the simplicial cover $X_*\to X$ and the spectrum object $A$:  
    \begin{itemize}
    \item[{--}] A map $A\to B$ of bounded above spectrum objects induces a map of spectral sequences \[ BKSS(X_*,A)\to BKSS(X_*,B)\] and a map of filtered graded abelian groups \[ \cH{\topo{X}}^*(X,A) \to \cH{\topo{X}}^*(X,B), \] and \begin{center}
        \begin{tikzcd}
            E^{p,q}_{\infty}(BKSS(X_*,A))\arrow[r,"\cong"] \arrow[d]& \gr^p\cH{\topo{X}}^{p+q}(X,A)\arrow[d]\\
            E^{p,q}_{\infty}(BKSS(X_*,B))\arrow[r,"\cong"] &\gr^p\cH{\topo{X}}^{p+q}(X,B).
        \end{tikzcd} 
    \end{center} commutes.
    \item[{--}] A map $f\colon (X_*\to X)\to (Y_*\to Y)$ of simplicial covers in $\topo{X}$ (i.e. of the underlying augmented simplicial objects) induces a map of spectral sequences \[ BKSS(Y_*,A)\to BKSS(X_*,A)\] and a map of filtered graded abelian groups \[ \cH{\topo{X}}^*(Y,A)\to \cH{\topo{X}}^*(X,A),\] and \begin{center}
        \begin{tikzcd}
            E^{p,q}_{\infty}(BKSS(Y_*,A))\arrow[r,"\cong"] \arrow[d]& \gr^p\cH{\topo{X}}^{p+q}(Y,A)\arrow[d]\\
            E^{p,q}_{\infty}(BKSS(X_*,A))\arrow[r,"\cong"] &\gr^p\cH{\topo{X}}^{p+q}(X,A)
        \end{tikzcd} 
    \end{center} commutes. 
    \end{itemize}
\end{romanenum}
\end{cor}
\begin{proof}
The functor $\map_{\stab{\topo{X}}}(-,-)\colon \stab{\topo{X}}^{\operatorname{op}}\times \stab{\topo{X}}\to\Sp$ preserves small limits in the first variable by \cref{spectralenrichmentcocontinuous}.  
In particular, if $X_*\to X \in \Fun(\Delta^{\operatorname{op}}\cup \{-1\}, \topo{X})$ is a simplicial cover, then $\cH{\topo{X}}(X,-)$ is the totalization of the cosimplicial object 
\[ \cH{\topo{X}}(X_*,-)\colon  \Delta\to \Fun(\stab{\topo{X}}, \Sp).\]

For $A\in\stab{\topo{X}}_{\leq N}$ and $X\in\topo{X}$, $\cH{\topo{X}}(X,A)\in \Sp_{\leq N}$ since for $k\geq N+1$, \[\Omega^{\infty-k}\cH{\topo{X}}(X,A)\cong \Omega^{\infty}\cH{\topo{X}}(X,\Sigma^kA)\cong \Map_{\topo{X}}(X, \Omega^{\infty-k}A)=0.\] 
The statement is now an immediate consequence of \cref{propositionbousfieldkanspectralsequence}. 
\end{proof}

In case $\stab{\mathcal X}$ is closed symmetric monoidal and admits countable products, we also obtain spectral sequences for internal cohomology: 
\begin{cor}\label{Bousfieldkanspectralsequenceinternal}
    Suppose $\mathcal X$ is a big topos such that the symmetric monoidal structure on $\stab{\mathcal X}$ is closed and $\stab{\topo{X}}$ has $\Delta^{\operatorname{op}}$-indexed colimits and $\Delta$-indexed limits.   
    If $X_*\to X$ is a simplicial cover in $\topo{X}$, there is a spectral sequence with $E_1$-page
    \[ E_1^{p,q}=\icH{\topo{X}}^q(X_p,A),\] which converges to $\icH{\topo{X}}^q(X,A)$ if $A\in \stab{\topo{X}}$ is bounded above.
    This is natural in the coefficient object $A$ and the simplicial cover $X_*\to X$. 
\end{cor}
The above corollary in particular applies to the category of accessible (hyper)\-sheaves on a (hyper)ac\-ces\-si\-ble explicit covering site by \cref{symmetricmonoidalstructureaccessiblesheavesofspectraclosed}. 

\begin{proof}
Since the symmetric monoidal structure on $\stab{\topo{X}}$ is closed, the tensor product preserves $\Delta^{\operatorname{op}}$-indexed colimits in both variables. Since $\stab{\topo{X}}$ has $\Delta$-indexed limits and $\Delta^{\operatorname{op}}$-indexed colimits, \cref{leftadjointsstableundercolimits} and \cite[Proposition 5.2.6.2]{highertopostheory} imply that $\icH{\topo{X}}(X,-)\cong \clim{\Delta}\icH{\topo{X}}(X_*,-)$. 
Suppose that $N\in\mathbb N_0$ such that $A\in\stab{\topo{X}}_{\leq N}$. By \cref{monoidalstructurecompatibletstructure}, for $T\in\stab{\topo{X}}_{>N}$ and $Y\in \topo{X}$, $T\otimes \Sigma^{\infty}_{+}Y\in \stab{\topo{X}}_{>N}$.
    This implies that 
    \[ \Map_{\stab{\topo{X}}}(T, \icH{\topo{X}}(Y,A))\cong \Map_{\stab{\topo{X}}}(T\otimes \Sigma^{\infty}_{+}Y,A)\cong \Map_{\stab{\topo{X}}}(\tau_{\leq  N}(Y\otimes \Sigma^{\infty}_{+}X),A)=0, \] i.e.\ $\icH{\topo{X}}(Y,A)\in \stab{\topo{X}}_{\leq N}$ for all $Y\in\topo{X}$. 
    The statement now follows from \cref{propositionbousfieldkanspectralsequence}. 
\end{proof}

\cref{constantsheafandtotalization} implies the following: 
\begin{cor}\label{constantsheafandtotalizationcohomology}Suppose $\topo{X}$ is a big topos such that $\stab{\topo{X}}$ admits countable products and denote by $c_{\Sp}\colon \Sp\to\stab{\topo{X}}$ the stabilization of the constant sheaf functor. 
If $X_*\to X$ is a simplicial cover in $\topo{X}$, for bounded above spectrum objects $A\in\stab{\topo{X}}_{<\infty}$, 
\[ c_{\Sp}\cH{\topo{X}}(X,A)\cong \clim{\Delta}(c_{\Sp}\cH{\topo{X}}(X_*,A))\] via the canonical map. 
    \end{cor} 
\begin{proof}For $A\in\stab{\topo{X}}_{\leq N}$ and $X\in\topo{X}$, $\cH{\topo{X}}(X,A)\in \Sp_{\leq N}$ since for $k\geq N+1$, \[\Omega^{\infty-k}\cH{\topo{X}}(X,A)=\Omega^{\infty}\cH{\topo{X}}(X,\Sigma^kA)=\Map_{\topo{X}}(X, \Omega^{\infty-k}A)=0.\] 
Hence by \cref{spectralenrichmentcocontinuous,constantsheafandtotalization}, 
\[ c_{\Sp}\cH{\topo{X}}(X,A)\cong c_{\Sp}(\clim{\Delta}\cH{\topo{X}}(X_*,A))\cong \clim{\Delta}(c_{\Sp}\cH{\topo{X}}(X_*,A)).\qedhere\] 
    \end{proof} 
    We will apply the following construction to compare continuous and condensed/solid group cohomology.  
\begin{definition}[Simplicial resolution]\label{definitionsimplicialcomplex}
    For a big topos $\topo{X}$, denote by $\mathbb Z[-]\colon \tau_{\leq 0}\topo{X}\to \Ab(\tau_{\leq 0}\topo{X})$ the free functor. 
    If $\pi \colon A\to B\in \tau_{\leq 0}\topo{X}$ is an epimorphism, denote by \[\mathbb Z[\check{C}(\pi)]\in \Fun(\Delta^{\operatorname{op}}, \Ab(\tau_{\leq 0}\topo{X}))\] the free simplicial abelian group on the \v{C}ech nerve of $\pi$. 
    By \cref{simplicialresolutionisresolution} below, its unnormalized Moore complex (\cite[Definition 1.2.3.8]{higheralgebra}) $S_*^{\pi}\in \Ch(\Ab(\tau_{\leq 0}\topo{X}))$ augments to  a resolution of $\mathbb Z[B]$ via $\mathbb Z[\pi]\colon S_0^{\pi}=\mathbb Z[A]\to \mathbb Z[B]$. 
    We call this the \emph{simplicial resolution} of $\mathbb Z[B]$ along $\pi$. 
    Explicitly, $S_*^{\pi}$ is given by \[ S_n^{\pi}\coloneqq \mathbb Z[\check{C}(p)_n]=\mathbb Z[\underbrace{A\times_{B}\times \ldots \times_B A}_{n+1 \text{ times}}]\, \text{ and boundaries } \partial_n\coloneqq \sum_{i=0}^{n} (-1)^i d_i\colon  S_n^{\pi}\to S_{n-1}^{\pi}\] the alternating sum of the face maps. If $B=*$, we write $S_*^{A}\coloneqq S_*^{A\to *}$ for the corresponding resolution of $\mathbb Z$. 
\end{definition}

\begin{lemma}\label{simplicialresolutionisresolution}The unnormalized Moore complex $S_*^{\pi}$ augments to a resolution of $\mathbb Z[B]$ via $\mathbb Z[\pi]\colon S_0^{\pi}=\mathbb Z[A]\to \mathbb Z[B]$.\end{lemma}

    \begin{proof}
        Choose an exhaustion $\topo{X}_*\colon\Lambda\to \Pr^L$ of $\topo{X}$ by topoi and $\lambda\in\Lambda$ with $A,B\in\topo{X}_{\lambda}\subseteq \topo{X}$.
        Since $\topo{X}_{\lambda}\subseteq \topo{X}$ is closed under finite limits, we can compute $\check{C}(\pi)$ in $\topo{X}_{\lambda}$.
        Since $\tau_{\leq 0}\topo{X}\to\Ab(\tau_{\leq 0}\topo{X})$ restricts to a functor $\tau_{\leq 0}\topo{X}_{\lambda}\to \Ab(\tau_{\leq 0}\topo{X}_{\lambda})$ and $\Ab(\tau_{\leq 0}\topo{X}_{\lambda})\subseteq \Ab(\tau_{\leq 0}\topo{X})$ is exact, it follows that the simplicial resolution of $\pi$ in $\Ab(\tau_{\leq 0}\topo{X})$ is the simplicial resolution of $\pi$ in $\Ab(\tau_{\leq 0}\topo{X}_{\lambda})$. By exactness of $\Ab(\tau_{\leq 0}\topo{X}_{\lambda})\subseteq \Ab(\tau_{\leq 0}\topo{X})$ we may therefore assume that $\topo{X}=\topo{X}_{\lambda}$ is a topos. Denote by \[H\mathbb Z[-]\colon\topo{X}\to \LMod{H\mathbb Z}{\stab{\topo{X}}}\] the composition of $\Sigma^{\infty}_{+}$ with the free $H\mathbb Z$-module functor.
    Since $H\mathbb Z[-]$ is a left adjoint, 
    \[\colim{\Delta^{\operatorname{op}}}H\mathbb Z[\check{C}(\pi)]\cong H\mathbb Z[\colim{\Delta^{\operatorname{op}}}\check{C}(\pi)]\cong H\mathbb Z[B].\] The second equivalence holds since $\pi$ is an effective epimorphism, cf.\ \cref{characterizationseffectiveepi}. 
    By \cref{propositionbousfieldkanspectralsequence} in the category $\stab{\topo{X}}^{\operatorname{op}}$ with the opposite $t$-structure/\cite[Proposition 1.2.4.5, Variant 1.2.4.9]{higheralgebra}, we have a spectral sequence with $E_1$-page $E_1^{p,q}=\pi_q(H\mathbb Z[\check{C}(\pi)_p])$ which converges completely to $\pi_{p+q}H\mathbb Z[B]$. 
    By \cref{freemoduloediscrete}, \[\pi_{q}H\mathbb Z[B]=\begin{cases}
    \mathbb Z[B] & q=0\\ 0 & \text{ else }
        \end{cases},\] and   \[E_{1}^{*,q}=\begin{cases}
            0 & q \neq 0 \\ 
            S_*^{\pi} & q=0
        \end{cases}.\] 
        
    This implies that \[E_{\infty}^{p,q}=E_2^{p,q}=\begin{cases}
            0 & q\neq 0 \\ 
                H_p(S_*^{\pi}) & q=0        \end{cases},\] and hence \[E_{\infty}^{p,0}=\gr^p\pi_{p}(H\mathbb Z[B])=\pi_{p}(H\mathbb Z[B])=\begin{cases}
        0 & p\neq 0 \\ 
        \mathbb Z[B] & p=0.
        \end{cases}\qedhere\]
    \end{proof}
The following spectral sequence is central for proving compactness properties of certain condensed abelian groups in \cref{section:compactcondensedabeliangroups} and for computations of derived internal Homs which enter our discussion of solid modules.
\begin{lemma}[{\cite[Corollary 4.8]{Scholzecondensed}}, Breen-Deligne spectral sequence]\label{Breendelignespectralsequence}
    Suppose $\topo{X}$ is a 1-topos. 
    There is a spectral sequence of functors 
    \[\Ab(\topo{X})^{\operatorname{op}}\times \mathcal D(\Ab(\topo{X}))\times \topo{X}^{\operatorname{op}}\to \Ab\] with $E_1^{p,q}(A,B,T)=\prod_{k=1}^{b_p}\Ext^{q}_{\Ab(\topo{X})}(\mathbb Z[A^{i_{k,p}}\times T],B)$, $b_p,i_{k,p}\in\mathbb N_0$.
    
    On $\topo{X}^{\operatorname{op}}\times \Ab(\topo{X})^{\operatorname{op}}\times \mathcal D(\Ab(\topo{X}))_{<\infty}$, this converges to \[ (A,B,T)\mapsto \Ext^{p+q}_{\Ab(\topo{X})}(\mathbb Z[T]\otimes A,B).\] 
    \end{lemma}
    \begin{proof}Choose a small 1-category $\mathcal C$ and a left-exact localization $L\colon \mathcal P(\mathcal C, \Set)\to\topo{X}$. 
        By \cite[Theorem 4.5]{Scholzecondensed}, every abelian group $A$ admits a resolution $S_*(A)\to A$ with terms \[S_n=\oplus_{k=1}^{b_n}\mathbb Z[A^{i_k}], i_k,b_n\in\mathbb N_0, n\in\mathbb N_0,\] and this is functorial in the abelian group $A$. (\cite{Scholzecondensed} attributes this to Eilenberg, MacLane, Breen and Deligne).
        Since $L$ induces a left-exact localization $\Fun(\mathcal C^{\operatorname{op}}, \Ab)\to \Ab(\mathcal X)$, for $A\in \Ab(\topo{X})$, the degreewise sheafification of 
        $\mathcal C^{\operatorname{op}}\ni c\mapsto S_*(A(c))$ is a resolution of $A$ with terms 
        $S_n=\oplus_{k=1}^{b_n}\mathbb Z[A^{i_k}]$. We denote this resolution by $S_*(A)$, it obviously depends functorially on $A$. 
    
        For $T\in\tau_{\leq 0}\topo{X}$, 
        $\mathbb Z[T]$ is the sheafification of \[\mathcal C^{\operatorname{op}}\ni c\mapsto \mathbb Z[T(c)].\] Since $\mathbb Z[T(c)]$ is flat for all $c\in\mathcal C$, $\mathbb Z[T]$ is flat in $\Ab(\topo{X})$ by \cite[\href{https://stacks.math.columbia.edu/tag}{Tag 03ET, Tag 03ES}]{stacks-project}. In particular,  $C_*(A,T)\coloneqq S_*(A)\otimes \mathbb Z[T]$ is a resolution of $A\otimes \mathbb Z[T]$ for all $T\in\tau_{\leq 0}\topo{X}$. 
        For $n\in\mathbb N_0$ denote by $C_*^{\leq n}(A,T)$ the stupid truncation of this resolution. 
        By \cite[Proposition 6.4.5.7]{highertopostheory}, there exists a hypercomplete 1-localic topos $\widehat{\topo{Y}}$ with $\tau_{\leq 0}\widehat{\topo{Y}}\cong \topo{X}$. \cref{derivedcategorymodules} and \cref{modulecategoriesareenriched} provide a $\mathcal D(\Ab)$-enrichment of $\mathcal D(\Ab(\topo{X}))\cong \LMod{H\mathbb Z}{\stab{\widehat{\topo{Y}}}}$. The $\mathcal D(\Ab)$-enriched mapping functor \[\map_{\mathcal D(\Ab(\topo{X}))}(-,-)\colon \mathcal D(\Ab(\topo{X}))^{\operatorname{op}}\times\mathcal D(\Ab(\topo{X}))\to\mathcal D(\Ab(\topo{X}))\] preserves small limits in both variables (\cref{continuityenrichmentorbits}), and in particular, \[\map_{\mathcal D(\Ab(\topo{X}))}(-,-)(C_*(-,-),-)\cong \clim{n\in\mathbb N_0}\map_{\mathcal D(\Ab(\topo{X}))}(C_*^{\leq n}(-,-),-).\]
        The long exact homology sequence implies that the canonical chain map \[\Cone(C_*^{\leq n}(-,-)\to C_*^{\leq n+1}(-,-))\to C_{n+1}(-,-)[n+1]\] is a quasi-isomorphism for all $n\in\mathbb N_0$. 
        We therefore obtain fiber sequences \begin{align*}\Omega^{n+1}\map_{\mathcal D(\Ab(\topo{X}))}(C_{n+1}(-_1,-_3),-_2)& \to \map_{\mathcal D(\Ab(\topo{X}))}(C_*^{\leq n+1}(-_1,-_3),-_2)\\ &\to \map_{\mathcal D(\Ab(\topo{X}))}(C_*^{\leq n}(-_1,-_3),-_2)\end{align*} in $\Fun(\Ab(\topo{X})^{\operatorname{op}}\times \mathcal D(\Ab(\topo{X}))\times \topo{X}^{\operatorname{op}}, \mathcal D(\Ab)).$
        By \cite[Definition 1.2.2.9]{higheralgebra}, the cofiltered object \[\ldots \to \map_{\mathcal D(\Ab(\topo{X}))}(C_*^{\leq n}(-,-),-)\to \map_{\mathcal D(\Ab(\topo{X}))}(C_*^{\leq n+1}(-,-),-)\to \ldots\] yields a spectral sequence with $E_1$-page \begin{align*}E_1^{p,q}(A,B,T)&=\pi_{-(p+q)}(\Fib\bigl(\map_{\mathcal D(\Ab(\topo{X}))}(C^{\leq p}(A,T),B)\to \map_{\mathcal D(\Ab(\topo{X}))}(C_*^{\leq p-1}(A,T),B)\bigr))\\&=\Ext^q_{\Ab(\topo{X})}(\oplus_{k=1}^{b_p}\mathbb Z[A^{i_{k,p}}\times T],B)\end{align*} in $\Fun(\Ab(\topo{X})^{\operatorname{op}}\times \mathcal D(\Ab(\topo{X}))\times \topo{X}^{\operatorname{op}}, \mathcal D(\Ab))$. 
        If $B\in\mathcal D(\Ab(\topo{X}))_{\leq N}$ is bounded above, then \[\map_{\mathcal D(\Ab(\topo{X}))}(C_*^{\leq n}(A,T),B)\in \mathcal D(\Ab)_{\leq N}\] for all $n\in\mathbb N_0$, whence by \cite[Proposition 1.2.4.5]{higheralgebra} the spectral sequence converges to \[\pi_{-(p+q)}\map_{\mathcal D(\Ab(\topo{X}))}(C_*(A,T),B)=\Ext^{p+q}_{\Ab(\topo{X})}(A\otimes \mathbb Z[T],B)\] with filtration 
        \[ \{F^s\coloneqq \ker(\Ext^*(C_*(A,T),B)\to \Ext^*(C_*^{\leq s}(A,T),B))\}_{s\geq -1}.\]
        This filtration is bounded, for $m\in\mathbb Z$: \[F^{m+N-1}\Ext^m(A\otimes \mathbb Z[T],B)=0, \text{ and } F^{-1}\Ext^m(A\otimes \mathbb Z[T],B)=\Ext^m(A\otimes\mathbb Z[T],B).\qedhere\] 
    \end{proof}
\newpage
\section{Condensed Mathematics}
In this chapter, we review cornerstones of condensed mathematics from \cite[Lectures 1-6]{Scholzecondensed}, paying particular attention to set-theoretic size issues. 
We begin by recalling basic results on the \textit{weight} (\cref{definitionweight}) of topological spaces and related categorical properties of the categories of compact Hausdorff spaces and profinite sets. A key fact is that the category of profinite sets is cocompactly generated and for every uncountable cardinal $\kappa$, the $\kappa$-cocompact objects are the $\kappa$-light profinite sets (\cref{profincoaccessible}). This implies that $(\Pro(\Fin),\text{surjections})$ forms a hyperaccessible explicit covering site (\cref{accessiblecoveringsitedefinition}, \cref{profiniteisexplicitcoveringsite}). The associated category of accessible hypersheaves is called condensed animae and denoted $\Cond{}(\an)$. 
By our discussion from the first chapter, the category of condensed animae is not presentable and in particular not a topos; however it retains many topos-like properties. We review this in \cref{section:condensedcategories}, it relies on the following description of condensed animae: 
For an uncountable cardinal $\kappa$, denote by $\Cond{\kappa}(\an)$ the category of condensed hypersheaves (\cref{kappacondensedtopology}) on the category of $\kappa$-light (\cref{definitionweight}) profinite sets $\Pro(\Fin)_{\kappa}$. 
Via sheafification of left Kan extension, these categories form a large filtered system with $\Cond{}(\an)\cong \colim{\kappa}\Cond{\kappa}(\an)$ (\cref{condensedasfilteredcolimitoverregular}).
Restricted to the poset of regular/strong limit cardinals, this is an exhaustion of $\Cond{}(\an)$ by topoi (\cref{definitionbigpresentable}, \cref{condensedfiltrationinfinitytopoi}). This description allows many statements about condensed categories to be deduced from their $\kappa$-condensed analogues. In particular, for a presentable category $\mathcal C$, many computations in $\Cond{}(\mathcal C)$, such as small limits and colimits, can be carried out in $\Cond{\kappa}(\mathcal C)\subseteq \Cond{}(\mathcal C)$ for some sufficiently large strong limit/regular cardinal $\kappa$. 
Some sources (e.g.\ \cite{Analyticstacks}, \cite{barwick2019pyknoticobjectsibasic}) work with $\kappa$-condensed instead of condensed categories. We therefore chose to formulate all results for both the condensed and the $\kappa$-condensed setting. 

The first part of this chapter is structured as follows: After recalling basic properties of the weight of topological spaces, we review $\kappa$-condensed categories and show that large classes of topological spaces embed fully faithfully into $\kappa$-condensed sets in \cref{section:kappacondensedcategories}. We apply this embedding to model $\kappa$-condensed animae as hypersheaves on other categories of topological spaces in  \cref{section:othermodelscondensedcat}.  
In \cref{section:condensedcategories} we specialize our discussion from the first chapter to establish basic structural properties of condensed categories and their relation to $\kappa$-condensed categories. This uses the results from \cref{section:propertiesweightsection}.
We then recall that there is a fully faithful functor from T1 compactly generated topological spaces to condensed sets which has a partially defined right adjoint. In \cref{section:condensedsheafcohomology}, we compare condensed cohomology of topological spaces with their sheaf cohomology. 
\cite{Scholzecondensed} constructed a comparison map from sheaf to condensed cohomology for compact Hausdorff spaces. We use the description of $\kappa$-condensed animae as hypersheaves on the category of all $\lambda$-small topological spaces to extend this to general ($\To$) topological spaces. Descent for condensed and sheaf cohomology along local-section covers then yields generalizations of the identifications between sheaf and condensed cohomology obtained in \cite{Scholzecondensed}.

In the second part of this chapter, we study the notion of solidity, the condensed analogue of non-archimedean completeness. As preparation, we record basic structural properties of derived and underived condensed module categories in \cref{section:condensedmodulecategories}, and establish some compactness properties of condensed abelian groups in \cref{section:compactcondensedabeliangroups}. 
In \cref{section:solidmodules}, we discuss the category of solid modules over a condensed ring $\algebra{R}$. We show that the category of solid $R$-modules is closed under limits, colimits and extensions and has enough projectives, and construct a (derived) solidification functor and a (derived) solid tensor product. 
To do so, we slightly adapt the discussion from \cite{Scholzecondensed} to obtain arguments that work for $\Cond{\kappa}(\Ab)$ for arbitrary uncountable cutoff cardinals $\kappa$. 
We compute the (derived and underived) solidification of free condensed abelian groups on compact Hausdorff spaces and CW-complexes. These computations imply that condensed cohomology with solid coefficients is homotopy invariant and lead to identifications of singular/sheaf and condensed cohomology with solid coefficients for large classes of topological spaces, see \cref{section:condensedcohomologywithsolidcoefficients}. 

\subsection{Some point-set topology}\label{section:propertiesweightsection}
In this section, we record some basic results on the \textit{weight} of topological spaces and related categorical properties of the categories of compact Hausdorff spaces and profinite sets. We will use this to relate condensed sets with topological spaces and to compare $\kappa$-condensed with condensed categories. 
 \begin{lemma}[{\cite[\href{https://stacks.math.columbia.edu/tag/08ZY}{Lemma 08ZY}]{stacks-project}}]\label{profinitetdch}
    The functor $\lim\colon \Pro(\Fin)\to \Top$ is fully faithful with essential image the totally disconnected compact Hausdorff spaces. 
 \end{lemma}
 We freely identify totally disconnected compact Hausdorff spaces with profinite sets in the sequel. 

 \begin{definition}
    A compact Hausdorff space $P$ is \emph{extremally disconnected} if it is projective in the category of compact Hausdorff spaces, i.e.\ if all continuous surjections $X\to P$ from a compact Hausdorff space $X$ admit a section. 
    Denote by $\edCH{}$ the (large) category of extremally disconnected compact Hausdorff spaces and continuous maps. 
 \end{definition}

 Recall that the forget functor from compact Hausdorff spaces to all topological spaces admits a left adjoint, called Stone-\v{C}ech compactification. 
 \begin{lemma}[{\cite{gleasonprojectives}}]\label{gleasonedch}
    \begin{romanenum}
    \item A compact Hausdorff space is extremally disconnected if and only if the closure of every open subset is open. 
    \item The Stone-\v{C}ech-compactification of a discrete topological space is extremally disconnected compact Hausdorff.
    \item Extremally disconnected compact Hausdorff spaces are totally disconnected. 
    \item Every compact Hausdorff space is a quotient of the Stone-\v{C}ech-compactification of its set of points. In particular, every extremally disconnected compact Hausdorff space is a retract of the Stone-\v{C}ech-compactification of its set of points. 
    \end{romanenum}
 \end{lemma}
 \begin{proof}
    The first statement is \cite[Definition 1.1, Theorem 2.5]{gleasonprojectives}. 
    Together with \cite[Theorem 4.1]{gleasonprojectives}, this implies that Stone-\v{C}ech-compactifications of discrete topological spaces are extremally disconnected. 

    Suppose that $X$ is extremally disconnected compact Hausdorff and $x,y\in X$ are distinct points. Choose disjoint open neighborhoods $x\in U_x,y\in U_y\subseteq X, U_x\cap U_y=\emptyset$. The closures $\overline{U_x}, \overline{U_y}$ are clopen by the first statement, and disjoint by \cite[Lemma 2.2]{gleasonprojectives}. In particular, $x$ and $y$ lie in different connected components of $X$, which shows that $X$ is totally disconnected. 

    By definition of the Stone-\v{C}ech-compactification, for a compact Hausdorff space $K$, the \textit{identity} $K^{\delta}\to K$ from the underlying set of $K$ to $K$ factors over a map $q\colon \beta(K^{\delta}) \to K$ which is surjective. 
    If $K$ is extremally disconnected compact Hausdorff, $q$ admits a section which exhibits $K$ as a retract of $\beta(K^{\delta})$. 
 \end{proof}
 The following definition will be central to us as it characterizes $\kappa$-cocompact profinite sets/compact Hausdorff spaces, see \cref{profincoaccessible}. 
\begin{definition}\label{definitionweight}
    The \emph{weight} $\wt(X)$ of a topological space $X$ is the smallest infinite cardinal such that $X$ admits a basis of cardinality $<\kappa$, i.e.\ a $\kappa$-small collection $\mathcal B$ of open subsets of $X$ such that every open subset of $X$ is a union of elements in $\mathcal B$. 
    We say that a topological space $X$ is $\kappa$-\emph{light} if $\wt(X)< \kappa$. 

    A topological space is \emph{light} if it admits a countable basis, i.e.\ is $\aleph_1$-light for $\aleph_1$ the smallest uncountable cardinal.
    
    For an uncountable cardinal $\kappa$ denote by \[\CH_{\kappa}\subseteq \CH, \, \Pro(\Fin)_{\kappa}\subseteq \Pro(\Fin), \text{ and }\edCH{\kappa}\subseteq \edCH{}\] the full subcategories on $\kappa$-light compact Hausdorff/profinite/extremally disconnected spaces, respectively.  
\end{definition}
\begin{rem}\label{weightsmallercardinality}
For a compact Hausdorff space $X$, $\wt(X)\leq \max\{|X|, \aleph_0\}$ by \cite[Theorem 3.1.21]{Engelkingtopology}. In particular, the categories $\edCH{\kappa}\subseteq \Pro(\Fin)_{\kappa}\subseteq \CH_{\kappa}$ are small for all small cardinals $\kappa$. 
\end{rem} 
\begin{lemma}\label{weightofsubspacesproducts}
    \begin{romanenum}
    \item For a subspace $A\subseteq X$, $\wt(A)\leq \wt(X)$.
    \item For a set of topological spaces $(X_i)_{i\in I}$,    \[ \wt(\prod_{i\in I}X_i)\leq \max\{|I|, \wt(X_i), i\in I\}.\] \end{romanenum}
\end{lemma}
\begin{proof}
This is immediate from the definition of the subspace and product topology. 
\end{proof}
\begin{lemma}[{\cite[Theorem 3.1.22]{Engelkingtopology}}]\label{weightofquotients}
If $X$ is a compact Hausdorff space and $Y\to X$ is a quotient map, then $\wt(X)\leq \wt(Y)$. 
\end{lemma}
\begin{lemma}[{\cite[Lemma 2.5]{aoki2024ktheoryringscontinuousfunctions}}]\label{embedintocubesofweightdimension}
    Suppose $\kappa$ is an uncountable cardinal.
    A compact Hausdorff space $X$ has $\wt(X)\leq \kappa$ if and only if there exists a continuous injection $X\hookrightarrow\prod_{\kappa}[0,1]$. 
\end{lemma}
\begin{proof}Suppose that $X$ is a compact Hausdorff space. 
    If $X\hookrightarrow\prod_{\kappa}[0,1]$ is a continuous injection, then $X$ is homeomorphic to its image, and hence every basis on $\prod_{\kappa}[0,1]$ pulls back to a basis on $X$, which shows that $\wt(X)\leq\wt(\prod_{\kappa}[0,1])=\kappa$.  
    Conversely, suppose that $\wt(X)\leq \kappa$. 
    For $f\in\mathcal C(X,[0,1])$ let $\mathcal U_f\coloneqq f^{-1}((0,1])\subseteq X$. We claim that these sets form a basis of the topology on $X$. 
    If $U\subsetneq X$ is open, choose $x\in U$. By Urysohn's lemma, there exists a continuous map $f_{x}\colon X\to [0,1]$ with $f_x(x)=1$ and $f_x|_{X\setminus U}=0$. Hence $U=\bigcup_{x\in X} U_{f_x}$, which shows that \[\mathcal B\coloneqq \{ U_f\, | \, f\in\mathcal C(X,[0,1])\}\] forms a basis of the topology on $X$.  By \cite[Theorem 1.1.15]{Engelkingtopology}, $\mathcal B$ has a subbasis of cardinality $\leq \wt(X)$, i.e.\ there exists a $\kappa$-small collection $\{f_i\}_{i\in I}\subseteq \mathcal C(X,[0,1])$ such that $(U_{f_i})_{i\in I}$ is a basis of the topology on $X$. It follows that \[\prod_{i\in I}f_i\colon X\to \prod_{i\in I}[0,1]\] is injective. 
    Indeed: Since $X$ is $\To$, for $x,y\in X$ with $x\neq y$ there exists an open neighborhood $x\in U_x$ such that $y\notin U_x$. Choose $i\in I$ with $x\in U_{f_i}\subseteq U_x$. Then $f_i(x)\neq f_i(y)$. 
\end{proof}
\begin{lemma}\label{profinitesetseparatedbyfunctions}
    For a profinite set $X$, the canonical map $X\to \prod_{f\in\mathcal C(X, \{0,1\})}\{0,1\}$ is injective.
\end{lemma}
\begin{proof}
    Suppose that $x,y\in X$ with $x\neq y$. 
    Then $X\setminus \{x\}, X\setminus \{y\}$ form an open cover of $X$. As $X$ is profinite, there exists a refinement to an open cover $X=\bigsqcup_{i\in I}U_i$ with \[U_i\subseteq X\setminus\{x\}\text{  or }U_i\subseteq X\setminus \{y\}\text{ for all } i \in I, \] see e.g.\ \cite[\href{https://stacks.math.columbia.edu/tag/08ZZ}{Lemma 08ZZ}]{stacks-project}.
    Choose $i_x\in I$ with $x\in U_i$. 
    Then \begin{align*}1_{U_{i_x}}\colon X& \to \{0,1\} \\t&\mapsto \begin{cases}
        1 & t \in U_{i_x} \\ 0 & \text{ else }
    \end{cases}\end{align*} is continuous, $1_{U_{i_x}}(x)=1$ and $1_{U_{i_x}}(y)=0$. 
    This shows that \[X\to \prod_{f\in\mathcal C(X, \{0,1\})}\{0,1\}\] is injective. 
\end{proof}
\begin{cor}\label{propertiesweight}
    For a profinite set $X$, the following hold: 
    \begin{romanenum}
    \item $\wt(X)=\max\{|\mathcal C(X, \{0,1\})|, \aleph_0\}$.  
    \item There exists a continuous injection 
    \[ X\hookrightarrow \prod_{\wt(X)}\{0,1\}.\] 
    \end{romanenum}

    In particular, $|X|\leq 2^{\wt(X)}$ and $\wt(X)\leq \max\{2^{|X|}, \aleph_0\}$. 
\end{cor}
\begin{proof}
    We first show $i)$. By \cref{profinitesetseparatedbyfunctions},  $X\hookrightarrow\prod_{\mathcal C(X, \{0,1\})}\{0,1\}$ is injective. 
    Since $X$ is compact and $\prod_{\mathcal C(X, \{0,1\})}\{0,1\}$ is Hausdorff, this is a homeomorphism onto its image, so in particular 
    \[\wt(X)\leq \wt(\prod_{\mathcal C(X, \{0,1\})}\{0,1\})\leq \max\{\aleph_0,|\mathcal C(X, \{0,1\})|\}.\] 
    We now show that $|\mathcal C(X, \{0,1\})|\leq \wt(X)$. 
    Suppose that $\mathcal B$ is a basis of the topology on $X$ with $|\mathcal B|\leq \wt(X)$ and $f\in \mathcal C(X, \{0,1\})$. Since $f^{-1}(\{0\})$ is open in $X$, there exist a subset $\mathcal B_f\subseteq \mathcal B$ with $f^{-1}(\{0\})=\cup_{U\in \mathcal B_f}U$. 
    The closed subspace $f^{-1}(\{0\})\subseteq X$ is compact, so we can choose $\mathcal{B}_f$ to be finite. 
    This determines an injection \[\mathcal C(X, \{0,1\})\hookrightarrow \{ F\subseteq \mathcal B \, | \, F \text{ finite }\}\] to the set of finite subsets of $\mathcal B$. 
    In particular, 
    \[|\mathcal C(X, \{0,1\})|\leq |\bigcup_{n\in\mathbb N_0}\mathcal B^n|\leq \wt(X)\times \aleph_0 \leq \wt(X)^2= \wt(X)\] since $\wt(X)$ is infinite.
    The second statement follows from the first and \cref{profinitesetseparatedbyfunctions}. 
    \cref{profinitesetseparatedbyfunctions} also implies that $|X|\leq 2^{|\mathcal C(X, \{0,1\})|}$ and obviously, $|\mathcal C(X, \{0,1\})|\leq 2^{|X|}$. 
\end{proof}

\begin{cor}\label{kappalightkappasmallimit}
Suppose that $\kappa$ is an uncountable regular cardinal. 
A profinite set $X$ is $\kappa$-light if and only if it can be realized as a $\kappa$-small filtered limit of finite sets. 
\end{cor}
\begin{proof}
    If $I$ is a $\kappa$-small filtered category, for $F\colon I\to\Fin$, 
    $\clim{I}F\hookrightarrow\prod_{i\in I}F(i)$, whence 
    \[\wt(\clim{I}F)\leq \wt(\prod_{i\in I}F)\leq |I|\times \aleph_0< \kappa.\] 
    Conversely, if $X$ is a $\kappa$-light profinite set, for $F\subseteq \mathcal C(X, \{0,1\})$ finite let \[X_F\coloneqq \im(\prod_{f\in F} f\colon X\to \prod_{f\in F}\{0,1\}).\] 
    By \cref{profinitesetseparatedbyfunctions}, the canonical map $l\colon X\to \clim{\substack{F\subseteq \mathcal C(X, \{0,1\})\\ F\text{ finite }}} X_F$ is injective. 
    As for $G\subseteq \mathcal C(X,\{0,1\})$ finite, 
    \[X\xrightarrow{l}\clim{\substack{F\subseteq \mathcal C(X, \{0,1\})\\ F\text{ finite }}} X_F\xrightarrow{\pi_G}X_G\] is surjective, it follows that $l$ is surjective and hence a homeomorphism since both sides are compact Hausdorff. 
    This exhibits $X$ as $\kappa$-small filtered limit of finite sets.
\end{proof}

\begin{lemma}\label{compacthausdorffquotientoftdch}
    Suppose $\kappa$ is an uncountable cardinal. A compact Hausdorff space is $\kappa$-light if and only if it is the quotient of a $\kappa$-light profinite set.  
    A compact Hausdorff space is light if and only if it is metrizable. 
\end{lemma}
\begin{proof}
Suppose $X$ is a $\kappa$-light compact Hausdorff space. By \cref{embedintocubesofweightdimension}, there exists a continuous injection $X\hookrightarrow \prod_{\wt(X)}[0,1]$. 
Since $[0,1]$ is compact metrizable, by the Hausdorff-Alexandroff theorem there exists a quotient map $\prod_{\mathbb N}\{0,1\}\to [0,1]$ from the Cantor set. 
Let $S\coloneqq \prod_{\wt(X)}\prod_{\mathbb N}\{0,1\}$.
Then $\wt(S)\leq \max\{\aleph_0, \wt(X)\}<\kappa$.  
Let $T\coloneqq S\times_{\prod_{\wt(X)}[0,1]}X$. 
This is a closed subspace of $S$ and in particular profinite of $\wt(T)\leq \wt(S)<\kappa$. 
The surjection $S\to \prod_{\kappa}[0,1]$ pulls back to a surjection $T\to X$ which is a quotient map since $T$ and $X$ are compact Hausdorff. The converse implication is \cref{weightofquotients}.

By \cref{profinitesetseparatedbyfunctions} and \cref{propertiesweight}, every light profinite set $X$ is homeomorphic to a subspace of the Cantor set $\prod_{\mathbb N}\{0,1\}$ and in particular metrizable.     
Conversely, by the Hausdorff-Alexandroff theorem, every (non-empty) compact metrizable space is a quotient of the Cantor set $\{0,1\}^{\mathbb N}$ which is light profinite. 
\cref{weightofquotients} now implies that compact metrizable spaces are light. 
\end{proof}

\begin{lemma}\label{kappalightforstronglimit}
    Suppose $\kappa$ is a strong limit cardinal.\footnote{This means that $\lambda<\kappa\Rightarrow 2^{\lambda}<\kappa$.} 
    For a compact Hausdorff space $X$, the following are equivalent: 
    \begin{romanenum}
        \item\label{wtsmallerkappa} $X$ is $\kappa$-light, i.e.\ $\wt(X)<\kappa$.
        \item\label{sizesmallerkappa} $X$ is $\kappa$-small, i.e.\ $|X|<\kappa$. 
        \item\label{quotientedchkappa} $X$ is the quotient of a $\kappa$-light extremally disconnected compact Hausdorff space. 
    \end{romanenum}
\end{lemma}
\begin{proof}
    Suppose that $\wt(X)<\kappa$. Then $X$ is the quotient of a $\kappa$-light profinite set $T$ by \cref{compacthausdorffquotientoftdch} and in particular, $|X|\leq |T|$. 
    By \cref{propertiesweight}, $|T|\leq 2^{\wt(T)}$ and since $\kappa$ is a strong limit cardinal, $2^{\wt(T)}<\kappa$. This shows that $\ref{wtsmallerkappa}\Rightarrow \ref{sizesmallerkappa}$. 
   As the collection of all open subsets of $X$ is a basis, $\wt(X)\leq \max\{2^{|X|}, \aleph_0\}$, which shows $\ref{sizesmallerkappa}\Rightarrow \ref{wtsmallerkappa}$. 
   
   The implication $\ref{quotientedchkappa}\Rightarrow \ref{wtsmallerkappa}$ holds by \cref{weightofquotients}. 
    Suppose now that $|X|<\kappa$. 
    By \cref{gleasonedch}, $X$ is the quotient of $\beta(X^{\delta})$, the Stone-\v{C}ech compactification of its set of points. 
    By construction of the Stone-\v{C}ech-compactification, 
    $|\beta(X^{\delta})|\leq 2^{|X|}<\kappa$ (since $\kappa$ is a strong limit cardinal), whence $\beta(X^{\delta})$ is $\kappa$-light by the equivalence of \ref{sizesmallerkappa} and \ref{wtsmallerkappa}. This proves $\ref{sizesmallerkappa}\Rightarrow \ref{quotientedchkappa}$. 
\end{proof}
The following result will be central for our discussion of condensed categories and motivated our discussion of weight of topological spaces. 
\begin{lemma}[{\cite[Lemma 2.1.5]{LucasMannthesis}, \cite[Proposition 2.11]{aoki2024ktheoryringscontinuousfunctions}}]\label{profincoaccessible}
\vspace{0pt}\noindent
    \begin{romanenum}
    \item The category $\Pro(\Fin)^{\operatorname{op}}$ is compactly generated. 
    For an uncountable regular cardinal $\kappa$, the $\kappa$-compact objects of $\Pro(\Fin)^{\operatorname{op}}$ are precisely the $\kappa$-light profinite sets.

    \item The category $\CH^{\operatorname{op}}$ is $\aleph_1$-compactly generated. 
    For an uncountable regular cardinal $\kappa$, the $\kappa$-compact objects are precisely the $\kappa$-light compact Hausdorff spaces. 
\end{romanenum}
\end{lemma}
\begin{proof}The first statement is essentially \cite[Lemma 2.1.5]{LucasMannthesis}.
    The category \[\Pro(\Fin)^{\operatorname{op}}=\Ind(\Fin^{\operatorname{op}})=\Ind_{\aleph_0}(\Fin^{\operatorname{op}})\] is $\aleph_0$-compactly generated by the objects of $\Fin^{\operatorname{op}}$.   
    By \cite[Lemma 2.1.5]{LucasMannthesis}, $X\in\Pro(\Fin)$ is $\kappa$-cocompact if and only if it is a $\kappa$-small filtered limit of finite sets, which is equivalent to being $\kappa$-light by \cref{kappalightkappasmallimit}. 

    By \cite[4.7, 6.5c]{GabrielUlmer}, $\CH^{\operatorname{op}}$ is $\aleph_1$-compactly generated by $[0,1]$ and by \cite[Proposition 2.11]{aoki2024ktheoryringscontinuousfunctions}, the $\kappa$-compact objects are precisely the $\kappa$-light compact Hausdorff spaces.  
\end{proof}
We will use the following variants of continuity to compare topological spaces with ($\kappa$-)condensed sets. 
\begin{definition}\label{definitionkcontinuous}
    Suppose $\kappa$ is an uncountable cardinal and $X,Y$ are topological spaces. 
    \begin{romanenum}\item A map $X\to Y$ is $\kappa$-\emph{continuous} if for all $\kappa$-light compact Hausdorff spaces $K$ and all continuous maps $K\to X$, the composition $K\to X\to Y$ is continuous.
    Denote by $\mathcal C_{\kappa}(X,Y)$ the set of $\kappa$-continuous functions $X\to Y$ and by $\kappacont$ the large category of all topological spaces and $\kappa$-continuous maps. 

    \item A map $X\to Y$ is $k$-\emph{continuous} if for all compact Hausdorff spaces $K$ and all continuous maps $K\to X$, the composition $K\to X\to Y$ is continuous. Denote by $\mathcal C_k(X,Y)$ the set of $k$-continuous maps $X\to Y$.
    
    \item A topological space $X$ is $\kappa$-\emph{compactly generated} if every $\kappa$-continuous map $X\to Y$ to a topological space $Y$ is continuous, i.e.\ $\mathcal C(X,-)=\mathcal C_{\kappa}(X,-)$.
     
    A topological space is \emph{compactly generated} if every $k$-continuous map $X\to Y$ to a topological space $Y$ is continuous, i.e.\ $\mathcal C(X,-)=\mathcal C_{k}(X,-)$.
    \end{romanenum}
\end{definition}

For cardinals $\lambda\geq\kappa$, we have inclusion 
\[ \mathcal C(-,-)\subseteq \mathcal C_k(-,-)\subseteq \mathcal C_{\lambda}(-,-)\subseteq \mathcal C_{\kappa}(-,-).\] 

\begin{lemma}
    A continuous map is $\aleph_1$-continuous if and only if it is sequentially continuous. 
\end{lemma}
\begin{proof}
    Suppose that $f\colon X\to Y$ is sequentially continuous. If $K\to X$ is a continuous map, the composite $K\to X\to Y$ is sequentially continuous as continuous maps are sequentially continuous and compositions of sequentially continuous maps are sequentially continuous. Since all light profinite sets are metrizable (\cref{propertiesweight}), this implies that $f$ is $\aleph_1$-continuous. Conversely, denote by $\mathbb N\cup\{\infty\}$ the one-point compactification of the natural numbers (with discrete topology). This is compact Hausdorff and totally disconnected, hence profinite and \[\{ \{ k\} \, |\, k\in\mathbb N_0\}\cup \{ \mathbb N_0\setminus \{1, \ldots,k\}\cup \{\infty\}, k\in\mathbb N_0\}\] is a countable basis for its topology. If $(a_n)_{n\in\mathbb N_0}$ is a convergent sequence in $X$ with limit $a_{\infty}$, then \begin{align*}a_*\colon \mathbb N\cup\{\infty\}&\to X\\ n &\mapsto a_n\end{align*} defines a continuous map. By definition of the one-point compactification, $f\circ a_*\colon \mathbb N\cup\{\infty\}\to Y$ is continuous if and only if $f(a_n)_n$ converges to $f(a_{\infty})$. This shows that every $\aleph_1$-continuous map is sequentially continuous. 
\end{proof}
\begin{ex}
If $X$ is Hausdorff and compactly generated, then $X$ is $\kappa$-compactly generated for all cardinals $\kappa>\wt(X)$ since \[X\cong\colim{\substack{K\to X\\ K \in \CH}}K \cong \colim{\substack{K\subseteq X\\ K\in \CH}}K \cong \colim{\substack{K\subseteq X\\ K \in\CH_{\kappa}}}K\cong\colim{\substack{K\to X\\ K\in \CH_{\kappa}}}K.\]
The first equivalence holds since $X$ is compactly generated, the second since the compact subspaces of $X$ form a cofinal subcategory of $\CH{}_{/X}$, the third since all subspaces of $X$ are $\kappa$-light, and the last since the compact subspaces form a cofinal subcategory of ${\CH_{\kappa}}_{/X}$. 
\end{ex}

We end this section with two results on the space of connected components of a compact Hausdorff space/topological group which we will use in our comparison of continuous and solid group cohomology in \cref{section:continuousandsolidgroupcohomology}. 
The following is a consequence of \cite[Theorem 6.1.23]{Engelkingtopology}. 
\begin{lemma}\label{connectedcomponentsofchisch}
For a compact Hausdorff space $X$, the space of connected components $\pi_0X$, endowed with the quotient topology from the projection $X\to\pi_0X$, is profinite. 
\end{lemma}
\begin{proof}
As quotient of a compact space, $\pi_0X$ is compact. 
Next, we show that $\pi_0X$ is Hausdorff.  

By \cite[Theorem 6.1.23]{Engelkingtopology}, the connected component of a point $x\in X$ is the intersection of all clopen subsets of $X$ containing $x$. In particular, if $x_1,x_2\in X$ lie in distinct connected components, there exist clopen subsets $x_1\in C_1\subseteq X,x_2\in C_2\subseteq X$ with $x_1\notin C_2, x_2\notin C_1$. 
Denote by $q\colon X\to\pi_0X$ the quotient map. For $i=1,2$: As $C_i$ is clopen, for a connected component $K\subseteq X$, $C_i\cap K=\emptyset$ or $K\subseteq C_i$, whence $q^{-1}(q(C_i))=C_i$. In particular, $q(C_i), q(C_2)\subseteq \pi_0X$ are clopen and disjoint, which shows that $\pi_0X$ is Hausdorff. 
It follows from \cite[Theorem 6.1.23]{Engelkingtopology} that $X$ is also totally disconnected. 
\end{proof}
\begin{cor}\label{weightofconnectedcomponents}
    Suppose $G$ is a topological group, and endow its set of connected components $\pi_0G$ with the quotient topology induced by $G\to \pi_0G$. 
    This is a totally disconnected and Hausdorff topological group.
    If $\pi_0G$ is locally compact, then $\wt(\pi_0G)\leq \wt(G)$.  
    \end{cor}
    \begin{proof}By continuity of the group operation, the connected component of the identity in a topological group is a normal subgroup, whence $\pi_0G$ with quotient topology is a topological group. As the connected component of the identity is closed, $\pi_0G$ is $\To$ and hence Hausdorff. $\pi_0G$ is totally disconnected by e.g.\ \cite[Lemma 2.9]{stroppel2006lcgrps}.  

    Suppose now that $\pi_0G$ is locally compact. By van Dantzig's theorem, $\pi_0G$ is a coproduct of profinite sets $\pi_0G=\sqcup_{i\in I}S_i$. Choose $I$ such that $S_i\neq \emptyset$ for all $i\in I$ and denote by $q\colon G\to \pi_0G$ the quotient map. 
    Then $G=\sqcup_{i\in I}q^{-1}(S_i)$. 
        Since $q^{-1}(S_i)\neq \emptyset$ for all $i\in I$, $|I|\leq \wt(G)$. 
        As $q^{-1}(S_i)\to S_i$ is a quotient map for all $i\in I$, it follows from \cref{weightofquotients} that \[\wt(\pi_0G)\leq \sum_{i\in I}\wt(S_i)\leq \sum_{i\in I}\wt(q^{-1}(S_i))\leq |I|\wt(G)\leq \wt(G).\qedhere\] 
    \end{proof}
\subsection{\texorpdfstring{$\kappa$}{k}-condensed categories}\label{section:kappacondensedcategories}
\begin{definition}\label{kappacondensedtopology}
    For an uncountable cardinal $\kappa$, the condensed topology on $\Pro(\Fin)_{\kappa}$ is the Grothendieck topology generated by finite, jointly surjective covering families, i.e.\ \[S\subseteq \oc{{\Pro(\Fin)_{\kappa}}}{X}\] is a covering sieve if and only if there exists a finite set $\{X_i\to X\}_{i=1}^n\subseteq S$ such that $\sqcup_{i=1}^n X_i\to X$ is onto. This defines a Grothendieck topology by \cref{condensedtopologyonprofinite}.

    The category of $\kappa$-\emph{condensed animae} \[\Cond{\kappa}(\an)\coloneqq \hypershv_{\condo}(\Pro(\Fin)_{\kappa}, \an)\] is the topos of hypercomplete \cite[p. 666]{highertopostheory} condensed sheaves on $\Pro(\Fin)_{\kappa}$. 

    For a presentable category $\mathcal C$, define the category of $\kappa$-condensed objects in $\mathcal C$ by 
    \[ \Cond{\kappa}(\mathcal C)\coloneqq \Cond{\kappa}(\an)\otimes_{\Pr^L}\mathcal C, \] where $\otimes_{\Pr^L}$ denotes the tensor product of presentable categories (\cite[Proposition 4.8.1.15]{higheralgebra}).

    For $\kappa=\aleph_1$ the smallest uncountable cardinal, 
    \[ \Cond{\light}(\mathcal C)\coloneqq \Cond{\aleph_1}(\mathcal C)\] is called category of light condensed objects in $\mathcal C$. 
\end{definition}
\begin{rem}
$\kappa$-condensed categories for strong limit cardinals $\kappa$ were introduced in \cite{Scholzecondensed}. Our definition is equivalent to theirs by \cref{kappalightforstronglimit}.  
\cite{Analyticstacks} discussed light condensed categories and  
\cite{LucasMannthesis} worked with $\kappa$-condensed categories for uncountable regular cardinals. $\kappa$-condensed categories for arbitrary uncountable cardinals $\kappa$ were introduced in \cite{aoki2024ktheoryringscontinuousfunctions}. 
By \cref{identifysheaveswithtensorproductinprl}, for a presentable category $\mathcal D$, \[ \Cond{\kappa}(\mathcal D)\cong \Fun^{\operatorname{lim}}(\Cond{\kappa}(\an)^{\operatorname{op}}, \mathcal D).\] One can more generally define  $\Cond{\kappa}(\mathcal D)\coloneqq \Fun^{\operatorname{lim}}(\Cond{\kappa}(\an)^{\operatorname{op}}, \mathcal D)$ for categories $\mathcal D$ with small limits. 
\end{rem}
\begin{ex}\label{explicitcharacterizationforcondensedset}
A functor $F\colon \Pro(\Fin)_{\kappa}^{\operatorname{op}}\to \Set\subseteq \an$ is a $\kappa$-condensed set/anima if and only if the following hold: 
\begin{romanenum}
    \item $F$ preserves finite products, i.e.\ $F(\emptyset)=*$ and for $A_1,A_2\in \Pro(\Fin)_{\kappa}$, the canonical maps $F(A_1\sqcup A_2)\to F(A_k)$ exhibit $F(A_1\sqcup A_2)=F(A_1)\times F(A_2)$. 
    \item For all continuous surjections $q\colon X\to Y\in\Pro(\Fin)_{\kappa}$, \[F(Y)\cong \Eq(F(X)\rightrightarrows F(X\times_Y X)).\]
\end{romanenum}
This follows from \cite[Proposition A.5.7.2]{SAG} since $\Set\subseteq \an$ is closed under limits and $\Delta_{s,\leq 1}\subseteq \Delta_{s}$ is left 1-cofinal, see \cite[Example 6.14]{du2023reshapinglimitdiagramscofinality}. 
In particular, for a topological space $X$, \[\underline{X}_{\kappa}\coloneqq \mathcal C(-,X)\colon\Pro(\Fin)_{\kappa}^{\operatorname{op}}\to\Set\] defines a condensed functor as every continuous surjection between compact Hausdorff spaces is a quotient map.  
\end{ex}
 More generally, \cref{sheafconditionexplicitcoveringsite} implies that $\kappa$-condensed objects in a presentable category $\mathcal C$ can be identified with functors $\Pro(\Fin)_{\kappa}^{\operatorname{op}}\to\mathcal C$ satisfying a descent condition for \textit{condensed hypercovers}: 
\begin{definition}
    \label{Definitioncondensedhypercover} 
    An augmented semi-simplicial object $U_*\colon \Delta_{s,+}^{\operatorname{op}}\to\Pro(\Fin)$ is a \emph{condensed hypercover} if for all $n\in\mathbb N_0$, the map 
    $U([n])\to M_nU([n])\in\Pro(\Fin)$ induced by the unit of $i^*\dashv i_*$ is surjective. 

    If $\mathcal C$ is a presentable category, denote by \[\Fun_{\condo}(\Pro(\Fin)_{\kappa}^{\operatorname{op}}, \mathcal C)\subseteq \Fun(\Pro(\Fin)_{\kappa}^{\operatorname{op}}, \mathcal C)\] the full subcategory on finite products-preserving functors $F\colon\Pro(\Fin)_{\kappa}^{\operatorname{op}}\to\mathcal C$ such that for all condensed hypercovers $U_*\colon \Delta_{s,+}^{\operatorname{op}}\to\Pro(\Fin)_{\kappa}$, the canonical map 
$F(U_{-1})\to \clim{\Delta_{s}}F(U_*)$ is an equivalence. We call such functors \emph{condensed functors}.
\end{definition}
\begin{cor}[{\cite[Proposition A.3.2.1, Proposition A.5.7.2]{SAG}}] \label{sheafconditionexplicitkappacondensed}
    For a presentable category $\mathcal C$, restriction along the Yoneda embedding \[\Pro(\Fin)_{\kappa}^{\operatorname{op}}\hookrightarrow \Cond{\kappa}(\an)^{\operatorname{op}}\] defines a fully faithful functor 
\[\Cond{\kappa}(\mathcal C)\hookrightarrow \Fun(\Pro(\Fin)_{\kappa}^{\operatorname{op}}, \mathcal C).\] Its essential image consists of the condensed functors. 
\end{cor}
\begin{proof}
Denote by $S\subseteq \Pro(\Fin)_{\kappa}$ the wide subcategory on continuous surjections $X\to Y$.
As $\Pro(\Fin)_{\kappa}$ has finite coproducts and pullbacks and the forget functor $\Pro(\Fin)_{\kappa}\subseteq \CH\to \Set$ preserves coproducts and pullbacks and is conservative, it is straightforward to check that $(\Pro(\Fin)_{\kappa},S)$ is an explicit covering site and the associated Grothendieck topology (\cite[Proposition A.3.2.1]{SAG}) is the condensed topology on $\Pro(\Fin)_{\kappa}$. 
\cref{sheafconditionexplicitcoveringsite} now implies that pullback along the \textit{hypersheafified Yoneda embedding} \[\Pro(\Fin)_{\kappa}\to \Fun(\Pro(\Fin)_{\kappa}^{\operatorname{op}}, \an)\to\Cond{\kappa}(\an)\] is fully faithful with essential image the condensed functors. 
By \cref{explicitcharacterizationforcondensedset}, the Yoneda embedding $\Pro(\Fin)_{\kappa}\to \Fun(\Pro(\Fin)_{\kappa}^{\operatorname{op}},\an)$ factors over \[\Cond{\kappa}(\an)\subseteq\Fun(\Pro(\Fin)_{\kappa}^{\operatorname{op}},\an),\] i.e. no hypersheafification is necessary. \end{proof}

\begin{lemma}\label{globalsectionsleftadjointcondensed}
The global sections functor $\Gamma\colon \Cond{\kappa}(\an)\to \an$ is a left adjoint. 
\end{lemma}
\begin{proof} By the adjoint functor theorem (\cite[Corollary 5.5.2.9]{highertopostheory}), it suffices to show that the global sections functor preserves small colimits. 
    Denote by $\mathcal P(\Pro(\Fin)_{\kappa})\xrightarrow{L_{\cond{\kappa}}}\Cond{\kappa}(\an)$ the left adjoint of the inclusion. 
We will show that the global sections functor \[\Gamma_0\colon \mathcal P(\Pro(\Fin)_{\kappa})\to \an\] of $\mathcal P(\Pro(\Fin)_{\kappa})$ factors over $\mathcal P(\Pro(\Fin)_{\kappa})\xrightarrow{L_{\cond{\kappa}}}\Cond{\kappa}(\an)\to \an$. Then the induced functor $\Cond{\kappa}(\an)\to\an$ preserves small colimits and is the global section functor: 
It is left adjoint to the composite \[ \an\xrightarrow{c}\mathcal P(\Pro(\Fin)_{\kappa})\xrightarrow{L_{\cond{\kappa}}}\Cond{\kappa}(\an)\] which equals the constant sheaf functor since it is cocontinuous and preserves the terminal object. 

We first show that $\Gamma_0$ factors over the category $\Shv_{\condo}(\Pro(\Fin)_{\kappa})$ of not necessarily hypercomplete sheaves. As $\Gamma_0$ is evaluation at $*$ and for $X\in\Pro(\Fin)_{\kappa}$, $\mathcal C(-,X)$ is a condensed sheaf (\cref{explicitcharacterizationforcondensedset}), the composite \[j\colon \Pro(\Fin)_{\kappa}\hookrightarrow \Shv_{\condo}(\Pro(\Fin)_{\kappa})\xrightarrow{\Gamma_0}\an\] of $\Gamma_0$ with the (hypersheafified) Yoneda embedding factors as $\Pro(\Fin)_{\kappa}\xrightarrow{f}\Set\subseteq \an$, where $f\colon\Pro(\Fin)\subseteq \CH\to \Set$ denotes the forget functor. In particular, $j$ preserves finite limits. 
If $\{X_i\to X\}_{i\in I}$ is a condensed cover, there exists a finite subset $F$ such that $\sqcup_{i\in F}X_i\to X$ is surjective. This implies that $\sqcup_{i\in I}j(F_i)\to j(X)$ is a surjection, i.e. an epimorphism in $\Set$ and hence an effective epimorphism in $\an$. It now follows from \cite[Proposition 6.2.3.20]{highertopostheory} that $\Gamma_0$ factors over $\Shv_{\condo}(\Pro(\Fin)_{\kappa})$. 
The induced functor $\Gamma_1\colon \Shv_{\condo}(\Pro(\Fin)_{\kappa})\to \an$ preserves small limits and colimits. Hence by \cite[Propositions 5.5.6.28, 6.5.1.12]{highertopostheory}, $\Gamma_1$ preserves $\infty$-connective morphisms. 
As $\an$ is hypercomplete, this implies that $\Gamma_1$ factors over a functor $\Gamma\colon \Cond{\kappa}(\an)=\hypershv_{\condo}(\Pro(\Fin)_{\kappa})\to \an$.
\end{proof}
We now recall that $\kappa$-condensed sets are a good approximation of ($\kappa$-compactly generated) topological spaces.  
For an uncountable cardinal $\kappa$, denote by $\kappacont$ the category of topological spaces and $\kappa$-continuous maps (\cref{definitionkcontinuous}). 
\begin{proposition}[{\cite[Proposition 1.7]{Scholzecondensed}}]\label{kappacontinuousfullyfaithfullyintocondensed}
    The functor 
    \begin{align*}\underline{(-)}_{\kappa}\colon\kappacont&\to \Cond{\kappa}(\Set)\\ X &\mapsto \mathcal C(-,X)\colon {\Pro(\Fin)^{\operatorname{op}}_{\kappa}}\to \Set\end{align*} is fully faithful. 
    It admits a left adjoint which sends a $\kappa$-condensed set $X$ to the topological space \[X(*)_{\kappa}\coloneqq \colim{S\in \Pro(\Fin)_{\kappa/X}} S, \] where the colimit is computed in the category of all topological spaces with continuous maps.  
    For a $\kappa$-condensed set $X$, the space $X(*)_{\kappa}$ is $\kappa$-compactly generated with set of points $X(*)$. 
\end{proposition}
\begin{rem}\label{underlinefunctorontop}
    Denote by $\Top\subseteq \kappacont$ the wide subcategory on topological spaces and continuous maps. 
    As the right adjoint $(*)_{\kappa}$ takes values in $\kappa$-compactly generated spaces, 
    the adjunction $\underline{(-)}_{\kappa}\dashv (*)_{\kappa}$ induces an adjoint pair
    \[ \underline{(-)}_{\kappa}\colon \Top\leftrightarrows \Cond{\kappa}(\Set)\colon (*)_{\kappa}.\] 
\end{rem}
\begin{proof}This is essentially \cite[Proposition 1.7]{Scholzecondensed}. By \cref{explicitcharacterizationforcondensedset}, for all topological spaces $X$, $\underline{X}_{\kappa}\in\Cond{\kappa}(\Set)$. 
    As every $\kappa$-condensed set is a colimit of representables, for every $X\in\Cond{\kappa}(\Set)$ and  $Y\in\kappacont$,  
    \begin{align*}
    \Hom_{\Cond{\kappa}(\Set)}(X, \underline{Y}_{\kappa})&\cong \clim{S\in (\Pro(\Fin)_{\kappa/X})^{\operatorname{op}}}\Hom_{\Cond{\kappa}(\Set)}(\underline{S}_{\kappa}, \underline{Y}_{\kappa})\\
    &\cong  \clim{S\in (\Pro(\Fin)_{\kappa/X})^{\operatorname{op}}}\mathcal C(S,Y)\\&\cong \mathcal C(\colim{S\in \Pro(\Fin)_{\kappa/X}}S,Y)
    \end{align*} where the colimit $X(*)_{\kappa}\coloneqq \colim{S\in\Pro(\Fin)_{\kappa/X}} S$ is computed in the category $\Top$ of topological spaces with continuous maps.
    This space is $\kappa$-compactly generated by \cref{compacthausdorffquotientoftdch}, so in particular, $\mathcal C(X(*)_{\kappa},Y)\cong\mathcal C_{\kappa}(X(*)_{\kappa},Y)$. As all the above identifications are natural in $X$ and $Y$, this shows that \[X\mapsto \colim{S\in \Pro(\Fin)_{\kappa/X}}S\] is left adjoint to $\underline{(-)}_{\kappa}$. 
    \cref{globalsectionsleftadjointcondensed} and \cite[Proposition 5.5.6.28]{highertopostheory} imply that the global section functor \[\Cond{\kappa}(\Set)\to \Set, \, F\mapsto F(*)\] preserves colimits, whence \[X(*)_{\kappa}=\colim{S\in \Pro(\Fin)_{\kappa/X}} S \to X(*)\] is a bijection of sets. 
    It remains to show that $\underline{(-)}_{\kappa}$ is fully faithful. Fix a topological space $T$. 
    For a profinite set $S\in \Pro(\Fin)_{\kappa}$, $\underline{T}_{\kappa}(S)=\mathcal C(S,T)\subseteq \prod_{s\in S}T$ is the subset on elements $(t_s)_{s\in S}$ such that $S\to T, s\mapsto t_s$ is continuous. If $X\to \underline{T}_{\kappa}$ is a map of $\kappa$-condensed sets, for all $S\in \Pro(\Fin)_{\kappa}$ and $s\in S$, the inclusion $\{i_s\}\hookrightarrow S$ induces a commutative diagram   \begin{center}
    \begin{tikzcd}
        X(S)\arrow[r]\arrow[d,"i_s^*"] & T(S)\arrow[r,hookrightarrow] &  \prod_{s\in S}T \arrow[d,"\pi_s"]\\
        X(\{s\})\arrow[rr,"f(*)"] & & T.
    \end{tikzcd}
    \end{center} A map of $\kappa$-condensed sets $X\to\underline{T}_{\kappa}$ is therefore the same data as a map of sets $f(*)\colon X(*)\to T$ such that for all profinite sets $S\in\Pro(\Fin)_{\kappa}$, and all $x\in X(S)$, the composition \[S\to X(*)\to T, \, s\mapsto f(*)(i_s^*(x))\] is continuous. Since $X(S)=\Map_{\Cond{\kappa}(\Set)}(\underline{S}_{\kappa},X)$, this is equivalent to requiring that for all $S\in\Pro(\Fin)_{\kappa}$ and all maps $\underline{S}_{\kappa}\xrightarrow{t} X\in \Cond{\kappa}(\Set)$, the composition $S\xrightarrow{t(*)}X(*)\xrightarrow{f(*)} T$ is continuous.     
    In particular, if $X=\underline{Y}_{\kappa}$ for a topological space $Y$, then $\Hom_{\Cond{\kappa}(\Set)}(\underline{Y}_{\kappa}, \underline{T}_{\kappa})$ is the set of maps $Y\to T$ such that for all $S\in \Pro(\Fin)_{\kappa}$ and all \[t\in\mathcal C(S,Y)=\underline{Y}_{\kappa}(S)\cong \Hom_{\Cond{\kappa}(\Set)}(\underline{S}_{\kappa}, \underline{Y}_{\kappa}),\] the composition $S\xrightarrow{t}Y \to T$ is continuous.
    By \cref{compacthausdorffquotientoftdch}, these are precisely the $\kappa$-continuous maps, which shows that $\underline{(-)}_{\kappa}$ is fully faithful. 
\end{proof}

\begin{lemma}[{\cite[Proposition 2.1.8]{LucasMannthesis}, \cite[Proposition 2.5]{aoki2024semitopologicalktheorysolidification}}]\label{leftKanextensionsfullyfaithfulcondensed}

    Suppose $\mathcal C$ is a presentable category.     
    For uncountable cardinals $\lambda\geq \kappa$, the restriction
    \[r^{\lambda}_{\kappa}\colon  \Cond{\lambda}(\mathcal C)\to \Cond{\kappa}(\mathcal C)\] has a left adjoint. 

    Suppose in addition that one of the following holds: 
    \begin{romanenum}
        \item $\lambda$ and $\kappa$ are uncountable regular cardinals and $\kappa$-filtered colimits commute with $\Delta$-indexed limits in $\mathcal C$.
        \item $\lambda, \kappa$ are strong limit cardinals and $\cof(\kappa)$-filtered colimits commute with finite products in $\mathcal C$. 
    \end{romanenum}
    Then the left adjoint of restriction is fully faithful and given by left Kan extension along \[\Pro(\Fin)_{\kappa}^{\operatorname{op}}\subseteq \Pro(\Fin)_{\lambda}^{\operatorname{op}}.\]  
\end{lemma}
\begin{proof}
    By construction of the restriction, 
    the composition \[\Cond{\kappa}(\mathcal C)\subseteq \Fun(\Pro(\Fin)_{\kappa}^{\operatorname{op}}, \mathcal C)\xrightarrow{i_{\kappa,!}^{\lambda}}\Fun(\Pro(\Fin)_{\lambda}^{\operatorname{op}}, \mathcal C)\xrightarrow{L_{\cond{\lambda}}}\Cond{\lambda}(\mathcal C)\] of left Kan extension along $\Pro(\Fin)_{\kappa}^{\operatorname{op}}\subseteq \Pro(\Fin)_{\lambda}^{\operatorname{op}}$ followed by $\lambda$-condensed hypersheafification is left adjoint to the restriction. 
    The left Kan extension is fully faithful since $\Pro(\Fin)_{\kappa}^{\operatorname{op}}\subseteq \Pro(\Fin)_{\lambda}^{\operatorname{op}}$ is. 
    By \cite[Proposition 2.1.8]{LucasMannthesis} or \cref{kappacondensedkappaacc} in the first and by \cite[Proposition 2.5]{aoki2024semitopologicalktheorysolidification} in the second case, 
    left Kan extension $i^{\lambda}_{\kappa,!}$ restricts to a functor $\Cond{\kappa}(\mathcal C)\to \Cond{\lambda}(\mathcal C)$ which is hence left adjoint to the restriction.  
\end{proof}
\subsection{Other models for \texorpdfstring{$\kappa$}{k}-condensed animae}
\label{section:othermodelscondensedcat}
In this section, we model $\Cond{\kappa}(\an)$ as hypercomplete sheaves on other classes of topological spaces. This will be central for our comparison of condensed with sheaf cohomology.  
\begin{definition}\label{definitioncondensedtopology}Suppose $\kappa$ is an uncountable cardinal and $\mathbb T\subseteq \Top$ is a small full subcategory. Denote by $\underline{(-)}_{\kappa}\colon \Top\to \Cond{\kappa}(\an)$ the restriction of the functor from \cref{kappacontinuousfullyfaithfullyintocondensed} to the wide subcategory $\Top\subseteq \kappacont$.  
    A collection of maps $\{T^i\to T\}_{i\in I}\subseteq \mathbb T$ is a $\kappa$-\emph{condensed cover} if \[p\colon \bigsqcup_{i\in I} \underline{T^i}_{\kappa}\to \underline{T}_{\kappa}\] is an epimorphism in $\Cond{\kappa}(\Set)$. 
\end{definition}
The main result of this section is the following: 
\begin{proposition}\label{condensedonothercategories} Suppose $\kappa$ is an uncountable cardinal and $\mathbb T\subseteq \Top$ is a small full subcategory. 
    \begin{romanenum}
    \item If $\mathbb T$ contains all $\kappa$-light profinite sets, then the $\kappa$-condensed covers form a coverage (\cref{definitioncoverage}) on $\mathbb T$ and right Kan extension  yields an equivalence 
    \[ \Cond{\kappa}(\an)\cong  \hypershv_{\cond{\kappa}}(\mathbb T, \an)\] between hypercomplete, $\kappa$-condensed sheaves on $\mathbb T$ (\cref{sheavesquasitopology}) and $\kappa$-condensed animae. 

    \item If $\kappa$ is a strong limit cardinal and $\mathbb T$ contains all $\kappa$-light extremally disconnected compact Hausdorff spaces, then the $\kappa$-condensed covers constitute a coverage on $\mathbb T$ and right Kan extension yields equivalences
    \[ \Cond{\kappa}(\an)\cong \mathcal P_{\Sigma}(\edCH{\kappa})\cong \hypershv_{\cond{\kappa}}(\mathbb T, \an)\] between hypercomplete sheaves in the $\kappa$-condensed topology on $\mathbb T$ and $\kappa$-condensed animae. 
    \end{romanenum}
    In both cases, the hypersheafified Yoneda embedding \[\mathbb T\xhookrightarrow{Y} \mathcal P(\mathbb T)\xrightarrow{L} \hypershv_{\cond{\kappa}}(\mathbb T, \an)\] sends a topological space $X$ to $\mathcal C_{\kappa}(-,X)$. 
\end{proposition}Before giving the proof, we make some preliminary observations about $\kappa$-condensed covers.
\begin{lemma}\label{condensedepisurjective}
If $\{ X^i\to X\}_{i\in I}\subseteq \Top$ is a $\kappa$-condensed cover, the map $\sqcup_{i\in I}X^i\to X$ is surjective.  
\end{lemma} 
\begin{proof}
    Suppose $\{X^i\to X\}_{i\in I}\subseteq \Top$ is a $\kappa$-condensed cover and let $Y\coloneqq \sqcup_{i\in I}\underline{X^i}_{\kappa}$. 
    As epimorphisms in the 1-topos $\Cond{\kappa}(\Set)$ are effective and $(*)_{\kappa}\colon \Cond{\kappa}(\Set)\to \Top$ is a left adjoint, \[\underline{X}_{\kappa}(*)_{\kappa}=\Coeq(Y(*)_{\kappa}\times_{\underline{X}_{\kappa}(*)_{\kappa}}Y(*)_{\kappa}\rightrightarrows Y(*)_{\kappa})\] in $\Top$. As $\Top\to \Set$ reflects colimits, this implies that $\sqcup_{i\in I}(\underline{X^i}_{\kappa}(*)_{\kappa})=Y(*)_{\kappa}\to\underline{X}_{\kappa}(*)_{\kappa}$ is surjective, whence $\sqcup_{i\in I}X^i\to X$ is surjective by \cref{kappacontinuousfullyfaithfullyintocondensed}. 
\end{proof}
\begin{lemma}\label{characterizationcondensedcovers}
    \begin{romanenum}
    \item\label{firstcharacterizationcondensedcovers} A family $\{ X^i\to X\}_{i\in I}\subseteq \Top$ is a $\kappa$-condensed cover if and only if for all $K\to X\in\oc{{\Pro(\Fin)_{\kappa}}}{X}$, there exists a finite family $\{ K^j\to K\}_{j=1}^n\subseteq \Pro(\Fin)_{\kappa}$ such that $\sqcup_{j=1}^n K^j\to K$ is surjective, and for all $j\in J$ there is $i(j)\in I$ such that $K_{j}\to K\to X$ factors over $K_j\to X^{i(j)}\to X$.  
    \item If $\kappa$ is a strong limit cardinal, a family $\{ X^i\to X\}_{i\in I}\subseteq \Top$ is a $\kappa$-condensed cover if and only if for all $K\to X\in\oc{{\edCH{\kappa}}}{X}$, there exists a finite family $\{ K^j\to K\}_{j=1}^n\subseteq \edCH{\kappa}$ such that $\sqcup_{j=1}^n K^j\to K$ is surjective and for all $j\in J$ there is $i(j)\in I$ such that $K_{j}\to K\to X$ factors over $K_j\to X^{i(j)}\to X$. 
    \end{romanenum}
\end{lemma}

\begin{proof}
    We first show the following: 
    \begin{enumerate}\label{characterizationcondensedcoversonprofin}
    \item[(*)] Suppose $K\in\Pro(\Fin)_{\kappa}$. A map $\sqcup_{i\in I}X^i\to \underline{K}_{\kappa}\in\Cond{\kappa}(\Set)$ is an epimorphism if and only if there exists a finite family $\{ K^j\to K\}_{j=1}^n\subseteq \Pro(\Fin)_{\kappa}$ such that $\sqcup_{j=1}^n K^j\to K$ is surjective and for all $j\in J$ there is $i(j)\in I$ such that $\underline{K^{j}}_{\kappa}\to \underline{K}_{\kappa}$ factors over a map $\underline{K^j}_{\kappa}\to X^{i(j)}$.
    \end{enumerate}
    The if-statement holds since $\sqcup_{j=1}^n\underline{K^j}_{\kappa}\to \underline{K}_{\kappa}$ is an epimorphism of $\kappa$-condensed sets by definition of the condensed topology on $\Pro(\Fin)_{\kappa}$ (see also \cite[Propositions 6.2.3.20 and 7.2.1.14]{highertopostheory}). 
    Conversely, if $\sqcup_{i\in I}X^i\to \underline{K}_{\kappa}$ is an epimorphism of $\kappa$-condensed sets, then 
    \[\sqcup_{i\in I}\sqcup_{Y\in \oc{\Pro(\Fin)_{\kappa}}{X^i}}\underline{Y}_{\kappa}\to \sqcup_{i\in I}X^i \to \underline{K}_{\kappa}\] is an epimorphism. As the condensed topology on $\Pro(\Fin)_{\kappa}$ is finitary, there exist a finite subset $F\subseteq I$ and for all $f\in F$ a finite subset $P_f\subseteq \oc{{\Pro(\Fin)_{\kappa}}}{X^f}$ such that \[\sqcup_{f\in F}\sqcup_{Y\in {P_f}}\underline{Y}_{\kappa}\to \underline{K}_{\kappa}\] is an epimorphism of $\kappa$-condensed sets, cf.\ \cite[\href{https://stacks.math.columbia.edu/tag/0D06}{Lemma 0D06}]{stacks-project}.  
    By \cref{kappacontinuousfullyfaithfullyintocondensed}, the maps $\underline{Y}_{\kappa}\to \underline{K}_{\kappa}$ are induced by continuous maps $Y\to K$ and by \cref{condensedepisurjective}, \[\sqcup_{f\in F}\sqcup_{Y\in P_f}Y\to K\] is surjective. This shows (*). 

    We now prove $(i)$. 
    As colimits in the 1-topos $\Cond{\kappa}(\Set)$ are universal, epimorphisms are effective and every $\kappa$-condensed set is a colimit of representables $\underline{K}_{\kappa}, K\in\Pro(\Fin)_{\kappa}$, a map $A\to B\in\Cond{\kappa}(\Set)$ is an epimorphism if and only if for all $\underline{K}_{\kappa}\to B$, $K\in\Pro(\Fin)_{\kappa}$, 
    \[A\times_B \underline{K}_{\kappa}\to \underline{K}_{\kappa}\] is an epimorphism of $\kappa$-condensed sets. 
    In particular, a family of continuous maps \[\{ X^i\to X\}_{i\in I}\subseteq \Top\] is a $\kappa$-condensed cover if and only if for all $K\in\Pro(\Fin)_{\kappa}$ and $\underline{K}_{\kappa}\to \underline{X}_{\kappa}\in\Cond{\kappa}(\Set)$, \[ \sqcup_{i\in I}\underline{X^i\times_X K}_{\kappa}=\sqcup_{i\in I}\underline{X^i}_{\kappa}\times_{\underline{X}_{\kappa}}\underline{K}_{\kappa}\to \underline{K}_{\kappa}\] is an epimorphism of $\kappa$-condensed sets.
    By (*), this holds if and only if for all $\underline{K}_{\kappa}\to X$, $K\in\Pro(\Fin)_{\kappa}$, there exists a finite family $\{ K^j\to K\}_{j\in J}\subseteq \Pro(\Fin)_{\kappa}$ such that 
    \begin{romanenum}
        \item $\sqcup_{j\in J}K^j\to K$ is onto, and 
        \item for all $j\in J$ there is $i(j)\in I$ such that 
    $K^j\to K$ factors as $\underline{K^j}_{\kappa}\to \underline{X^{i(j)}\times_X K}_{\kappa}\to \underline{K}_{\kappa}$.\end{romanenum} 
    By \cref{kappacontinuousfullyfaithfullyintocondensed}, $(ii)$ holds if and only if for all $j\in J$, $K^j\to K$ factors over a $\kappa$-continuous (and hence continuous) map $K^j\to X^{i(j)}\times_{X}K$. 
    This proves the first statement. 

    Suppose now that $\kappa$ is a strong limit cardinal. 
        By \cref{kappalightforstronglimit}, every $\kappa$-light profinite space is a quotient of a $\kappa$-light extremally disconnected compact Hausdorff space. 
        For $X\in\Pro(\Fin)_{\kappa}$ choose a quotient map $Y\to X, X\in\edCH{\kappa}$. 
        As $Y\times_XY\in\Pro(\Fin)_{\kappa}$, there exists a quotient map $Z\to Y\times_X Y$ with $Z\in\edCH{\kappa}$. 
        By definition of the condensed topology on $\Pro(\Fin)_{\kappa}$, 
        \[\Coeq(\underline{Z}_{\kappa}\to \underline{Y\times_X Y}_{\kappa}\rightrightarrows \underline{Y}_{\kappa})=\underline{X}_{\kappa}\in\Cond{\kappa}(\Set).\] 
        This shows that $\underline{X}_{\kappa}, X\in\edCH{\kappa}$ generates $\Cond{\kappa}(\Set)$ under colimits. 
         As epimorphisms in $\Cond{\kappa}(\Set)$ are effective and colimits are universal, this implies that a map $A\to B\in\Cond{\kappa}(\Set)$ is an epimorphism if and only if for all $\underline{K}_{\kappa}\to B, K\in\edCH{\kappa}$, $A\times_B \underline{K}_{\kappa}\to \underline{K}_{\kappa}$ is an epimorphism in $\Cond{\kappa}(\Set)$. 
         The characterization of $\kappa$-condensed covers now follows from \cref{kappalightforstronglimit} and the above. 
\end{proof}

\begin{lemma}\label{kappacondensedcoversquasicovering}
    Suppose $\mathbb T\subseteq \Top$ is a small full subcategory.  
    \begin{romanenum}
    \item If $\mathbb T$ contains all $\kappa$-light profinite sets, then the $\kappa$-condensed covers form a coverage on $\mathbb T$. 
    \item If $\kappa$ is a strong limit cardinal and $\mathbb T$ contains all $\kappa$-light extremally disconnected compact Hausdorff spaces, then the $\kappa$-condensed covers form a coverage on $\mathbb T$.
    \end{romanenum}
    \end{lemma}
    \begin{proof}
        Suppose first that $\Pro(\Fin)_{\kappa}\subseteq \mathbb T$. 
        Suppose that $\{T^i\to T\}_{i\in I}$ is a condensed cover and $S\to T\in \mathbb T$ is a continuous map. 
        Since epimorphisms are effective and colimits are universal, epimorphisms in a 1-topos are stable under pullback. In particular, 
        \[\sqcup_{i\in I}\underline{T^i\times_T S}_{\kappa}\cong (\sqcup_{i\in I}\underline{T^i}_{\kappa})\times_{\underline{T}_{\kappa}}\underline{S}_{\kappa}\to \underline{S}_{\kappa}\] is an epimorphism of $\kappa$-condensed sets. 
        As for $i\in I$, \[ \bigsqcup_{{X\in {\oc{{\Pro(\Fin)_{\kappa}}}{T^i\times_T S}}}}\underline{X}_{\kappa}\to \underline{T^i\times_T S}_{\kappa}\] is an epimorphism of $\kappa$-condensed sets, this implies that 
        \[ \{X\to T^i\times_T S\to S\}_{X\in \oc{{\Pro(\Fin)_{\kappa}}}{T^i\times_T S}, i\in I}\subseteq \mathbb T\] is a $\kappa$-condensed cover of $S$, which shows that $\kappa$-condensed covers define a coverage on $\mathbb T$.
                
        Suppose now that $\kappa$ is a strong limit cardinal and $\edCH{\kappa}\subseteq \mathbb T$.  
        Suppose that $\{T_i\to T\}_{i\in I}$ is a $\kappa$-condensed cover and $S\to T\in \mathbb T$. 
        Above, we constructed a family of continuous maps $\{X_j\to S\}_{j\in J}$ such that
        \begin{romanenum}
            \item for all $j\in J$, $X_j\in\Pro(\Fin)_{\kappa}$, 
            \item $\sqcup_{j\in J}\underline{X_j}_{\kappa}\to \underline{S}_{\kappa}$ is an epimorphism of $\kappa$-condensed sets, and 
            \item for all $j\in J$, there exists $i(j)\in I$ such that $X_j\to S\to T$ factors as $X_j\to T_{i(j)}\to T$. 
        \end{romanenum}
        Since $\kappa$ is a strong limit cardinal, for all $j\in J$, $X_j$ is the quotient of a $\kappa$-light extremally disconnected compact Hausdorff space $\tilde{X_j}$ (\cref{kappalightforstronglimit}), and by definition of the condensed topology, the quotient map $\underline{\tilde{X}_j}_{\kappa}\to \underline{X_j}_{\kappa}$ is an epimorphism of $\kappa$-condensed sets for all $j\in J$. This implies that $\{\tilde{X}_j\to X_j\to S\}_{j\in J}\subseteq \mathbb T$ is a $\kappa$-condensed cover, which shows that $\kappa$-condensed covers define a coverage on $\mathbb T$. 
    \end{proof} 
\begin{cor}\label{kappacontinuousfunctionscondensedsheaf}
        Suppose that $\mathbb T\subseteq \Top$ is a small full subcategory and $\kappa$ is an uncountable cardinal such that the $\kappa$-condensed coverings form a coverage on $\mathbb T$. 
        For a topological space $S$, \[\mathcal C_{\kappa}(-,S)\colon \mathbb T^{\operatorname{op}}\to \Set\subseteq \an\] defines is a $\kappa$-condensed sheaf, i.e.\ $\mathcal C_{\kappa}(-,S)\in \hypershv_{\cond{\kappa}}(\mathbb T)$.  
\end{cor}
\begin{proof}
        Suppose $\{X^i\to X\}_{i\in I}\in \mathbb T$ is a $\kappa$-condensed cover. As $\sqcup_{i\in I}X_i\to X$ is surjective (\cref{condensedepisurjective}), for every topological space $S$, \[\Map_{\Set}(X,S)\cong \Eq(\prod_{i\in I}\Map_{\Set}(X^i,S)\rightrightarrows \prod_{i,j\in I}\Map_{\Set}(X^i\times_X X^j,S)).\] 
        We claim that this isomorphism restricts to an isomorphism      
    \[\mathcal C_{\kappa}(X,S)\cong \Eq(\prod_{i\in I}\mathcal C_{\kappa}(X^i,S)\rightrightarrows \prod_{i,j\in I}\mathcal C_{\kappa}(X^i\times_X X^j,S)), \] i.e.\ that a map of sets $X\to S$ is $\kappa$-continuous if and only if the composite $\sqcup_{i\in I}X^i\to X\to S$ is $\kappa$-continuous. The only-if statement holds by continuity of $\sqcup_{i\in I}X^i\to X$. 
    Suppose now that $X\to S$ is a map of sets such that $\sqcup_{i\in I}X^i\to X\to S$ is $\kappa$-continuous. 
    By \cref{compacthausdorffquotientoftdch}, it suffices to show that for all $K\to X\in \oc{{\Pro(\Fin)_{\kappa}}}{X}$, the composite $K\to X\to S$ is continuous. 
    Fix $K\to X\in \oc{{\Pro(\Fin)_{\kappa}}}{X}$. 
    By \cref{characterizationcondensedcovers}, there exists a $\kappa$-light profinite space $\tilde K$ with a continuous surjection $\tilde K\to K$ and a continuous map $\tilde K\to \sqcup_{i\in I}X^i$ such that $\tilde K\to K\to X$ factors as $\tilde K\to \sqcup_{i\in I}X^i\to X$. 
    In particular, $\tilde K\to K\to X$ is continuous. 
    As $\tilde K$ and $K$ are $\kappa$-light compact Hausdorff, $K\to \tilde K$ is a quotient map, whence $K\to X\to S$ is continuous. 
\end{proof}
\cref{condensedonothercategories} is now a straightforward consequence:  
\begin{proof}[Proof of \cref{condensedonothercategories}]\label{proofofcondensedonothercats}
   By \cref{kappacondensedcoversquasicovering}, the $\kappa$-condensed covers form a coverage $\tau_{\cond{\kappa}}$ on $\mathbb T$ in both cases. 
   Suppose first that $\Pro(\Fin)_{\kappa}\subseteq \mathbb T$ and denote by $\tau_{\condo}$ the coverage on $\Pro(\Fin)_{\kappa}$ given be finite, jointly surjective covering families. 

    By \cite[Propositions 6.2.3.20 and 7.2.1.14]{highertopostheory}, the inclusion $(\Pro(\Fin)_{\kappa}, \tau_{\condo})\to (\mathbb T, \tau_{\cond{\kappa}})$ is a continuous functor, and for $T\in \mathbb T$, \[\{ X\to T\}_{X\to T\in\Pro(\Fin)_{\kappa/T}}\subseteq \mathbb T\] is a $\kappa$-condensed cover in $\mathbb T$. 
    By \cref{characterizationcondensedcovers}, $(\Pro(\Fin)_{\kappa}, \tau_{\condo})\to (\mathbb T, \tau_{\cond{\kappa}})$ has the covering lifting property (\cref{definitioncoveringliftingproperty}). 
    It now follows from 
    \cref{criterionrightkanextensionisocovers} that 
    \[ \Cond{\kappa}(\an)\cong \hypershv_{\cond{\kappa}}(\mathbb T)\] via right Kan extension and restriction. 
    By \cref{kappacontinuousfunctionscondensedsheaf}, for $T\in \mathbb T$: $\mathcal C_{\kappa}(-,T)\colon \mathbb T^{\operatorname{op}}\to \Set$ is a $\kappa$-condensed sheaf. 
    The inclusion $\mathcal C(-,T)\subseteq \mathcal C_{\kappa}(-,T)$ therefore induces a monomorphism \[L_{\cond{\kappa}}(\mathcal C(-,T))\to \mathcal C_{\kappa}(-,T)\in \hypershv_{\cond{\kappa}}(\mathbb T),\] where $L_{\cond{\kappa}}\colon \Fun(\mathbb T^{\operatorname{op}},\an)\to\hypershv_{\cond{\kappa}}(\mathbb T)$ denotes condensed hypersheafification which is left-exact by \cite[Propositions 6.2.2.7, 6.2.1.1 and 6.5.1.16]{highertopostheory}.  
    As the restriction \[\hypershv_{\cond{\kappa}}(\mathbb T)\to \Cond{\kappa}(\an)\] is an equivalence and for $X\in\Pro(\Fin)_{\kappa}$, 
    \[ \mathcal C(X,T)=\mathcal C_{\kappa}(X,T), \] it follows that $L_{\cond{\kappa}}(\mathcal C(-,T))\cong \mathcal C_{\kappa}(-,T)$. 

    Suppose now that $\edCH{\kappa}\subseteq \mathbb T$ and equip both sides with the $\kappa$-condensed coverages $\tau_{\cond{\kappa}}$. The inclusion $(\edCH{\kappa}, \tau_{\cond{\kappa}})\subseteq (\mathbb T, \tau_{\cond{\kappa}})$ is continuous and by \cref{characterizationcondensedcovers}, it has the covering lifting property. 
    Since $\kappa$ is a strong limit cardinal, \cref{kappalightforstronglimit} implies that for $T\in\mathbb T$, \[\{ X\to T\}_{X\in\edCH{\kappa}/T}\subseteq \mathbb T\] is a $\kappa$-condensed cover. 
    It now follows from \cref{criterionrightkanextensionisocovers} that right Kan extension induces an equivalence 
    \[ \hypershv_{\cond{\kappa}}(\edCH{\kappa})\cong \hypershv_{\cond{\kappa}}(\mathbb T).\] 
    By setting $\mathbb T=\Pro(\Fin)_{\kappa}$, we obtain that \[ \hypershv_{\cond{\kappa}}(\edCH{\kappa}^{\operatorname{op}})\cong \hypershv_{\cond{\kappa}}(\Pro(\Fin)_{\kappa})\cong \Cond{\kappa}(\an), \] where the right equivalence is the first statement of the proposition for $\mathbb T=\Pro(\Fin)_{\kappa}$. 

    It remains to show that $\hypershv_{\cond{\kappa}}(\edCH{\kappa})\cong \mathcal P_{\Sigma}(\edCH{\kappa})$. As for $X^0,X^1\in \edCH{\kappa}$, \[\{X^i\to X^0\sqcup X^1\}_{i=0,1}\] is a $\kappa$-condensed cover, and the empty cover of the empty set is a $\kappa$-condensed cover, \[{\Shv}_{\cond{\kappa}}(\edCH{\kappa})\subseteq \mathcal P_{\Sigma}(\edCH{\kappa}).\] 
    Suppose that $\{ X^i\to X\}_{i\in I}\subseteq \edCH{\kappa}$ is such that $\sqcup_{i\in I}X^i\to X$ is an epimorphism of $\kappa$-condensed sets. 
    Since $\underline{X}_{\kappa}$ is representable and the condensed topology is finitary, there exists a finite subset $F\subseteq I$ such that $\sqcup_{i\in F}\underline{X^i}_{\kappa}\to \underline{X}_{\kappa}$ is an epimorphism of condensed sets (\cite[\href{https://stacks.math.columbia.edu/tag/0D06}{Lemma 0D06}]{stacks-project}). 
    By \cref{condensedepisurjective}, \[q\colon \sqcup_{i\in F}X^i\to X\] is surjective, and hence admits a section since $X$ is extremally disconnected. This implies that every presheaf \[F\colon \edCH{\kappa}^{\operatorname{op}}\to \an\] is $\{ \sqcup_{i\in F}X^i\to X\}$-descendable, and hence every additive sheaf $F\in\mathcal P_{\Sigma}(\edCH{\kappa})$ is $\{ X^i\to X\}_{i\in I}$-descendable, which shows that \[{\Shv}_{\cond{\kappa}}(\edCH{\kappa})=\mathcal P_{\Sigma}(\edCH{\kappa}).\] 
    By \cite[Lemma 6.5.1.2]{highertopostheory}, the homotopy group functors $\pi_n\colon\mathcal P(\mathcal C)\to\tau_{\leq 0}\mathcal P(\mathcal C)$ restrict to functors \[\mathcal P_{\Sigma}(\edCH{\kappa})\to \tau_{\leq 0}\mathcal P_{\Sigma}(\edCH{\kappa}), \] which implies that $\mathcal P_{\Sigma}(\edCH{\kappa})\cong \hypershv_{\cond{\kappa}}(\edCH{\kappa})$ is hypercomplete.
    The same argument as above shows that for $T\in\mathbb T$, $L_{\cond{\kappa}}(\mathcal C(-,T))\cong \mathcal C_{\kappa}(-,T)$. 
\end{proof}
In particular, we obtain the following models for $\Cond{\kappa}(\an)$:  
\begin{cor}[{\cite[Proposition 2.3, Proposition 2.7]{Scholzecondensed}}]\label{condensedoncompactextremallydisconnected}

    Suppose $\kappa$ is an uncountable cardinal.
    \begin{romanenum} 
    \item The finite, jointly surjective covers form a coverage $\tau$ on the category $\CH_{\kappa}$ of $\kappa$-light compact Hausdorff spaces which generates (\cref{definitiongrothendiecktopologyviacovers}) the $\kappa$-condensed topology on $\CH_{\kappa}$. In particular, restriction along $\Pro(\Fin)_{\kappa}\subseteq \CH_{\kappa}$ is an equivalence \[ \hypershv_{\tau}(\CH_{\kappa}, \an)\cong \Cond{\kappa}(\an).\]  

    \item The restriction $\Cond{\kappa}(\an)\to \mathcal P_{\Sigma}(\edCH{\kappa})$ has a fully faithful right adjoint given by right Kan extension. 
    If $\kappa$ is a strong limit cardinal, then this is an equivalence.
    \end{romanenum} 
\end{cor}
\begin{proof}
    We first prove the first statement.  
    If $\{ X_i\to X\}_{i=1}^n\subseteq \CH_{\kappa}$ is such that $\sqcup_{i=1}^n X_i\to X$ is onto, for all $Y\to X\in\CH_{\kappa}$: $\sqcup_{i=1}^nY\times_{X}X_i\to X$ is onto. Since $X$ is Hausdorff, $Y\times_{X}X_i$ is a closed subspace of $Y\times X$ and in particular $\kappa$-light compact Hausdorff. 
    This shows that finite, jointly surjective covers form a coverage on $\CH_{\kappa}$.
    By \cref{compacthausdorffquotientoftdch}, every $\kappa$-light compact Hausdorff space is the quotient of a $\kappa$-light profinite set. Together with \cref{characterizationcondensedcovers}, this implies that a family $\{ X^i\to X\}_{i\in I}\subseteq \CH_{\kappa}$ is a $\kappa$-condensed cover if and only if there exists a finite subset $F\subseteq I$ such that $\sqcup_{i\in F}X^i\to X$ is onto.
    This shows that $\kappa$-condensed covers and finite, jointly surjective maps generate the same Grothendieck topology on $\CH_{\kappa}$. 
    \cref{condensedonothercategories} now implies that 
    \[ \hypershv_{\tau}(\CH_{\kappa}, \an)\cong \Cond{\kappa}(\an)\] via restriction.   

    Equip $\edCH{\kappa}$ with the coverage given by covers of the form $\{ X_j\to \sqcup_{j=1}^nX^j\}_{j=1, \ldots, n}$ and the empty cover of the empty set, and $\Pro(\Fin)_{\kappa}$ with the coverage defined by finite, jointly surjective covers. We claim that $\edCH{\kappa}\subseteq \Pro(\Fin)_{\kappa}$ has the covering lifting property. Suppose that \[\{ X_i\to X\}_{i=1}^n\subseteq \Pro(\Fin)_{\kappa}\] is such that $q\colon \sqcup_{i=1}^n X_i\to X$ is onto and $X\in\edCH{\kappa}$. 
    Since $\sqcup_{i=1}^n X_i$ is compact Hausdorff, the map $q$ has a section $s\colon X\to \sqcup_{i=1}^n X_i$. As $X=\sqcup_{i=1}^n s^{-1}(X_i)$, \[\{ s^{-1}(X_i)\to X\}_{i=1}^n\subseteq \edCH{\kappa}\] is a cover for $X$, which shows that $\edCH{\kappa}\subseteq \Pro(\Fin)_{\kappa}$ has the covering lifting property. 
    By \cite[Lemma 6.5.1.2]{highertopostheory}, the homotopy group functors $\pi_n\colon\mathcal P(\edCH{\kappa})\to\tau_{\leq 0}\mathcal P(\edCH{\kappa})$ restrict to functors \[\mathcal P_{\Sigma}(\edCH{\kappa})\to \tau_{\leq 0}\mathcal P_{\Sigma}(\edCH{\kappa}), \] which implies that $\mathcal P_{\Sigma}(\edCH{\kappa})$ is hypercomplete. It now follows from \cref{rightkanextensionsheaves} that right Kan extension $\mathcal P(\edCH{\kappa})\to\mathcal P(\Pro(\Fin)_{\kappa})$ restricts to a functor 
    $\mathcal P_{\Sigma}(\edCH{\kappa})\to \Cond{\kappa}(\an)$ which is fully faithful since right Kan extension is.  
    By \cref{condensedonothercategories}, it is an equivalence if $\kappa$ is a strong limit cardinal. 
\end{proof}

Recall that an object $c$ in a category $\mathcal C$ is compact if $\Map_{\mathcal C}(c,-)\colon \mathcal C\to\an$ preserves filtered colimits and compact projective if 
$\Map_{\mathcal C}(c,-)$ preserves sifted colimits. 
\cref{condensedoncompactextremallydisconnected} implies that the following families of \textit{representables} are compact respectively projective in $\Cond{\kappa}(\an)$ and $\Cond{\kappa}(\Set)$: 
\begin{cor}\label{compactprojectivegeneratorsforcondensedsets}
    Suppose $\kappa$ is an uncountable cardinal. 
    \begin{romanenum}
    \item For a $\kappa$-light compact Hausdorff space $X\in\CH_{\kappa}$, $\underline{X}_{\kappa}$ is compact in $\tau_{\leq n}\Cond{\kappa}(\an)$ for all $n\in\mathbb N_0$ and in particular in $\Cond{\kappa}(\Set)=\tau_{\leq 0}\Cond{\kappa}(\an)$. 
    \item For a $\kappa$-light extremally disconnected compact Hausdorff space $X\in\edCH{\kappa}$, $\underline{X}_{\kappa}$ is compact projective in $\Cond{\kappa}(\an)$. If $\kappa$ is a strong limit cardinal, then \[\underline{X}_{\kappa},\,  X\in\edCH{\kappa}\] generate $\Cond{\kappa}(\an)$ under colimits. 
    \end{romanenum} 
\end{cor}

\begin{proof}
    By \cref{condensedoncompactextremallydisconnected}, the $\kappa$-condensed topology on $\CH_{\kappa}$ is finitary (i.e.\ every cover can be refined by a finite cover). Hence the representable sheaves $\underline{X}_{\kappa},X\in\CH_{\kappa}$ in the associated hypercomplete topos $\Cond{\kappa}(\an)$ are coherent by \cite[Corollary A.2.3.2, Proposition A.2.2.2]{SAG}. By \cite[Corollary A.2.3.2]{SAG}, this implies that for all $X\in\CH_{\kappa}$, $\underline{X}_{\kappa}$ is compact in $\tau_{\leq n}\Cond{\kappa}(\an)$ for all $n\in\mathbb N_0$. 
    
    By \cref{condensedoncompactextremallydisconnected}, the restriction $i^*\colon \Cond{\kappa}(\an)\to \mathcal P_{\Sigma}(\edCH{\kappa})$ has a right adjoint given by right Kan extension. In particular, $i^*$ preserves colimits. Since colimits in $\an$ are universal, sifted colimits in $\mathcal P_{\Sigma}(\edCH{\kappa})$ are computed pointwise, i.e.\ for representables $h_X\in\mathcal P_{\Sigma}(\edCH{\kappa}), X\in\edCH{\kappa}$, \[\Map_{\mathcal P_{\Sigma}(\edCH{\kappa})}(h_X,-)\colon \mathcal P_{\Sigma}(\edCH{\kappa})\to \an\] preserves sifted colimits, which implies that \[\Map_{\Cond{\kappa}(\an)}(\underline{X}_{\kappa},-)\cong \Map_{\mathcal P_{\Sigma}(\edCH{\kappa})}(h_X,-)\circ i^*\] does as well.      
    If $\kappa$ is a strong limit cardinal, then \[\Cond{\kappa}(\an)\cong \mathcal P_{\Sigma}(\edCH{\kappa})\] and hence the represented condensed sets $\underline{X}_{\kappa}$, for $X\in \edCH{\kappa}$ generate $\Cond{\kappa}(\an)$ under colimits. 
    \end{proof}

    \begin{cor}\label{cantorsetgenerates}The Cantor set $\underline{\{0,1\}^{\mathbb N}}_{\kappa}$ generates $\Cond{\light}(\an)$ under colimits. 
    \end{cor}
    \begin{proof}
        By the Hausdorff-Alexandroff theorem, every non-empty metrizable compact Hausdorff space is the quotient of the Cantor set $\{0,1\}^{\mathbb N}$. 
        It follows from \cite[Proposition 20.4.5.1]{SAG} that $\underline{\{0,1\}}^{\mathbb N}_{\kappa}$ generates $\Cond{\light}(\an)$ under small colimits. \end{proof}
\subsection{Condensed categories}\label{section:condensedcategories}
We now review condensed categories and record basic categorical properties thereof which follow from our  discussion in the first chapter and \cref{section:propertiesweightsection}.
\begin{lemma}\label{profiniteisexplicitcoveringsite}
$(\Pro(\Fin), \text{surjective maps})$ is an explicit covering site which is $\kappa$-hyperaccessible for every uncountable regular cardinal $\kappa$. 
\end{lemma}
\begin{proof}
    As the forget functor $\Pro(\Fin)\subseteq \CH\to \Set$ reflects limits and finite colimits, it is straightforward to check that 
$(\Pro(\Fin), \text{surjective maps})$ is an explicit covering site. 
The category $\Pro(\Fin)$ is $\kappa$-coaccessible for all uncountable regular cardinals $\kappa$ by \cite[Lemma 5.4.2.10, Remark 5.4.2.9]{highertopostheory}. 
It now follows from \cite[Lemma 2.1.7]{LucasMannthesis} that 
$(\Pro(\Fin), \text{surjective maps})$ is a $\kappa$-hyperaccessible explicit covering site for all uncountable regular cardinals $\kappa$. 
\end{proof}
\begin{definition}\label{definitioncondensed}
    Denote by $S\subseteq \Pro(\Fin)$ the wide subcategory on surjections. 
    For a presentable category $\mathcal C$ denote by 
    \[ \Cond{}(\mathcal C)\coloneqq \hypershvacc_{S}(\Pro(\Fin), \mathcal C)\] the category of hypercomplete, accessible $S$-sheaves valued in $\mathcal C$ (\cref{accessiblecoveringsitedefinition}). 
\end{definition}
This is equivalent to the definition of \cite{Scholzecondensed}, see \cref{condensedascolimit} below. 

More concretely, a functor $F\colon \Pro(\Fin)^{\operatorname{op}}\to\mathcal C$ is a condensed object in $\mathcal C$ if and only if the following are satisfied:
\begin{romanenum}
    \item $F$ is accessible,
    \item $F$ preserves finite products, and 
    \item for every condensed hypercover $U_*\colon \Delta_{s,+}^{\operatorname{op}}\to\Pro(\Fin)$ (\cref{Definitioncondensedhypercover}), 
\[F(U_{-1})=\clim{\Delta_{s}}F(U_*).\] 
\end{romanenum}
\begin{ex}\label{sheafconditioncondensedsetexplicit}
A functor $F\colon \Pro(\Fin)^{\operatorname{op}}\to \Set$ is a condensed set if and only if 
\begin{romanenum}
    \item $F$ is accessible, 
    \item $F$ preserves finite products, and 
    \item for all continuous surjections $q\colon X\to Y\in\Pro(\Fin)_{\kappa}$, \[F(Y)\cong \Eq(F(X)\rightrightarrows F(X\times_Y X)).\]
\end{romanenum}
This holds since $\Set\subseteq \an$ reflects filtered colimits and small limits and $\Delta_{s,\leq 1}\subseteq \Delta_{s}$ is left 1-cofinal, see \cite[Example 6.14]{du2023reshapinglimitdiagramscofinality}. 
\end{ex}
\begin{rem}[Pyknotic versus condensed]\label{pyknoticversuscondensed}
If $\kappa$ is an uncountable, small strong limit regular cardinal and $\mathcal U_{\kappa}\subseteq \mathcal U_0$ is the universe associated to ${\kappa}$, i.e. the set of $\kappa$-small sets, then $\Cond{\kappa}(\an)$ is the category of pyknotic animae corresponding to $\mathcal U_{\kappa}\in\mathcal U_0$.
As the inclusion $i\colon \an_{\mathcal U_0}\subseteq \an_{\mathcal U_1}$ from $\mathcal U_0$-small animae into $\mathcal U_1$-small animae preserves $\mathcal U_0$-small limits, the functor \[i_*\colon \Fun(\Pro(\Fin)_{\mathcal U_0}^{\operatorname{op}},\an_{\mathcal U_0})\hookrightarrow \Fun(\Pro(\Fin)_{\mathcal U_0}^{\operatorname{op}},\an_{\mathcal U_1})\] restricts to a fully faithful functor $\Cond{}(\an)_{\mathcal U_0}\hookrightarrow \operatorname{Pyk}(\an)_{\mathcal U_0}^{\mathcal U_1}$ from condensed anima (defined in $\mathcal U_0$) to pyknotic animae (defined with respect to $\mathcal U_0\in \mathcal U_1$). Here $\Pro(\Fin)_{\mathcal U_0}$ denotes the category of $\mathcal U_0$-small profinite sets. By \cref{kappalightforstronglimit}, these are precisely the $\mathcal U_0$-light profinite sets. 
As $(\Pro(\Fin)_{\mathcal U_2}, \text{ epimorphisms})$ is $\mathcal U_0$- and $\mathcal U_1$-hyperaccessible (\cref{profiniteisexplicitcoveringsite}), left Kan extension along 
$\Pro(\Fin)_{\mathcal U_0}\subseteq \Pro(\Fin)_{\mathcal U_1}$ restricts to a functor \[\operatorname{Pyk}(\an)_{\mathcal U_0}^{\mathcal U_2}\hookrightarrow \operatorname{Pyk}(\an)_{\mathcal U_1}^{\mathcal U_2}\] (\cref{kappacondensedkappaacc} in $\mathcal U_1$). By definition of condensed anima (in $\mathcal U_1$), this factors over \[\Cond{}(\an)_{\mathcal U_1}\subseteq \operatorname{Pyk}(\an)_{\mathcal U_1},\] so we obtain fully faithful functors 
\[ \Cond{\kappa}(\an)_{\mathcal U_0}\subseteq \operatorname{Pyk}(\an)_{\mathcal U_0}^{\mathcal U_1}\subseteq \Cond{\kappa}(\an)_{\mathcal U_1}\subseteq \operatorname{Pyk}(\an)_{\mathcal U_1}^{\mathcal U_2}.\] 
\end{rem}
As $\Pro(\Fin)$ is not an essentially small category and the condensed topology is subcanonical,  $\Cond{}(\an)$ is not an accessible category by \cref{accessiblesheavesnottopos}, and in particular not a topos (\cite[Proposition 6.1.3.19]{highertopostheory}).
We recall here from the first chapter that $\Cond{}(\an)$ however resembles a topos in many ways. \cref{Giraudsaxiomsaccessiblesheaves} implies:
    \begin{cor}[Giraud's axioms]\label{Girauds axioms}
    Condensed animae are a locally small category and satisfy all of Giraud's axioms except accessibility, i.e.\ 
\begin{romanenum}
    \item All small colimits and all small limits exist.
    \item Small colimits are universal.
    \item Small coproducts are disjoint.
    \item Every groupoid object of $\Cond{}(\an)$ is effective. 
\end{romanenum}
\end{cor}
  \cref{closednessmonoidalstructureaccessiblesheaves} implies:  
\begin{cor}[{\cite[Proposition 2.1.11]{LucasMannthesis}}]\label{closedmonoidalstructurecondensedcategories}
    If $\mathcal C$ is a presentably symmetric monoidal category, $\Cond{}(\mathcal C)$ inherits the structure of a closed monoidal category such that for all regular cardinals $\kappa$ such that $\mathcal C$ is $\kappa$-accessible, $\Cond{\kappa}(\mathcal C)\subseteq \Cond{}(\mathcal C)$ is a symmetric monoidal subcategory. 
\end{cor}
We have also shown that the category of spectrum objects in $\Cond{}(\an)$ behaves like the stabilization of a presentable category in many
aspects:
\begin{romanenum}\label{condensedspectrabasicproperties}
    \item $\stab{\Cond{}(\an)}\cong \Cond{}(\Sp)$ is a stable category
(\cref{spectrumobjectsinbigcatsstable}, \cref{stabilizationisspectrumvaluedsheaves}).
    \item The functor $\Omega^{\infty}\colon \Cond{}(\Sp)\to \Cond{}(\an)$
admits a left adjoint $\Sigma^{\infty}_{+}\colon \Cond{}(\an)\to
\Cond{}(\Sp)$, which factors as a composite of left adjoints
$\Cond{}(\an)\to \CGrp(\Cond{}(\an))\to \Cond{}(\Sp)$
(\cref{existencesuspension}).
    \item The essential image of $\CGrp(\Cond{}(\an))\to \Cond{}(\Sp)$ is
the connective part of a $t$-structure on $\Cond{}(\Sp)$ with
$\Cond{}(\Sp)^{\heart}\cong \Cond{}(\Ab)$.
(\cref{tstructurespectrumobjects})
    \item There exists an essentially unique cocontinuous symmetric
monoidal structure on $\Cond{}(\Sp)$ with a symmetric monoidal enhancement of $\Sigma^{\infty}_{+}$ (\cref{monoidalstructurespectrumobjectscocontinuous}).
This symmetric monoidal structure is closed (\cref{symmetricmonoidalstructureaccessiblesheavesofspectraclosed}).
    \item The symmetric monoidal structure on $\Cond{}(\Sp)$ is compatible
with the $t$-structure, i.e.\ $\Cond{}(\Sp)_{\geq 0}\subseteq
\Cond{}(\Sp)$ is a symmetric monoidal subcategory
(\cref{monoidalstructurecompatibletstructure}).
\end{romanenum}

\subsubsection{Condensed and $\kappa$-condensed categories}
We now recall that for a presentable category $\mathcal C$, $\Cond{}(\mathcal C)$ can be exhausted by $\kappa$-condensed categories $\Cond{\kappa}(\mathcal C)$, as $\kappa$ ranges over a suitable family of uncountable cardinals. We will apply this to deduce statements on condensed categories from the corresponding $\kappa$-condensed statements for sufficiently large $\kappa$. 
The identification from \cref{sheafconditionexplicitkappacondensed}, yields for every uncountable cardinal $\kappa$ a restriction functor \[\Cond{}(\mathcal C)\to \hypershv_{\condo}(\Pro(\Fin)_{\kappa}^{\operatorname{op}}, \mathcal C)\cong \Cond{\kappa}(\mathcal C).\] 
\begin{lemma}\label{leftadjointrestrictions}
Suppose $\mathcal C$ is a presentable category. 
For every uncountable cardinal $\kappa$, the restriction 
$\Cond{}(\mathcal  C)\to \Cond{\kappa}(\mathcal C)$ has a left adjoint. 
\end{lemma}
\begin{proof}
Chose a regular cardinal $\mu$ such that $\mathcal C$ is $\mu$-accessible. 
Suppose $\kappa$ is an uncountable cardinal and choose a regular cardinal $\lambda\geq \kappa, \lambda\gg\mu$. 
Then the restriction $r^{\kappa}\colon \Cond{}(\mathcal C)\to\Cond{\kappa}(\mathcal C)$ factors as \[\Cond{}(\mathcal C)\xrightarrow{r^{\lambda}}\Cond{\lambda}(\mathcal C)\xrightarrow{r^{\lambda}_{\kappa}}\Cond{\kappa}(\mathcal C).\]
By \cref{profincoaccessible} and \cref{kappacondensedkappaacc}, $r^{\lambda}$ has a left adjoint.
The functor $r^{\lambda}_{\kappa}$ has a left adjoint by \cref{leftKanextensionsfullyfaithfulcondensed}. 
\end{proof}
Denote by $\Lambda$ the poset of all small uncountable cardinals and by $\Lambda^{\triangleright}\coloneqq \Lambda\cup \{\infty\}$ its cone. 
The \textit{restrictions} yield a functor 
\begin{align*}
(\Lambda^{\triangleright})^{\operatorname{op}}& \to \vlCat\\ 
\lambda&\mapsto \begin{cases}\Cond{\lambda}(\mathcal C) & \lambda\neq \infty\\
\Cond{}(\mathcal C) & \lambda=\infty.\end{cases} 
\end{align*}
By \cref{leftadjointrestrictions} and \cref{leftKanextensionsfullyfaithfulcondensed}, this enhances to a functor \[(\Lambda^{\triangleright})^{\operatorname{op}}\to \vlCat^R.\]  
Denote by $\Cond{*}(\mathcal C)\colon\Lambda^{\triangleright}\to \vlCat^L$ its opposite. (Recall that $\vlCat^L\cong (\vlCat^R)^{\operatorname{op}}$, see e.g.\ \cite[Theorem B]{HaugsengHebestreitLinskensNuiten}.)
\cref{kappacondensedkappaacc} implies: 
\begin{cor}\label{condensedasfilteredcolimitoverregular}
    Suppose $\mathcal C$ is a presentable category. 
    Then \[\Cond{*}(\mathcal C)\colon\Lambda^{\triangleright}\to \vlCat^L\to \vlCat\] is a colimit diagram. 
\end{cor}
\begin{proof}
    Choose a regular cardinal $\mu$ such that $\mathcal C$ is $\mu$-accessible. Since the uncountable regular cardinals $\kappa\geq \mu$ form a cofinal subcategory of $\Lambda$, the statement follows from \cref{accessiblesheavesisbigtopos} and \cref{profincoaccessible}.  
\end{proof}
\begin{cor}\label{condensedascolimit} Suppose that $\mathcal C$ is a presentable category. 
    Denote by $\Lambda^{r}, \Lambda^s$ the large posets of small regular/strong limit cardinals, respectively. 
    Then \[ \colim{\kappa\in\Lambda^r}\Cond{\kappa}(\mathcal C)\cong \colim{\kappa\in\Lambda^s}\Cond{\kappa}(\mathcal C)\cong \Cond{}(\mathcal C)\] where the colimits are computed in $\vlCat$ and the transition maps are the left adjoints of the restriction.
\end{cor}
\begin{proof}This follows from \cref{condensedasfilteredcolimitoverregular} as the regular/strong limit cardinals form a cofinal system in the poset of all small cardinals. 
\end{proof}
\begin{rem}
\cite{Scholzecondensed} defined $\Cond{}(\mathcal C)$ as $\colim{\kappa\in\Lambda^s}\Cond{\kappa}(\mathcal C)$. 
The identification \[\Cond{}(\mathcal C)\cong \colim{\kappa\in\Lambda^s}\Cond{\kappa}(\mathcal C)\] is instructive as for all strong limit cardinals $\kappa$, $\Cond{\kappa}(\mathcal C)\cong\mathcal P_{\Sigma}(\edCH{\kappa})$ is a particularly nice category. For example, $\Cond{\kappa}(\Ab)$ has enough compact projectives for all strong limit cardinals $\kappa$, which implies that the same holds for $\Cond{}(\Ab)$, see \cref{enoughprojectives}.  
\end{rem}
\begin{cor}\label{condensedstronglimitfullyfaithful}
    Suppose that $\mathcal C$ is a presentable category such that filtered colimits commute with finite products in $\mathcal C$ (e.g.\ $\mathcal C$ is compactly generated or stable). 
    Then the left adjoint \[\Cond{\kappa}(\mathcal C)\to \Cond{}(\mathcal C)\] of restriction is fully faithful for all strong limit cardinals $\kappa$.
\end{cor}
\begin{proof}
    By \cref{leftKanextensionsfullyfaithfulcondensed}, for all strong limit cardinals $\kappa\leq \lambda$, the left adjoint $\Cond{\kappa}(\mathcal C)\to \Cond{\lambda}(\mathcal C)$ of restriction is fully faithful. 
    \cref{filteredcolimitsofcategoriesappendix} now implies that \[\Cond{\kappa}(\mathcal C)\to \colim{\substack{\lambda\\ \lambda\text{ str. limit}}}\Cond{\lambda}(\mathcal C)\cong \Cond{}(\mathcal C)\] is fully faithful for all strong limit cardinals $\kappa$. 
\end{proof}
\cref{kappacondensedkappaacc} and \cref{limitscolimitscanbecomputedfinitestagesheaves} imply: 
\begin{cor}\label{condensedextremallydisconnected1}Suppose $\mathcal C$ is a presentable category and $\mu$ is a regular cardinal such that $\mathcal C$ is $\mu$-accessible and $\kappa\geq \mu$ is an uncountable regular cardinal. 
    \begin{romanenum}
    \item The left adjoint $\Cond{\kappa}(\mathcal C)\to\Cond{}(\mathcal C)$ of restriction is fully faithful. Its essential image consists of the condensed hypersheaves $\Pro(\Fin)^{\operatorname{op}}\to \mathcal C$ which preserve $\kappa$-filtered colimits. 
    \item In particular, $\Cond{\kappa}(\mathcal C)\subseteq \Cond{}(\mathcal C)$ is closed under $\mu$-small limits and small colimits. 
    \item If $\kappa>>\mu$, then $\Cond{\kappa}(\mathcal C)\subseteq \Cond{}(\mathcal C)$ is closed under $\kappa$-small limits.
\item The category $\Cond{}(\mathcal C)$ has small limits and colimits and small limits and colimits can always be computed in some stage $\Cond{\lambda}(\mathcal C)\subseteq \Cond{}(\mathcal C)$ for an uncountable regular cardinal $\lambda\geq \mu$.
    \end{romanenum}
\end{cor}
\begin{proof}This is a special case of \cref{kappacondensedkappaacc} and \cref{limitscolimitscanbecomputedfinitestagesheaves}, note that $\mathcal C$ is $\lambda$-accessible for all $\lambda\gg \mu$ by \cite[Lemma 5.4.2.10]{highertopostheory}.
\end{proof} 
\begin{cor}\label{condensedfiltrationinfinitytopoi}
Denote by $\Lambda^{r/s}$ the posets of small regular/strong limit cardinals, respectively. 
The functor \[ \Cond{*}(\an)\colon \Lambda^{r/s}\to \Pr^L, \lambda\mapsto \Cond{\lambda}(\an)\] with transition maps given by the left adjoints of the restriction is an exhaustion of $\Cond{}(\an)$ by topoi (\cref{definitionbigpresentable}). 
\end{cor}
\begin{proof}
    By \cref{condensedascolimit}, $\Cond{}(\an)$ is the colimit of both functors in $\vlCat$, and by \cref{leftKanextensionsfullyfaithfulcondensed}, the transition functors $\Cond{\kappa}(\an)\to\Cond{\lambda}(\an)$ are fully faithful left adjoints for all pairs of regular or strong limit cardinals $\kappa\leq \lambda$, respectively. 
    As for all uncountable cardinals $\kappa$, \[\Pro(\Fin)_{\kappa}\subseteq \Pro(\Fin)\] is closed under finite limits, the left adjoint $\Cond{\kappa}(\an)\to\Cond{\lambda}(\an)$ preserves finite limits by \cref{criterioncoveringflat}. 
\end{proof}
As the condensed topology is finitary, the following holds: 
\begin{lemma}\label{filteredcolimitsincondensedncategorycomputedpointwise}If $\mathcal C$ is a presentable $(n,1)$-category so that filtered colimits commute with finite limits in $\mathcal C$, then filtered colimits in $\Cond{(\kappa)}(\mathcal C)\subseteq \Fun(\Pro(\Fin)_{(\kappa)}^{\operatorname{op}},\mathcal C)$ can be computed pointwise. 
\end{lemma}
\begin{proof}
    Suppose that $F\colon I\to \Cond{\kappa}(\mathcal C)$ is a filtered diagram and denote by $\colim{I}F$ its colimit in $\Fun(\Pro(\Fin)_{(\kappa)}^{\operatorname{op}},\mathcal C)$. As $\mathcal C$ has filtered colimits, this colimit is computed pointwise. 
    Since finite limits commute with filtered colimits in $\mathcal C$, $\colim{I}F$ preserves finite products. 
    Suppose that $U_*\colon \Delta_{s,+}^{\operatorname{op}}\to \Pro(\Fin)_{(\kappa)}$ is a condensed hypercover. 
    As $\Delta_{s,\leq n}\subseteq \Delta_{s}$ is left $n$-cofinal (\cite[Example 6.14]{du2023reshapinglimitdiagramscofinality}), for a functor $G\colon \Pro(\Fin)_{(\kappa)}^{\operatorname{op}}\to\mathcal C$, 
\[\clim{\Delta_{s}}G(U_*)\cong \clim{\Delta_{s,\leq n}}G(U_*).\] 
    As $\Delta_{s,\leq n}$ is finite and filtered colimits commute with finite limits in $\mathcal C$, \cref{sheafconditionexplicitkappacondensed}/\cref{definitioncondensed} imply that \[\clim{\Delta_s}\colim{I}F(U_*)\cong \clim{\Delta_{s,\leq n}}\colim{I}F(U_*)\cong \colim{I}\clim{\Delta_{s,\leq n}}F(U_*)\cong \colim{I}F(U_{-1}),\] whence $\colim{I}F\in \Cond{(\kappa)}(\mathcal C)$ by \cref{sheafconditionexplicitkappacondensed}/\cref{definitioncondensed}.  
\end{proof}
\subsubsection{Condensed sets and topological spaces}\label{section:condensedsetsandtopologicalspaces}
In this section, we recall how $\To$ topological spaces embed into condensed sets.
\begin{proposition}[{\cite[Proposition 2.15]{Scholzecondensed}}]\label{underlinexaccessible}
    Suppose $X$ is a $\To$ topological space. 
    If $\kappa>|X|$ is an uncountable regular cardinal, then $\mathcal C(-,X)\colon\Pro(\Fin)^{\operatorname{op}}\to \Set$ is $\kappa$-accessible. 
    In particular, \[\underline{X}\coloneqq \mathcal C(-,X)\colon \Pro(\Fin)^{\operatorname{op}}\to \Set\] defines a condensed set.  
\end{proposition}

        \begin{rem}
        If $X$ is a topological space which is not $\To$, then $\cont(-,X)$ is not accessible and in particular not a condensed set, see\ \cite[Warning 2.14]{Scholzecondensed}.
    \end{rem}
The proof of \cref{underlinexaccessible} uses the following observation. 
    \begin{lemma}[{\cite[Theorem 3]{STONE1979203}}]\label{inverselimitnonempty}
        Suppose that $F\colon I\to \CH$ is a cofiltered diagram with $\clim{I}F=\emptyset$. 
        Then there exists $i\in I$ such that $F(i)=\emptyset$. 
    \end{lemma}
    \begin{proof}[Proof of \cref{underlinexaccessible}]We follow the proof of \cite[Proposition 2.15]{Scholzecondensed}.
        As every continuous surjection between totally disconnected compact Hausdorff spaces is a quotient map, for every topological space $X$, \[\mathcal C(-,X)\colon \Pro(\Fin)^{\operatorname{op}}\to \Set\subseteq \an\] is a condensed hypersheaf. 
        Suppose now that $X$ is a $\To$ topological space and $\kappa>|X|$ is an uncountable regular cardinal. By \cref{profincoaccessible}, we have to show that for $Y\in \Pro(\Fin)$, the canonical map 
    \[c\colon \colim{K\in \oc{{\Pro(\Fin)_{\kappa}^{\operatorname{op}}}}{Y}}\mathcal C(K,X)\to \mathcal C(Y,X)\] is a bijection.
    Fix $Y\in\Pro(\Fin)$ and suppose that $(f_i\colon K_i\to X, t_i\colon Y\to K_i)_{i=1,2}$ are continuous maps with $K_i\in \Pro(\Fin)_{\kappa}$ such that $f_1\circ t_1=f_2\circ t_2$. 
    Since $Y$ is compact, \[\tilde K\coloneqq t_1(Y)\times t_2(Y)\subseteq K_1\times K_2\] is a compact subspace, and hence profinite of weight $<\kappa$. 
    By construction, the two maps \[f_i\circ \pi_i\colon K_1\times K_2\to X, i=1,2\] restrict to the same map $\tilde K\to K_i\to X$, whence $c$ is injective.
    We now show surjectivity. Fix $f\colon Y\to X\in \mathcal C(Y,X)$.
    Below, we will construct $K\in\oc{{\Pro(\Fin)_{\kappa}^{\operatorname{op}}}}{Y}$ such that the two maps \[Y\times_K Y\rightrightarrows Y\xrightarrow{f} X\] agree. 
    If such $K$ is found, denote by $\tilde K\subseteq K$ the image of the map $Y\to K$. 
    As $Y$ is compact, $\tilde K$ is profinite (as closed subspace of the profinite space $K$). It has $\wt(\tilde K)\leq \wt(K)< \kappa$. 
    Since $Y\times_K Y\rightrightarrows Y\to \tilde K$ is a coequalizer diagram, $f$ factors over $\tilde f\colon\tilde K\to X$, which then proves surjectivity of the map $c$. 
    We now construct a space $K$ as above. 
    The diagonal induces an injective continuous map \[l\colon Y\to \clim{K\in \oc{{\Pro(\Fin)_{\kappa}^{\operatorname{op}}}}{Y}} Y\times_{K} Y.\] 
    \cref{profinitesetseparatedbyfunctions} implies that $l$ is surjective: Suppose $x,y\in Y, x\neq y$. 
    By \cref{profinitesetseparatedbyfunctions}, there exists a continuous map $s\colon Y\to \{0,1\}$ with $s(x)\neq s(y)$.
    Fix $(f\colon Y\to K) \in \oc{{\Pro(\Fin)_{\kappa}^{\operatorname{op}}}}{Y}$ and denote by \[\pi_f\colon \clim{K\in \oc{{\Pro(\Fin)_{\kappa}^{\operatorname{op}}}}{Y}} Y\times_{K} Y\to Y\times_K Y\] the projection. Since $f$ factors as 
    \[Y\xrightarrow{s\times f}\{0,1\}\times K\xrightarrow{\pi_K}K,\] $(x,y)\notin \im(\pi_f)$, which shows that $\im(\pi_f)=\Delta_Y$ for all $f\colon Y\to K\in  \oc{{\Pro(\Fin)_{\kappa}^{\operatorname{op}}}}{Y}$. 
    As $\oc{{\Pro(\Fin)_{\kappa}^{\operatorname{op}}}}{Y}$ is cofiltered, this shows that $l$ is surjective, and hence a homeomorphism as its target and source are compact Hausdorff. 
    For $K\in\oc{{\Pro(\Fin)_{\kappa}^{\operatorname{op}}}}{Y}$ let \[d_K\colon Y\times_K Y\to Y\times Y\to X\times X.\] 
    Suppose $x,y\in X$, $x\neq y$. Since $\oc{{\Pro(\Fin)_{\kappa}^{\operatorname{op}}}}{Y}$ is cofiltered and
    \[\clim{K\in \oc{{\Pro(\Fin)_{\kappa}^{\operatorname{op}}}}{Y}}d_K^{-1}(x,y)=\clim{K\in \oc{{\Pro(\Fin)_{\kappa}^{\operatorname{op}}}}{Y}} Y\times_K Y\times_{X\times X}\{(x,y)\}\cong \Delta_{Y}\times_{X\times X}\{(x,y)\}=\emptyset, \] by \cref{inverselimitnonempty} there exists $K_{(x,y)}\in \oc{{\Pro(\Fin)_{\kappa}^{\operatorname{op}}}}{Y}$ with $d_{K_{(x,y)}}^{-1}(x,y)=\emptyset$. 
    As $|X^2|<\kappa^2= \kappa$, 
    \[K_{\infty}\coloneqq \prod_{\substack{(x,y)\in X^2\\ x\neq y}}K_{(x,y)}\in \Pro(\Fin)_{\kappa}\] (\cref{weightofsubspacesproducts}). Since for all $(x,y)\in X\times X, x\neq y$, \[t\colon Y\times_{K_{\infty}}Y\to Y\times Y\xrightarrow{f\times f} X\times X\] factors as 
    \[Y\times_{K_{\infty}}Y\to Y\times_{K_{(x,y)}}Y\to Y\times Y\to X\times X, \] $\im(t)\subseteq \Delta_{X}$, i.e.\ the two maps $Y\times_K Y\rightrightarrows Y\xrightarrow{f} X$ agree. 
    This shows that $c$ is an isomorphism, i.e.\ \[\underline{X}\colon\Pro(\Fin)^{\operatorname{op}}\to \Set\] preserves $\kappa$-filtered colimits.
    \end{proof}
\begin{cor}[{\cite[Proposition 1.7]{Scholzecondensed}}]\label{underlinefunctor}
    Denote by $\Tkcont$ the category of $\To$ topological spaces with $k$-continuous maps (\cref{definitionkcontinuous}). 
    The assignment $X\mapsto \underline{X}$ defines a fully faithful functor 
    \begin{align*}
        \underline{(-)}\colon \Tkcont\to \Cond{}(\Set).
    \end{align*}
    \end{cor}
    \begin{proof}
    For $\To$ topological spaces $X,Y$, choose an uncountable regular cardinal $\kappa>|X|+|Y|$. 
    Then $\underline{X}, \underline{Y}$ are $\kappa$-accessible by \cref{underlinexaccessible}, whence 
    \[\Hom_{\Cond{}(\Set)}(\underline{X}, \underline{Y})\cong \Hom_{\Cond{\kappa}(\Set)}(\underline{X}_{\kappa}, \underline{Y}_{\kappa})\] by \cref{kappacondensedkappaacc}.  
    By \cref{kappacontinuousfullyfaithfullyintocondensed}, 
    $\Hom_{\Cond{\kappa}(\Set)}(\underline{X}_{\kappa}, \underline{Y}_{\kappa})
    \cong \mathcal C_{\kappa}(X,Y)$. 
    The same holds for regular cardinals $\lambda\geq \kappa$, which implies that $\mathcal C_{\kappa}(X,Y)= \mathcal C_{k}(X,Y)$. This shows that $\underline{(-)}$ is fully faithful.  
    \end{proof}
    We will show below (\cref{Underlinefunctorrightadjointeithoutkappa}) that $\underline{(-)}$ has a partially defined left adjoint. 

\cref{profincoaccessible,underlinexaccessible} imply the following: 
\begin{cor}\label{topspacecomefromregular}
    Suppose $X$ is a $\To$ topological space and $\kappa$ is a regular or strong limit cardinal such that one of the following is satisfied: 
    \begin{romanenum}
        \item $\cof(\kappa)>|X|$.\footnote{The cofinality $\cof(\kappa)$ of a cardinal $\kappa$ is the minimal cardinal $\lambda$ with a cofinal functor $\beta\to\kappa$. A map $f\colon \beta\to \kappa$ is cofinal if for $\alpha_0<\alpha$, there exists $\beta_0<\beta$ with $\alpha_0<f(\beta_0)$. More concretely, \[\operatorname{cof}(\kappa)=\min\{ |I|\, |\, \kappa=\sum_{i\in I}\lambda_i, \forall i\in I,|\lambda_i|<\kappa\}.\] }
        \item $X$ is compact Hausdorff and $\kappa>\wt(X)$. 
    \end{romanenum}
    Then the left adjoint 
    $\Cond{\kappa}(\an)\to\Cond{}(\an)$ maps $\underline{X}_{\kappa}$ to $\underline{X}$. 
\end{cor}
\begin{proof}
    By \cref{condensedfiltrationinfinitytopoi} and \cref{fullyfaithfulnessandpreservationoflimits}, for a regular and strong limit cardinals $\kappa$, the left adjoint $\Cond{\kappa}(\an)\to\Cond{}(\an)$ of the restriction is fully faithful. 
    As $\underline{X}$ restricts to $\underline{X}_{\kappa}$, it therefore suffices to show that $\underline{X}$ lies in the essential image of $\Cond{\kappa}(\an)\hookrightarrow \Cond{}(\an)$. 
    Since $\kappa\geq \cof(\kappa)$, \[\Cond{\operatorname{cof}(\kappa)}(\an)\to \Cond{}(\an)\] factors over left adjoints \[\Cond{\operatorname{cof}(\kappa)}(\an)\to \Cond{\kappa}(\an)\to \Cond{}(\an)\] of the restrictions, so it is enough to show that $\underline{X}$ lies in the essential image of \[\Cond{\operatorname{cof}(\kappa)}(\an)\to \Cond{}(\an), \] i.e.\ that $\underline{X}$ preserves $\cof(\kappa)$-filtered colimits (\cref{condensedextremallydisconnected1}). 
    If $\cof(\kappa)>|X|$, this holds by \cref{underlinexaccessible}. 
    If $X$ is compact Hausdorff and $\kappa>\wt(X)$, then \[\mathcal C(-,X)\colon \CH^{\operatorname{op}}\to \Set\] preserves $\kappa$-filtered colimits by \cref{profincoaccessible}. 
    As $\Pro(\Fin)\subseteq \CH$ is closed under limits, this implies that $\underline{X}\colon\Pro(\Fin)^{\operatorname{op}}\to\an$ preserves $\kappa$-filtered colimits. 
\end{proof}
\cref{condensedascolimit} and \cref{compactprojectivegeneratorsforcondensedsets} imply: 
\begin{cor}\label{compactprojectiverepresentables}
    \begin{romanenum}
    \item For $X\in \edCH{}$, $\underline{X}$ is a compact projective condensed anima. 
    Every condensed anima is a small colimit of \textit{representables} $\underline{X}, \, X\in\edCH{}$. 
    
    \item If $K$ is compact Hausdorff, then $\underline{K}$ is compact in $\tau_{\leq n}\Cond{}(\an)$ for all $n\in\mathbb N_0$.
    \end{romanenum}
\end{cor}
\begin{proof}
    We first show that $\underline{X}, X\in\edCH{}$ generate $\Cond{}(\an)$ under small colimits.
    For a condensed anima $T\in \Cond{}(\an)$ choose a strong limit cardinal $\kappa$ with $T\in\Cond{\kappa}(\an)\subseteq \Cond{}(\an)$. 
    By \cref{compactprojectivegeneratorsforcondensedsets}, there exists a small diagram $X_*\colon I\to \edCH{\kappa}$ such that $\underline{T}_{\kappa}\cong \colim{I}\underline{X_i}_{\kappa}$. 
    As $\Cond{\kappa}(\an)\subseteq\Cond{}(\an)$ is closed under small colimits (\cref{condensedfiltrationinfinitytopoi,fullyfaithfulnessandpreservationoflimits}) and maps $\underline{X}_{\kappa}$ to $\underline{X}$ for $X\in\edCH{\kappa}$ (\cref{topspacecomefromregular}), this shows that $\underline{X},X\in\edCH{}$ generates $\Cond{}(\an)$ under small colimits. 

    We now show that for $X$ extremally disconnected compact Hausdorff, $\underline{X}$ is compact projective in $\Cond{}(\an)$. If $F\colon I\to \Cond{}(\an)$ is a small sifted diagram, choose a regular cardinal $\mu>\wt(X)$ such that $F$ factors over $\Cond{\mu}(\an)\subseteq \Cond{}(\an)$. As $\Cond{\mu}(\an)\subseteq \Cond{}(\an)$ is closed under small colimits (\cref{condensedextremallydisconnected1}) and $\underline{X}_{\mu}\in\Cond{\mu}(\an)$ is compact projective (\cref{compactprojectivegeneratorsforcondensedsets}), \cref{topspacecomefromregular} implies that \begin{align*}\colim{I}\Map_{\Cond{}(\an)}(\underline{X},F)&\cong\colim{I}\Map_{\Cond{\mu}(\an)}(\underline{X}_{\mu},F)\\&\cong\Map_{\Cond{\mu}(\an)}(\underline{X}_{\mu}, \colim{I}F)\cong\Map_{\Cond{}(\an)}(\underline{X}, \colim{I}F), \end{align*} which shows that $\underline{X}$ is projective in $\Cond{}(\an)$. 

    Suppose now that $K$ is a compact Hausdorff space and $F\colon I\to \tau_{\leq n}\Cond{}(\an)$ is a filtered diagram. Choose a regular cardinal $\mu>\wt(K)$ such that $F$ factors over $\Cond{\mu}(\an)\subseteq \Cond{}(\an)$. 
    Since $\underline{K}_{\mu}\in \tau_{\leq n}\Cond{\mu}(\an)$ is compact (\cref{compactprojectivegeneratorsforcondensedsets}),
    \[ \Cond{\mu}(\an)\subseteq \Cond{}(\an)\] is closed under small colimits (\cref{condensedextremallydisconnected1}), and \[\Cond{\mu}(\an)\cap\tau_{\leq n}\Cond{}(\an)=\tau_{\leq n}\Cond{\mu}(\an) \text{ (\cref{truncationbigtopos})},\]it follows from \cref{topspacecomefromregular} that  
    \begin{align*}\colim{I}\Map_{\Cond{}(\an)}(\underline{K},F)&\cong\colim{I}\Map_{\Cond{\mu}(\an)}(\underline{K}_{\mu},F)\\&\cong\Map_{\Cond{\mu}(\an)}(\underline{K}_{\mu}, \colim{I}F)\cong\Map_{\Cond{}(\an)}(\underline{K}, \colim{I}F),\end{align*} which proves that $\underline{K}$ is compact in $\tau_{\leq n}\Cond{}(\an)$ for all $n\in\mathbb N_0$. 
\end{proof}

\begin{definition}
    \begin{romanenum}
    \item A condensed set $X$ is \emph{quasi-compact} if for every effective epimorphism $\sqcup_{i\in I}U_i\to P\in \Cond{}(\Set)$, there exists a finite subset $F\subseteq I$ such that $\sqcup_{i\in F}U_i\to P$ is an effective epimorphism. 
    \item A condensed set $X$ is $\To$ if for all $*\to X$ and $\underline{T}\to X\in\oc{\Pro(\Fin)}{X}$, the pullback $P\coloneqq  *\times_X \underline{T}$ is a quasi-compact condensed set. 
    \item Denote by $\To\Cond{}(\Set)\subseteq \Cond{}(\Set)$ the full subcategory on $\To$ condensed sets.
    \end{romanenum} 
\end{definition}
\cref{compactprojectiverepresentables} implies: 
\begin{cor}\label{compacthausdorffisquasicompact}
For a compact Hausdorff space $K$, $\underline{K}$ is a quasi-compact condensed set. 
\end{cor}
\begin{proof}
    Suppose that $K$ is a compact Hausdorff space and $\{X_i\to \underline{K}\}_{i\in I}$ is a small collection of maps so that $\sqcup_{i\in I}X_i\to \underline{K}$ is an effective epimorphism. 
    Choose a regular cardinal $\kappa>\wt(K)$ such that $X_i\in\Cond{\kappa}(\an)$ for all $i\in I$. 
    By \cref{topspacecomefromregular}, $\underline{K}=\underline{K}_{\kappa}\in\Cond{\kappa}(\an)\subseteq \Cond{}(\an)$. 
    As $\Cond{\kappa}(\an)\subseteq \Cond{}(\an)$ is closed under small colimits and finite limits (\cref{condensedextremallydisconnected1}), $\sqcup_{i\in I}X_i\to K$ is an effective epimorphism in $\Cond{\kappa}(\an)$, i.e. $\sqcup_{i\in I}\tau_{\leq 0}X_i\to \underline{K}$ is an epimorphism in $\Cond{\kappa}(\Set)$. By effectivity of epimorphisms and universality of colimits in the 1-topos $\Cond{\kappa}(\Set)$, this means that \[\underline{K}=\underline{K}_{\kappa}\cong \Coeq(\sqcup_{i,j\in I}\tau_{\leq 0}X_i\times_{\underline{K}}\tau_{\leq 0}X_j\rightrightarrows \sqcup_{i\in I}\tau_{\leq 0}X_i).\]  
    The right-hand side is equivalent to \[\colim{\substack{F\subseteq I\\ F\text{ finite }}}\Coeq(\sqcup_{i,j\in F}\tau_{\leq 0}X_i\times_{\underline{K}}\tau_{\leq 0}X_j\rightrightarrows \sqcup_{i\in F}\tau_{\leq 0}X_i).\] 
    Hence by compactness of $\underline{K}=\underline{K}_{\kappa}$ in $\Cond{\kappa}(\Set)$ (\cref{compactprojectivegeneratorsforcondensedsets}), there exists $F\subseteq I$ finite so that $\id_{\underline{K}_{\kappa}}$ factors over \[\Coeq(\sqcup_{i,j\in F}\tau_{\leq 0}X_i\times_{K}\tau_{\leq 0}X_j\rightrightarrows \sqcup_{i\in I}\tau_{\leq 0}X_i)\to \underline{K}_{\kappa}.\] This implies that $\sqcup_{i\in F}\tau_{\leq 0}X_i\to \underline{K}$ is an epimorphism in $\Cond{\kappa}(\Set)$, i.e. $\sqcup_{i\in F}X_i\to \underline{K}_{\kappa}$ is an effective epimorphism in $\Cond{\kappa}(\an)$. As $\Cond{\kappa}(\an)\to\Cond{}(\an)$ preserves finite limits and small colimits, it preserves effective epimorphisms, which proves that $\underline{K}$ is quasi-compact. 
\end{proof}
\begin{proposition}[{\cite[Proposition 1.7, Proposition 2.15]{Scholzecondensed}}]\label{Underlinefunctorrightadjointeithoutkappa}
    The functor \[\underline{(-)}\colon \Tkcont\to \Cond{}(\Set)\] factors over $\To\Cond{}(\Set)\subseteq \Cond{}(\Set)$. 
    The induced functor \[\underline{(-)}\colon \Tkcont\to \To\Cond{}(\Set)\] admits a left adjoint which sends a $\To$ condensed set $X$ to $X(*)_{\operatorname{top}}\coloneqq \colim{\substack{\underline{K}\to X\\ K \in \Pro(\Fin)}} K$. This is a compactly generated $\To$ topological space with set of points $X(*)$. 
\end{proposition}
\begin{proof}
    Suppose $X$ is a $\To$ topological space. 
    By fully faithfulness of $\underline{(-)}$, an element of $\oc{{\Pro(\Fin)_{\kappa}}}{\underline{X}}$ consists of a $\kappa$-light profinite set $T$ and a $k$-continuous and hence continuous map $T\to X$. As $X$ is $\To$, the pullback $T\times_{X}*$ in topological spaces is a closed subspace of $T$ and in particular compact Hausdorff.
    Since for all $K\in\Pro(\Fin)_{\kappa}$, \[\underline{T\times_X *}(K)=\mathcal C(K,T\times_X *)\cong \mathcal C(K,T)\times_{\mathcal C(K,X)}\mathcal C(K,*)\cong (\underline{T}\times_{\underline{X}}*)(K), \] \[\underline{T\times_X *}\cong \underline{T}\times_{\underline{X}}*.\]
    By \cref{compacthausdorffisquasicompact}, $\underline{T\times_X*}$ is quasi-compact, which shows that $\underline{X}$ is $\To$. 
    
    We now show that $\To\Cond{}(\Set)\ni X\mapsto X(*)_{\operatorname{top}}$ is well-defined. Suppose $X\in \To\Cond{}(\Set)$ and choose a regular cardinal such that $X$ is $\kappa$-accessible, i.e. $X\in\Cond{\kappa}(\Set)\subseteq \Cond{}(\Set)$. As $X$ is left Kan extended from \[\Pro(\Fin)_{\kappa}^{\operatorname{op}}\subseteq \Pro(\Fin)^{\operatorname{op}}, \] the colimit \[X(*)_{\operatorname{top}}=\colim{K\in\oc{\Pro(\Fin)}{X}}K=\colim{K\in \oc{{\Pro(\Fin)_{\kappa}}}{X}} K=X(*)_{\kappa}\] exists. It is $\kappa$-compactly generated and the underlying set of points is $X(*)$ by \cref{kappacontinuousfullyfaithfullyintocondensed}.
    We now want to show that $X(*)_{\operatorname{top}}$ is $\To$. 
    As $X(*)_{\operatorname{top}}$ is compactly generated, by \cref{compacthausdorffquotientoftdch} it suffices to show that for all $x\in X(*)_{\operatorname{top}}$ and all $K\in\oc{{\Pro(\Fin)_{\kappa}}}{X}$, $K\times_{X(*)_{\operatorname{top}}} \{x\}\subseteq K$ is a closed subspace, i.e.\ $K\times_{X(*)_{\operatorname{top}}}\{x\}$ is compact. 
    Fix $K\in\Pro(\Fin)_{\kappa}$ with a map $\underline{K}\to X$. 
    As \[\Cond{}(\Set)\to \Set, \,  F\mapsto F(*)\] preserves limits, the  canonical map 
    \[ (\underline{K}\times_{X}\{x\})(*)_{\operatorname{top}}\to K\times_{X(*)_{\operatorname{top}}}\{x\}\] is a continuous bijection, hence it is enough to show that $(\underline{K}\times_{X}\{x\})(*)_{\operatorname{top}}$ is compact. 
    As $\underline{K}\times_{X}\{x\}\in\Cond{\kappa}(\Set)\subseteq \Cond{}(\Set)$,  \[q\colon \sqcup_{Y\in\oc{\Pro(\Fin)_{\kappa}}{\underline{K}\times_{X}\{x\}}}\underline{Y}_{\kappa}\to \underline{K}\times_{X}\{x\}\] is an effective epimorphism in $\Cond{\kappa}(\an)$ and hence in $\Cond{}(\an)$ by \cref{condensedextremallydisconnected1}. 
    Since $X$ is $\To$, $\underline{K}\times_{X}\{x\}$ is quasi-compact, i.e. there exists a finite collection of $\kappa$-light profinite sets $(Y_i)_{i=1}^n$ with an effective epimorphism $\sqcup_{i=1}^n\underline{Y_i}_{\kappa}\to  \underline{K}\times_{X}\{x\}$. 
    By \cref{condensedextremallydisconnected1}, this is also an effective epimorphism in $\Cond{\kappa}(\an)$, i.e. 
    \[\underline{K}\times_{X}\{x\}=\Coeq(\sqcup_{i,j=1}^n \underline{Y_i}_{\kappa}\times_{\underline{K}\times_{X}\{x\}}\underline{Y_j}_{\kappa}\rightrightarrows \sqcup_{i=1}^n\underline{Y_i}_{\kappa})\] by effectivity of epimorphisms and universality of colimits in the 1-topos $\Cond{\kappa}(\Set)$. 
    As $(-)(*)_{\kappa}$ is a left adjoint, this implies that \[\sqcup_{i=1}^nY_i=(\sqcup_{i=1}^n\underline{Y_i}_{\kappa})(*)_{\kappa}\to (\underline{K}\times_{X}\{x\})(*)_{\kappa}\] is a continuous surjection, and in particular, $ (\underline{K}\times_{X}\{x\})(*)_{\kappa}$ is compact, which proves that $X(*)_{\operatorname{top}}$ is $\To$. 
    It follows from \cref{kappacontinuousfullyfaithfullyintocondensed} that $X\mapsto X(*)_{\operatorname{top}}$ it is left adjoint to $\underline{(-)}$. 
\end{proof}
\cref{condensedoncompactextremallydisconnected} implies the following which we will use in the proofs of \cref{examplessflatrings,simplicialresolutionisresolution}. 
\begin{cor}\label{approximationbycompacty}
For a Hausdorff space $X$, 
\[\underline{X}_{(\kappa)}\cong  \colim{\substack{K\to X\\ K\in\operatorname{CH}_{(\kappa)}}}\underline{K}_{(\kappa)}\cong \colim{\substack{K\subseteq X\\ K\CH_{(\kappa)}}}\underline{K}_{(\kappa)}\] where the colimit in the middle is over the category $\oc{\CH_{(\kappa)}}{X}$ of ($\kappa$-light) compact Hausdorff spaces with a continuous map to $X$, and the colimit on the right runs over all ($\kappa$-light) compact subspaces of $X$. 
\end{cor}
\begin{proof}
By \cref{condensedoncompactextremallydisconnected,kappacontinuousfullyfaithfullyintocondensed}, 
\[\underline{X}_{\kappa}\cong \colim{\substack{\underline{K}_{\kappa}\to X\\ K\CH_{\kappa}}}\underline{K}_{\kappa}\cong{\colim{\substack{K\to X\\ K\CH_{\kappa}}}}\underline{K}_{\kappa}.\] 
The image of a continuous map $f\colon K\to X$ map is a compact subspace of $X$. As $X$ is Hausdorff, $\im(f)$ is compact Hausdorff, whence $K\to \im(f)$ is a quotient map and in particular, \[\wt(\im(f))\leq \wt(K)\] by \cref{weightofquotients}. 
This shows that the $\kappa$-light compact subspaces of $X$ are cofinal in $\oc{\CH_{\kappa}}{X}$ and hence 
\[\underline{X}_{\kappa}\cong \colim{\substack{K\subseteq X\\ K\CH_{\kappa}}}\underline{K}_{\kappa}\] for all uncountable cardinals $\kappa$. 
Choose a regular cardinal $\kappa>\wt(X)+|X|$. 
As \[\Cond{\kappa}(\an)\subseteq \Cond{}(\an)\] is closed under small colimits, \cref{topspacecomefromregular} and the above imply that \[ \underline{X}=\colim{\substack{K\subseteq X\\ K\CH_{\kappa}}}\underline{K}.\] By \cref{weightofsubspacesproducts}, \[ \colim{\substack{K\subseteq X\\ K\CH_{\kappa}}}\underline{K}=\colim{\substack{K\subseteq X\\ K\CH{}}}\underline{K}.\] 
As the compact subspaces of $X$ are cofinal in compact spaces with a continuous map to $X$, \[\colim{\substack{K\subseteq X\\ K\CH{}}}\underline{K}\cong \colim{\substack{{K}\to X\\ K\CH}}\underline{K}.\qedhere\]  
\end{proof}
We will need the following observation to prove \cref{tstructureleftcomplete}, which enters our discussion of solid modules. 
\begin{lemma}\label{condensedanimaPostnikovcomplete}
    $\Cond{(\kappa)}(\an)$ is Postnikov complete, i.e. $\Cond{(\kappa)}(\an)\cong \clim{n}\left(\tau_{\leq n}\Cond{(\kappa)}(\an)\right)$.  
\end{lemma}
\begin{proof}
    We first show that for all uncountable cardinals $\kappa$, $\Cond{\kappa}(\Set)$ is a replete topos (\cite[Definition 3.1.1]{Bhatt2013ThePT}), i.e.\ that the following holds: If $F\colon\mathbb N^{\geq}\to \Cond{\kappa}(\Set)$ is a diagram such that for all $n\leq m\in\mathbb N_0$, $F_m\to F_n$ is an epimorphism, then the projection $\pi_n\colon \clim{m}F_m\to F_n$ is an epimorphism for all $n\in\mathbb N_0$. 
    Fix $X\in\Pro(\Fin)_{\kappa}$ and $x\in F_n(X)$. We want to construct a sequence of quotient maps of profinite sets $(q_{k}\colon X_{k+1}\to X_k\in\Pro(\Fin)_{\kappa})_{k\in\mathbb N_0}$ and elements $(x_k\in F_{k+n}(X_k))_{k\in\mathbb N_0}$ such that $X_0=X, x_0=x$ and for all $k\in\mathbb N_0$, \[(F_{k+1+n}(X_{k+1})\to F_{k+n}(X_{k+1})\to F_{k+n}(X_{k}))(x_{k+1})=(F_{k+n}(q_{k+1}))(x_k).\]
  
    Suppose that we have constructed a finite such sequence $X_{k}\to X_{k-1}\to \ldots \to X$ and $(x_k, \ldots, x_{1},x)$. 
    By surjectivity of $F_{n+k+1}\to F_{n+k}$, there exists a cover 
    $\{X^i\to X_k\}_{i\in I}$ such that the image of $x_k$ under $F_{n+k}(X)\hookrightarrow \prod_{I}F_{n+k}(X^i)$ lies in the image of $\prod_{I}F_{n+k+1}(X^i)\to \prod_{I}F_{n+k}(X^i)$.  
    Since the condensed topology is finitary, we can assume that $I$ is finite and since $\Pro(\Fin)_{\kappa}$ has finite coproducts and condensed sets send finite coproducts to products, we can reduce to the case that $I$ consists of a single element $X_{k+1}\coloneqq \sqcup_{i\in I}X^i$. Denote by $q_{k+1}\colon X_{k+1}\to X_k$ the induced map and choose $x_{k+1}\in F_{n+k+1}(X_{k+1})$ with \[(F_{n+k+1}(X_{k+1})\to F_{n+k}(X_{k+1}))(x_{k+1})=(F_{n+k}(q_{k+1}))(x_k).\] By repeating this process we obtain a sequence with the claimed properties.
    Let \[X_{\infty}\coloneqq \clim{k\in\mathbb N_0}X_k.\] 
    As $\kappa$-small filtered limit of $\kappa$-light profinite sets, $X_{\infty}$ is $\kappa$-light profinite. 
    For $k\in\mathbb N_0$ the projection $p_k\colon X_{\infty}\to X_k$ is surjective since all transition maps are surjective. 
    By construction, \[((F_{k+n}(p_k))(x_k))_{k\in\mathbb N_0}\in \clim{k\in\mathbb N_0}F_{k+n}(X_{\infty})\subseteq \prod_{k}F_k(X_{\infty}).\] 
    This shows that 
    $(F_n(q_0))(x)$ lies in the image of the projection $\clim{m}F_m(X_{\infty})\to F_n(X_{\infty})$, which proves that $\clim{m}F_m\to F_n$ is an epimorphism of $\kappa$-condensed sets, i.e.\ $\Cond{\kappa}(\Set)$ is replete. 
    
    It now follows from  \cite[Theorem A]{mondal2024postnikovcompletenessrepletetopoi} that the associated hypercomplete topos $\Cond{\kappa}(\an)$ is Postnikov complete. 
    As the large poset of small regular cardinals $\Lambda^r$ is $\aleph_1$-filtered and for $n\in\mathbb N_0$, \[\tau_{\leq n}\Cond{(\an)}\cong \colim{\kappa\in\Lambda^r}\tau_{\leq n}\Cond{\kappa}(\an)\text{ (\cref{truncationbigtopos})},\] \[\colim{\kappa\in\Lambda^r}\clim{n\in\mathbb N_0}(\Cond{\kappa}(\an)_{\leq n})\cong \clim{n\in\mathbb N_0}\Cond{}(\an)_{\leq n},\] see e.g.\ \cite[Corollary A.2.11]{LucasMannthesis} in the large universe. Hence by the above, \[\Cond{}(\an)\cong  \clim{n\in\mathbb N_0}\tau_{\leq n}\Cond{}(\an).\qedhere\] 
\end{proof}
\begin{notation}
We will in the sequel prove results about condensed and $\kappa$-condensed categories in parallel and use $\Cond{(\kappa)}(\mathcal C)$ as a placeholder for the category of condensed (or $\kappa$-condensed) objects in a presentable category $\mathcal C$. \end{notation}

\cref{globalsectionsleftadjointcondensed} implies: 
\begin{cor}\label{globalsectionstexactcocontinuous}
 The global sections functor stabilizes to a $t$-exact, colimits preserving functor
\[\Gamma^{(\kappa)}_{\Sp}\colon\Cond{(\kappa)}(\Sp)\to \Sp.\] 
\end{cor}
\begin{proof}
By \cref{globalsectionsleftadjointcondensed}, $\Gamma^{\kappa}\colon \Cond{\kappa}(\an)\to\an$ is a left adjoint for all uncountable cardinals $\kappa$. Hence its stabilization $\Gamma^k_{\Sp}$ is a left adjoint and $t$-exact by \cref{geometricmorphismstabilization,texactnessgeometricmorphism}. 
As for all regular cardinals $\kappa$, $i_{\kappa}\colon\Cond{\kappa}(\an)\to \Cond{}(\an)$ is fully faithful, \[\Gamma^{\kappa}=\Map_{\Cond{\kappa}(\an)}(*,-)\cong \Map_{\Cond{}(\an)}(*,i_{\kappa}-)\cong \Gamma\circ i_{\kappa},\] where $\Gamma\colon \Cond{}(\an)\to\an$ is the global sections functor for $\Cond{}(\an)$. 
It now follows from \cref{condensedextremallydisconnected1,condensedfiltrationinfinitytopoi,tstructurespectrumobjects} that 
\[ \Gamma_{\Sp}\colon\Cond{}(\Sp)\to \Sp\] is $t$-exact and preserves small colimits as well. 
\end{proof}

\subsection{Condensed and sheaf cohomology}\label{section:condensedsheafcohomology}
In this section, we compare condensed and sheaf cohomology. \cite{Scholzecondensed} constructed a comparison map $\cH{\sheaf}(X,-)\to \ccH^*(\underline{X}, \underline{-})$ from sheaf to condensed cohomology for compact Hausdorff spaces $X$. We apply the identification of $\kappa$-condensed animae with hypersheaves on the category of all $\kappa$-small topological spaces from \cref{section:othermodelscondensedcat} to extend this comparison map to arbitrary topological spaces and then use descent of sheaf and condensed cohomology along \textit{local section covers} to generalize comparison results from \cite{Scholzecondensed}.
\begin{definition}
    If $X$ is a topological space, denote by $\Op(X)$ the poset of open subsets of $X$, ordered by inclusion, and equip it with the Zariski topology, where a sieve $\{U_i\to U\}_{i\in I}\subseteq \oc{\Op(X)}{U}$ is a covering if and only if $U=\cup_{i\in I}U_i$. 
    Denote by $\Shv_{\Zar}(X)$ the topos of $\an$-valued Zariski sheaves on $\Op(X)$. 
    Sheaf cohomology of $X$ is defined as \[\cH{\sheaf}(X,-)\coloneqq \cH{{\Shv}_{\Zar}(X)}(X,-)\colon \stab{{\Shv}_{\Zar}(X)}\to \Sp.\]  
\end{definition}
\begin{rem}
By \cref{derivedcategorymodules}/\cite[Theorem 2.1.2.3]{SAG}, \[\mathcal D(\Shv_{\Zar}(X,\Ab))\cong \LMod{H\mathbb Z}{\stab{\hypershv_{\Zar}(X)}}\subseteq \LMod{H\mathbb Z}{\stab{\Shv_{\Zar}(X)}}\] is equivalent to the category of hypercomplete $H\mathbb Z$-modules in $\stab{\Shv_{\Zar}(X)}$ which is a full subcategory of $\LMod{H\mathbb Z}{\stab{\Shv_{\Zar}(X)}}$. 
In particular, the above definition of sheaf cohomology recovers more classical notions of sheaf cohomology (after taking homotopy groups). 
\end{rem}
\begin{definition}\label{definitioncondensedcohomology}
    We denote by \[\cckH{(\kappa\text{-})}(-,-)\coloneqq \cH{\Cond{(\kappa)}(\an)}(-,-)\colon  \Cond{(\kappa)}(\an)^{\operatorname{op}}\times\Cond{(\kappa)}(\Sp)\to \Sp\] the cohomology functor of $\Cond{(\kappa)}(\an)$, and by \[\icckH{(\kappa\text{-})}(-,-)\coloneqq \icH{\Cond{(\kappa)}(\an)}(-,-)\colon  \Cond{(\kappa)}(\an)^{\operatorname{op}}\times\Cond{(\kappa)}(\Sp)\to \Cond{(\kappa)}(\Sp)\] the internal cohomology, cf.\ \cref{definitioncohomologyinatopos}. 
\end{definition} 
\cref{cohomologycanbecomputedonfinitelelvel} and \cref{condensedfiltrationinfinitytopoi} imply:  
\begin{cor}\label{condensedcohomologycanbecomputedonfinitestage}
    For all strong limit/regular cardinals $\kappa$, 
    \[\ccH(-,-)|_{\Cond{\kappa}(\an)^{\operatorname{op}}\times \Cond{\kappa}(\Sp)}\cong\ckH(-,-).\] 
\end{cor}
\begin{cor}\label{condensedcohomologycanbecomputedonfinitestagetop} Suppose $X$ is a T1 topological space. 
    For all strong limit/regular cardinals $\kappa$ with $\cof(\kappa)>|X|$,
    \[\ccH(\underline{X},-)|_{\Cond{\kappa}(\Sp)}\cong\ckH(\underline{X}_{\kappa},-).\] 
    In particular, if $M$ is a $\To$ topological abelian group, for all strong limit/regular cardinals $\kappa$ with $\cof(\kappa)>|X|,|M|$, 
    \[\ccH(\underline{X}, \underline{M})\cong \ckH(\underline{X}_{\kappa}, \underline{M}_{\kappa}).\] 
\end{cor}
\begin{proof}This follows from \cref{condensedcohomologycanbecomputedonfinitestage} and \cref{topspacecomefromregular}. 
\end{proof}
We now describe a comparison map from sheaf to condensed cohomology. 
\begin{construction}\label{comparisonmapkappacondensedsheaf}
    Suppose $X$ is a topological space and $\lambda\geq\max\{\kappa,|X|\}$ is an uncountable cardinal. 
    Denote by $\Top^{\lambda}$ the category of $\lambda$-small topological spaces and continuous maps. 
    By \cite[Theorem 3.1.21]{Engelkingtopology}, 
    $\Pro(\Fin)_{\kappa}\subseteq \Pro(\Fin)_{\lambda}\subseteq \Top^{\lambda}$, whence the $\kappa$-condensed topology $\tau_{\cond{\kappa}}$ (\cref{definitioncondensedtopology}) defines a coverage on $\Top^{\lambda}$ by \cref{condensedonothercategories}. 
    
    Equip $\oc{\Top^{\lambda}}{X}$ with the Grothendieck topology pulled back from the $\kappa$-condensed topology (\cref{topologyonslice}).
    By \cref{sheavesonslice} and \cref{condensedonothercategories}, \[\hypershv_{\cond{\kappa}}(\oc{\Top^{\lambda}}{X})\cong \oc{\hypershv_{\cond{\kappa}}(\Top^{\lambda})}{\underline{X}_{\kappa}}\cong \oc{\Cond{\kappa}(\an)}{\underline{X}_{\kappa}}.\] The inclusion $\Op(X)\to \oc{\Top^{\lambda}}{X}$ preserves finite limits and is continuous and therefore defines a morphism of sites \[i\colon (\Op(X), \Zar)\to (\oc{\Top^{\lambda}}{X},[\tau_{\cond{\kappa}}]), \] cf.\ \cref{criterioncoveringflat}. 
    We obtain geometric morphisms \[ \oc{\Cond{\kappa}(\an)}{\underline{X}_{\kappa}}\cong \hypershv_{\cond{\kappa}}(\oc{\Top^{\lambda}}{X}, \an)\leftrightarrows \Shv_{\cond{\kappa}}(\oc{\Top^{\lambda}}{X}, \an)\leftrightarrows \Shv_{\Zar}(X), \]  where the left adjoint pair is hypersheafification and forgetting, and the right adjoint pair is induced by the morphism of sites $i$ (\cref{morphismofquasisitesinducesfunctoronsheaves}). We denote by \[\xradg{i}{X}\colon \oc{\Cond{\kappa}(\an)}{\underline{X}_{\kappa}}\leftrightarrows \Shv_{\Zar}(X)\colon \xladg{i}{X}, \xradg{i}{X}\vdash \xladg{i}{X} \] their composite.  
    Composition with the forget functor \[\oc{\Cond{\kappa}(\an)}{\underline{X}_{\kappa}}\to\Cond{\kappa}(\an), \] which is left adjoint to $-\times\underline{X}_{\kappa}$, yields an adjoint pair 
\[ \ladg{i}\colon \Shv_{\Zar}(X)\rightleftarrows \Cond{\kappa}(\an)\colon \radg{i}, \, \ladg{i}\dashv \radg{i}.\] 
As the left-exact functor $\Shv_{\Zar}(X)\to \oc{\Cond{\kappa}(\an)}{\underline{X}_{\kappa}}$ preserves the terminal object $X$, \[\ladg{i}\mathcal C(-,X)=\underline{X}_{\kappa}.\] The adjunction $\ladg{i}\dashv \radg{i}$ therefore induces an equivalence 
\[ \cH{\Zar}(X, \stradg{i}-)\cong \ckH(\underline{X}_{\kappa},-)\in \Fun(\Cond{\kappa}(\Sp), \Sp), \] cf.\ \cref{existenceleftadjointmoregeneral}.
The functor $i_{*, \Sp}\colon \Cond{\kappa}(\Sp)\to \stab{\Shv_{\Zar}(X)}$ is left t-exact as stabilization of a right adjoint (\cref{texactnessgeometricmorphism}), and hence the counit for $\tau_{\geq 0}\vdash(\Cond{\kappa}(\Sp)_{\geq 0}\subseteq \Cond{\kappa}(\Sp))$ yields a natural transformation \[\pi_0i_{*, \Sp}|_{\Cond{\kappa}(\Ab)}\to i_{*, \Sp}|_{\Cond{\kappa}(\Ab)}.\]  
This induces a comparison map \begin{align}\label{comparisonmap}\cH{\sheaf}(X, \pi_0\stradg{i}-)\to  \cH{\sheaf}(X, \stradg{i}-)\cong \ckH(\underline{X}_{\kappa},-)\in\Fun(\Cond{\kappa}(\Ab), \Sp).\end{align}
By construction of the equivalence $\Cond{\kappa}(\an)\cong \hypershv_{\cond{\kappa}}(\Top^{\lambda})$, for $A\in\Cond{\kappa}(\an)$, 
\[\radg{i}A\colon \Op(X)^{\operatorname{op}}\to \an, U\mapsto \Map_{\Cond{\kappa}(\an)}(\underline{U}_{\kappa},A).\] 
This implies that for $A\in \Cond{\kappa}(\Ab)$, \[\pi_0\stradg{i}A\colon \Op(X)^{\operatorname{op}}\to \Ab, U\mapsto\Map_{\Cond{\kappa}(\Set)}(\underline{U}_{\kappa},A).\]
We will denote this sheaf by $A|_{\Op(X)}\coloneqq \pi_0\stradg{i}A$ in the sequel. 
If $A$ is a topological abelian group, then $\pi_0\stradg{i}\underline{A}_{\kappa}=\mathcal C_{\kappa}(-,A)$ is the sheaf of $\kappa$-continuous functions into $A$ (\cref{condensedonothercategories}). 

As $\mathcal C(-,A)\subseteq \mathcal C_{\kappa}(-,A)=:A_{\kappa}$ is a subsheaf, we obtain a map 
\begin{align}\label{comparisonmapkappacondensedsheaftopological}c_{\kappa,A}\colon \cH{\sheaf}(X,A)\to \cH{\sheaf}(X,A_{\kappa})\to \ckH(\underline{X}_{\kappa},\underline{A}_{\kappa}).\end{align}
\end{construction}

The above construction also yields a comparison map from sheaf to condensed cohomology for $\To$ topological spaces: 
\begin{construction}\label{comparisonmapcondensedinfinitysheaf}
 Suppose that $X$ is a $\To$ topological space. By \cref{condensedcohomologycanbecomputedonfinitestagetop}, for regular/strong limit cardinals $\kappa$ with $\cof(\kappa)> |X|$, 
    \[\ccH(\underline{X},-)\cong \ckH(\underline{X},-)\circ (\Cond{}(\Sp)\xrightarrow{r^{\kappa}_{\Sp}}\Cond{\kappa}(\Sp)).\] 
    \ref{comparisonmap} therefore induces a comparison map 
    \[\cH{\sheaf}(X, \pi_0j^*_{\Sp}r^{\kappa}_{\Sp}-)\to \cH{\sheaf}(X,j^*_{\Sp}r^{\kappa}_{\Sp}-)\cong  \ccH(\underline{X},-).\] 
     For $U\in\Op(X)$, $\underline{U}=\underline{U}_{\kappa}\in\Cond{\kappa}(\Set)\subseteq \Cond{}(\Set)$ by \cref{underlinexaccessible}. 
    This implies that for $A\in\Cond{}(\Ab)$, \[A|_{\Op(X)}\coloneqq \pi_0j^*_{\Sp}r^{\kappa}_{\Sp}A\in \Shv_{\Zar}(X,\Ab)\] equals \[\Op(X)^{\operatorname{op}}\to \Ab, U\mapsto \Map_{\Cond{}(\Set)}(\underline{U},A).\] 
In particular, if $A$ is a Hausdorff topological group, then $\underline{A}|_{\Op(X)}=\mathcal C_k(-,A)$.

As $\mathcal C(-,A)\subseteq \mathcal C_{k}(-,A)=:A_{k}$ is a subsheaf, we obtain a map 
\begin{align}\label{comparisonmapkappacondensedsheaftopologicalwithoutkappa}c_{k,A}\colon \cH{\sheaf}(X,A)\to \cH{\sheaf}(X,A_{\kappa})\to \ccH(\underline{X},\underline{A}).\end{align}
\end{construction}
\begin{thm}[{\cite[Theorem 3.2, 3.3]{Scholzecondensed}}]\label{cohomologyclausenscholze}
    Suppose $X$ is a ($\kappa$-light) compact Hausdorff space. 
    
    If $A$ is a product of a discrete abelian group and a finite-dimensional normed $\mathbb R$-vector space, then 
    \[\cH{\sheaf}(X,A)\cong \cckH{(\kappa\text{-})}(\underline{X}_{(\kappa)}, \underline{A}_{(\kappa)}) \] via the comparison map described above. 
\end{thm}
\begin{proof}
This is an immediate consequence of \cite[Theorem 3.2, 3.3]{Scholzecondensed}, we recall their proof and spell this out in \cref{section:proofofcohomologyclausenscholze}.
\end{proof} 
\begin{rem}
    In general, condensed and sheaf cohomology are not isomorphic: 
    If $X$ is a profinite topological space, for all $\mathcal F\in \Shv_{\Zar}(X, \Ab)$ sheaves of abelian groups on $X$, \[H^p_{\operatorname{\sheaf}}(X, \mathcal F)=0\text{ for }p>0, \] see e.g.\ \cite[\href{https://stacks.math.columbia.edu/tag/0A3F}{Tag 0A3F}]{stacks-project}. 
    In contrast, if $X$ is profinite and $\mathbb Z[\underline X]$ is not a projective condensed abelian group, e.g.\ if $X=Y^2$ for a profinite but not discrete space $Y$, see \cite[Proposition 4.8]{Scholzeanalytic}, then there exists a condensed abelian group $M\in \Cond{}(\Ab)$ such that $\Ext^1_{\Cond{}(\Ab)}(\mathbb Z[\underline X],M)\neq 0$. 
\end{rem}

We now want to generalize \cref{cohomologyclausenscholze}. 
\begin{notation}\label{localsectioncover}
    \begin{romanenum}
    \item Suppose $A\in \Cond{\kappa}(\Ab)$ is a $\kappa$-condensed abelian group. 
    We say that a topological space $X$ is $\kappa$-$A$-\emph{exact} if the comparison map 
    \[ \cH{\sheaf}(X,A|_{\Op(X)})\to \ckH(\underline X_{\kappa},A)\] is an isomorphism.
    \item Suppose $A\in\Cond{}(\Ab)$. A $\To$ topological space $X$ is $A$-\emph{exact} if the comparison map 
    \[ \cH{\sheaf}(X,A|_{\Op(X)})\to \ccH(\underline X,A)\] is an isomorphism.
    \item A \emph{local section cover} is a collection of continuous maps $\{X_i\to X\}_{i\in I}$ such that the map $p\colon\bigsqcup_{i\in I} X_i\to X$ admits local sections, i.e.\ there exists an open cover $X=\bigcup_j U_j$ and continuous maps $s_j\colon U_j\to \bigsqcup_{i\in I}X_i$ such that $p\circ s_j\colon U_j\hookrightarrow X$ is the inclusion for all $j\in J$.

    \item For $A\in\Cond{\kappa}(\Ab)$, a local section cover $\{X_i\to X\}_{i\in I}$ of a topological space is $\kappa$-$A$-\emph{exact} if for all $(i_1,\ldots,i_n)\in I^n$, $n\in\mathbb N_1$,  
    $X_{i_1}\times_X \times \ldots \times_{X}X_{i_n}$ is $\kappa$-$A$-exact.
    \item For $A\in\Cond{}(\Ab)$, a local section cover $\{X_i\to X\}_{i\in I}$ of $\To$ topological spaces is $A$-\emph{exact} if for all $(i_1,\ldots,i_n)\in I^n$, $n\in\mathbb N_1$,  
    $X_{i_1}\times_X \times \ldots \times_{X}X_{i_n}$ is $A$-exact.
    \end{romanenum} 
\end{notation}
 Our next goal is to prove the following: 
\begin{proposition}\label{goodcovers}
    \vspace{0pt}\noindent
    \begin{romanenum}
    \item For a $\kappa$-condensed abelian group $A\in\Cond{\kappa}(\Ab)$, a topological space is $\kappa$-$A$-exact if and only if it admits a $\kappa$-$A$-exact local section cover.
    \item For a condensed abelian group $A$, a $\To$ topological space is $A$-exact if and only if it admits an $A$-exact local section cover by $\To$ topological spaces. 
    \end{romanenum}
\end{proposition}
The only-if statement is trivial: 
If $X$ is ($\kappa$-)$A$-exact, then $\{id_X\colon X\to X\}$ is an $A$-exact local section cover. 
To prove the converse, we use the \textit{gros topos} to interpolate between $\Shv_{\Zar}(X)$ and $\Cond{\kappa}(\an)$, the proof will be completed on page \pageref{proofgoodcovers}. 
\begin{definition}\label{grostoposdefinition}Fix an uncountable cardinal $\lambda$ and denote
     by $\Top^{\lambda}$ the category of $\lambda$-small topological spaces. 
     The local section covers constitute a coverage on $\Top^{\lambda}$ (\cref{grostopos}). 
    Denote by $(\Top^{\lambda},LS)$ the associated site. 
    The corresponding topos $\grostop\coloneqq {\Shv}_{LS}(\Top^{\lambda})$ is called the \emph{gros topos} of ($\lambda$-small) topological spaces. 
\end{definition} 
\begin{rem}The local section topology is generated by open covers which form a coverage on $\Top^{\lambda}$.
    The local section topology is subcanonical, i.e.\ for all $x\in\Top^{\lambda}$, $h_x\coloneqq \mathcal C(-,x)\colon \Top^{\lambda}\to \Set$ is an $LS$-sheaf. 

    The category of topological spaces is not accessible. For example, the Sierpi\'nsky space is not $\kappa$-compact for any regular cardinal $\kappa$, see e.g.\ \cite[p. 49]{Hovey}. We can therefore not apply the formalism from \cref{explicitcoveringssitessection} to study sheaves on the large category $\Top$.
\end{rem}
Suppose $\lambda\geq \kappa$ are uncountable cardinals.
For a $\lambda$-small topological space $X\in \Top^{\lambda}$ equip $\oc{\Top^{\lambda}}{X}$ with the topology pulled back from the local section topology (\cref{topologyonslice}). Then ${\Shv}_{LS}(\oc{\Top^{\lambda}}{X})\cong \oc{{\grostop}}{{h_X}}$ by \cref{sheavesonslice}. 
The inclusion \[t_{X}\colon \Op(X)\to \oc{\Top^{\lambda}}{X}\] obviously defines a continuous functor $(\Op(X), \Zar)\to (\oc{\Top^{\lambda}}{X},LS)$.  
Since $\Op(X)$ has finite limits and $t_X$ preserves finite limits, $t_{X}$ is a morphism of sites by \cref{criterioncoveringflat}.
Since the local section topology is coarser than the $\kappa$-condensed topology (\cref{grostopos}), the identity defines morphisms of sites 
\[(\Top^{\lambda},LS)\to (\Top^{\lambda},[\tau_{\cond{\kappa}}]).\] \label{grostoposgeometricmorphismcondensed}
\[ (\oc{\Top^{\lambda}}{X},LS)\to (\oc{\Top^{\lambda}}{X},[\tau_{\cond{\kappa}}]).\]
This yields a commutative diagram of left adjoint functors  
\begin{center}\label{sheavesgrostoposgeometricmorphism}
    \begin{tikzcd}
        \Shv_{\Zar}(X)\arrow[r,"\xladg{t}{X}"]\arrow[dr,"\ladg{t}"] & \oc{{\grostop}}{{h_X}}\arrow[d]\arrow[r,"\xladg{j}{X}"] & \oc{{\Cond{\kappa}(\an)}}{\underline{X}_{\kappa}}\arrow[d]\\ 
        & \grostop\arrow[r,"\ladg{j}"]& \Cond{\kappa}(\an),
    \end{tikzcd}
\end{center}
    where the horizontal arrows are induced by the corresponding morphisms of sites (\cref{morphismofquasisitesinducesfunctoronsheaves}) and condensed hypersheafification, and the vertical arrows are the forget functors, which are left adjoint to $h_X\times -$ and $\underline{X}_{\kappa}\times-$, respectively. 
    The composite $\Shv_{\Zar}(X)\to \Cond{\kappa}(\an)$  is the functor $\ladg{i}$ from above. 
    \begin{lemma}\label{grostoposcondensedzariskirepresentables}
    For a $\lambda$-small topological space $X\in\Top^{\lambda}$, 
    $\ladg{t}(\mathcal C(-,X))=h_X$ and $\ladg{j}{{h_X}}=\underline{X}_{\kappa}$. 
    \end{lemma}
    \begin{proof}
    This holds since $\xladg{t}{X}$ and $\xladg{j}{X}$ are left-exact and in particular preserve the terminal object.
    \end{proof}
    By \cref{existenceleftadjointmoregeneral}, this implies the following: 
\begin{cor}\label{grostoposequalssheaf}
    For $X\in\Top^{\lambda}$, \begin{equation*}\cH{\grostop}(h_X,-)\cong \cH{\sheaf}(X, \stradg{t}-),\end{equation*}
and in particular, \[\cH{\grostop}(h_X, \stradg{j}-)\cong \cH{\sheaf}(X, \stradg{t}\stradg{j}-)=\cH{\sheaf}(X, \stradg{i}-).\]
\end{cor}
\begin{lemma}\label{restrictionhasrightadjoint}
The right adjoint $\radg{t}\colon \grostop\to\Shv_{\Zar}(X)$ has a further right adjoint. 

\end{lemma}
\begin{proof}
    Since $-\times h_X\colon\grostop\to\oc{{\grostop}}{{h_X}}$ has a right adjoint, it suffices to show that the restriction \[\xradg{t}{X}\colon \oc{{\grostop}}{{h_X}}\to \Shv_{\Zar}(X)\] has a right adjoint. The morphism of sites \[(\Op(X), \Zar)\to (\oc{\Top^{\lambda}}{X},LS)\] has the covering lifting property (\cref{definitioncoveringliftingproperty}), hence right Kan extension restricts to a functor 
    \[{\Shv}_{\Zar}(\Op(X))\to \Shv_{LS}(\oc{\Top^{\lambda}}{X})\cong \oc{\grostop}{{h_X}}\] which is right adjoint to $\xradg{t}{X}$ by \cref{rightkanextensionsheaves}. 
\end{proof}
\begin{cor}\label{grostoposrecoversheafcohomology} For $A\in \Ab(\tau_{\leq 0}\grostop)$ and a $\lambda$-small topological space $X\in \Top^{\lambda}$, 
\[ \cH{\grostop}(h_X,A)\cong \cH{\sheaf}(X,A|_{\Op(X)}),\] where \[ A|_{\Op(X)}\colon\Op(X)^{\operatorname{op}}\to \Ab, U\mapsto A(U)\] is the restriction of $A$ along $t\colon\Op(X)\to \Top^{\lambda}$. 
\end{cor}
\begin{proof}
\cref{texactnessgeometricmorphism,restrictionhasrightadjoint} imply that $t_{*,\Sp}$ is $t$-exact, whence $t_{*,\Sp}A=\pi_0t_{*,\Sp}A=A\circ t$. The statement now follows from \cref{grostoposequalssheaf}. 
\end{proof}
\begin{cor}\label{characterizationexactnessgrostopos}Suppose $X$ is a $\lambda$-small topological space. 
    \begin{romanenum}
        \item For $A\in \stab{\grostop}_{\leq 0}$, the counit $\pi_0A\to A$ induces an equivalence 
\[ \cH{\grostop}(h_X, \pi_0A)\cong \cH{\grostop}(h_X,A)\] if and only if the counit
$\pi_0\stradg{t}A\to\stradg{t}A$ induces an equivalence \[\cH{\sheaf}(X, \pi_0\stradg{t}A)\cong \cH{\sheaf}(X, \stradg{t}A).\]

\item In particular, for $\kappa\leq \lambda$ and $A\in \Cond{\kappa}(\Ab)$, \[ \cH{\sheaf}(X, \pi_0 \stradg{i}A)\cong \ckH(\underline{X}_{\kappa},A)\] via the comparison map (\ref{comparisonmap}) if and only if $\pi_0\stradg{j}A\to \stradg{j}A$ induces an equivalence \[\cH{{\grostop}}(h_X, \pi_0\stradg{j}A)\cong \cH{{\grostop}}(h_X, \stradg{j}A).\]
    \end{romanenum}
\end{cor} 
\begin{proof}The second statement follows from the first since $\radg{t}\circ\radg{j}=\radg{i}$.  
Since $\radg{t}$ has a left and a right adjoint, $\stradg{t}$ is $t$-exact by \cref{texactnessgeometricmorphism}.  
As the adjunction $\ladg{t}\vdash \radg{t}$ induces an equivalence \[\cH{\grostop}(h_X, -)\cong \cH{\sheaf}(X, \stradg{t}-)\in\Fun(\stab{\grostop},\Sp)\] (\cref{geometricmorphismcohomology}), this implies the first statement.
\end{proof}
The remaining implication of \cref{goodcovers} is now straightforward:
\begin{proof}[Proof of \cref{goodcovers}]\label{proofgoodcovers} Suppose first that $\{X_i\to X\}$ is a $\kappa$-$A$-exact local section cover and choose a cardinal $\lambda\geq \kappa,|\sqcup_{i\in I}X_i\sqcup X|$. By definition of the local section topology, $p\colon \sqcup_{i\in I} h_{X_i}\to h_{X}$ is an effective epimorphism in ${\grostop}$ (\cite[Proposition 6.2.3.20]{highertopostheory}), and hence \[\cH{\grostop}(h_X,-)=\clim{\Delta}(\cH{\grostop}(h_{\check{C}(p)},-))\] by \cref{spectralenrichmentcocontinuous}.  
For all $n\in\Delta$,  \[\check{C}(p)_n=\sqcup_{(i_0, \ldots,i_n)\in I^{n+1}}X_{i_0}\times_X\ldots\times_X X_{i_n}\] is ($\kappa$)-$A$-exact as coproduct of ($\kappa$-)$A$-exact spaces and $\lambda$-small since \[|\sqcup_{(i_0, \ldots,i_n)\in I^{n+1}}X_{i_0}\times_X\ldots\times_X X_{i_n}|\leq |\sqcup_{i\in I}X_i|^{n+1}<\lambda^{n+1}=\lambda.\] 
It now follows from \cref{characterizationexactnessgrostopos} that the map \[\pi_0\stradg{j}A\to\stradg{j}A\] induces an isomorphism of cosimplicial spectra
\[ \cH{\grostop}(h_{\check{C}(p)_*}, \pi_0\stradg{j}A)\cong \cH{\grostop}(h_{\check{C}(p)_*}, \stradg{j}A).\] 
In particular, \[\cH{\grostop}(h_X, \pi_0\stradg{j}A)\cong \cH{\grostop}(h_X, \stradg{j}A), \] whence $X$ is $\kappa$-$A$-exact by \cref{characterizationexactnessgrostopos}. 

Suppose now that $A\in\Cond{}(\Ab)$ and $\{ X_i\to X\}_{i\in I}$ is an $A$-exact local section cover. 
Choose a regular cardinal $\kappa$ such that $|X\sqcup \bigsqcup_{i\in I}X_i|<\kappa$ and $A\in\Cond{\kappa}(\Ab)\subseteq \Cond{}(\Ab)$. 
For a $\kappa$-small $\To$ topological space $Y$, 
\[\cckH{\kappa-}(\underline{Y}_{\kappa},A)\cong \ccH(\underline{Y},A)\] by \cref{condensedcohomologycanbecomputedonfinitestagetop}. 
In particular, by construction of the comparison map (\ref{comparisonmapcondensedinfinitysheaf}), a $\kappa$-small $\To$-topological space $Y$ is $A$-exact if and only if it is $\kappa$-$A$-exact.
Since for all $n\in\mathbb N_0$, \[\check{C}(p)_n=\sqcup_{(i_0, \ldots,i_n)\in I^{n+1}}X_{i_0}\times_X\ldots\times_X X_{i_n}\leq |\sqcup_{i\in I}X_i|^{n+1}<\kappa, \] it follows that $\{ X_i\to X\}_{i\in I}$ is a $\kappa$-$A$-exact local section cover. We have shown above that this implies that $X$ is $\kappa$-$A$-exact, and whence $A$-exact since $|X|<\kappa$. 
\end{proof}
\begin{cor}\label{condensedandsheafcohomology1}
    Suppose $X$ is a ($\kappa$-light), locally compact Hausdorff space. 
    \begin{romanenum}
    \item For a topological abelian group $A$, $\mathcal C(-,A)=\mathcal C_{k}(-,A)=\mathcal C_{\kappa}(-,A)\colon \Op(X)^{\operatorname{op}}\to \Ab$. 
    
    \item If $A$ is a product of a discrete abelian group and a finite-dimensional normed $\mathbb R$-vector space, then 
    \[ \cH{\sheaf}(X,A)\cong \cckH{(\kappa\text{-})}(\underline{X}_{(\kappa)}, \underline{A}_{(\kappa)})\] via the comparison map (\ref{comparisonmapkappacondensedsheaf}/\ref{comparisonmapcondensedinfinitysheaf}). 
    \end{romanenum}
\end{cor}
\begin{rem}
    The statement for discrete abelian groups (and more generally bounded above spectra) was obtained in \cite[Corollary 4.12]{haine2022descentsheavescompacthausdorff} using shape theory. \cite[Proposition 1.3.5, Corollary 1.3.6]{CatrinMairPhDthesis} identified $\cH{\sheaf}(X,A)\cong \cH{\condo}(\underline{X},A)$ for a locally contractible or paracompact, compactly generated Hausdorff space $X$ and a discrete abelian group $A$. 
    We give an alternative proof for locally contractible spaces below, see \cref{condensedandsheafcohomologylocallycontractible}.   
\end{rem}
\begin{proof}
    Since $X$ is locally compact and $\kappa$-light, every point of $x$ has a neighborhood basis by $\kappa$-light compact Hausdorff spaces. 
    This implies that the monomorphisms \[\mathcal C_{\kappa}(-,A)\subseteq\mathcal C_k(-,A)\subseteq\mathcal C(-,A)\] are isomorphisms on stalks and hence isomorphisms.  
    For $x\in X$ choose a compact neighborhood $x\in K_x\subseteq X$. Then $\{ K_x\to X\}_{x\in X}$ is a local section cover and for $x\in X$, $\wt(K_x)\leq \wt(X)$.
    For $(x_1, \ldots,x_n)\in X^n$, $K_{x_1}\times_{X}\ldots \times_{X}K_{x_n}=\cap_{i=1}^n K_{x_i}$ is compact Hausdorff, whence $\{ K_x\to X\}_{x\in X}$ is ($\kappa$-)$A$-exact for $A$ a product of a discrete abelian group and a finite-dimensional normed $\mathbb R$-vector space by \cref{cohomologyclausenscholze}.
\end{proof}

We will prove further identifications of sheaf and condensed cohomology at the end of this chapter, see \cref{section:condensedcohomologywithsolidcoefficients}. 

Recall also the following vanishing result for sheaf cohomology which will be central for our discussion of solid modules.  
\begin{lemma}[{\cite[Theorem 9.11, Theorem 9.16]{Bredon}}]\label{sheafcohomologywvspacecoeffvanonparcomplc}
Suppose $X$ is a paracompact, locally compact Hausdorff space. 
For all finite-dimensional normed $\mathbb R$-vector spaces $V$, \[ \cH{\sheaf}^*(X,V)=\mathcal C(X,V)\in \Sp^{\heart}\] is concentrated in degree $0$.
\end{lemma}

For $(\To)$ topological abelian groups $A$, it is natural to consider the following variant of \cref{localsectioncover}:
\begin{notation}\label{deltakappaexact}Suppose $A$ is a topological abelian group.  
    \begin{romanenum}
\item We say that a topological space $X$ is  $\epsilon$-$\kappa$-$A$-\emph{exact} if the map \[ c_{\kappa,A}(X)\colon \cH{\sheaf}(X,A)\to \cH{\sheaf}(X,A_{\kappa})\to \ckH(\underline{X}_{\kappa},\underline{A}_{\kappa}) \, \, (\ref{comparisonmapkappacondensedsheaftopological})\]  is an equivalence. 

\item If $A$ is $\To$, we say that a $\To$ topological space $X$ is $\epsilon$-$A$-\emph{exact} if the map 
\[c_{k,A}\colon \cH{\sheaf}(X,A)\to \cH{\sheaf}(X,A_{k})\to \ccH(\underline{X},\underline{A})\, \, (\ref{comparisonmapkappacondensedsheaftopologicalwithoutkappa})\] is an equivalence.
    \end{romanenum}  
\end{notation}
\begin{cor}\label{deltakappaexactgoodcover}
    Suppose $A$ is a topological abelian group.
    \begin{romanenum}
    \item A topological space $X$ is $\epsilon$-$\kappa$-$A$-exact if and only if there is a local section cover $\{ X_i\to X\}_{i\in I}$ such that for all $(i_1,\ldots,i_n)\in I^n$, $n\in\mathbb N_1$,  
    \[ X_{i_1}\times_{X}\ldots \times_X X_{i_n}\] is $\epsilon$-$\kappa$-$A$-exact. 
    \item If $A$ is $\To$, a $\To$ topological space is $\epsilon$-$A$-exact if and only if there exists a local section cover $\{ X_i\to X\}_{i\in I}$ by $\To$ topological spaces such that for all $(i_1,\ldots,i_n)\in I^n$, $n\in\mathbb N_1$,  
    \[ X_{i_1}\times_{X}\ldots \times_X X_{i_n}\] is $\epsilon$-$A$-exact.  
    \end{romanenum}
\end{cor}
\begin{proof}
The only-if statements are clear as $\{\id_X\colon X\to X\}$ is a local section cover. 

We show the converse implication for $\epsilon$-$\kappa$-$A$-exactness.
\cref{grostoposequalssheaf} implies that for $\lambda\geq \kappa$, a $\lambda$-small topological space $Y$ is $\epsilon$-$\kappa$-$A$-exact if and only if the map \[\mathcal C(-,A)\to \mathcal C_{\kappa}(-,A)=\pi_0\stradg{j}\underline{A}_{\kappa}\to \stradg{j}\underline{A}_{\kappa}\] induces an equivalence \begin{align*}\cH{\grostop}(h_Y,\mathcal C(-,A))\cong \underbrace{\cH{\grostop}(h_Y,\stradg{j}\underline{A}_{\kappa})}_{\cong \ckH(\underline{Y}_{\kappa},\underline{A}_{\kappa})}.\end{align*}
Suppose that $\{ X_i\to X\}_{i\in I}$ is a local section cover such that for all $(i_1,\ldots,i_n)\in I^n$, $n\in\mathbb N_1$,  
    \[ X_{i_1}\times_{X}\ldots \times_X X_{i_n}\] is $\epsilon$-$\kappa$-$A$-exact and choose $\lambda>|\sqcup_{i\in I}X^i\sqcup X|$. 
    As $\sqcup_{i\in I}h_{X^i}\to h_{X}$ is an effective epimorphism in $\grostop$ (\cite[Proposition 6.2.3.20]{highertopostheory}), by \cref{spectralenrichmentcocontinuous}, 
\[ \cH{\grostop}(h_{X},-)\cong \clim{\Delta}\cH{\grostop}(\check{C}(\sqcup_{i\in I}h_{X_i}\to h_X),-).\] 
For $n\in\mathbb N_0$, \[\check{C}(\sqcup_{i\in I}h_{X^i}\to h_X)([n])=\sqcup_{(i_1,\ldots,i_n)\in I^n}h_{X_{i_1}\times_XX_{i_2}\times_{X}\ldots \times_X X_{i_n}}\] as the Yoneda embedding $h\colon \Top^{\lambda}\to\grostop$ preserves limits. 
It follows that $\mathcal C(-,A)\to \stradg{j}\underline{A}_{\kappa}$ induces an equivalence of cosimplicial spectra 
\[ \cH{\grostop}(\check{C}(\sqcup_{i\in I}h_{X_i}\to h_X),\mathcal C(-,A))\to \cH{\grostop}(\check{C}(\sqcup_{i\in I}h_{X_i}\to h_X),\stradg{j}\underline{A}_{\kappa}),\] and in particular, $X$ is $\epsilon$-$\kappa$-$A$-exact.  

Suppose now that $X$ is $\To$ and $\{X^i\to X\}_{i\in I}$ is an $\epsilon$-$A$-exact local section cover by $\To$ topological spaces. Choose a regular cardinal $\kappa>|X|+|\sqcup_{i\in I}X_i|$. 
By \cref{condensedcohomologycanbecomputedonfinitestage} and construction of the comparison map (\ref{comparisonmapkappacondensedsheaftopologicalwithoutkappa}), a $\kappa$-small $\To$ topological space is $\epsilon$-$A$-exact if and only if it is $\epsilon$-$\kappa$-$A$-exact. 
It now follows from the above that $X$ is $\epsilon$-$A$-exact. 
\end{proof}

We will use the following observation to identify sheaf and condensed cohomology of locally contractible spaces with coefficients in discrete abelian groups (\cref{condensedandsheafcohomologylocallycontractible}):
    \begin{lemma}\label{comparisonmaponbasis}
    Suppose that $X$ is a topological space, $\mathcal B\subseteq \Op(X)$ is a basis for the topology on $X$, i.e. every open subset of $X$ is a union of elements in $\mathcal B$. 
    
    \begin{romanenum}
    \item If $A\in \Cond{\kappa}(\Ab)$ is such that for all $U\in\mathcal B$, the comparison map 
    \[\cH{\sheaf}(U,\pi_0\stradg{i}A)\to \ckH(\underline{U}_{\kappa},A)\] is an equivalence, then the comparison map \[\cH{\sheaf}(X,\pi_0\stradg{i}A)\to \ckH(\underline{X}_{\kappa},A)\] is an equivalence. 
    \item If $A$ is a topological abelian group such that for all $U\in\mathcal B$, the comparison map 
    \[ \cH{\sheaf}(U,A)\to \ckH(\underline{U}_{\kappa},\underline{A}_{\kappa})\] is an equivalence, then the comparison map \[\cH{\sheaf}(X,A)\to \ckH(\underline{X}_{\kappa},\underline{A}_{\kappa})\] is an equivalence. 
    \item If $\phi\colon A\to B\in\Cond{\kappa}(\Ab)$ is such that for all $U\in\mathcal B$, 
    \[\stradg{j}(\phi)\colon \cH{\sheaf}(U,\stradg{i}A)\to \cH{\sheaf}(U,\stradg{i}B)\] is an equivalence, then \[ \stradg{j}(\phi)\colon \ckH(\underline{X}_{\kappa},A)\to \ckH(\underline{X}_{\kappa},B)\] is an equivalence.  
    \end{romanenum}
    \end{lemma}
    \begin{rem}
By \cref{condensedcohomologycanbecomputedonfinitestage}, the corresponding results for condensed cohomology hold as well (for $\To$ topological spaces). 
    \end{rem}

    \begin{proof}
        Suppose first that $A\in\Cond{\kappa}(\Ab)$ is such that for all $U\in\mathcal B$, the comparison map 
    \[\cH{\sheaf}(U,\pi_0\stradg{j}A)\to \ckH(\underline{U}_{\kappa},A)\] is an equivalence, then the comparison map \[\cH{\sheaf}(X,\pi_0\stradg{j}A)\to \ckH(\underline{X}_{\kappa},A)\] is an equivalence. 
    We claim that $\pi_0\stradg{i}A\to\stradg{i}A$ is an equivalence. 
    As right adjoint of a geometric morphism, $\radg{i}\colon \Cond{\kappa}(\an)\to \Shv(X)$ preserves hypercomplete objects, i.e. factors over $\widehat{\Shv}(X)$.
    This implies that $\stradg{i}$ factors over a morphism $\Cond{\kappa}(\Sp)\to \stab{\widehat{\Shv}(X)}\subseteq \stab{\Shv(X)}$. 
    It therefore suffices to show that $\epsilon\colon \pi_0\stradg{i}A\to \stradg{i}A$ is an isomorphism on stalks (cf. \cref{stalkfunctors}). By our assumption, 
        \[\eta(U)\colon \cH{\sheaf}(U,\pi_0\stradg{j}A)\to \cH{\sheaf}(U,\stradg{j}A)=\ckH(\underline{U}_{\kappa},A)\] is an equivalence for all $U\in\mathcal B$. 
        As $\mathcal B$ forms a basis of the topology, for $x\in X$, $\{ B\in\mathcal B\, |\, x\in B\}$ is cofinal in the open neighbourhoods of $x$. 
        It now follows that $\pi_0\stradg{i}A\to \stradg{i}A$ is an isomorphism on stalks and hence an isomorphism. 
        In particular, 
        \[ \cH{\sheaf}(X,\pi_0\stradg{i}A)\cong \cH{\sheaf}(X,\stradg{i}A),\] whence \[\cH{\sheaf}(X,\pi_0\stradg{i}A)\cong \ckH(\underline{X}_{\kappa},\underline{A}_{\kappa})\] via the comparison map.  

        The same argument shows that if $A$ is a topological abelian group such that for all $U\in\mathcal B$, the comparison map 
        \[ \cH{\sheaf}(U,A)\to\ckH(\underline{U}_{\kappa},A)\] is an equivalence, then $A\to\pi_0\stradg{i}A\to \stradg{i}A$ is an equivalence on stalks and hence an equivalence and in particular, 
        \[ \cH{\sheaf}(X,A)\cong \cH{\sheaf}(X,\stradg{j}\underline{A}_{\kappa}),\] i.e. \[\cH{\sheaf}(X,A)\cong \ckH(\underline{X}_{\kappa},\underline{A}_{\kappa})\] via the comparison map.  
    
        The same argument also shows that if $\phi\colon A\to B\in\Cond{\kappa}(\Ab)$ is such that for all $U\in\mathcal B$, 
        \[ \ckH(\underline{U}_{\kappa},A)\to \ckH(\underline{U}_{\kappa},B)\] is an equivalence, then $\stradg{j}(\phi)\colon \stradg{j}A\to \stradg{j}(B)$ is an equivalence on stalks and hence an equivalence. 
        As $\ckH(\underline{X}_{\kappa},-)\cong \cH{\sheaf}(X,\stradg{j}-)$, this implies that 
        \[ \ckH(\underline{X}_{\kappa},A)\to \ckH(\underline{Y}_{\kappa},B)\] is an equivalence.
    \end{proof}

We now deduce from \cref{cohomologyclausenscholze} that internal condensed cohomology of a compact Hausdorff space with discrete coefficients has a trivial condensed structure (\cref{internalcohomologydiscretecompacthausdorff}). This will be central for the proof of \cref{solidificationoftopologicalspaces}. 
\begin{lemma}\label{internalcohomologyprofinite}Suppose $X$ is ($\kappa$-light) profinite set. 
For all discrete abelian groups $A$, 
        \[\icckH{(\kappa\text{-})}(\underline{X}_{(\kappa)}, \underline{A}_{(\kappa)})\cong \underline{\mathcal C(X,A)}_{(\kappa)}\] is represented by the discrete abelian group $\mathcal C(X,A)$ and concentrated in degree $0$. 
\end{lemma}
\begin{proof}
    By \cref{condensedcohomologycanbecomputedonfinitestagetop} and \cref{forkappalargeenoughclosedmonoidal}, it suffices to show the $\kappa$-condensed statement. 
For $T\in\Cond{\kappa}(\an)$,  
\[\ickH(T, \underline{A}_{\kappa})\in\Cond{\kappa}(\Sp)\cong\Fun^{\lim}(\Cond{\kappa}(\an)^{\operatorname{op}}, \Sp)\] is the small limits-preserving functor \[\ckH(T\times -,A)\colon\Cond{\kappa}(\an)^{\operatorname{op}}\to \Sp, \] cf.\ \cref{enrichmentrecoversinternalhom}. 
By \cite[Theorem 3.2]{Scholzecondensed}/\cref{cohomologyclausenscholzeprofinite}, for $X,S\in\Pro(\Fin)_{\kappa}$, \[ \cckH{(\kappa\text{-})}(\underline{X}_{\kappa}\times \underline{S}_{\kappa}, \underline{A}_{\kappa})\cong \cH{\sheaf}(X\times S,A)\cong \mathcal C(X\times S,A)\] is concentrated in degree $0$. 
This shows that for $X\in\Pro(\Fin)_{\kappa}$, \[\ickH(\underline{X}_{\kappa},A)\cong \underline{\mathcal C(X,A)}_{\kappa}, \] where $\mathcal C(X,A)$ is endowed with the compact open topology which is discrete since $X$ is compact Hausdorff and $A$ is discrete.
\end{proof}
    Denote by $c_{\Sp}\colon \Sp\to\Cond{(\kappa)}(\Sp)$ the stabilization of the constant sheaf functor. 
\begin{cor}\label{internalcohomologydiscretecompacthausdorff}
    Suppose that $X$ is a ($\kappa$-light) compact Hausdorff space. 
        For all discrete abelian groups $A$, the counit 
        \[c_{\Sp}\cckH{(\kappa\text{-})}(X, \underline{A}_{(\kappa)})\to\icckH{(\kappa\text{-})}(\underline{X}_{(\kappa)}, \underline{A}_{(\kappa)})\] is an equivalence. 
\end{cor}
\begin{proof}Choose a quotient map $p\colon P\to X$ from a ($\kappa$-light) profinite set $P$. 
This exists by \cref{compacthausdorffquotientoftdch}.
Then $\underline{X}_{(\kappa)}\cong \colim{\Delta^{\operatorname{op}}}\underline{\check{C}(p)}_{(\kappa)}$ (\cite[Proposition 6.2.3.20]{highertopostheory}), and hence \[ \cckH{(\kappa\text{-})}(\underline{X}_{\kappa}, \underline{A}_{\kappa})\cong \clim{\Delta}\left(\cckH{(\kappa\text{-})}(\underline{\check{C}({p})}_{(\kappa)}, \underline{A}_{(\kappa)})\right) \] and \[\icckH{(\kappa\text{-})}(\underline{X}_{(\kappa)}, \underline{A}_{(\kappa)})\cong  
        \clim{\Delta}\left(\icckH{(\kappa\text{-})}(\underline{\check{C}({p})}_{(\kappa)}, \underline{A}_{(\kappa)})\right)\] by \cref{spectralenrichmentcocontinuous} and \cref{Bousfieldkanspectralsequenceinternal}.
        By \cref{constantsheafandtotalizationcohomology}, \[c_{\Sp}\left(\cckH{\kappa\text{-}}(\underline{X}_{(\kappa)},\underline{A}_{(\kappa)})\right)\cong c_{\Sp}\left(\clim{\Delta}\cckH{(\kappa\text{-})}(\underline{\check{C}({p})}_{(\kappa)}, \underline{A}_{(\kappa)})\right)\cong \clim{\Delta}\left(c_{\Sp}\cckH{(\kappa\text{-})}(\underline{\check{C}({p})}_{(\kappa)},A)\right).\]

        As the \v{C}ech nerve $\check{C}(p)$ consists of ($\kappa$-light) profinite spaces, and $c_{\Sp}$ is $t$-exact (\cref{texactnessgeometricmorphism}), 
        \cref{internalcohomologyprofinite} implies that for all discrete abelian groups $A$, the counit for $c_{\Sp}\dashv \Gamma_{\Sp}$ induces an equivalence of cosimplicial spectra \[c_{\Sp}\left(\cckH{(\kappa\text{-})}(\underline{\check{C}({p})}_{(\kappa)}, \underline{A}_{(\kappa)})\right)\cong \icckH{(\kappa\text{-})}(\underline{\check{C}({p})}_{(\kappa)}, \underline{A}_{(\kappa)}).\qedhere\] 
\end{proof}
    \begin{rem}
        Internal condensed cohomology has a non-constant condensed structure in many cases: For a topological space $X$ and a topological abelian group $A$, 
        $c_{\Sp}\ckH^0(\underline{X}_{\kappa},\underline{A}_{\kappa})$ is the ($\kappa$-)condensed abelian group represented by the discrete abelian group $\ckH^0(X,A)=\mathcal C_{\kappa}(X,A)$. 
        On the other hand, $\ickH^0(\underline{X}_{\kappa}, \underline{A}_{\kappa})$ is represented by the topological abelian group $\mathcal C_{\kappa}(X,A)$ with compact open topology. 
        For example, if $X$ is an infinite set and $A$ is a compact abelian group, then 
        \[\ickH^0(X,A)\cong \prod_{X}\underline{A}_{(\kappa)}\] is represented by the compact topological group $\prod_{X}A$, whereas $c_{\Sp}\ckH^0(X,A)$ is represented by the discrete abelian group $(\prod_{X}A)^{\delta}$. 
    \end{rem}

\subsection{Condensed module categories}\label{section:condensedmodulecategories}
In this section, we record basic properties of derived and underived $(\kappa\text{-})$condensed module categories relevant to the discussion of solid modules in \cref{section:solidmodules}. Many of the statements below are special cases of the results established in \cref{section:modulecategories}, we spell them out here for the convenience of the reader. 

\begin{notation}
Endow $\Cond{(\kappa)}(\Ab)$ with the symmetric monoidal structure arising from the tensor product of abelian groups via \cref{constructionmonoidalstructureaccessiblesheaves}. 

For a $(\kappa\text{-})$condensed ring $\algebra{R}\in\Alg(\Cond{(\kappa)}(\Ab))$, denote by \[\Cond{(\kappa)}(\algebra{R})\coloneqq \LMod{\algebra{R}}{\Cond{(\kappa)}(\Ab)}\] the (underived) category of left $\algebra{R}$-modules.  
\end{notation}
\begin{rems}\label{remarkstensorproductcondensedmodules}
This has the following properties: 
    \begin{romanenum}
            \item \label{forgetfunctorpreserveslimitscolimitscondensedmodules} By \cref{closednessmonoidalstructureaccessiblesheaves}, the symmetric monoidal structure on $\Cond{(\kappa)}(\Ab)$ is closed. 
            
            As $\Cond{(\kappa)}(\Ab)$ has all small limits and colimits (\cref{condensedextremallydisconnected1}), this implies that for a ($\kappa$-)condensed ring $\algebra{R}$, $\Cond{(\kappa)}(\algebra{R})$ has small limits and colimits and the forget functor \[\Cond{(\kappa)}(\algebra{R})\to \Cond{(\kappa)}(\Ab)\] reflects small limits and colimits (\cref{forgetfreeadjunctionmodules}). 
    \item Moreover, if $\algebra{R}$ is commutative, then $\Cond{(\kappa)}(\algebra{R})$ inherits a closed symmetric monoidal structure such that the free $\algebra{R}$-module functor enhances to a symmetric monoidal functor. (\cref{symmetricmonoidalstructuremodulecategories,internalhommodules}).  
         \item For uncountable cardinals $\kappa\leq\lambda$, the left-exact left adjoints \[\Cond{\kappa}(\Set)\to \Cond{\lambda}(\Set)\to \Cond{}(\Set)\] induce symmetric monoidal left adjoints \[\Cond{\kappa}(\Ab)\to\Cond{\lambda}(\Ab)\to \Cond{}(\Ab).\] For $\algebra{R}\in \Alg(\Cond{\kappa}(\Ab))$, these induce left adjoints \[\Cond{\kappa}(\algebra{R})\to\Cond{\lambda}(\algebra{R})\to\Cond{}(\algebra{R}),\] which enhance to symmetric monoidal functors if $\algebra{R}$ is commutative (\cref{localisationinduceslocalisationonmodulecategories,symmetricmonoidalstructuremodulesnatural}). 
    If $\kappa,\lambda$ are regular or strong limit cardinals, then 
    \[\Cond{\kappa}(\algebra{R})\to\Cond{\lambda}(\algebra{R})\to\Cond{}(\algebra{R})\] are fully faithful. (\cref{leftKanextensionsfullyfaithfulcondensed,condensedstronglimitfullyfaithful,condensedextremallydisconnected1,naturalitylocalizationmodulecategories}).
    By \cref{filteredcolimitsmodules,condensedascolimit}, \[\Cond{}(\algebra{R})\cong \colim{\kappa}\Cond{\kappa}(\algebra{R}),\] where the colimit is computed in $\vlCat$ and runs over all small regular/strong limit cardinals $\kappa$ with $\algebra{R}\in\Alg(\Cond{\kappa}(\Ab))$. 
    \item For a regular cardinal $\kappa$, $\Cond{\kappa}(\algebra{R})\subseteq \Cond{}(\algebra{R})$ is closed under small colimits and $\kappa$-small limits. This follows from \cref{condensedextremallydisconnected1,forgetfreeadjunctionmodules}. 
\item For a commutative, ($\kappa$-)condensed ring $\algebra{R}$, the tensor product on $\Cond{(\kappa)}(\algebra{R})$ can be described as follows: 
Symmetric monoidality of the free functor \[\algebra{R}[-]\colon \Cond{(\kappa)}(\Ab)\to \Cond{(\kappa)}(\algebra{R})\] and cocontinuity of the symmetric monoidal structure on $\Cond{\kappa}(\algebra{R})$ imply that for $\kappa$-condensed $R$-modules $M,N$, $M\otimes_{\algebra{R}}N$ is the $\kappa$-condensed sheafification of \[\Pro(\Fin)_{\kappa}^{\operatorname{op}}\to \Ab, \,  T\mapsto M(T)\otimes_{\algebra{R}(T)}N(T)\] with the canonical condensed $\algebra{R}$-module structure.  

If $\mu$ is a regular cardinal with \[\algebra{R}\in \CAlg(\Cond{\mu}(\Ab))\subseteq \CAlg(\Cond{}(\Ab)),\] for regular cardinals $\kappa\geq \mu$, \[\Cond{\kappa}(\algebra{R})\hookrightarrow\Cond{}(\algebra{R})\] is symmetric monoidal with essential image \[\Cond{}(\algebra{R})\times_{\Cond{}(\Set)}\Cond{\kappa}(\Set), \] so we can always compute the tensor product of to condensed $R$-modules in $\Cond{\kappa}(\algebra{R})$ for a sufficiently large regular/strong limit cardinal $\kappa$. This implies that for $M,N\in\Cond{}(\algebra{R})$, 
$M\otimes_{\algebra{R}}N$ is the condensed sheafification (\cref{accessiblesheafification}) of 
\[\Pro(\Fin)^{\operatorname{op}}\to \Ab, \, T\mapsto M(T)\otimes_{\algebra{R}(T)}N(T), \] equipped with the canonical condensed $\algebra{R}$-module structure.
    \item Denote by $c\colon\Ab\to \Cond{(\kappa)}(\Ab)$ the constant sheaf functor. If $\algebra{R}=c\algebra{S}$ for some ring $S$, then \[\Cond{(\kappa)}(\LMod{S}{\Ab})\cong \Cond{(\kappa)}(\algebra{R})\] (\cref{modulesoverconstantrings}). If $S$ is commutative, this enhances to a symmetric monoidal equivalence where the right-hand side is endowed with the symmetric monoidal structure induced from $\LMod{S}{\Ab}^{\otimes}$ via \cref{constructionmonoidalstructureaccessiblesheaves}.  
    \item The symmetric monoidal structure of $\Sp$ induces a symmetric monoidal structure on $\Cond{(\kappa)}{\Sp}$ which is compatible with the $t$-structure (\cref{monoidalstructurecompatibletstructure}), whence ${\Cond{(\kappa)}(\Sp)}^{\heart}$ inherits a symmetric monoidal structure. The equivalence $\Cond{(\kappa)}(\Ab)\cong {\Cond{(\kappa)}(\Sp)}^{\heart}$ enhances to a symmetric monoidal functor by \cref{symmetricmonoidalstructureonspheartisinducedbyab}.
\end{romanenum} 
\end{rems}
\begin{lemma}[{\cite[Theorem 1.10]{Scholzecondensed}}]\label{condensedabeliangroupsgrothendieckaxioms}
For a ($\kappa$-)condensed ring $\algebra{R}\in \Alg(\Cond{(\kappa)}(\Ab))$, the category $\Cond{(\kappa)}(\algebra{R})$ is abelian and satisfies Grothendieck's axioms $AB 3$, $AB 3^*$, $ AB 4$, $ AB 5$ and $AB 6$, to wit: 
\begin{romanenum}
\item All small limits ($AB 3$) and small colimits ($AB 3^*$) exist,  
\item Small coproducts ($AB 4$) and small filtered colimits ($AB 5$) are exact.
\item ($AB 6$): For any small index set $I$ and $I$-indexed family of small filtered diagrams \[(F_i\colon J_i\to \BMod{\algebra{R}}{\infty}^{\heart})_{i \in I},\] 
\[ \colim{(j_i\in J_i)_i}\prod_{i\in I} F_{i}(j_i)\cong \prod_{i\in I}\colim{J_i}F_i.\]
\end{romanenum}
Moreover, $\kappa$-small products are exact in $\Cond{\kappa}(\algebra{R})$. 
If $\kappa$ is a strong limit cardinal, then small products in $\Cond{\kappa}(\algebra{R})$ are exact. Small products in $\Cond{}(\algebra{R})$ are exact.
\end{lemma}
\begin{proof} 
    Since the forget functor $\Cond{(\kappa)}(\algebra{R})\to \Cond{(\kappa)}(\Ab)$ is conservative and creates small limits and colimits, it suffices to prove the statements for $\Cond{(\kappa)}(\Ab)$. We first do this for $\Cond{\kappa}(\Ab)$. 
    The category $\Cond{\kappa}(\Ab)\cong \Shv_{\condo}(\Pro(\Fin)_{\kappa}, \Ab)$ has all small limits and colimits, and limits and filtered colimits can be computed pointwise: 
    For limits this is clear, for filtered colimits this holds by \cref{filteredcolimitsincondensedncategorycomputedpointwise}. Since the category of abelian groups $\Ab$ satisfies $AB 3$, $ AB 3^*$,  $ AB 4$, $AB 5$ and $AB 6$, it follows that $\Cond{\kappa}(\Ab)$ does as well.
    Next, we show that $\kappa$-small products are exact in $\Cond{\kappa}(\Ab)$. 
    Suppose that  $(X_i\to Y_i)_{i\in I}$ is a family of epimorphisms in $\Cond{\kappa}(\Ab)$ with $|I|<\kappa$. 
    Fix $T\in\Pro(\Fin)_{\kappa}$ and $(y_i)_{i\in I}\in\prod_{i\in I}Y_i(T)$. For all $i\in I$: Since $X_i\to Y_i$ is an epimorphism, there exists a cover 
    $\{q_{j}\colon T_{j}^i\to T\}_{j\in J_i}$ such that \[q_j^*(y_i)\in \im(X_i(T_j^i)\to Y_i(T_j^i))\text{ for all } j\in J_i.\] 
    \begin{center}
        \begin{tikzcd}
            X_i(T)\arrow[d]\arrow[r, hookrightarrow,"\prod_jq_j^*"] & \prod_{j\in J_i}X_i(T_j^i)\arrow[d]\\ 
            Y_i(T)\arrow[r, hookrightarrow,"\prod_jq_j^*"] & \prod_{j\in J_i} Y_i(T_j^i).
        \end{tikzcd}
    \end{center} Since the condensed topology is finitary, we can assume that $J_i$ is finite for all $i\in I$ and by taking the coproduct $T_i\coloneqq \sqcup_{j\in J_i}T_{j_i}$ that it consists of a single element $T_i$. 
    Let \[T^y=\{(t_i)_{i\in I}\in\prod_{i\in I}T_i\, | \, q_i(t_i)=q_k(t_k)\, \forall i,k\in I\}\subseteq \prod_{i\in I}T_i.\] 
    The space $T^y$ is profinite as closed subspace of the profinite set $\prod_{i\in I}T_i$ and since $|I|<\kappa$, \[\wt(T^y)\leq\max\{|I|, \wt(T_i)\}<\kappa.\] 
    The induced map $q_y\colon T^y\to T$ is onto as for all $i\in I$, $T_i\to T$ is onto. 
    We obtain a commutative diagram 
    \begin{center}
        \begin{tikzcd}
           \prod_{i\in I} X_i(T)\arrow[d]\arrow[r, hookrightarrow, "q_y^*"] & \prod_{i\in I}X_i(T^y)\arrow[d]\\ 
           \prod_{i\in I} Y_i(T)\arrow[r, hookrightarrow,"q_y^*"] & \prod_{i\in I}Y_i(T^y).
        \end{tikzcd}
    \end{center} 
    The pullbacks along $q_y$ factor as \[\prod_{i\in I}Y_i(T)\hookrightarrow\prod_{i\in I}Y_i(T_i)\to\prod_{i\in I}Y_i(T^y), \] respectively \[\prod_{i\in I}X_i(T_i)\hookrightarrow\prod_{i\in i}X_i(T_i)\to \prod_{i\in I}X_i(T^y), \] 
    whence \[q_y^*((y_i)_{i\in I})\in \im(\prod_{i\in I}X_i(T^y)\to \prod_{i\in I}Y_i(T^y)).\] 
    This shows that $\prod_{i\in I}X_i\to\prod_{i\in I}Y_i$ is an epimorphism of condensed abelian groups, which proves that $\kappa$-small products are exact. 
    \cref{condensedoncompactextremallydisconnected} implies that if $\kappa$ is a strong limit cardinal, then \[\Cond{\kappa}(\Ab)\cong \Ab(\tau_{\leq 0}\Cond{\kappa}(\an))\cong \mathcal P_{\Sigma}(\edCH{\kappa}, \Ab).\]
    As finite finite products commute with all limits and colimits in $\Ab$, small limits and colimits in $ \mathcal P_{\Sigma}(\edCH{\kappa}, \Ab)$ can be computed pointwise, which implies that arbitrary products are exact in $\Cond{\kappa}(\Ab)$ for strong limit cardinals $\kappa$. 

    For every uncountable regular cardinal $\kappa$, $\Cond{\kappa}(\Ab)\subseteq \Cond{}(\Ab)$ is closed under $\kappa$-small limits by \cref{condensedextremallydisconnected1}.
    As $\Cond{\kappa}(\Ab)\subseteq \Cond{}(\Ab)$ is also closed under small colimits and every small diagram $I\to \Cond{}(\Ab)$ factors over $\Cond{\kappa}(\Ab)$ for some uncountable regular cardinal $\kappa$, it follows from the above that products in $\Cond{}(\Ab)$ are exact and that $\Cond{}(\Ab)$ satisfies $AB 3$, $AB 3^*$, $ AB 4$, $ AB 5$ and $AB 6$.
\end{proof}
\begin{notation}
For a $(\kappa\text{-})$condensed ring $\algebra{R}\in\Alg(\Cond{(\kappa)}(\Ab))$, denote by 
\[\algebra{R}[-]\colon \Cond{(\kappa)}(\Set)\to \Cond{(\kappa)}(\algebra{R})\] the left adjoint of the forget functor. 
\end{notation}\cref{compactprojectiverepresentables}, \cref{compactprojectivegeneratorsforcondensedsets} and \cref{cantorsetgenerates} imply: 
\begin{cor}\label{enoughprojectives}
    Suppose $\algebra{R}\in\Alg(\Cond{(\kappa)}(\Ab))$ is a $(\kappa)$-condensed ring. 
    \begin{romanenum}
    \item For $K\in\CH_{(\kappa)}$, $\algebra{R}[\underline{K}_{(\kappa)}]$ is compact in $\Cond{(\kappa)}(\algebra{R})$. 
    \item For $X\in\edCH{\kappa}$, $\algebra{R}[\underline{X}_{(\kappa)}]$ is compact projective in $\Cond{(\kappa)}(\algebra{R})$. 
    \item If $\kappa$ is a strong limit cardinal, then $\Cond{\kappa}(\algebra{R})$ has enough projectives. 
    \item For every condensed $\algebra{R}$-module $M$, there exists a small family $(X_i)_{i\in I}$ of extremally disconnected compact Hausdorff spaces with a quotient map $\oplus_{i\in I}\algebra{R}[\underline{X_i}]\to M$. In particular, $\Cond{}(\algebra{R})$ has enough projectives. 
    \end{romanenum}
\end{cor} 
\begin{proof}
    As the forget functor $\Cond{(\kappa)}(\algebra{R})\to\Cond{(\kappa)}(\Ab)\to \Cond{(\kappa)}(\Set)$ commutes with filtered colimits (\cref{remarkstensorproductcondensedmodules,filteredcolimitsincondensedncategorycomputedpointwise}), its left adjoint $\algebra{R}[-]$ preserves compact objects. 
    In particular, for a ($\kappa$-light) compact Hausdorff space $K$, $\algebra{R}[\underline{K}_{(\kappa)}]$ is compact in $\Cond{(\kappa)}(\algebra{R})$ by \cref{compactprojectiverepresentables}/\cref{compactprojectivegeneratorsforcondensedsets}.

    Next, we show that for a ($\kappa$-light) extremally disconnected compact Hausdorff space $X$, $\algebra{R}[\underline{X}_{(\kappa)}]$ is projective in $\Cond{(\kappa)}(\algebra{R})$.
    For a $\kappa$-light extremally disconnected compact Hausdorff space $X$, 
    \[ \Hom_{\Cond{\kappa}(\algebra{R})}(R[\underline{X}_{\kappa}],-)\] factors as 
\[ \Cond{(\kappa)}(\algebra{R})\xrightarrow{f}\Cond{(\kappa)}(\Ab)\xrightarrow{r^*}\mathcal P_{\Sigma}(\edCH{\kappa},\Ab)\xrightarrow{\operatorname{ev}_X}\Ab,\] where $r^*$ is the restriction along $\edCH{\kappa}\subseteq \Pro(\Fin)_{\kappa}$. \cref{condensedoncompactextremallydisconnected} implies that $r^*$ has a right adjoint given by right Kan extension, whence $r^*$ preserves colimits. 
As $f$ and $\operatorname{ev}_X$ are exact, it follows that $\algebra{R}[\underline{X}_{\kappa}]$ is projective in $\Sol{(\kappa)}(\algebra{R})$. 
Suppose now that $X$ is extremally disconnected compact Hausdorff and $0\to A\to B\to C\to 0$ is a short exact sequence of condensed $\algebra{R}$-modules. Choose a regular cardinal $\kappa$ such that $R\in\Alg(\Cond{\kappa}(\Ab))$, $A,B,C\in\Cond{\kappa}(\algebra{R})$ and $\wt(X)<\kappa$.
Denote by $r^{\kappa}$ the right adjoint of $i_{\kappa}\colon \Cond{\kappa}(\algebra{R})\subseteq \Cond{}(\algebra{R})$. As $i_{\kappa}$ is fully faithful and exact, \[ 0\to r^{\kappa}A\to r^{\kappa}B\to r^{\kappa}C\to 0\] is exact. 
By \cref{topspacecomefromregular}, $i_{\kappa}\algebra{R}[\underline{X}_{\kappa}]=\algebra{R}[\underline{X}]$. It now follows from the above that \[0\to \Hom_{\Cond{}(\algebra{R})}(\algebra{R}[\underline{X}],A)\to \Hom_{\Cond{}(\algebra{R})}(\algebra{R}[\underline{X}],B)\to \Hom_{\Cond{}(\algebra{R})}(\algebra{R}[\underline{X}],C)\to 0\] is exact, which shows that $R[\underline{X}]$ is projective. 

Suppose now that $\kappa$ is a strong limit cardinal. As the forget functor $\Cond{\kappa}(\algebra{R})\to \Cond{\kappa}(\Set)$ is faithful, \cref{compactprojectivegeneratorsforcondensedsets} implies that for every $\kappa$-condensed $\algebra{R}$-module $M$, the canonical map \[\oplus_{S\in{\edCH{\kappa}}{/M}}\algebra{R}[\underline{X}_{\kappa}]\to M\] is an epimorphism in $\Cond{\kappa}(\algebra{R})$, where we consider $M$ as a $\kappa$-condensed set in the indexing category. As $\edCH{\kappa}$ is small and for all $S\in\edCH{\kappa}$, $M(S)$ is small, $\edCH{\kappa/M}$ is small. 
This shows in particular that $\Cond{\kappa}(\algebra{R})$ has enough projectives for strong limit cardinals $\kappa$. 

Suppose now that $M\in\Cond{}(\algebra{R})$ and choose a strong limit cardinal $\kappa$ with $\algebra{R}\in\Alg(\Cond{\kappa}(\Ab))$ and $M\in\Cond{\kappa}(\algebra{R})\subseteq \Cond{}(\algebra{R})$. As $\Cond{\kappa}(\algebra{R})\subseteq \Cond{}(\algebra{R})$ is a left adjoint, it follows from the above and \cref{topspacecomefromregular} that $\oplus_{S\in{\edCH{\kappa}}{/M}}\algebra{R}[\underline{X}]\to M$ is an epimorphism in $\Cond{}(\algebra{R})$. In particular, $\Cond{}(\algebra{R})$ has enough projectives.  
\end{proof}

\begin{rems}\label{condensedmodulesgrothendieckabelian}
    \begin{enumerate}\item As the forget functor $\Cond{}(\algebra{R})\to\Cond{}(\Set)$ is conservative, for all uncountable cardinals $\kappa$, $\oplus_{X\in\Pro(\Fin)_{\kappa}}\algebra{R}[\underline{X}_{\kappa}]$ is a generator for $\Cond{\kappa}(\algebra{R})$. 
    Together with \cref{condensedabeliangroupsgrothendieckaxioms}, this implies that $\Cond{\kappa}(\algebra{R})$ is Grothendieck abelian for all uncountable cardinals $\kappa$. 
    
   \item The category $\Cond{}(\algebra{R})$ is not Grothendieck abelian (if $\algebra{R}\neq 0$) since it does not admit a generator. 
The generating family $\{\algebra{R}[\underline{X}],X\in\edCH{}\}$ is large and the direct sum 
$\oplus_{X\in\edCH{}}\algebra{R}[\underline{X}]$ does NOT exist in $\Cond{}({\algebra{R}})$ (unless $\algebra{R}=0$). 

\item The author is unaware whether for arbitrary uncountable cardinals $\kappa$, $\Cond{\kappa}(\algebra{R})$ admits enough projectives, but expects that this is not the case.  

\item If $\topo{X}$ is a topos, the category $\Ab(\tau_{\leq 0}\topo{X})$ of abelian group objects in the underlying 1-topos has enough injectives, cf.\ \cite[\href{https://stacks.math.columbia.edu/tag/01DP}{Tag 01DP}]{stacks-project}. 
In contrast to that, the category of condensed abelian groups has no non-zero injectives, cf.\  \cite{noinjectivesincondensed}.
    \end{enumerate} 
\end{rems}
We now collect some basic properties about derived/stable condensed module categories which follow from our discussion in \cref{section:modulecategories}.
\cref{derivedcategorymodules} implies: 
\begin{cor}\label{condensedmodulesarederivedcats}
    If $\algebra{R}\in\Alg(\Cond{(\kappa)}(\Ab))$ is a $(\kappa\text{-})$condensed ring, \[\Cond{(\kappa)}(\algebra{R})\cong\CMod{\algebra{R}}{(\kappa)}^{\heart}\subseteq \CMod{\algebra{R}}{(\kappa)}\] extends to a $t$-exact equivalence 
    \[\mathcal D(\Cond{(\kappa)}(\algebra{R}))\cong \LMod{\algebra{R}}{\Cond{(\kappa)}(\Sp)},\] where on the right-hand side we consider $\algebra{R}$ as an algebra in $\Cond{(\kappa)}(\Sp)$ via the Eilenberg-Mac Lane functor (\cref{Eilenbergmaclanealgebra}).
\end{cor}
This is instructive as by our discussion from \cref{section:modulecategories}, the right-hand side has good categorical properties. 
In particular, we have shown the following: 
\begin{cor}\label{derivedcondensedmodulesallsmallimitsandcolimits}
For $\algebra{R}\in\Alg(\Cond{(\kappa)}(\Ab))$, the category $\mathcal D(\Cond{(\kappa)}(\algebra{R}))$ is big presentable and has small limits and colimits.
\end{cor}
\begin{proof}\cref{condensedmodulesarederivedcats} and \cref{filteredcolimitsmodules} imply that $\mathcal D(\Cond{(\kappa)}(\algebra{R}))$ is big presentable. 
As $\Cond{(\kappa)}(\Sp)$ has small colimits and limits (\cref{condensedextremallydisconnected1}) and the symmetric monoidal structure on $\Cond{(\kappa)}(\Sp)$ is closed (\cref{symmetricmonoidalstructureaccessiblesheavesofspectraclosed}), (co)com\-ple\-te\-ness follows from \cref{forgetfreeadjunctionmodules} and the identification of \cref{condensedmodulesarederivedcats}. 
\end{proof}
\begin{cor}\label{condensedmodulesaremodulesincondensed}
    Denote by $c_{\Sp}\colon \Alg(\Sp)\to \Alg(\Cond{(\kappa)}(\Sp))$ the functor induced by the constant sheaf functor. 
    \begin{romanenum}
    \item For $\algebra{R}\in\Alg(\Sp)$, \[ \Cond{(\kappa)}(\LMod{\algebra{R}}{\Sp})\cong \CMod{c_{\Sp}\algebra{R}}{(\kappa)}.\]

    \item In particular, for a ring $\algebra{R}\in\Alg(\Ab)$, 
    \[  \Cond{(\kappa)}(\mathcal D(\algebra{R}))\cong \CMod{c_{\Sp}\algebra{R}}{(\kappa)}\cong \mathcal D(\Cond{(\kappa)}(c_{\Sp}\algebra{R})).\]
    \item For a commutative ring spectrum $\algebra{R}\in \CAlg(\Sp)$, the above equivalence enhances to a symmetric monoidal equivalence 
\[ \Cond{(\kappa)}(\LMod{\algebra{R}}{\Sp})\cong \CMod{c_{\Sp}\algebra{R}}{(\kappa)},\] 
where the left-hand side is endowed with the symmetric monoidal structure induced from $\LMod{\Sp}{\algebra{R}}$ via \cref{constructionmonoidalstructureaccessiblesheaves}, and the right-hand side with the symmetric monoidal structure provided by \cref{symmetricmonoidalstructurebigtoposmodules}. 
    \end{romanenum}
\end{cor}
\begin{proof}The first identification is \cref{modulesoverconstantrings}, the second statement follows from the first and \cref{derivedcategorymodules,condensedmodulesarederivedcats}, and third identification holds by \cref{modulesoverconstantringssymmetricmonoidalaccessibleshaves}. 
\end{proof}
 \cref{modulesaccessiblesheavesareclosed} implies: 
\begin{cor}\label{condensedmodulesclosed}
    If $\algebra{R}\in\CAlg(\Cond{(\kappa)}(\Sp))$ is a commutative condensed ring spectrum, the symmetric monoidal structure on $\LMod{\algebra{R}}{\Cond{(\kappa)}{(\Sp)}}$ described in \cref{symmetricmonoidalstructurebigtoposmodules} is closed. 
\end{cor}
\begin{proof}For $\Cond{}(\Sp)$, this is \cref{modulesaccessiblesheavesareclosed}.
    For $\Cond{\kappa}(\Sp)$, this follows from \cite[Theorem 3.4.4.2]{higheralgebra}. 
\end{proof}

\begin{notation}\label{internalhomandenrichmentderivedcondensedcategories}
    For a commutative condensed ring $\algebra{R}\in\CAlg(\Cond{(\kappa)}(\Ab))$, endow $\mathcal D(\Cond{(\kappa)}(\algebra{R}))$ with the closed symmetric monoidal structure induced from the symmetric monoidal structure on $\LMod{\algebra{R}}{\Cond{(\kappa)}(\Sp)}$ from \cref{condensedmodulesclosed} via the equivalence \[\mathcal D(\Cond{(\kappa)}(\algebra{R}))\cong \LMod{\algebra{R}}{\Cond{(\kappa)}(\Sp)}\] from \cref{condensedmodulesarederivedcats}. Denote by $\imap_{\mathcal D(\Cond{(\kappa)}(\algebra{R}))}(-,-)$ the internal Hom. 
    By \cref{modulecategoriesareenriched}, \[\mathcal D(\Cond{(\kappa)}(\algebra{R}))\cong \LMod{\algebra{R}}{\Cond{(\kappa)}(\Sp)}\] is enriched in $\mathcal D(\algebra{R}(*))$. We denote by $\map_{\mathcal D(\Cond{(\kappa)}(\algebra{R}))}(-,-)$ the enriched mapping functor. 
\end{notation}
\begin{rems}\label{enrichment}

    \begin{enumerate}
\item Denote by $c\colon\CAlg(\Ab)\to\CAlg(\Cond{}(\Ab))$ the functor induced by the constant sheaf functor. By \cref{condensedmodulesaremodulesincondensed}, for a ring $\algebra{R}$, the symmetric monoidal structure on $\mathcal D(\Cond{(\kappa)}(c\algebra{R}))$ from above is equivalent to the one induced by the symmetric monoidal structure on $\mathcal D(\algebra{R})\cong \LMod{\algebra{R}}{\Sp}$.
    In particular, for $M,N\in\mathcal D(\Cond{(\kappa)}(\algebra{R}))$, 
$M\otimes_{\mathcal D(\Cond{(\kappa)}(\algebra{R}))}N$ is the ($\kappa$-)condensed sheafification of 
\[ \Pro(\Fin)_{(\kappa)}^{\operatorname{op}}\ni K\mapsto M(K)\otimes_{\mathcal D(\algebra{R})}N(K).\] 
\item By \cref{monoidalstructurecompatibletstructure}, the symmetric monoidal structure on $\mathcal D(\Cond{(\kappa)}(\algebra{R}))$ is compatible with the $t$-structure. The induced symmetric monoidal structure on the heart $\Cond{(\kappa)}(\algebra{R})$ is the one induced by the tensor product of abelian groups, cf.\ \cref{symmetricmonoidalstructureonspheartisinducedbyab}. 
    \end{enumerate}
\end{rems}
Recall from \cref{condensedmodulesarederivedcats} that $\mathcal D(\Cond{(\kappa)}(\Ab))\cong \LMod{H\mathbb Z}{\Cond{(\kappa)}(\Sp)}$. 
    By \cite[Corollary 3.4.1.9, Corollary 4.5.1.6]{higheralgebra}, this implies that for all algebras $R\in\Alg(\Cond{(\kappa)}(\Ab))$, 
    \[\LMod{\algebra{R}}{\mathcal D(\Cond{(\kappa)}(\Ab))}\cong \LMod{\algebra{R}}{\Cond{(\kappa)}{\Sp}}.\]  
    In particular, $\LMod{\algebra{R}}{\mathcal D(\Cond{(\kappa)}(\Ab))}$ inherits a $t$-structure via \cref{tstructureonmodules}. \cref{condensedmodulesarederivedcats} now implies: 
\begin{cor}\label{modulesinderivedcategoryisderivedcategoryofmodules} For $R\in \Alg(\Cond{(\kappa)}(\Ab))$, the inclusion \[\Cond{(\kappa)}(\algebra{R})\cong  \LMod{R}{\mathcal D(\Cond{(\kappa)}(\Ab))}^{\heart}\hookrightarrow \LMod{R}{\mathcal D(\Cond{(\kappa)}(\Ab))}\] extends to an equivalence \[\mathcal D({\Cond{(\kappa)}(\algebra{R})})\cong \LMod{R}{\mathcal D(\Cond{(\kappa)}(\Ab))}.\] 
\end{cor}
\cref{condensedanimaPostnikovcomplete} implies the following: 
\begin{cor}\label{tstructureleftcomplete}If $\algebra{R}\in\Alg(\Cond{(\kappa)}(\Sp)_{\geq 0})$ is a connective $(\kappa\text{-})$condensed ring, the $t$-structure on $\LMod{\algebra{R}}{\Cond{(\kappa)}(\Sp)}$ from \cref{tstructureonmodules} is left and right complete. 
\end{cor}
\begin{proof}
    \cite[Proposition 1.3.2.7]{SAG} implies that for all uncountable cardinals $\kappa$, \[\Cond{\kappa}(\Sp)_{\leq 0}\subseteq \Cond{\kappa}(\Sp)\] is closed under countable coproducts and that the $t$-structure on $\Cond{\kappa}(\Sp)$ is right-separated. It follows from \cref{tstructurespectrumobjects,condensedextremallydisconnected1,tstructurespectrumobjects} that this also holds for $\Cond{}(\Sp)$.
    As the forget functor \[\LMod{\algebra{R}}{\Cond{(\kappa)}(\Sp)}\to \Cond{(\kappa)}(\Sp)\] is $t$-exact, conservative and reflects colimits (\cref{forgetfreeadjunctionmodules}), it follows from \cite[Proposition 1.2.1.19]{higheralgebra} that the $t$-structure on $\LMod{\algebra{R}}{\Cond{(\kappa)}(\Sp)}$ is right-complete. 

    To show left-completeness, we proceed as in \cite[Corollary 3.26]{mondal2024postnikovcompletenessrepletetopoi}. 
    As $\LMod{\algebra{R}}{\Cond{(\kappa)}(\Sp)}$ has small limits, the functor 
    \[\clim{n\in\mathbb Z}\tau_{\leq n}-\colon \LMod{\algebra{R}}{\Cond{(\kappa)}(\Sp)}\to \clim{n\in\mathbb Z}\left(\LMod{\algebra{R}}{\Cond{(\kappa)}(\Sp)}_{\leq n}\right), X\mapsto (\tau_{\leq n}X)\] has a right adjoint 
    \[ \clim{}\colon \clim{n\in\mathbb Z}\left(\LMod{\algebra{R}}{\Cond{(\kappa)}(\Sp)}_{\leq n}\right)\to \LMod{\algebra{R}}{\Cond{(\kappa)}(\Sp)}, (X_n)_n\mapsto \clim{n\in\mathbb Z}X_n\] by   
    \cite[Proposition 5.1.10]{landinfinity}. 
    The unit and counit of this adjunction are the canonical maps \[\eta\colon X\to \clim{n\in\mathbb Z}\tau_{\leq n}X\text{ and } (\epsilon_n\colon (\tau_{\leq n}\clim{k\in\mathbb Z}X_k\to \tau_{\leq n}X_n=X_n))_{n\in\mathbb Z},\] respectively. 
    So it remains to show that these maps are equivalences, then left-completeness follows. 
    By \cref{tstructureonmodules,forgetfreeadjunctionmodules} it suffices to check this for $\Cond{(\kappa)}(\Sp)$. 
    By right separatedness of the $t$-structure and \cref{recognitionprinciplebigtopoi}, the functors \[\Omega^{\infty-m}\colon\Cond{(\kappa)}(\Sp)\to \Cond{(\kappa)}(\an), m\in\mathbb N_0\] are jointly conservative. 
    As $\Omega^{\infty-m}$ preserves limits and $\Omega^{\infty-m}\tau_{\leq n}=\tau_{\leq n+m}\Omega^{\infty-m}$, \[\Omega^{\infty-m}(\clim{n}\tau_{\leq n}X)\cong (\clim{n}\Omega^{\infty-m}\tau_{\leq n}X)\cong \clim{n}\tau_{\leq n+m}\Omega^{\infty-m}X\cong \clim{n}\tau_{\leq n}\Omega^{\infty-m}X.\]
    Since Postnikov towers in $\Cond{(\kappa)}(\an)$ converge (\cref{condensedanimaPostnikovcomplete}), this implies that \[\Omega^{\infty-m}\eta\colon \Omega^{\infty-m}X\to \Omega^{\infty-m}\clim{n}\tau_{\leq n}X\] is an equivalence for all $m\in\mathbb N_0$.  
    Suppose now that $(X_k)_k\in \clim{k\in\mathbb Z}\Cond{(\kappa)}(\Sp)_{\leq k}$. Then \[(\Omega^{\infty-m}X_{k-m})_k\in \clim{k\in\mathbb Z}\Cond{(\kappa)}(\an)_{\leq k}.\] 
    As $\lim{}\colon \clim{k\in\mathbb Z}\Cond{(\kappa)}(\an)_{\leq k}\to \Cond{(\kappa)}(\an)$ is an equivalence (\cref{condensedanimaPostnikovcomplete}), for all $k\in\mathbb Z$, \[\Omega^{\infty-m}X_{k-m}\cong \tau_{\leq k}\clim{n}(\Omega^{\infty-m}X_{n-m})\cong \Omega^{\infty-m}\tau_{\leq k-m}\clim{n}X_n.\] This shows that $\Omega^{\infty-m}\epsilon$ is an equivalence for all $m\in\mathbb N_0$. 
\end{proof}

\subsection{Compact condensed abelian groups}\label{section:compactcondensedabeliangroups}
In this section, we collect some compactness properties of condensed abelian groups from \cite[Lecture 4]{Scholzecondensed} which we will use in our discussion of solid abelian groups in the next section. 
Throughout, we use \cref{internalhomandenrichmentderivedcondensedcategories}. 
\begin{lemma}\label{suspensionspectraonchcompact}
    Fix $n\in\mathbb N_0$. 
\begin{romanenum}
    \item Suppose $\kappa$ is an uncountable cardinal. For a $\kappa$-light compact Hausdorff space $S$, 
    \[ \ickH(\underline{S}_{\kappa},-)\colon \Cond{\kappa}(\Sp)_{\leq n}\to\Cond{\kappa}(\Sp)\] and 
    \[ \ckH(\underline{S}_{\kappa},-)\colon \Cond{\kappa}(\Sp)_{\leq n}\to \Sp\] commute with filtered colimits. 

    \item For a compact Hausdorff space $S$,  
    \[ \iccH(\underline{S},-)\colon \Cond{}(\Sp)_{\leq n}\to\Cond{}(\Sp)\] and 
    \[ \ccH(\underline{S},-)\colon \Cond{}(\Sp)_{\leq n}\to \Sp\] commute with filtered colimits.
    \end{romanenum}
\end{lemma}
\begin{proof}
    It is enough to prove the statement for internal cohomology, then the statement for cohomology follows from \cref{globalsectionstexactcocontinuous}. By \cref{infinitelooppreservesfilteredcolimits}, \[\Omega^{\infty-i}\colon\Cond{(\kappa)}(\Sp)\to \Cond{(\kappa)}(\an)\] preserves filtered colimits for all $i\in\mathbb N_1$. As the functors $\Omega^{\infty-i}, i\in\mathbb N_0$ are jointly conservative (\cite[Remark 1.4.2.25]{higheralgebra}),
it suffices to show that \[\Omega^{\infty-i}\icckH{(\kappa\text{-})}(\underline{S}_{(\kappa)},-)|_{\Cond{(\kappa)}(\Sp)_{\leq n}}\cong\iMap_{\Cond{(\kappa)}(\an)}(\underline{S}_{(\kappa)}, \Omega^{\infty}\Sigma^i-)|_{\Cond{(\kappa)}(\Sp)_{\leq n}}\] preserves filtered colimits for all $i\in\mathbb N_0$. The equivalence of these functors was explained in \cref{propertiesinternalhomspectra}. 
By \cref{condensedextremallydisconnected1,tstructurespectrumobjects}, it suffices to show this for $\Cond{\kappa}(\Sp)$. 
As $\Omega^{\infty}\Sigma^i\colon \Cond{\kappa}(\Sp)_{\leq n}\to \Cond{\kappa}(\an)_{\leq i+n}$ preserves filtered colimits and \[\iMap_{\Cond{(\kappa)}(\an)}(\underline{S}_{\kappa},-)=(\Pro(\Fin)_{\kappa}^{\operatorname{op}}\ni T\mapsto \Map_{\Cond{\kappa}(\an)}(\underline{S}_{\kappa}\times \underline{T}_{\kappa},-))\] (\cref{enrichmentrecoversinternalhom}), the statement now follows from \cref{compactprojectiverepresentables,filteredcolimitsincondensedncategorycomputedpointwise}. 
\end{proof}
\begin{cor}\label{freemodulesonchpseudocoherent}
    Suppose $\kappa$ is an uncountable cardinal and $R\in \CAlg(\Cond{\kappa}(\Sp))$ is a commutative condensed ring spectrum. 
    For $n\in\mathbb N_0$ and $S\in \CH_{(\kappa)}$, 
    \[ \imap_{\CMod{\algebra{R}}{(\kappa)}}(\algebra{R}[\underline{S}_{\kappa}],-)\colon \CMod{\algebra{R}}{(\kappa)}_{\leq n}\to\CMod{\algebra{R}}{(\kappa)}\] commutes with filtered colimits. 
\end{cor}
\begin{proof}
    Since the symmetric monoidal structure on $\Cond{(\kappa)}(\Sp)$ is compatible with small colimits, the forget functor $f\colon\CMod{\algebra{R}}{(\kappa)}\to \Cond{(\kappa)}(\Sp)$ reflects colimits by \cref{forgetfreeadjunctionmodules}.
    As \[f\circ \imap_{\CMod{\algebra{R}}{(\kappa)}}(\algebra{R}[\underline{K}_{(\kappa)}],-)\cong \imap_{\Cond{(\kappa)}(\Sp)}(\Sigma^{\infty}_{+}\underline{K}_{(\kappa)},-)\circ f\] (both are right adjoint to $\algebra{R}[\underline{K}_{(\kappa)}]\otimes_{\algebra{R}}\algebra{R}[-]\cong \algebra{R}[-]\circ (\underline{K}_{(\kappa)}\times-)$) and $f$ is $t$-exact,
    \[ f\circ \imap_{\CMod{\algebra{R}}{(\kappa)}}(\algebra{R}[\underline{K}_{(\kappa)}],-)\colon \CMod{\algebra{R}}{(\kappa)}_{\leq n}\to \Cond{(\kappa)}(\Sp)\] preserves filtered colimits for all $n\in\mathbb N_0$ by  
    \cref{suspensionspectraonchcompact}, which implies that \[ \imap_{\CMod{\algebra{R}}{(\kappa)}}(\algebra{R}[\underline{K}_{(\kappa)}],-)\colon \CMod{\algebra{R}}{(\kappa)}_{\leq n}\to \CMod{\algebra{R}}{(\kappa)}\] does as well. 
\end{proof}

\paragraph{Some computations using Breen-Deligne resolutions}
We end this section with a recollection of computations using Breen-Deligne resolutions from \cite[Lecture 4]{Scholzecondensed}. 
\begin{cor}\label{freemodulespseudocoherentdabenriched}
    Suppose $\algebra{R}\in \CAlg(\Cond{(\kappa)}(\Ab))$. 
    For all $S\in\CH_{(\kappa)}$ and $n\in\mathbb N_0$, \[\map_{\mathcal D(\Cond{(\kappa)}(\algebra{R}))}(\algebra{R}[\underline{S}_{(\kappa)}],-)\colon \mathcal D(\Cond{(\kappa)}(\algebra{R}))_{\leq n}\to \mathcal D(\algebra{R}(*))\] commutes with filtered colimits. 
\end{cor}
\begin{proof}
By \cref{globalsectionstexactcocontinuous} and \cref{forgetfreeadjunctionmodules}, the global sections functor \[\Gamma_{\algebra{R}}\colon\mathcal D(\Cond{(\kappa)}(\algebra{R}))\to\mathcal D(\Cond{}(\algebra{R}(*)))\] preserves small colimits. 
 As \[\map_{\mathcal D(\Cond{(\kappa)}(\algebra{R}))}(\algebra{R}[\underline{K}_{(\kappa)}],-)=\Gamma_{\algebra{R}}\imap_{\mathcal D(\Cond{(\kappa)}(\algebra{R}))}(\algebra{R}[\underline{K}_{(\kappa)}],-), \] the statement now follows from \cref{freemodulesonchpseudocoherent}.\end{proof}
Together with \cref{Breendelignespectralsequence}, this implies: 
\begin{cor}\label{compacttopologicalgroupspseudocoherent}
Suppose $A$ is a ($\kappa$-light) compact Hausdorff topological abelian group and $n\in\mathbb Z$. The internal Hom (\ref{internalhomandenrichmentderivedcondensedcategories})
\[ \imap_{\mathcal D(\Cond{(\kappa)}(\Ab))}(\underline{A}_{(\kappa)},-)\colon \mathcal D(\Cond{(\kappa)}(\Ab))_{\leq n}\to \mathcal D(\Cond{(\kappa)}(\Ab))\] and the $\mathcal D(\Ab)$-enriched mapping functor (\cref{internalhomandenrichmentderivedcondensedcategories})\[ \map_{\mathcal D(\Cond{(\kappa)}(\Ab))}(\underline{A}_{(\kappa)},-)\colon \mathcal D(\Cond{(\kappa)}(\Ab))_{\leq n}\to \mathcal D(\Ab)\]commute with filtered colimits.
\end{cor}

\begin{proof}We first prove the $\kappa$-condensed statement. Under the equivalence \[\mathcal D(\Cond{\kappa}(\Ab))\cong\Cond{\kappa}(\mathcal D(\Ab))\] provided by \cref{condensedmodulesaremodulesincondensed}, \[\imap_{\mathcal D(\Cond{\kappa}(\Ab))}(\underline{A}_{(\kappa)},-)=(\Pro(\Fin)_{\kappa}^{\operatorname{op}}\ni T\mapsto \map_{\mathcal D(\Cond{\kappa}(\Ab))}(A\otimes \mathbb Z[\underline{T}_{\kappa}],-)), \] cf.\ \cref{enrichmentrecoversinternalhom}. Here $A\otimes \mathbb Z[\underline{T}_{\kappa}]$ denotes the tensor product in $\Cond{\kappa}(\Ab)$, which agrees with the tensor product in $\mathcal D(\Cond{\kappa}(\Ab))$ since $\mathbb Z[\underline{T}_{\kappa}]$ is flat. 
    It therefore suffices to show prove the statement for $\map_{\mathcal D(\Cond{\kappa}(\Ab))}(\underline{A}_{\kappa},-)$, then the statement for $\imap_{\mathcal D(\Cond{\kappa}(\Ab))}(\underline{A}_{\kappa},-)$ follows. 
    \cref{Breendelignespectralsequence} provides a convergent spectral sequence of functors \[\mathcal D(\Cond{\kappa}(\Ab))_{\leq N}\times \CH_{\kappa}^{\operatorname{op}}\to\mathcal D(\Ab),\]\[ E_1^{p,q}(B,T)=\prod_{k=1}^{b_p}\pi_{-q}\map_{\mathcal D(\Cond{\kappa}(\Ab))}(\mathbb Z[\underline{A^{i_k}\times T}_{\kappa}],B)\Rightarrow \pi_{-(p+q)}\map_{\mathcal D(\Cond{\kappa}(\Ab))}(\underline{A}_{\kappa}\otimes\mathbb Z[\underline{T}_{\kappa}],B).\]
     As for all $n\in\mathbb N_0$, $\sqcup_{k=1}^{b_n} A^{i_k,n}\times T\in \CH_{\kappa}$, \cref{freemodulespseudocoherentdabenriched} implies that the $E_1$-page preserves filtered colimits in the variable from $\mathcal D(\Cond{\kappa}(\Ab))_{\leq N}$, hence the limit \[\pi_{-i}\map_{\mathcal D(\Cond{\kappa}(\Ab))}(A\otimes\mathbb Z[\underline{T}_{\kappa}],-)\colon \mathcal D(\Cond{\kappa}(\Ab))_{\leq N}\to\Ab\] of the spectral sequence commutes with filtered colimits for all $T\in\CH_{\kappa}$ and all $i\in\mathbb Z$. 
    As the functors \[\pi_{-i}\colon\mathcal D(\Ab)\to \Ab,i\in\mathbb Z\] are jointly conservative and commute with filtered colimits, this implies the statement for $\Cond{\kappa}(\Ab)$. 

    We now deduce the statement for $\Cond{}(\Ab)$ from this. Fix $n\in\mathbb N_0$ and suppose that \[F\colon I\to \mathcal D(\Cond{}(\Ab))_{\leq n}\] is a small filtered diagram. Choose a regular cardinal $\kappa\geq |A|$ such that $F$ factors over \[\mathcal D(\Cond{\kappa}(\Ab))_{\leq n}\subseteq \mathcal D(\Cond{}(\Ab))_{\leq n}.\] 
    By \cref{forkappalargeenoughclosedmonoidal}, there exists a regular cardinal $\lambda\geq \kappa$ such that for all $i\in I$, \begin{align*}\imap_{\mathcal D(\Cond{}(\Ab))}(A,F(i))& \cong \imap_{\mathcal D(\Cond{\lambda}(\Ab))}(A,F(i)), \\ \text{ and }\imap_{\mathcal D(\Cond{}(\Ab))}(A, \colim{I}F)& \cong \imap_{\mathcal D(\Cond{\lambda}(\Ab))}(A, \colim{I}F)\end{align*} in $\mathcal D(\Cond{\lambda}(\Ab))\subseteq \mathcal D(\Cond{}(\Ab))$.  
    Since $\mathcal D(\Cond{\lambda}(\Ab))\subseteq \mathcal D(\Cond{}(\Ab))$ is closed under small colimits, it follows from the above that 
    \[ \colim{I}\, \imap_{\mathcal D(\Cond{}(\Ab))}(A,F)\cong \imap_{\mathcal D(\Cond{}(\Ab))}(A, \colim{I}F).\] By \cref{globalsectionstexactcocontinuous} and \cref{forgetfreeadjunctionmodules}, the global sections functor $\Gamma\colon \mathcal D(\Cond{}(\Ab))\to\mathcal D(\Ab)$ preserves colimits. 
    It now follows that \[\map_{\mathcal D(\Cond{}(\Ab))}(A,-)|_{\mathcal D(\Cond{}(\Ab))_{\leq n}}=\Gamma \imap_{\mathcal D(\Cond{}(\Ab))}(A,-)|_{\mathcal D(\Cond{}(\Ab))_{\leq n}}\] preserves filtered colimits for all $n\in\mathbb N_0$. 
\end{proof}
\begin{lemma}\label{cohomologyofrwithdiscretecoefficientsvanishes}
Suppose $M$ is a discrete abelian group. 
For $A\in \mathcal D(\Cond{(\kappa)}(\underline{\mathbb R}_{(\kappa)}))$, \[\imap_{\mathcal D(\Cond{(\kappa)}(\Ab))}(A,\underline{M}_{(\kappa)})=0.\] 
\end{lemma}
\begin{proof}
    As $\mathcal D(\Cond{(\kappa)}(\underline{\mathbb{R}}_{(\kappa)}))\cong \LMod{\underline{\mathbb R}_{(\kappa)}}{\Cond{(\kappa)}(\Sp)}$ has small limits and colimits (\cref{derivedcondensedmodulesallsmallimitsandcolimits}), $\imap_{\mathcal D(\Cond{(\kappa)}(\Ab))}(-,\underline{\mathbb R}_{(\kappa)})$ commutes with small limits in the first variable (\cref{leftadjointsstableundercolimits}). 
    As the free $\underline{\mathbb R}_{(\kappa)}$-modules generate $\mathcal D(\Cond{(\kappa)}(\underline{\mathbb R}_{(\kappa)}))$ under small colimits (\cite[Proposition 4.7.3.14]{higheralgebra}), it suffices to show that for $T\in\mathcal D(\Cond{(\kappa)}(\Ab))$, \[\imap_{\mathcal D(\Cond{(\kappa)}(\Ab))}(\underline{\mathbb R}_{(\kappa)}[T],\underline{M}_{\kappa})=0.\]  
    For $T\in\mathcal D(\Cond{(\kappa)}(\Ab))$, 
    \[ \imap_{\mathcal D(\Cond{(\kappa)}(\Ab))}(\underline{\mathbb R}_{(\kappa)}[T],\underline{M}_{\kappa})\cong \imap_{\mathcal D(\Cond{(\kappa)}(\Ab))}(T,\imap_{\mathcal D(\Cond{(\kappa)}(\Ab))}(\underline{\mathbb R}_{(\kappa)},\underline{M}_{\kappa}))\] (\cref{forgetfreeadjunctionmodules}). It is therefore enough to show that 
    \[ \imap_{\mathcal D(\Cond{(\kappa)}(\Ab))}(\underline{\mathbb R}_{(\kappa)},\underline{M}_{\kappa})=0.\] 

    By \cref{forkappalargeenoughclosedmonoidal}, it suffices to show this for $\Cond{\kappa}(\Ab)$.
    Under the identification \[\mathcal D(\Cond{\kappa}(\Ab))\cong \Cond{\kappa}(\mathcal D(\Ab))\subseteq \Fun(\Pro(\Fin)_{\kappa}^{\operatorname{op}},\mathcal D(\Ab)),\]  
    $\imap_{\mathcal D(\Cond{\kappa}(\Ab))}(\underline{\mathbb R}_{\kappa},\underline{M}_{\kappa})$ becomes  
    \[ \Pro(\Fin)_{\kappa}^{\operatorname{op}}\ni T\mapsto \map_{\mathcal D(\Cond{\kappa}(\Ab))}(\underline{\mathbb R}_{\kappa}\otimes_{\mathcal D(\Cond{\kappa}(\Ab))}\mathbb Z[\underline{T}_{\kappa}],\underline{M}_{\kappa}),\] cf. \cref{enrichmentrecoversinternalhom}. 
    It therefore suffices to show that $\map_{\mathcal D(\Cond{\kappa}(\Ab))}(\underline{\mathbb R}_{\kappa}\otimes\mathbb Z[\underline{T}_{\kappa}],\underline{M}_{\kappa})=0$ for all $S\in\Pro(\Fin)_{\kappa}$. 

    For a $\kappa$-condensed abelian group $A$ denote by $S_*(A)\to A\in\Cond{\kappa}(\Ab)$ its Breen-Deligne resolution (\cite[Theorem 4.5]{Scholzecondensed}). As filtered colimits in $\Cond{\kappa}(\Ab)$ are exact (\cref{condensedabeliangroupsgrothendieckaxioms}), \[A\cong S_*(A)=\colim{n\to \infty}S_*(A)^{\leq n}\] is the colimit over its stupid truncations. This implies that for  $T\in\Pro(\Fin)_{\kappa}$, 
    \[ \map_{\mathcal D(\Cond{\kappa}(\Ab))}(A\otimes \mathbb Z[\underline{T}_{\kappa}],\underline{\mathbb R}_{\kappa})\cong \clim{n\to \infty}\map_{\mathcal D(\Cond{\kappa}(\Ab))}(S_*(A)^{\leq n}\otimes \mathbb \mathbb Z[\underline{T}_{\kappa}],\underline{\mathbb R}_{\kappa}).\] 
    For $n\in\mathbb N_0$, we have cofiber sequences $S_*(A)^{\leq n}\to S_*(A)^{\leq n+1}\to S_{n+1}(A)[n+1]=\oplus_{k=1}^{b_{n+1}}\mathbb Z[A^{i_{k,n+1}}]$. 
    By \cref{condensedandsheafcohomology1}, 
    \[ \map_{\mathcal D(\Cond{\kappa}(\Ab))}(S_{n+1}(\underline{\mathbb R}_{\kappa})\otimes_{\mathcal D(\Cond{\kappa}(\Ab))}\mathbb Z[\underline{T}_{\kappa}],\underline{M}_{\kappa})\cong \prod_{k=1}^{b_{n+1}}\cH{\sheaf}(\mathbb R^{i_{k,n+1}}\times T,M)\] and  
    \[ \map_{\mathcal D(\Cond{\kappa}(\Ab))}(S_{n+1}(0)\otimes_{\mathcal D(\Cond{\kappa}(\Ab))}\mathbb Z[\underline{T}_{\kappa}],\underline{M}_{\kappa})\cong \prod_{k=1}^{b_{n+1}}\cH{\sheaf}(T,M).\]  
    By \cite[Proposition 11.13]{Bredon}, sheaf cohomology with constant coefficients is a homotopy invariant, see also \cref{grostoposcohomologyomotopyinvariant}. 
    This implies that the chain map $\pi_*\colon S_*(\underline{\mathbb R}_{\kappa})\to S_*(0)$ induced by $\underline{\mathbb R}_{\kappa}\to 0$ (\cite[Theorem 4.10]{Scholzecondensed}) defines an equivalence 
    \[  \map_{\mathcal D(\Cond{\kappa}(\Ab))}(S_{n+1}(0)\otimes\mathbb Z[\underline{T}_{\kappa}],\underline{M}_{\kappa})\cong \map_{\mathcal D(\Cond{\kappa}(\Ab))}(S_{n+1}(\underline{\mathbb R}_{\kappa})\otimes\mathbb Z[\underline{T}_{\kappa}],\underline{M}_{\kappa}).\] 
    It follows inductively that
    \[ \pi^*\colon \map_{\mathcal D(\Cond{\kappa}(\Ab))}(S_{*}(0)_{\leq n}\otimes\mathbb Z[\underline{T}_{\kappa}],\underline{M}_{\kappa})\to\map_{\mathcal D(\Cond{\kappa}(\Ab))}(S_{*}(\underline{\mathbb R}_{\kappa})_{\leq n}\otimes\mathbb Z[\underline{T}_{\kappa}],M)\] is an equivalence for all $n\in\mathbb N_0$, and hence 
    \begin{align*}0&= \map_{\mathcal D(\Cond{\kappa}(\Ab))}(S_*(0)\otimes_{\mathcal D(\Cond{\kappa}(\Ab))}\mathbb Z[\underline{T}_{\kappa}],\underline{M}_{\kappa})\\&\cong \map_{\mathcal D(\Cond{\kappa}(\Ab))}(S_{*}(\underline{\mathbb R}_{\kappa})\otimes_{\mathcal D(\Cond{\kappa}(\Ab))}\mathbb Z[\underline{T}_{\kappa}],\underline{M}_{\kappa})\\ &\cong \map_{\mathcal D(\Cond{\kappa}(\Ab))}(\underline{\mathbb R}_{\kappa}\otimes\mathbb Z[T],\underline{M}_{\kappa}).\qedhere\end{align*}
\end{proof}
\begin{lemma}[{\cite[Theorem 4.3]{Scholzecondensed}}]\label{cohomologyofcompacthausdorffabeliangroupwithcoefficientsinrvanishes}
If $A$ is a ($\kappa$-light) compact Hausdorff abelian topological group, then 
$\imap_{\mathcal D(\Cond{(\kappa)}(\Ab))}(\underline{A}_{(\kappa)},\underline{\mathbb R}_{(\kappa)})=0$. 
\end{lemma}
\begin{proof}By \cref{forkappalargeenoughclosedmonoidal}, it suffices to  show this for $\Cond{\kappa}(\Ab)$. 
    Under the identification \[\mathcal D(\Cond{\kappa}(\Ab))\cong \Cond{\kappa}(\mathcal D(\Ab))\subseteq \Fun(\Pro(\Fin)_{\kappa}^{\operatorname{op}},\mathcal D(\Ab)),\] $\imap_{\mathcal D(\Cond{\kappa}(\algebra{R}))}(A,\underline{\mathbb R}_{\kappa})$ becomes 
    \[ \Pro(\Fin)_{\kappa}^{\operatorname{op}}\ni T\mapsto \map_{\mathcal D(\Cond{\kappa}(\Ab))}(\underline{A}_{\kappa}\otimes_{\mathcal D(\Cond{\kappa}(\Ab))} \mathbb Z[\underline{T}_{\kappa}],\underline{\mathbb R}_{\kappa}),\] cf. \cref{enrichmentrecoversinternalhom}. 
    By \cref{Breendelignespectralsequence}, we have a spectral sequence 
    \[E_1^{p,q}=\oplus_{k=1}^{b_p}\Ext^q_{\Cond{\kappa}(\Ab)}(\mathbb Z[\underline{A^{n_k}\times T}_{\kappa}],\underline{\mathbb R}_{\kappa})\Rightarrow \pi_{-(p+q)}\map_{\mathcal D(\Cond{\kappa}(\Ab))}(\underline{A}_{\kappa}\otimes_{\mathcal D(\Cond{\kappa}(\algebra{R}))} \mathbb Z[\underline{T}_{\kappa}],\underline{\mathbb R}_{\kappa}).\] 
    By \cref{condensedandsheafcohomology1,sheafcohomologywvspacecoeffvanonparcomplc}, 
    \[ E_1^{*,q}(T)=\begin{cases}
    0 & q\neq 0\\ 
    \Hom_{\Cond{\kappa}(\Ab)}(S_*(A)\otimes \mathbb Z[\underline{T}_{\kappa}],\underline{\mathbb R}_{\kappa}) & q=0,
    \end{cases}\] where $S_*(A)$ denotes the Breen-Deligne resolution of $A$. 
    We are therefore reduced to showing that 
    $\Hom_{\Cond{\kappa}(\Ab)}(S_*(A)\otimes\mathbb Z[\underline{T}_{\kappa}],\underline{\mathbb R}_{\kappa})$ is exact.

    \cref{enrichmentrecoversinternalhom,kappacontinuousfullyfaithfullyintocondensed} imply that for $n\in\mathbb N_0$, \[\iHom_{\Cond{\kappa}(\Ab)}(S_n(A),\underline{\mathbb R}_{\kappa})\cong\prod_{k=1}^{b_p}\underline{\mathcal C(A^{i_{k,p}},\mathbb R)}_{\kappa},\] where $\mathcal C(A^{i_{k,p}},\mathbb R)$ is endowed with the compact open topology and $b_p, i_{k,p}$ are the coefficients of the Breen-Deligne resolution. 
    The sup-norm makes this a Banach space. 
    \cref{kappacontinuousfullyfaithfullyintocondensed} implies that the differentials $S_n(A)\to S_{n-1}(A)$ induce $\kappa$-continuous group homomorphisms \[ \prod_{k=1}^{b_{p+1}}\mathcal C(A^{i_{k,p+1}},\mathbb R)\to \prod_{k=1}^{b_p}\mathcal C(A^{i_{k,p}},\mathbb R).\]  As both terms are Banach spaces, these are continuous, i.e. have bounded norm. 
    This implies that $\Hom_{\Cond{\kappa}(\Ab)}(S_*(A)\otimes \mathbb Z[\underline{T}_{\kappa}],\underline{\mathbb R}_{\kappa})=\iHom_{\Cond{\kappa}(\Ab)}(S_n(A),\underline{\mathbb R}_{\kappa})(\mathbb Z[\underline{T}_{\kappa}])$ is represented by the cochain complex of Banach spaces 
    \[\ldots \to  \prod_{k=1}^{b_{p}}\mathcal C(A^{i_{k,p}}\times T,\mathbb R)\xrightarrow{d_{p+1}}  \prod_{k=1}^{b_{p+1}}\mathcal C(A^{i_{k,p+1}}\times T,\mathbb R)\ldots \to \] and the boundary maps have bounded norm. 
    Choose a chain homotopy $h_*\colon S(A)_*\to S(A)_{*+1}$ in $\Ch(\Cond{\kappa}(\Ab))$ from the multiplication by $2\colon S(A)_*\to S(A)_*$ to the map $S(2)\colon S(A)_*\to S(A)_*$ induced by multiplication by $2$ on $A$. This exists by \cite[Proposition 4.17]{Scholzecondensed}. 
    The induced map \[h_{n}^*\colon \prod_{k=1}^{b_{n+1}}\mathcal C(A^{i_{k,n+1}},\mathbb R)_{\kappa}=\iHom_{\Cond{\kappa}(\Ab)}(S(A)_{n+1},\mathbb R)\to \iHom_{\Cond{\kappa}(\Ab)}(S(A)_n,\mathbb R)=\prod_{k=1}^{b_{n}}\underline{\mathcal C(A^{i_{k,n}},\mathbb R)}_{\kappa}\] corresponds to a $\kappa$-continuous and hence continuous map 
    \[ \prod_{k=1}^{b_{n+1}}\mathcal C(A^{i_{k,n+1}},\mathbb R)\to \prod_{k=1}^{b_n}\mathcal C(A^{i_{k,p}},\mathbb R),\] 
 which induces a continuous map 
    \[ h_{n}^*(T)\colon \prod_{k=1}^{b_{n+1}}\mathcal C(A^{i_{k,n+1}}\times T,\mathbb R)\to \prod_{k=1}^{b_{n}}\mathcal C(A^{i_{k,n}}\times T,\mathbb R).\] This is a map of Banach spaces and hence has bounded norm. 

    Assume that \[f\in \Hom_{\Cond{\kappa}(\Ab)}(S_n(A)\otimes\mathbb Z[T],\mathbb R)=\prod_{k=1}^{b_n}\mathcal C(A^{i_{k,n}}\times T,\mathbb R)\] with $d_n(f)=0$. 
    Then $2f=S(2)^*f+d_{n-1}(h_{n-1}^*f)$, whence 
    \[ f=\frac{1}{2}(S(2)^*f+d_{n-1}(h_{n-1}^*f)).\] 
    In particular, $S(2)^*f\in\ker(d_n)$. 
    By iterating this process, we find that \[ f=\frac{1}{2^m}S(2^m)^*f+d_{n-1}(\sum_{k=1}^m\frac{1}{2^{k}}(h_{n-1}^*S(2^{k-1})^*f).\] 
    As $||\frac{1}{2^n}S(2^n)^*f||\leq\frac{1}{2^n}||f||$ and $d_{n-1}$ is $\mathbb R$-linear, it follows that \[f=\sum_{k\in\mathbb N}d_{n-1}(\frac{1}{2^{k}}(h_{n-1}^*S(2^{k-1})^*f)).\] 
    As $h_{n-1}^*(T)$ has bounded norm (see above), for $k\in\mathbb N_0$, \[||h_{n-1}^*(T)S(2^{k-1})^*f||\leq ||h_{n-1}^*(T)||||S(2^{k-1})^*f||\leq ||h_{n-1}^*(T)||||f||<\infty,\] and hence the limit \[\sum_{k\in\mathbb N}\frac{1}{2^{k}}(h_{n-1}^*(T)S(2^{k-1})^*f)\in \mathcal C(\sqcup_{k=1}^{b_n}A^{i_{k,n}}\times T,\mathbb R)\] exists. 
    As $d_{n-1}$ is continuous and $\mathbb R$-linear, this implies that \[d_{n-1}(\sum_{k\in\mathbb N}\frac{1}{2^{k}}(h_{n-1}^*S(2^{k-1})^*f))=f.\qedhere\] 
\end{proof}

\begin{lemma}[{\cite[Theorem 4.3]{Scholzecondensed}}]\label{cohomologyoftori}
    Suppose $M$ is a discrete abelian group. 
    \begin{romanenum}
        \item The connecting homomorphism of the exponential sequence \[0\to \mathbb Z\to\mathbb R\to\sphere{1}\to 0\] is an equivalence  
        \[ M[-1]\cong \imap_{\mathcal D(\Cond{(\kappa)}(\Ab))}(\underline{\sphere{1}}_{(\kappa)},\underline{M}_{(\kappa)}).\]
        \item For a ($\kappa$-small) set $I$, the projections $\prod_{I}\sphere{1}\to\sphere{1}$ induce an isomorphism  
    \[ \bigoplus_{I}\imap_{\mathcal D(\Cond{(\kappa)}(\Ab))}(\underline{\sphere{1}}_{(\kappa)},\underline{M}_{(\kappa)})\cong \imap_{\mathcal D(\Cond{\kappa}(\Ab))}(\prod_{I}\underline{\sphere{1}}_{\kappa},\underline{M}_{(\kappa)}).\]
    \end{romanenum}
\end{lemma}
\begin{proof}
    The first statement holds by \cref{cohomologyofrwithdiscretecoefficientsvanishes}. 
    The first statement implies the second statement fo $I$ finite. It remains to show that
\[ \colim{\substack{F\subseteq I\\ F\text{ finite }}}\imap_{\mathcal D(\Cond{\kappa}(\Ab))}(\prod_{F}\underline{\sphere{1}}_{\kappa},\underline{M}_{\kappa})\to \imap_{\mathcal D(\Cond{\kappa}(\Ab))}(\prod_{I}\underline{\sphere{1}}_{\kappa},\underline{M}_{\kappa})\] is an equivalence. 
By \cref{forkappalargeenoughclosedmonoidal}, it suffices to prove this for $\Cond{\kappa}(\Ab)$. 
By \cref{enrichmentrecoversinternalhom}, under the identification $\mathcal D(\Cond{\kappa}(\Ab))\cong \Cond{\kappa}(\mathcal D(\Ab))\subseteq \Fun(\Pro(\Fin)_{\kappa}^{\operatorname{op}},\mathcal D(\Ab))$, 
\[ \colim{\substack{F\subseteq I\\ F\text{ finite }}}\imap_{\mathcal D(\Cond{\kappa}(\Ab))}(\prod_{F}\underline{\sphere{1}}_{\kappa},\underline{M}_{\kappa})\] is the $\kappa$-condensed sheafification of 
\[\Pro(\Fin)_{\kappa}^{\operatorname{op}}\ni T\mapsto \colim{\substack{F\subseteq I\\ F\text{ finite }}}\map_{\mathcal D(\Cond{\kappa}(\Ab))}(\prod_{F}\underline{\sphere{1}}_{\kappa}\otimes\mathbb Z[\underline{T}_{\kappa}],\underline{M}_{\kappa}),\] and \[ \imap_{\mathcal D(\Cond{\kappa}(\Ab))}(\prod_{I}\underline{\sphere{1}}_{\kappa},\underline{M}_{\kappa})\] is
\[\Pro(\Fin)_{\kappa}^{\operatorname{op}}\ni T\mapsto \map_{\mathcal D(\Cond{\kappa}(\Ab))}(\prod_{I}\underline{\sphere{1}}_{\kappa}\otimes\mathbb Z[\underline{T}_{\kappa}],\underline{M}_{\kappa}).\] 

It therefore suffices to show that for $T\in\Pro(\Fin)_{\kappa}$, 
\[ \colim{\substack{F\subseteq I\\ F\text{ finite }}}\map_{\mathcal D(\Cond{\kappa}(\Ab))}(\prod_{F}\underline{\sphere{1}}_{\kappa}\otimes\mathbb Z[\underline{T}_{\kappa}],\underline{M}_{\kappa})\cong \map_{\mathcal D(\Cond{\kappa}(\Ab))}(\prod_{I}\underline{\sphere{1}}_{\kappa}\otimes\mathbb Z[\underline{T}_{\kappa}],\underline{M}_{\kappa}).\] 
As filtered colimits in $\Cond{\kappa}(\Ab)$ are exact, it is enough to show that for all $i\in\mathbb Z$ and $T\in\Pro(\Fin)_{\kappa}$, 
\[\colim{\substack{F\subseteq I\\ F\text{ finite }}} \pi_{i}\map_{\mathcal D(\Cond{\kappa}(\Ab))}(\prod_F\underline{\sphere{1}}_{\kappa}\otimes \mathbb Z[\underline{T}_{\kappa}],\underline{M}_{\kappa})\to  \pi_{i}\map_{\mathcal D(\Cond{\kappa}(\Ab))}(\prod_I\underline{\sphere{1}}_{\kappa}\otimes \mathbb Z[\underline{T}_{\kappa}],\underline{M}_{\kappa})\] is an equivalence. 

By \cref{Breendelignespectralsequence}, we have spectral sequences 
\begin{align*} E_1^{p,q}=\colim{\substack{F\subseteq I\\ F\text{ finite }}}\prod_{k=1}^{b_p}\pi_{-q}\map_{\mathcal D(\Cond{\kappa}(\Ab))}(\mathbb Z[\prod_F\underline{\sphere{i_{n_k,p}}\times T}_{\kappa}],\underline{M}_{\kappa}) \end{align*} and 
\[ E_1^{p,q}=\prod_{k=1}^{b_p}\pi_{-q}\map_{\mathcal D(\Cond{\kappa}(\Ab))}(\mathbb Z[\prod_I\underline{\sphere{i_{n_k,b}}\times T}_{\kappa}],\underline{M}_{\kappa})\Rightarrow \pi_{-(p+q)}\map_{\mathcal D(\Cond{\kappa}(\Ab))}(\prod_I\underline{\sphere{1}}_{\kappa}\otimes \mathbb Z[\underline{T}_{\kappa}],\underline{M}_{\kappa})\] which converge completely to 
\[  \colim{\substack{F\subseteq I\\ F\text{ finite }}}\pi_{-(p+q)}\map_{\mathcal D(\Cond{\kappa}(\Ab))}(\prod_F\underline{\sphere{1}}_{\kappa}\otimes \mathbb Z[\underline{T}_{\kappa}],\underline{M}_{\kappa})\] and \[  \pi_{-(p+q)}\map_{\mathcal D(\Cond{\kappa}(\Ab))}(\prod_I\underline{\sphere{1}}_{\kappa}\otimes \mathbb Z[\underline{T}_{\kappa}],\underline{M}_{\kappa}),\] respectively. 

Naturality of these spectral sequences (\cref{Breendelignespectralsequence}) implies that the map on limit terms is induced by the canonical map of $E_1$-pages 
\[ \colim{\substack{F\subseteq I\\ F\text{ finite }}}\prod_{k=1}^{b_p}\pi_{-q}\map_{\mathcal D(\Cond{\kappa}(\Ab))}(\mathbb Z[\prod_F\underline{\sphere{i_{n_k,b}}\times T}_{\kappa}],\underline{M}_{\kappa})\to \prod_{k=1}^{b_p}\pi_{-q}\map_{\mathcal D(\Cond{\kappa}(\Ab))}(\mathbb Z[\prod_I\underline{\sphere{i_{n_k,p}}\times T}_{\kappa}],\underline{M}_{\kappa}).\] 
By \cref{cohomologyclausenscholze}, \[ \colim{\substack{F\subseteq I\\ F\text{ finite }}}\prod_{k=1}^{b_p}\pi_{-q}\map_{\mathcal D(\Cond{\kappa}(\Ab))}(\mathbb Z[\prod_F\underline{\sphere{i_{n_k,b}}\times T}_{\kappa}],\underline{M}_{\kappa})\cong \prod_{k=1}^{b_p}\colim{\substack{F\subseteq I\\ F\text{ finite}}}\cH{\sheaf}^q(\prod_{F}\sphere{1}\times T,M)\] and \[\prod_{k=1}^{b_p}\pi_{-q}\map_{\mathcal D(\Cond{\kappa}(\Ab))}(\mathbb Z[\prod_I\underline{\sphere{i_{n_k,p}}\times T}_{\kappa}],\underline{M}_{\kappa})\cong \prod_{k=1}^{b_p}\cH{\sheaf}^{q}(\prod_{I}{\sphere{1}}\times T,M).\] 
As $\prod_{I,F}\mathbb S^1\times T$ are paracompact and Hausdorff, their sheaf cohomology agrees with their \v{C}ech cohomology (see e.g.\ \cite[Theorem 5.1.10]{Godement}), and since $\prod_{I}\mathbb S^1\times T\cong \clim{\substack{F\subseteq I\\ F\text{ finite }}}\prod_{F}\sphere{1}\times T$ is an inverse limit of compact Hausdorff spaces, the canonical map \[ \colim{\substack{F\subseteq I\\ F\text{ finite}}}\check{H}^*(\prod_{F}\sphere{1}\times T,M)\to \check{H}^*(\prod_{I}\sphere{1}\times T,M)\] on \v{C}ech cohomology is an isomorphism, cf. \cite[page 231]{Godement}. 
\end{proof}
\begin{cor}\label{duality}
    For a ($\kappa$-small) set $I$ and a discrete abelian group $M$, the canonical map 
    \[\oplus_{I}\underline{M}_{(\kappa)}\to \imap_{\mathcal D(\Cond{\kappa}(\Ab))}(\prod_I\underline{\mathbb Z}_{(\kappa)},\underline{M}_{(\kappa)})\] is an equivalence. 
\end{cor}
\begin{proof}
    As ($\kappa$-small) products in $\Cond{(\kappa)}(\Ab)$ are exact (\cref{condensedabeliangroupsgrothendieckaxioms}), the exponential exact sequence \[ 0\to \underline{\mathbb Z}_{(\kappa)}\to \underline{\mathbb R}_{(\kappa)}\to \underline{\sphere{1}}_{(\kappa)}\to 0\] yields a map of fiber sequences 
\begin{center}\begin{tikzcd}[cramped, sep=small]
    \oplus_{I}\imap_{\mathcal D(\Cond{(\kappa)}(\Ab))}(\underline{\sphere{1}}_{(\kappa)},\underline{M}_{(\kappa)})\arrow[r] \arrow[d] &  \imap_{\mathcal D(\Cond{(\kappa)}(\Ab))}(\prod_{I}\underline{\sphere{1}}_{(\kappa)},\underline{M}_{(\kappa)})\arrow[d]
\\\oplus_{I}\imap_{\mathcal D(\Cond{}(\Ab))}(\underline{\mathbb R}_{(\kappa)},\underline{M}_{(\kappa)}) \arrow[d]\arrow[r] &  \imap_{\mathcal D(\Cond{(\kappa)}(\Ab))}(\prod_{I}\underline{\mathbb R}_{(\kappa)},\underline{M}_{(\kappa)})\arrow[d]\\ 
\oplus_{I}\imap_{\mathcal D(\Cond{(\kappa)}(\Ab))}(\underline{\mathbb Z}_{(\kappa)},\underline{M}_{(\kappa)})\arrow[r] & \imap_{\mathcal D(\Cond{(\kappa)}(\Ab))}(\prod_I\underline{\mathbb Z}_{(\kappa)},\underline{M}_{(\kappa)})
\end{tikzcd}
\end{center}
By \cref{cohomologyofrwithdiscretecoefficientsvanishes}, the two middle terms vanish, and by \cref{cohomologyoftori}, the upper horizontal map is an equivalence, whence so is the bottom horizontal map. 
\end{proof}
\subsection{Solid modules}
\label{section:solidmodules}
A central feature of the condensed formalism is the existence of an extremely well-behaved notion of (non-archimedean) completeness for condensed abelian groups, called solidity, which we review in this section. 
Denote by $\iHom_{\Cond{(\kappa)}(\Ab)}(-,-)$ the internal Hom of ($\kappa$-)condensed abelian groups. 
\begin{definition}[{\cite[Def. 5.1]{Scholzecondensed}}]\label{definitionsolid}
    For a profinite set $S$ let 
    \[ \mathbb Z[\underline{S}_{(\kappa)}]^{\blacksquare}\coloneqq \iHom_{\Cond{(\kappa)}(\Ab)}(\iHom_{\Cond{(\kappa)}(\Ab)}(\mathbb Z[\underline{S}_{(\kappa)}], \mathbb Z), \mathbb Z).\]
    The counit for $\mathbb Z[\underline{S}_{(\kappa)}]\otimes \dashv \iHom_{\Cond{(\kappa)}(\Ab)}(\mathbb Z[\underline{S}_{(\kappa)}],-)$ evaluates to a map \[\mathbb Z[\underline{S}_{(\kappa)}]\otimes_{\Cond{(\kappa)}(\Ab)}\iHom_{\Cond{(\kappa)}(\Ab)}(\mathbb Z[\underline{S}_{(\kappa)}], \mathbb Z)\to \mathbb Z\] which is adjoint to a map $\mathbb Z[\underline{S}_{(\kappa)}]\to \mathbb Z[\underline{S}_{(\kappa)}]^{\blacksquare}$. 

    A ($\kappa$-)condensed abelian group $A$ is \emph{solid} if for all $S\in\Pro(\Fin)_{(\kappa)}$, pullback along the map \[\mathbb Z[\underline{S}_{(\kappa)}]\to\mathbb Z[\underline{S}_{(\kappa)}]^{\blacksquare}\] is an isomorphism 
    \[ \Hom_{\Cond{(\kappa)}(\Ab)}(\mathbb Z[\underline{S}_{(\kappa)}]^{\blacksquare},A)\cong \Hom_{\Cond{(\kappa)}(\Ab)}(\mathbb Z[\underline{S}_{(\kappa)}],A).\] 
    Denote by $\Sol{(\kappa)}\subseteq \Cond{(\kappa)}(\Ab)$ the full subcategories of solid abelian groups. More generally, if $\algebra{R}\in \Alg(\Cond{(\kappa)}(\Ab))$ is a discrete condensed ring, denote by $\Sol{(\kappa)}(\algebra{R})\subseteq \Cond{(\kappa)}(\algebra{R})$ the full subcategory on $\algebra{R}$-modules whose underlying condensed abelian group is solid. 
\end{definition}  
We will show below (\cref{nullsequencessummable}) that if $M$ is a $(\To)$ topological abelian group such that $\underline{M}_{\kappa}/\underline{M}$ is solid, then every null-sequence in $M$ is summable, and that locally profinite abelian groups represent solid abelian groups (\cref{profinitesolid}).
This justifies to consider solidity as a condensed analogue of (non-archimedean) completeness. 

The category of solid modules enjoys excellent formal properties, including the existence of a completion functor and a completed tensor product. We now give an overview of the results we prove in this section. 
    \begin{romanenum}
        \item \label{solidclosedunnderlimitscolimitslist}$\Sol{(\kappa)}(\algebra{R})\subseteq \Cond{(\kappa)}(\algebra{R})$ is an abelian category closed under small limits, small colimits and extensions. (\cref{solidclosedunderlimitscolimitskappa}, \cref{underivedsolidificationwithoutkappa})
        \item The inclusion $\Sol{(\kappa)}(\algebra{R})\subseteq\Cond{(\kappa)}(\algebra{R})$ has a left adjoint $(-)^{\solid\algebra{R}}$, called solidification (\cref{underivedsolidification}, \cref{underivedsolidificationwithoutkappa}). 
        \item \label{solidificationfreeonprofinitelist}For $S\in\Pro(\Fin)_{(\kappa)}$, $\mathbb Z[\underline{S}_{(\kappa)}]^{\solid\mathbb Z}=\mathbb Z[\underline{S}_{(\kappa)}]^{\blacksquare}\cong \prod_{I}\mathbb Z$ for a set $I$ of cardinality $|I|\leq \wt(S)$ (\cref{underivedsolidification}, \cref{solidificationoffreemodules}, \cref{Nobelingspecker}). 
        \item For a commutative ring $\algebra{R}\in\CAlg(\Cond{(\kappa)}(\Ab))$, the localisation $(-)^{\solid\algebra{R}}$ is symmetric monoidal, i.e. there exists cocontinuous symmetric monoidal structure on $\Sol{(\kappa)}(\algebra{R})$ such that $(-)^{\solid\algebra{R}}$ enhances to a symmetric monoidal functor. (\cref{solidmodulesismodulesinsolid})
        \item \label{profinitestableundersolidtensorproductlist} For $S\in\Pro(\Fin)_{(\kappa)}$, $\algebra{R}[\underline{S}_{(\kappa)}]^{\solid\algebra{R}}$ is projective in $\Sol{(\kappa)}(\algebra{R})$. In particular, $\Sol{(\kappa)}(\algebra{R})$ has enough projectives (\cref{solidenoughprojectiveskappa,solidenoughprojectiveswithoutkappa}).
        If $\algebra{R}$ is commutative and $P,Q\in\Sol{(\kappa)}(\algebra{R})$ are projective, so is their solid tensor product $P\otimes_{\Sol{(\kappa)}(\algebra{R})}Q$. (\cref{solidprojectivestensorproduct})
    \end{romanenum}
    We also discuss derived categories of solid modules. The inclusion $\Sol{(\kappa)}(\algebra{R})\subseteq\Cond{(\kappa)}(\algebra{R})$ extends uniquely to a $t$-exact, small colimits preserving functor \[f\colon \mathcal D(\Sol{(\kappa)}(\algebra{R}))\to \mathcal D(\Cond{(\kappa)}(\algebra{R})).\] 
    In many cases, this functor admits a left adjoint, for instance for $\Cond{}(\algebra{R})$, (\cref{derivedsolidificationwithoutkappa}) as well as for $\Cond{\kappa}(\algebra{R})$ for strong limit cardinals $\kappa$ (\cref{derivedsolificationexistsstronglimitcardinal}) or $\kappa=\aleph_1$ (\cref{derivedsolidificationexistssflatrings,lightringssflat}), and for $\algebra{R}=\mathbb Z$ (\cref{derivedsolidificationexistsabeliangroupskappa}).   
    \begin{romanenum}[resume]
    \item The forget functor $\mathcal D(\Sol{(\kappa)})\to\mathcal D(\Cond{(\kappa)}(\Ab))$ is fully faithful with essential image 
    \[ \{ M\in\mathcal D(\Cond{(\kappa)}(\Ab))\,|\, H_i(M)\in\Sol{(\kappa)}\text{ for all } i\in\mathbb Z\}.\] (\cref{derivedsolidissolidderived,derivedsolidissolidderivedwithoutkappa}).
    \item The localisation $\mathcal D(\Cond{(\kappa)}(\Ab))\to\mathcal D(\Sol{(\kappa)}(\Ab))$ is symmetric monoidal, and the induced symmetric monoidal structure on $\mathcal D(\Sol{(\kappa)})$ is closed. (\cref{symmetricmonoidalstructureonderivedsolidabeliangroups}). 
    \end{romanenum}
    
    We will prove the above statements in a peculiar order, with later proofs building on earlier results. 
    As a first step, we apply a theorem of Nöbeling and Specker to compute $\mathbb Z[\underline{S}_{\kappa}]^{\blacksquare}$ for ($\kappa$-light) profinite sets $S$ (\cref{Nobelingspecker}). 
    Next, we show that $\kappa$-condensed solid abelian groups satisfy a stronger locality condition than the one assumed in their definition (\cref{solidderived}).
    This serves as the essential ingredient for the proofs of the above claims.  
    We will first prove the $\kappa$-condensed statements and then use that $\Sol{}(\algebra{R})=\colim{\substack{\kappa\\ \kappa\text{ regular }}}\Sol{\kappa}(\algebra{R})$ (\cref{colimitdescriprionsolid}) to deduce the corresponding condensed results.  
    We then show that if $\Cond{(\kappa)}(\algebra{R})$ has enough projectives, then the forget functor $\mathcal D(\Sol{}(\algebra{R}))\to\mathcal D(\Cond{}(\algebra{R}))$ admits a left adjoint $(-)^{L\solid\algebra{R}}$ (\cref{derivedsolificationexistsstronglimitcardinal}, \cref{derivedsolidificationwithoutkappa}). This in particular holds for $\Cond{}(\algebra{R})$ as well as for $\Cond{\kappa}(\algebra{R})$ for strong limit cardinals $\kappa$. 
    Our arguments are minor modifications of the discussion in \cite[Lectures 5,6]{Scholzecondensed} who proved the statements for $\Sol{}$ and $\Sol{\kappa}$ for strong limit cardinals, and worked with the homotopy category of the derived category. Small modifications are necessary since the proofs of \cite[Lectures 5,6]{Scholzecondensed} used that for strong limit cardinals $\kappa$, $\Cond{(\kappa)}(\Ab)$ has enough projectives. For light solid abelian groups, the proofs can be simplified significantly, see \cite[Lectures 5,6]{Analyticstacks} or \cite[section 3]{Juansolidgeometry}. The above results rely crucially on the fact that ($\kappa$-light) profinite sets generate $\Cond{(\kappa)}(\an)$ under colimits. The availability of a well-behaved theory of (non-archimedean) complete abelian group objects is a key advantage of the condensed \textit{topos} over the gros topos. 

    In \cref{section:computationssolidification}, we give examples of solid abelian groups and compute the solidification of free condensed abelian groups on compact Hausdorff spaces and CW-complexes (\cref{solidificationcompacthausdorffspace}, \cref{solidcw}). We will use these computations in \cref{section:condensedcohomologywithsolidcoefficients} to obtain further identifications of condensed with sheaf cohomology. Combined with \ref{profinitestableundersolidtensorproductlist}, they serve as the essential ingredient for our comparison of solid group cohomology with continuous group cohomology. 
    
    In \cref{section:sflatrings}, we characterize ($\kappa$-)condensed rings for which the forget functor \[\mathcal D(\Sol{(\kappa)}(\algebra{R}))\to\mathcal D(\Sol{(\kappa)})\] factors over an equivalence \begin{align}\label{modulesinsolidaresolidmodulesderivedforsflatrings} \mathcal D(\Sol{(\kappa)}(\algebra{R}))\cong \LMod{\algebra{R}^{\solid}}{\mathcal D(\Sol{(\kappa)}(\algebra{R}))}.\end{align}
    We call such rings ($\kappa$-)$s$-flat and give many examples (\cref{lightringssflat,examplessflatrings,profiniteringssflat}). 
     The identification \ref{modulesinsolidaresolidmodulesderivedforsflatrings} implies that the following holds: 
    Suppose that $\algebra{R}$ is a ($\kappa$-)$s$-flat condensed ring. 
    \begin{romanenum}
        \item The forget functor \[f\colon \mathcal D(\Sol{(\kappa)}(\algebra{R}))\to \mathcal D(\Cond{(\kappa)}(\algebra{R}))\] has a left adjoint $(-)^{L\solid\algebra{R}}$. (\cref{derivedsolidificationexistssflatrings})
 
        \item  If $\algebra{R}$ is commutative, there exists a cocontinuous symmetric monoidal structure on $\mathcal D(\Sol{\kappa}(\algebra{R}))$ such that $(-)^{L\solid\algebra{R}}$ enhances to a symmetric monoidal functor. This symmetric monoidal structure is closed and compatible with the $t$-structure on $\mathcal D(\Sol{(\kappa)}(\algebra{R}))$ (\cref{symmetricmonoidalstructuresolidrmodules,symmetricmonoidalstructuresolidrmodulescompatibletstructure}). 
        The induced monoidal structure on the heart $\Sol{(\kappa)}(\algebra{R})$ is the symmetric monoidal structure from \cref{solidmodulesismodulesinsolid}.
        
        \item If $\algebra{R}^{L\solid\mathbb Z}\cong \algebra{R}^{\solid\mathbb Z}$, then the functor $f\colon \mathcal D(\Sol{(\kappa)}(\algebra{R}))\to\mathcal D(\Cond{(\kappa)}(\algebra{R}))$ is fully faithful with essential image \[\{ M\in \mathcal D(\Cond{(\kappa)}(\algebra{R}))\, |\, H_i(M)\in\Sol{(\kappa)}(\algebra{R}) \, \text{ for all } i\in\mathbb Z\}.\] (\cref{derivedsolidificationunderlyingforacyclicrings}). This in particular holds if $\algebra{R}\in\Alg(\Sol{(\kappa)})\subseteq \Alg(\Cond{(\kappa)}(\Ab))$. 
    \end{romanenum}

We now work towards the proof of the statements listed above.
Recall from \cref{internalcohomologyprofinite}, that for a ($\kappa$-light) profinite space $S$, 
\[\iHom_{\Cond{(\kappa)}(\Ab)}(\mathbb Z[\underline{S}_{(\kappa)}], \mathbb Z)\cong \underline{\mathcal C(S, \mathbb Z)}_{(\kappa)}, \] where $\mathcal C(S, \mathbb Z)$ is endowed with the discrete topology. 
In particular, \[\mathbb Z[\underline{S}_{(\kappa)}]^{\blacksquare}\cong \iHom_{\Cond{(\kappa)}(\Ab)}(\underline{\mathcal C(S, \mathbb Z)}_{(\kappa)}, \mathbb Z).\] 
\begin{lemma}[Nöbeling, Specker]\label{Nobelingspecker}
    For a profinite set $S$, there exists a set $I$ of cardinality $|I|\leq \wt(S)$ such that $\mathcal C(S, \mathbb Z)\cong \oplus_{i\in I}\mathbb Z$. 
    In particular, for $S\in\Pro(\Fin)_{(\kappa)}$, there exists a ($\kappa$-small) set $I$ such that 
    $\mathbb Z[\underline{S}_{(\kappa)}]^{\blacksquare}\cong \prod_{I}\mathbb Z$. 
\end{lemma}
\begin{proof}This is an immediate consequence of a theorem of Nöbeling \cite{nobeling} generalizing work of Specker, see also \cite[Theorem 5.4]{Scholzecondensed} for a proof. 
They show that for a profinite set $S$, there exists a set $I$ with $\mathcal C(S, \mathbb Z)=\oplus_I \mathbb Z$. 
If $J$ is a set such that $X\hookrightarrow \prod_{J}\{0,1\}$, they choose $I$ to consists of a collection of finite products of elements of $J$. In particular, $|I|\leq |\bigcup_{n\in\mathbb N} J^n|\leq \max\{\aleph_0,|J|\}$. 
By \cref{propertiesweight}, we can choose $|J|\leq \wt(X)$. 
In particular, if $\kappa$ is an uncountable cardinal, for $S\in\Pro(\Fin)_{\kappa}$ there exists a set $I$ of cardinality $|I|< \kappa$ with $\mathcal C(S, \mathbb Z)\cong \oplus_{I}\mathbb Z$.
\end{proof}
Our next goal is to prove the following characterisation of solid abelian groups which will be central for proving the results on solid modules listed above.  
\begin{proposition}\label{solidderived}
    \begin{romanenum}
    \item A condensed abelian group $M\in\Cond{\kappa}(\Ab)$ is solid if and only if for all $S\in \Pro(\Fin)_{\kappa}$, 
    \[ \imap_{\mathcal D(\Cond{\kappa}(\Ab))}(\mathbb Z[\underline{S}_{\kappa}]^{\blacksquare},M)\cong \imap_{\mathcal D(\Cond{\kappa}(\Ab))}(\mathbb Z[\underline{S}_{\kappa}],M)\] via the pullback along $\mathbb Z[\underline{S}_{\kappa}]\to\mathbb Z[\underline{S}_{\kappa}]^{\blacksquare}$. 
    
    \item For a solid abelian group $M$ and $S\in\Pro(\Fin)_{\kappa}$, 
    \[\imap_{\mathcal D(\Cond{\kappa}(\Ab))}(\mathbb Z[\underline{S}]^{\blacksquare},M)\cong \imap_{\mathcal D(\Cond{\kappa}(\Ab))}(\mathbb Z[\underline{S}_{\kappa}],M)\in \mathcal D(\Cond{\kappa}(\Ab))^{\heart}, \] and in particular, 
      \[\map_{\mathcal D(\Cond{\kappa}(\Ab))}(\mathbb Z[\underline{S}_{\kappa}]^{\blacksquare},M)\cong \map_{\mathcal D(\Cond{\kappa}(\Ab))}(\mathbb Z[\underline{S}_{\kappa}],M)\in \mathcal D(\Ab)^{\heart}.\]

    \item The category $\Sol{\kappa}\subseteq \Cond{\kappa}(\Ab)$ is closed under filtered colimits.
    \end{romanenum} 
\end{proposition}
This was shown for strong limit cardinals in \cite[Lectures 5,6]{Scholzecondensed} and for $\Cond{\light}(\Ab)$ in \cite[Lectures 5,6]{Analyticstacks}, see also \cite{Juansolidgeometry}. 
Our proof of \cref{solidderived} is a minor modification of the argument from \cite[p. 38ff]{Scholzecondensed} and will be given on page \pageref{proofsolidderived} after the necessary preparation. 
The analogous results for $\Sol{}$ hold as well (\cref{underivedsolidificationwithoutkappa,derivedsolidissolidderivedwithoutkappa,solidenoughprojectiveswithoutkappa}), we will deduce this from the above statement. 
\begin{notation}
    For $S\in\Pro(\Fin)_{\kappa}$ and $A\in\Cond{\kappa}(\Ab)$ let \[\mathcal M(S,A)\coloneqq \iHom_{\Cond{\kappa}(\Ab)}(\underline{\mathcal C(S, \mathbb Z)}_{\kappa},A).\]
\end{notation}
\begin{ex}
    For $S\in\Pro(\Fin)_{\kappa}$, $\mathcal M(S, \mathbb Z)\cong \mathbb Z[\underline{S}_{\kappa}]^{\blacksquare}$ by \cref{Nobelingspecker}. 
\end{ex}
By \cref{Nobelingspecker}, for all $S\in\Pro(\Fin)_{\kappa}$ there exists a $\kappa$-small set $I$ such that $\mathcal M(S,A)=\prod_{I}A$ for all $A\in\Cond{\kappa}(\Ab)$. 
As $\kappa$-small products are exact in $\Cond{\kappa}(\Ab)$ (\cref{condensedabeliangroupsgrothendieckaxioms}), this implies that \[\mathcal M(S,-)\colon\Cond{\kappa}(\Ab)\to\Cond{\kappa}(\Ab)\] is an exact functor. 
In particular, the exponential sequence $0\to\mathbb Z\to\mathbb R\to\sphere{1}\to 0$ yields exact sequences 
\[ 0\to \mathbb Z[\underline{S}]^{\blacksquare}\to\mathcal M(S, \mathbb R)\to\mathcal M(S, \sphere{1})\to 0, \,  S\in\Pro(\Fin)_{\kappa}.\] 

\begin{lemma}[{\cite[Lecture 6]{Scholzecondensed}}]\label{extendalongexponentialexactsequence}
 Suppose $(J_i)_{i\in I}, (U_k)_{k\in K}$ are small collections of $\kappa$-small sets. 
   Every group homomorphism \[ \oplus_I\prod_{J_i}\underline{\mathbb Z}_{\kappa}\to \oplus_{K}\prod_{U_k}\underline{\mathbb Z}_{\kappa}\text{ in }\Cond{\kappa}(\Ab)\] extends uniquely to a map \[ \oplus_I\prod_{J_i}\underline{\mathbb R}_{\kappa}\to \oplus_{K}\prod_{U_k}\underline{\mathbb R}_{\kappa}\text{ in } \Cond{\kappa}(\Ab).\] 
\end{lemma}
\begin{proof}
    It of course suffices to show this for $I=*$. We show that pullback along $\prod_{J_i}\underline{\mathbb Z}_{\kappa}\to \prod_{J_i}\underline{\mathbb R}_{\kappa}$ is an isomorphism
    \begin{align}\label{equation1} \Hom_{\Cond{\kappa}(\Ab)}(\prod_{J_i}\underline{\mathbb R}_{\kappa}, \oplus_{K}\prod_{U_k}\underline{\mathbb R}_{\kappa})\cong \Hom_{\Cond{\kappa}(\Ab)}(\prod_{J_i}\underline{\mathbb Z}_{\kappa}, \oplus_{K}\prod_{U_k}\underline{\mathbb R}_{\kappa}).\end{align} 

    Since $|J_i|<\kappa$, the sequence of $\kappa$-condensed abelian groups \[0\to \prod_{J_i}\underline{\mathbb Z}_{\kappa}\to\prod_{J_i}\underline{\mathbb R}_{\kappa}\to\prod_{J_i}\underline{\sphere{1}}_{\kappa}\to 0\in \Cond{\kappa}(\Ab)\] is exact (\cref{condensedabeliangroupsgrothendieckaxioms}), whence it is enough to show that for $k=0,1$,
    \begin{align}\label{claimvanishingcohomology}\pi_{-k}\map_{\mathcal D(\Cond{\kappa}(\Ab))}(\prod_{J_i}\underline{\sphere{1}}_{\kappa}, \oplus_{K}\prod_{U_k}\underline{\mathbb R}_{\kappa})=\Ext^{k}_{\Cond{\kappa}(\Ab)}(\prod_{J_i}\underline{\sphere{1}}_{\kappa}, \oplus_{K}\prod_{U_k}\underline{\mathbb R}_{\kappa})=0, \end{align} then \ref{equation1} follows.  
    By \cref{compacttopologicalgroupspseudocoherent} and exactness of filtered colimits in $\Cond{\kappa}(\Ab)$ (\cref{condensedabeliangroupsgrothendieckaxioms}), \[\map_{\mathcal D(\Cond{\kappa}(\Ab))}(\prod_{J_i}\underline{\sphere{1}}_{\kappa},-)\colon\Cond{\kappa}(\Ab)\to \mathcal D(\Ab)\] commutes with filtered colimits. 
    Since $\kappa$-small products in $\Cond{\kappa}(\Ab)$ are exact (\cref{condensedabeliangroupsgrothendieckaxioms}), it follows that
    \begin{align*}\map_{\mathcal D(\Cond{\kappa}(\Ab))}(\prod_{J_i}\underline{\sphere{1}}_{\kappa}, \oplus_{K}\prod_{U_k}\underline{\mathbb R}_{\kappa})& \cong \colim{\substack{F\subseteq K\\ F\text{ finite }}}\prod_{f\in F}\prod_{U_f}\map_{\mathcal D(\Cond{\kappa}(\Ab))}(\prod_{J_i}\underline{\sphere{1}}_{\kappa}, \underline{\mathbb R}_{\kappa})\\& \cong \colim{\substack{F\subseteq K\\ F\text{ finite }}}\prod_{f\in F}\prod_{U_f} * \end{align*} is contractible by \cref{cohomologyofcompacthausdorffabeliangroupwithcoefficientsinrvanishes}. 
    This proves (\ref{claimvanishingcohomology}). 
\end{proof}
\begin{proposition}[{\cite[Lecture 6]{Scholzecondensed}}]\label{goodcomplexcohomologicaldimension}
    Suppose $C_*\in \Ch(\Cond{\kappa}(\Ab))_{\geq 0}$ is a chain complex such that \begin{itemize}\label{goodcomplex}
        \item[(*)] For all $i\in\mathbb N_0$, there exists a set $J_i$ and a collection of $\kappa$-light profinite sets $(S^{j})_{j\in J_i}$ such that $C_i=\oplus_{j\in J_i}\mathbb Z[\underline{S^j}_{\kappa}]^{\blacksquare}$. 
        \end{itemize}
    For all $S\in\Pro(\Fin)_{\kappa}$, 
    \begin{align}\label{claimconnectivity}\map_{\mathcal D(\Cond{\kappa}(\Ab))}(\mathbb Z[\underline{S}_{\kappa}],C_*)\in \mathcal D(\Ab)_{\geq -1}\end{align} and \[ \imap_{\mathcal D(\Cond{\kappa}(\Ab))}(\mathcal M(S, \sphere{1}),C_*)\in \mathcal D(\Cond{\kappa}(\Ab))_{\geq -1}.\] 
    \end{proposition}
    \begin{proof} 
    Suppose that $(C_*, \partial)$ is a chain complex as in the lemma and for $n\in\mathbb N_0$, $j\in J_n$ choose $I_j$, $|I_j|<\kappa$ such that $\mathbb Z[\underline{S^j}_{\kappa}]^{\blacksquare}\cong \prod_{I_j}\mathbb Z$.     
    By \cref{extendalongexponentialexactsequence}, we obtain a chain complex $(C_*^{\mathbb R}, \partial^{\mathbb R}_*)$ with $C_n=\oplus_{j\in J_n}\prod_{I_j}\underline{\mathbb R}_{\kappa}$ and $\partial_{n}^{\mathbb R}$ the unique extension of $\partial_n$. 
    In particular, for all $i\in\mathbb N_1$, $\partial_{n}^{\mathbb R}$ descends to a group homomorphism 
    \[\partial_{n}^{\sphere{1}}\colon  \oplus_{j\in J_n}\prod_{I_j}\underline{\sphere{1}}_{\kappa}\to \oplus_{j\in J_{n-1}}\prod_{I_j}\underline{\sphere{1}}_{\kappa}.\]     
    This defines a chain complex $(C_*^{\sphere{1}}, \partial^{\sphere{1}})$. 
    The short exact sequence of chain complexes 
    \[ 0\to C_*\to C_*^{\mathbb R}\to C_*^{\sphere{1}}\to 0\] yields fiber sequences
    \begin{align*} \imap_{\mathcal D(\Cond{\kappa}(\Ab))}(\mathcal M(S, \sphere{1}),C_*)\to &\imap_{\mathcal D(\Cond{\kappa}(\Ab))}(\mathcal M(S, \sphere{1}),C_*^{\mathbb R})\\& \to \imap_{\mathcal D(\Cond{\kappa}(\Ab))}(\mathcal M(S, \sphere{1}),C_*^{\sphere{1}})\end{align*} and 
    \begin{align*} \map_{\mathcal D(\Cond{\kappa}(\Ab))}(\mathbb Z[\underline{S}_{\kappa}],C_*)\to &\map_{\mathcal D(\Cond{\kappa}(\Ab))}(\mathbb Z[\underline{S}_{\kappa}],C_*^{\mathbb R})\\ &\to \map_{\mathcal D(\Cond{\kappa}(\Ab))}(\mathbb Z[\underline{S}_{\kappa}],C_*^{\sphere{1}}).\end{align*}
    It is therefore enough to show that for $ A=\underline{\mathbb R}_{\kappa}, \underline{\sphere{1}}_{\kappa}$ and $S\in\Pro(\Fin)_{\kappa}$, \[\imap_{\mathcal D(\Cond{\kappa}(\Ab))}(\mathcal M(S, \sphere{1}),C_*^{A})\in \mathcal D(\Cond{\kappa}(\Ab))_{\geq 0}\] and \[\map_{\mathcal D(\Cond{\kappa}(\Ab))}(\mathbb Z[\underline{S}_{\kappa}],C_*^{A})\in \mathcal D(\Ab)_{\geq 0}, \] then it follows that \[\imap_{\mathcal D(\Cond{\kappa}(\Ab))}(\mathcal M(S, \mathbb Z),C_*)\in \mathcal D(\Cond{\kappa}(\Ab))_{\geq -1}\] and \[\map_{\mathcal D(\Cond{\kappa}(\Ab))}(\mathbb Z[\underline{S}_{\kappa}],C_*)\in \mathcal D(\Ab)_{\geq -1}.\] 
    
   Denote by $\mathcal H\subseteq \mathcal D(\Cond{\kappa}(\Ab))$ the full subcategory on objects $H\in\mathcal D(\Cond{\kappa}(\Ab))$ such that \[\imap_{\mathcal D(\Cond{\kappa}(\Ab))}(\mathcal M(S, \sphere{1}),H)\in \mathcal D(\Cond{\kappa}(\Ab))_{\geq 0}\] and \[\map_{\mathcal D(\Cond{\kappa}(\Ab))}(\mathbb Z[\underline{S}_{\kappa}],H)\in \mathcal D(\Ab)_{\geq 0}\] for all $S\in\Pro(\Fin)_{\kappa}$. 

    We first show that $\underline{\mathbb R}_{\kappa},\underline{\sphere{1}}_{\kappa}\in\mathcal H$. 
    \cref{cohomologyisext1} and \cite[Theorem 3.2, 3.3]{Scholzecondensed}/\cref{cohomologyclausenscholzeprofinite} imply that
    \[ \map_{\mathcal D(\Cond{\kappa}(\Ab))}(\mathbb Z[\underline{S}_{\kappa}], \underline{\mathbb R}_{\kappa}),\,  \map_{\mathcal D(\Cond{\kappa}(\Ab))}(\mathbb Z[\underline{S}_{\kappa}], \mathbb Z)\in \mathcal D(\Ab)^{\heart}\] are concentrated in degree $0$.  
    It now follows from the exponential exact sequence that for $S\in\Pro(\Fin)_{\kappa}$, 
    \[ \map_{\mathcal D(\Cond{\kappa}(\Ab))}(\mathbb Z[\underline{S}_{\kappa}], \underline{\sphere{1}}_{\kappa})\in \mathcal D(\Ab)_{\geq 0}.\] 
    By \cref{cohomologyoftori,cohomologyofcompacthausdorffabeliangroupwithcoefficientsinrvanishes,Nobelingspecker}, for all $S\in\Pro(\Fin)_{\kappa}$, \
    \[ \imap_{\mathcal D(\Cond{\kappa}(\Ab))}(\mathcal M(S, \sphere{1}), \underline{\mathbb R}_{\kappa})=0,\]
    and \[\imap_{\mathcal D(\Cond{\kappa}(\Ab))}(\mathcal M(S,\sphere{1}),\mathbb Z)\in \Cond{\kappa}(\Ab)[-1].\] 
    Hence the exponential exact sequence implies that \[ \imap_{\mathcal D(\Cond{\kappa}(\Ab))}(\mathcal M(S,\sphere{1}),\underline{\sphere{1}}_{\kappa})\cong \Sigma \imap_{\mathcal D(\Cond{\kappa}(\Ab))}(\mathcal M(S,\sphere{1}),\mathbb Z)\in\mathcal D(\Cond{\kappa}(\Ab))^{\heart},\] which shows that $\underline{\mathbb R}_{\kappa},\underline{\sphere{1}}_{\kappa}\in\mathcal H$.

    As filtered colimits in $\Cond{\kappa}(\Ab)$ are exact (\cref{condensedabeliangroupsgrothendieckaxioms}), \cref{freemodulespseudocoherentdabenriched} and \cref{compacttopologicalgroupspseudocoherent} imply that $\mathcal H\cap \Cond{\kappa}(\Ab)$ is closed under filtered colimits. 
    Since $\kappa$-small products in $\Cond{\kappa}(\Ab)$ are exact (\cref{condensedabeliangroupsgrothendieckaxioms}), $\mathcal H\cap \Cond{\kappa}(\Ab)$ is also closed under $\kappa$-small products. 
    It now follows from the above that if $(I_j)_{j\in J}$ is a small family of $\kappa$-small sets, then \begin{align}\label{levelsoffreesolutioninhs}\oplus_{j\in J}\prod_{i\in I_j}\underline{\mathbb R}_{\kappa}, \oplus_{j\in J}\prod_{i\in I_j}\underline{\sphere{1}}_{\kappa}\in\mathcal H.\end{align} 

    We now want to show that for $A=\mathbb R,\sphere{1}$ and all complexes $C_*^{A}$ of the above form, $C_*^{A}\in \mathcal H$. 
    We first explain that it suffices to show that $\tau_{\leq n}C_*^{A}\in\mathcal H$ for all $n\in\mathbb N_0$. 
    For $n\in\mathbb Z$, the cofiber sequence $\tau_{\leq n}C_*^{A}\to \tau_{\leq n+1}C_*^{A}\to H_{n+1}(C_*^{A})[n+1]$ implies that 
    for $X\in \Cond{\kappa}(\Ab))$, 
    \[ \tau_{\geq -n}\imap_{\mathcal D(\Cond{\kappa}(\Ab))}(X,\tau_{\leq n}C_*^{A})\cong \tau_{\geq -n}\imap_{\mathcal D(\Cond{\kappa}(\Ab))}(X,\tau_{\leq n+1}C_*^{A})\] and 
    \[ \tau_{\geq -n}\map_{\mathcal D(\Cond{\kappa}(\Ab))}(X,\tau_{\leq n}C_*^{A})\cong \tau_{\geq -n}\map_{\mathcal D(\Cond{\kappa}(\Ab))}(X,\tau_{\leq n+1}C_*^{A}).\]
    As the $t$-structure on $\mathcal D(\Cond{\kappa}(\Ab))$ is left-complete (\cref{tstructureleftcomplete}), $C_*^{A}=\clim{n}\tau_{\leq n}C_*^{A}$, and by the above, for all $n\in\mathbb N_0$ and $X\in\Cond{\kappa}(\Ab)$, 
    \[\tau_{\geq -n}\imap_{\mathcal D(\Cond{\kappa}(\Ab))}(X,C_*^A)\cong \imap_{\mathcal D(\Cond{\kappa}(\Ab))}(X,\tau_{\leq n}C_*^A)\] and 
    \[\tau_{\geq -n}\map_{\mathcal D(\Cond{\kappa}(\Ab))}(X,C_*^A)\cong \map_{\mathcal D(\Cond{\kappa}(\Ab))}(X,\tau_{\leq n}C_*^A).\]
    It therefore suffices to show that $\tau_{\leq n}C_*^{A}\in\mathcal H$ for all $n\in\mathbb Z_{\geq -1}$. 
    We prove this by induction on $n$, the case $n=-1$ is trivial. 
    Suppose now that $n\in\mathbb Z_{\geq 0}$ such that for all complexes $C_*^A$ of the above form, $\tau_{\leq n-1}C_*^A\in\mathcal H$.  
    Fix such a complex $C_*^{A}$ and denote by $C_*^{A, \leq n}$ its stupid truncations.
    We have fiber sequences
    \[ C_*^{A, \leq n}\to \tau_{\leq n}C_*^A\to \im(\partial_{n+1}^A)[n+1],\]
    \[ \ker(\partial_{n+1}^A)[n+1]\to C_{n+1}^{A}[n+1]\to \im(\partial_{n+1}^A)[n+1],\] 
    and \[ C_*^{A,\leq n-1}\to C_*^{A, \leq n}\to C_n^A[n].\] 
    By induction hypothesis, $\tau_{\leq n-1}C_*^{A, \leq n-1}=C_*^{A,\leq n-1}\in \mathcal H$, and by \ref{levelsoffreesolutioninhs}, $C_n^{A}[n], C_{n+1}^{A}[n+1]\in\mathcal H$ for all $n\in\mathbb N_0$. 
    Since $\mathcal H$ is closed under extensions and cofibers and for $m\in\mathbb N_0$, $\mathcal H[m]\subseteq \mathcal H$, it suffices to show that $\ker(\partial^A_{n+1})[1]\in \mathcal H$, then it follows that $\tau_{\leq n}C_*^A\in\mathcal H$.

    We now show that $\ker(\partial^A_{n+1})[1]\in \mathcal H$. 
    Suppose that $C_*^{A}$ is constructed from a chain complex $C_*$ with \[C_n=\oplus_{j\in J_n}\mathbb Z[\underline{S^j}_{\kappa}]^{\blacksquare}, S_j\in\Pro(\Fin)_{\kappa} \text{ for all }j\in J_n,n\in\mathbb N_0.\] 
    For $j\in J_{n}$: Since $\mathbb Z[\underline{S^j}_{\kappa}]$ is compact in $\Cond{\kappa}(\Ab)$, there exists a finite subset $F(j)\subseteq J_{n-1}$ such that 
    \[\mathbb Z[\underline{S^j}_{\kappa}]\to \mathbb Z[\underline{S^j}_{\kappa}]^{\blacksquare}\xrightarrow{\partial_j|_{\mathbb Z[\underline{S^j}_{\kappa}]^{\blacksquare}}}\oplus_{k\in J_{i-1}}\mathbb Z[\underline{S^k}_{\kappa}]^{\blacksquare}\] factors over $\oplus_{k\in F(j)}\mathbb Z[\underline{S^k}_{\kappa}]^{\blacksquare}$. 
    This implies that the induced map \[\partial_{n}|_{\mathbb Z[\underline{S^j}_{\kappa}]^{\blacksquare}}\colon \mathbb Z[\underline{S^j}_{\kappa}]^{\blacksquare}=\prod_{I_j}\mathbb Z\to \oplus_{k\in J_{n-1}}\mathbb Z[\underline{S^k}_{\kappa}]^{\blacksquare}=\oplus_{k\in J_{n-1}}\prod_{U_k}\mathbb Z\] factors over $\oplus_{k\in F(j)}\mathbb Z[\underline{S^k}_{\kappa}]^{\blacksquare}$ as well. By construction of $\partial_i^A$, it follows that \[\partial_i^{A}(\prod_{I_j}A)\subseteq \oplus_{k\in F(j)}\prod_{U_k}A,\] where $I_j,U_k$ are such that $\mathbb Z[\underline{S^j}_{\kappa}]^{\blacksquare}=\prod_{I_j}\mathbb Z, \mathbb Z[\underline{S^k}_{\kappa}]^{\blacksquare}=\prod_{U_k}\mathbb Z$. 
    In particular, $\ker(\partial_n^{A})$ is a filtered colimit of condensed abelian groups of the form \[\ker(\prod_I \underline{A}_{\kappa}\to\prod_U \underline{A}_{\kappa})\] for homomorphisms $f_A\colon\prod_I \underline{A}_{\kappa}\to \prod_U \underline{A}_{\kappa}\in\Cond{\kappa}(\Ab)$ with $|I|,|U|<\kappa$. 

    Since filtered colimits in $\Cond{\kappa}(\Ab)$ are exact, \cref{freemodulespseudocoherentdabenriched} and \cref{compacttopologicalgroupspseudocoherent} imply that
    \[ \imap_{\mathcal D(\Cond{\kappa}(\Ab))}(\mathcal M(S, \mathbb Z),-)\colon \Cond{\kappa}(\Ab)[1]\to \mathcal D(\Cond{\kappa}(\Ab))\] and \[ \map_{\mathcal D(\Cond{\kappa}(\Ab))}(\mathbb Z[\underline{S}_{\kappa}],-)\colon \Cond{\kappa}(\Ab)[1]\to \mathcal D(\Ab)\] commute with filtered colimits, so in particular, 
    \[\mathcal H\cap \Cond{\kappa}(\Ab)[1]\subseteq \mathcal D(\Cond{\kappa}(\Ab))\] is closed under filtered colimits. 
    We are therefore reduced to showing that for all $\kappa$-small sets $I,U$ and all homomorphisms 
    $f_A\colon \prod_{I}\underline{A}_{\kappa}\to\prod_{U}\underline{A}_{\kappa}\in\Cond{\kappa}(\Ab),$ $\ker(f_A)[1]\in \mathcal H$. 
    
    We first show this for $A=\underline{\sphere{1}}_{\kappa}$.
    Since for all $\kappa$-light compact Hausdorff spaces $X$, $\underline{X}_{\kappa}(*)_{\kappa}=X$, the map $f_{\sphere{1}}$ induces a homomorphism of topological abelian groups $f_{\sphere{1}}(*)_{\kappa}\colon \prod_{I}{\sphere{1}}\to\prod_{J}\underline{\sphere{1}}$. By fully faithfulness of $\underline{-}_{\kappa}$ (\cref{kappacontinuousfullyfaithfullyintocondensed}), $f_{\sphere{1}}=\underline{f_{\sphere{1}}(*)}_{\kappa}$. 
    In particular, \[\ker(f_{\sphere{1}})=\underline{\ker(f_{\sphere{1}}(*))}_{\kappa}\] is represented by a compact abelian group of weight $\leq\max\{|I|, \mathbb N\}<\kappa$. 
    Since Pontryagin-duality restricts to a contravariant equivalence between compact Hausdorff and discrete abelian groups, taking Pontryagin duals yields an epimorphism \[ q\colon \oplus_{k}\mathbb Z\to \ker(f_{\sphere{1}}(*))^{v}\] to the Pontryagin-dual $\ker(f_{\sphere{1}}(*))^{v}$ of $\ker(f_{\sphere{1}}(*))$. As subgroup of a free abelian group, $\ker(q)$ is free i.e.\ $\ker(q)=\oplus_K\mathbb Z$ for $|K|\leq |I|$. 
    By taking Pontryagin duals again we obtain a resolution \[0\to \ker(f_{\sphere{1}}(*))\to \prod_{I}\sphere{1}\to \prod_{K}\sphere{1}\to 0.\] As $\prod_{I}\sphere{1}\to \prod_{K}\sphere{1}$ is a continuous surjection between $\kappa$-light compact Hausdorff spaces, \[\prod_{I}\underline{\sphere{1}}_{\kappa}\to \prod_{K}\underline{\sphere{1}}_{\kappa}\] is an epimorphism of $\kappa$-condensed sets by \cref{condensedoncompactextremallydisconnected}.
    This implies that the induced sequence of $\kappa$-condensed abelian groups 
    \[0\to \ker(f_{\sphere{1}})\to \prod_{I}\underline{\sphere{1}}_{\kappa}\to \prod_{K}\underline{\sphere{1}}_{\kappa}\to 0\] is exact. 
    As $|K|,|I|<\kappa$, $\prod_{K}\underline{\sphere{1}}_{\kappa}, \prod_{I}\underline{\sphere{1}}_{\kappa}\in\mathcal H$ by \cref{levelsoffreesolutioninhs}, which shows that $\ker(f_{\sphere{1}})[1]\in\mathcal H$.

    By \cref{duality}, for $|K|<\kappa$, $\oplus_{K}\mathbb Z\cong \imap_{\mathcal D(\Cond{\kappa}(\Ab))}(\prod_{K}\mathbb Z,\mathbb Z)$, whence $f_{\mathbb Z}\colon\prod_{I}\mathbb Z\to\prod_{J}\mathbb Z$ is dual to a map $f^{v}\colon \oplus_{J}\mathbb Z\to \oplus_I \mathbb Z$. 
    This implies that $f_{\mathbb R}$ is dual to a map $f^{v}_{\mathbb R}\colon \oplus_{J}\underline{\mathbb R}_{\kappa}\to\oplus_{I}\underline{\mathbb R}_{\kappa}$. 
    The map $\oplus_J\underline{\mathbb R}_{\kappa}\to \Coker(f^{v}_{\mathbb R}(*)_{\kappa})$ admits a continuous $\mathbb R$-linear section, whence \[\Coker(f^{v}_{\mathbb R}(*)_{\kappa})=\oplus_{K}\mathbb R\] for $|K|\leq |I|$. 
    As $\oplus_{J/I}\underline{\mathbb R}_{\kappa}\cong \underline{(-)}_{\kappa}\circ (*)_{\kappa}(\oplus_{J/I}\underline{\mathbb R}_{\kappa})$ and $ \oplus_{J}\mathbb R(*)_{\kappa}\to \Coker(f^{v}_{\mathbb R})(*)_{\kappa}$ admits a continuous section, it follows that
    \[ \Coker(f^{v}_{\mathbb R})=\underline{(-)}_{\kappa}(\Coker(f^v_{\mathbb R}(*)_{\kappa}))=\oplus_K \underline{\mathbb R}_{\kappa}.\]   Consequently, $\ker(f_{\mathbb R})=\prod_K\underline{\mathbb R}_{\kappa}$. Since $|K|< \kappa$, $\ker(f_{\mathbb R})\in\mathcal H$ by \ref{levelsoffreesolutioninhs}. 
    \end{proof}

\begin{notation}\label{Gcategorysoliddefinition}
    For a $\kappa$-light profinite set $S$, the connecting homomorphism of the short exact sequence \[ 0\to \mathbb Z[\underline{S}_{\kappa}]^{\blacksquare}\to\mathcal M(S, \mathbb R)\to\mathcal M(S, \sphere{1})\to 0\] yields a natural transformation 
    \[ \partial^*\colon \imap_{\mathcal D(\Cond{\kappa}(\Ab))}(\mathcal M(S, \sphere{1}),-)\to \imap_{\mathcal D(\Cond{\kappa}(\Ab))}(\mathbb Z[\underline{S}_{\kappa}]^{\blacksquare},-)[-1].\] 
    Denote by \[ c_S\colon \imap_{\mathcal D(\Cond{\kappa}(\Ab))}(\mathcal M(S, \sphere{1}),-)\to \imap_{\mathcal D(\Cond{\kappa}(\Ab))}(\mathbb Z[\underline{S}_{\kappa}],-)[-1]\] the composition of $\partial^*$ with pullback along $\mathbb Z[\underline{S}_{\kappa}]\to\mathbb Z[\underline{S}_{\kappa}]^{\blacksquare}$. 

    Let $\mathcal G\subseteq \mathcal D(\Cond{\kappa}(\Ab))$ be the full subcategory on objects $G$ such that for all $S\in\Pro(\Fin)_{\kappa}$, 
    \begin{align}\label{claimsolidinternalprojective} c_S(G)\colon \imap_{\mathcal D(\Cond{\kappa}(\Ab))}(\mathcal M(S, \sphere{1}),G)\to \imap_{\mathcal D(\Cond{\kappa}(\Ab))}(\mathbb Z[\underline{S}_{\kappa}],G)[-1]\end{align} is an equivalence. 
\end{notation}
We will identify $\mathcal G$ with the derived category of solid abelian groups $\mathcal D(\Sol{\kappa})$ in \cref{derivedsolidissolidderived} below, but first prove some results on $\mathcal G$ which enter the proof of \cref{solidderived}. 
\begin{lemma}\label{zingcategorysolid}
    $\mathbb Z\in\mathcal G$. 
\end{lemma}
\begin{proof}Fix $S\in\Pro(\Fin)_{\kappa}$.
    By \cref{Nobelingspecker}, there exists $|I|<\kappa$ with $\mathcal C(S, \mathbb Z)\cong \oplus_{I}\mathbb Z$. This implies that $\mathcal M(S, \mathbb R)=\prod_{I}\mathbb R$.
    \cref{cohomologyoftori,cohomologyofcompacthausdorffabeliangroupwithcoefficientsinrvanishes} imply that for all discrete abelian groups $M$,     
    \[ \oplus_{I}M[-1]\cong \imap_{\mathcal D(\Cond{\kappa}(\Ab))}(\mathcal M(S, \sphere{1}),M)\cong \imap_{\mathcal D(\Cond{\kappa}(\Ab))}(\mathbb Z[\underline{S}_{\kappa}]^{\blacksquare},M)[-1].\] 
    It therefore suffices to show that pullback along $\eta\colon \mathbb Z[\underline{S}_{\kappa}]\to\mathbb Z[\underline{S}_{\kappa}]^{\blacksquare}$ yields an equivalence  
    \[  \imap_{\mathcal D(\Cond{\kappa}(\Ab))}(\mathbb Z[\underline{S}_{\kappa}]^{\blacksquare},\mathbb Z)\cong \imap_{\mathcal D(\Cond{\kappa}(\Ab))}(\mathbb Z[\underline{S}_{\kappa}], \mathbb Z).\] 
    Since $\mathbb Z[\underline{S}_{\kappa}]$ is flat, under the identification \[\mathcal D(\Cond{\kappa}(\Ab))\cong\Cond{\kappa}(\mathcal D(\Ab))\subseteq \Fun(\Pro(\Fin)_{\kappa}^{\operatorname{op}}, \mathcal D(\Ab)), \] $\imap_{\mathcal D(\Cond{\kappa}(\Ab))}(\mathbb Z[\underline{S}_{\kappa}], \mathbb Z)$ becomes 
    \begin{align*}\Pro(\Fin)_{\kappa}& \to \mathcal D(\Ab), \\ T &\mapsto \map_{\mathcal D(\Cond{\kappa}(\Ab))}(\mathbb Z[\underline{S\times T}_{\kappa}], \mathbb Z), \end{align*} cf.\ \cref{enrichmentrecoversinternalhom}. 
    By \cite[Theorem 3.2]{Scholzecondensed}/\cref{cohomologyclausenscholzeprofinite}, for all $T\in\Pro(\Fin)_{\kappa}$, \[\ckH(\underline{S\times T}_{\kappa}, \mathbb Z)\cong \mathcal C(S\times T, \mathbb Z)\in \mathcal D(\Cond{\kappa}(\Ab))^{\heart}.\] 
    It therefore suffices to show that $\eta$ induces an isomorphism \[\oplus_{I}\mathbb Z\cong \pi_0\imap_{\mathcal D(\Cond{\kappa}(\Ab))}(\mathbb Z[\underline{S}_{\kappa}]^{\blacksquare}, \mathbb Z)\xrightarrow{\eta^*} \pi_0\imap_{\mathcal D(\Cond{\kappa}(\Ab))}(\mathbb Z[\underline{S}_{\kappa}], \mathbb Z)\cong \underline{\mathcal C(S, \mathbb Z)}_{\kappa}.\] 
    This holds by construction of the identification of \cref{cohomologyoftori}. 
\end{proof}
\begin{cor}\label{Gcategorysolidissolid}
    \begin{romanenum}
        \item If $G\in\mathcal G$ and $M\in\mathcal D(\Cond{\kappa}(\underline{\mathbb R}_{\kappa}))\cong \CMod{\underline{\mathbb R}_{\kappa}}{\kappa}$ is a $\kappa$-condensed $\mathbb R$-module, then \[\imap_{\mathcal D(\Cond{\kappa}(\Ab))}(M,G)=0.\]
        \item For $G\in\mathcal G$ and $S\in\Pro(\Fin)_{\kappa}$, pullback along $\mathbb Z[\underline{S}_{\kappa}]\to\mathbb Z[\underline{S}_{\kappa}]^{\blacksquare}$ induces an isomorphism \[\imap_{\mathcal D(\Cond{\kappa}(\Ab))}(\mathbb Z[\underline{S}]^{\blacksquare},G)\cong\imap_{\mathcal D(\Cond{\kappa}(\Ab))}(\mathbb Z[\underline{S}_{\kappa}],G).\]
    \end{romanenum}
\end{cor} 
\begin{proof}
    Suppose $G\in\mathcal G$. Evaluating \ref{claimsolidinternalprojective} at $S=*$ implies that \[\imap_{\mathcal D(\Cond{\kappa}(\Ab))}(\sphere{1},G)\cong \imap_{\mathcal D(\Cond{\kappa}(\Ab))}(\mathbb Z,G)[-1]\] via the connecting homomorphism. Since
    \[ \imap_{\mathcal D(\Cond{\kappa}(\Ab))}(\sphere{1},G)\to\imap_{\mathcal D(\Cond{\kappa}(\Ab))}(\mathbb R,G)\to \imap_{\mathcal D(\Cond{\kappa}(\Ab))}(\mathbb Z,G)\] is a cofiber sequence, this implies that $\imap_{\mathcal D(\Cond{\kappa}(\Ab))}(\mathbb R,G)=0.$ 
    It follows that for every $A\in\mathcal D(\Cond{\kappa}(\Ab))$, 
    \[ \imap_{\mathcal D(\Cond{\kappa}(\Ab))}(\mathbb R\otimes_{\mathcal D(\Cond{\kappa}(\Ab))} A,M)\cong \imap_{\mathcal D(\Cond{\kappa}(\Ab))}(A, \imap_{\mathcal D(\Cond{\kappa}(\Ab))}(\mathbb R,M))=0.\] 
    Since the free modules $\mathbb R[A], A\in \mathcal D(\Cond{\kappa}(\Ab))$ generate $\LMod{\mathbb R}{\mathcal D(\Cond{\kappa}(\Ab))}$ under $\Delta^{\operatorname{op}}$-indexed colimits (\cite[Proposition 4.7.3.14]{higheralgebra}), \cref{forgetfreeadjunctionmodules} implies that \[\imap_{\mathcal D(\Cond{\kappa}(\Ab))}(V,G)=0\] for all condensed $\mathbb R$-modules $V$. 
    In particular, the connecting homomorphism \[\imap_{\mathcal D(\Cond{\kappa}(\Ab))}(\mathcal M(S, \sphere{1}),G)\to \imap_{\mathcal D(\Cond{\kappa}(\Ab))}(\mathbb Z[\underline{S}_{\kappa}]^{\blacksquare},G)[-1]\] from the exponential exact sequence is an equivalence for all $G\in\mathcal G$, which implies that pullback along $\mathbb Z[\underline{S}_{\kappa}]\to\mathbb Z[\underline{S}_{\kappa}]^{\blacksquare}$ induces an isomorphism \[\imap_{\mathcal D(\Cond{\kappa}(\Ab))}(\mathbb Z[\underline{S}_{\kappa}]^{\blacksquare},G)\cong\imap_{\mathcal D(\Cond{\kappa}(\Ab))}(\mathbb Z[\underline{S}_{\kappa}],G)\] for all $G\in\mathcal G$ and all $S\in\Pro(\Fin)_{\kappa}$. 
\end{proof}
\begin{cor}\label{Gcategorysolidclosedunderfilteredcolimits}
    \begin{romanenum}
\item The category $\mathcal G\subseteq\mathcal D(\Cond{\kappa}(\Ab))$ is a stable subcategory closed under extensions. 
\item For all $n\in\mathbb N_0$, $\mathcal G\cap \mathcal D(\Cond{\kappa}(\Ab))_{\leq n}\subseteq \mathcal D(\Cond{\kappa}(\Ab))$ is closed under filtered colimits. 
In particular, for $n\leq m\in\mathbb Z$, 
$\mathcal G\cap \mathcal D(\Cond{\kappa}(\Ab))_{[n,m]}\subseteq \mathcal D(\Cond{\kappa}(\Ab))$ is closed under filtered colimits.
    \end{romanenum} 
\end{cor} 
\begin{proof} Since for $S\in\Pro(\Fin)_{\kappa}$, \[\imap_{\mathcal D(\Cond{\kappa}(\Ab))}(\mathcal M(S, \sphere{1}),-)\text{ and }\imap_{\mathcal D(\Cond{\kappa}(\Ab))}(\mathbb Z[\underline{S}_{\kappa}],-)\] are exact functors, 
$\mathcal G$ is a stable subcategory closed under extensions. 
For $S\in\Pro(\Fin)_{\kappa}$ and $n\in\mathbb N_0$, \[ \imap_{\mathcal D(\Cond{\kappa}(\Ab))}(\mathcal M(S, \sphere{1}),-)|_{\mathcal D(\Cond{\kappa}(\Ab))_{\leq n}}\text{ and }\imap_{\mathcal D(\Cond{\kappa}(\Ab))}(\mathbb Z[\underline{S}_{\kappa}],-)|_{\mathcal D(\Cond{\kappa}(\Ab))_{\leq n}}\] preserve filtered colimits by \cref{Nobelingspecker}, \cref{compacttopologicalgroupspseudocoherent} and \cref{freemodulesonchpseudocoherent}. 
Since filtered colimits in $\Cond{\kappa}(\Ab)$ are exact, this implies that \[\mathcal G\cap \mathcal D(\Cond{\kappa}(\Ab))_{\leq n}\subseteq \mathcal D(\Cond{\kappa}(\Ab))\] and $\mathcal G\cap\mathcal D(\Cond{\kappa}(\Ab))_{[m,n]}$ are closed under filtered colimits for all $m\leq n\in\mathbb Z$.     
\end{proof}

\begin{cor}\label{Gcategorysolidcontainsgoodcomplexes}
If $C_*$ is a chain complex as in \cref{goodcomplex}, then $C_*$ represents an element in $\mathcal G$. 
\end{cor}

\begin{proof}
    By definition, $\mathcal G\subseteq \mathcal D(\Cond{\kappa}(\Ab))$ is closed under limits. Since $\kappa$-small products in $\Cond{\kappa}(\Ab)$ are exact (\cref{condensedabeliangroupsgrothendieckaxioms}), it follows that $\mathcal G\cap \Cond{\kappa}(\Ab)\subseteq \Cond{\kappa}(\Ab)$ is closed under $\kappa$-small products. \cref{Nobelingspecker} and \cref{zingcategorysolid} now imply that $\mathbb Z[\underline{S}_{\kappa}]^{\blacksquare}\in \mathcal G$ for all $S\in\Pro(\Fin)_{\kappa}$.
    Since $\mathcal G\cap\Cond{\kappa}(\Ab)\subseteq \Cond{\kappa}(\Ab)$ is closed under filtered colimits (\cref{Gcategorysolidclosedunderfilteredcolimits}) and $\kappa$-small products, $\mathcal G\cap \Cond{\kappa}(\Ab)\subseteq \Cond{\kappa}(\Ab)$ is closed under direct sums. In particular, if $(S_i)_{i\in I}\subseteq \Pro(\Fin)_{\kappa}$ is a small family of profinite sets, then $\oplus_{i\in I}\mathbb Z[\underline{S_i}_{\kappa}]^{\blacksquare}\in\mathcal G$. 
    
    As $\mathcal G\subseteq \mathcal D(\Cond{\kappa}(\Ab))$ is closed under extensions (\cref{Gcategorysolidclosedunderfilteredcolimits}), it follows by induction on the length that every finite length complex \[ 0\to C_n\to  C_{n-1}\to \ldots\to C_0\] with $C_i=\oplus_{j\in J_i}\mathbb Z[\underline{S_j}_{\kappa}]^{\blacksquare}$, $S_{j}\in\Pro(\Fin)_{\kappa}$ represents an element in $\mathcal G$.
    Suppose now that $C_*$ is an unbounded complex as in \cref{goodcomplexcohomologicaldimension}. 
    Let \[ C_*^{>n}\coloneqq (\ldots\to  C_{k+n}\to C_{k+n-1}\to\ldots\to C_{n+1}\to 0)\] the stupid connective cover and denote by $C_*^{\leq n}$ the stupid truncation.   
    The canonical map \[\Cone(C_*^{\leq n}\to C_*)\to C_*^{>n}\] is a quasi-isomorphism. 
    By \cref{goodcomplexcohomologicaldimension}, for $S\in\Pro(\Fin)_{\kappa}$, \[ \map_{\mathcal D(\Cond{\kappa}(\Ab))}(\mathbb Z[\underline{S}_{\kappa}],C_*^{>n})[-n-1]\in \mathcal D(\Cond{\kappa}(\Ab))_{\geq -1}, \] i.e.\ \[ \map_{\mathcal D(\Cond{\kappa}(\Ab))}(\mathbb Z[\underline{S}_{\kappa}],C_*^{>n})\in \mathcal D(\Cond{\kappa}(\Ab))_{\geq n}\] and \[ \imap_{\mathcal D(\Cond{\kappa}(\Ab))}(\mathcal M(S, \sphere{1}),C_*^{>n})[-n-1]\in \mathcal D(\Cond{\kappa}(\Ab))_{\geq -1}, \] i.e.\ \[ \imap_{\mathcal D(\Cond{\kappa}(\Ab))}(\mathcal M(S, \sphere{1}),C_*^{>n})\in \mathcal D(\Cond{\kappa}(\Ab))_{\geq n}.\]
    This implies that for $i < n$, $C_*^{\leq n}\to C_*$ induces isomorphisms \[ \pi_i(\map_{\mathcal D(\Cond{\kappa}(\Ab))}(\mathbb Z[\underline{S}_{\kappa}],C_*^{\leq n}))\cong\pi_i(\map_{\mathcal D(\Cond{\kappa}(\Ab))})(\mathbb Z[\underline{S}_{\kappa}],C_*)\] and \[ \pi_i(\imap_{\mathcal D(\Cond{\kappa}(\Ab))}(\mathcal M(S, \sphere{1}),C_*^{\leq n}))\cong \pi_i(\imap_{\mathcal D(\Cond{\kappa}(\Ab))}(\mathcal M(S, \sphere{1}),C_*)).\] 
    The first isomorphism implies that for $i<n$, \[ \pi_i\imap_{\mathcal D(\Cond{\kappa}(\Ab))}(\mathbb Z[\underline{S}_{\kappa}],C_*^{\leq n})\cong \pi_i(\imap_{\mathcal D(\Cond{\kappa}(\Ab))}(\mathbb Z[\underline{S}_{\kappa}],C_*)).\]
    Indeed, by t-exactness of sheafification and \cref{enrichmentrecoversinternalhom}, both sides are the $\kappa$-condensed sheafification of \[\Pro(\Fin)_{\kappa}^{\operatorname{op}}\ni T\mapsto \pi_i(\map_{\mathcal D(\Cond{\kappa}(\Ab))}(\mathbb Z[\underline{S\times T}_{\kappa}],C^{\leq n}_*))\cong \pi_i(\map_{\mathcal D(\Cond{\kappa}(\Ab))}(\mathbb Z[\underline{S\times T}_{\kappa}],C_*)).\] 
    Since for all $n\in\mathbb N_0$, $C^{\leq n}_*\in\mathcal G$, and the functors $\pi_i, i\in\mathbb Z$ are jointly conservative, this shows that $C_*$ represents an element in $\mathcal G$.
\end{proof}

\begin{lemma}\label{Gcategoryheart}Denote by $\mathcal G\subseteq \mathcal D(\Cond{\kappa}(\Ab))$ the subcategory from \cref{Gcategorysoliddefinition}. 
Then \[\mathcal G\cap \mathcal D(\Cond{\kappa}(\Ab))^{\heart}=\Sol{\kappa}.\]
\end{lemma}
\begin{proof}   
By \cref{Gcategorysolidissolid}, for $S\in\Pro(\Fin)_{\kappa}$, and $G\in\mathcal G$,
\[\imap_{\mathcal D(\Cond{\kappa}(\Ab))}(\mathbb Z[\underline{S}_{\kappa}]^{\blacksquare},G)\cong \imap_{\mathcal D(\Cond{\kappa}(\Ab))}(\mathbb Z[\underline{S}_{\kappa}],G), \] whence $\mathcal G\cap \mathcal D(\Cond{\kappa}(\Ab))^{\heart}\subseteq \Sol{\kappa}$. 

We claim that every solid abelian group $M$ admits a resolution $C_*\to M$ with $C_*$ a complex as in \cref{goodcomplex}, then \cref{Gcategorysolidcontainsgoodcomplexes} implies that $\Sol{\kappa}\subseteq \mathcal G$. 
Since $\mathbb Z[\underline{S}_{\kappa}],S\in\Pro(\Fin)_{\kappa}$ generate $\Cond{\kappa}(\Ab)$ under colimits, for $M\in\Cond{\kappa}(\Ab)$ there exists a set of $\kappa$-light profinite sets $(S_i)_{i\in I}\subseteq \Pro(\Fin)_{\kappa}$ with an epimorphism $\oplus_{i\in I}\mathbb Z[\underline{S_i}_{\kappa}]\to M$. 
If $M$ is solid, this factors over an epimorphism
$\oplus_{i\in I}\mathbb Z[\underline{S_i}_{\kappa}]^{\blacksquare}\to M$.
Suppose that a partial resolution 
\[ C_n\xrightarrow{\partial_n} \ldots \to C_{n-1}\to M\] as in \cref{goodcomplex} has been constructed. 
Choose a small collection of $\kappa$-light profinite sets $(S_i)_{i\in I}$ with a surjection $\oplus_{i\in I}\mathbb Z[\underline{S_i}_{\kappa}]\to \ker(\partial_n)$. 
Since $C_n\in\mathcal G\cap\mathcal D(\Cond{\kappa}(\Ab))^{\heart}$ is solid by the above, $\oplus_{i\in I}\mathbb Z[\underline{S^i}_{\kappa}]\to C_n$ extends uniquely to a group homomorphism \[\partial_{n+1}\colon \oplus_{i\in I}\mathbb Z[\underline{S_i}_{\kappa}]^{\blacksquare}\to C_n.\] 
By construction, $\ker(\partial_n)\subseteq \im(\partial_{n+1})$.
Since for $i\in I$, $\partial_{n}\circ \partial_{n+1}|_{\mathbb Z[\underline{S_i}_{\kappa}]^{\blacksquare}}$ is an extension of the zero morphism $\mathbb Z[\underline{S_i}_{\kappa}]\xrightarrow{0}C_{n-1}$ and $C_{n-1}$ is solid, $\partial_n\circ \partial_{n+1}=0$, i.e.\ $\im(\partial_{n+1})\subseteq \ker(\partial_n)$.
By repeating this construction we obtain the desired resolution $C_*$ of $M$.
\end{proof}
We can now finally deduce \cref{solidderived} from this. 
\begin{proof}[Proof of \cref{solidderived}]\label{proofsolidderived}
By \cref{Gcategorysolidissolid} and \cref{Gcategoryheart}, for $M\in\Sol{\kappa}$ and $S\in\Pro(\Fin)_{\kappa}$, 
pullback along $\mathbb Z[\underline{S}_{\kappa}]\to\mathbb Z[\underline{S}_{\kappa}]^{\blacksquare}$ is an isomorphism
\[\imap_{\mathcal D(\Cond{\kappa}(\Ab))}(\mathbb Z[\underline{S}_{\kappa}]^{\blacksquare},M)\cong \imap_{\mathcal D(\Cond{\kappa}(\Ab))}(\mathbb Z[\underline{S}_{\kappa}],M).\] 

Suppose that $M\in\Sol{\kappa}$ and $C_*\to M$ is a resolution of $M$ as in \cref{goodcomplexcohomologicaldimension}. It was shown in the proof of \cref{Gcategorysolidcontainsgoodcomplexes} that such a resolution exists. 
By \cref{goodcomplexcohomologicaldimension}, \[\imap_{\mathcal D(\Cond{\kappa}(\Ab))}(\mathcal M(S, \sphere{1}),M)\cong \imap_{\mathcal D(\Cond{\kappa}(\Ab))}(\mathcal M(S, \sphere{1}),C_*)\in \mathcal D(\Cond{\kappa}(\Ab))_{\geq -1}.\] 
Since $M\in\mathcal G$, it follows that 
\[\imap_{\mathcal D(\Cond{\kappa}(\Ab))}(\mathbb Z[\underline{S}_{\kappa}],M)\cong \imap_{\mathcal D(\Cond{\kappa}(\Ab))}(\mathcal M(S, \sphere{1}),M)[1]\in \mathcal D(\Cond{\kappa}(\Ab))_{\geq 0},\] which shows that $\imap_{\mathcal D(\Cond{\kappa}(\Ab))}(\mathbb Z[\underline{S}_{\kappa}],M)\in \mathcal D(\Cond{\kappa}(\Ab))^{\heart}$ is concentrated in degree $0$.

\cref{Gcategoryheart} and \cref{Gcategorysolidclosedunderfilteredcolimits} imply that $\mathcal D(\Cond{\kappa}(\Ab))^{\heart}\cap\mathcal G=\Sol{\kappa}\subseteq \Cond{\kappa}(\Ab)$ is closed under filtered colimits. 
\end{proof}

\begin{cor}\label{solidificationissolid}
    \begin{romanenum}\item For $S\in\Pro(\Fin)_{\kappa}$, $\mathbb Z[\underline{S}_{\kappa}]^{\blacksquare}\in\Sol{\kappa}$. 
    
    \item The family $\mathbb Z[\underline{S}_{\kappa}]^{\blacksquare}, S\in\Pro(\Fin)_{\kappa}$ generates $\Sol{\kappa}$ under small colimits.
    \end{romanenum} 
\end{cor}
\begin{proof}
    By \cref{Gcategorysolidcontainsgoodcomplexes}, $\mathbb Z[\underline{S}_{\kappa}]^{\blacksquare}\in\mathcal G$. 
    We showed in the proof of \cref{solidderived} that $\mathcal G\cap\mathcal D(\Cond{\kappa}(\Ab))^{\heart}=\Sol{\kappa}$. 
    Since every solid abelian group admits a resolution as in \cref{goodcomplex}, $\mathbb Z[\underline{S}_{\kappa}]^{\blacksquare}, S\in\Pro(\Fin)_{\kappa}$ generates $\Sol{\kappa}$ under colimits. 
\end{proof}

\begin{lemma}\label{solidclosedunderlimitscolimitskappa}
    Suppose $\kappa$ is an uncountable cardinal and $\algebra{R}\in\Alg(\Cond{(\kappa)}(\Ab))$. 
    The category \[\Sol{\kappa}(\algebra{R})\subseteq \Cond{\kappa}(\algebra{R})\] is  stable under small limits, colimits and extensions.
\end{lemma}
\begin{proof}
Since the tensor product on $\Cond{\kappa}(\Ab)$ is cocontinuous in both variables, the forget functor $\Cond{\kappa}(\algebra{R})\to \Cond{\kappa}(\Ab)$ reflects limits and colimits, hence it suffices to show that \[\Sol{\kappa}\subseteq \Cond{\kappa}(\Ab)\] is stable under small limits, colimits and extensions. It is immediate from the definition that  $\Sol{\kappa}(\mathbb Z)\subseteq \Cond{\kappa}(\Ab)$ is closed under limits. 
By \cref{solidderived}, $\Sol{\kappa}\subseteq \Cond{\kappa}(\Ab)$ is closed under filtered colimits. 
This implies that $\Sol{\kappa}$ is closed under direct sums (which are filtered colimits of finite direct sums).

If $0\to A\to B\to C\to 0$ is a short exact sequence in $\Hom_{\Cond{\kappa}(\Ab)}$ with $A\in\Sol{\kappa}$ for $S\in\Pro(\Fin)_{\kappa}$ we obtain a commutative diagram 
\begin{center}
        \begin{tikzcd}[cramped, sep=small]
            0\arrow[r]& \Hom(\mathbb Z[\underline{S}_{\kappa}]^{\blacksquare},A)\arrow[r]\arrow[d,"\cong"] & \Hom(\mathbb Z[\underline{S}_{\kappa}]^{\blacksquare},B) \arrow[r]\arrow[d]& \Hom(\mathbb Z[\underline{S}_{\kappa}]^{\blacksquare},C)\arrow[r]\arrow[d] & \Ext^1_{\Cond{\kappa}(\Ab)}(\mathbb Z[\underline{S}_{\kappa}]^{\blacksquare},A)\arrow[r]\arrow[d,"\cong"] & \ldots\\
            0\arrow[r]& \Hom(\mathbb Z[\underline{S}_{\kappa}],A)\arrow[r] & \Hom(\mathbb Z[\underline{S}_{\kappa}],B) \arrow[r]& \Hom(\mathbb Z[\underline{S}_{\kappa}],C)\arrow[r] & \underbrace{\Ext^1_{\Cond{\kappa}(\Ab)}(\mathbb Z[\underline{S}_{\kappa}],A)}_{=0}\arrow[r] & \ldots 
        \end{tikzcd}
\end{center}
By \cref{solidderived}, the right vertical map is an isomorphism and the bottom-right corner is $0$. 
Hence by the Snake lemma, $B\in\Sol{\kappa}$ if and only if $C\in\Sol{\kappa}$, which shows that $\Sol{\kappa}\subseteq \Cond{\kappa}(\Ab)$ is closed under extensions and cofibers.  We explained above that is also closed under filtered colimits and finite coproducts, and hence under arbitrary small colimits, cf.\ \cite[Proposition 4.4.2.6, 4.4.2.7]{highertopostheory}. 
\end{proof}
\begin{cor}\label{underivedsolidification}Suppose $\kappa$ is an uncountable cardinal and $\algebra{R}\in\Alg(\Cond{\kappa}(\Ab))$. 
    \begin{romanenum}
        \item The category $\Sol{\kappa}(\algebra{R})$ is presentable and Grothendieck abelian. 
        \item The forget functor $\Sol{\kappa}(\algebra{R})\to \Cond{\kappa}(\algebra{R})$ admits a left adjoint $(-)^{\solid\algebra{R}}$. 
        \item For $S\in\Pro(\Fin)_{\kappa}$, 
    \[ \mathbb Z[\underline{S}_{\kappa}]^{\solid\mathbb Z}\cong \mathbb Z[\underline{S}_{\kappa}]^{\blacksquare}.\] \item \label{solidificationunderlying}For a condensed $\algebra{R}$-module $M$, the condensed abelian group underlying $M^{\solid\algebra{R}}$ is $M^{\solid\mathbb Z}$ and the $\algebra{R}$-module structure is given by \[ \algebra{R}\times M^{\solid\mathbb Z}\to \algebra{R}^{\solid\mathbb Z}\times M^{\solid\mathbb Z}=(\algebra{R}\times M)^{\solid\mathbb Z}\to M^{\solid\mathbb Z}, \] where the left two maps are induced by the unit of the solidification adjunction on condensed abelian groups and the right map is the solidification of the multiplication map.
    \end{romanenum} 
\end{cor}
\begin{proof}
    By \cref{solidclosedunderlimitscolimitskappa}, $\Sol{\kappa}(\algebra{R})\subseteq \Cond{\kappa}(\algebra{R})$ is closed under limits and colimits. 
    Hence by \cite{Adamek1989}, $\Sol{\kappa}(\algebra{R})$ is presentable and the forget functor $\Sol{\kappa}(\algebra{R})\to\Cond{\kappa}(\algebra{R})$ admits a left adjoint $(-)^{\solid\algebra{R}}$. 
    By \cref{condensedabeliangroupsgrothendieckaxioms}, filtered colimits in $\Cond{\kappa}(\algebra{R})$ are exact and $\oplus_{X\in \Pro(\Fin)_{\kappa}}\algebra{R}[\underline{X}_{\kappa}]$ is a generator. By \cref{solidclosedunderlimitscolimitskappa} and the above, this implies that filtered colimits in $\Sol{\kappa}(\algebra{R})$ are exact and that $(\oplus_{X\in \Pro(\Fin)_{\kappa}}\algebra{R}[\underline{X}_{\kappa}])^{\solid\algebra{R}}$ is a generator, which proves that $\Sol{\kappa}(\algebra{R})$ is Grothendieck abelian. For $S\in\Pro(\Fin)_{\kappa}$, $\mathbb Z[\underline{S}_{\kappa}]^{\blacksquare}\in \Sol{\kappa}$ by \cref{solidificationissolid}. Hence by definition of $\Sol{\kappa}$ and $(-)^{\solid\mathbb Z}$, $\mathbb Z[\underline{S}_{\kappa}]^{\blacksquare}=\mathbb Z[\underline{S}_{\kappa}]^{\solid\mathbb Z}$. 
    It is straightforward to check that the left adjoint $(-)^{\solid\algebra{R}}$ satisfies the description from \ref{solidificationunderlying}.  
\end{proof}

\begin{cor}\label{profinitesolidifiestoprojective}
    For $S\in\Pro(\Fin)_{\kappa}$, $\mathbb Z[\underline{S}_{\kappa}]^{\solid}$ is projective in $\Sol{\kappa}$. 
\end{cor}
\begin{proof}Denote by $f\colon \Sol{\kappa}\to\Cond{\kappa}(\Ab)$ the forget functor. By \cref{solidderived,underivedsolidification}, for $S\in\Pro(\Fin)_{\kappa}$, \[ \Hom_{\Cond{\kappa}(\Ab)}(f\mathbb Z[\underline{S}_{\kappa}]^{\solid},-)\circ f\colon\Sol{\kappa}\to \Ab\] is exact. 
    As $f$ is fully faithful, this shows that $\mathbb Z[\underline{S}_{\kappa}]^{\solid}$ is projective in $\Sol{\kappa}$. 
\end{proof}

\begin{cor}\label{solidenoughprojectiveskappa}
Suppose $\algebra{R}\in \Alg(\Cond{\kappa}(\Ab))$ is a $\kappa$-condensed ring.

For $S\in\Pro(\Fin)_{\kappa}$, $\algebra{R}[\underline{S}_{\kappa}]^{\solid\algebra{R}}$ is projective in $\Sol{\kappa}(\algebra{R})$. 
For every solid $R$-module  $M$, there exists a small family of $\kappa$-light profinite sets $S^i,i\in I$ with an epimorphism $\oplus_{i\in I}\algebra{R}[\underline{S^i}_{\kappa}]^{\solid\algebra{R}}\to M$. 
In particular, $\Sol{\kappa}(\algebra{R})$ has enough projectives. 
\end{cor}
\begin{proof}
    As the forget functor $\Cond{\kappa}(\algebra{R})\to \Cond{}(\Ab)$ reflects small limits and colimits, it follows from \cref{solidclosedunderlimitscolimitskappa} that the forget functor $\Sol{\kappa}(\algebra{R})\to\Sol{\kappa}$ preserves small limits and colimits. 
    Whence it admits a left adjoint $\algebra{R}^{\solid}[-]$ by \cite{Adamek1989}. As $f$ is exact, $\algebra{R}^{\solid}[-]$ preserves projectives, and since both are left adjoint to the forget functor $\Sol{\kappa}(\algebra{R})\to\Cond{\kappa}(\Ab)$, \[(-)^{\solid\algebra{R}}\circ \algebra{R}[-]\cong \algebra{R}^{\solid}[-]\circ (-)^{\solid\mathbb Z}.\] 
    It now follows from \cref{profiniteprojectiveinsolidkappa} that for all $S\in\Pro(\Fin)_{\kappa}$, \[\algebra{R}[\underline{S}_{\kappa}]^{\solid\algebra{R}}=\algebra{R}^{\solid}[(\mathbb Z[\underline{S}_{\kappa}]^{\solid})]\] is projective in $\Sol{\kappa}(\algebra{R})$. 

    Suppose now that $P\in\Sol{\kappa}(\algebra{R})$. Since $\oc{{\Pro(\Fin)_{\kappa}}}{P}$ is a small category, there exists a small family of $\kappa$-light profinite sets $(S^i)_{i\in I}$ and an epimorphism $\sqcup_{i\in I}\underline{S^i}_{\kappa}\to P$ of $\kappa$-condensed sets. This induces an epimorphism $\oplus_{i\in I}\algebra{R}[\underline{S^i}_{\kappa}]\to P$ of $\kappa$-condensed $\algebra{R}$-modules which factors over an epimorphism \[ \oplus_{i\in I}\algebra{R}[\underline{S^i}_{\kappa}]^{\solid\algebra{R}}\to P\] since $P$ is solid.
\end{proof}

\begin{cor}\label{derivedkappasolidmodulespresentable}
    Suppose $\kappa$ is an uncountable cardinal. For $\algebra{R}\in\Alg(\Cond{\kappa}(\Ab))$, 
    $\mathcal D(\Sol{\kappa}(\algebra{R}))$ is presentable and in particular admits small limits and colimits.  
\end{cor}
\begin{proof}
    As $\Sol{\kappa}(\algebra{R})$ is Grothendieck abelian (\cref{underivedsolidification}), its derived category $\mathcal D(\Sol{\kappa}(\algebra{R}))$ is presentable by \cite[Proposition 1.3.5.21]{higheralgebra}. 
    In particular, it admits small limits and colimits by \cite[Definition 5.5.0.1, Corollary 5.5.2.4]{highertopostheory}. 
\end{proof}

\begin{lemma}\label{solidleftcomplete}
The $t$-structure on $\mathcal D(\Sol{\kappa}(\algebra{R}))$ is left and right-complete. 
\end{lemma}
\begin{proof}
    Since $\Sol{\kappa}(\algebra{R})$ is Grothendieck abelian (\cref{underivedsolidification}), the $t$-structure on $\mathcal D(\Sol{\kappa}(\algebra{R}))$ is right-complete by \cite[Proposition 1.3.5.21]{higheralgebra}. 

    By \cref{solidclosedunderlimitscolimitskappa} and \cref{condensedabeliangroupsgrothendieckaxioms}, countable products in $\Sol{\kappa}(\algebra{R})$ are exact. This implies that $\mathcal D(\Sol{\kappa}(\algebra{R}))$ has countable products, and that they can be computed on representing chain complexes. In particular, $\mathcal D(\Sol{\kappa}(\algebra{R}))_{\geq 0}\subseteq\mathcal D(\Sol{\kappa}(\algebra{R}))$ is stable under countable coproducts. 
    Since the $t$-structure on $\mathcal D(\Sol{\kappa}(\algebra{R}))$ is left-separated, it now follows from \cite[Proposition 1.2.1.19]{higheralgebra} that it is left-complete. 
\end{proof}

Since $\Sol{\kappa}(\algebra{R})\to\Cond{\kappa}(\algebra{R})$ is exact, the induced functor $\Ch(\Sol{\kappa}(\algebra{R}))\to \Ch(\Cond{\kappa}(\algebra{R}))$ descends to a functor $f\colon \mathcal D(\Sol{\kappa}(\algebra{R}))\to\mathcal D(\Cond{\kappa}(\algebra{R}))$. We now record basic properties of this functor. 

\begin{cor}\label{derivedforgetfunctorsolidpreservescolimits}
    The functor $f\colon \mathcal D(\Sol{\kappa}(\algebra{R}))\to\mathcal D(\Cond{\kappa}(\algebra{R}))$ preserves small colimits.
\end{cor}
\begin{proof}
    Since coproducts in $\Cond{\kappa}(\algebra{R})$ are exact (\cref{condensedabeliangroupsgrothendieckaxioms}), \cref{solidclosedunderlimitscolimitskappa} implies that they are exact in $\Sol{\kappa}(\algebra{R})$ as well. 
    It follows that coproducts in $\mathcal D(\Sol{\kappa}(\algebra{R}))$ and $\mathcal D(\Cond{\kappa}(\algebra{R}))$ can be computed on representing chain complexes (as degreewise coproducts), and in particular, $f$ preserves coproducts by \cref{solidclosedunderlimitscolimitskappa}. 
    \cref{exactfunctorsinduceexactfunctorsonderivedcats} and \cref{solidclosedunderlimitscolimitskappa} imply that $f$ preserves finite colimits, and hence $f$ preserves small colimits by \cite[Proposition 4.4.2.6]{highertopostheory}. 
\end{proof}
\begin{rem}As $\Cond{\kappa}(\algebra{R})$ and $\Sol{\kappa}(\algebra{R})$ are Grothendieck abelian, by \cite[Proposition C.3.1.1, C.3.2.1, Theorem C.5.4.9]{SAG}, the functor $f$ is determined up to contractible choice by being $t$-exact, cocontinuous and restricting to the forget functor \[\Sol{\kappa}(\algebra{R})\to\Cond{\kappa}(\algebra{R})\] on hearts, see also \cite[Proposition A.2]{cartiermodulesmattisweiss}. \end{rem}
\begin{lemma}\label{derivedsolificationexistsstronglimitcardinal}
If $\algebra{R}$ is a $\kappa$-condensed ring such that $\Cond{\kappa}(\algebra{R})$ has enough projectives, then the forget functor 
\[ \mathcal D(\Sol{\kappa}(\algebra{R}))\to \mathcal D(\Cond{\kappa}(\algebra{R}))\] has a left adjoint. 
\end{lemma}
This in particular applies if $\kappa$ is a strong limit cardinal, cf. \cref{enoughprojectives}. 
\begin{proof}
    By \cref{solidleftcomplete}, the $t$-structure on $\mathcal D(\Sol{(\kappa)}(\algebra{R}))$ is left and right complete and by \cref{solidclosedunderlimitscolimitskappa}, filtered colimits in $\Sol{\kappa}(\algebra{R})$ are exact. 
    \cref{existenceunboundedderivedfunctors} now implies that there exists a unique colimits preserving, right t-exact functor $(-)^{L\solid\algebra{R}}\colon \mathcal D(\Cond{(\kappa)}(\algebra{R}))\to\mathcal D(\Sol{\kappa}(\algebra{R}))$ such that \[H_0\circ (-)^{L\solid\algebra{R}}|_{\Cond{\kappa}(\algebra{R})}=(-)^{\solid\algebra{R}}\] and for all projective $R$-modules $P$, \[P^{L\solid\algebra{R}}\in\mathcal D(\Sol{\kappa}(\algebra{R}))^{\heart}.\] 
    The composite $f\circ (-)^{L\solid\algebra{R}}$ is cocontinuous, right t-exact and carries projectives of $\Cond{\kappa}(\algebra{R})$ to the heart. 
    As $\Cond{\kappa}(\algebra{R})$ has enough projectives (\cref{enoughprojectives}), filtered colimits in $\Cond{\kappa}(\algebra{R})$ are exact (\cref{condensedabeliangroupsgrothendieckaxioms}) and the $t$-structure  $\mathcal D(\Cond{\kappa}(\algebra{R}))$ is left and right complete (\cref{tstructureleftcomplete,modulesinderivedcategoryisderivedcategoryofmodules}),  \cref{existenceunboundedderivedfunctors} implies that the unit $\id\to f_0\circ (-)^{\solid\algebra{R}}$ enhances essentially uniquely to a natural transformation $\id \to f\circ (-)^{L\solid\algebra{R}}$.     
    We claim that this exhibits $f$ as right adjoint to $(-)^{L\solid\algebra{R}}$.
    Denote by $\mathcal D\subseteq \mathcal D(\Cond{\kappa}(\algebra{R}))$ the full subcategory on objects $X$ such that \[ c(\eta,X)\colon \Map_{\mathcal D(\Sol{(\kappa)}(\algebra{R}))}(X^{L\solid\algebra{R}},-)\xrightarrow{f_*} \Map_{\mathcal D(\Cond{(\kappa)}(\algebra{R}))}(fX^{L\solid\algebra{R}},f-)\xrightarrow{\eta^*}\Map_{\mathcal D(\Cond{\kappa}(\algebra{R}))}(X,f-)\] is an equivalence. 
    We want to show that $\mathcal D=\mathcal D(\Cond{(\kappa)}(\algebra{R}))$. 
    As $f$ and $(-)^{L\solid\algebra{R}}$ preserve colimits, $\mathcal D$ is closed under colimits. 
    Since $f$ and $(-)^{L\solid\algebra{R}}$ are exact, $X\in\mathcal D$ if and only if $\Sigma^iX\in\mathcal D$ for some $i\in\mathbb Z$. 
    By right-completeness of the $t$-structure on $\mathcal D(\Cond{(\kappa)}(\Ab))$, it therefore suffices to show that $\mathcal D(\Cond{(\kappa)}(\algebra{R}))_{\geq 0}\subseteq \mathcal D$. 
    Fix $X\in\mathcal D(\Cond{(\kappa)}(\algebra{R}))_{\geq 0}$. For $i\in\mathbb N_0$, \begin{align*}\pi_i\Map_{\mathcal D(\Sol{(\kappa)}(\algebra{R}))}(X^{L\solid\algebra{R}},-)& \cong \pi_0\Map_{\mathcal D(\Sol{(\kappa)}(\algebra{R}))}(X^{L\solid\algebra{R}},\Omega^i-)\\ &\cong \pi_0\Map_{\mathcal D(\Sol{(\kappa)}(\algebra{R}))}(X^{L\solid\algebra{R}},\tau_{\geq 0}\Omega^i-)\end{align*} by right $t$-exactness of $(-)^{L\solid\algebra{R}}$, and 
    \begin{align*}\pi_i\Map_{\mathcal D(\Cond{(\kappa)}(\algebra{R}))}(X,f-)&\cong \pi_0\Map_{\mathcal D(\Cond{(\kappa)}(\algebra{R}))}(X,\Omega^if-)\\ &\cong \pi_0\Map_{\mathcal D(\Cond{(\kappa)}(\algebra{R}))}(X,\tau_{\geq 0}\Omega^if-) \\ 
    &\cong \pi_0\Map_{\mathcal D(\Cond{(\kappa)}(\algebra{R}))}(X,f\tau_{\geq 0}\Omega^i-),\end{align*} by $t$-exactness of $f$.
    As $\tau_{\geq 0}\Omega^i=\Omega^i\tau_{\geq -i}$ and $\Omega^i\colon\mathcal D(\Sol{\kappa}(\algebra{R}))_{>-\infty}\to \mathcal D(\Sol{\kappa}(\algebra{R}))_{>-\infty}$ is an equivalence, this shows that $X\in\mathcal D$ if and only if $\pi_0c(\eta,X)(Y)$ is an equivalence for all $Y\in\mathcal D(\Sol{(\kappa)}(\algebra{R}))_{>-\infty}$. 
    Choose a chain complex $C_*$ of projective $\kappa$-condensed $\algebra{R}$-modules representing $X$ which is concentrated in degrees $\geq 0$. This exists since $\Cond{\kappa}(\algebra{R})$ has enough projectives. 
    Then $X^{L\solid\algebra{R}}$ is computed by applying $(-)^{\solid\algebra{R}}$ degreewise to $C_*$. 
    Indeed: Denote by $C_*^{\leq n}$ the stupid truncations of $C_*$. Exactness of filtered colimits in $\Cond{(\kappa)}(\algebra{R})$ and $\Sol{\kappa}(\algebra{R})$ implies that $X^{L\solid\algebra{R}}\cong \colim{k}(C_*^{\leq k})^{L\solid\algebra{R}}$.
    Using the cofiber sequences \[ C_*^{\leq k}\to C_*^{\leq k+1}\to C_{k+1}[k+1],\] it follows by induction on $k$ that $(C_*^{\leq k})^{L\solid\algebra{R}}$ is computed by applying $(-)^{\solid\algebra{R}}$ degreewise to $C_*^{\leq k}$.   
    The unit for the localisation $\Cond{(\kappa)}(\algebra{R})\to \Sol{(\kappa)}(\algebra{R})$ defines a chain map $C_*\to C_*^{\solid\algebra{R}}$, and the same inductive argument shows that this represents the homotopy class of the unit $\eta_X$. 
    The argument from the proof of \cite[\href{https://stacks.math.columbia.edu/tag/0FNC}{Tag 0FNC}]{stacks-project} now shows that $\pi_0c(\eta,X)(Y)$ is an equivalence for all $Y\in\mathcal D(\Sol{\kappa}(\algebra{R}))_{>-\infty}$. 
\end{proof}
We will show below (\cref{derivedsolidificationwithoutkappa}, \cref{derivedsolidificationexistssflatrings,lightringssflat}), that the left adjoint also exists over $\Cond{}(\algebra{R})$ and for $\kappa=\aleph_1$. 

\begin{lemma}\label{fpreservescounatblelimits}
    The functor $f\colon \mathcal D(\Sol{\kappa}(\algebra{R}))\to\mathcal D(\Cond{\kappa}(\algebra{R}))$ preserves $\kappa$-small limits. 
\end{lemma}
\begin{proof}
Since $\kappa$-small products in $\Cond{\kappa}(\algebra{R})$ are exact (\cref{condensedabeliangroupsgrothendieckaxioms}), they are exact in $\Sol{\kappa}(\algebra{R})$ by \cref{solidclosedunderlimitscolimitskappa}. This implies that $\kappa$-small products in $\mathcal D(\Sol{\kappa}(\algebra{R}))$ can be computed on the level of chain complexes (as degreewise product). 
Since the same holds for $\mathcal D(\Cond{\kappa}(\Ab))$, \cref{solidclosedunderlimitscolimitskappa} implies that $f$ preserves all $\kappa$-small products. Since $f$ is exact, this shows that $f$ preserves all $\kappa$-small limits (\cite[Proposition 4.4.2.7]{highertopostheory}).
\end{proof}

\begin{lemma}\label{derivedsolidforgetfunctorfullyfaithful}
    The forget functor $f\colon \mathcal D(\Sol{\kappa})\to \mathcal D(\Cond{\kappa}(\Ab))$ is fully faithful. 
\end{lemma}

\begin{proof}We want to deduce this from \cref{fullyfaithfulnessderivedfunctors}. 
    By \cref{solidclosedunderlimitscolimitskappa,underivedsolidification} and \cref{condensedabeliangroupsgrothendieckaxioms}, $\Cond{\kappa}(\Ab)$ and $\Sol{\kappa}$ are Grothendieck abelian, countable products in both categories are exact, and the functor $\Sol{\kappa}\to \Cond{\kappa}(\Ab)$ is fully faithful and preserves small limits and colimits. 

    By \cref{solidenoughprojectiveskappa}, $\Sol{\kappa}$ has enough projectives, and every projective object $p\in\Sol{\kappa}$ is a direct summand of $\oplus_{S\in\Pro(\Fin)_{\kappa}}\mathbb Z[\underline{S}_{\kappa}]^{\solid}$. \cref{solidderived} therefore implies that for all solid abelian groups $M\in\Sol{\kappa}$ and for all projectives $p\in\Sol{\kappa}$, \[\Ext^*_{\Cond{\kappa}(\Ab)}(p,M)\cong \Ext^0_{\Cond{\kappa}(\Ab)}(p,M)\] is concentrated in degree $0$. 
    The statement now follows from \cref{fullyfaithfulnessderivedfunctors}. 
\end{proof}
The essential image of the forget functor $\mathcal D(\Sol{\kappa})\hookrightarrow\mathcal D(\Cond{\kappa}(\Ab))$ can be characterized as follows: 
\begin{cor}\label{derivedsolidissolidderived}
    For $M\in\mathcal D(\Cond{\kappa}(\Ab))$, the following are equivalent: 
    \begin{romanenum}
        \item\label{essentialimageforgetfunctorsolid} $M$ lies in the essential image of the forget functor $ \mathcal D(\Sol{\kappa})\to\mathcal D(\Cond{\kappa}(\Ab))$. 
    \item\label{derivedsolid} For all $S\in\Pro(\Fin)_{\kappa}$, the map $\mathbb Z[\underline{S}_{\kappa}]\to\mathbb Z[\underline{S}_{\kappa}]^{\blacksquare}$ induces an equivalence \[\map_{\mathcal D(\Cond{\kappa}(\Ab))}(\mathbb Z[\underline{S}_{\kappa}]^{\blacksquare},M)\cong \map_{\mathcal D(\Cond{\kappa}(\Ab))}(\mathbb Z[\underline{S}_{\kappa}],M).\] 
    \item\label{internallyderivedsolid} For all $S\in\Pro(\Fin)_{\kappa}$, the map $\mathbb Z[\underline{S}_{\kappa}]\to\mathbb Z[\underline{S}_{\kappa}]^{\blacksquare}$ induces an equivalence \[\imap_{\mathcal D(\Cond{\kappa}(\Ab))}(\mathbb Z[\underline{S}_{\kappa}]^{\blacksquare},M)\cong \imap_{\mathcal D(\Cond{\kappa}(\Ab))}(\mathbb Z[\underline{S}_{\kappa}],M).\] 
    \item \label{Gcategoryidentifywithsolid}$M$ lies in the category $\mathcal G$ from \cref{Gcategorysoliddefinition}.
    \item \label{essentialimageforgethomology} For all $i\in\mathbb Z$, $H_i(M)\in\Sol{\kappa}$. 
    \end{romanenum}
    \end{cor}
    \begin{proof}
        Since $\Sol{\kappa}\subseteq \Cond{\kappa}(\Ab)$ preserves limits and colimits (\cref{solidclosedunderlimitscolimitskappa}), $\ref{essentialimageforgetfunctorsolid}\Rightarrow \ref{essentialimageforgethomology}$. 
        Clearly, $\ref{internallyderivedsolid}\Rightarrow \ref{derivedsolid}$, and by \cref{Gcategorysolidissolid}, $\ref{Gcategoryidentifywithsolid}\Rightarrow \ref{internallyderivedsolid}$.

        We now show that $\ref{derivedsolid}\Rightarrow\ref{essentialimageforgetfunctorsolid}$. 
        Denote by $\mathcal D(\Cond{\kappa}(\Ab))_{\solid}\subseteq \mathcal D(\Cond{\kappa}(\Ab))$ the full subcategory on objects satisfying $\ref{derivedsolid}$. 
    Since $\mathcal D(\Ab)\to \Sp$ is conservative, $M\in\mathcal D(\Cond{\kappa}(\Ab))_{\solid}$ if and only if for all $S\in\Pro(\Fin)_{\kappa}, i\in\mathbb N_0$, 
    $\mathbb Z[\underline{S}_{\kappa}]\to\mathbb Z[\underline{S}_{\kappa}]^{\blacksquare}$ induces an equivalence
    \[ \Map_{\mathcal D(\Cond{(\kappa)}(\Ab))}(\Sigma^i\mathbb Z[\underline{S}_{\kappa}]^{\blacksquare},M)\cong \Map_{\mathcal D(\Cond{(\kappa)}(\Ab))}(\Sigma^i\mathbb Z[\underline{S}_{\kappa}],M).\] 
    As $\mathcal D(\Cond{\kappa}(\Ab))$ is presentable, by \cite[Proposition 5.5.4.15]{highertopostheory}, the inclusion \[\mathcal D(\Cond{\kappa}(\Ab))_{\solid}\subseteq \mathcal D(\Cond{\kappa}(\Ab))\] has a left adjoint $L$, and for $S\in\Pro(\Fin)_{\kappa}$, 
    $L(\mathbb Z[\underline{S}_{\kappa}])\cong L(\mathbb Z[\underline{S}_{\kappa}]^{\blacksquare})$. 
    By \cref{Gcategorysolidissolid} and \cref{goodcomplexcohomologicaldimension}, for all $S\in\Pro(\Fin)_{\kappa}$, $\mathbb Z[\underline{S}_{\kappa}]^{\blacksquare}\in\mathcal D(\Cond{\kappa}(\Ab))_{\solid}$, whence $L(\mathbb Z[\underline{S}_{\kappa}])\cong\mathbb Z[\underline{S}_{\kappa}]^{\blacksquare}$. In particular, \[L(\mathbb Z[\underline{S}_{\kappa}])\in \essim(f)\subseteq \mathcal D(\Cond{\kappa}(\Ab))\] by \cref{underivedsolidification}. As the essential image of $f$ is a stable subcategory closed under small colimits (\cref{derivedforgetfunctorsolidpreservescolimits}, \cref{fpreservescounatblelimits}), and $\Sigma^i\mathbb Z[\underline{S}_{\kappa}], i\in\mathbb Z, S\in\Pro(\Fin)_{\kappa}$ generates $\mathcal D(\Cond{\kappa}(\Ab))$ under colimits, it follows that \[\im(L)=\mathcal D(\Cond{\kappa}(\Ab))_{\solid}\subseteq \essim(f).\]

    We now show that $\ref{essentialimageforgetfunctorsolid}\Rightarrow \ref{Gcategoryidentifywithsolid}$. 
    By \cref{Gcategoryheart}, $f(\Sol{\kappa})\subseteq \mathcal G.$ Since $\mathcal G\subseteq \mathcal D(\Cond{\kappa}(\Ab))$ is a stable subcategory closed under extensions and $f$ is $t$-exact, this implies that for $n\leq m\in\mathbb Z$, \[f(\mathcal D(\Sol{\kappa})_{[n,m]})\subseteq \mathcal G.\] In particular, for $n\leq m\in\mathbb Z$ and $A\in\mathcal D(\Sol{\kappa})$, $\tau_{\leq m}\tau_{\geq n}fA\in\mathcal D(\Cond{\kappa}(\Ab))_{\solid}$ since $f$ is $t$-exact.
    Right completeness of the $t$-structure on $\mathcal D(\Cond{\kappa}(\Ab))$ and \cref{Gcategorysolidclosedunderfilteredcolimits} imply that $\tau_{\leq m}A=\tau_{\leq m}\colim{n}\tau_{\geq n}A\in \mathcal G$ for all $m\in\mathbb N_0$. 
    As the $t$-structure on $\mathcal D(\Cond{}(\Ab))$ is left-complete (\cref{tstructureleftcomplete}), and $\mathcal G$ is closed under limits, this implies that $\essim(f)\subseteq \mathcal G$. 

    It remains to show that $\ref{essentialimageforgethomology}\Rightarrow \ref{essentialimageforgetfunctorsolid}$. Suppose $A\in\mathcal D(\Cond{\kappa}(\Ab))$ with $H_i(A)\in\Sol{\kappa}$ for all $i\in\mathbb Z$. By \cref{tstructureleftcomplete}, \[A=\colim{k}\tau_{\geq k}A=\colim{k}\clim{i}\tau_{\leq i}\tau_{\geq k}A.\] 
    As $\essim(f)\subseteq\mathcal D(\Cond{\kappa}(\Ab))$ is closed under colimits and countable limits (\cref{derivedforgetfunctorsolidpreservescolimits}, \cref{fpreservescounatblelimits}), it is enough to show that for all $n\leq m\in\mathbb Z$, $\tau_{\geq n}\tau_{\leq m}A\in \essim(f).$
    Since $\essim(f)$ is a stable subcategory, we can assume that $n=0$.
    The claim now follows by induction on $m$: The case $n=m=0$ is trivial, and the induction step follows since $\essim(f)\subseteq \mathcal D(\Cond{\kappa}(\Ab))$ is a stable subcategory. 
    \end{proof}
    We will from now on freely identify $\mathcal D(\Sol{\kappa})$ with the essential image of the forget functor $\mathcal D(\Sol{\kappa})\hookrightarrow \mathcal D(\Cond{\kappa}(\Ab))$. 
\begin{cor}\label{derivedsolidificationexistsabeliangroupskappa}
    For every uncountable cardinal $\kappa$, $\mathcal D(\Sol{\kappa})\subseteq \mathcal D(\Cond{\kappa}(\Ab))$ is closed under small limits and colimits. 
    The inclusion $\mathcal D(\Sol{\kappa})\subseteq \mathcal D(\Cond{\kappa}(\Ab))$ admits a left adjoint $(-)^{L\solid}$. 
    \end{cor}
    \begin{proof}The third characterization of the essential image of the forget functor $f\colon \mathcal D(\Sol{\kappa})\to\mathcal D(\Cond{\kappa}(\Ab))$ from \cref{derivedsolidissolidderived} implies that $\mathcal D(\Sol{\kappa})\subseteq \mathcal D(\Cond{\kappa}(\Ab))$ is closed under small limits. By \cref{derivedforgetfunctorsolidpreservescolimits}, it is closed under small colimits. As $\Cond{\kappa}(\Ab), \Sol{\kappa}$ are Grothendieck abelian (\cref{underivedsolidification}, \cref{condensedmodulesgrothendieckabelian}), their derived categories are presentable by \cite[Proposition 1.3.5.21]{higheralgebra}. 
        It now follows from the adjoint functor theorem (\cite[Corollary 5.5.2.9]{highertopostheory}) that $\mathcal D(\Sol{\kappa})\subseteq \mathcal D(\Cond{\kappa}(\Ab))$ has a left adjoint $(-)^{L\solid}$. 
    \end{proof}

\begin{cor}\label{profiniteprojectiveinsolidkappa}
\begin{romanenum}
\item For $S\in \Pro(\Fin)_{\kappa}$, $\mathbb Z[\underline{S}_{\kappa}]^{L\solid}\cong\mathbb Z[\underline{S}_{\kappa}]^{\solid}$ and $\mathbb Z[\underline{S}_{\kappa}]^{\solid}$ is projective in $\Sol{\kappa}$. 
\item $\underline{\mathbb R}_{\kappa}^{L\solid}=0$.
\end{romanenum}
\end{cor} 
\begin{proof}
By \cref{Gcategorysolidissolid} and \cref{Gcategorysolidcontainsgoodcomplexes}, for $G\in\mathcal G$, and $S\in\Pro(\Fin)_{\kappa}$, \[ \map_{\mathcal D(\Cond{\kappa}(\Ab))}(\mathbb Z[\underline{S}_{\kappa}],G)\cong \map_{\mathcal D(\Cond{\kappa}(\Ab))}(\mathbb Z[\underline{S}_{\kappa}]^{\blacksquare},G).\] 
As the forget functor factors over an equivalence $\mathcal D(\Sol{\kappa})\cong \mathcal G$ (\cref{derivedsolidissolidderived}) and $\mathbb Z[\underline{S}_{\kappa}]^{\blacksquare}\in\mathcal G$ (\cref{Gcategorysolidcontainsgoodcomplexes}), this implies that $\mathbb Z[\underline{S}_{\kappa}]^{L\solid}\cong\mathbb Z[\underline{S}_{\kappa}]^{\blacksquare}\cong\mathbb Z[\underline{S}_{\kappa}]^{\solid}$. 
By \cref{profinitesolidifiestoprojective}, $\mathbb Z[\underline{S}_{\kappa}]^{\solid}$ is projective in $\Sol{\kappa}$. 
By \cref{Gcategorysolidissolid}, for $G\in\mathcal G$, $\Map_{\mathcal D(\Cond{\kappa}(\Ab))}(\underline{\mathbb R}_{\kappa},G)=*,$ whence $\underline{\mathbb R^{L\solid}}_{\kappa}=0$ by \cref{derivedsolidissolidderived}. 
 \end{proof}

We now deduce from the above results their condensed analogues. This relies on the following observation. 
\begin{lemma}\label{colimitdescriprionsolid}
    Suppose $\algebra{R}\in\Alg(\Cond{}(\Ab))$ and choose a regular cardinal $\mu$ with $\algebra{R}\in \Alg(\Cond{\mu}(\Ab))$. 
    \begin{romanenum}
        \item For regular cardinals $\lambda\geq\kappa\geq \mu$, 
    \[ \Cond{\kappa}(\algebra{R})\hookrightarrow \Cond{\lambda}(\algebra{R})\] restricts to a functor 
    \[\Sol{\kappa}(\algebra{R})\hookrightarrow \Sol{\lambda}(\algebra{R}).\] 
    \item This exhibits
    \[\colim{\substack{\lambda\geq \mu\\ \lambda\text{ regular }}}\Sol{\lambda}(\algebra{R})\cong \Sol{}(\algebra{R})\] as colimit in $\vlCat$.  
    \item For $M\in\Cond{}(\algebra{R})$, $M\in\Sol{}(\algebra{R})$ if and only if for all regular cardinals $\kappa\geq \mu$, $r^{\kappa}M\in \Sol{\kappa}(\algebra{R})$.
    \item For $\kappa\geq \mu$, 
    \[\Sol{\kappa}(\algebra{R})\cong \Cond{\kappa}(\algebra{R})\times_{\Cond{}(\algebra{R})}\Sol{}(\algebra{R})\cong \Cond{\kappa}(\Ab)\times_{\Cond{}(\Ab)}\Sol{}(\algebra{R}).\] 
    \item $\Sol{}(\algebra{R})$ has small limits and colimits and for $\kappa\geq \mu$, 
     $\Sol{\kappa}(\algebra{R})\subseteq \Sol{}(\algebra{R})$ is closed under small colimits and $\kappa$-small limits. 
    \end{romanenum}  
\end{lemma}

\begin{proof}
    We first show that for regular cardinals $\kappa\leq \lambda$, the left Kan extension \[i_{\kappa}^{\lambda}\colon \Cond{\kappa}(\Ab)\to \Cond{\lambda}(\Ab)\] restricts to $\Sol{\kappa}\to \Sol{\lambda}$.
By \cref{solidclosedunderlimitscolimitskappa}, $\Sol{\kappa}\subseteq \Cond{\kappa}(\Ab)$ is closed under small limits and colimits. In particular, $\prod_{I}\mathbb Z\in\Sol{\kappa}$ for all small sets $I$. \cref{solidificationissolid} and \cref{Nobelingspecker} imply that $\Sol{\kappa}$ is generated by $\prod_{I}\mathbb Z, |I|< \kappa$ under small colimits. 
For all regular cardinals $\lambda\geq\kappa$, the functor $i^{\lambda}_{\kappa}\colon \Cond{\kappa}(\Ab)\to\Cond{\lambda}(\Ab)$ preserves $\kappa$-small limits by \cref{condensedextremallydisconnected1}. 
As $i^{\lambda}_{\kappa}$ also preserves colimits and $\Sol{\kappa/\lambda}\subseteq \Cond{\kappa/\lambda}(\Ab)$ is closed under small colimits and limits (\cref{solidclosedunderlimitscolimitskappa}), this shows that $i^{\lambda}_{\kappa}(\Sol{\kappa})\subseteq \Sol{\lambda}$. 
It now follows from \cref{filteredcolimitsofcategoriesappendix} that \[\colim{\substack{\kappa\\ \kappa\text{ regular }}}\Sol{\kappa}(\mathbb Z)\subseteq \Cond{}(\Ab)\] is a full subcategory.
We claim that it consists precisely of the solid abelian groups. 
Fix a regular cardinal $\kappa$ and denote by $i_{\kappa}\colon\Cond{\kappa}(\Ab)\hookrightarrow\Cond{}(\Ab)$ the left adjoint. \cref{topspacecomefromregular} implies that for $S\in\Pro(\Fin)_{\kappa}$, $i_{\kappa}\mathbb Z[\underline{S}_{\kappa}]=\mathbb Z[\underline{S}]$. 
Since $i_{\kappa}$ preserves $\kappa$-small limits (\cref{condensedextremallydisconnected1}),
$(i_{\kappa}\mathbb Z[\underline{S}_{\kappa}])^{\blacksquare}\cong i_{\kappa}(\mathbb Z[\underline{S}_{\kappa}])^{\blacksquare}$ by \cref{Nobelingspecker}. 
This implies that $M\in\Sol{\kappa}$ if and only if for all regular cardinals $\kappa$, $r^{\kappa}M\in\Sol{\kappa}$, and hence $\colim{\substack{\kappa\\ \kappa\text{ regular }}}\Sol{\kappa}=\Sol{}\subseteq \Cond{}(\Ab)$. 

Suppose now that $\algebra{R}\in \Alg(\Cond{}(\Ab))$ and choose an uncountable regular cardinal $\mu$ such that $\algebra{R}\in\Alg(\Cond{\mu}(\Ab))\subseteq \Alg(\Cond{}(\Ab))$. 
As for all cardinals $\kappa\geq \mu$, 
\[\Sol{\kappa}(\algebra{R})\cong \LMod{\algebra{R}}{\Cond{\kappa}(\Ab)}\times_{\Cond{\kappa}(\Ab)}\Sol{\kappa}\cong \LMod{\algebra{R}}{\Cond{}(\Ab)}\times_{\Cond{}(\Ab)}\Sol{\kappa}, \] it follows from the above that for all $\lambda\geq \kappa\geq \mu$, $\Cond{\kappa}(\algebra{R})\hookrightarrow \Cond{\lambda}(\algebra{R})$ restricts to a functor $\Sol{\kappa}(\algebra{R})\to \Sol{\lambda}(\algebra{R})$ and that \[\colim{\substack{\kappa\geq \mu\\ \kappa\text{ regular}}}\Sol{\kappa}(\algebra{R})=\Sol{}(\algebra{R})\subseteq \LMod{\algebra{R}}{\Cond{}(\Ab)}.\] 
In particular, the induced functors $\Sol{\kappa}(\algebra{R})\to \Sol{}(\algebra{R})$ are fully faithful with essential image \[\Cond{\kappa}(\Ab)\times_{\Cond{}(\Ab)}\Sol{}(\algebra{R})\cong \Cond{\kappa}(\algebra{R})\times_{\Cond{}(\algebra{R})}\Sol{}(\algebra{R}).\] 

We have shown above that for $M\in\Cond{}(\algebra{R})$, $M\in\Sol{}(\algebra{R})$ if and only if for all regular cardinals $\kappa\geq \mu$, $r^{\kappa}M\in \Sol{\kappa}(\algebra{R})$. 

Finally,
\cref{solidclosedunderlimitscolimitskappa} and \cref{condensedextremallydisconnected1,forgetfreeadjunctionmodules} imply that for all regular cardinals $\mu\leq \kappa\leq \lambda$,\[\Sol{\kappa}(\algebra{R})\to\Sol{\lambda}(\algebra{R})\] preserves small colimits and $\kappa$-small limits. It follows from \cref{filteredcolimitsofcategoriesappendix} that $\Sol{}(\algebra{R})$ has small limits and colimits and that $\Sol{\kappa}(\algebra{R})\subseteq \Sol{}(\algebra{R})$ is closed under $\kappa$-small limits and small colimits. 
\end{proof}

\begin{cor}\label{underivedsolidificationwithoutkappa}
        Suppose $\algebra{R}\in \Alg(\Cond{}(\Ab))$ is a discrete condensed ring. 
    \begin{romanenum}
    \item \label{solidclosedunderlimitswithoutkappa}The category \[\Sol{}(\algebra{R})\subseteq \Cond{}(\algebra{R})\] is an abelian subcategory stable under small limits, colimits and extensions. In particular, $\Sol{}(\algebra{R})$ has all small colimits and limits.
    \item The forget functor $\Sol{}(\algebra{R})\to \Cond{}(\algebra{R})$ admits a left adjoint $(-)^{\solid\algebra{R}}$.
    \item For a condensed $\algebra{R}$-module $M$, the condensed abelian group underlying $M^{\solid\algebra{R}}$ is $M^{\solid\mathbb Z}$ and the $\algebra{R}$-module structure is given by \[ \algebra{R}\times M^{\solid\mathbb Z}\to \algebra{R}^{\solid\mathbb Z}\times M^{\solid\mathbb Z}=(\algebra{R}\times M)^{\solid\mathbb Z}\to M^{\solid\mathbb Z}, \] where the left two maps are induced by the unit of the solidification adjunction on for condensed abelian groups and the right map is the solidification of the scalar multiplication. 
    \end{romanenum}
\end{cor}
\begin{proof}
    \cref{solidclosedunderlimitscolimitskappa,colimitdescriprionsolid,filteredcolimitsofcategoriesappendix} imply that $\Sol{}(\algebra{R})\subseteq \Cond{}(\algebra{R})$ is closed under small limits, colimits and extensions and has small limits and colimits. 

    We now show that the forget functor $\Sol{}(\algebra{R})\to \Cond{}(\algebra{R})$ has a left adjoint. Fix a regular cardinal $\mu$ with $\algebra{R}\in\Alg(\Cond{\mu}(\Ab))$. 
As \[\Sol{}(\algebra{R})\cong \colim{\substack{\kappa\geq \mu\\ \kappa\text{ regular}}}\Sol{\kappa}(\algebra{R}), \] 
it suffices to show that for all regular cardinals $\lambda\geq\kappa\geq \mu$, the mate of the commutative diagram \begin{center}\begin{tikzcd}
    \Sol{\kappa}(\algebra{R})\arrow[d,"i_{\kappa,R}^{\lambda,\solid}"]\arrow[r] & \Cond{\kappa}(\algebra{R})\arrow[d,"i_{\kappa,R}^{\lambda}"]\\ 
    \Sol{\lambda}(\algebra{R})\arrow[r] & \Cond{\lambda}(\algebra{R})
\end{tikzcd}\end{center}commutes, i.e.\ 
\begin{align}\label{Beckchevalley}(-)^{\solid \algebra{R},\lambda}\circ i_{\kappa,R}^{\lambda}\cong i_{\kappa,R}^{\lambda,\solid}\circ (-)^{\solid\algebra{R},\kappa}, \end{align} is an equivalence, then it follows from \cref{adjunctionsbigpresentable} that $\colim{\lambda}(-)^{\solid\algebra{R},\lambda}$ is left adjoint to the forget functor.
As $\Sol{\kappa}(\algebra{R})\hookrightarrow \Cond{\kappa}(\algebra{R})$ is fully faithful, (\ref{Beckchevalley}) holds if and only if  the unit \[\id\to f_{\kappa}\circ (-)^{\solid\algebra{R},\kappa}\] induces an equivalence \begin{align}\label{replacementbeckchevalley}(-)^{\solid\algebra{R},\lambda}\circ i_{\kappa,R}^{\lambda}\cong (-)^{\solid\algebra{R},\lambda}\circ i_{\kappa,R}^{\lambda}\circ f_{\kappa}\circ (-)^{\solid\algebra{R},\kappa}.\end{align}
 As the forget functors $g_{\kappa/\lambda}^{\solid}\circ \Sol{\kappa/\lambda}(\algebra{R})\to \Sol{\kappa/\lambda}$ are conservative, \[g_{\kappa/\lambda}\circ (-)^{\solid\algebra{R},\kappa/\lambda}=(-)^{\solid\algebra{\mathbb Z},\kappa/\lambda}\ \text{ (by \cref{underivedsolidification})},\]
 \[g_{\lambda}\circ i_{\kappa,R}^{\lambda}=i_{\kappa,\mathbb Z}^{\lambda}\circ g_{\kappa}\, (\text{by monoidality of } i_{\kappa}^{\lambda}),\] and \[g_{\kappa}^{\solid}\circ f_{\kappa}=f_{\kappa}\circ g_{\kappa},\] it suffices to show \ref{replacementbeckchevalley} for $\algebra{R}=\mathbb Z$. 
 At $\mathbb Z[\underline{S}_{\kappa}], S\in\Pro(\Fin)_{\kappa}$, \ref{replacementbeckchevalley} holds by \cref{profiniteprojectiveinsolidkappa,underivedsolidification}. 
 As all functors in (\ref{replacementbeckchevalley}) are cocontinuous, this implies the statement for $\algebra{R}=\mathbb Z$, and hence for general $R$.  

By construction, the left adjoint $(-)^{\solid\algebra{R}}$ restricts to $(-)^{\solid\algebra{R}, \kappa}\colon \Cond{\kappa}(\algebra{R})\to \Sol{\kappa}(\algebra{R})$ for all regular cardinals $\kappa\geq \mu$, whence it follows from \cref{underivedsolidification}, that for a condensed $\algebra{R}$-module $M$, 
the abelian group underlying $M^{\solid\algebra{R}}$ is $M^{\solid\mathbb Z}$ with scalar multiplication given by \[ \algebra{R}\times M^{\solid\mathbb Z}\to \algebra{R}^{\solid\mathbb Z}\times M^{\solid\mathbb Z}=(\algebra{R}\times M)^{\solid\mathbb Z}\to M^{\solid\mathbb Z}.\qedhere\] 
\end{proof}
Together with \cref{topspacecomefromregular}, this implies: 
\begin{cor}\label{solidificationoffreemodules}
    Suppose $R\in \Alg(\Cond{}(\Ab))$ and $\kappa$ is a regular cardinal with \[R\in \Alg(\Cond{\kappa}(\Ab))\subseteq \Alg(\Cond{}(\Ab)).\] 
    
    For $S\in\Pro(\Fin)_{\kappa}$, the left adjoint $i_{\kappa}^{\solid}\colon \Sol{\kappa}(\algebra{R})\to \Sol{}(\algebra{R})$ sends $\algebra{R}[\underline{S}_{\kappa}]^{\solid\algebra{R}}$ to $\algebra{R}[\underline{S}]^{\solid\algebra{R}}$. 
\end{cor}
\begin{proof}We have shown in the proof of \cref{underivedsolidificationwithoutkappa} that $i_{\kappa}^{\solid}(\algebra{R}[\underline{S}_{\kappa}]^{\solid\algebra{R}})=(i_{\kappa}\algebra{R}[\underline{S}_{\kappa}])^{\solid\algebra{R}}$, where $i_{\kappa}\colon\Cond{\kappa}(\algebra{R})\to \Cond{}(\algebra{R})$ is the functor induced by the symmetric monoidal functor $\Cond{\kappa}(\Ab)\to \Cond{}(\Ab)$. 
As their right adjoints are equivalent, \[i_{\kappa}\circ R[-]\cong R[-]\circ (\Cond{\kappa}(\Set)\to \Cond{}(\Set)),\] whence by \cref{topspacecomefromregular}, $i_{\kappa}\algebra{R}[\underline{S}_{\kappa}]=\algebra{R}[\underline{S}]$. 
\end{proof}
\begin{cor}\label{solkappafullyfaithfulonbderivedandpreservescolimits}
    For $\algebra{R}\in\Alg(\Cond{}(\Ab))$ choose a regular cardinal $\mu$ with $\algebra{R}\in \Alg(\Cond{\mu}(\Ab))$. 
    \begin{romanenum}
    \item Then \[ \colim{\substack{\kappa\geq \mu\\ \kappa\text{ regular }}}\mathcal D(\Sol{\kappa}(\algebra{R}))\cong \mathcal D(\Sol{}(\algebra{R})).\] 
    \item For regular cardinals $\kappa\geq \mu$, 
    \[ \mathcal D(\Sol{\kappa}(\algebra{R}))\to \mathcal D(\Sol{}(\algebra{R}))\] is fully faithful and preserves colimits and $\kappa$-small limits. 
    \item The category $\mathcal D(\Sol{}(\algebra{R}))$ is big presentable and has small colimits and small limits. 
    \end{romanenum}
\end{cor}
\begin{proof}The first statement follows from  \cref{colimitdescriprionsolid} and 
\cref{Filteredcolimitderivedcategories}.
    Recall from \cref{derivedkappasolidmodulespresentable} that $\mathcal D(\Sol{\kappa}(\algebra{R}))$ is presentable for all regular cardinals $\kappa\geq \mu$. 
Since \[\mathcal D(\Sol{}(\algebra{R}))\cong \colim{\substack{\kappa\geq \mu\\ \kappa\text{ regular }}}\mathcal D(\Sol{\kappa}(\algebra{R})), \] it suffices to show that for all regular cardinals $\lambda\geq\kappa\geq \mu$, the derived functor \[\mathcal D(\Sol{\kappa}(\algebra{R}))\to\mathcal D(\Sol{\lambda}(\algebra{R}))\] of $i^{\lambda}_{\kappa}\colon \Sol{\kappa}(\algebra{R})\hookrightarrow \Sol{\lambda}(\algebra{R})$ is fully faithful and preserves colimits and $\kappa$-small limits, then the second and third statement follow from \cref{filteredcolimitsofcategoriesappendix}. 

We first deduce fully faithfulness from \cref{fullyfaithfulnessderivedfunctors}. 
The functor $i^{\lambda}_{\kappa}\colon \Sol{\kappa}(\algebra{R})\hookrightarrow \Sol{\lambda}(\algebra{R})$ is fully faithful and preserves $\kappa$-small limits and colimits by \cref{solidclosedunderlimitscolimitskappa,forgetfreeadjunctionmodules,condensedextremallydisconnected1}. 
We claim that it also preserves projective objects. 
By \cref{solidenoughprojectiveskappa}, every projective in $\Sol{\kappa}(\algebra{R})$ is a direct summand of $\oplus_{S\in\Pro(\Fin)_{\kappa}}\algebra{R}[\underline{S}_{\kappa}]^{\solid\algebra{R}}$. As $i^{\lambda}_{\kappa}$ preserves colimits, it therefore suffices to show that $i^{\lambda}_{\kappa}(\algebra{R}[\underline{S}_{\kappa}]^{\solid\algebra{R}})$ is projective in $\Sol{\kappa}(\algebra{R})$. 
\cref{solidificationoffreemodules} implies that for all $S\in\Pro(\Fin)_{\kappa}$, $i^{\lambda,\solid}_{\kappa}(\algebra{R}[\underline{S}_{\kappa}]^{\solid\algebra{R}})=\algebra{R}[\underline{S}_{\lambda}]^{\solid\algebra{R}}$, which is projective in $\Sol{\lambda}(\algebra{R})$ by \cref{profinitesolidifiestoprojective}.
As $\Sol{\kappa}(\algebra{R})$ has enough projectives (\cref{solidenoughprojectiveskappa}), it now follows from \cref{fullyfaithfulnessderivedfunctors} that $\mathcal D(i^{\lambda}_{\kappa})$ is fully faithful. 

By \cref{solidclosedunderlimitscolimitskappa,condensedabeliangroupsgrothendieckaxioms}, $\kappa$-small products and small coproducts exist and are exact in $\Sol{\kappa}(\algebra{R})$ and $\Sol{\lambda}(\algebra{R})$. 
In particular, $\kappa$-small products and small coproducts in their derived categories exist and can be computed as degreewise products/coproducts of representing chain complexes, respectively. 
\cref{colimitdescriprionsolid} implies that $\Sol{\kappa}(\algebra{R})\to\Sol{\lambda}(\algebra{R})$ preserves $\kappa$-small products and coproducts, whence $\mathcal D(i^{\lambda}_{\kappa})$ preserves $\kappa$-small products and small coproducts. As $\mathcal D(i^{\lambda}_{\kappa})$ is exact, it follows from \cite[Proposition 4.4.2.7]{highertopostheory} that $\mathcal D(i^{\lambda}_{\kappa})$ preserves $\kappa$-small limits and small colimits. 
\end{proof}
\begin{cor}\label{derivedsolidificationwithoutkappa}
    Suppose $\algebra{R}\in \Alg(\Cond{}(\Ab))$ is a condensed ring. 
    \begin{romanenum}
    \item \label{derivedforgetfunctor}The forget functor extends uniquely to a $t$-exact, small colimits preserving functor \[f\colon\mathcal D(\Sol{}(\algebra{R}))\to\mathcal D(\Cond{}(\algebra{R})).\] 

    \item The functor $f\colon \mathcal D(\Sol{}(\algebra{R}))\to \mathcal D(\Cond{}(\algebra{R}))$ admits a left adjoint $(-)^{L\solid\algebra{R}}$. \end{romanenum} 
\end{cor} 

\begin{proof}We first prove the first statement. Choose a regular cardinal $\mu$ with $\algebra{R}\in\Alg(\Cond{\mu}(\Ab))$. \cref{solkappafullyfaithfulonbderivedandpreservescolimits} implies that \[\Fun^{\operatorname{colim}, \operatorname{t-ex}}(\mathcal D(\Sol{}(\algebra{R})), \mathcal D(\Cond{}(\algebra{R})))\cong \clim{\substack{\kappa\geq \mu\\ \kappa\text{ regular}}} \Fun^{\operatorname{colim}, \operatorname{t-ex}}(\mathcal D(\Sol{\kappa}(\algebra{R})), \mathcal D(\Cond{}(\algebra{R})))\]  where $\Fun^{\operatorname{colim}, \operatorname{t-ex}}$ refers to small colimits preserving, $t$-exact functors. By \cref{colimitdescriprionsolid} and \ref{solidclosedunderlimitscolimitskappa}, 
\[ \Fun^{\operatorname{colim}, \operatorname{ex}}(\Sol{}(\algebra{R}), \Cond{}(\algebra{R}))\cong \clim{\substack{\kappa\geq \mu\\ \kappa\text{ regular}}}\Fun^{\operatorname{colim}, \operatorname{ex}}(\Sol{\kappa}(\algebra{R}), \Cond{}(\algebra{R})), \] where $\Fun^{\operatorname{colim}, \operatorname{t-ex}}$ refers to small colimits preserving, exact functors. 
This implies that  \begin{align*}&\Fun^{\operatorname{colim}, \operatorname{t-ex}}(\mathcal D(\Sol{}(\algebra{R})), \mathcal D(\Cond{}(\algebra{R})))\times_{\Fun^{\operatorname{colim}, \operatorname{t-ex}}(\Sol{}(\algebra{R}), \Cond{}(\algebra{R}))}\{ f\} \\ &\cong \clim{\lambda\geq \mu} \Bigl(\Fun^{\operatorname{colim}, \operatorname{t-ex}}(\mathcal D(\Sol{\lambda}(\algebra{R})), \mathcal D(\Cond{}(\algebra{R})))\times_{\Fun^{\operatorname{colim}, \operatorname{ex}}(\Sol{\lambda}(\algebra{R}), \Cond{}(\algebra{R}))}\{ f|_{\Sol{\lambda}(\algebra{R})}\}\Bigr)\\ & \cong \clim{\lambda \geq \mu}\Bigl(\Fun^{\operatorname{colim}, \operatorname{t-ex}}(\mathcal D(\Sol{\lambda}(\algebra{R})), \mathcal D(\Cond{}(\algebra{R})))\times_{\Fun^{\operatorname{colim}, \operatorname{ex}}(\Sol{\lambda}(\algebra{R}), \Cond{\lambda}(\algebra{R}))}\{ f_{\lambda}\}\Bigr), \end{align*} where the limits run over all regular cardinals $\lambda\geq\mu$ and $f_{(\lambda)}\colon \Sol{(\lambda)}(\algebra{R})\to \Cond{(\lambda)}(\algebra{R})$ denotes the forget functor. 
Since $\Cond{\lambda}(\algebra{R})$ and $\Sol{\lambda}(\algebra{R})$ are Grothendieck abelian (\cref{underivedsolidification}, \cref{condensedmodulesgrothendieckabelian}), the right-hand side is contractible by \cite[Proposition C.3.1.1, C.3.2.1, Theorem C.5.4.9]{SAG}/\cite[Proposition A.2]{cartiermodulesmattisweiss}. This proves \ref{derivedforgetfunctor}. 
 
As $\Cond{}(\algebra{R})$ has enough projectives, \cref{existencetotalderivedfunctorsonaccessiblesheaves1cat} implies that 
$F\mapsto H_0\circ F|_{\Sol{}(\algebra{R})}$ defines an equivalence \begin{align}\label{existencetotalderivedfunctorscondensed}\Fun^{\operatorname{colim},\operatorname{p},\operatorname{r-t-ex}}(\mathcal D(\Cond{}(\algebra{R})),\mathcal D(\Sol{}(\algebra{R})))\cong \Fun^{\operatorname{colim}}(\Cond{}(\algebra{R}),\Sol{}(\algebra{R})),\end{align} between small colimits preserving, right $t$-exact functors $\mathcal D(\Cond{}(\algebra{R}))\to \mathcal D(\Sol{}(\algebra{R}))$ which map projectives to the heart, and small colimits preserving functors $\Cond{}(\algebra{R})\to \Sol{}(\algebra{R})$. 
($\mathcal D(\Sol{}(\algebra{R}))$ satisfies the conditions of \cref{existencetotalderivedfunctorsonaccessiblesheaves1cat} by \cref{solidleftcomplete,solkappafullyfaithfulonbderivedandpreservescolimits,underivedsolidification} and \cite[Theorem 1.3.5.21]{higheralgebra}.) 

Denote by $(-)^{L\solid\algebra{R}}$ a functor in the fiber of $(-)^{\solid\algebra{R}}$. We claim that this is left adjoint to $f$.  

 We showed above that $f$ is $t$-exact and preserves small colimits.
As $\mathcal D(\Cond{}(\Ab))$ is stable, big presentable, has small colimits and its $t$-structure is left and right complete (\cref{tstructureleftcomplete}) and compatible with filtered colimits (\cref{condensedabeliangroupsgrothendieckaxioms}),  \cref{existencetotalderivedfunctorsonaccessiblesheaves1cat} implies that the unit of the localisation $\Sol{}(\algebra{R})\to\Cond{}(\algebra{R})$ enhances essentially uniquely to a natural transformation $\id \to f\circ (-)^{L\solid\algebra{R}}$.
The same argument as in the proof of \cref{derivedsolificationexistsstronglimitcardinal} shows that this exhibits $(-)^{L\solid\algebra{R}}$ as left adjoint to $f$. 
\end{proof}
\begin{definition}Suppose $\algebra{R}\in\Alg(\Cond{(\kappa)}(\Ab))$. 
    If it exists, we call the left adjoint \[(-)^{L\solid\algebra{R}}\colon\mathcal D(\Cond{}(\algebra{R}))\to \mathcal D(\Sol{}(\algebra{R}))\] of the forget functor \emph{derived solidification}.
We write $(-)^{L\solid}\coloneqq (-)^{L\solid\mathbb Z}$ for the derived solidification over abelian groups, which exists by \cref{derivedsolidificationexistsabeliangroupskappa,derivedsolidificationwithoutkappa}.
\end{definition}
\begin{rem}\label{derivedsolidificationisderivedfunctorofsolidification}
\begin{enumerate}    
\item As the forget functor $\mathcal D(\Sol{(\kappa)}(\algebra{R}))\to\mathcal D(\Cond{(\kappa)}(\algebra{R}))$ is the derived functor of $\Sol{(\kappa)}(\algebra{R})\subseteq \Cond{(\kappa)}(\algebra{R})$, the derived solidification defined above agrees with the derived functor of solidification at points where the latter is defined by \cite[\href{https://stacks.math.columbia.edu/tag/0FND}{Tag 0FND}]{stacks-project}.
 
\item Suppose that $C_*$ is a bounded-below chain complex of $(-)^{L\solid\algebra{R}}$-acyclic ($\kappa$-)condensed $\algebra{R}$-modules, i.e. $C_n^{L\solid\algebra{R}}\cong C_n^{\solid\algebra{R}}$ for all $n\in\mathbb Z$. 
Then $C_*^{L\solid\algebra{R}}\cong C_*^{\solid\algebra{R}}$ is computed by applying $(-)^{\solid\algebra{R}}$ degreewise to $C_*$. 
Indeed: Using the cofiber sequences $C_*^{\leq m}\to C_*^{\leq m+1}\to C_{m+1}[m+1]$, it follows by induction on $m$ that $(C_*^{\leq m})^{L\solid\algebra{R}}=(C_*^{\solid\algebra{R}})^{\leq m}$ is computed by applying $(-)^{\solid{\algebra{R}}}$ degreewise to $C_*^{\leq m}$. Hence by exactness of filtered colimits in $\Cond{(\kappa)}(\algebra{R})$ and $\Sol{(\kappa)}(\algebra{R})$, \[(C_*)^{L\solid\algebra{R}}\cong \colim{m}(C_*^{\leq m})^{L\solid\algebra{R}}\cong \colim{m}(C_*^{\solid\algebra{R}})^{\leq m}\cong C_*^{\solid\algebra{R}}\] is computed by applying $(-)^{\solid\algebra{R}}$ degreewise to $C_*$. 

\item This implies the following: 
Suppose that $R$ is such that for all $M\in\Cond{(\kappa)}(\algebra{R})$, there exists a $(-)^{L\solid}$-acyclic $R$-module $A$ with an epimorphism $A\to M$. Then the derived functor of solidification $(-)^{\solid\algebra{R}}$ is everywhere defined by \cite[\href{https://stacks.math.columbia.edu/tag/0794}{Tag 0794}]{stacks-project} (in the large universe), applied to the class of $(-)^{L\solid\algebra{R}}$-acyclic objects \[\mathcal P=\{ A\in\Cond{(\kappa)}(\algebra{R})\, |\, A^{L\solid\algebra{R}} =A^{\solid\algebra{R}}\in\mathcal D(\Sol{(\kappa)})^{\heart}\}.\]\end{enumerate} 
\end{rem}

\begin{lemma}\label{freeabeliangroupsacyclicforsolidificationwithoutkappa}
    \vspace{0pt}\noindent
    \begin{romanenum}
        \item The forget functor $\mathcal D(\Sol{})\to\mathcal D(\Cond{}(\Ab))$ is fully faithful. 
        \item \label{freeabeliangrouponprofiniteacyclicforsolidfication}For $S\in\Pro(\Fin)$, $\mathbb Z[\underline{S}]^{L\solid}\cong \mathbb Z[\underline{S}]^{\solid}\cong \prod_{I}\mathbb Z$ for some small set $I$ of cardinality $|I|\leq \wt(S)$.
    \item \label{solidificationfactors}For all uncountable regular cardinals $\kappa$, the derived solidification  \[(-)^{L\solid}\colon\mathcal D(\Cond{}(\Ab))\to\mathcal D(\Sol{})\] restricts to \[(-)^{L\solid,\kappa}\colon\mathcal D(\Cond{\kappa}(\Ab))\to\mathcal D(\Sol{\kappa})\subseteq \mathcal D(\Sol{}),\] where $(-)^{L\solid,\kappa}$ denotes the derived solidification on the $\kappa$-condensed level.
    \item In particular, for all uncountable regular cardinals $\kappa$, \[\mathcal D(\Sol{\kappa})\cong \mathcal D(\Sol{})\times_{\mathcal D(\Cond{}(\Ab))}\mathcal D(\Cond{\kappa}(\Ab)).\]
\end{romanenum}
\end{lemma}
\begin{proof}
     By construction (see the proof of \cref{derivedsolidificationwithoutkappa}), the forget functor is the colimit of the forget functors $\mathcal D(\Sol{\kappa})\to\mathcal D(\Cond{\kappa}(\Ab))$ over all small regular cardinals, and hence fully faithful by \cref{derivedsolidforgetfunctorfullyfaithful,solkappafullyfaithfulonbderivedandpreservescolimits,filteredcolimitsofcategoriesappendix}. 

     To prove \ref{freeabeliangrouponprofiniteacyclicforsolidfication}, we show that for $M\in\mathcal D(\Sol{})$, pullback along the map $\mathbb Z[\underline{S}]^{L\solid}\to \mathbb Z[\underline{S}]^{\solid}$ adjoint to $\mathbb Z[\underline{S}]\to \mathbb Z[\underline{S}]^{\solid}$ gives an isomorphism   
     \[\Map_{\mathcal D(\Sol{})}(\mathbb Z[\underline{S}]^{\solid},M)\cong \Map_{\mathcal D(\Sol{})}(\mathbb Z[\underline{S}]^{L\solid},M).\]  
     Fix $M\in\mathcal D(\Sol{})$, $S\in\Pro(\Fin)$ and choose a regular cardinal $\kappa>\wt(S)$ with $M\in\mathcal D(\Sol{\kappa})$. 
     Then $fM\in\mathcal D(\Cond{\kappa}(\Ab))\subseteq \mathcal D(\Ab)$ and \[\Map_{\mathcal D(\Cond{}(\Ab))}(\mathbb Z[\underline{S}],fM)\cong \Map_{\mathcal D(\Cond{\kappa}(\Ab))}(\mathbb Z[\underline{S}_{\kappa}],fM).\] 
     By \cref{derivedsolidissolidderived}, \cref{solidificationoffreemodules} and fully faithfulness of $f$,
     \begin{align*}\Map_{\mathcal D(\Cond{\kappa}(\Ab))}(\mathbb Z[\underline{S}_{\kappa}],fM)&\cong \Map_{\mathcal D(\Cond{\kappa}(\Ab))}(\mathbb Z[\underline{S}_{\kappa}]^{\solid},fM)\\ &\cong \Map_{\mathcal D(\Cond{}(\Ab))}(\mathbb Z[\underline{S}]^{\solid},fM)\\ &\cong \Map_{\mathcal D(\Sol{})}(\mathbb Z[\underline{S}]^{\solid},M).\end{align*}
     This shows that the map $\mathbb Z[\underline{S}]^{L\solid}\to \mathbb Z[\underline{S}]^{\solid}$ adjoint to $\mathbb Z[\underline{S}]\to\mathbb Z[\underline{S}]^{\solid}$ is an equivalence for all $S\in\Pro(\Fin)$.
     By \cref{solidificationoffreemodules}, $\mathbb Z[\underline{S}]^{\solid}\cong \prod_{I}\mathbb Z$ for a small set $I$ of cardinality $|I|\leq \wt(S)$. In particular, $\mathbb Z[\underline{S}]^{\solid}\in\mathcal D(\Sol{\kappa})$ for all $\kappa>\wt(S)$. 

    As for an uncountable regular cardinal $\kappa$, $\mathcal D(\Cond{\kappa}(\Ab))\subseteq \mathcal D(\Cond{}(\Ab))$ is closed under small colimits and generated under small colimits by $\mathbb Z[\underline{S}], S\in\Pro(\Fin)_{\kappa}$ and $\mathcal D(\Sol{\kappa})\subseteq \mathcal D(\Sol{})$ is closed under small colimits, this shows that  
    \[\mathcal D(\Cond{\kappa}(\Ab))\subseteq \mathcal D(\Cond{}(\Ab))\xrightarrow{(-)^{L\solid}}\mathcal D(\Sol{})\] factors over \[\mathcal D(\Sol{\kappa})\subseteq \mathcal D(\Sol{}).\] 
    The induced functor \[\mathcal D(\Cond{\kappa}(\Ab))\to\mathcal D(\Sol{\kappa})\] is obviously left adjoint to the forget functor $\mathcal D(\Sol{\kappa})\to\mathcal D(\Cond{\kappa}(\Ab))$, and hence equivalent to $(-)^{\solid,\kappa}$. 

    For all uncountable regular cardinals $\kappa$, the composite \[f_{\kappa}\colon \mathcal D(\Sol{\kappa})\hookrightarrow \mathcal D(\Sol{})\xrightarrow{f}\mathcal D(\Cond{}(\Ab))\] is fully faithful by \cref{solkappafullyfaithfulonbderivedandpreservescolimits} and the above. By construction of the forget functor $f$, $f_{\kappa}$ factors over $\mathcal D(\Cond{\kappa}(\Ab))$. $f_{\kappa}$ therefore factors over a functor
    \[\tilde{f}_{\kappa}\colon \mathcal D(\Sol{\kappa})\to \mathcal D(\Sol{})\times_{\mathcal D(\Cond{}(\Ab))}\mathcal D(\Cond{\kappa}(\Ab)).\]
    As \[ \mathcal D(\Sol{})\to \mathcal D(\Cond{}(\Ab)), \, \mathcal D(\Cond{\kappa}(\Ab))\to \mathcal D(\Cond{}(\Ab))\] are fully faithful, $\tilde{f}_{\kappa}$ is fully faithful. 
    By \ref{solidificationfactors}, for $M\in\mathcal D(\Sol{})$ with $fM\in\mathcal D(\Cond{\kappa}(\Ab))$, $(fM)^{L\solid}\in\mathcal D(\Sol{\kappa})$. As $(fM)^{L\solid}\cong M$, this shows that $\tilde{f}_{\kappa}$ is essentially surjective. 
\end{proof}
 
\begin{cor}\label{derivedsolidissolidderivedwithoutkappa}
    For $M\in \mathcal D(\Cond{}(\Ab))$, the following are equivalent: 
\begin{romanenum}
    \item $M$ lies in the essential image of the forget functor $\mathcal D(\Sol{})\to\mathcal D(\Cond{}(\Ab))$. 
    \item \label{solidenrichedwithoutkappa} For all $S\in\Pro(\Fin)$, the unit $\mathbb Z[\underline{S}]\to\mathbb Z[\underline{S}]^{\solid}$ induces an equivalence \[\map_{\mathcal D(\Cond{}(\Ab))}(\mathbb Z[\underline{S}]^{\solid},M)\cong \map_{\mathcal D(\Cond{}(\Ab))}(\mathbb Z[\underline{S}],M).\] 
    \item \label{solidinternallywithoutkappa} For all $S\in\Pro(\Fin)$, the unit $\mathbb Z[\underline{S}]\to\mathbb Z[\underline{S}]^{\solid}$ induces an equivalence \[\imap_{\mathcal D(\Cond{}(\Ab))}(\mathbb Z[\underline{S}]^{\solid},M)\cong \imap_{\mathcal D(\Cond{}(\Ab))}(\mathbb Z[\underline{S}],M).\] 
    \item For all $i\in\mathbb Z$, $H_i(M)\in\Sol{}$. 
   \end{romanenum}
\end{cor}
\begin{proof}
    By construction (see the proof of \cref{derivedsolidificationwithoutkappa}), the forget functor is the colimit of the forget functors $\mathcal D(\Sol{\kappa})\to\mathcal D(\Cond{\kappa}(\Ab))$ over all small regular cardinals. 
    As for all regular cardinals $\kappa$, $\Sol{\kappa}\to \Sol{}$ is exact (\cref{colimitdescriprionsolid}), it follows from \cref{derivedsolidissolidderived} that the essential image of $\mathcal D(\Sol{})\subseteq \mathcal D(\Cond{}(\Ab))$ consists precisely of those $M$ such that $H_i(M)\in \Sol{}$ for all $i\in\mathbb Z$. 

    Suppose now that $M\in\mathcal D(\Cond{}(\Ab))$ satisfies \ref{solidenrichedwithoutkappa} and choose an uncountable regular cardinal $\kappa$ with $M\in\mathcal D(\Cond{\kappa}(\Ab))$. 
    \cref{solidificationoffreemodules,topspacecomefromregular,adjunctionsspectrallyenriched} and our assumption on $M$ imply that for all $S\in\Pro(\Fin)_{\kappa}$, the unit $\mathbb Z[\underline{S}_{\kappa}]\to\mathbb Z[\underline{S}_{\kappa}]^{\solid}$ induces an equivalence \[\map_{\mathcal D(\Cond{\kappa}(\Ab))}(\mathbb Z[\underline{S}_{\kappa}]^{\solid},M)\cong \map_{\mathcal D(\Cond{\kappa}(\Ab))}(\mathbb Z[\underline{S}],M).\] Whence by \cref{underivedsolidification,derivedsolidissolidderived}, $M$ lies in the essential image of \[\mathcal D(\Sol{\kappa})\to\mathcal D(\Cond{\kappa}(\Ab))\subseteq \mathcal D(\Cond{}(\Ab)).\] It now follows from \cref{freeabeliangroupsacyclicforsolidificationwithoutkappa} that $M$ lies in the essential image of the forget functor \[\mathcal D(\Sol{})\to\mathcal D(\Cond{}(\Ab)).\] 
    Conversely, suppose $M\in\mathcal D(\Sol{})$ and $S\in\Pro(\Fin)$. 
    For regular cardinals $\kappa>\wt(S)$, $\Cond{\kappa}(\Ab)\subseteq \Cond{}(\Ab)$ sends $\mathbb Z[\underline{S}_{\kappa}]$ to $\mathbb Z[\underline{S}]$ and $\mathbb Z[\underline{S}_{\kappa}]^{\solid}$ to $\mathbb Z[\underline{S}]^{\solid}$ by \cref{topspacecomefromregular,solidificationoffreemodules}. 
    Hence by \cref{forkappalargeenoughclosedmonoidal}, there exists a regular cardinal $\kappa>\wt(S)$ such that $fM\in\mathcal D(\Cond{\kappa}(\Ab))$ and \[ \imap_{\mathcal D(\Cond{\kappa}(\Ab))}(\mathbb Z[\underline{S}_{\kappa}],fM)\cong \imap_{\mathcal D(\Cond{}(\Ab))}(\mathbb Z[\underline{S}],fM),\]  \[ \imap_{\mathcal D(\Cond{\kappa}(\Ab))}(\mathbb Z[\underline{S}_{\kappa}]^{\solid},fM)\cong \imap_{\mathcal D(\Cond{}(\Ab))}(\mathbb Z[\underline{S}]^{\solid},fM),\] 
    It now follows from \cref{derivedsolidissolidderived} that $fM$ satisfies \ref{solidinternallywithoutkappa}, and in particular also \ref{solidenrichedwithoutkappa}.
\end{proof}
We now explain that $\mathcal D(\Sol{(\kappa)})$ inherits a closed symmetric monoidal structure from $\mathcal D(\Cond{(\kappa)}(\Ab))$. 
\begin{definition}
    Suppose that $\mathcal C^{\otimes}$ is a symmetric monoidal category. A localization \[L\colon \mathcal C\to \mathcal B\subseteq \mathcal C\] is symmetric monoidal if there exists a symmetric monoidal structure on $\mathcal B$ such that $L$ enhances to a symmetric monoidal functor $L^{\otimes}\colon \mathcal C^{\otimes}\to\mathcal B^{\otimes}$ with the following universal property: 
    If $\mathcal D^{\otimes}$ is a symmetric monoidal category, pullback along $L^{\otimes}$ defines a fully faithful functor 
\[ \Fun^{\otimes}(\mathcal B^{\otimes}, \mathcal D^{\otimes})\to \Fun^{\otimes}(\mathcal C^{\otimes}, \mathcal D^{\otimes})\] with essential image the symmetric monoidal functors $\mathcal C^{\otimes}\to\mathcal D^{\otimes}$ whose underlying functor $\mathcal C\to\mathcal D$ inverts $L$-equivalences. 
\end{definition}
The universal property determines the symmetric monoidal structure on $\mathcal B$ up to contractible choice. 
We now show that $(-)^{L\solid}$ enhances to a symmetric monoidal localisation. 
\begin{lemma}\label{derivedmappingspacesaresolid}
For $A\in\mathcal D(\Cond{(\kappa)}(\Ab))$ and $B\in\mathcal D(\Sol{(\kappa)})$, \[\imap_{\mathcal D(\Cond{(\kappa)}(\Ab))}(A,B)\in\mathcal D(\Sol{(\kappa)})\subseteq\mathcal D(\Cond{(\kappa)}(\Ab)).\] 
\end{lemma}
\begin{proof}For $A,B\in \mathcal D(\Cond{\kappa}(\Ab))$, the symmetric monoidal structure on $\mathcal D(\Cond{(\kappa)}(\Ab))$ defines an equivalence \[\imap_{\mathcal D(\Cond{\kappa}(\Ab))}(A\otimes B,-)\cong \imap_{\mathcal D(\Cond{\kappa}(\Ab))}(A, \imap_{\mathcal D(\Cond{\kappa}(\Ab))}(B,-)).\] 
Recall from \cref{derivedsolidissolidderived} and \cref{Gcategorysoliddefinition} that for $C\in\mathcal D(\Sol{\kappa})$ and $S\in\Pro(\Fin)_{\kappa}$, 
\[\imap_{\mathcal D(\Cond{\kappa}(\Ab))}(\mathcal M(S,\sphere{1}),C)\cong \imap_{\mathcal D(\Cond{\kappa}(\Ab))}(\mathbb Z[\underline{S}_{\kappa}],C)[-1].\] 
This implies that for $S\in\Pro(\Fin)_{\kappa}$ and $C\in\mathcal D(\Sol{\kappa})$, \begin{align*}&\imap_{\mathcal D(\Cond{\kappa}(\Ab))}(\mathcal M(S, \sphere{1}), \imap_{\mathcal D(\Cond{\kappa}(\Ab))}(B,C))\\&\cong \imap_{\mathcal D(\Cond{\kappa}(\Ab))}(B, \imap_{\mathcal D(\Cond{\kappa}(\Ab))}(\mathcal M(S, \sphere{1}),C))\\&\cong  \imap_{\mathcal D(\Cond{\kappa}(\Ab))}(B, \imap_{\mathcal D(\Cond{\kappa}(\Ab))}(\mathbb Z[\underline{S}_{\kappa}],C)[-1])\\ &\cong \imap_{\mathcal D(\Cond{\kappa}(\Ab))}(B, \imap_{\mathcal D(\Cond{\kappa}(\Ab))}(\mathbb Z[\underline{S}_{\kappa}],C))[-1]\\&\cong \imap_{\mathcal D(\Cond{\kappa}(\Ab))}(\mathbb Z[\underline{S}_{\kappa}], \imap_{\mathcal D(\Cond{\kappa}(\Ab))}(B,C))[-1],\end{align*}
which shows that
\[ \imap_{\mathcal D(\Cond{\kappa}(\Ab))}(B,C)\in \mathcal G\cong \mathcal D(\Sol{\kappa}).\qedhere\] 
  
By \cref{forkappalargeenoughclosedmonoidal}, for $A,B\in\mathcal D(\Cond{}(\Ab))$ there exists a regular cardinal $\kappa$ with \[A,B\in\mathcal D(\Cond{\kappa}(\Ab))\subseteq \mathcal D(\Cond{}(\Ab))\] and \[\imap_{\mathcal D(\Cond{\kappa}(\Ab))}(A,B)=\imap_{\mathcal D(\Cond{}(\Ab))}(A,B).\]
Since $\mathcal D(\Sol{})\times_{\mathcal D(\Cond{}(\Ab))}\mathcal D(\Cond{\kappa}(\Ab))\cong \mathcal D(\Sol{\kappa})$ (\cref{freeabeliangroupsacyclicforsolidificationwithoutkappa}), the condensed statement follows from the above. 
\end{proof}

\begin{cor}\label{symmetricmonoidalstructureonderivedsolidabeliangroups}
    The localization $(-)^{L\solid}\colon\mathcal D(\Cond{(\kappa)}(\Ab))\to\mathcal D(\Sol{(\kappa)})$ is symmetric monoidal. 
    The induced symmetric monoidal structure on $\mathcal D(\Sol{(\kappa)})$ is closed.
\end{cor}
\begin{proof}
We first show that the fiber of $(-)^{L\solid}$ is a tensor ideal, i.e.\ that for $A\in \mathcal D(\Cond{(\kappa)}(\Ab))$ with $A^{L\solid}=0$ and $B\in\mathcal D(\Cond{(\kappa)}(\Ab))$, $(A\otimes B)^{L\solid}=0$. 
Suppose that $A^{L\solid}=0$, and fix $B\in\mathcal D(\Cond{\kappa}(\Ab))$ and $C\in\mathcal D(\Sol{\kappa})$. 
By \cref{derivedmappingspacesaresolid}, $\imap_{\mathcal D(\Cond{(\kappa)}(\Ab))}(B,C)\in\mathcal D(\Sol{(\kappa)})$, and hence \begin{align*}\Map_{\mathcal D(\Sol{(\kappa)})}((A\otimes B)^{L\solid},C)& \cong \Map_{\mathcal D(\Cond{(\kappa)}(\Ab))}(A\otimes B, C)\\& \cong \Map_{\mathcal D(\Cond{(\kappa)}(\Ab))}(A, \imap_{\mathcal D(\Cond{(\kappa)}(\Ab))}(B,C))\\ &\cong \Map_{\mathcal D(\Cond{(\kappa)}(\Ab))}(A^{L\solid}, \imap_{\mathcal D(\Cond{(\kappa)}(\Ab))}(B,C))\\ &=0, \end{align*} which shows that $(A\otimes B)^{L\solid}=0$. 
Since $\mathcal D(\Cond{(\kappa)}(\Ab))$ is stable and the symmetric monoidal structure on $\mathcal D(\Cond{(\kappa)}(\Ab))$ is cocontinuous, it now follows from \cite[Proposition 4.1.7.4]{higheralgebra} that $(-)^{L\solid}$ is a symmetric monoidal localization. 

Since $\mathcal D(\Sol{(\kappa)})\subseteq \mathcal D(\Cond{(\kappa)}(\Ab))$ is a full subcategory and for $A\in\mathcal D(\Sol{(\kappa)})$, the internal Hom $\imap_{\mathcal D(\Cond{(\kappa)}(\Ab))}(A,-)$ restricts to a functor $\mathcal D(\Sol{(\kappa)})\to\mathcal D(\Sol{(\kappa)})$ (\cref{derivedmappingspacesaresolid}), this is a right adjoint to $A\otimes_{\mathcal D(\Sol{(\kappa)})}-$.  
\end{proof} 

\begin{lemma}\label{tstructuresolidcompatiblewithsymmetricmonoidalstructurewithoutkappa}
The symmetric monoidal structure on $\mathcal D(\Sol{(\kappa)})$ is compatible with the $t$-structure. 
\end{lemma}
\begin{proof}The unit $\mathbb Z\in\mathcal D(\Sol{\kappa})$ is connective.  
    By construction of the symmetric monoidal structure, for $M,N\in\mathcal D(\Sol{(\kappa)})$, \[M\otimes_{\mathcal D(\Sol{\kappa})} N\cong (fM\otimes_{\mathcal D(\Cond{(\kappa)}(\Ab))} fN)^{L\solid}.\] 
    Since $(-)^{L\solid}$ is left adjoint to the $t$-exact functor $f$, it is right $t$-exact. 
    Since $f$ is $t$-exact and the $t$-structure on $\mathcal D(\Cond{(\kappa)}(\Ab))$ is compatible with the symmetric monoidal structure, this implies that for $M,N\in\mathcal D(\Sol{(\kappa)})_{\geq 0}$, $M\otimes_{\mathcal D(\Sol{(\kappa)}(\Ab))} N\in\mathcal D(\Sol{(\kappa)}(\Ab))_{\geq 0}$. 
\end{proof}
\begin{rem}\label{monoidalstructureonsolidabeliangroupsclosed}
 In particular, $\Sol{(\kappa)}$ inherits a symmetric monoidal structure such that \[H_0\colon \mathcal D(\Sol{(\kappa)})_{\geq 0}\to \Sol{(\kappa)}\] enhances to a symmetric monoidal functor, see \cite[Example 2.2.1.10]{higheralgebra}. 
This symmetric monoidal structure is closed with internal Hom 
\[ \pi_0\imap_{\mathcal D(\Sol{(\kappa)})}(-,-).\] \end{rem}
\begin{rem}\label{monoidalstructurederivedtensorprduct}
    As for $T,S\in\Pro(\Fin)_{\kappa}$, \begin{align}\label{acyclicityprofinitesolidtensorproduct}\mathbb Z[\underline{S}_{(\kappa)}]^{\solid}\otimes_{\mathcal D(\Sol{(\kappa)})}\mathbb Z[\underline{T}_{(\kappa)}]^{\solid}\cong \mathbb Z[\underline{T\times S}_{(\kappa)}]^{\solid}\in\mathcal D(\Sol{(\kappa)})^{\heart},\end{align} the  tensor product $-\otimes_{\mathcal D(\Sol{\kappa})}-$ is the derived functor of $-\otimes_{\Sol{\kappa}}-$. 
In more detail: Since the $t$-structure on $\mathcal D(\Sol{(\kappa)})$ is right-complete and $-\otimes_{\mathcal D(\Sol{(\kappa)})}-$ is cocontinuous in both variables, $-\otimes_{\mathcal D(\Sol{(\kappa)})}-$ is determined by its valued on bounded-below complexes.
\cref{acyclicityprofinitesolidtensorproduct} and exactness of filtered colimits in $\Sol{(\kappa)}(\Ab)$ imply that $-\otimes_{\mathcal D(\Sol{(\kappa)})}-$ can be computed on bounded complexes by resolving both variables with bounded below complexes $C_*$ with \[C_n=\oplus_{i\in I_n}\mathbb Z[\underline{S^i}_{(\kappa)}]^{\solid},S^i\in\Pro(\Fin)_{(\kappa)}\] for all $n\in\mathbb Z$. 

By \cref{Nobelingspecker} and \cite[Corollary 3.4.3]{Juansolidgeometry}, for $S\in\Pro(\Fin)_{\light}$, $\mathbb Z[\underline{S}_{\light}]^{\solid}$ is flat. In particular, the  category of light solid abelian groups has enough flat objects, which implies that for $M\in\Sol{\light}$, $M\otimes_{\mathcal D(\Sol{\light})}-$ is the derived functor of $M\otimes_{\Sol{\light}}-$. 

An example of Efimov shows that that $\Sol{\kappa}$ in general does not have enough flat objects, so in particular, $M\otimes_{\mathcal D(\Sol{\kappa})}-$ is not the derived functor of $M\otimes_{\Sol{\kappa}}-$ for general $\kappa$ and $M\in\Sol{\kappa}$. 
More concretely, $M\otimes_{\mathcal D(\Sol{\kappa})}-$ is the derived functor of $M\otimes_{\Sol{\kappa}}-$ if and only if for all $S\in\Pro(\Fin)_{\kappa}$, $M\otimes_{\mathcal D(\Sol{\kappa})}\mathbb Z[\underline{S}_{\kappa}]^{\solid\mathbb Z}\in\mathcal D(\Sol{\kappa})^{\heart}$.  
This holds for example if the condensed abelian group underlying $M$ is pseudo-coherent, see \cref{examplessflatrings} below. 
\end{rem}

\begin{lemma}\label{solidmodulesismodulesinsolid}
Suppose $R$ is a ($\kappa$-)condensed ring. 
\begin{romanenum} 
\item The forget functor \[\Sol{(\kappa)}(\algebra{R})\to\Sol{(\kappa)}\] factors over an equivalence $\Sol{(\kappa)}(\algebra{R})\cong \LMod{\algebra{R}^{\solid}}{\Sol{(\kappa)}}$. 

\item If $\algebra{R}$ is commutative, then $(-)^{\solid\algebra{R}}$ is a symmetric monoidal localisation. In particular, $\Sol{(\kappa)}(\algebra{R})$ inherits a cocontinuous symmetric monoidal structure such that $(-)^{\solid\algebra{R}}$ enhances to a symmetric monoidal functor.
The enhancement $(-)^{\solid\algebra{R}}\colon \Cond{(\kappa)}(\algebra{R})\to \Sol{(\kappa)}(\algebra{R})\in\CAlg(\vlCat)$ of $(-)^{\solid\algebra{R}}$ is essentially unique. 
\end{romanenum}
\end{lemma}
\begin{proof}
The symmetric monoidal localisation $\Sol{(\kappa)}\subseteq \Cond{(\kappa)}(\Ab)$ induces a localisation 
\[\LMod{\algebra{R}^{\solid}}{\Sol{(\kappa)}}\subseteq \LMod{\algebra{R}}{\Cond{(\kappa)}(\Ab)}\] with essential image 
$\LMod{\algebra{R}}{\Cond{(\kappa)}(\Ab)}\times_{\Cond{(\kappa)}(\Ab)}\Sol{(\kappa)}$, cf.\ \cref{localisationinduceslocalisationonmodulecategories}. The essential image are precisely the solid $\algebra{R}$-modules, which shows that $\Sol{(\kappa)}(\algebra{R})\cong \LMod{\algebra{R}^{\solid}}{\Sol{(\kappa)}}$. 

As $\Sol{(\kappa)}(\Ab)$ has small colimits (\cref{solidclosedunderlimitscolimitskappa,underivedsolidificationwithoutkappa}) and the tensor product is cocontinuous in both variables (the symmetric monoidal structure is closed), for a commutative condensed ring $\algebra{R}$, \cref{symmetricmonoidalstructure,symmetricmonoidalstructuremodulesnatural} provide a symmetric monoidal structure on $\LMod{\algebra{R}^{\solid}}{\Sol{(\kappa)}}$ and a symmetric monoidal enhancement of $(-)^{\solid\algebra{R}}$.  
Suppose that $w\colon M\to N\in\Cond{(\kappa)}(\algebra{R})$ is such that $w^{\solid\algebra{R}}$ is an equivalence. 
As $(-)^{\solid\algebra{R}}$ enhances to a symmetric monoidal functor, $N\in\Cond{(\kappa)}(\algebra{R})$, $(w\otimes_{\Cond{(\kappa)}(\algebra{R})}N)^{\solid\algebra{R}}$ is an equivalence. It now follows from \cite[Proposition 4.1.7.4]{higheralgebra} (in the large universe) that the localisation $(-)^{\solid\algebra{R}}$ is symmetric monoidal and that the symmetric monoidal enhancement of $(-)^{\solid\algebra{R}}$ is essentially unique. 
\end{proof}

\begin{lemma}\label{solidenoughprojectiveswithoutkappa}
    For a profinite set $S$, $\algebra{R}[\underline{S}]^{\solid\algebra{R}}$ is projective in $\Sol{}(\algebra{R})$. 
    Every solid $\algebra{R}$-module is a quotient of $\oplus_{i\in I}\algebra{R}[\underline{S^i}]^{\solid\algebra{R}}$ for some small family of profinite sets $(S^i)_{i\in I}\in\Pro(\Fin)$.
    In particular, $\Sol{}(\algebra{R})$ has enough projectives. 
\end{lemma}
\begin{proof}
    By \cref{enoughprojectives}, for every condensed $\algebra{R}$-module $M$, there exists a small collection $S_i,i\in I$ of extremally disconnected compact Hausdorff spaces with a quotient map $\oplus_{i\in I}\algebra{R}[\underline{S_i}]\to M$. 
If $M$ is solid, this factors over a quotient map $\oplus_{i\in I}\algebra{R}[\underline{S_i}]^{\solid\algebra{R}}\to M$. 
It remains to show that for a profinite set $X$, $R[\underline{X}]^{\solid\algebra{R}}$ is projective. 
    
    As its right adjoint is exact, the free $R^{\solid}$-module functor $\Sol{}\to \Sol{}(\algebra{R})$ preserves projectives. 
    It therefore suffices to show that for $S\in\Pro(\Fin)$, $\mathbb Z[\underline{S}]^{\solid\mathbb Z}$ is projective in $\Sol{}$. 
Fix $S\in\Pro(\Fin)$ and choose a regular cardinal $\kappa>\wt(S)$. Then $\mathbb Z[\underline{S}]^{\solid\mathbb Z}\cong i_{\kappa}(\mathbb Z[\underline{S}_{\lambda}])^{\solid\mathbb Z,\kappa}\cong \prod_J\mathbb Z$ for a small set of cardinality $|J|\leq \wt(X)$ by \cref{Nobelingspecker,solidificationoffreemodules}. We now show that for a small set $I$, $\prod_{I}\mathbb Z$ is projective in $\Sol{}$. 
\cref{Nobelingspecker} and the universal property of Stone-\v{C}ech-compactification imply that if $S\coloneqq \beta(I)$ is the Stone-\v{C}ech compactification of a set $I$, then \[\mathcal C(S,\mathbb Z)\cong \colim{n\in\mathbb N_0}\mathcal C(S,[-n,n]\cap\mathbb Z)\cong \colim{n\in\mathbb N_0}\{-n,n\}^{I}.\] 
There is a monomorphism $\oplus_{I}\mathbb Z\hookrightarrow  \colim{n\in\mathbb N_0}\{-n,n\}^{I} $, whence by \cref{Nobelingspecker}, $\mathbb Z[\underline{S}]^{\solid}\cong \prod_{K}\mathbb Z$ for a small set $K$ with $|K|\geq |I|$. 
This implies that $\prod_{I}\mathbb Z$ is a direct summand of $\mathbb Z[\underline{S}]^{\solid}\cong \prod_{J}\mathbb Z$.
As $(-)^{\solid\mathbb Z}$ preserves projectives (its right adjoint is exact), $\mathbb Z[\underline{S}]^{\solid}$ is projective, which shows that $\prod_{I}\mathbb Z$ is projective. 
\end{proof}
\cref{solidenoughprojectiveswithoutkappa,solidenoughprojectiveskappa} imply the following, which will be essential for our identification of continuous and solid group cohomology (\cref{solidequalscontinuousgoodgroups}). 
\begin{cor}\label{solidprojectivestensorproduct}
    Suppose $\algebra{R}\in \CAlg(\Cond{(\kappa)}(\Ab))$. If $P,Q$ are projective solid $\algebra{R}$-modules, then $P\otimes_{\Sol{(\kappa)}(\algebra{R})}Q$ is projective in $\Sol{(\kappa)}(\algebra{R})$. 
\end{cor}

\begin{rem}
The analogue of \cref{solidprojectivestensorproduct} is very much not satisfied in $\Cond{(\kappa)}(\Ab)$:
For a strong limit cardinal $\kappa$, and a ($\kappa$-light) profinite space $X$, $\mathbb Z[\underline{X}_{(\kappa)}]\otimes_{\Cond{(\kappa)}(\Ab)}\mathbb Z[\underline{X}_{(\kappa)}]\cong\mathbb Z[\underline{X^2}_{(\kappa)}]$ is projective in $\Sol{(\kappa)}$ if and only if $X$ is discrete by \cite[Proposition 4.8]{Scholzeanalytic}.  
\end{rem}
\begin{proof}
    Suppose $P,Q$ are projective in $\Sol{(\kappa)}(\algebra{R})$. \cref{solidenoughprojectiveswithoutkappa,solidenoughprojectiveskappa} imply that there exist small families of ($\kappa$-light) profinite sets $(S^i)_{i\in I}, (T^j)_{j\in J}$ so that $P$ is a direct summand of $\oplus_{i\in I}\algebra{R}[\underline{S^i}_{(\kappa)}]^{\solid\algebra{R}}$, and $Q$ is a direct summand of $\oplus_{j\in J}\algebra{R}[\underline{T^j}_{(\kappa)}]^{\solid\algebra{R}}$. 
    Since $-\otimes_{\Sol{(\kappa)}(\algebra{R})}-$ is cocontinuous in both variables, this implies that $P\otimes_{\Sol{(\kappa)}(\algebra{R})}Q$ is a direct summand of \[\oplus_{j\in J}\oplus_{i\in I}\algebra{R}[\underline{S^i}_{(\kappa)}]^{\solid\algebra{R}}\otimes_{\Sol{(\kappa)}(\algebra{R})}\algebra{R}[\underline{T^j}_{(\kappa)}]^{\solid\algebra{R}}\cong \oplus_{i\in I}\oplus_{j\in J}\algebra{R}[\underline{S^i\times T^j}_{(\kappa)}]^{\solid\algebra{R}}.\] 
    As for all $i,j$, $S^i\times T^j\in\Pro(\Fin)_{(\kappa)}$, 
    $\oplus_{i\in I}\oplus_{j\in J}\algebra{R}[\underline{S^i\times T^j}_{(\kappa)}]^{\solid\algebra{R}}$ is projective in $\Sol{(\kappa)}(\algebra{R})$ by \cref{solidenoughprojectiveskappa,solidenoughprojectiveswithoutkappa}, respectively. 
\end{proof}

\subsubsection{Examples and computations}\label{section:computationssolidification}
In this section, we give some examples of solid modules and compute the derived and underived solidification of free abelian groups on compact Hausdorff spaces and CW-complexes. We will use these results to obtain further identifications of condensed with sheaf cohomology in to compute the condensed cohomology (\cref{section:condensedcohomologywithsolidcoefficients}), and to compare continuous with solid group cohomology (\cref{section:continuousandsolidgroupcohomology}). 
\begin{lemma}\label{nullsequencessummable}
    Let $\mathbb N\cup\{\infty\}$ denote one-point compactification of the set of natural numbers, let \[P\coloneqq \mathbb Z[\underline{\mathbb N\cup\{\infty\}}_{(\kappa)}]/\{\infty\}\coloneqq \Coker(\mathbb Z[\{\infty\}]\to\mathbb Z[\mathbb N\cup \{\infty\}]),\] and denote 
    by \[s\colon P\to P\in\Cond{\kappa}(\Ab)\] the homomorphism induced by \[\mathbb N_0\cup\{\infty\}\ni n\mapsto \begin{cases}
        [n]-[n+1] & n\in\mathbb N_0\\
        [\infty] & n=\infty
    \end{cases}.\] 
    The induced map $s^{\solid}\colon P^{\solid}\to P^{\solid}$ is an isomorphism. 
\end{lemma}
\begin{proof}
    As $\mathcal C(\mathbb N\cup\{\infty\}, \mathbb Z)\to \mathbb Z\oplus \bigoplus_{n\in\mathbb N_0}\mathbb Z, \, f\mapsto (f(\infty), (f(n)-f(\infty))_{n})$ is an isomorphism, 
     \[P^{\solid}=\Coker(\mathbb Z\to \mathbb Z[\mathbb N\cup \{\infty\}]^{\solid})\cong \iHom_{\Cond{(\kappa)}(\Ab)}(\oplus_{n\in\mathbb N_0}\mathbb Z, \mathbb Z).\] 
    Under this identification, $s^{\solid}$ is pullback along the isomorphism \[\oplus_{n\in\mathbb N_0}\mathbb Z\to \oplus_{n\in\mathbb N_0}\mathbb Z, (a_n)_n\mapsto (a_n-a_{n+1})_n, \] and hence an isomorphism. 
\end{proof}
\begin{rems}
\begin{romanenum}
    \item If $M$ is a $(\To)$ topological abelian group, then \[\Hom_{\Cond{(\kappa)}(\Ab)}(P, \underline{M}_{(\kappa)})=\text{Null}(M)\] is the abelian group of nullsequences in $M$, and the lemma implies that if $\underline{M}_{(\kappa)}$ is solid, then 
    $\text{Null}(M)\to \text{Null}(M), \, (a_n)_n\mapsto (a_n-a_{n+1})_n$ is an isomorphism, i.e.\ every nullsequence in $M$ is summable.
    \item In light condensed abelian groups, summability of all nullsequences characterizes solid abelian groups, i.e.\ a light condensed abelian group $M$ is solid if and only if 
\[ s^*\colon \Hom_{\Cond{\light}(\Ab)}(P,M)\to \Hom_{\Cond{\light}(\Ab)}(P,M)\] is an isomorphism, see e.g.\ \cite[Lecture 5]{Analyticstacks} or \cite[Proposition 3.2.3]{Juansolidgeometry}.
    \end{romanenum}    
\end{rems}
\begin{lemma}[{\cite{384945}}]\label{solidificationoftopologicalspaces}
    Denote by \[ \cshv{-}\colon \mathcal D(\Ab)\cong \LMod{H\mathbb Z}{\Sp}\to  \LMod{cH\mathbb Z}{\Cond{(\kappa)}(\Sp)}\cong \mathcal D(\Cond{(\kappa)}(\Ab))\] the functor induced by the constant sheaf functor. 
    \begin{romanenum}
    \item For $M\in\mathcal D(\Cond{(\kappa)}(\Ab))$, there is a natural map \[M^{L\solid}\to \imap_{\mathcal D(\Cond{(\kappa)}(\Ab))}(\cshv{\map_{\mathcal D(\Cond{(\kappa)}(\Ab))}(M, \mathbb Z)}, \mathbb Z).\] 
    \item This is an equivalence for $M=\mathbb Z[\underline{X}_{(\kappa)}], X\in\CH_{(\kappa)}$. 
    \end{romanenum}
\end{lemma}  
\begin{proof}
    We first construct the map \[M^{L\solid}\to \imap_{\mathcal D(\Cond{(\kappa)}(\Ab))}(\cshv{\map_{\mathcal D(\Cond{(\kappa)}(\Ab))}(M, \mathbb Z)}, \mathbb Z).\] 
    As $\mathbb Z$ is solid, for all $M\in\mathcal D(\Cond{(\kappa)}(\Ab))$, \[\imap_{\mathcal D(\Cond{(\kappa)}(\Ab))}(\cshv{\map_{\mathcal D(\Cond{(\kappa)}(\Ab))}(M, \mathbb Z)}, \mathbb Z)\in\mathcal D(\Sol{(\kappa)}(\Ab))\] by \cref{derivedmappingspacesaresolid}, so it suffices to construct a map \[M\to  \imap_{\mathcal D(\Cond{(\kappa)}(\Ab))}(\cshv{\map_{\mathcal D(\Cond{(\kappa)}(\Ab))}(M, \mathbb Z)}, \mathbb Z), \] then this factors uniquely over a map 
    \[\eta_{M}\colon M^{L\solid}\to  \imap_{\mathcal D(\Cond{(\kappa)}(\Ab))}(\cshv{\map_{\mathcal D(\Cond{(\kappa)}(\Ab))}(M, \mathbb Z)}, \mathbb Z).\] 

    The counit of the constant sheaf/global sections adjunction induces a natural transformation
    \[\epsilon\colon \cshv{-}\circ \map_{\mathcal D(\Cond{(\kappa)}(\Ab))}(-, \mathbb Z)\to \imap_{\mathcal D(\Cond{(\kappa)}(\Ab))}(-, \mathbb Z), \] and the counit for  $M\otimes -\dashv \imap_{\mathcal D(\Cond{(\kappa)}(\Ab))}(M,-)$ yields a map \[E\colon M\otimes \imap_{\mathcal D(\Cond{(\kappa)}(\Ab))}(M, \mathbb Z)\to\mathbb Z.\] 
    By definition of $\imap_{\mathcal D(\Cond{(\kappa)}(\Ab))}(-,-)$, the composite \[E\circ (\id\otimes \epsilon_*)\colon M\otimes_{\mathcal D(\Cond{(\kappa)}(\Ab))}\cshv{\map_{\mathcal D(\Cond{(\kappa)}(\Ab))}(\mathbb Z[X], \mathbb Z)}\to  \mathbb Z\] determines a map  
    \[ M\to \imap_{\mathcal D(\Cond{(\kappa)}(\Ab))}(\cshv{\map_{\mathcal D(\Cond{(\kappa)}(\Ab))}(M, \mathbb Z)}, \mathbb Z).\] 
    Denote by 
    \[\eta_{M} \colon M^{L\solid}\to \imap_{\mathcal D(\Cond{(\kappa)}(\Ab))}(\cshv{\map_{\mathcal D(\Cond{(\kappa)}(\Ab))}(M, \mathbb Z)}, \mathbb Z)\] the induced map. This is by construction natural in the derived $(\kappa)$-condensed abelian group $M$. 

    We now want to show that for $Y\in\Pro(\Fin)_{(\kappa)}$, 
    $\eta_{\mathbb Z[\underline{Y}_{(\kappa)}]}$ is an isomorphism. 
    By \cref{profiniteprojectiveinsolidkappa}/ \cref{freeabeliangroupsacyclicforsolidificationwithoutkappa}, $\mathbb Z[\underline{Y}_{(\kappa)}]^{L\solid}\cong \mathbb Z[\underline{Y}_{(\kappa)}]^{\solid}$ is concentrated in degree $0$. 
    We first show that the right-hand side is also concentrated in degree $0$.  
    By \cite[Theorem 3.2]{Scholzecondensed}/\cref{cohomologyclausenscholzeprofinite}, for $Y\in\Pro(\Fin)_{(\kappa)}$, \[\map_{\mathcal D(\Cond{(\kappa)(\Ab)})}(\mathbb Z[\underline{Y}_{(\kappa)}], \mathbb Z)\cong \mathcal C(Y, \mathbb Z)\] is concentrated in degree $0$. 
    As the constant sheaf functor is $t$-exact, this implies that 
    \[\cshv{\map_{\mathcal D(\Cond{(\kappa)})}(\mathbb Z[\underline{Y}_{(\kappa)}], \mathbb Z)}\cong \underline{\mathcal C(Y, \mathbb Z)}^{\delta}_{(\kappa)}\] is concentrated in degree $0$ and represented by the discrete abelian group $\mathcal C(Y,\mathbb Z)^{\delta}$.
    By \cref{Nobelingspecker}, $\mathcal C(Y, \mathbb Z)\cong \oplus_{I}\mathbb Z$ for some ($\kappa$-small) set $I$. 
    Since coproducts and ($\kappa$-small) products in $\Cond{(\kappa)}(\Ab)$ are exact and the constant sheaf functor preserves coproducts, it follows that \[\imap_{\mathcal D(\Cond{(\kappa)}(\Ab))}(\cshv{\map_{\mathcal D(\Cond{(\kappa)})}(\mathbb Z[\underline{Y}_{(\kappa)}], \mathbb Z)}, \mathbb Z)\cong \prod_{I}\mathbb Z\] is concentrated in degree $0$. 
    We are therefore reduced to showing that $H_0(\eta_{\mathbb Z[\underline{Y}_{(\kappa)}]})$ is an isomorphism. 
    By construction, 
    \[c\colon \mathbb Z[\underline{Y}_{(\kappa)}]\to \mathbb Z[\underline{Y}_{(\kappa)}]^{\solid}\xrightarrow{H_0(\eta_{\mathbb Z[\underline{Y}_{(\kappa)}]})} \iHom_{\Cond{(\kappa)}(\Ab)}(\underline{\mathcal C(Y, \mathbb Z)}^{\delta}_{(\kappa)}, \mathbb Z)\] is the map adjoint to 
    \[\mathbb Z[\underline{Y}_{(\kappa)}]\otimes\underline{\mathcal C(Y, \mathbb Z)}^{\delta}_{(\kappa)}\to \mathbb Z[\underline{Y}_{(\kappa)}]\otimes  \iHom_{\Cond{(\kappa)}(\Ab)}(\mathbb Z[\underline{Y}_{(\kappa)}], \mathbb Z)\to\mathbb Z, \] where the first map is induced by the counit of $c(-)\dashv \Gamma$, and the second map by the adjunction counit for $\mathbb Z[\underline{Y}_{(\kappa)}]\otimes-\vdash \iHom_{\Cond{(\kappa)}(\Ab)}(\mathbb Z[\underline{Y}_{(\kappa)}],-)$ at $\mathbb Z$.  
    By \cref{internalcohomologydiscretecompacthausdorff}, the first map is an isomorphism. The composite $\mathbb Z[\underline{Y}_{(\kappa)}]\xrightarrow{c}  \iHom_{\Cond{(\kappa)}(\Ab)}(\underline{\mathcal C(Y, \mathbb Z)}^{\delta}_{(\kappa)}, \mathbb Z)\cong \mathbb Z[\underline{Y}_{(\kappa)}]^{\blacksquare}$ is the map from the definition of solidity, whence $H_0(\eta_{\mathbb Z[\underline{Y}_{(\kappa)}]})$ is an isomorphism.
    
    Denote by $\mathcal C\subseteq\mathcal D(\Cond{(\kappa)}(\Ab))$ the full subcategory on objects $M$ for which $\eta_{M}$ is an equivalence. 
    As $(-)^{L\solid}$ and \[R\coloneqq  \imap_{\mathcal D(\Cond{(\kappa)}(\Ab))}(-, \mathbb Z)\circ \cshv{-}\circ\map_{\mathcal D(\Cond{(\kappa)}(\Ab))}(-, \mathbb Z)\] are exact, $\mathcal C$ is a stable subcategory of $\mathcal D(\Cond{(\kappa)}(\Ab))$.  
    Suppose now that $M\in\Cond{(\kappa)}(\Ab)$ and there exists a resolution $C_*\to M$ with $C_n=\mathbb Z[\underline{Y^n}_{(\kappa)}], \, Y^n\in\Pro(\Fin)_{(\kappa)}$ for all $n\in\mathbb N_0$. 
    We claim that this implies that $M\in\mathcal C$. 
    For $n\in\mathbb N_0$ denote by $c_n\colon C_*^{\leq n}\to C$ the map to from its stupid truncation.
    As $(-)^{L\solid}$ and $R$ are exact and right t-exact functors and \[\Cofib(C_*^{\leq n}\to C_*)\in \mathcal D(\Cond{(\kappa)}(\Ab))_{>n},\] \[H_k(c_n^{L\solid})\colon H_k\left((C_*^{\leq n})^{L\solid}\right)\to H_k(C_*^{L\solid})\text{  and } H_k(R(c_n))\colon H_k(R(C_*^{\leq n}))\to H_k(R(C_*))\] is an isomorphism for $k< n$. 
    It is therefore enough to show that $C_*^{\leq n}\in\mathcal C$ for all  $n\in\mathbb N_0$, then it follows that $M\cong C_*\in\mathcal C$. 
    We prove this by induction on $n$. 
    The case $n=0$ was shown above. The inductive statement follows since $C_*^{\leq n}\to C_*^{\leq n+1}\to C_{n+1}[n+1]$ is a cofiber sequence, $C_*^{\leq n}, C_{n+1}\in\mathcal C$, and $\mathcal C$ is a stable subcategory.

    Suppose now that $X$ is a compact Hausdorff space. 
    By \cref{compacthausdorffquotientoftdch}, there exists a quotient map $\pi\colon T\to X$ from a $\kappa$-light profinite space $T\in\Pro(\Fin)_{(\kappa)}$. The simplicial resolution $S_*^{\pi}$ of $\pi$ (\cref{definitionsimplicialcomplex}) is a resolution of $\mathbb Z[\underline{X}_{(\kappa)}]$ as above, which shows that $\mathbb Z[\underline{X}_{(\kappa)}]\in\mathcal C$.     
\end{proof}
\begin{cor}\label{solidificationcompacthausdorffspace}
    Suppose $X$ is a ($\kappa$-light) compact Hausdorff space and endow $\mathcal C(X, \mathbb Z)$ with the compact open topology which is discrete.
    \begin{romanenum}\item Then  \[\mathbb Z[\underline{X}_{(\kappa)}]^{\solid}\cong \iHom_{\Cond{(\kappa)}(\Ab)}(\underline{\mathcal C(X, \mathbb Z)}_{(\kappa)}, \mathbb Z)\] naturally in $X$, where $\iHom_{\Cond{(\kappa)}(\Ab)}(-,-)$ denotes the internal Hom of $\Cond{(\kappa)}(\Ab)$.
    \item Denote by $\pi_0X$ the set of connected components of $X$ and endow it with the quotient topology of the projection $X\to\pi_0X$. 
    The projection induces an isomorphism $\mathbb Z[\underline{X}_{(\kappa)}]^{\solid}\cong \mathbb Z[\underline{\pi_0 X}_{(\kappa)}]^{\solid}$. 
    \item Endow the set $\pi_0^{\gamma}X$ of path-connected components of $X$ with the quotient topology from the projection $X\to \pi_0^{\gamma}X$. If $\pi_0^{\gamma}X$ is Hausdorff, the projection induces an isomorphism $\mathbb Z[\underline{X}_{(\kappa)}]^{\solid}\cong \mathbb Z[\underline{\pi_0^{\gamma} X}_{(\kappa)}]^{\solid}$. 
    \end{romanenum}
\end{cor}
\begin{proof}
    We first recall that for $A\in\mathcal D(\Cond{(\kappa)}(\Ab))_{\leq 0}$ and $B\in\Cond{\kappa}(\Ab)$, the counit $\pi_0A\to A$ induces an equivalence \[\pi_0\imap_{\mathcal D(\Cond{(\kappa)}(\Ab))}(A,B)\cong \pi_0\imap_{\mathcal D(\Cond{(\kappa)}(\Ab))}(\pi_0A,B).\]
Fix $A\in\mathcal D(\Cond{(\kappa)}(\Ab))_{\leq 0}, B\in\Cond{\kappa}(\Ab)$. 
The cofiber sequences \[\tau_{\geq n}A\to \tau_{\geq n-1}A\to H_{n-1}(A)[n-1]\] imply that for $n\in\mathbb Z$, \[\tau_{\geq n}\imap_{\mathcal D(\Cond{\kappa}(\Ab))}(\tau_{\geq n-1}A,B)\cong \tau_{\geq n}\imap_{\mathcal D(\Cond{(\kappa)}(\Ab))}(\tau_{\geq n}A,B).\] 
As the $t$-structure on $\mathcal D(\Cond{\kappa}(\Ab))$ is right-complete (\cref{tstructureleftcomplete}), ${A}=\colim{n\to -\infty}\tau_{\geq n}{A}$, whence \[\tau_{\geq n}\imap_{\mathcal D(\Cond{\kappa}(\Ab))}({A},B)\cong \tau_{\geq n}\imap_{\mathcal D(\Cond{\kappa}(\Ab))}(\tau_{\geq n}A,B)\] by the above. 

Since $\cshv{-}$ is $t$-exact, \cref{solidificationoftopologicalspaces} now implies that 
\begin{align*}
\mathbb Z[\underline{X}_{(\kappa)}]^{\solid}&\cong \pi_0\imap_{\mathcal D(\Cond{(\kappa)}(\Ab))}(\cshv{\map_{\mathcal D(\Cond{(\kappa)}(\Ab))}(\mathbb Z[\underline{X}_{(\kappa)}], \mathbb Z)}, \mathbb Z)\\ 
&\cong \iHom_{\Cond{(\kappa)}(\Ab)}(\pi_0\cshv{\map_{\mathcal D(\Cond{(\kappa)}(\Ab))}(\mathbb Z[\underline{X}_{(\kappa)}], \mathbb Z)}, \mathbb Z)\\ 
& \cong \iHom_{\Cond{(\kappa)}(\Ab)}(\cshv{\pi_0\map_{\mathcal D(\Cond{(\kappa)}(\Ab))}(\mathbb Z[\underline{X}_{(\kappa)}], \mathbb Z)}, \mathbb Z)\\ 
& \cong \iHom_{\Cond{(\kappa)}(\Ab)}(\underline{\mathcal C(X, \mathbb Z)^{\delta}}_{(\kappa)}, \mathbb Z)
\end{align*} naturally in the ($\kappa$-light) compact Hausdorff space $X$.

As $\pi_0X$ is profinite (\cref{connectedcomponentsofchisch}) and $\wt(\pi_0X)\leq \wt(X)$ (\cref{weightofquotients}), this implies that $X\to\pi_0X$ induces an isomorphism \[\mathbb Z[\underline{X}_{(\kappa)}]^{\solid}\cong \mathbb Z[\underline{\pi_0 X}_{(\kappa)}]^{\solid}\] and analogously for $\pi_0^{\gamma}X$, provided that this space is Hausdorff. 
\end{proof}

\begin{cor}\label{profinitesolid}
    Suppose $M$ is a topological abelian group. 
    \begin{romanenum}\item If $M$ is locally profinite (i.e.\ locally compact Hausdorff and totally disconnected), then $\underline M_{(\kappa)}$ is solid. This in particular applies to discrete abelian groups. 
    \item If there exists an uncountable cardinal $\kappa$ such that $\underline{M}_{\kappa}$ is solid, then $M$ is totally path-disconnected.
    \item If $M$ is Hausdorff and $\underline{M}$ is solid, then $M$ is totally path-disconnected. 
    \end{romanenum}
\end{cor}
\begin{proof}We only prove the statement for condensed abelian groups, the statement for $\kappa$-condensed abelian groups can be shown completely analogously.
If $A$ is a discrete abelian group, choose a presentation 
    $A=\Coker(\oplus_{i\in I}\mathbb Z\to \oplus_{j\in J}\mathbb Z)$. As $\text{const}\colon \Ab\to \Cond{}(\Ab)$ preserves colimits, \[\underline{A}=\Coker(\oplus_{i\in I}\underline{\mathbb Z}\to \oplus_{j\in J}\underline{\mathbb Z}).\] Since $\underline{\mathbb Z}$ is a solid abelian group (\cref{solidificationissolid}) and $\Sol{}\subseteq \Cond{}(\Ab)$ is closed under small colimits, this implies that $\underline{A}$ is solid. As $\underline{(-)}\colon \Ab(\To\Top)\to \Cond{}(\Ab)$ preserves limits and $\Sol{}$ is closed under limits, it follows that limits of discrete abelian groups and in particular profinite abelian groups are solid. By van Dantzig's theorem, every locally profinite abelian group $G$ has a profinite open normal subgroup $N\subseteq G$.
    Since $G/N$ is discrete, the projection $G\to G/N$ admits a continuous section (which need not be a group homomorphism) and in particular, $\underline{G}\to\underline{G/N}$ is an epimorphism of condensed sets. 
    This implies that the sequence of condensed abelian groups 
    \[0\to\underline{N}\to\underline{G}\to\underline{G/N}\to 0\] is exact. 
    Since $N$ is profinite and $G/N$ is discrete, $\underline{N}$ and $\underline{G/N}$ are solid abelian groups by the above. 
    As $\Sol{}\subseteq \Ab$ is closed under extensions, it follows that $\underline{G}$ is a solid abelian group. 

Since $\mathbb Z[\underline{[0,1]}]^{\solid}\cong \mathbb Z$ (\cref{solidcw}), if $M$ is a $\To$ topological abelian group such that $\underline{M}$ is solid, then
\[\cont([0,1],M)\cong \Hom_{\Sol{}}(\mathbb Z[\underline{[0,1]}_{}]^{\solid}, \underline M)\cong \Hom_{\Sol{}}(\mathbb Z, \underline M)\cong M, \] i.e.\ every continuous path $[0,1]\to M$ is constant. 
\end{proof}

For a topological space $X$, denote by $\mathbb Z[\Sigma^{\infty}_{+}X]\in \mathcal D(\Ab)\cong \LMod{H\mathbb Z}{\Sp}$ the free $\mathbb Z$-module on $\Sigma^{\infty}_{+}X\in\Sp$. 
Denote by \[ \cshv{-}\colon\mathcal D(\Ab)\to\mathcal D(\Cond{(\kappa)}(\Ab))\] the functor induced by the constant sheaf functor (\cref{constantsheafsymmetricmonoidal,localisationinduceslocalisationonmodulecategories,condensedmodulesaremodulesincondensed}).  
\begin{cor}[{\cite[Example 6.5]{Scholzecondensed}}]\label{solidcw}
    If $X$ is a CW-complex, then \[\mathbb Z[\underline X_{(\kappa)}]^{L\solid}\cong \cshv{\mathbb Z[\Sigma^{\infty}_{+}X]}\in \mathcal D(\Sol{(\kappa)}{}).\]
\end{cor}
\begin{rem}
The singular chains functor $C_*^{\operatorname{sing}}\colon \operatorname{Top}\to \Ch(\Ab)\to \mathcal D(\Ab)$ descends to a functor $C_*^{\operatorname{sing}}(-)\colon \an\to \mathcal D(\Ab)$. This functor preserves colimits and is therefore equivalent to $ \mathbb Z[-]\circ \Sigma^{\infty}_{+}$. 
\end{rem}
    \begin{proof}
        As every compact subspace of a CW-complex is contained in a finite subcomplex, \[\underline{X}_{(\kappa)}=\colim{F\subseteq X}\underline{F}_{(\kappa)}\] is the colimit of its finite subcomplexes. 
        It therefore suffices to show that for finite CW-complexes $F$, there is an equivalence $\mathbb Z[\underline{F}_{(\kappa)}]^{L\solid}\cong \cshv{\mathbb Z[\Sigma^{\infty}_{+}F]}$ which is natural with respect to continuous maps. 
        \cref{weightofquotients} implies that finite CW-complexes are light compact Hausdorff spaces. Fix a finite CW-complex $F$.
        By \cref{cohomologyclausenscholze}, 
        \[ \map_{\mathcal D(\Cond{(\kappa)}(\Ab))}(\mathbb Z[\underline{F}_{(\kappa)}], \mathbb Z)\cong \cH{\sheaf}(F, \mathbb Z)\] naturally in $F$. 
        Since $F$ is locally contractible, the functor \[ C^{*}_{\sing}(-)\colon \Op(F)^{\operatorname{op}}\to \Ch(\Ab)\] which sends an open subset to its complex of singular cochains is a resolution of $\mathbb Z$. For all $k\in\mathbb N_0$, $C^k_{\operatorname{sing}}(-)$ is flasque, whence \[\cH{\sheaf}(F, \mathbb Z)\cong C^*_{\sing}(F)=\Hom_{\Ab}(C_*^{\operatorname{sing}(F)}, \mathbb Z)\] is represented by the singular cochain complex of $F$. This identification is natural in $F$. 
        Since $C_*^{\operatorname{sing}}(F)\cong \mathbb Z[\Sigma^{\infty}_{+}F]$, it follows that \[ \Hom_{\Ab}(C_*^{\operatorname{sing}(F)}, \mathbb Z) \cong \map_{\mathcal D(\Ab)}(\mathbb Z[\Sigma^{\infty}_{+}F], \mathbb Z)\] naturally in $F$.

        By \cref{constantsheafsymmetricmonoidal,symmetricmonoidalstructuremodulesnatural,condensedmodulesaremodulesincondensed}, $c(-)$ enhances to a symmetric monoidal functor. 
        In particular, \[\cshv{\mathbb Z[\Sigma^{\infty}_{+}(-)]}\otimes \cshv{-}\cong \cshv{-}\circ (\mathbb Z[\Sigma^{\infty}_{+}-]\otimes-).\] This determines an equivalence \[ \map_{\mathcal D(\Ab)}(\mathbb Z[\Sigma^{\infty}_{+}F],-)\circ \Gamma\cong \map_{\mathcal D(\Cond{\kappa}(\Ab))}(\cshv{\mathbb Z[\Sigma^{\infty}_{+}F]},-)\] between their respective right adjoints, which is natural in $F$. As $\Gamma$ is $t$-exact (\cref{globalsectionstexactcocontinuous}), $\Gamma(\mathbb Z)=\mathbb Z$, whence \[  \map_{\mathcal D(\Cond{(\kappa)}(\Ab))}(\mathbb Z[\underline{F}_{(\kappa)}], \mathbb Z)\cong \map_{\mathcal D(\Cond{(\kappa)}(\Ab))}(\cshv{\mathbb Z[\Sigma^{\infty}_{+}F]}, \mathbb Z)\] naturally in $F$. 
        \cref{solidificationoftopologicalspaces} now implies that  
        \[ \mathbb Z[\underline{F}_{(\kappa)}]^{L\solid}\cong \imap_{\mathcal D(\Cond{(\kappa)}(\Ab))}(\cshv{\map_{\mathcal D(\Cond{(\kappa)}(\Ab))}(\cshv{H\mathbb Z[\Sigma^{\infty}_{+}F]}, \mathbb Z)}, \mathbb Z)\] naturally in $F$. 
        Since $\cshv{-}\colon\mathcal D(\Ab)\to\mathcal D(\Cond{(\kappa)}(\Ab))$ is $t$-exact, \[H_i(\cshv{\mathbb Z[\Sigma^{\infty}_{+}F]})\cong \cshv{H_i(\mathbb Z[\Sigma^{\infty}_{+}F])}\] is solid for all $i\in\mathbb Z$, whence $\cshv{\mathbb Z[\Sigma^{\infty}_{+}F]}\in\mathcal D(\Sol{})$ by \cref{derivedsolidificationunderlyingforacyclicrings}. 
        We show by induction on $\operatorname{dim}(F)$ that the map \[\eta_{\cshv{\mathbb Z[\Sigma^{\infty}_{+}F]}}\colon \underbrace{\left(\cshv{\mathbb Z[\Sigma^{\infty}_{+}F]}\right)^{L\solid}}_{=\cshv{\mathbb Z[\Sigma^{\infty}_{+}F]}}\to \imap_{\mathcal D(\Cond{(\kappa)}(\Ab))}(\cshv{\map_{\mathcal D(\Cond{(\kappa)}(\Ab))}(\cshv{\mathbb Z[\Sigma^{\infty}_{+}F]}, \mathbb Z)}, \mathbb Z)\] from \cref{solidificationoftopologicalspaces} is an equivalence.  
        This then implies the statement since the target is naturally equivalent to $\mathbb Z[\underline{F}_{(\kappa)}]^{L\solid}$ by the above. 
        If $\operatorname{dim}(F)=0$, this holds by \cref{solidificationoftopologicalspaces}. 
        The inductive step follows as for an $(n+1)$-dimensional finite CW-complex $F$ with $n$-sceleton $F^n$, \[\mathbb Z[\Sigma^{\infty}_{+}F^n]\to \mathbb Z[\Sigma^{\infty}_{+}F]\to \oplus_{n+1\text{-cells}}\mathbb Z[n+1]\] is a cofiber sequence in $\mathcal D(\Ab)$ and the full subcategory $\mathcal C\subseteq \mathcal D(\Cond{}(\Ab))$ on objects $M$ for which $\eta_M$ is an equivalence is a stable subcategory. 
    \end{proof} 
\subsubsection{Derived solid modules over $s$-flat rings}\label{section:sflatrings}
Since the category $\Sol{(\kappa)}(\algebra{R})$ generally does not have enough flat objects, it is not true for general condensed rings $\algebra{R}$ that 
\[ \mathcal D(\Sol{(\kappa)}(\algebra{R}))\cong \LMod{\algebra{R}^{\solid\mathbb Z}}{\mathcal D(\Sol{(\kappa)})}.\] 
We now characterize the rings for which this holds and give many examples below. 
\begin{definition}\label{definitionsflatrings}
A ($\kappa$-)condensed ring $\algebra{R}$ is ($\kappa$)-$s$-\emph{flat} if for all $S\in\Pro(\Fin)_{(\kappa)}$, 
\[ \algebra{R}^{\solid\mathbb Z}\otimes_{\mathcal D(\Sol{(\kappa)})}\mathbb Z[\underline{S}_{(\kappa)}]^{\solid\mathbb Z}\in \mathcal D(\Sol{(\kappa)})^{\heart},\] where $(-)^{\solid\mathbb Z}$ denotes the left adjoint of $\Sol{(\kappa)}\subseteq \Cond{(\kappa)}(\Ab)$. 
\end{definition}

This definition is motivated by the following observation: 
\begin{proposition}\label{derivedsolidrmodulesismodulecategory}
    For $\algebra{R}\in\Alg(\Cond{(\kappa)}(\Ab))$, the forget functor $\Sol{(\kappa)}(\algebra{R})\to\Sol{(\kappa)}$ extends uniquely to a small colimits preserving, $t$-exact functor \[\mathcal D(\Sol{(\kappa)}(\algebra{R}))\to \mathcal D(\Sol{(\kappa)}).\]
    
   This functor factors over an equivalence \[ \mathcal D(\Sol{(\kappa)}(\algebra{R}))\cong \LMod{\algebra{R}^{\solid}}{\mathcal D(\Sol{(\kappa)})}\to\mathcal D(\Sol{(\kappa)})\] if and only if $R$ is ($\kappa$)-$s$-flat. 
\end{proposition}
\begin{proof}We first prove the $\kappa$-condensed statement. The existence and uniqueness of the functor $\mathcal D(\Sol{\kappa}(\algebra{R}))\to\mathcal D(\Cond{\kappa}(\algebra{R}))$ follows from \cite[Proposition C.3.1.1, C.3.2.1, Theorem C.5.4.9]{SAG} since $\Sol{\kappa}(\algebra{R}), \Cond{\kappa}(\algebra{R})$ are Grothendieck abelian (\cref{underivedsolidification,condensedmodulesgrothendieckabelian}), see also \cite[Proposition A.2]{cartiermodulesmattisweiss}. 

Suppose now that  $\algebra{R}$ is $\kappa\text{-}s$-flat. 
We want to apply \cite[Remark C.5.4.11]{SAG} to identify \[\mathcal D(\Sol{\kappa}(\algebra{R}))\cong \LMod{\algebra{R}^{\solid}}{\mathcal D(\Sol{\kappa})}\to\mathcal D(\Sol{\kappa}),\] and we thank Ko Aoki for pointing out this reference to us. 

Since $\Sol{\kappa}$ is Grothendieck abelian, the $t$-structure on $\mathcal D(\Sol{\kappa})_{\geq 0}$ is accessible by \cite[Proposition 1.3.5.21]{higheralgebra}. Hence by \cref{tstructureonmodulespresentable}, the $t$-structure on $\mathcal D(\Sol{\kappa})$ pulls back to a $t$-structure on $\LMod{\algebra{R}^{\solid}}{\mathcal D(\Sol{\kappa})}$ with \[\LMod{\algebra{R}^{\solid}}{\mathcal D(\Sol{\kappa})}_{\geq 0}\cong \LMod{\algebra{R}^{\solid}}{\mathcal D(\Sol{\kappa})_{\geq 0}}.\] 

As $\mathcal D(\Sol{\kappa})$ is presentable (\cref{derivedkappasolidmodulespresentable}) and its symmetric monoidal structure is closed (\cref{symmetricmonoidalstructureonderivedsolidabeliangroups}), $\LMod{\algebra{R}^{\solid}}{\mathcal D(\Sol{\kappa})}$ is presentable by \cite[Corollary 4.2.3.7]{higheralgebra}. 
By \cref{forgetfreeadjunctionmodules}, the forget functor $\LMod{\algebra{R}^{\solid}}{\mathcal D(\Sol{\kappa})}\to \mathcal D(\Sol{(\kappa)})$ reflects small limits and colimits. In particular, $\LMod{\algebra{R}^{\solid}}{\mathcal D(\Sol{(\kappa)})}$ is stable by \cref{derivedcategoryisstable}. 
As countable products and coproducts and filtered colimits in $\Sol{\kappa}$ are exact, ${\mathcal D(\Sol{\kappa})}_{\geq 0/\leq 0}\subseteq \mathcal D(\Sol{\kappa})$ are stable under countable products and coproducts and the $t$-structure is compatible with filtered colimits. 
As the forget functor $\LMod{\algebra{R}}{\mathcal D(\Sol{\kappa})}\to \mathcal D(\Sol{\kappa})$ reflects limits and colimits and is $t$-exact, this implies that \[\LMod{\algebra{R}}{\mathcal D(\Sol{\kappa})}_{\geq 0/\leq 0}\subseteq \LMod{\algebra{R}}{\mathcal D(\Sol{\kappa})}\] are stable under countable products and coproducts and the $t$-structure on $\LMod{\algebra{R}}{\mathcal D(\Sol{\kappa})}$ is compatible with filtered colimits. It now follows from \cite[Proposition 1.2.1.19]{higheralgebra} that the $t$-structure on $\mathcal D(\Sol{\kappa})$ is right and left-complete. The category 
$\LMod{\algebra{R}^{\solid}}{\mathcal D(\Sol{\kappa})}_{\geq 0}$ is prestable by \cite[Corollary C.1.2.3]{SAG}. 
As the forget functor \[\LMod{\algebra{R}^{\solid}}{\mathcal D(\Sol{\kappa})}_{\geq 0}\to\mathcal D(\Sol{\kappa})_{\geq 0}\] reflects small limits and colimits, \cite[Example C.1.4.5]{SAG} implies that $\LMod{\algebra{R}^{\solid}}{\mathcal D(\Sol{\kappa})}_{\geq 0}$ is Grothendieck prestable (\cite[Definition C.1.4.2]{SAG}). 

Next, we show that $\LMod{\algebra{R}^{\solid}}{\mathcal D(\Sol{\kappa})}_{\geq 0}$ is $0$-complicial (\cite[Definition C.5.3.1]{SAG}). 
Fix $M\in\LMod{\algebra{R}^{\solid}}{\mathcal D(\Sol{\kappa})}_{\geq 0}$. By \cite[Proposition C.5.3.2]{SAG}, there exists $Q\in\Sol{\kappa}$ with a map $q\colon Q\to fM\in\mathcal D(\Sol{\kappa})$ such that $H_0(q)$ is an epimorphism. By \cref{solidenoughprojectiveskappa}, we can assume that $Q=\oplus_{i\in I}\mathbb Z[\underline{S}_i]^{\solid}$ for a small family of $\kappa$-light profinite sets $S^i,i\in I$. Denote by $p\colon \algebra{R}^{\solid}[Q]\to M\in\LMod{\algebra{R}^{\solid}}{\mathcal D(\Sol{\kappa})}$ the adjoint map. 
As $R$ is ($\kappa$-)$s$-flat, by definition of the $t$-structure on $\LMod{\algebra{R}^{\solid}}{\mathcal D(\Sol{\kappa})}$, \[R^{\solid}[Q]\in \LMod{\algebra{R}^{\solid}}{\mathcal D(\Sol{\kappa})}^{\heart}.\] 
The map $H_0(p)$ is an epimorphism. Indeed: $H_0(\algebra{R}^{\solid}[Q])$ is the free $\algebra{R}^{\solid}$-module on $Q$ in $\Sol{\kappa}(\algebra{R})$, and $H_0(p)\colon H_0(\algebra{R}^{\solid}[Q])\to H_0(M)\in \LMod{\algebra{R}^{\solid}}{\Sol{\kappa}}$ is the map adjoint to $H_0(q)$.
As the forget functor $f^{\heart}\colon \LMod{\algebra{R}^{\solid}}{\Sol{\kappa}}\to \Sol{\kappa}$ preserves colimits and $H_0(q)$ factors as \[Q\xrightarrow{\eta}H_0(\algebra{R}^{\solid}[Q])\xrightarrow{f^{\heart}H_0(p)} H_0(M),\] it follows that \[0=\Coker(f^{\heart}H_0(p))=f^{\heart}\Coker(H_0(p)),\] which shows that $H_0(p)$ is an epimorphism. 
This shows that $\LMod{\algebra{R}^{\solid}}{\mathcal D(\Sol{\kappa})}_{\geq 0}$ is $0$-complicial. 

By \cref{solidmodulesismodulesinsolid}, \[\LMod{\algebra{R}^{\solid}}{\mathcal D(\Sol{\kappa})}^{\heart}\cong \LMod{\algebra{R}^{\solid}}{\Sol{\kappa}}.\] 
It now follows from \cite[Remark C.5.4.11]{SAG} that the equivalences of their hearts extends to a $t$-exact, small colimits preserving equivalence 
\[\mathcal D({\Sol{\kappa}(\algebra{R})})\cong \LMod{\algebra{R}^{\solid}}{\mathcal D(\Sol{\kappa})}.\] By \cite[Proposition C.3.1.1, C.3.2.1, Theorem C.5.4.9]{SAG}/\cite[Proposition A.2]{cartiermodulesmattisweiss}, this is essentially unique and the forget functor \[\mathcal D({\Sol{\kappa}(\algebra{R})})\to \mathcal D(\Sol{\kappa})\] factors over the above equivalence. 

Conversely, suppose that $\algebra{R}$ is a $\kappa$-condensed ring such that the forget functor factors over an equivalence \[ \mathcal D(\Sol{\kappa}(\algebra{R}))\cong \LMod{\algebra{R}^{\solid}}{\mathcal D(\Sol{\kappa})}\to\mathcal D(\Sol{\kappa}).\] 
    By \cref{existenceunboundedderivedfunctors}, there exists a unique small colimits preserving, right $t$-exact functor \[L\colon\mathcal D(\Sol{\kappa})\to\mathcal D(\Sol{\kappa}(\algebra{R}))\] such that for all projectives $P\in\Sol{\kappa}$, $L(P)\in\mathcal D(\Sol{\kappa}(\algebra{R}))^{\heart}$ and $H_0\circ L|_{\Sol{\kappa}}$ is the free $R^{\solid}$-module functor. 
    The same argument as in the proof of \cref{derivedsolificationexistsstronglimitcardinal} implies that this is a left adjoint to the forget functor $f\colon \mathcal D(\Sol{\kappa}(\algebra{R}))\to\mathcal D(\Sol{\kappa})$. 
    In particular, for $S\in\Pro(\Fin)_{\kappa}$ \[fL(\mathbb Z[\underline{S}_{\kappa}]^{\solid})\cong fR^{\solid}[\mathbb Z[\underline{S}_{\kappa}]^{\solid}]\cong \algebra{R}^{\solid}\otimes_{\mathcal D(\Sol{\kappa})}\mathbb Z[\underline{S}_{\kappa}]^{\solid}\in\mathcal D(\Sol{\kappa})^{\heart},\] which shows that $R$ is $\kappa$-$s$-flat.

We now deduce from this the condensed statements. 
Fix a condensed ring $\algebra{R}$ and choose a regular cardinal $\mu$ with $\algebra{R}\in\Alg(\Cond{\mu}(\Ab))$. 
By \cref{colimitdescriprionsolid} and \ref{solidclosedunderlimitscolimitskappa},
\[ \Fun^{\operatorname{colim}, \operatorname{ex}}(\Sol{}(\algebra{R}), \Sol{})\cong \clim{\substack{\kappa\geq \mu\\ \kappa\text{ regular}}}\Fun^{\operatorname{colim}, \operatorname{ex}}(\Sol{\kappa}(\algebra{R}), \Sol{}), \] where $\Fun^{\operatorname{colim}, \operatorname{ex}}$ refers to small colimits preserving, exact functors. 
By \cref{solkappafullyfaithfulonbderivedandpreservescolimits}, 
\[ \Fun^{\operatorname{colim}, \operatorname{t-ex}}(\Sol{}(\algebra{R}), \Sol{})\cong \clim{\substack{\kappa\geq \mu\\ \kappa\text{ regular}}}\Fun^{\operatorname{colim}, \operatorname{t-ex}}(\Sol{\kappa}(\algebra{R}), \Sol{}), \] where $\Fun^{\operatorname{colim}, \operatorname{t-ex}}$ refers to small colimits preserving, $t$-exact functors. This implies that 
\begin{align*}&\Fun^{\operatorname{colim}, \operatorname{t-ex}}(\mathcal D(\Sol{}(\algebra{R})), \mathcal D(\Sol{}))\times_{\Fun^{\operatorname{colim}, \operatorname{ex}}(\Sol{}(\algebra{R}), \Cond{}(\algebra{R}))}\{ f\} \\ &\cong \clim{\lambda\geq \mu} \Bigl(\Fun^{\operatorname{colim}, \operatorname{t-ex}}(\mathcal D(\Sol{\lambda}(\algebra{R})), \mathcal D(\Sol{}))\times_{\Fun^{\operatorname{colim}, \operatorname{ex}}(\Sol{\lambda}(\algebra{R}), \Sol{})}\{ f|_{\Sol{\lambda}(\algebra{R})}\}\Bigr)\\ & \cong \clim{\lambda \geq \mu}\Bigl(\Fun^{\operatorname{colim}, \operatorname{t-ex}}(\mathcal D(\Sol{\lambda}(\algebra{R})), \mathcal D(\Sol{\lambda}))\times_{\Fun^{\operatorname{colim}, \operatorname{ex}}(\Sol{\lambda}(\algebra{R}), \Sol{\lambda})}\{ f_{\lambda}\}\Bigr), \end{align*} where the limits run over all regular cardinals $\lambda\geq\mu$. By \cite[Proposition C.3.1.1, C.3.2.1, Theorem C.5.4.9]{SAG}/\cite[Proposition A.2]{cartiermodulesmattisweiss}, the right-hand side is contractible, i.e. there exists an essentially unique $t$-exact, cocontinuous extension $\mathcal D(\Sol{}(\algebra{R}))\to \mathcal D(\Sol{})$ of the forget functor. 

Suppose now that $\algebra{R}$ is $s$-flat. \cref{freeabeliangroupsacyclicforsolidificationwithoutkappa} implies that $\algebra{R}$ is $\kappa$-$s$-flat for all $\kappa\geq \mu$.  We claim that the forget functor \[ f\colon \mathcal D(\Sol{}(\algebra{R}))\to\mathcal D(\Sol{})\] admits a left-adjoint $L$, such that for all regular cardinals $\kappa\geq \mu$ , $L(\mathcal D(\Sol{\kappa}))\subseteq \mathcal D(\Sol{\kappa}(\algebra{R}))$. Then it follows from the above and the Barr-Beck theorem \cite[Theorem 4.7.0.3/Theorem 4.8.5.8]{higheralgebra} that $f$ factors over an equivalence \[\mathcal D(\Sol{}(\algebra{R}))\cong \LMod{\algebra{R}}{\mathcal D(\Sol{})}.\] 
Denote by $\Lambda^r$ the poset of uncountable regular cardinals. 
By the above, the forget functors assemble into a natural transformation 
\[f_*\colon\mathcal D(\Sol{*}(\algebra{R}))\to\mathcal D(\Sol{*})\in \Fun(\Lambda^r_{\geq \mu},\Pr^L)\] and for all $\kappa\geq \mu$, $f_{\kappa}$ admits a left ajoint $L_{\kappa}$. 
Hence by \cref{adjunctionsbigpresentable} applied to the opposite categories, it suffices to show that for all $\mu\leq \kappa\leq \lambda$, the mate \[ L_{\lambda}\circ i_{\kappa}^{\lambda}\to i_{\kappa,R}^{\lambda}\circ L_{\kappa}\] 
of the diagram 
\begin{center}
\begin{tikzcd}
    \mathcal D(\Sol{\kappa}(\algebra{R}))\arrow[d,"f_{\kappa}"]\arrow[r,"i_{\kappa,R}^{\lambda}"] & \mathcal D(\Sol{\lambda}(\algebra{R}))\arrow[d,"f_{\lambda}"]\\ 
    \mathcal D(\Sol{\kappa})\arrow[r,"i_{\kappa}^{\lambda}"] & \mathcal D(\Sol{\lambda})
\end{tikzcd}
\end{center} commutes, then it follows that $L\coloneqq \colim{\kappa}L_{\kappa}$ is a left adjoint of $f$ with \[L(\mathcal D(\Sol{\kappa}))=L_{\kappa}(\mathcal D(\Sol{\kappa}))\subseteq \mathcal D(\Sol{\kappa}(\algebra{R}))\] for all $\kappa\geq \mu$. As $L_{\lambda},i_{\kappa}^{\lambda},L_{\kappa},i_{\kappa,R}^{\lambda}$ are right $t$-exact, cocontinuous and carry projectives to the heart (since $\algebra{R}$ is $\kappa/\lambda$-$s$-flat), by \cref{existenceunboundedderivedfunctors}, it suffices to check that that for $M\in\Sol{\kappa}$, 
\[H_0(\beta)\colon H_0(L_{\lambda}\circ i_{\kappa}^{\lambda})(M)\to H_0(i_{\kappa,R}^{\lambda}\circ L_{\kappa})(M)\] is an equivalence. This holds since $\Sol{\kappa}\to \Sol{\lambda}$ is symmetric monoidal. 
This shows the if-statement for $\Sol{}$. 

Suppose now that $\algebra{R}$ is a condensed ring such that the forget functor factors over an equivalence \[ \mathcal D(\Sol{}(\algebra{R}))\cong \LMod{\algebra{R}^{\solid}}{\mathcal D(\Sol{})}\to\mathcal D(\Sol{}).\] By \cref{filteredcolimitsmodules,Filteredcolimitderivedcategories,colimitdescriprionsolid}, this implies that for all regular cardinals $\kappa$ with $\algebra{R}\in \Alg(\Cond{\kappa}(\Ab))$, the forget functor factors over an equivalence \[ \mathcal D(\Sol{\kappa}(\algebra{R}))\cong \LMod{\algebra{R}^{\solid}}{\mathcal D(\Sol{\kappa})},\] whence $\algebra{R}$ is $\kappa$-$s$-flat by the above statement on $\kappa$-solid modules. 
 As the functors \[i_{\kappa}\colon \mathcal D(\Sol{\kappa})\to\mathcal D(\Sol{})\] are symmetric monoidal and for $S\in\Pro(\Fin)_{\kappa}$, $i_{\kappa}(\mathbb Z[\underline{S}_{\kappa}]^{\solid})\cong \mathbb Z[\underline{S}]^{\solid}$ (\cref{freeabeliangroupsacyclicforsolidificationwithoutkappa}), this implies that $\algebra{R}$ is $s$-flat. 
\end{proof}

\begin{rem} By \cref{monoidalstructurederivedtensorprduct},$\kappa$-$s$-flatness is equivalent to $\algebra{R}^{\solid}\otimes_{\mathcal D(\Sol{(\kappa)})}-$ being the derived functor of $\algebra{R}^{\solid}\otimes_{\Sol{(\kappa)}}-$.
\end{rem}
 \begin{ex}\label{lightringssflat}If $\kappa=\aleph_1$, then $\mathbb Z[\underline{S}_{\kappa}]^{\solid}$ is flat for all $S\in\Pro(\Fin)_{\kappa}$ by \cref{Nobelingspecker} and \cite[Corollary 3.4.3]{Juansolidgeometry}, and hence all $\aleph_1$-condensed rings are $\aleph_1\text{-}s$-flat. 
 \end{ex}

\begin{lemma}\label{examplessflatrings}
    \begin{romanenum}\item If $\algebra{R}$ is a ($\kappa$-)condensed ring such that the underlying condensed abelian group $\underline{R}_{(\kappa)}^{\solid}$ admits a resolution $C_*\to \underline{R}_{(\kappa)}^{\solid}$ with $C_n=\mathbb Z[\underline{S^n}_{(\kappa)}], S^n\in\Pro(\Fin)_{(\kappa)}$ for all $n\in\mathbb N_0$, then $\algebra{R}$ is ($\kappa$)-$s$-flat.
    \item If $G$ is a Hausdorff topological group, the group ring $\mathbb Z[\underline{G}_{\kappa}]^{\solid}$ is $\kappa\text{-}s$-flat. 
    \end{romanenum}
\end{lemma}
\begin{proof}Suppose first that  $\algebra{R}$ is a ($\kappa$-)condensed ring such that the underlying ($\kappa$-)condensed abelian group $\underline{R}_{(\kappa)}^{\solid}$ admits a resolution $C_*\to \underline{R}_{(\kappa)}^{\solid}$ with $C_n=\mathbb Z[\underline{S^n}_{(\kappa)}], S^n\in\Pro(\Fin)_{(\kappa)}$.
    By exactness of filtered colimits in $\Sol{(\kappa)}$ and $\Cond{(\kappa)}(\Ab)$, for $T\in\Pro(\Fin)_{(\kappa)}$, \[\algebra{R}^{\solid}\otimes_{\mathcal D(\Sol{(\kappa)})}\mathbb Z[\underline{T}_{(\kappa)}]\cong \colim{m} (C_*^{\leq m}\otimes_{\mathcal D(\Sol{(\kappa)})}\mathbb Z[\underline{T}_{(\kappa)}]),\] where $C_*^{\leq m}$ denotes the stupid truncation.
    Using the cofiber sequence \[ C_*^{\leq m}\to C_*^{\leq m+1}\to C_{m+1}[m+1],\] it follows by induction on $m$ that $(C_*^{\leq m}\otimes_{\mathcal D(\Sol{(\kappa)})}\mathbb Z[\underline{T}_{(\kappa)}])$ is computed by degreewise applying $-\otimes_{\mathcal D(\Sol{(\kappa)})}\mathbb Z[\underline{T}_{(\kappa)}]^{\solid}$ to $C_*^{\leq m}$, and that this is isomorphic to the degreewise $I$-fold product $(C_*^{\leq m})^I$ where $I$ is a set with $\mathbb Z[\underline{T}_{(\kappa)}]^{\solid}\cong \prod_I \mathbb Z$.
    Hence by exactness of filtered colimits in $\Sol{(\kappa)}$, \[C_*\otimes_{\mathcal D(\Sol{(\kappa)})}\mathbb Z[\underline{T}_{(\kappa)}]^{\solid}\cong C_*^{I}.\] As ($\kappa$-small) products in $\Sol{(\kappa)}$ are exact, $C_*^{I}$ is a resolution of $(\algebra{R^{\solid}})^I$, and in particular \[\mathbb Z[\underline{T}_{(\kappa)}]^{\solid}\otimes_{\mathcal D(\Sol{\kappa})}\algebra{R}^{\solid}\cong (\algebra{R}^{\solid})^{I}\] is concentrated in degree $0$.  

    Suppose now that $G$ is a Hausdorff topological group.
    By \cref{approximationbycompacty}, \[ \colim{\substack{K\subseteq G\\ K\in \operatorname{CH}_{(\kappa)}}}\underline{K}_{(\kappa)}\cong \underline{G}_{(\kappa)},\] and hence \[\mathbb Z[\underline{G}_{(\kappa)}]^{\solid\mathbb Z}\cong \colim{\substack{K\subseteq G\\K\in\operatorname{CH}_{(\kappa)}}}\mathbb Z[\underline{K}_{(\kappa)}]^{\solid}.\] As the colimit is filtered and filtered colimits in $\Sol{(\kappa)}$ are exact (\cref{solidclosedunderlimitscolimitskappa,underivedsolidificationwithoutkappa,condensedabeliangroupsgrothendieckaxioms}), 
    \[ \mathbb Z[\underline{G}_{(\kappa)}]^{\solid\mathbb Z}\otimes_{\mathcal D(\Sol{(\kappa)})}\mathbb Z[\underline{T}_{(\kappa)}]^{\solid}\cong \colim{\substack{K\subseteq G\\K\in\operatorname{CH}_{(\kappa)}}} (\mathbb Z[\underline{K}_{(\kappa)}]^{\solid}\otimes_{\mathcal D(\Sol{(\kappa)})}\mathbb Z[\underline{T}_{(\kappa)}]^{\solid}).\] 
    By \cref{connectedcomponentsofchisch}, for a compact Hausdorff space $K$, the set of connected components $\pi_0K$ (with quotient topology from $K\to \pi_0K$) is profinite and by \cref{weightofquotients}, $\wt(\pi_0K)\leq \wt(K)$. Hence by \cref{solidificationcompacthausdorffspace}, for $K\in\CH_{(\kappa)}$, \[\mathbb Z[\underline{K}_{(\kappa)}]^{\solid}\otimes_{\mathcal D(\Sol{(\kappa)})}\mathbb Z[\underline{T}_{(\kappa)}]^{\solid}\cong \mathbb Z[\underline{\pi_0K\times T}_{(\kappa)}]^{\solid}\in \mathcal D(\Sol{(\kappa)})^{\heart},\] which shows that $\mathbb Z[\underline{G}_{(\kappa)}]$ is ($\kappa$-)$s$-flat. 
 \end{proof}
 \begin{ex}\label{profiniteringssflat}
 If $\algebra{R}$ is a ($\kappa$-light) profinite ring, then the Breen-Deligne resolution (\cite[Theorem 4.5]{Scholzecondensed}) for $\underline{\algebra{R}}_{(\kappa)}^{\solid}=\underline{\algebra{R}}_{(\kappa)}$ is a resolution as in $i)$ of the above lemma. In particular, $\underline{R}_{(\kappa)}$ is ($\kappa$-)$s$-flat. 
 \end{ex}
 \cref{derivedsolidrmodulesismodulecategory} implies the following: 
\begin{cor}\label{derivedsolidificationexistssflatrings}
    Suppose that $\algebra{R}$ is a ($\kappa$-)$s$-flat ($\kappa$-)condensed ring. 
The functor \[\mathcal D(\Sol{(\kappa)}(\algebra{R}))\to \mathcal D(\Cond{(\kappa)}(\algebra{R}))\] preserves small limits and colimits and admits a left adjoint $(-)^{L\solid\algebra{R}}$.
\end{cor} 
\begin{proof}
    We first show that the forget functor preserves small limits and colimits. 
    As the tensor product $-\otimes_{\mathcal D(\Sol{(\kappa)})}-$ commutes with colimits in both variables (\cref{symmetricmonoidalstructureonderivedsolidabeliangroups}) and $\mathcal D(\Sol{(\kappa)})$ has small limits and colimits (\cref{derivedkappasolidmodulespresentable,solkappafullyfaithfulonbderivedandpreservescolimits}), the forget functor \[ \mathcal D(\Sol{(\kappa)}(\algebra{R}))\cong \LMod{\algebra{R}^{\solid}}{\mathcal D(\Sol{(\kappa)})}\to \mathcal D(\Sol{(\kappa)})\] reflects small limits and colimits (\cref{forgetfreeadjunctionmodules}). 
    By \cref{derivedsolidificationexistsabeliangroupskappa,derivedsolidificationwithoutkappa},
    \[\mathcal D(\Sol{(\kappa)})\subseteq \mathcal D(\Cond{(\kappa)}(\Ab))\] is closed under small limits and colimits. 
    By \cref{modulesinderivedcategoryisderivedcategoryofmodules}, \[\mathcal D(\Cond{(\kappa)}(\algebra{R}))\cong \LMod{\algebra{R}}{\mathcal D(\Cond{(\kappa)}(\Ab))}\] and by \cref{derivedcondensedmodulesallsmallimitsandcolimits,closedmonoidalstructurecondensedcategories,forgetfreeadjunctionmodules}, the forget functor \[\mathcal D(\Cond{(\kappa)}(\algebra{R}))\cong \LMod{\algebra{R}}{\mathcal D(\Cond{(\kappa)}(\Ab))}\to \mathcal D(\Cond{(\kappa)}(\Ab))\] reflects small limits and colimits.
    As \begin{center}
        \begin{tikzcd}
        \mathcal D(\Sol{\kappa}(\algebra{R}))\arrow[r]\arrow[d] & \mathcal D(\Sol{\kappa})\arrow[d]\\ \mathcal D(\Cond{\kappa}(\algebra{R}))\arrow[r] & \mathcal D(\Cond{\kappa}(\Ab))
        \end{tikzcd}
    \end{center} commutes, this implies that 
    \[ f\colon \mathcal D(\Sol{(\kappa)}(\algebra{R}))\to \mathcal D(\Cond{(\kappa)}(\algebra{R}))\] preserves small limits and colimits. 
    The functor $\mathcal D(\Sol{}(\algebra{R}))\to\mathcal D(\Cond{}(\algebra{R}))$ has a left-adjoint by \cref{derivedsolidificationwithoutkappa}. 
    Since $\Sol{\kappa}(\algebra{R})$ and $\Cond{\kappa}(\algebra{R})$ are Grothendieck abelian, their derived categories $\mathcal D(\Cond{\kappa}(\Ab))$ and $\mathcal D(\Sol{\kappa}(\algebra{R}))$ are presentable by \cite[Proposition 1.3.5.21]{higheralgebra}.
    As $f$ preserves small limits and colimits, it follows from the adjoint functor theorem \cite[Corollary 5.5.2.9]{highertopostheory} that $f$ admits a left adjoint \[(-)^{L\solid\algebra{R}}\colon\mathcal D(\Cond{\kappa}(\algebra{R}))\to\mathcal D(\Sol{\kappa}(\algebra{R})).\qedhere\] 
    \end{proof}
\begin{cor}\label{derivedsolidmodulesarenriched}
Suppose $\algebra{R}\in\Alg(\Cond{(\kappa)}(\Ab))$ is a ($\kappa$)-$s$-flat condensed ring. 
The left-tensoring of \[\mathcal D(\Sol{(\kappa)}(\algebra{R}))\cong \LMod{\algebra{R}^{\solid}}{\mathcal D(\Sol{(\kappa)})}\] over $\mathcal D(\Sol{(\kappa)})$ described in \cref{modulecategorieslefttensored} exhibits 
$\LMod{\algebra{R}^{\solid}}{\mathcal D(\Sol{(\kappa)}(\algebra{R}))}$ as $\mathcal D(\Sol{(\kappa)})$-enriched. 
\end{cor}
\begin{proof}
    By \cref{derivedcondensedmodulesallsmallimitsandcolimits}/\cref{solkappafullyfaithfulonbderivedandpreservescolimits}, $\mathcal D(\Sol{(\kappa)}(\Ab))$ has all small limits and colimits, and by \cref{symmetricmonoidalstructureonderivedsolidabeliangroups}, the symmetric monoidal structure is closed. The statement now follows from \cref{modulecategoryisagainenriched}. 
\end{proof}
\begin{rem}\label{moregeneralenrichmentssolid}
The solidification and the constant sheaf functor yield symmetric monoidal left adjoints 
\[ \mathcal D(\Ab)\to \mathcal D(\Cond{(\kappa)})\to \mathcal D(\Sol{(\kappa)}).\] 
By the above, these functors exhibit $\mathcal D(\Sol{(\kappa)}(\algebra{R}))$ as enriched in $\mathcal D(\Cond{(\kappa)}(\Ab))$ and $\mathcal D(\Ab)$. 
Analogously, $\Sol{(\kappa)}(\algebra{R})$ is enriched in $\Sol{(\kappa)}$, $\Cond{(\kappa)}(\Ab)$ and $\Ab$. 
\end{rem}
\begin{lemma}\label{factorizationforgetfunctor2}
    For a ($\kappa$-$s$)-flat $\kappa$-condensed ring $\algebra{R}\in \Alg(\Cond{(\kappa)}(\Ab))$, the forget functor \[\mathcal D(\Sol{(\kappa)}(\algebra{R}))\to\mathcal D(\Cond{(\kappa)}(\algebra{R}))\] factors as 
    \begin{align*}\mathcal D(\Sol{(\kappa)}(\algebra{R}))\cong\LMod{\algebra{R}^{\solid}}{\mathcal D(\Sol{(\kappa)})}\xrightarrow{q^*}\LMod{\algebra{R}^{L\solid}}{\mathcal D(\Sol{(\kappa)})}\xhookrightarrow{f_{LM}^R} &\LMod{\algebra{R}}{\mathcal D(\Cond{(\kappa)}(\Ab))}\\ &\cong\mathcal D(\Cond{(\kappa)}(\algebra{R})),
\end{align*} 
    where 
    \begin{romanenum}
    \item the left equivalence is supplied by \cref{derivedsolidrmodulesismodulecategory} and the right equivalence by \cref{modulesinderivedcategoryisderivedcategoryofmodules},
    \item $q^*$ is restriction of scalars along the algebra map \[q\colon \algebra{R}^{L\solid}\to\algebra{R}^{\solid}\in \Alg(\mathcal D(\Sol{(\kappa)})_{\geq 0})\subseteq \Alg(\mathcal D(\Sol{(\kappa)}))\] provided by the unit of the localisation $ \Alg(\mathcal D(\Sol{(\kappa)}(\Ab))_{\geq 0})\to \Alg(\Sol{(\kappa)})$ induced by the symmetric monoidal left adjoint \[H_0\colon \mathcal D(\Sol{(\kappa)})_{\geq 0}\to \Sol{(\kappa)}\] (\cref{adjuncttioninducesadjunctiononalgebraobjects}).    
    \item $f_{LM}^R$ is the right adjoint induced by the symmetric monoidal localization \[ f\colon \mathcal D(\Sol{(\kappa)})\rightleftarrows \mathcal D(\Cond{(\kappa)}(\Ab))\colon (-)^{L\solid}  \] cf.\ \cref{localisationinduceslocalisationonmodulecategories}.
    \end{romanenum}
\end{lemma}
    \begin{proof}All functors above are $t$-exact and on hearts, the composite \begin{align*}\Sol{(\kappa)}(\algebra{R})\cong\LMod{\algebra{R}^{\solid}}{\Sol{(\kappa)}}\xrightarrow{q^*}\LMod{\algebra{R}^{L\solid}}{\mathcal D(\Sol{(\kappa)})}^{\heart}\xhookrightarrow{f_{LM}^R} &\LMod{\algebra{R}}{\mathcal D(\Cond{(\kappa)})}^{\heart}\\&\cong\Cond{(\kappa)}(\algebra{R})\end{align*} is the forget functor $\Sol{(\kappa)}(\algebra{R})\subseteq \Cond{(\kappa)}(\algebra{R})$.  
    Hence by \cref{derivedsolidificationwithoutkappa} or \cite[Proposition C.3.1.1, C.3.2.1, Theorem C.5.4.9]{SAG}/\cite[Proposition A.2]{cartiermodulesmattisweiss} respectively, it suffices to show that $q^*$ and $f_{\LM}$ preserve small colimits. 
    This follows from \cref{forgetfreeadjunctionmodules} since the symmetric monoidal structures on $\mathcal D(\Sol{(\kappa)})$ and $\mathcal D(\Cond{(\kappa)}(\Ab))$ are cocontinuous (they are closed), $\mathcal D(\Cond{(\kappa)}(\Ab))$ is cocomplete and $\mathcal D(\Sol{(\kappa)})\subseteq \mathcal D(\Cond{(\kappa)}(\Ab))$ is closed under small colimits by \cref{derivedforgetfunctorsolidpreservescolimits}, \cref{derivedsolidificationwithoutkappa}.  \end{proof} 

    For a ($\kappa$-)$s$-flat commutative condensed ring $\algebra{R}\in\CAlg(\Cond{(\kappa)}(\Ab))$, endow \[\mathcal D(\Sol{(\kappa)}(\algebra{R}))\cong \LMod{\algebra{R}^{\solid}}{\mathcal D(\Sol{(\kappa)})}\] with the symmetric monoidal structure described in \cref{symmetricmonoidalstructure}. 
\begin{cor}\label{symmetricmonoidalstructuresolidrmodules}
    Suppose that $\algebra{R}\in\CAlg(\Cond{(\kappa)}(\Ab))$ is a ($\kappa$-)$s$-flat commutative condensed ring. 
    The symmetric monoidal structure on $\mathcal D(\Sol{(\kappa)}(\algebra{R}))$ is closed and \[(-)^{L\solid\algebra{R}}\colon\mathcal D(\Cond{(\kappa)}(\algebra{R}))\to\mathcal D(\Sol{(\kappa)}(\algebra{R}))\] enhances to a symmetric monoidal functor. 
\end{cor}
\begin{proof}
    By \cref{derivedkappasolidmodulespresentable} and \cref{solkappafullyfaithfulonbderivedandpreservescolimits}, respectively, the category $\mathcal D(\Sol{(\kappa)}(\algebra{R}))$ has all small limits and colimits. 
    As the symmetric monoidal structure on $\mathcal D(\Sol{(\kappa)})$ is closed (\cref{symmetricmonoidalstructureonderivedsolidabeliangroups}), \cref{internalhommodules} now implies that the symmetric monoidal structure on \[\mathcal D(\Sol{(\kappa)}(\algebra{R}))\cong \LMod{\algebra{R}^{\solid}}{\mathcal D(\Sol{(\kappa)}(\Ab))}\] is closed. 
    Since $(-)^{\solid}=H_0\circ (-)^{L\solid}$ is symmetric monoidal, the adjunction unit enhances to a morphism of commutative algebras \[\eta\colon \algebra{R}\to\algebra{R}^{\solid}\in\CAlg(\Cond{(\kappa)}(\Ab))\subseteq \CAlg(\mathcal D(\Cond{(\kappa)}(\Ab))).\] This exhibits $\algebra{R}^{\solid}$ as commutative algebra in $\LMod{\algebra{R}}{\mathcal D(\Cond{(\kappa)}(\Ab))}$ and by \cite[Corollary 3.4.1.9]{higheralgebra}, \[\LMod{\algebra{R}^{\solid}}{\LMod{\algebra{R}}{\mathcal D(\Cond{(\kappa)}(\Ab))}}^{\otimes}\cong \LMod{\algebra{R}^{\solid}}{\mathcal D(\Cond{(\kappa)}(\Ab))}^{\otimes}\] as symmetric monoidal categories. 
    By \cite[Corollary 4.8.5.22]{higheralgebra}, the free $\algebra{R}$-module functor enhances to a symmetric monoidal functor \begin{align*}\algebra{R}^{\solid}[-]^{\otimes} \colon \mathcal D(\Cond{(\kappa)}(\algebra{R}))\cong \LMod{\algebra{R}}{\mathcal D(\Cond{(\kappa)}(\Ab))}\to &\LMod{\algebra{R}^{\solid}}{\LMod{\algebra{R}}{\mathcal D(\Cond{(\kappa)}(\Ab))}}\\ & \LMod{\algebra{R}^{\solid}}{\mathcal D(\Cond{\kappa}(\Ab))}\\&\cong \mathcal D(\Cond{(\kappa)}(\algebra{R}^{\solid})), \end{align*} where we used the identifications from \cref{condensedmodulesarederivedcats} (which we also used to define symmetric monoidal structures on the derived category).  
    
    Let $S\coloneqq \algebra{R}^{\solid}$. 
    Then \[(-)^{L\solid S}\circ S[-]\colon \mathcal D(\Cond{(\kappa)}(\algebra{R}))\to \mathcal D(\Cond{(\kappa)}(\algebra{R}^{\solid}))\] is left adjoint to \begin{align*}\mathcal D(\Sol{(\kappa)}(\algebra{R}))\cong\LMod{\algebra{R}^{\solid}}{\mathcal D(\Sol{(\kappa)})}\hookrightarrow&\LMod{\algebra{R}^{\solid}}{\mathcal D(\Cond{(\kappa)}(\Ab))}\\ &\xrightarrow{\eta^*}\LMod{\algebra{R}}{\mathcal D(\Cond{(\kappa)})}(\Ab).\end{align*} This functor is $t$-exact, preserves small colimits and restricts to the forget functor on hearts, and is therefore equivalent to the forget functor $\mathcal D(\Sol{(\kappa)}(\algebra{R}))\to\mathcal D(\Cond{(\kappa)}(\algebra{R}))$ by \cref{derivedsolidificationwithoutkappa} and \cite[Proposition C.3.1.1, C.3.2.1, Theorem C.5.4.9]{SAG}/\cite[Proposition A.2]{cartiermodulesmattisweiss}, respectively. Hence, $(-)^{L\solid\algebra{R}}\cong (-)^{L\solid S}\circ S[-]$, and we are reduced to showing that $(-)^{L\solid S}$ enhances to a symmetric monoidal functor.
Since $S^{L\solid\mathbb Z}\cong S$,  by \cref{factorizationforgetfunctor2}, \[(-)^{L\solid S}\colon \mathcal D(\Cond{(\kappa)}(S))\cong \LMod{S}{\mathcal D(\Cond{(\kappa)}(\Ab))}\to \LMod{S^{L\solid}}{\mathcal D(\Sol{(\kappa)})}\cong \mathcal D(\Sol{(\kappa)}(S))\] is the functor induced by the symmetric monoidal functor \[(-)^{L\solid\mathbb Z}\colon\mathcal D(\Cond{(\kappa)}(\Ab))\to\mathcal D(\Sol{(\kappa)}(\Ab)).\]
As $\mathcal D(\Sol{(\kappa)}), \mathcal D(\Cond{(\kappa)}(\Ab))$ have small colimits, and the tensor products on both categories preserves small colimits in both variables (\cref{symmetricmonoidalstructureonderivedsolidabeliangroups,internalhomandenrichmentderivedcondensedcategories}), $(-)^{L\solid\mathbb Z}$ enhances to a symmetric monoidal functor \[\LMod{S}{\mathcal D(\Cond{(\kappa)}(\Ab))}^{\otimes}\to\LMod{S}{\mathcal D(\Sol{(\kappa)})}^{\otimes}\] by \cref{symmetricmonoidalstructuremodulesnatural}. 
\end{proof}

\begin{lemma}\label{symmetricmonoidalstructuresolidrmodulescompatibletstructure}Suppose $\algebra{R}\in \CAlg(\Cond{(\kappa)}(\Ab))$ is ($\kappa$-)$s$-flat. The symmetric monoidal structure on $\mathcal D(\Sol{(\kappa)}(\algebra{R}))$ provided by the above lemma is compatible with the t-structure, i.e.\ \[\mathcal D(\Sol{(\kappa)}(\algebra{R}))_{\geq 0}\subseteq\mathcal D(\Sol{(\kappa)}(\algebra{R}))\] is a symmetric monoidal subcategory. 
    
The induced symmetric monoidal structure on the heart $\Sol{(\kappa)}(\algebra{R})$ is the one described in \cref{solidmodulesismodulesinsolid}. 
\end{lemma}
\begin{proof}
    By construction of the symmetric monoidal structure, restriction along the commutative algebra map $\algebra{R}\to\algebra{R}^{\solid}$ defines a symmetric monoidal equivalence $\mathcal D(\Sol{(\kappa)}(\algebra{R}^{\solid}))\cong \mathcal D(\Sol{(\kappa)}(\algebra{R}))$. This functor is $t$-exact, we can therefore assume that $\algebra{R}=\algebra{R}^{\solid}\in \CAlg(\Sol{(\kappa)})$. Denote by \[f_R\colon\mathcal D(\Sol{(\kappa)}(R))\rightleftarrows \mathcal D(\Cond{(\kappa)}(R))\colon (-)^{L\solid R}\] the solidification adjunction. 
    Since $R^{L\solid}\cong R$, $f_{\algebra{R}}$ is fully faithful by \cref{derivedsolidificationunderlyingforacyclicrings} below, i.e.\ $(-)^{L\solid \algebra{R}}f_{\algebra{R}}=1$. As $(-)^{L\solid \algebra{R}}$ is symmetric monoidal, for $M,N\in\mathcal D(\Sol{(\kappa)}(\algebra{R}))$,
\[ (M\otimes_{\mathcal D(\Sol{(\kappa)}(\algebra{R}))}N)\cong (f_{\algebra{R}}M)^{L\solid \algebra{R}}\otimes_{\mathcal D(\Sol{(\kappa)}(\algebra{R}))}(f_{\algebra{R}}N)^{L\solid \algebra{R}}\cong (f_{\algebra{R}}M\otimes_{\mathcal D(\Cond{(\kappa)}(\algebra{R}))}f_{\algebra{R}}N)^{L\solid \algebra{R}}.\]
Since $f_{\algebra{R}}$ is $t$-exact, for $M,N\in\mathcal D(\Sol{(\kappa)}(\algebra{R}))_{\geq 0}$, 
\[ (f_{\algebra{R}}M\otimes_{\algebra{R}}f_{\algebra{R}}N)\in\mathcal D(\Cond{(\kappa)}(\algebra{R}))_{\geq 0}\] by \cref{monoidalstructurecompatibletstructure}.  As left adjoint of a $t$-exact functor, $(-)^{L\solid \algebra{R}}$ is right $t$-exact, so in particular 
\[ (M\otimes_{\mathcal D(\Sol{(\kappa)}(\algebra{R}))}N)=(fM\otimes_{\mathcal D(\Sol{(\kappa)}(\algebra{R}))}fN)^{L\solid \algebra{R}}\in \mathcal D(\Sol{(\kappa)}(\algebra{R}))_{\geq 0}.\] 
Since $\mathcal D(\Sol{(\kappa)}(\algebra{R}))_{\geq 0}\subseteq \mathcal D(\Sol{(\kappa)}(\algebra{R}))$ contains the unit $\algebra{R}^{\solid}$, it now follows from \cite[Remark 2.2.1.2]{higheralgebra} that $\mathcal D(\Sol{(\kappa)}(\algebra{R}))_{\geq 0}\subseteq \mathcal D(\Sol{(\kappa)}(\algebra{R}))$ is a symmetric monoidal subcategory. 
In particular, by \cite[Example 2.2.1.10]{higheralgebra}, $\Sol{(\kappa)}(\algebra{R})$ inherits a symmetric monoidal structure such that $\pi_0\colon \mathcal D(\Sol{(\kappa)}(\algebra{R}))_{\geq 0}\to \Sol{(\kappa)}(\algebra{R})$ is symmetric monoidal. This implies that $(-)^{\solid\algebra{R}}$ is symmetric monoidal as well, which shows that the induced symmetric monoidal structure on $\Sol{(\kappa)}(\algebra{R})$ is the one described in \cref{solidmodulesismodulesinsolid}. 
\end{proof}
\begin{rem}
In the above situation, the tensor product $-\otimes_{\mathcal D(\Sol{(\kappa)}(\algebra{R}))}-$ is the derived functor of the tensor product $-\otimes_{\Sol{(\kappa)}(\algebra{R})}-$: \cref{existenceunboundedderivedfunctors} and \cref{solkappafullyfaithfulonbderivedandpreservescolimits} imply that for a general commutative ($\kappa$-)condensed ring $\algebra{R}$, there is an essentially unique small colimits preserving, right $t$-exact functor \[(-\otimes_{\Sol{(\kappa)}(\algebra{R})}^L-)\colon \mathcal D(\Sol{(\kappa)}(\algebra{R}))\times \mathcal D(\Sol{(\kappa)}(\algebra{R}))\to\mathcal D(\Sol{(\kappa)}(\algebra{R}))\] such that 
\[H_0\circ (-\otimes_{\Sol{(\kappa)}(\algebra{R})}^L-)|_{\Sol{(\kappa)}(\algebra{R})\times \Sol{(\kappa)}(\algebra{R})}=-\otimes_{\Sol{(\kappa)}(\algebra{R})}-\] and for projectives $P,Q\in\Sol{(\kappa)}(\algebra{R})$, $P\otimes_{\Sol{(\kappa)}(\algebra{R})}^L Q\in \mathcal D(\Sol{(\kappa)}(\algebra{R}))^{\heart}$. 
We call this the derived functor of $-\otimes_{\Sol{(\kappa)}(\algebra{R})}-$. 

If $R$ is ($\kappa$)-$s$-flat, then $-\otimes_{\mathcal D(\Sol{(\kappa)}(\algebra{R}))}-$ is such a functor. 
For a general commutative, ($\kappa$-)condensed ring $\algebra{R}$, we are not aware of a way to construct a symmetric monoidal structure on $\mathcal D(\Sol{(\kappa)}(\algebra{R}))$ which has $-\otimes_{\Sol{(\kappa)}(\algebra{R})}^L-$ as underlying tensor product. 
\end{rem}

\cref{factorizationforgetfunctor2} implies: 
 \begin{lemma}\label{derivedsolidificationunderlyingforacyclicrings}
     Suppose $\algebra{R}\in \Alg(\Cond{(\kappa)}(\Ab))$ is a ($\kappa$-)$s$-flat ring which is $(-)^{\solid}$-acyclic, i.e.\ $\algebra{R}^{L\solid}\cong\algebra{R}^{\solid}$. 
     \begin{romanenum}
     \item\label{cohomologysolidacyclicrings} The functor $f\colon \mathcal D(\Sol{(\kappa)}(\algebra{R}))\to\mathcal D(\Cond{(\kappa)}(\algebra{R}))$ is fully faithful with essential image
     \[ \{M\in\mathcal D(\Cond{(\kappa)}(\algebra{R}))\, |\, H_i(M)\in\Sol{(\kappa)}(\algebra{R}) \text{ for all }i\in\mathbb Z\}.\] 
     \item Denote by $g\colon\mathcal D(\Cond{\kappa}(\algebra{R}))\to\mathcal D(\Cond{\kappa}(\Ab))$ and $g^{\solid}\colon\mathcal D(\Sol{\kappa}(\algebra{R}))\to\mathcal D(\Sol{\kappa})$ the forget functors. Then $g^{\solid}\circ (-)^{L\solid\algebra{R}}\cong (-)^{L\solid\mathbb Z}\circ g$.
     \item \label{internalhomsolidrmodulesissolid} For $A\in\mathcal D(\Cond{(\kappa)}(\algebra{R})), B\in\mathcal D(\Sol{(\kappa)}(\algebra{R}))$, 
    \[\imap_{\mathcal D(\Cond{(\kappa)}(\algebra{R}))}(A,B)\in\mathcal D(\Sol{(\kappa)}(\algebra{R}))\subseteq \mathcal D(\Cond{(\kappa)}(\algebra{R})).\]
    In particular, $f\colon \mathcal D(\Sol{(\kappa)}(\algebra{R}))\to\mathcal D(\Cond{(\kappa)}(\algebra{R}))$ preserves the internal Hom.
     \end{romanenum} 
 \end{lemma} 

 \begin{proof} We first show \ref{cohomologysolidacyclicrings}. By \cref{factorizationforgetfunctor2}, $f\colon \mathcal D(\Sol{(\kappa)}(\algebra{R}))\to\mathcal D(\Cond{(\kappa)}(\algebra{R}))$ is fully faithful and its essential image consists of those $M\in\mathcal D(\Cond{(\kappa)}(\algebra{R}))$ such that \[gM\in \im(\mathcal D(\Sol{(\kappa)})\to\mathcal D(\Cond{(\kappa)}(\Ab))).\] 
    By \cref{derivedsolidissolidderived,derivedsolidissolidderivedwithoutkappa}, these are precisely those $M$ with $H_i(gM)\in\Sol{(\kappa)}$ for all $i\in\mathbb Z$. As $g$ is $t$-exact, this shows \ref{cohomologysolidacyclicrings}. 
Since $\algebra{R}^{L\solid\mathbb Z}\cong \algebra{R}^{\solid\mathbb Z}$, \cref{factorizationforgetfunctor2} implies that \[(-)^{L\solid\algebra{R}}\colon \mathcal D(\Cond{(\kappa)}(\algebra{R}))\cong \LMod{\algebra{R}}{\mathcal D(\Cond{(\kappa)}(\Ab))}\to\mathcal D(\Sol{(\kappa)}(\algebra{R}))\cong \LMod{\algebra{R}^{L\solid\mathbb Z}}{\mathcal D(\Sol{(\kappa)})}\] is the functor induced by the symmetric monoidal functor $(-)^{L\solid\mathbb Z}$ (\cref{localisationinduceslocalisationonmodulecategories}), whence \[g^{\solid}\circ (-)^{L\solid\algebra{R}}\cong (-)^{L\solid\mathbb Z}\circ g,\] cf.\ \cref{naturalitylocalizationmodulecategories}. 
    It remains to show that for $B\in\mathcal D(\Sol{(\kappa)}(\algebra{R}))$, $A\in\mathcal D(\Cond{(\kappa)}(\algebra{R}))$, \[\imap_{\mathcal D(\Cond{(\kappa)}(\algebra{R}))}(A,B)\in\mathcal D(\Sol{(\kappa)}(\algebra{R}))\subseteq \mathcal D(\Cond{(\kappa)}(\algebra{R})).\] 
    Since $\mathcal D(\Sol{(\kappa)}(\algebra{R}))\subseteq\mathcal D(\Cond{(\kappa)}(\algebra{R}))$ is closed under limits, 
    it suffices to show this for free modules $A=\algebra{R}[c],c\in\mathcal D(\Cond{(\kappa)}(\Ab))$ which generate \[\mathcal D(\Cond{(\kappa)}(\algebra{R}))\cong \LMod{\algebra{R}}{\mathcal D(\Cond{(\kappa)}(\Ab))}\] under $\Delta^{\operatorname{op}}$-indexed colimits (\cref{modulesinderivedcategoryisderivedcategoryofmodules} and \cite[Proposition 4.7.3.14]{higheralgebra}). 
    By the above description of the essential image of $f$, it is enough to show that \begin{align}\label{conditionderivedmappingspacesolidfreemodules}g(\imap_{\mathcal D(\Cond{(\kappa)}(\algebra{R}))}(\algebra{R}[c],B))\in\mathcal D(\Sol{(\kappa)})\end{align} for all $c\in\mathcal D(\Cond{(\kappa)}(\Ab))$. 
    The symmetric monoidal structure of the free $\algebra{R}$-module functor yields an equivalence \[ \imap_{\mathcal D(\Cond{(\kappa)}(\Ab))}(c,-)\circ g\cong g\circ \imap_{\mathcal D(\Cond{(\kappa)}(\algebra{R}))}(\algebra{R}[c],-).\] 
    Since $gB\in\mathcal D(\Sol{(\kappa)})$, (\ref{conditionderivedmappingspacesolidfreemodules}) now follows from \cref{derivedmappingspacesaresolid}.
\end{proof}
For a $(\kappa)$-$s$-flat, $(-)^{\solid\mathbb Z}$-acyclic ring $\algebra{R}$, the essential image of the forget functor \[\mathcal D(\Sol{(\kappa)}(\algebra{R}))\hookrightarrow \mathcal D(\Cond{(\kappa)}(\algebra{R}))\] can be characterized as follows: 
    \begin{cor}\label{solidificationonderivedmappingspacesforacyclicrings}
        Suppose $\algebra{R}\in\Alg(\Cond{(\kappa)}(\Ab))$ is ($\kappa$-)$s$-flat and $\algebra{R}^{L\solid\mathbb Z}\cong\algebra{R}^{\solid\mathbb Z}$. 

        For $M\in\mathcal D(\Cond{\kappa}(\Ab))$, the following are equivalent: 
        \begin{romanenum}
            \item \label{essentialimagesolidrmodules}$M$ lies in the essential image of the forget functor $\mathcal D(\Sol{(\kappa)}(\algebra{R}))\hookrightarrow \mathcal D(\Cond{(\kappa)}(\algebra{R}))$. 
            \item \label{solidderivedrmodules}For all $S\in\Pro(\Fin)_{(\kappa)}$, the adjunction unit $\algebra{R}[\underline{S}_{(\kappa)}]\to\algebra{R}[\underline{S}_{(\kappa)}]^{L\solid\algebra{R}}$ induces an equivalence \[\map_{\mathcal D(\Cond{(\kappa)}(\algebra{R}))}(\algebra{R}[\underline{S}_{(\kappa)}]^{L\solid\algebra{R}},M)\cong \map_{\mathcal D(\Cond{(\kappa)}(\algebra{R}))}(\algebra{R}[\underline{S}_{(\kappa)}],M).\] 
        \item \label{solidderivedrmodulesinternally} For all $S\in\Pro(\Fin)_{(\kappa)}$, the adjunction unit $\algebra{R}[\underline{S}_{(\kappa)}]\to\algebra{R}[\underline{S}_{(\kappa)}]^{L\solid\algebra{R}}$ induces an equivalence \[\imap_{\mathcal D(\Cond{(\kappa)}(\algebra{R}))}(\algebra{R}[\underline{S}_{(\kappa)}]^{L\solid\algebra{R}},M)\cong \imap_{\mathcal D(\Cond{(\kappa)}(\algebra{R}))}(\algebra{R}[\underline{S}_{(\kappa)}],M).\] 
        \end{romanenum}
    \end{cor}
    \begin{proof}
        Clearly, $\ref{solidderivedrmodulesinternally}\Rightarrow \ref{solidderivedrmodules}$.

        We now show that $\ref{essentialimagesolidrmodules}\Rightarrow \ref{solidderivedrmodulesinternally}$. 
        As $(-)^{L\solid\algebra{R}}$ is symmetric monoidal, \cref{derivedsolidificationunderlyingforacyclicrings} implies that for  $A\in\mathcal D(\Cond{(\kappa)}(\algebra{R}))$, \[\imap_{\mathcal D(\Cond{(\kappa)}(\algebra{R}))}(A^{L\solid\algebra{R}},-)\circ f\] factors over $\mathcal D(\Sol{(\kappa)}\algebra{R})\subseteq \mathcal D(\Cond{(\kappa)}(\algebra{R}))$ and is right adjoint to \[A^{L\solid\algebra{R}}\otimes_{\mathcal D(\Sol{(\kappa)}(\algebra{R}))}-\cong (-)^{L\solid\algebra{R}}\circ (A\otimes_{\mathcal D(\Cond{(\kappa)}(\algebra{R}))}-).\] 
        This implies that \[\imap_{\mathcal D(\Cond{(\kappa)}(\algebra{R}))}(A,-)\circ f\cong \imap_{\mathcal D(\Cond{(\kappa)}(\algebra{R}))}(A^{L\solid\algebra{R}},-)\circ f, \] and in particular $\ref{essentialimagesolidrmodules}\Rightarrow \ref{solidderivedrmodulesinternally}$. 

        We now show that $ \ref{essentialimagesolidrmodules}\Leftrightarrow \ref{solidderivedrmodules}$. 
        Since $\algebra{R}[\underline{S}_{(\kappa)}], S\in\Pro(\Fin)_{(\kappa)}$ generates $\mathcal D(\Cond{(\kappa)}(\algebra{R}))$ under small colimits and shifts,
        \[\mathcal D(\Sol{(\kappa)}(\algebra{R}))\subseteq\mathcal D(\Cond{(\kappa)}(\algebra{R}))\] is the localization at \[\Sigma^i\algebra{R}[\underline{S}_{(\kappa)}]\to\Sigma^i\algebra{R}[\underline{S}_{(\kappa)}]^{L\solid\algebra{R}}, S\in\Pro(\Fin)_{(\kappa)}, i\in\mathbb Z, \] i.e.\ the essential image of $f$ consists of those $M\in\mathcal D(\Sol{(\kappa)}(\algebra{R}))$ such that \[\Map_{\mathcal D(\Cond{(\kappa)}(\Ab))}(\Sigma^i\algebra{R}[\underline{S}]^{L\solid\algebra{R}},M)\to\Map_{\mathcal D(\Cond{(\kappa)}(\Ab))}(\Sigma^i\algebra{R} [\underline{S}_{(\kappa)}],M)\] is an equivalence for all $S\in\Pro(\Fin)_{(\kappa)}, i\in\mathbb Z$. This shows that $ \ref{essentialimagesolidrmodules}\Leftrightarrow \ref{solidderivedrmodules}$. 
    \end{proof}

\subsection{Cohomology with solid coefficients}\label{section:condensedcohomologywithsolidcoefficients}
    We now apply the above computations to compare $(\kappa)$-condensed cohomology with coefficients in solid abelian groups with sheaf and singular cohomology. 
\begin{notation}
    Since $\mathcal D(\Cond{(\kappa)}(\Ab))\cong \LMod{H\mathbb Z}{\Cond{(\kappa)}(\Sp)}$ (\cref{condensedmodulesaremodulesincondensed}), condensed  cohomology enhances to a functor \[\cckH{(\kappa\text{-})}(-,-)\colon \Cond{\kappa}(\an)^{\operatorname{op}}\times \mathcal D(\Cond{\kappa}(\Ab))\to \mathcal D(\Ab), \] cf.\ \cref{cohomologyisext1}.
    In this section, we denote by \[\cckH{(\kappa\text{-})}(-,-)\colon \Cond{(\kappa)}(\an)^{\operatorname{op}}\times \mathcal D(\Sol{(\kappa)})\to \mathcal D(\Ab)\] the restriction of this functor to $\Cond{(\kappa)}(\an)^{\operatorname{op}}\times \mathcal D(\Sol{(\kappa)})\subseteq \Cond{(\kappa)}(\an)^{\operatorname{op}}\times\mathcal D(\Cond{(\kappa)}(\Ab))$. 
\end{notation}
 
    \cref{solidcw} implies that cohomology with solid coefficients is homotopy invariant:
    \begin{cor}\label{solidcohomologyhomotopyinvariant}
        Suppose that $h\colon X\to Y$ is a homotopy equivalence. 
        \begin{romanenum}
        \item For all uncountable cardinals $\kappa$, $h$ induces an equivalence  
        \[ \mathbb Z[\underline{X}_{\kappa}]^{L\solid}\cong\mathbb Z[\underline{Y}_{\kappa}]^{L\solid}\in\mathcal D(\Sol{\kappa}).\] 
        In particular, for all $A\in \mathcal D(\Sol{\kappa})$ 
        \[h^*\colon \ckH(\underline{Y},A)\cong \ckH(\underline{X},A)\] is an equivalence.  

        \item If $X$ and $Y$ are $\To$, then $h$ induces an equivalence \[\mathbb Z[\underline{X}]^{L\solid}\cong\mathbb Z[\underline{Y}]^{L\solid}\in\mathcal D(\Sol{}).\] In particular, for all $A\in\mathcal D(\Sol{})$,
        \[h^*\colon \ccH(\underline{Y},A)\cong \ccH(\underline{X},A)\] is an equivalence.
        \end{romanenum}
    \end{cor}
    \begin{proof}
        As $(-)^{L\solid}$ is symmetric monoidal and $\mathbb Z[\underline{[0,1]}_{(\kappa)}]^{L\solid}\cong\mathbb Z$ (\cref{solidcw}), for a ($\kappa$-)condensed set $X$, \[ \mathbb Z[X\times\underline{[0,1]}_{(\kappa)}]^{L\solid}=\mathbb Z[X]^{L\solid}\otimes_{\mathcal D(\Sol{(\kappa)})}\mathbb Z\cong \mathbb Z[X]^{L\solid}.\] This implies that $\mathbb Z[-]^{L\solid}\circ \underline{(-)}_{(\kappa)}$ inverts homotopy equivalences. 
        By \cref{cohomologyisext1} and \cref{derivedsolidissolidderived}, \cref{freeabeliangroupsacyclicforsolidificationwithoutkappa}, for $A\in\Sol{(\kappa)}$,
        \[ \cckH{(\kappa\text{-})}(-,A)|_{\tau_{\leq 0}\Cond{(\kappa)}(\an)^{\operatorname{op}}}\cong \map_{\mathcal D(\Cond{(\kappa)}(\Ab))}(\mathbb Z[-],A)\cong \map_{\mathcal D(\Cond{(\kappa)}(\Ab))}(\mathbb Z[-]^{L\solid},A).\] This implies the statement on cohomology.  
    \end{proof}
    
    In the appendix, we show that gros topos cohomology with coefficients in constant sheaves inverts homotopy equivalences. Together with \cref{solidcohomologyhomotopyinvariant}, this implies the following:  
        \begin{cor}\label{deltakappaexactnesshomotopyinvariant}
    Suppose that $A$ is a discrete abelian group.
    \begin{romanenum}
    \item For a topological space $X$, the comparison map
    \[ \cH{\sheaf}(X,A)\to \ckH(\underline{X}_{\kappa},\underline{A}_{\kappa}) \, \,(\ref{comparisonmapkappacondensedsheaftopological})\] is an equivalence if and only if $X$ is homotopy equivalent to a topological space $Y$ such that the comparison map 
    \[ \cH{\sheaf}(Y,A)\to \ckH(\underline{Y}_{\kappa},\underline{A}_{\kappa}) \,\, (\ref{comparisonmapkappacondensedsheaftopological})\] is an equivalence. 

    \item For a $\To$ topological space $X$, the comparison map
    \[\cH{\sheaf}(X,A)\to \ccH(\underline{X},\underline{A})\,  (\ref{comparisonmapkappacondensedsheaftopologicalwithoutkappa})\] is an equivalence if and only if $X$ is homotopy equivalent to a $\To$ topological space $Y$ such that the comparison map 
    \[ \cH{\sheaf}(Y,A)\to \ccH(\underline{Y},\underline{A}) \,\, \,(\ref{comparisonmapkappacondensedsheaftopologicalwithoutkappa})\] is an equivalence. 
    \end{romanenum}    
\end{cor}

    \begin{proof}
        It suffices to prove the $\kappa$-condensed statement, then the condensed statement  follows from \cref{condensedcohomologycanbecomputedonfinitestage}. 
        Fix $\lambda\geq \kappa>|X|+|Y|$.  
        By construction of the comparison map, for a $\lambda$-small topological space $T$, the map  
        \[ \cH{\sheaf}(T,A)\to \cH{\sheaf}(T,\mathcal C_{\kappa}(-,A))\to \ckH(\underline{T}_{\kappa},\underline{A}_{\kappa})\] is an equivalence if and only if 
        $\mathcal C(-,A)\to \mathcal C_{k}(-,A)\to \mathcal C_{\kappa}(-,A)=\pi_{0}\stradg{j}\underline{A}_{\kappa}\to \stradg{j}\underline{A}_{\kappa}$ induces an equivalence 
        \[ \cH{\grostop}(h_T,A)\cong \cH{\grostop}(h_T,\stradg{j}\underline{A}_{\kappa}).\]
        As $\underline{A}_{\kappa}$ is solid (\cref{profinitesolid}), \[\cH{\grostop}(h_{-},\stradg{j}\underline{A}_{\kappa})\cong \ckH(\underline{(-)}_{\kappa},\underline{A}_{\kappa})\colon (\Top^{\lambda})^{\operatorname{op}}\to \Ab\] inverts homotopy equivalences by \cref{solidcohomologyhomotopyinvariant}. 

        By \cref{grostoposcohomologyomotopyinvariant},
        \[ \cH{\grostop}(h_{-},h_A)\colon (\Top^{\lambda})^{\operatorname{op}}\to \Sp\] inverts homotopy equivalences as well.
    \end{proof}
    \cref{homotopyinvariancesheafcohomologyconstantcoefficients,solidcohomologyhomotopyinvariant} and \cref{comparisonmaponbasis} imply the following: 
    \begin{lemma}\label{condensedandsheafcohomologylocallycontractible}
        Suppose $X$ is a topological space which is homotopy equivalent to a locally contractible topological space. 
    \begin{romanenum}
        \item For all uncountable cardinals $\kappa$ and all discrete abelian groups $A$, 
        \[\cH{\sheaf}(X,A)\cong \ckH(\underline{X}_{\kappa},\underline{A}_{\kappa})\] via the comparison map (\ref{comparisonmapkappacondensedsheaftopological}). 
        \item For a solid abelian groups $A\in\Sol{\kappa}$, denote by $A(*)^{\delta}$ the discrete abelian group $A(*)$. 
        The counit of the global sections geometric morphism defines a map $\underline{A(*)^{\delta}}_{\kappa}\to A$ which induces an isomorphism 
        \[ \ckH(\underline{X}_{\kappa},\underline{A(*)^{\delta}})\cong  \ckH^*(\underline{X}_{\kappa},A).\]  
        \item If $X$ is $\To$, for all solid abelian groups $A\in\Sol{}$, \[ \cH{\sheaf}^*(X,A(*)^{\delta})\cong \ccH(\underline{X},\underline{A(*)^{\delta}})\cong \ccH(\underline{X},A),\] 
        where the first equivalence is the comparison map (\ref{comparisonmapkappacondensedsheaftopologicalwithoutkappa}) and the second map is the counit $\underline{A(*)^{\delta}}\to A$.  
    \end{romanenum}
    \end{lemma}
    \begin{proof}Recall from \cref{colimitdescriprionsolid} that a condensed abelian group $M$ is solid if and only if $r^kM\in\Sol{\kappa}$ for all regular cardinals $\kappa$. By \cref{condensedcohomologycanbecomputedonfinitestage}, it therefore suffices to prove the $\kappa$-condensed statements.  

    As for a discrete abelian group $M$, $\underline{M}_{\kappa}\in\Sol{\kappa}$ (\cref{profinitesolid}), by \cref{deltakappaexactnesshomotopyinvariant,solidcohomologyhomotopyinvariant}, we can assume that $X$ is locally contractible. 
    By \cref{comparisonmaponbasis}, it is enough to prove the statement for $X$ contractible, and using \cref{deltakappaexactnesshomotopyinvariant,solidcohomologyhomotopyinvariant}, we can reduce to the case $X=*$. 

    For $X=*$, $\Shv(X)\cong\an$, and under this identification, the geometric morphism 
        \[\ladg{i}\colon \an\cong \Shv(X)\rightleftarrows \Cond{\kappa}(\an)\colon \radg{i}\] is the global sections geometric morphism. 
        The global sections functor stabilizes to a $t$-exact functor (\cref{globalsectionstexactcocontinuous}), whence $\stradg{i}B=\pi_0\stradg{i}B=\mathcal C(-,B(*)^{\delta})$ for all $B\in\Cond{\kappa}(\Ab)$. 
        In particular, for all discrete abelian groups $B$, \[ \cH{\sheaf}(*,B)\cong \ckH(*,\underline{B}_{\kappa})\] via the comparison map.
        By \cref{cohomologyglobalsections}, $\ckH(*,-)$ is the stabilization of the global sections functor.  
        As the global sections functor $\Gamma$ is cocontinuous (\cref{globalsectionstexactcocontinuous}), the unit $\id_{\an}\to \Gamma\circ c$ is an equivalence. It follows from \cref{geometricmorphismstabilization} that the unit $\id_{\Sp}\to \Gamma_{\Sp}\circ c_{\Sp}$ is an equivalence as well, whence the counit $\underline{A(*)^{\delta}}_{\kappa}\to A$ induces an equivalence \[\ckH(*,\underline{A(*)^{\delta}}_{\kappa})\cong \ckH(*,A).\qedhere\] 
    \end{proof}

    \begin{rem}For a locally contractible space $X$ and a discrete abelian group $A$, \[\cH{\sheaf}(X,A)\cong \cH{\sing}(X,A)\] by e.g.\ \cite{sella2016comparison}. 
    \end{rem}

  \cref{deltakappaexactnesshomotopyinvariant,condensedandsheafcohomology1} imply the following: 
    \begin{cor}\label{condensedandsheafcohomologylocallycompacthausdorffhomotopy}
        Suppose $X$ is a topological space which is homotopy equivalent to a locally compact Hausdorff space.
    \begin{romanenum}
        \item  For all discrete abelian groups $A$, \[ \cH{\sheaf}(X,A)\cong \ckH(\underline{X}_{\kappa},\underline{A}_{\kappa})\] via the comparison map (\ref{comparisonmapkappacondensedsheaftopological}).  
        \item If $X$ is a $\To$ topological space, for all discrete abelian groups $A$, 
        \[\cH{\sheaf}(X,A)\cong \ccH(\underline{X},\underline{A})\] via the comparison map (\ref{comparisonmapkappacondensedsheaftopologicalwithoutkappa}). 
    \end{romanenum}
    \end{cor}

     Denote by $c_{\Sp}\colon \Sp\to\Cond{(\kappa)}(\Sp)$ the stabilization of the constant sheaf functor. 
    \cref{internalcohomologydiscretecompacthausdorff} implies that the condensed structure on internal cohomology is constant in the following situation: 
    \begin{cor}\label{internalcohomologydiscrete}
        Suppose that $X$ is a ($\kappa$-)condensed set and there exists a ($\kappa$-light) compact Hausdorff space $K$ such that $\mathbb Z[X]^{L\solid}\cong \mathbb Z[\underline{K}_{(\kappa)}]^{L\solid}$. 

        For all discrete abelian groups $A$, the counit
        \[c_{\Sp}\left(\cckH{(\kappa\text{-})}(X, \underline{A}_{(\kappa)})\right)\to\icckH{(\kappa\text{-})}(X, \underline{A}_{(\kappa)})\] is an equivalence.
    \end{cor}    
    \begin{ex}By \cref{solidcohomologyhomotopyinvariant}, \cref{internalcohomologydiscrete} in particular applies if $X=\underline{Y}_{(\kappa)}$ for a ($\To$) topological space $Y$ which is homotopy equivalent to a ($\kappa$-light) compact Hausdorff space. \end{ex}
    \begin{proof}
        Denote by $\imap_{\mathcal D(\Cond{(\kappa)}(\Ab))}(-,-)$ the internal Hom, by \[f\colon \mathcal D(\Cond{\kappa}(\Ab))\cong \LMod{H\mathbb Z}{\Cond{\kappa}(\Sp)}\leftrightarrows\Cond{\kappa}(\Sp)\colon H\mathbb Z[-]\] the free-forget adjunction and let \[\mathbb Z[-]\coloneqq H\mathbb Z[-]\circ \Sigma^{\infty}_{+}\colon \Cond{(\kappa)}(\an)\to \mathcal D(\Cond{(\kappa)}(\Ab)).\] By \cref{freemoduloediscrete}, this restricts to a functor $\Cond{(\kappa)}(\Set)\to \Cond{(\kappa)}(\Ab)$, and this is left adjoint to the forget functor. 
        The symmetric monoidal structure of the free module functor $H\mathbb Z[-]$ yields an equivalence 
        \[\icckH{(\kappa\text{-})}(-,f-)\cong f\circ \imap_{\mathcal D(\Cond{(\kappa)}(\Ab))}(\mathbb Z[-],-).\] Since discrete abelian groups are solid (\cref{profinitesolid}), \cref{derivedsolidificationunderlyingforacyclicrings} implies that \[ \icckH{(\kappa\text{-})}(X, \underline{A}_{(\kappa)})\cong f\imap_{\mathcal D(\Cond{(\kappa)}(\Ab))}(\mathbb Z[X], \underline{A}_{(\kappa)})\cong f\imap_{\mathcal D(\Cond{(\kappa)}(\Ab))}(\mathbb Z[X]^{L\solid}, \underline{A}_{(\kappa)}).\] 
        We can therefore assume that $X=\underline{K}_{(\kappa)}$ is represented by a ($\kappa$-light) compact Hausdorff space $K$. 
        The statement now follows from \cref{internalcohomologydiscretecompacthausdorff}.
    \end{proof}
    \begin{rem}
        As the forget functor $f\colon\mathcal D(\Cond{(\kappa)}(\Ab))\to\Cond{(\kappa)}(\Sp)$ is conservative, \cref{internalhommodules} implies that for $X$ as in the above corollary, \[\cshv{\map_{\mathcal D(\Cond{(\kappa)}(\Ab))}(\mathbb Z[X], \underline{A}_{(\kappa)})}\cong \imap_{\mathcal D(\Cond{(\kappa)}(\Ab))}(\mathbb Z[X], \underline{A}_{\kappa})\] via the counit for the adjoint pair 
        \[\cshv{-}\colon \mathcal D(\Ab)\rightleftarrows \mathcal D(\Cond{(\kappa)}(\Ab))\colon\Gamma .\] 
    \end{rem}
\newpage
\section{Group Cohomology}
The condensed formalism yields a notion of group cohomology for topological groups which we study in this chapter. Our main focus is on its relation to continuous group cohomology and the condensed/singular/sheaf cohomology of classifying spaces. 

First, we review the category of $G$-objects for a group object $G$ in a big topos. This category is itself a big topos, and we define group cohomology as the cohomology of its terminal object and record some of its basic properties. Under mild conditions on a big topos $\topo{X}$, there is a \textit{derived fixed-point functor} \[(-)^G\colon \LMod{\mathbb S[G]}{\stab{\topo{X}}}\to \stab{\topo{X}}\] refining group cohomology (\cref{internalgroupcohomologytotalisaton}), which we also discuss very briefly. 
Then we specialize to the condensed setting. 
For a (Hausdorff) topological group $G$, the \v{C}ech-to-cohomology spectral sequence (\cref{Bousfieldkanspectralsequencehomology}) for the cover $\underline{G}_{(\kappa)}\to *$ yields a natural comparison map \[ \contgrpcoh{}^*(G,-)\to \condgrpcoh{\kappa\text{-}}^*(\underline{G}_{(\kappa)}, \underline{(-)}_{(\kappa)})\] from continuous to $(\kappa\text{-})$condensed group cohomology of the associated ($\kappa$-)condensed group. 
Our comparison of sheaf and condensed cohomology implies that this map is an isomorphism in two important cases: For locally profinite groups $G$ and solid continuous $G$-modules as coefficients (\cref{condensedequalscontinuousprofinite}), as well as for locally compact groups $G$ and finite-dimensional, continuous real $G$-representations as coefficients (\cref{continuousequalscondensedonvectorspacecoefficients}).  
However, we will see that condensed group cohomology is in general a much more refined invariant than continuous group cohomology. 
In \cref{section:condensedcohomologyandcohomologyofclassifyingspaces}, we deduce from homotopy invariance of condensed cohomology with solid coefficients, that for a (Hausdorff) topological group $G$, condensed group cohomology with coefficients in a solid abelian group with trivial $\underline{G}_{(\kappa)}$-action is isomorphic to the condensed cohomology of a classifying space $BG$ of numerable principal $G$-bundles. 
For many topological groups $G$, the condensed cohomology of a classifying space $BG$ is isomorphic to its singular and sheaf cohomology (\cref{section:Cohomologyclassifyingspaces}). 
This leads to plenty examples where condensed group cohomology differs from continuous group cohomology (\cref{condensedequalscontinuousimpliescontractible}). 
Finally, we show that for a large class of topological groups, continuous group cohomology with solid coefficients can be identified with $\Ext^*_{\Sol{}(\mathbb Z[\underline{G}_{(\kappa)}])}(\mathbb Z, \underline{(-)}_{(\kappa)})$. 
This relies on the existence of enough projectives in $\Sol{(\kappa)}(\mathbb Z[\underline{G}_{(\kappa)}])$ (\cref{solidenoughprojectiveskappa}, \cref{solidenoughprojectiveswithoutkappa}) and the computations from \cref{section:computationssolidification}. 
\subsection{\texorpdfstring{$G$}{G}-objects in a big topos}

\begin{definition}\label{definitiongobjectsinbigtopos}
    Suppose $\topo{X}$ is a big topos and equip it with the cartesian monoidal structure.  
    A group object in $\topo{X}$ is an algebra $G\in \Alg(\topo{X}^{\times})$ such that $\tau_{\leq 0}G\in\Alg(\tau_{\leq 0}\topo{X}^{\times})$ is a group object, i.e.\ the shear maps 
    \begin{align*} \tau_{\leq 0}G\times\tau_{\leq 0}G&\to \tau_{\leq 0}G\times\tau_{\leq 0}G\\ (g,h)&\mapsto (gh,h)\\ (g,h)&\mapsto (hg,h)\end{align*} are isomorphisms. 
    Denote by $\Grp(\topo{X})\subseteq \Alg(\topo{X}^{\times})$ the full subcategory on group objects.
    For a group object $G\in \Grp(\topo{X})$ denote by $\Act{G}{\topo{X}}\coloneqq \LMod{G}{\mathcal X}$ the category of left $G$-modules in $\mathcal X$.
\end{definition}In this section, we record basic properties of the category of $G$-objects in a big topos. We first show that this category is itself a big topos. Then we recall from \cite{Principalinftybundles} that the \textit{homotopy orbit functor }\[-//G\colon \Act{G}{\topo{X}}\to\topo{X}\] (\cref{deftrivfixedpoints}) classifies principal $\infty$-$G$-bundles (\cref{definitionprincipalinfinitybundle}) in $\topo{X}$. This will be central for our comparison of condensed group cohomology with cohomology of classifying spaces. To identify $G$-objects with principal $\infty$-$G$-bundles, we review simplicial models for $G$-objects. This also yields handy description of restriction and extension of scalars functors, which we apply to identify stable module categories in $\Act{G}{\topo{X}}$ with modules over the corresponding group rings (\cref{modulesingroupobjectsaremodulesovergrouprings}).
\begin{lemma}If $G$ is a group object in a topos $\topo{X}$, then $\Act{G}{\topo{X}}$ is a topos.
    \end{lemma} 
    \begin{proof}By (\cite[Proposition 6.1.0.6]{highertopostheory}, $\topo{X}^{\times}$ is presentably symmetric monoidal. (Cocontinuity of the cartesian monoidal structure follows from universality of colimits). By \cite[Corollary 4.2.3.7]{higheralgebra}, this implies that $\Act{G}{\topo{X}}$ is presentable and that the forget functor $\Act{G}{\topo{X}}\to\topo{X}$ preserves small limits and colimits. 
    Using that the forget functor is conservative, it is now  straightforward to check that $\Act{G}{\topo{X}}$ satisfies Giraud's axioms (\cite[Proposition 6.1.0.6]{highertopostheory}). 
    \end{proof}
    
\begin{lemma}\label{adjunctioninducesadjunctionongobjects}Suppose that $R\colon \mathcal X\leftrightarrows \mathcal Y \colon L$ is a geometric morphism (\cref{definitiongeometricmorphism}) between big topoi.
    \begin{romanenum}
        \item $L$ and $R$ lift to an adjoint pair \[L^{\Alg}\colon \Alg(\topo{X}^{\times})\rightleftarrows \Alg(\topo{Y}^{\times})\colon R^{\Alg}, \] and $L^{\Alg}$ preserves group objects. 
        \item If $G\in \operatorname{Grp}(\mathcal X)$ is a group object, then $L\dashv R$ lift to an adjoint pair \[L^G\colon \Act{G}{\topo{X}}\rightleftarrows \Act{L^{\Alg}G}{\topo{Y}}\colon R^G\] and $L^G$ preserves finite limits. 
        If $L$ is fully faithful, so is $L^G$. 
    \end{romanenum}
            \end{lemma}
            \begin{proof}
                Since $L$ and $R$ preserve finite products, they promote uniquely to symmetric monoidal functors (\cite[Corollary 2.4.1.9]{higheralgebra}) and in particular lift to an adjoint pair \[L^{\Alg}\colon\Alg(\topo{X}^{\times})\rightleftarrows \Alg(\topo{Y}^{\times})\colon R^{\Alg}, \] cf.\ \cref{adjuncttioninducesadjunctiononalgebraobjects}.  Next, we show that $L^{\Alg}$ preserves group objects.  
                By \cite[Proposition 5.5.6.16]{highertopostheory}, $L$ and $R$ restrict to functors $l\colon \tau_{\leq 0}\topo{X}\leftrightarrows \tau_{\leq 0}\topo{Y}\colon r$, and $\tau_{\leq 0}L\cong l\tau_{\leq 0}$ since both are left adjoint to \[(\tau_{\leq 0}\topo{X}\hookrightarrow \topo{X})\circ r.\] The functor $l$ preserves products and therefore induces a functor $\Alg(\tau_{\leq 0}\topo{X}^{\times})\to\Alg(\tau_{\leq 0}\topo{Y}^{\times})$ which preserves group objects. This shows that $L^{\Alg}$ preserves group objects. 

                As explained in \cref{localisationinduceslocalisationonmodulecategories}, $L$ and $R$ lift to an adjoint pair \[L^G\colon \Act{G}{\topo{X}}\rightleftarrows \Act{L^{\Alg}G}{\topo{Y}}\colon R^G,\] and if $L$ is fully faithful, so is $L^G$. 
                Since the forget functors $\Act{G}{\topo{X}}\to\topo{X}$ and $\Act{L^{\Alg}G}{\topo{Y}}\to\topo{Y}$ create limits (\cref{forgetfreeadjunctionmodules}) and $L$ preserves finite limits, it follows that $L^G$ preserves finite limits. 
            \end{proof}
            Together with \cref{filteredcolimitsmodules}, this implies: 
\begin{cor}\label{Gobjectsisbigtopos}Suppose $\topo{X}$ is a big topos and $\mathcal X_*\colon \Lambda\to \Pr^L$ is an exhaustion of $\topo{X}$ by topoi. 
    \begin{romanenum}
\item The equivalence $\colim{\lambda\in\Lambda}\Alg(\topo{X}^{\times}_{\lambda})\cong \Alg(\topo{X}^{\times})$ restricts to an equivalence \[\colim{\lambda\in\Lambda}\Grp(\topo{X}_{\lambda})\cong \Grp(\topo{X}).\] 
\item For $G\in\Grp(\topo{X})$, $\Act{G}{\topo{X}}$ is a big topos. If $\mu\in\Lambda$ is such that $G\in \Grp(\mathcal X_{\mu})\subseteq \Grp(\mathcal X)$, then \begin{align*} \phantom{}_{\mu\backslash}{\Lambda}&\to \Pr^L\\\mu\to\lambda&\mapsto \lAct{G}{\topo{X}}{\lambda}\end{align*} is an exhaustion for $\Act{G}{\topo{X}}$ by topoi. 
\item In particular, for all $\kappa\to\lambda\in \phantom{}_{\mu\backslash}\Lambda$, \[\lAct{G}{\topo{X}}{\kappa}\to \lAct{G}{\topo{X}}{\lambda}\to \Act{G}{\topo{X}}\] are fully faithful, left-exact left adjoints.
    \end{romanenum}
\end{cor}
\begin{proof}By \cref{adjunctioninducesadjunctionongobjects}, the equivalence $\colim{\lambda\in\Lambda}\Alg(\topo{X}^{\times}_{\lambda})\cong \Alg(\topo{X}^{\times})$ restricts to a fully faithful functor \[\colim{\lambda\in\Lambda}\Grp(\topo{X}_{\lambda})\to \Grp(\topo{X}).\] 
Since $\tau_{\leq 0}^{\topo{X}}$ restricts to $\tau_{\leq 0}^{\topo{X}_{\lambda}}$ for all $\lambda\in\Lambda$ (\cref{truncationbigtopos}), this is essentially surjective. 
The remaining statements hold by \cref{filteredcolimitsmodules} and \cref{fullyfaithfulnessandpreservationoflimits}.  
\end{proof}

\cref{existenceleftadjointpullback} implies: 
\begin{cor}\label{existenceinduction}
    Suppose that $\topo{X}$ is a big topos and 
    $\phi\colon H\to G\in\Grp(\topo{X})$ is a group homomorphism. 
   Restriction of scalars
    $\phi^*\colon \Act{G}{\topo{X}}\to \Act{H}{\topo{X}}$ along $\phi$ admits a left adjoint $\phi_{!}\dashv \phi^*$. 
\end{cor}
\begin{rem}
If $\topo{X}$ is presentable, i.e.\ a topos, then $\phi^*$ also has a right adjoint $\phi_*$ by \cite[Corollary 5.5.2.9]{highertopostheory}.
Being a right adjoint, $\phi^*$ preserves finite products and hence enhances essentially uniquely to a symmetric monoidal functor 
\[ (\phi^*)^{\times}\colon (\Act{G}{\topo{X}})^{\times}\to (\Act{H}{\topo{X}})^{\times}\] by \cite[Corollary 2.4.1.9]{higheralgebra}. 
\end{rem}
\cref{universalpropertyspectrumobjectsbigpresentablemonoidal} implies the following:  
\begin{cor}\label{monoidalstructurestabilizationpullback}
    Suppose that $\topo{X}$ is a big topos and $\phi\colon H\to G\in\Grp(\topo{X})$ is a group homomorphism. 
    Denote by $\phi^*\colon \Act{G}{\topo{X}}\to\Act{H}{\topo{X}}$ restriction of scalars along $\phi$ and by $\Sp(\phi^*)\colon\stab{\Act{G}{\topo{X}}}\to \stab{\Act{H}{\topo{X}}}$ its stabilization.  
    There is an essentially unique symmetric monoidal enhancement $\Sp(\phi^*)^{\otimes}$ of $\Sp(\phi^*)$ with \[\Sp(\phi^*)^{\otimes}\circ \Sigma^{\infty,\otimes}_{+}\cong \Sigma^{\infty,\otimes}_{+}\circ (\phi^*)^{\times}.\] 
\end{cor}
\begin{proof}
    Choose an exhaustion $\topo{X}_*\colon\Lambda\to \Pr^L$ of $\topo{X}$ by topoi such that $G,H\in\Grp(\topo{X}_{\lambda})$ for all $\lambda\in\Lambda$. For all $\lambda\in\Lambda$, $\phi^*|_{\topo{X}_{\lambda}}$ factors as \[ \lAct{G}{\topo{X}}{\lambda}\xrightarrow{\phi^*_{\lambda}}\Act{H}{\topo{X}_{\lambda}}\subseteq {\Act{H}{\topo{X}}}, \] where $\phi^*_{\lambda}$ denotes restriction of scalars along $H\to G\in\Grp(\topo{X}_{\lambda})$, and by universality of colimits in $\topo{X}_{\lambda}$ and \cref{forgetfreeadjunctionmodules}, $\phi^*_{\lambda}$ preserves small colimits. 

    Denote by \[\Fun^{\operatorname{colim},\Lambda,\otimes}(\stab{\Act{G}{\topo{X}}}, \stab{\Act{H}{\topo{X}}})\subseteq \Fun^{\otimes}(\stab{\Act{G}{\topo{X}}}^{\otimes}, \stab{\Act{H}{\topo{X}}}^{\otimes})\] the full subcategory on symmetric monoidal functors $F^{\otimes}\colon \stab{\Act{G}{\topo{X}}}^{\otimes}\to \stab{\Act{H}{\topo{X}}}^{\otimes}$ such that for all $\lambda\in\Lambda$, \[F|_{\stab{\lAct{G}{\topo{X}}{\lambda}}}\colon \stab{\lAct{G}{\topo{X}}{\lambda}}\to \stab{\Act{H}{\topo{X}}}\] preserves small colimits. 
    By \cref{Gobjectsisbigtopos}, $(\lAct{G}{\topo{X}}{*})^{\times}$ is an exhaustion of $(\Act{G}{\topo{X}})^{\times}$ by presentably symmetric monoidal categories. 
    By construction of the symmetric monoidal enhancement \[\Sigma^{\infty,\otimes}_{+}\colon \Act{G}{\topo{X}}\to \stab{\Act{G}{\topo{X}}}^{\otimes}\] (see the proof of \cref{symmetricmonoidalstructureonspectrumobjects}), this implies that 
    \begin{align*}&\Fun^{\operatorname{colim}, \Lambda, \otimes}(\stab{\Act{G}{\topo{X}}}^{\otimes}, \stab{\Act{H}{\topo{X}}}^{\otimes})\times_{\Fun^{\otimes}((\Act{G}{\topo{X}})^{\times}, \stab{\Act{H}{\topo{X}}}^{\otimes})}\{ \Sigma^{\infty,\otimes}_{+,H}\circ \phi^*\}\\ &\cong \clim{\lambda\in\Lambda}(\Fun^{\operatorname{colim}, \otimes}(\stab{\lAct{G}{\topo{X}}{\lambda}}^{\otimes}, \stab{\Act{H}{\topo{X}}}^{\otimes})\times_{\Fun^{\otimes}((\lAct{G}{\topo{X}}{\lambda})^{\times}, \stab{\Act{H}{\topo{X}}}^{\otimes})}\{\Sigma^{\infty,\otimes}_{+,H}\circ \phi^*\}.\end{align*}
    As for all $\lambda\in\Lambda$, $\phi^*|_{\lAct{G}{\topo{X}}{\lambda}}$ factors as \[\lAct{G}{\topo{X}}{\lambda}\xrightarrow{\phi^*_{\lambda}}\Act{H}{\topo{X}_{\lambda}}\subseteq \Act{H}{\topo{X}}\] and in particular preserves small colimits (\cref{forgetfreeadjunctionmodules}), the right-hand side is contractible by \cref{universalpropertyspectrumobjectsbigpresentablemonoidal}, i.e.\ there exists an essentially unique filler \begin{center}\begin{tikzcd}
       (\Act{G}{\topo{X}})^{\times}\arrow[d,"(\phi^*)^{\times}"'] \arrow[r,"\Sigma^{\infty, \otimes}_{+,G}"] & \stab{\Act{G}{\topo{X}}}^{\otimes}\arrow[d,dashed, "(\phi_{\Sp})^{\otimes}", "\exists!"']\\ 
        {\Act{H}{\topo{X}}}^{\times}\arrow[r,"\Sigma^{\infty, \otimes}_{+,H}"] & \stab{\Act{H}{\topo{X}}}^{\otimes}.
        \end{tikzcd}
        \end{center} such that for all $\lambda\in\Lambda$, $(\phi_{\Sp})|_{\stab{\lAct{G}{\topo{X}}{\lambda}}}$ preserves small colimits. 
        It remains to show that this is a symmetric monoidal enhancement of $\Sp(\phi^*)$, i.e. that there is an equivalence $\phi_{\Sp}\cong \Sp(\phi^*)$. 
        By \cref{Gobjectsisbigtopos,stabilizationbigpresentablecategoriesbigpresentable}, \begin{align*}&\Fun^{\operatorname{colim},\Lambda}(\stab{\Act{G}{\topo{X}}}, \stab{\Act{H}{\topo{X}}})\times_{\Fun^{\operatorname{colim},\Lambda}(\Act{G}{\topo{X}}, \Act{H}{\topo{X}})}\{  \Sigma^{\infty}_{+,H}\circ \phi^*\}\\ &\cong \clim{\lambda\in\Lambda}(\Fun^{\operatorname{colim}}(\stab{\lAct{G}{\topo{X}}{\lambda}}, \stab{\Act{H}{\topo{X}}})\times_{\Fun^{\operatorname{colim}}(\lAct{G}{\topo{X}}{\lambda}, \stab{\Act{H}{\topo{X}}})}\{ \Sigma^{\infty}_{+,H}\circ \phi^*\}),\end{align*} and the right-hand side is contractible by \cref{universalpropertyspectrumobjectsbigpresentable}. 
        It therefore suffices to show that $\stab{\phi^*}$ enhances to an element of the left-hand side. 
        By construction, for all $\lambda\in\Lambda$, $\stab{\phi^*}|_{\stab{\lAct{G}{\topo{X}}{\lambda}}}$ factors as $\stab{\lAct{G}{\topo{X}}{\lambda}}\xrightarrow{\stab{\phi^*_{\lambda}}}\stab{\Act{H}{\topo{X_{\lambda}}}}\subseteq \stab{\Act{H}{\topo{X}}}$ and in particular preserves small colimits. It remains to describe an equivalence $\stab{\phi^*}\circ \Sigma^{\infty}_{+,G}\cong \Sigma^{\infty}_{+,H}\circ \phi^*$. 

        As for all $\lambda\in\Lambda$, $\phi^{*}_{\lambda}$ and $\lAct{G}{\topo{X}}{\lambda}\subseteq \Act{G}{\topo{X}}$ are left-exact left adjoints, the mate \[ \Sigma^{\infty}_{+,H}\circ \phi^*\circ (\lAct{G}{\topo{X}}{\lambda}\subseteq \Act{G}{\topo{X}}) \to \stab{\phi^*}\circ \Sigma^{\infty}_{+,G}\circ (\lAct{G}{\topo{X}}{\lambda}\subseteq \Act{G}{\topo{X}})\] of the outer square in the commutative diagram 
        \begin{center}
 \begin{tikzcd}
    \lAct{G}{\topo{X}}{\lambda}\arrow[r, hookrightarrow ]&\Act{G}{\topo{X}}\arrow[r,"\phi^*"]& \Act{H}{\topo{X}}\\ 
    \stab{\lAct{G}{\topo{X}}{\lambda}}\arrow[r, hookrightarrow ]\arrow[u,"\Omega^{\infty}"]&\stab{\Act{G}{\topo{X}}}\arrow[r,"\stab{\phi^*}"]\arrow[u,"\Omega^{\infty}"] & \stab{\Act{H}{\topo{X}}}\arrow[u,"\Omega^{\infty}"']
    \end{tikzcd}
\end{center} is an equivalence. As the mate of the left square is also an equivalence (\cref{stabilizationcommuteswithsuspension}), it follows from the pasting law for mates (\cite[Lemma 2.2.4]{carmeli2022ambidexterity}) that the mate \[\beta\colon \Sigma^{\infty}_{+,H}\circ \phi^*\to \stab{\phi^*}\circ \Sigma^{\infty}_{+,G}\] of the right square is an equivalence on $\lAct{G}{\topo{X}}{\lambda}\subseteq \Act{G}{\topo{X}}$. As $\lambda$ was arbitrary, this implies that $\beta$ is an equivalence. 
\end{proof} 
\begin{definition}\label{deftrivfixedpoints}
    Suppose $\topo{X}$ is a big topos.
    If $i\colon H\subseteq G\in\Grp(\topo{X})$ is a monomorphism, we denote by $\res{H}{G}\colon \Act{G}{\topo{X}}\to\Act{H}{\topo{X}}$ restriction of scalars along $i$ and by $\ind{G}{H}$ its left adjoint.

    For a group object $G\in\topo{X}$, we denote by $\triv\colon \Act{G}{\topo{X}}\to \topo{X}$ restriction of scalars along the group homomorphism $G\to *$. The left adjoint of $\triv$ is called \emph{homotopy orbit functor} and is denoted $-//G$. 
    If it exists, the right adjoint of $\triv$ is called $G$-\textit{fixed-point functor} and is denoted $(-)^G$. 
\end{definition}
We will show below (\cref{fixedpointfunctorexistsforcartesianclosed}) that if $\topo{X}$ is a cartesian closed big topos, $G$-fixed-point functors exist for all groups $G\in\Grp(\topo{X})$. 

\paragraph{Simplicial models for G-objects}
We now recall from \cite[section 4.2.2]{higheralgebra} that $G$-modules in a big topos can be described as simplicial objects encoding a coherent $G$-action on an object $X\in\mathcal X$. 
\begin{recollection}Suppose $\topo{X}$ is a big topos.  
By \cite[Propositions 2.4.1.7, 4.1.2.10]{higheralgebra} and \cref{characterizationgroupobjects}, there is a fully faithul functor
\[\Grp(\topo{X})\hookrightarrow \Fun(\Delta^{\operatorname{op}}, \topo{X})\] which sends a group object $G$ to the simplicial object $[n]\mapsto G^n$ with \begin{align*}d_i\colon G^n&\to G^{n-1}, \\(g_1, \ldots,g_n)&\mapsto (g_1, \ldots, g_{i}g_{i+1},g_{i+2}, \ldots,g_n)\end{align*} and \begin{align*}s_i\colon G^{n}& \to G^{n+1}, \\ (g_1, \ldots,g_n)& \mapsto (g_1, \ldots, g_{i-1},1,g_i, \ldots,g_n).\end{align*} 
This is natural in the cartesian monoidal category $\topo{X}$. 
The essential image of this functor consists of those simplicial objects $G\colon \Delta^{\operatorname{op}}\to \topo{X}$ satisfying the following two conditions: 
\begin{romanenum}\item For all $n\in\mathbb N_0$, the face maps 
    \[\{ G([n])\to G(\{i-1,i\})\}_{1\leq i\leq n}\] exhibit $G([n])\cong \prod_{i=1}^n G(\{i-1,i\})$. 
    \item The map \[\tau_{\leq 0}G([1])\times \tau_{\leq 0}G([1])\cong \tau_{\leq 0}G(\{0,1\})\times \tau_{\leq 0}G(\{1,2\})\cong \tau_{\leq 0}G([2])\xrightarrow{d_1} \tau_{\leq 0}G([1])\] defines a group structure on $\tau_{\leq 0}G([1])$ in $\tau_{\leq 0}\topo{X}.$ 
\end{romanenum}
\end{recollection} 
\begin{lemma}\label{simplicialmodelgobjects}
    Suppose $\topo{X}$ is a big topos, $G\in\Grp(\topo{X})$ is a group object and denote by \[(*//G)_*\in\Fun(\Delta^{\operatorname{op}}, \topo{X})\] the corresponding simplicial object. 
    There is a fully faithful functor 
    \[\Act{G}{\topo{X}}\hookrightarrow \oc{\Fun(\Delta^{\operatorname{op}}, \topo{X})}{(*//G)_*}.\] 
    Its essential image consists of the simplicial objects $X\to (*//G)_*$ such that for all $n\in\mathbb N_0$, the maps \[X([n])\to X(\{n\})\cong X([0]),\, \text{ and }\, X([n])\to (*//G)([n])\cong G^n\] exhibit $X([n])$ as $X([n])\cong G^n\times X([0])$. 
    This is natural in the group object $G$ and with respect to finite products-preserving functors of big topoi. 
\end{lemma}
\begin{proof}This follows from the discussion in \cite[Chapter 4]{higheralgebra}, see \cref{modulecategoriesincartesian}.
\end{proof}
Note for future reference that the simplicial objects in the essential image of \[ \Act{G}{\topo{X}}\hookrightarrow \oc{\Fun(\Delta^{\operatorname{op}}, \topo{X})}{(*//G)_*}\to \Fun(\Delta^{\operatorname{op}}, \topo{X})\] satisfy \cite[Proposition 6.1.2.6, 4{{''}}]{highertopostheory}, i.e.\ are groupoid objects (\cite[Definition 6.1.2.7]{highertopostheory}). 
\begin{definition}\label{definitionclassifyingspaceofgbundles}
    Suppose that $\topo{X}$ is a big topos and $G\in\Grp(\topo{X})$ is a group object. 
    Choose an exhaustion $\topo{X}_*\colon\Lambda\to \Pr^L$ of $\topo{X}$ by topoi and $\mu\in \Lambda$ with $G\in\Grp(\topo{X}_{\mu})$. 
    Let \[-//G\colon \Act{G}{\topo{X}}\cong \colim{\lambda\in \phantom{}_{\mu\backslash}\Lambda}\lAct{G}{\topo{X}}{\lambda}\subseteq \colim{\lambda\in\phantom{}_{\mu\backslash}\Lambda}\oc{\Fun(\Delta^{\operatorname{op}}, \topo{X_{\lambda}})}{(*//G)_*}\to \colim{\lambda\in\phantom{}_{\mu\backslash\Lambda}}\Fun(\Delta^{\operatorname{op}}, \topo{X}_{\lambda})\xrightarrow{\colim{\Delta^{\operatorname{op}}}}\colim{\lambda\in\phantom{}_{\mu\backslash\Lambda}}\topo{X}_{\lambda}\cong \topo{X}\] and $\mathbb BG\coloneqq (-//G)(*)\coloneqq \colim{\Delta^{\operatorname{op}}}(*//G)_*$.
    
    Since for all $\lambda\in\Lambda$, $\topo{X}_{\lambda}\subseteq \topo{X}$ is closed under finite limits and small colimits and $\Lambda$ is filtered, this does not depend on the chosen exhaustion of $\topo{X}$ by topoi nor the choice of $\mu$.
    
    Denote by $-//G\colon \Act{G}{\topo{X}}\to \oc{\topo{X}}{\mathbb BG}$ the induced functor \[\Act{G}{\topo{X}}\cong \oc{\Act{G}{\topo{X}}}{*}\xrightarrow{-//G}\oc{\topo{X}}{\mathbb BG}.\] 
    We call $\mathbb BG$ the \textit{classifying space of principal} $\infty$-$G$-\textit{bundles} in $\topo{X}$. This terminology is motivated by \cref{characterizationbundles} below. 
\end{definition}
\begin{rem}
By \cref{basechange} below, the homotopy orbits functor (\cref{deftrivfixedpoints}) factors as \[\Act{G}{\topo{X}}\xrightarrow{-//G}\oc{\topo{X}}{\mathbb BG}\xrightarrow{(\mathbb BG\to *)_!} \topo{X},\] where $(BG\to *)_{!}$ denotes the forget functor. The  above notation therefore does not clash with \cref{deftrivfixedpoints}.  
\end{rem}
\begin{proposition}[{\cite[Proposition 3.2.76]{sati2022equivariantprincipalinfinitybundles}}]\label{gobjectsequalsgbundles}
    The functor \[-//G\colon  \Act{G}{\topo{X}}\to \oc{\topo{X}}{\mathbb BG}\] is an equivalence. 
\end{proposition}
\begin{proof}
    Choose an exhaustion $\topo{X}_*\colon\Lambda\to \Pr^L$ of $\topo{X}$ by topoi and $\mu\in\Lambda$ with $G\in\Grp(\topo{X}_{\mu})$. 
    Upon replacing $\Lambda$ by $_{\mu\backslash }\Lambda$ we can assume that $G\in\Grp(\topo{X}_{\lambda})$ for all $\lambda\in\Lambda$.
    Since $\mathcal{X}_{\lambda}\subseteq \mathcal{X}_{\infty}$ is closed under small colimits, for all $\lambda\in\Lambda$, $\mathbb BG\in\topo{X}_{\lambda}$, $\oc{{\topo{X}_{\lambda}}}{\mathbb BG}\subseteq \oc{\topo{X}}{\mathbb BG}$ is a full subcategory and 
    \[ \Act{G}{\topo{X}}\cong \oc{\Act{G}{\topo{X}}}{*}\xrightarrow{-//G} \oc{\topo{X}}{\mathbb BG}\] restricts to \[\lAct{G}{\topo{X}}{\lambda}\to \oc{{\topo{X}_{\lambda}}}{\mathbb BG}, \] which is an equivalence by \cite[Proposition 3.2.63]{sati2022equivariantprincipalinfinitybundles}. 
    This implies that \[-//G\colon \Act{G}{\topo{X}}\cong \oc{\Act{G}{\topo{X}}}{*}\xrightarrow{-//G} \oc{\topo{X}}{\mathbb BG}\] is an equivalence. 
\end{proof}
We recall below that $\mathbb BG=*//G$ classifies principal $\infty$-$G$-bundles. 
In preparation for that, we describe extension and restriction of scalars along a group homomorphism in the above model. A group homomorphism $\phi\colon H\to G\in \Grp(\topo{X})\hookrightarrow \Fun(\Delta^{\operatorname{op}}, \topo{X})$ induces a map $\mathbb  B\phi\colon \mathbb BH\to\mathbb BG$. 
\begin{lemma}\label{basechange}Suppose $\phi\colon H\to G\in\Grp(\topo{X})$ is a group homomorphism in a big topos $\topo{X}$. 
    Under the equivalences $\Act{G}{\topo{X}}\cong \oc{\topo{X}}{\mathbb BG}, \Act{H}{\topo{X}}\cong \oc{\topo{X}}{\mathbb BH}$, restriction along $\phi$ equals $\mathbb BH\times_{\mathbb BG}-$ and its left adjoint $\phi_{!}\colon \oc{\topo{X}}{\mathbb BH}\to\oc{\topo{X}}{\mathbb BG}$ is postcomposition with $\mathbb B\phi$. 
\end{lemma} 
\begin{proof}
    Since $\phi_{!}\dashv \phi^*$, the description of $\phi_{!}$ follows from the description of $\phi^*$.
    Choose an exhaustion $\topo{X}_*\colon\Lambda\to \Pr^L$ of $\topo{X}$ by topoi and $\mu\in\Lambda$ with $H,G\in\Grp(\topo{X}_{\mu})$. Upon replacing $\Lambda$ by $_{\mu\backslash }\Lambda$ we can assume that $H,G\in\Grp(\topo{X}_{\lambda})$ for all $\lambda\in\Lambda$. 
    Denote by \[\phi_{\Delta}\colon (*//H)_*\to (*//G)_*\in \colim{\lambda\in\Lambda}\Fun(\Delta^{\operatorname{op}}, \topo{X}_{\lambda})\] the image of $\phi$ under $\Grp(\topo{X}_{\lambda})\subseteq \colim{\lambda\in\Lambda}\Fun(\Delta^{\operatorname{op}}, \topo{X}_{\lambda})$. 
    This induces a map \[\mathbb B\phi\coloneqq \colim{\Delta^{\operatorname{op}}}\phi_{\Delta}\colon\mathbb BH\to\mathbb BG\in\topo{X}_{\lambda}\subseteq \topo{X}.\]
    The pullback \[\oc{\Fun(\Delta^{\operatorname{op}}, \topo{X})}{(*//G)_*}\xrightarrow{(*//H)_*\times_{(*//G)_*}-}\oc{\Fun(\Delta^{\operatorname{op}}, \topo{X})}{(*//H)_*}\] exists since $\topo{X}$ has finite limits and restricts to 
    \[\phi^*\colon \Act{G}{\topo{X}}\to \Act{H}{\topo{X}}.\]  
    For $\lambda\in\Lambda$, the equivalence $-//G\colon \Act{G}{\topo{X}}\to \oc{\topo{X}}{\mathbb BG}$ restricts to an equivalence $-//^{\lambda}G\colon \lAct{G}{\topo{X}}{\lambda}\to \oc{\topo{X}_{\lambda}}{\mathbb BG}$, and in $\lAct{G}{\topo{X}}{\lambda}$, \[-//^{\lambda}G\circ \bigl((*//G)_*\times_{\mathbb BG}(-//^{\lambda}G)\bigr)\cong (-//^{\lambda}G)\] by universality of colimits in $\topo{X}_{\lambda}$.
    This implies that for all $\lambda\in\Lambda$, \[\id_{\lAct{G}{\topo{X}}{\lambda}}\cong (*//G)_*\times_{\mathbb BG}(-//^{\lambda}G), \] and hence 
   \[\id_{\Act{G}{\topo{X}}}\cong (*//G)_*\times_{\mathbb BG}(-//G)\] since $\topo{X}_{\lambda}\subseteq \topo{X}$.
   In particular,  
    \[\phi^*=(*//H)_*\times_{(*//G)_*}-\cong (*//H)_*\times_{(*//G)_*}((*//G)_*\times_{\mathbb BG}(-//G))\cong (*//H)_*\times_{\mathbb BG}(-//G),\]
    so we are reduced to showing that the canonical map 
    \begin{align*} (-//H)\circ ( (*//H)_*\times_{\mathbb BG}-)\circ (-//G))&=\colim{\Delta^{\operatorname{op}}}\left(((*//H)_*\times_{\mathbb BG}-)\circ (-//G)\right)\\& \to \left((\colim{\Delta^{\operatorname{op}}}(*//H)_*)\times_{\mathbb BG}-\right)\circ (-//G)=\mathbb BH\times_{\mathbb BG}(-//G)\end{align*} is an equivalence. 
    As for all $\lambda\in\Lambda$, $\lAct{G}{\topo{X}}{\lambda}\subseteq \Act{G}{\topo{X}}$ is closed under finite limits and small colimits, it suffices to show the statement for $\topo{X}_{\lambda}$, where it holds by universality of colimits. 
\end{proof}
The above implies that for a group object $G$ in a big topos $\topo{X}$, the homotopy orbits functor $(-//G)$, $\triv$ and the free functor $\ind{G}{*}$ satisfy the following relations: 
\begin{cor}\label{orbitsoftrivialandinducedobjects}Suppose $G\in\Grp(\topo{X})$ is a group object in a big topos $\topo{X}$.
    \begin{romanenum}\item  Then \[(-//G)\circ \triv\cong -\times \mathbb BG \text{ and } (-//G)\circ \ind{G}{*}\cong \id_{\topo{X}}.\] 
    \item For $X,Y\in \topo{X}$, $\ind{G}{*}(X\times Y)\cong \ind{G}{*}(X)\times \triv(Y)$ naturally in $X$ and $Y$.
    \end{romanenum} 
\end{cor}
\begin{proof}
    By \cref{basechange}, under the equivalence $\Act{G}{\topo{X}}\cong \oc{\topo{X}}{\mathbb BG}$, 
    \[\triv(-)\cong (-\times \mathbb BG\xrightarrow{\pi_{\mathbb BG}}\mathbb BG), \]
    $\ind{G}{*}\colon\topo{X}\cong \oc{\topo{X}}{*}\to\oc{\topo{X}}{\mathbb BG}$ is postcomposition with $*\to \mathbb BG$, and  $-//G\colon\Act{G}{\topo{X}}\cong \oc{\topo{X}}{\mathbb BG}\to \topo{X}$ is the forget functor (i.e. postcomposition with $\mathbb BG\to *$).  
    This implies that \[(-//G)\circ \triv\cong -\times \mathbb BG, \, \, (-//G)\circ \ind{G}{*}\cong \id_{\topo{X}},\] and for $X,Y\in \topo{X}$, \begin{align*}\ind{G}{*}(X\times Y)\cong (X\times Y\to *\to\mathbb BG)&\cong \left((X\times \mathbb BG)\times_{\mathbb BG}Y\to Y\to *\to\mathbb BG \right)\\ &\cong \triv(X)\times_{\oc{\topo{X}}{\mathbb BG}} \ind{G}{*}(Y)\end{align*} naturally in $X$ and $Y$. 
\end{proof}
This implies that if $\topo{X}$ is a cartesian closed big topos, there exists a homotopy fixed-point functor \[(-)^G\colon \Act{G}{\topo{X}}\to\topo{X}\] for all group objects $G\in\Grp(\topo{X})$: 
\begin{cor}\label{fixedpointfunctorexistsforcartesianclosed}
    Suppose that $\topo{X}$ is a big topos such that the cartesian monoidal structure on $\topo{X}$ is closed with internal mapping spaces $\iMap_{\topo{X}}(-,-)$. 
    Then \[ \Act{G}{\topo{X}}\cong \oc{\topo{X}}{\mathbb BG}\xrightarrow{\iMap_{\topo{X}}(\mathbb BG,-)}\oc{\topo{X}}{\iMap_{\topo{X}}(\mathbb BG, \mathbb BG)}\xrightarrow{\{\id_{\mathbb BG}\}\times_{\iMap_{\topo{X}}(\mathbb BG, \mathbb BG)}-} \topo{X}\] is right adjoint to $\triv$. 
\end{cor}
\begin{proof}Denote by $q_{!}\colon \oc{\topo{X}}{\mathbb BG}\to\topo{X}$ the forget functor, i.e. postcomposition with $\mathbb BG\to *$. 
    We have the following equivalences of functors 
    $\topo{X}^{\operatorname{op}}\times\oc{\topo{X}}{\mathbb BG}\to \an$: 
    \begin{align*}\Map_{\topo{X}}(-_1, \iMap_{\topo{X}}&(\mathbb BG,q_{!}-_2))\times_{\iMap_{\topo{X}}(\mathbb BG, \mathbb BG)}\{\id_{\mathbb BG}\}\\& \cong \Map_{\topo{X}}(-_1, \iMap_{\topo{X}}(\mathbb BG,q_{!}-_2))\times_{\Map_{\topo{X}}(-_1, \iMap_{\topo{X}}(\mathbb BG, \mathbb BG))}\{ \operatorname{const}_{\id_{\mathbb BG}}\})\\ & \cong \Map_{\topo{X}}(-_1\times \mathbb BG,q_{!}-_2)\times_{\Map_{\topo{X}}(-_1\times \mathbb BG, \mathbb BG)}\{\pi_{\mathbb BG}\}\\ & \cong \Map_{\oc{\topo{X}}{\mathbb BG}}(-_1\times \mathbb BG,-_2).\end{align*}
    By \cref{orbitsoftrivialandinducedobjects}, the right-hand side is equivalent to $\Map_{\oc{\topo{X}}{\mathbb BG}}(\triv(-_1),-_2)$. 
\end{proof}

\paragraph{Principal $\infty$-$G$-bundles}
We now recall from \cite{Principalinftybundles} that the forget functor $\oc{\topo{X}}{\mathbb BG}\to\topo{X}$ classifies \textit{principal} $\infty$-$G$-\textit{bundles}. 
\begin{definition}\label{definitionprincipalinfinitybundle}
    Suppose $\topo{X}$ is a big topos topos and $G\in\Grp(\topo{X})$. 
    For $B\in\topo{X}$, a map $E\to \triv B\in \Act{G}{\topo{X}}$ is a \emph{principal} $\infty$-$G$-\emph{bundle} if there exists an effective epimorphism $p\colon U\to B\in \topo{X}$ such that \[E\times_{\triv B} \triv U\cong \ind{G}{*}(U)\in \oc{\Act{G}{\topo{X}}}{\triv B}.\]
\end{definition}
\begin{lemma}[{\cite[section 3.2]{Principalinftybundles}}]\label{homotopyquoptientbundle}Suppose $G$ is a group object in a big topos $\topo{X}$. 

For $E\in\Act{G}{\topo{X}}\subseteq \oc{\Fun(\Delta^{\operatorname{op}}, \topo{X})}{(*//G)_*}$, \[E\cong \check{C}(E([0])\to E//G)\in \oc{\Fun(\Delta^{\operatorname{op}}, \topo{X})}{(*//G)_*},\] where $E([0])\to \colim{\Delta^{\operatorname{op}}}E=E//G$ is the canonical map.  
\end{lemma}
\begin{proof}
    Choose an exhaustion $\topo{X}_*\colon\Lambda\to \Pr^L$ of $\topo{X}$ by topoi such that $G\in\Grp(\topo{X}_{\lambda})$ for all $\lambda\in\Lambda$. 
    Then $\lAct{G}{\topo{X}}{*}\colon \Lambda\to \Pr^L$ is an exhaustion of $\Act{G}{\topo{X}}$ by topoi (\cref{Gobjectsisbigtopos}). 
    For $E\in \Act{G}{\topo{X}}$ choose $\lambda\in\Lambda$ with $E\in\lAct{G}{\topo{X}}{\lambda}$.
    The associated element $E\in\oc{\Fun(\Delta^{\operatorname{op}}, \topo{X}_{\lambda})}{(*//G)_*}\subseteq \oc{\Fun(\Delta^{\operatorname{op}}, \topo{X})}{(*//G)_*}$ is a groupoid over $(*//G)_*$, whence 
    \[E\cong \check{C}\bigl(E([0])\to \colim{\Delta^{\operatorname{op}}}^{\topo{X}_{\lambda}}E\bigr)\in \oc{\Fun(\Delta^{\operatorname{op}}, \topo{X}_{\lambda})}{(*//G)_*}\] by effectivity of groupoids in the topos $\topo{X}_{\lambda}$.  
    Since $\topo{X}_{\lambda}\subseteq \topo{X}$ is closed under finite limits and small colimits, 
    this implies that \[E\cong \check{C}\bigl(E([0])\to \colim{\Delta^{\operatorname{op}}}^{\topo{X}}E\bigr)\in \oc{\Fun(\Delta^{\operatorname{op}}, \topo{X})}{(*//G)_*}.\qedhere\]
\end{proof}
\begin{cor}[{\cite[section 3.2]{Principalinftybundles}}]\label{characterizationbundles}Suppose $G$ is a group object in a big topos $\topo{X}$. 
    \begin{romanenum}
    \item \label{characterisationbundles1} A map $p\colon E\to \triv B\in \Act{G}{\topo{X}}\subseteq \oc{\Fun(\Delta^{\operatorname{op}}, \topo{X})}{(*//G)_*}$ is a principal $\infty$-$G$-bundle if and only if \[E//G\xrightarrow{p//G}\triv(B)//G\cong \mathbb BG\times B\xrightarrow{\pi_B} B\] is an equivalence. 

    \item For $B\in\topo{X}$ denote by $\Bun{G}{\topo{X}}{B}\subseteq \oc{\Act{G}{\topo{X}}}{\triv B}$ the full subcategory on principal $\infty$-$G$-bundles and by $(\triv B\to *)_{!}\colon \oc{\Act{G}{\topo{X}}}{\triv B}\xrightarrow{(\triv B\to *)_{!}}\oc{\Act{G}{\topo{X}}}{*}$ postcomposition with $\triv B\to *$. 
    The functor \[\oc{\Act{G}{\topo{X}}}{\triv B}\xrightarrow{(\triv B\to *)_{!}}\oc{\Act{G}{\topo{X}}}{*}\xrightarrow{-//G}\oc{\topo{X}}{\mathbb BG}\] induces an equivalence \[ \Bun{G}{B}{\topo{X}}\to \oc{\topo{X}}{\mathbb BG}\times_{\topo{X}}\{B\}\cong \Map_{\topo{X}}(B, \mathbb BG).\]
    \end{romanenum} 
\end{cor}
\begin{proof}We begin with the first statement. 
    Suppose that $E\to \triv B\in\Act{G}{\topo{X}}$ is a principal $\infty$-$G$-bundle and $p\colon U\to B\in \topo{X}$ such that \[E\times_{\triv B} \triv U\cong \ind{G}{*}(U)\in \oc{\Act{G}{\topo{X}}}{\triv B}.\]
    Choose an exhaustion $\topo{X}_*\colon\Lambda\to \Pr^L$ of $\topo{X}$ by topoi such that $G\in\Grp(\topo{X}_{\lambda})$ for all $\lambda\in\Lambda$, and choose $\lambda\in\Lambda$ with $E\in\lAct{G}{\topo{X}}{\lambda}$ and $U,X\in\topo{X}_{\lambda}$.
    As $\topo{X}_{\lambda}\subseteq \topo{X}$ is closed under finite limits and small colimits (\cref{fullyfaithfulnessandpreservationoflimits}), $p\colon U\to X$ is an effective epimorphism in $\topo{X_{\lambda}}$. 
    The functor $\triv\colon\topo{X}\to\Act{G}{\topo{X}}$ restricts to $\triv^{\lambda}\colon \topo{X_{\lambda}}\to\lAct{G}{\topo{X}}{\lambda}$, the restriction along $G\to *\in\Grp(\topo{X}_{\lambda})$.  
    Since $\triv^{\lambda}$ preserves small limits and colimits, $\triv^{\lambda}(p)$ is an effective epimorphism in $\lAct{G}{\topo{X}}{\lambda}$.
    As colimits in $\lAct{G}{\topo{X}}{\lambda}$ are universal and $\lAct{G}{\topo{X}}{\lambda}\subseteq \Act{G}{\topo{X}}$ is closed under finite limits and small colimits (\cref{Gobjectsisbigtopos}, \cref{fullyfaithfulnessandpreservationoflimits}), this implies that 
    \[ E\cong \colim{\Delta^{\operatorname{op}}}E\times_{\triv B}\check{C}(\triv p)\] in $\lAct{G}{\topo{X}}{\lambda}$ and $\Act{G}{\topo{X}}$.  
    By \cref{orbitsoftrivialandinducedobjects}, for $n\in\mathbb N_0$,
    \begin{align*}E\times_{\triv B}  \check{C}(\triv p) ([n])& =E\times_{\triv B} \triv U\times_{\triv B} \triv\underbrace{(U\times_B\ldots \times_BU)}_{n\text{ times }}\\
    &\cong (\ind{G}{*}(U))\times_{\triv B}\triv\underbrace{(U\times_B\ldots \times_BU)}_{n\text{ times }}\\&\cong \ind{G}{*}(*)\times \triv U\times_{\triv B}\triv\underbrace{(U\times_B\ldots \times_BU)}_{n\text{ times }}
    \\&\cong \ind{G}{*}(*)\times \triv\check{C}(p)([n])\\&\cong \ind{G}{*}(\check{C}(\triv p)[n]).\end{align*} This defines an isomorphism of simplicial objects $E\times_{\triv B}\check{C}(\triv p)\cong \ind{G}{*}(\check{C}(\triv(p)))$. 
    As $\Act{G}{\topo{X}}\xrightarrow{-//G} \topo{X}$ preserves colimits, it follows from \cref{orbitsoftrivialandinducedobjects} that \[E//G\cong \left(\colim{\Delta^{\operatorname{op}}}\ind{G}{*}\check{C}(p)\right)//G\cong \colim{\Delta^{\operatorname{op}}}\check{C}(p)\cong B.\]
    By construction, this equivalence is the composition \[E//G\to \triv B//G\cong \mathbb BG\times B\xrightarrow{\pi_B} B.\] 

    For the converse implication of \ref{characterisationbundles1}, it suffices to show that for all $E\in\Act{G}{\topo{X}}$, the adjunction unit of $-//G\dashv \triv$ defines a bundle map $E\to \triv E//G$. 
    By \cref{homotopyquoptientbundle}, \[E\cong \check{C}(E([0])\to E//G)\in \oc{\Fun(\Delta^{\operatorname{op}}, \topo{X})}{(*//G)_*}, \] whence $E([0])\to E//G$ is an effective epimorphism in $\topo{X}$.
    For $X\in\topo{X}$, \[\triv X\in \Act{G}{\topo{X}}\subseteq \oc{\Fun(\Delta^{\operatorname{op}}, \topo{X})}{(*//G)_*}\] is the simplicial object \[(*//G)_*\times X\xrightarrow{\pi} (*//G)_*\in \oc{\Fun(\Delta^{\operatorname{op}}, \topo{X})}{(*//G)_*},\] where $\pi$ denotes the projection.  
    \cref{simplicialmodelgobjects} implies that \[\Act{G}{\topo{X}}\subseteq  \oc{\Fun(\Delta^{\operatorname{op}}, \topo{X})}{(*//G)}\] is closed under pullbacks which are computed pointwise. This implies that $E\times_{\triv E//G}\triv E([0])$ is represented by the simplicial object \[[n]\mapsto E([n])\times_{E//G}E([0]), \] where the pullback is computed with respect to the canonical maps $E([n])\to \colim{\Delta^{\operatorname{op}}}E=E//G$.
    By \cref{homotopyquoptientbundle}, this is equivalent to \[\check{C}(E[0]\to E//G)\times_{E//G}E([0])\cong \check{C}(E([0])\times_{E//G}E([0])\xrightarrow{\pi_1}E([0])), \]
    and \[ (E([0])\times_{E//G}E([0])\xrightarrow{\pi_1}E([0]))\cong (d_1\colon E([1])\to E([0]))\] is the face map of the simplicial object $E$. By \cref{simplicialmodelgobjects}, \[(d_1\colon E([1])\to E([0]))\cong (G\times E([0])\xrightarrow{\pi_{E([0])}}E([0])).\]
    This shows that \[E\times_{\triv E//G}\triv E([0])\cong \check{C}(G\to *)\times E([0])\cong \ind{G}{*}(E[0])\in \oc{\Fun(\Delta^{\operatorname{op}}, \topo{X})}{(*//G)_*}, \] which proves the characterization of bundles. 

    For $B\in \topo{X}$, the equivalence $-//G\colon \Act{G}{\topo{X}}\to \oc{\topo{X}}{\mathbb BG}$ induces an equivalence \[ \oc{\Act{G}{\topo{X}}}{\triv B}\cong \oc{\topo{X}}{B\times \mathbb BG}.\] Denote by \[\pi\colon \oc{\Act{G}{\topo{X}}}{\triv B}\cong \oc{\topo{X}}{B\times\mathbb BG}\xrightarrow{(B\times\mathbb BG\to B)_{!}} \oc{\topo{X}}{B}\] the composition of this equivalence with postcomposition with the projection $B\times\mathbb BG\to B$. 
    By the above characterization of bundles, \begin{center}
    \begin{tikzcd}
    \Bun{G}{B}{\topo{X}}\arrow[r]\arrow[d] & \oc{\Act{G}{\topo{X}}}{\triv B}\arrow[d,"\pi"]\\ 
    \{ \id_B\} \arrow[r] & \oc{\topo{X}}{B}
    \end{tikzcd}
    \end{center} is a pullback. 
    As $\oc{\Act{G}{\topo{X}}}{\triv B}\cong \oc{\topo{X}}{B\times\mathbb BG}\cong \oc{\topo{X}}{B}\times_{\topo{X}}\oc{\topo{X}}{\mathbb BG}$, it follows that 
    \[ \Bun{G}{B}{\topo{X}}\cong \oc{\topo{X}}{\mathbb BG}\times_{\topo{X}}\{B\}=\Map_{\topo{X}}(B,\mathbb BG).\qedhere\] 
\end{proof}
We have shown in the above proof that for all $E\in\Act{G}{\topo{X}}$, the adjunction unit for $(-//G)\dashv \triv$ evaluates to a principal $\infty$-$G$-bundle $E\to \triv E//G$. 
In particular, there is no distinction between objects with an $\infty$-$G$-action and (total spaces) of principal $\infty$-$G$-bundles. 
The unit evaluates to a principal $\infty$-$G$-bundle $*\to \triv \mathbb BG$, where $*$ denotes the terminal object of $\Act{G}{\topo{X}}$.
One should think of this as the universal principal $\infty$-$G$-bundle: The above proof shows that for $E\in\Act{G}{\topo{X}}$, 
\[E\cong E//G\times_{\mathbb BG}*\in \Act{G}{\topo{X}},\] i.e. $E$ is \textit{classified} by the map $(E\to *)//G\colon E//G\to \mathbb BG$. 

\begin{cor}\label{deloopingclassifyingspace}
For a group object $G\in\Grp(\topo{X})$ in a big topos $\topo{X}$, \[\Omega \mathbb BG\cong G\in\topo{X}.\] 
\end{cor}
\begin{proof}
    By \cref{homotopyquoptientbundle}, \[ G\cong (*//G)_*([1])\cong \check{C}(*\to \mathbb BG)([1])\cong *\times_{\mathbb BG}*\cong \Omega \mathbb BG.\qedhere\] 
\end{proof}

\paragraph{Modules over group rings}
\begin{definition}
    Suppose $\topo{X}$ is a big topos. The symmetric monoidal functor \[\Sigma^{\infty}_{+}\colon \topo{X}^{\times}\to\stab{\topo{X}}^{\otimes} \text{ (\cref{symmetricmonoidalstructureonspectrumobjects})}\] induces a functor $\mathbb S[-]\colon \Alg(\topo{X}^{\times})\to \Alg(\stab{\topo{X}})$. 
    Recall from \cite[Proposition 3.2.4.3/Example 3.2.4.4]{higheralgebra} that for a symmetric monoidal category $\mathcal C^{\otimes}$, $\Alg(\mathcal C)$ admits a monoidal structure such that the forget functor $\Alg(\mathcal C)\to\mathcal C$ enhances to a monoidal functor.
    For a ring spectrum $\algebra{R}\in\Alg(\stab{\topo{X}})$ and a group object $G\in\Grp(\topo{X})$, the $\algebra{R}$-linear \emph{group ring} of $G$ is 
    \[\algebra{R}[G]\coloneqq \algebra{R}\otimes_{\Alg(\stab{\topo{X}})} \mathbb S[G]\in\Alg(\stab{\topo{X}}).\]  
\end{definition}
Our next goal is to show that under mild conditions on a big topos $\topo{X}$, there is an equivalence \[\LMod{\triv_{\Sp}\algebra{R}}{\stab{\Act{G}{\topo{X}}}}\cong \LMod{\algebra{R}[G]}{\stab{\topo{X}}},\] see \cref{modulesingroupobjectsaremodulesovergrouprings} below. 
\begin{lemma}\label{GobjectsstabilisetoGmodules}Suppose that $G$ is a group object in a big topos $\topo{X}$.
    The stabilization (\cref{definitionstabilization}) of $\res{G}{*}\colon \Act{G}{\topo{X}}\to \mathcal X$ factors over an equivalence \[\stab{\Act{G}{\topo{X}}}\cong \LMod{\mathbb S[G]}{\stab{\topo{X}}}\xrightarrow{f}\stab{\topo{X}}, \] where $f$ denotes the forget functor. 
\end{lemma}
\begin{proof}
    For a symmetric monoidal category $\mathcal C$, $\mathcal C\hookrightarrow \Fun(\mathcal C, \mathcal C), c\mapsto c\otimes-$ enhances to a symmetric monoidal functor, where $\Fun(\mathcal C, \mathcal C)$ is endowed with the symmetric monoidal structure given by composition, cf.\ \cite[page 656]{higheralgebra}/\cref{compositionmonoidalstructure}. 

    The forget functor $\Act{G}{\topo{X}}=\LMod{G}{\topo{X}}\to\topo{X}$ exhibits $\Act{G}{\topo{X}}$ as monadic over $\topo{X}$ with monad \[G\in \Alg(\topo{X}^{\times})\subseteq \Alg(\Fun(\topo{X}, \topo{X})).\] Hence by \cite[Example 4.7.3.10/Theorem 4.7.3.5]{higheralgebra}, its stabilization exhibits $\stab{\Act{G}{\topo{X}}}$ as monadic over $\stab{\topo{X}}$ with monad \[\stab{G\times -}=\mathbb S[G]\otimes-\in \Alg(\stab{\topo{X}})\subseteq \Alg(\Fun(\stab{\topo{X}}, \stab{\topo{X}})), \] i.e.\ $\stab{\Act{G}{\topo{X}}}\cong \LMod{\mathbb S[G]}{\stab{\topo{X}}}$.
\end{proof}
\begin{rem}
If it exists, a commutative enhancement of $G$ determines a commutative enhancement of $\mathbb S[G]$ which defines a symmetric monoidal structure on $\LMod{\mathbb S[G]}{\stab{\topo{X}}}$. For $G\neq 0$, this is different from the structure induced by the cartesian monoidal structure on $\Act{G}{\topo{X}}$: The forget functor $\stab{\Act{G}{\topo{X}}}\to\stab{\topo{X}}$ is symmetric monoidal, which fails for the symmetric monoidal structure on $\LMod{\mathbb S[G]}{\stab{\topo{X}}}$ induced by the commutative algebra structure of $\mathbb S[G]$ (if $G\neq 0$).
\end{rem}
\begin{definition}
    If $G$ is a group object in a big topos $\topo{X}$, denote by $\res{G}{*, \Sp}\colon \stab{\Act{G}{\topo{X}}}\to\stab{\topo{X}}$ the stabilization of the restriction $\res{G}{*}\colon \Act{G}{\topo{X}}\to \topo{X}$.
    By \cref{GobjectsstabilisetoGmodules} and \cref{forgetfreeadjunctionmodules},
    this admits a left adjoint which we denote by $\ind{G}{*, \Sp}$.  
\end{definition}
\begin{rem}\label{indgspisstabilization}
Since $\Omega^{\infty}\res{G}{*, \Sp}\cong \res{G}{*}\Omega^{\infty}$, $\Sigma^{\infty}_{+}\ind{G}{*}\cong \ind{G}{*, \Sp}\Sigma^{\infty}_{+}$.
\end{rem}
\cref{universalpropertyspectrumobjectsbigpresentable} and \cref{orbitsoftrivialandinducedobjects} imply the following: 
\begin{lemma}\label{inductiontrivialstable}
Suppose $G$ is a group object in a big topos $\topo{X}$. 
There is a preferred equivalence \[\ind{G}{*, \Sp}(-_1)\otimes \triv_{\Sp}(-_1)\cong \ind{G}{*, \Sp}(-_1\otimes -_2)\in\Fun(\stab{\topo{X}}\times\stab{\topo{X}}, \stab{\Act{G}{\topo{X}}}).\] 
\end{lemma}
\begin{proof}
Choose an exhaustion $\topo{X}_*\colon\Lambda\to\CAlg(\Pr^L)$ of $\topo{X}$ by topoi such that $G\in\Grp(\topo{X}_{\lambda})\subseteq \Grp(\topo{X})$ for all $\lambda\in\Lambda$. 
Denote by \[\Fun^{\operatorname{colim},\Lambda,v}(\stab{\topo{X}}\times\stab{\topo{X}}, \stab{\Act{G}{\topo{X}}})\subseteq \Fun(\stab{\topo{X}}\times\stab{\topo{X}}, \stab{\Act{G}{\topo{X}}})\] the full subcategory on functors $F$ such that for all $\lambda\in\Lambda$, $F|_{\stab{\topo{X}_{\lambda}}\times\stab{\topo{X}_{\lambda}}}$ preserves small colimits separately in both variables. 
Recall from \cref{stabilizationbigpresentablecategoriesbigpresentable} that $\stab{\topo{X}_*}$ is an exhaustion of $\stab{\topo{X}}$ by presentable categories. 
By construction (see the proof of \cref{existencesuspension}), \[\Sigma^{\infty}_{\topo{X},+}=\colim{\lambda\in\Lambda}\left((\stab{\topo{X}_{\lambda}}\hookrightarrow \stab{\topo{X}})\circ \Sigma^{\infty}_{+,\topo{X}_{\lambda}}\right).\] Alternatively, this follows from \cref{stabilizationcommuteswithsuspension,fullyfaithfulnessandpreservationoflimits}. 
\cref{universalpropertyspectrumobjectsbigpresentable,limitsincat} now imply that  
\begin{align*} \Fun^{\operatorname{colim},\Lambda,v}(\stab{\topo{X}}\times\stab{\topo{X}}, \stab{\Act{G}{\topo{X}}}) &\to \clim{\lambda\in\Lambda}\Fun(\topo{X}_{\lambda}\times \topo{X}_{\lambda}, \stab{\Act{G}{\topo{X}}})\cong \Fun(\topo{X}\times \topo{X}, \stab{\Act{G}{\topo{X}}}), \\ F&\mapsto F\circ (\Sigma^{\infty}_{+}\times \Sigma^{\infty}_{+})\end{align*} is fully faithful. We now explain that \[ \ind{G}{*, \Sp}(-_1)\otimes \triv_{\Sp}(-_1), \ind{G}{*, \Sp}(-_1\otimes -_2)\in \Fun^{\operatorname{colim},\Lambda,v}(\stab{\topo{X}}\times\stab{\topo{X}}, \stab{\Act{G}{\topo{X}}}).\]
By construction (\cref{existenceleftadjointpullback}), $\triv$ restricts to $\triv^{\lambda}\colon \topo{X}_{\lambda}\to\lAct{G}{\topo{X}}{\lambda}$ for all $\lambda\in\Lambda$, and hence $\triv_{\Sp}$ restricts to $\triv^{\lambda}_{\Sp}\colon \stab{\topo{X}_{\lambda}}\to \stab{\lAct{G}{\topo{X}}{\lambda}}$.
This functor is cocontinuous. 
As $\stab{\topo{X}_*}^{\otimes}$ is an exhaustion of $\stab{\topo{X}}^{\otimes}$ by presentably symmetric monoidal categories, the proof of \cref{existenceleftadjointpullback} shows that $\ind{G}{*,\Sp}$ restricts to left adjoints 
$\stab{\topo{X}_{\lambda}}\to\stab{\lAct{G}{\topo{X}}{\lambda}}$ for all $\lambda\in\Lambda$. 
Since $\stab{\lAct{G}{\topo{X}}{\lambda}}\subseteq \stab{\Act{G}{\topo{X}}}$ and $\stab{\topo{X}_{\lambda}}\subseteq \stab{\topo{X}}$ are presentably symmetric monoidal subcategories (\cref{bigpresentablymonoidalcatdef,symmetricmonoidalstructureonspectrumobjects}), it follows that
\[ \ind{G}{*, \Sp}(-_1)\otimes \triv_{\Sp}(-_1), \ind{G}{*, \Sp}(-_1\otimes -_2)\in \Fun^{\operatorname{colim},\Lambda,v}(\stab{\topo{X}}\times\stab{\topo{X}}, \stab{\Act{G}{\topo{X}}}).\]
We are therefore reduced to describing an equivalence 
\[\ind{G}{*, \Sp}\Sigma^{\infty}_{+}(-_1)\otimes \triv_{\Sp}\Sigma^{\infty}_{+}(-_2)\cong \ind{G}{*, \Sp}(\Sigma^{\infty}_{+}(-1)\otimes \Sigma^{\infty}_{+}(-_2)).\] 

Since $\ind{G}{*, \Sp}\Sigma^{\infty}_{+}\cong \Sigma^{\infty}_{+}\ind{G}{*}$ and $\triv_{\Sp}\Sigma^{\infty}_{+}\cong \Sigma^{\infty}_{+}\triv$ (\cref{indgspisstabilization,monoidalstructurestabilizationpullback}),
symmetric monoidality of $\Sigma^{\infty}_{+}$ yields an equivalence 
\[ \ind{G}{*, \Sp}\Sigma^{\infty}_{+}(-_1)\otimes \triv_{\Sp}\Sigma^{\infty}_{+}(-_2)\cong \Sigma^{\infty}_{+}(\ind{G}{*}(-_1)\times \triv(-_1)).\] 
By \cref{orbitsoftrivialandinducedobjects}, $\ind{G}{*}(-_1)\times \triv(-_2)\cong \ind{G}{*}(-_1\times -_2)$, whence \[ \Sigma^{\infty}_{+}(\ind{G}{*}(-_1)\times\triv(-_2))\cong \Sigma^{\infty}_{+}\ind{G}{*}(-_1\times-_2).\] 
As $\ind{G}{*, \Sp}\Sigma^{\infty}_{+}\cong \Sigma^{\infty}_{+}\ind{G}{*}$ and $\Sigma^{\infty}_{+}$ is symmetric monoidal, the right-hand side is equivalent to $\ind{G}{*, \Sp}(\Sigma^{\infty}_{+}(-_1)\otimes \Sigma^{\infty}_{+}(-_2))$. 
\end{proof}

\begin{cor}\label{modulesingroupobjectsaremodulesovergrouprings} 
Suppose that $\topo{X}$ is a big topos such that $\stab{\topo{X}}$ admits $\Delta^{\operatorname{op}}$-indexed colimits and the tensor product on $\stab{\topo{X}}$ preserves $\Delta^{\operatorname{op}}$-indexed colimits in both variables. 

For $\algebra{R}\in \Alg(\stab{\topo{X}})$ and $G\in\Grp(\topo{X})$, 
\[\LMod{\triv_{\Sp}\algebra{R}}{\stab{\Act{G}{\topo{X}}}}\cong \LMod{\algebra{R}[G]}{\stab{\topo{X}}}.\] 
\end{cor}
\begin{proof}The symmetric monoidal functor $\triv_{\Sp}^{\otimes}$ (\cref{monoidalstructurestabilizationpullback}) exhibits $\LMod{\mathbb S[G]}{\stab{\topo{X}}}$ as left-tensored over $\stab{\topo{X}}$. We claim that this is the left tensoring induced by the cocartesian fibration of operads \[\stab{\topo{X}}^{\otimes}\times_{\Comm^{\otimes}}\BM^{\otimes}\to \BM^{\otimes}\] as described in \cite[Remark 4.3.3.7]{higheralgebra}, then the statement follows from \cref{modulesinmodulecategories}. We will deduce this from \cite[Proposition 4.8.5.8]{higheralgebra}. 
By \cref{GobjectsstabilisetoGmodules}, \cref{forgetfreeadjunctionmodules}, $\res{G}{*, \Sp}\colon \stab{\Act{G}{\topo{X}}}\to \stab{\topo{X}}$ is conservative and preserves $\Delta^{\operatorname{op}}$-indexed colimits. 
By \cref{monoidalstructurestabilizationpullback}, it enhances to a symmetric monoidal functor, which implies that the tensor product on $\stab{\Act{G}{\topo{X}}}$ preserves $\Delta^{\operatorname{op}}$-indexed colimits in both variables. 
Since \[\res{G}{*, \Sp}\circ \triv_{\Sp}=\id_{\stab{\topo{X}}}, \] it follows that $\triv_{\Sp}$ preserves $\Delta^{\operatorname{op}}$-indexed colimits. This shows that the action map \[-\otimes -\colon \stab{\topo{X}}\times \LMod{\mathbb S[G]}{\stab{\topo{X}}}\to \LMod{\mathbb S[G]}{\stab{\topo{X}}}\] preserves $\Delta^{\operatorname{op}}$-indexed colimits in both variables. 
By \cref{inductiontrivialstable}, \[L\coloneqq \mathbb S[G]\otimes -\colon \stab{\topo{X}}\to \stab{\Act{G}{\topo{X}}}\] is equivalent to $\ind{G}{*, \Sp}(-)$, and hence right adjoint to $\res{G}{*, \Sp}$. 
For $m\in\LMod{\mathbb S[G]}{\stab{\topo{X}}}$ and $x\in\stab{\topo{X}}$, the counit $\epsilon\colon L\res{G}{*, \Sp}\to \id$ and the identification $\ind{G}{*, \Sp}(-\otimes -)\cong \triv_{\Sp}(-\otimes-)$ from the proof of \cref{inductiontrivialstable} induce a map \[ c_{x,m}\colon L(x\otimes \res{G}{*, \Sp}m)\cong \triv_{\Sp}(x)\otimes L(\res{G}{*, \Sp}m)\xrightarrow{\id\otimes \epsilon}\triv_{\Sp}(x)\otimes m.\] 
By construction of the identification $L(x\otimes \res{G}{*, \Sp}m)\cong \triv_{\Sp}(x)\otimes L(\res{G}{*, \Sp}m)$, the adjoint map 
\[ x\otimes \res{G}{*, \Sp}m\to \res{G}{*, \Sp}(\triv_{\Sp}x\otimes m)\] is an equivalence.
\cite[Proposition 4.8.5.8]{higheralgebra} now implies that there exists a $\stab{\topo{X}}$-linear equivalence 
\[\RMod{A}{\stab{\topo{X}}}\cong \LMod{\mathbb S[G]}{\stab{\topo{X}}},\] where $A=\operatorname{Mor}_{\stab{\topo{X}}}(\mathbb S[G], \mathbb S[G])$ the endomorphism algebra of $\mathbb S[G]$ and $\LMod{\mathbb S[G]}{\stab{\topo{X}}}$ is left-tensored over $\stab{\topo{X}}$ via $\triv_{\Sp}^{\otimes}$.  
The algebra structure on $A$ is induced by the counit for $L\dashv \res{G}{*, \Sp}$ and the equivalence $c_{\mathbb S[G],S[G]}$, i.e.\ under the fully faithful functor \[\Alg(\stab{\topo{X}}^{\otimes})\cong \Alg(\stab{\topo{X}}^{\otimes}_{\operatorname{rev}})\hookrightarrow \Alg(\Fun(\stab{\topo{X}}, \stab{\topo{X}}))\] induced by the symmetric monoidal functor \begin{align*}\stab{\topo{X}}\cong \stab{\topo{X}}_{\operatorname{rev}}& \hookrightarrow \Fun(\stab{\topo{X}}, \stab{\topo{X}})\\  c&\mapsto c\otimes -, \end{align*} $A$ becomes the endomorphism monad of the adjunction $L\dashv \res{G}{*, \Sp}$. By \cref{GobjectsstabilisetoGmodules}, the equivalence \[L\cong \ind{G}{*, \Sp}(-)\] yields an equivalence $A\cong \mathbb S[G]^{\operatorname{rev}}\in\Alg(\stab{\topo{X}}).$ 
\end{proof}

\cref{modulesingroupobjectsaremodulesovergrouprings} and \cref{derivedcategorymodules} imply the following: 
\begin{cor}\label{derivedgmodules}
    Suppose that $\topo{X}$ is a big topos with the following properties: 
    \begin{romanenum}\item  $\topo{X}$ admits an exhaustion (\cref{definitionbigpresentable}) by hypercomplete, 1-localic topoi, 
    \item $\Ab(\tau_{\leq 0}\topo{X})$ has countable coproducts, 
    \item $\stab{\topo{X}}$ has $\Delta^{\operatorname{op}}$-indexed colimits and the tensor product on $\stab{\topo{X}}$ preserves them in both variables. 
    \end{romanenum}
    For a discrete group $G\in\Grp(\tau_{\leq 0}\topo{X})$ and $\algebra{R}\in \Alg(\Ab(\tau_{\leq 0}\topo{X}))$,   
    \[\mathcal D(\LMod{\algebra{R}[G]}{\Ab(\tau_{\leq 0}\topo{X})})\cong \LMod{\triv_{\Sp}R}{\stab{\Act{G}{\topo{X}}}}.\] 
\end{cor} 
\cref{derivedgmodules} in particular applies to $\Cond{(\kappa)}(\an)$, or more generally to categories of accessible hypersheaves on a hyperaccessible explicit covering site $(\mathcal C,S)$ with $\mathcal C$ is a 1-category. 
\begin{proof}By \cref{modulesingroupobjectsaremodulesovergrouprings}, 
\[\LMod{\algebra{R}[G]}{\stab{\topo{X}}}\cong \LMod{\triv_{\Sp}\algebra{R}}{\stab{\Act{G}{\topo{X}}}}.\]
The statement now follows from \cref{derivedcategorymodules}. 
\end{proof}
We will use the following observation to describe $\algebra{R}$-linear enhancements of group cohomology, see \cref{remarkgroupcohomologyrefines} below. 
\begin{lemma}\label{globalsectionsforgobjectsfactors}
    Suppose $G$ is a group object in a big topos $\topo{X}$. Denote by $c_{\topo{X}}\colon \an\to \topo{X}$ and $c_{\Act{G}{\topo{X}}}\colon \an\to\Act{G}{\topo{X}}$ the constant sheaf functors (\cref{constantsheaffunctordefinition}) for $\topo{X}$ and $\Act{G}{\topo{X}}$, respectively. 

    There is an essentially unique equivalence $\triv\circ c_{\topo{X}}=c_{\Act{G}{\topo{X}}}$. 
\end{lemma}
\begin{proof}
    By \cref{existenceconstantsheaf}, it suffices to show that \[ \an\xrightarrow{c_{\topo{X}}}\topo{X}\xrightarrow{\triv}\Act{G}{\topo{X}}\] preserves colimits and the terminal object. Choose an exhaustion $\topo{X}_*\colon\Lambda\to \Pr^L$ of $\topo{X}$ by topoi and $\lambda\in\Lambda$ with $G\in\Grp(\topo{X}_{\lambda})$. 
    Then $\topo{X}_{\lambda}\to\topo{X}\xrightarrow{\triv}\Act{G}{\topo{X}}_{\lambda}$ factors as $\topo{X}_{\lambda}\xrightarrow{\triv^{\lambda}}\lAct{G}{\topo{X}}{\lambda}\to \Act{G}{\topo{X}}$, where $\triv^{\lambda}$ denotes restriction of scalars along $G\to *\in \Grp(\topo{X}_{\lambda})$. By \cite[Corollary 5.5.2.9]{highertopostheory} and \cref{existenceinduction} $\triv^{\lambda}$ is a left-exact left adjoint, and by \cref{Gobjectsisbigtopos,fullyfaithfulnessandpreservationoflimits}, $\lAct{G}{\topo{X}}{\lambda}\to\Act{G}{\topo{X}}$ is a left-exact left adjoint. 
    \cref{fullyfaithfulnessandpreservationoflimits,existenceconstantsheaf} imply that the constant sheaf functor for $\topo{X}$ factors as $\an\xrightarrow{c_{\lambda}}\topo{X}_{\lambda}\to \topo{X}$, where $c_{\lambda}$ denotes the constant sheaf functor for $\topo{X}_{\lambda}$. 
    This shows that $\triv\circ c_{\topo{X}}$ preserves colimits and the terminal object, and hence there is an essentially unique equivalence $\triv\circ c_{\topo{X}}\cong c_{\Act{G}{\topo{X}}}$.  
\end{proof}
\begin{rem}\label{remarkglobalsectionsforgobjectsfactors}
    The essentally unique equivalence $\triv\circ c_{\topo{X}}\cong c_{\Act{G}{\topo{X}}}$ defines a mate transformation \[ \beta\colon \Gamma_{\Act{G}{\topo{X}}}\circ \triv\to \Gamma_{\topo{X}},\] where $\Gamma_{\Act{G}{\topo{X}}}$, $\Gamma_{\topo{X}}$ denote the global sections functors of $\Act{G}{\topo{X}}$ and $\topo{X}$, respectively. This stabilizes to a natural transformation 
    \[\beta_{\Sp}\colon \Gamma_{\Act{G}{\topo{X}},\Sp}\circ \triv_{\Sp}\to \Gamma_{\topo{X},\Sp}\in \Fun(\stab{\topo{X}},\Sp).\] 
By \cref{geometricmorphismstabilization}, $\beta_{\Sp}$ is the Beck-Chevalley transformation associated to the commutative diagram 
\begin{center}
\begin{tikzcd}
\Sp\arrow[r,"c_{\stab{\topo{X}}}"]\arrow[d,"\id"]& \stab{\topo{X}}\arrow[d,"\triv_{\Sp}"]\\
\Sp\arrow[r,"c_{\stab{\Act{G}{\topo{X}}}}"]& \stab{\Act{G}{\topo{X}}}.
\end{tikzcd}
\end{center}
The symmetric monoidal structure of $c_{\stab{\topo{X}}},c_{\stab{\Act{G}{\topo{X}}}}$ (\cref{constantsheafsymmetricmonoidal}) and $\triv_{\Sp}$ (\cref{monoidalstructurestabilizationpullback}) determine an enhancement of $\beta_{\Sp}$ to a transformation of lax symmetric monoidal functors. This follows from \cref{adjunctionfiberwise}. In particular, we obtain a natural transformation
    \[ \beta_{\Alg,\Sp}\colon  \Gamma_{\Act{G}{\topo{X}},\Sp}\circ \triv_{\Sp}\to \Gamma_{\topo{X},\Sp}\in \Fun\left(\aCAlg(\stab{\topo{X}}),\aCAlg(\Sp)\right)\] between the induced functors on (commutative) algebra objects. 
\end{rem}
\subsection{Group cohomology in a big topos}\label{section:groupcohomologybigtopos}
We are now prepared to discuss two notions of group cohomology in a big topos.
\begin{definition}[Group cohomology]\label{definitiongroupcohomologybigtopos}
    Suppose $\topo{X}$ is a big topos and $G\in\Grp(\topo{X})$ is a group object. 

    \begin{romanenum}
    \item Define \emph{group cohomology} of $G$ in $\topo{X}$ as 
    \[ \grpcoh{\mathcal X}(G,-)\coloneqq \cH{\Act{G}{\topo{X}}}(*,-)\colon \stab{\Act{G}{\topo{X}}}\to\stab{\topo{X}}\] the cohomology (\cref{definitioncohomologyinatopos}) of the terminal object in the big topos $\Act{G}{\topo{X}}$. 
  
    \item If $\triv_{\Sp}\colon \stab{\topo{X}}\to\stab{\Act{G}{\topo{X}}}$ admits a right adjoint, we denote it by \[\intgrpcoh{\topo{X}}(G,-)\colon \stab{\Act{G}{\topo{X}}}\to\stab{{\topo{X}}}\] and refer to it as \emph{internal group cohomology} of $G$ in $\topo{X}$.
    \end{romanenum}
\end{definition}
\begin{rems}\label{groupcohomologyinternalrecovers}
    \begin{enumerate}
\item Group cohomology is the stabilization (\cref{definitionstabilization}) of the global sections functor $\Gamma_{\Act{G}{\topo{X}}}\colon \Act{G}{\topo{X}}\to\an$ (\cref{constantsheaffunctordefinition}), cf. \cref{cohomologyglobalsections}. 
\item Denote by $\Gamma_{\Sp,\topo{X}}\colon \stab{\topo{X}}\to \Sp$ the stabilization of the global sections functor in $\topo{X}$. 
If internal group cohomology exists, then \[\Gamma_{\Sp}\circ \intgrpcoh{\topo{X}}(G,-)\cong\grpcoh{\topo{X}}(G,-)\] since both are right adjoint to the stabilization of the constant sheaf functor $c_{\Act{G}{\topo{X}}}\colon \an\to \Act{G}{\topo{X}}$, which factors as $\an\xrightarrow{c_{\topo{X}}} \topo{X}\xrightarrow{\triv}\Act{G}{\topo{X}}$ by \cref{globalsectionsforgobjectsfactors}.
\item Internal group cohomology is not internal cohomology of the terminal object in the topos $\Act{G}{\topo{X}}$ (\cref{definitioncohomologyinatopos}), the latter is the identity $\id_{\stab{\Act{G}{\topo{X}}}}$. 
\item If $\topo{X}$ is cartesian closed, then $\triv$ has a right adjoint by \cref{fixedpointfunctorexistsforcartesianclosed}, and hence internal group cohomology $\intgrpcoh{\topo{X}}(G,-)=(-)^G_{\Sp}$ exists for all groups $G\in\Grp(\topo{X})$ by \cref{geometricmorphismstabilization}. 
    \end{enumerate}\end{rems}
\begin{rem}\label{groupcohomologyisext}
Suppose that $\topo{X}$ is a big topos which admits an exhaustion by hypercomplete, 1-localic topoi, $\Ab(\tau_{\leq 0}\topo{X})$ admits countable coproducts, $\stab{\topo{X}}$ has $\Delta^{\operatorname{op}}$-indexed colimits and the tensor product on $\stab{\topo{X}}$ preserves them in both variables. This holds for example if $\topo{X}=\hypershvacc_S(\mathcal C)$ for a hyperaccessible explicit covering site $(\mathcal C,S)$ with $\mathcal C$ a $1$-category, e.g.\ $\topo{X}=\Cond{}(\an)$. 

By \cref{derivedgmodules}, for $R\in\Alg(\stab{\topo{X}}^{\heart})$ and $G\in\Grp(\tau_{\leq 0}\topo{X})$, \[\LMod{\triv_{\Sp}\algebra{R}}{\stab{\Act{G}{\topo{X}}}}\cong \LMod{\algebra{R}[G]}{\stab{\topo{X}}}\cong \mathcal D(\LMod{\algebra{R}[G]}{\Ab(\tau_{\leq 0}\topo{X})}).\] Denote by $f\colon\mathcal D(\LMod{\algebra{R}[G]}{\Ab(\tau_{\leq 0}\topo{X})})\to\stab{\topo{X}}$ the forget functor. 
    Then \[\grpcoh{\topo{X}}^*(G,f-)\cong \pi_{-*}\map_{\LMod{\algebra{R}[G]}{\stab{\topo{X}}}}(\algebra{R},-)\cong \Ext^*_{\LMod{\algebra{R}[G]} {\Ab(\tau_{\leq 0}\topo{X})}}(\algebra{R},-)\] by \cref{adjunctionsspectrallyenriched}.  
\end{rem}
\begin{rems}\label{remarkgroupcohomologyrefines}
    It is natural to consider the following enhancements of group cohomology: 
\begin{enumerate}
    \item As group cohomology is right adjoint to the stabilization of the constant sheaf functor $c_{\stab{\Act{G}{\topo{X}}}}\colon\Sp\to \stab{\Act{G}{\topo{X}}}$ (\cref{cohomologyglobalsections}), \cref{constantsheafsymmetricmonoidal} determines an enhancement of group cohomology to a lax symmetric monoidal functor \[\Gamma_{\stab{\Act{G}{\topo{X}}}}\colon\stab{\Act{G}{\topo{X}}}^{\otimes}\to\Sp^{\otimes},\] cf. \cref{adjunctionfiberwise}. In particular, for $\algebra{R}\in\Alg(\stab{\topo{X}})$, group cohomology enhances to a functor 
    \[\grpcoh{\topo{X}}(G,-)_R\colon \LMod{\algebra{R}}{\stab{\Act{G}{\topo{X}}}}\to \LMod{\Gamma_{\stab{\Act{G}{\topo{X}}}}\algebra{R}}{\Sp}.\] 
    For $\algebra{R}\in\stab{\topo{X}}$, restriction of scalars along the algebra map $\beta_R\colon\Gamma_{\stab{\topo{X}}}\algebra{R}\to \Gamma_{\stab{\Act{G}{\topo{X}}}}\triv_{\Sp}\algebra{R}$ decribed in \cref{remarkglobalsectionsforgobjectsfactors} yields a functor 
    \[ \grpcoh{\topo{X}}(G,-)_R\coloneqq \beta_R^*\circ \grpcoh{\topo{X}}(G,-)_{\triv_{\Sp}\algebra{R}}\colon \LMod{\algebra{R}}{\stab{\Act{G}{\topo{X}}}}\to \LMod{\Gamma_{\Sp}\algebra{R}}{\Sp}\] enhancing group cohomology. 
    \item Suppose $G$ is a group object in a big topos such that internal group cohomology of $G$ in $\topo{X}$ exists. For every $\algebra{R}\in\Alg(\stab{\topo{X}})$, internal group cohomology enhances to a functor \[ \LMod{\triv_{\Sp}R}{\stab{\Act{G}{\topo{X}}}}\to \LMod{\algebra{R}}{\stab{\topo{X}}}.\] 
    Indeed, the symmetric monoial structure of $\triv_{\Sp}\colon \stab{\topo{X}}\to \stab{\Act{G}{\topo{X}}}$ determines an enhancement of $\intgrpcoh{\topo{X}}(G,-)$ to a lax symmetric monoidal functor. This implies that $\triv_{\Sp}\dashv \intgrpcoh{\topo{X}}(G,-)$ enhances to an adjoint pair
    \[ \triv_{\Sp}^R\colon \LMod{\algebra{R}}{\stab{{\topo{X}}}}\leftrightarrows \LMod{\triv \algebra{R}}{\stab{\Act{G}{\topo{X}}}}\colon \intgrpcoh{\topo{X}}(G,-)_R\] such that 
    \begin{center}
    \begin{tikzcd}
        \LMod{\triv \algebra{R}}{\stab{\Act{G}{\topo{X}}}}\arrow[d]\arrow[rr,"\intgrpcoh{\topo{X}}(G{,}-)_R"']&& \LMod{\algebra{R}}{\stab{\topo{X}}}\arrow[d] & \\ 
        {\stab{\Act{G}{\topo{X}}}}\arrow[rr, "\intgrpcoh{\topo{X}}(G{,}-)"] && {\stab{\topo{X}}} 
    \end{tikzcd}
    \end{center} commutes, cf.\ \cref{localisationinduceslocalisationonmodulecategories}. 
\end{enumerate} 
\end{rems}
\cref{adjunctioninducesadjunctionongobjects} and \cref{geometricmorphismcohomology} imply that group cohomology is invariant under geometric morphisms: 
\begin{cor}
    Suppose that $R\colon \topo{Y}\leftrightarrows \topo{X} \colon L$ is a geometric morphism and $G\in\Grp(\topo{X})$. 
    By \cref{adjunctioninducesadjunctionongobjects}, $K\coloneqq L^{\Alg}G\in in\Grp(\topo{Y})\subseteq \Alg(\topo{Y})$.  Denote by \[ R^{G}\colon \Act{K}{\topo{Y}}\leftrightarrows \Act{G}{\topo{X}}\colon L^G\] the induced adjunction (\cref{adjunctioninducesadjunctionongobjects}).
    
    Then \[\grpcoh{\topo{X}}(G,-)\cong \grpcoh{\topo{Y}}(K,R^G_{\Sp}-).\] 
\end{cor}
In particular, group cohomology in a big topos can always be computed in a subcategory which is a topos: 
\begin{cor}\label{groupcohomologyonfinitelevel}
    Suppose $G$ is a group object in a big topos $\mt{X}_{\infty}$ is a big topos and $\mt{X}_*\colon \Lambda\to \Pr^L$ is an exhaustion of $\mt{X}_{\infty}$ by topoi. 
    Suppose $\lambda\in\Lambda$ is such that $G\in\Grp(\mt{X}_{\lambda})\subseteq \Grp(\mt{X}_{\infty})$, denote by \[r_{\lambda}\colon\mt{X_{\infty}}\to \mt{X_{\lambda}}\] the right adjoint of $\mt{X_{\lambda}}\to\mt{X_{\infty}}$ (\cref{fullyfaithfulnessandpreservationoflimits}), and by $r^G_{\lambda}\colon \lAct{G}{\topo{X}}{\infty}\to \lAct{G}{\topo{X}}{{\lambda}}$ the induced right adjoint (\cref{adjunctioninducesadjunctionongobjects}).  
    Then \[\grpcoh{\mt{X}_{\infty}}(G,-)\cong \grpcoh{\mt{X}_{\lambda}}(G,-)\circ (r_{\lambda}^G)_{\Sp},\] and in particular, 
    \[ \grpcoh{\mt{X}_{\infty}}(G,-)|_{\stab{\lAct{G}{\topo{X}}{\lambda}}}=\grpcoh{\mt{X}_{\lambda}}(G,-).\] 
\end{cor}
If $G$ is a discrete group, group cohomology with coefficients in a trivial $\mathbb Z[G]$-module is isomorphic to the singular/sheaf cohomology of the classifying space $BG$ of numerable principal $G$-bundles. 
This has the following generalization: 
\begin{proposition}\label{grpcohomologyiscohomologyonclassifyingspaces}
    Suppose $G\in\Grp(\topo{X})$ is a group object in a big topos $\topo{X}$ and denote by $\mathbb BG$ the classifying space of principal $\infty$-$G$-bundles in $\topo{X}$ (\cref{definitionclassifyingspaceofgbundles}). 
    The adjunction $-//G\dashv \triv$ induces an equivalence
    \[ \grpcoh{\topo{X}}(G,-)\circ \triv_{\Sp}\cong \cH{\topo{X}}(\mathbb BG,-)\] of functors $\stab{\topo{X}}\to \Sp$. 
\end{proposition}
\begin{rem}
    For a group object $G$ in a big topos $\topo{X}$, $\triv_{\Sp}\colon \stab{\topo{X}}\to\stab{\Act{G}{\topo{X}}}$ is $t$-exact. 
    Indeed: Choose an exhaustion $\topo{X}_*\colon\Lambda\to \Pr^L$ of $\topo{X}$ by topoi such that $G\in\Grp(\topo{X}_{\lambda})$ for all $\lambda\in\Lambda$. By construction, $\triv_{\Sp}$ restricts to $\triv^{\lambda}_{\Sp}\colon \stab{\topo{X}_{\lambda}}\to \stab{\lAct{G}{\topo{X}}{\lambda}}$ for all $\lambda\in\Lambda$ (the stabilization of restriction of scalars along $G\to *\in\Grp(\topo{X}_{\lambda})$), and since $\triv^{\lambda}$ admits a right adjoint, $\triv^{\lambda}_{\Sp}$ is $t$-exact by \cref{texactnessgeometricmorphism}. By construction of the $t$-structure on $\stab{\topo{X}}$ (\cref{tstructurespectrumobjects}), this implies that $\triv_{\Sp}$ is $t$-exact. In particular, if $G\in \Grp(\tau_{\leq 0}\topo{X})$, then $\triv_{\Sp}$ restricts to a functor \[\triv_{\Sp}^{\heart}\colon \Ab(\tau_{\leq 0}\topo{X})=\stab{\topo{X}}^{\heart}\to \stab{\Act{G}{\topo{X}}}^{\heart}\cong \LMod{\mathbb Z[G]}{\Ab(\tau_{\leq 0}\topo{X})}\] which endows an abelian group object $M\in\Ab(\tau_{\leq 0}\topo{X})$ with the trivial $G$-action.  
\end{rem}
\begin{proof}[Proof of \cref{grpcohomologyiscohomologyonclassifyingspaces}]
    Choose an exhaustion $\topo{X}_*\colon\Lambda\to \Pr^L$ by topoi and $\lambda\in\Lambda$ with \[G\in\Grp(\topo{X}_{\lambda})\subseteq \Grp(\topo{X}),\] and denote by $r^{\lambda,(G)}\colon \Act{(G)}{\topo{X}}\to \Act{(G)}{\topo{X}}$ the right adjoint of $\lAct{(G)}{\topo{X}}{\lambda}\hookrightarrow \Act{(G)}{\topo{X}}$. Since $\topo{X}_{\lambda}\subseteq\topo{X}$ is closed under colimits and finite limits, $\mathbb BG\in\topo{X}_{\lambda}$ is the classifying space of $G$ in $\topo{X}_{\lambda}$. 
By construction, $\triv$ restricts to $\triv^{\lambda}\colon\topo{X}\to \lAct{G}{\topo{X}}{\lambda}$ (restriction of scalars along $G\to *\in\Grp(\topo{X}_{\lambda})$), and \begin{align}\label{restrictiontrivialisationlambda} r^{\lambda,G}\circ \triv=\triv^{\lambda}\circ\ r^{\lambda},\end{align} cf.\ \cref{naturalitylocalizationmodulecategories}.
By \cref{existenceleftadjointmoregeneral}, $\triv^{\lambda}_{\Sp}\colon \stab{\topo{X}_{\lambda}}\to \stab{\lAct{G}{\topo{X}}{\lambda}}$ admits a left adjoint $L$ with \[L\circ \Sigma^{\infty}_{+, \stab{\lAct{G}{\topo{X}}{\lambda}}}\cong \Sigma^{\infty}_{+, \topo{X}_{\lambda}}\circ (-//G), \] where $-//G\colon\lAct{G}{\topo{X}}{\lambda}\to \topo{X}_{\lambda}$ denotes the homotopy orbits functor. \cref{adjunctionsspectrallyenriched} implies that \[\cH{\topo{X}_{\lambda}}(\mathbb BG,-)\cong \map_{\stab{\topo{X}_{\lambda}}}(L\Sigma^{\infty}_{+, \stab{\Act{G}{\topo{X}}}}(*),-)\cong \grpcoh{\topo{X}_{\lambda}}(G,-)\circ \triv^{\lambda}_{\Sp}.\] 
By \ref{restrictiontrivialisationlambda}, \[r^{\lambda,G}_{\Sp}\circ \triv_{\Sp}=\triv^{\lambda}_{\Sp}\circ\ r^{\lambda}_{\Sp}, \] 
and hence \[\cH{\topo{X}_{\lambda}}(\mathbb BG,-)\circ r^{\lambda}_{\Sp}\cong \grpcoh{\topo{X}_{\lambda}}(G,-)\circ \triv^{\lambda}_{\Sp}\circ r^{\lambda}_{\Sp}\cong \grpcoh{\topo{X}_{\lambda}}(G,-)\circ r^{\lambda,G}_{\Sp}\circ \triv_{\Sp}.\]
By \cref{groupcohomologyonfinitelevel} and \cref{cohomologycanbecomputedonfinitelelvel}, \[\cH{\topo{X}}(\mathbb BG,-)\cong \cH{\topo{X}_{\lambda}}(\mathbb BG,-)\circ r^{\lambda}_{\Sp}\] and 
\[\grpcoh{\topo{X}}(G,-)\cong \grpcoh{\topo{X}_{\lambda}}(G,-)\circ r^{\lambda,G}_{\Sp} .\qedhere\]
\end{proof}
If $G\in\Grp(\topo{X})$ is a group object in a big topos $\topo{X}$ and $E\to \triv B$ is a principal $\infty$-$G$-bundle in $\topo{X}$, the classifying map $c\colon B\cong E//G\to \mathbb BG$ (\cref{characterizationbundles}) induces a natural transformation
\[c^*\colon \grpcoh{\topo{X}}(G, \triv_{\Sp}-)\cong \cH{\topo{X}}(\mathbb BG,-)\to \cH{\topo{X}}(B,-)\] of functors $\stab{\topo{X}}\to \Sp$. 
\begin{cor}\label{groupcohomologyofbundles}
    Suppose $G\in \Grp(\topo{X})$ is a group object in a big topos $\topo{X}$ and $E\to \triv B$ is a principal $\infty$-$G$-bundle. 
    Suppose that $M\in \stab{\topo{X}}$ is such that for all $p\in\mathbb N_0$, the projection $G^p\times \res{G}{*}E\to G^p$ induces an equivalence \[\cH{\topo{X}}(G^p,M)\cong \cH{\topo{X}}(G^p\times \res{G}{*}E,M).\] 
    Then pullback along the classifying map $B\to\mathbb BG$ is an equivalence \[ \grpcoh{\mathcal X}(G, \triv_{\Sp}(M))\cong \cH{\topo{X}}(B,M).\]
\end{cor}
\begin{proof}
Denote by $E_*$ the action groupoid of $E$, i.e.\ the image of 
$E$ under \[\Act{G}{\topo{X}}\hookrightarrow \oc{\Fun(\Delta^{\operatorname{op}}, \topo{X})}{*//G}\to\Fun(\Delta^{\operatorname{op}}, \topo{X}).\]  
By construction of the identification $\Act{G}{\topo{X}}\cong \oc{\topo{X}}{\mathbb BG}$ and \cref{spectralenrichmentcocontinuous}, 
$c^*$ is the totalization of the morphism of cosimplicial objects
\[ \cH{\topo{X}}((*//G)_*,-)\to \cH{\topo{X}}(E_*,-)\in \Fun(\Delta, \Fun(\stab{\topo{X}}, \Sp))\] induced by the map $E_*\to (*//G)_*$. 
This implies the statement as for $[n]\in\Delta^{\operatorname{op}}$, \[E([n])\cong (*//G)([n])\times E([0])\cong G^n\times \res{G}{*}E\] and $E([n])\to (*//G)([n])$ is the projection.   
\end{proof}

For a group object $G$ in a big topos $\topo{X}$, we obtain the following spectral sequence by applying \cref{Bousfieldkanspectralsequencehomology} to the effective epimorphism $G\to *\in\Act{G}{\topo{X}}$:
\begin{cor}\label{chechtocohomologyspectralsequencegroupcohomology}
    Suppose $G$ is a group object in a big topos $\mathcal X$. There is a spectral sequence of functors $\stab{\Act{G}{\topo{X}}}\to \Ab$ with $E_1$-page \[E_1^{p,q}=\cH{\topo{X}}^q(G^p, \res{G}{*, \Sp}-).\] 
    On bounded above objects $\bigcup_{n\in\mathbb N_0}\LMod{\mathbb S[G]}{\stab{\topo{X}}}_{\leq n}$, this converges completely to $\grpcoh{\topo{X}}^{p+q}(G,-)$.
\end{cor}
\begin{proof}
    By \cref{Bousfieldkanspectralsequencehomology}, there is a spectral sequence with $E_1$-page \[ E_1^{p,q}=\cH{\Act{G}{\topo{X}}}^q(G^{p+1},-)\] which converges completely to $\grpcoh{\topo{X}}^{p+q}(G,-)\cong \cH{\Act{G}{\topo{X}}}(*,-)$ on bounded above spectrum objects. 
    For all $m\in\mathbb N_0$, the map $G^m\to G^{m+1}, g\mapsto (1,g)$ is adjoint to an equivalence $\ind{G}{*}(G^m)\cong G^{m+1}$. Hence by \cref{indgspisstabilization}, \[\Sigma^{\infty}_{+, \Act{G}{\topo{X}}}\check{C}(G\to *)([m])\cong \Sigma^{\infty}_{+,\Act{G}{\topo{X}}}G^{m+1}\cong \Sigma^{\infty}_{+,\Act{G}{\topo{X}}}\ind{G}{*}(G^m)\cong\ind{G}{*, \Sp}\Sigma^{\infty}_{+, \topo{X}}G^m\] for all $m\in\mathbb N_0$. It follows from \cref{adjunctionsspectrallyenriched} that for all $p,q\in\mathbb Z$,  
    \[  \cH{\Act{G}{\topo{X}}}^q(G^{p+1},-)\cong \cH{\topo{X}}^q(G^p,\res{G}{*,\Sp}-).\qedhere\]
\end{proof}
\begin{rem}\label{edgehomomorphismgroupcohomology}
For $M\in \stab{\Act{G}{\topo{X}}}^{\heart}\cong \LMod{\mathbb Z[\pi_0^{\topo{X}}G]}{\Ab(\tau_{\leq 0}\topo{X})}$, the spectral sequence (\cref{chechtocohomologyspectralsequencegroupcohomology}) is concentrated in the upper right quadrant, whence we obtain an edge homomorphism \[E_2^{*,0}(M)\twoheadrightarrow E_2^{*, \infty}(M)\hookrightarrow \grpcoh{\topo{X}}^*(G,M).\]  
Denote by $S_*^{G}\in \Ch(\LMod{\mathbb Z[\pi_0^{\topo{X}}G]}{\Ab(\tau_{\leq 0}\topo{X})})$ the simplicial resolution of $\pi_0^{\topo{X}}G\to *\in\tau_{\leq 0}\Act{G}{\topo{X}}$ (\cref{definitionsimplicialcomplex}). Then \[E_1^{*,0}|_{\LMod{\mathbb Z[\pi_0^{\topo{X}}G]}{\Ab(\tau_{\leq 0}\topo{X})}}\cong \Hom_{\LMod{\mathbb Z[\pi_0^{\topo{X}}G]}{\Ab(\tau_{\leq 0}\topo{X})}}(S_*^{G},-), \] i.e.\ the edge homomorphism is a natural transformation \[E_2^{*,0}|_{\LMod{\mathbb Z[\pi_0^{\topo{X}}G]}{\Ab(\tau_{\leq 0}\topo{X})}}\cong H^*(\Hom_{\LMod{\mathbb Z[\pi_0^{\topo{X}}G]}{\Ab(\tau_{\leq 0}\topo{X})}}(S_*^{G},-))\Rightarrow \grpcoh{\topo{X}}^*(G,-)|_{\LMod{\mathbb Z[\pi_0^{\topo{X}}G]}{\Ab(\tau_{\leq 0}\topo{X})}}\] of functors $\LMod{\mathbb Z[\pi_0^{\topo{X}}G]}{\Ab(\tau_{\leq 0}\topo{X})}\to \grAb$. 
\end{rem}
We will apply this below to compare condensed with continuous group cohomology.
\subsubsection{Internal Group Cohomology}\label{section:internalgroupcohomology}
We now record some basic results on internal group cohomology. 
\begin{lemma}\label{internalgroupcohomologytotalisaton}
    Suppose that $\topo{X}$ is a big topos such that $\stab{\topo{X}}$ has all $\Delta$-indexed limits and all $\Delta^{\operatorname{op}}$-indexed colimits and the symmetric monoidal structure on $\stab{\topo{X}}$ is closed. 
    
    For all group objects $G\in\Grp(\topo{X})$, internal group cohomology $\intgrpcoh{\topo{X}}(G,-)$ exists and   
    \[\intgrpcoh{\topo{X}}(G,-)\cong \clim{ p\in \Delta}\icH{\topo{X}}(G^p,-)\circ \res{G}{*, \Sp}(-).\]   
\end{lemma}
\begin{rem}The conditions of the lemma are satisfied if $\mathcal X$ is a topos or the category of (hyper)ac\-ces\-si\-ble sheaves on a (hyper)ac\-ces\-si\-ble explicit covering site (\cref{symmetricmonoidalstructureaccessiblesheavesofspectraclosed}).\end{rem}  
\begin{proof}
    Recall from \cref{GobjectsstabilisetoGmodules} that $\stab{\Act{G}{\topo{X}}}\cong \LMod{\mathbb S[G]}{\stab{\topo{X}}}$. 
As $\stab{\topo{X}}$ has all $\Delta^{\operatorname{op}}$-indexed colimits and its monoidal structure is closed, by \cref{forgetfreeadjunctionmodules} $\stab{\Act{G}{\topo{X}}}$ has all $\Delta^{\operatorname{op}}$-indexed colimits and the forget functor \[\res{G}{*, \Sp}\colon \stab{\Act{G}{\topo{X}}}\to\stab{\topo{X}}\] reflects them. 
Since $\res{G}{*, \Sp}$ enhances to a symmetric monoidal functor, it follows that the tensor product on $\stab{\Act{G}{\topo{X}}}$ preserves $\Delta^{\operatorname{op}}$-indexed colimits in both variables.
In particular, \[\triv_{\Sp}(-)\cong \triv_{\Sp}(-)\otimes \colim{\Delta^{\operatorname{op}}}\Sigma^{\infty}_{+, \Act{G}{\topo{X}}}\check{C}(G\to *)\cong \colim{\Delta^{\operatorname{op}}} (\triv_{\Sp}(-)\otimes \Sigma^{\infty}_{+, \Act{G}{\topo{X}}}\check{C}(G\to *)),\] where we used that $\colim{\Delta^{\operatorname{op}}}\Sigma^{\infty}_{+, \Act{G}{\topo{X}}}\check{C}(G\to *)\cong \Sigma^{\infty}_{+, \Act{G}{\topo{X}}}(*)$ is the monoidal unit of $\stab{\Act{G}{\topo{X}}}$. By \cref{indgspisstabilization} and \cref{inductiontrivialstable}, for  $[n]\in\Delta^{\operatorname{op}}$, \begin{align*}\triv_{\Sp}(-)\otimes_{\Sp(\Act{G}{\topo{X}})} \Sigma^{\infty}_{+, \Act{G}{\topo{X}}}\check{C}(G\to *)([n])&\cong \triv_{\Sp}(-)\otimes_{\Sp(\Act{G}{\topo{X}})}\Sigma^{\infty}_{+, \Act{G}{\topo{X}}}G^{n+1}\\& \cong \triv_{\Sp}(-)\otimes_{\Sp(\Act{G}{\topo{X}})} \ind{G}{*, \Sp}(\Sigma^{\infty}_{+, \topo{X}}G^n)\\  &\cong \ind{G}{*, \Sp}(-\otimes_{\stab{\topo{X}}}\Sigma^{\infty}_{+, \topo{X}}G^n)\end{align*} is left adjoint to \[\imap_{\stab{\topo{X}}}(\Sigma^{\infty}_{+, \topo{X}}G^n,-)\circ \res{G}{*, \Sp}=\icH{\topo{X}}(G^n, \res{G}{*, \Sp}-).\]
Since $\stab{\topo{X}}$ has all $\Delta^{\operatorname{op}}$-indexed limits and $\stab{\Act{G}{\topo{X}}}$ has all $\Delta^{\operatorname{op}}$-indexed colimits, 
\[ \Fun^{L}(\stab{\topo{X}}, \stab{\Act{G}{\topo{X}}})\subseteq \Fun(\stab{\topo{X}}, \stab{\topo{X}})\] is closed under $\Delta^{\operatorname{op}}$-indexed colimits and they are computed pointwise (\cref{leftadjointsstableundercolimits}). Analogously, \[ \Fun^{R}(\stab{\Act{G}{\topo{X}}}, \stab{\topo{X}})\subseteq \Fun(\stab{\Act{G}{\topo{X}}}, \stab{\topo{X}})\] is closed under $\Delta$-indexed limits and they can be computed pointwise. 
\cite[Proposition 5.2.6.2]{highertopostheory} now implies that $\clim{n\in \Delta}\icH{\topo{X}}(G^n,-)\circ \res{G}{*, \Sp}$ is right adjoint to $\triv_{\Sp}(-)$. 
\end{proof}

Suppose $G$ is a group object in a big topos $\topo{X}$ such that internal group cohomology $\intgrpcoh{\topo{X}}(G,-)$ exists. As $\grpcoh{\topo{X}}(G,-)\cong \Gamma_{\Sp}\circ \intgrpcoh{\topo{X}}(G,-)$ (\cref{groupcohomologyinternalrecovers}), the adjunction counit of the stabilization of the global sections geometric morphism $c_{\Sp}(-)\colon \Sp\rightleftarrows \stab{\topo{X}}\colon  \Gamma_{\Sp}$ provides a natural transformation \[c_{\Sp}\left(\grpcoh{\topo{X}}(G,-)\right)\to \intgrpcoh{\topo{X}}(G,-).\] 
In rare cases, this is an equivalence: 
\begin{cor}\label{internalgroupcohomologyconstant}
    Suppose that $\topo{X}$ is a big topos such that $\stab{\topo{X}}$ has all countable products and all $\Delta^{\operatorname{op}}$-indexed colimits and the symmetric monoidal structure on $\stab{\topo{X}}$ is closed. 
    Suppose that $G\in\Grp(\topo{X})$ and $M\in\LMod{\mathbb S[G]}{\stab{\topo{X}}}$ is such that for all $p\in\mathbb N_0$, 
    \[c_{\Sp}\left(\cH{\topo{X}}(G^p, \res{G}{*,\Sp}M)\right)\cong \icH{\topo{X}}(G^p, \res{G}{*,\Sp}M)\] via the counit of $c(-)\dashv \Gamma_{\Sp}$. 
    Then $c_{\Sp}\left(\grpcoh{\topo{X}}(G,M)\right)\cong \intgrpcoh{\topo{X}}(G,M)$ via the counit. 
\end{cor}
\begin{proof}As $\stab{\topo{X}}$ has countable products (by assumption) and finite limits (\cref{spectrumobjectsinbigcatsstable}), $\stab{\topo{X}}$ has all countable limits (\cite[Proposition 4.4.2.6]{highertopostheory}), and in particular all $\Delta$-indexed limits. By \cref{internalgroupcohomologytotalisaton}, 
    \[ \intgrpcoh{\topo{X}}(G,-)\cong \clim{p\in\Delta}\icH{\topo{X}}(G^p,\res{G}{*,\Sp}-),\] 
    whence \[\grpcoh{\topo{X}}(G,-)\cong\Gamma_{\Sp}\intgrpcoh{\topo{X}}(G,-)\cong \clim{p\in\Delta}\Gamma_{\Sp}\icH{\topo{X}}(G^p,\res{G}{*,\Sp}-)\cong\clim{p\in\Delta}\cH{\topo{X}}(G^p,-).\]  
    The statement now follows from \cref{constantsheafandtotalization}.
\end{proof}

\subsubsection{Condensed Group Cohomology}\label{section:condensedgroupcohomology}
We now specialize the above discussion to condensed group cohomology. 
\begin{definition}
    For a ($\kappa$-)condensed group $G\in \Grp(\Cond{(\kappa)}(\an))$, denote by 
    \[ \condgrpcoh{(\kappa\text{-})}(G,-)\coloneqq \grpcoh{\Cond{(\kappa)}(\an)}(G,-)\] its ($\kappa$-condensed) group cohomology, and by 
    \[ \icondgrpcoh{(\kappa\text{-})}(G,-)\coloneqq \intgrpcoh{\Cond{(\kappa)}(\an)}(G,-)\] the internal ($\kappa$-)condensed group cohomology. 
    The latter exists by \cref{groupcohomologyinternalrecovers}. 
\end{definition}
\begin{rem}\label{condensedgroupcohomologyisext}For a discrete ($\kappa$-)condensed group $G\in\Grp(\Cond{(\kappa)}(\Set))$,  \[\LMod{H\mathbb Z}{\stab{\Act{G}{\Cond{(\kappa)}(\an)}}}\cong \mathcal D(\Cond{(\kappa)}(\mathbb Z[G])) (\text{ \cref{derivedgmodules}}).\] Denote by $f\colon\mathcal D(\Cond{(\kappa)}(\mathbb Z[G]))\to\stab{\Act{G}{\Cond{(\kappa)}(\an)}}$ the forget functor. Then \[ \condgrpcoh{(\kappa\text{-})}^*(G,f-)\cong \Ext^*_{\Cond{(\kappa)}(\mathbb Z[G])}(\mathbb Z,-) \text{ (\cref{groupcohomologyisext})}.\] 

    Moreover, ($\kappa$-)condensed group cohomology   
    $\condgrpcoh{(\kappa\text{-})}(G,-)$ refines to a functor 
    \[\condgrpcoh{(\kappa\text{-})}(G,-)_{\mathbb Z}\colon \mathcal D(\Cond{(\kappa)}(\LMod{\mathbb Z[G]}{\Cond{(\kappa)}(\Ab)}))\to \mathcal D(\Ab), \] and internal ($\kappa$-)condensed group cohomology enhances to a functor \[\mathcal D(\Cond{(\kappa)}(\mathbb Z[G]))\to \mathcal D(\Cond{(\kappa)}(\Ab)),\text{ cf. \cref{remarkgroupcohomologyrefines}}.\] 
\end{rem}
\cref{groupcohomologyonfinitelevel} implies that condensed group cohomology can always be computed on some level \[\stab{\Act{G}{\Cond{\kappa}(\an)}}\subseteq \stab{\Act{G}{\Cond{}(\an)}}:\] 
\begin{cor}Suppose $G$ is a Hausdorff topological group. 
    \begin{romanenum}
    \item If $\kappa$ is an uncountable regular or strong limit cardinal with $\cof(\kappa)>|G|$, then \[ \condgrpcoh{}(\underline{G},-)|_{\Sp(\Act{\underline{G}_{\kappa}}{\Cond{\kappa}(\an)})}\cong \condgrpcoh{\kappa\text{-}}(\underline{G}_{\kappa},-).\] 
    \item In particular, if $M$ is a $\To$ continuous $G$-module, for all uncountable regular/strong limit cardinals $\kappa$ with $\cof(\kappa)>\max\{|M|,|G|\}$, 
    \[\condgrpcoh{}(\underline{G}, \underline{M})\cong \condgrpcoh{\kappa\text{-}}(\underline{G}_{\kappa}, \underline{M}_{\kappa}).\] 
    \end{romanenum}
\end{cor}
\begin{proof}
    By \cref{topspacecomefromregular}, for all regular cardinals and all strong limit cardinals $\kappa$ with $\cof(\kappa)>|G|$, \[\underline{G}=\underline{G}_{\kappa}\in \Cond{\kappa}(\Set)\subseteq \Cond{}(\Set).\] The statement now follows from \cref{groupcohomologyonfinitelevel} and \cref{condensedfiltrationinfinitytopoi} since
    \[ \Cond{\kappa}(\Grp)\cong \Cond{}(\Grp)\times_{\Cond{}(\Set)}\Cond{\kappa}(\Set)\] and 
    \[\Cond{\kappa}(\Ab)\cong \Cond{}(\Ab)\times_{\Cond{}(\Set)}\Cond{\kappa}(\Set).\qedhere\] 
\end{proof}
\cref{internalgroupcohomologyconstant} implies that internal ($\kappa$-)condensed group cohomology is constant in some cases:
    Denote by $c_{\Sp}\colon \Sp\to \Cond{(\kappa)}(\Sp)$ the stabilization of the constant sheaf functor. 
\begin{cor}\label{internalgroupcohomologyconstantcondensed}
    Suppose $G\in \Grp(\Cond{(\kappa)}(\Set))$ and there exists a ($\kappa$-light) compact Hausdorff space $K$ with $\mathbb Z[G]^{L\solid}\cong \mathbb Z[\underline{K}_{(\kappa)}]^{L\solid}$ in $\mathcal D(\Sol{(\kappa)})$.  
    If $A\in\Cond{(\kappa)}(\mathbb Z[G])$ is such that the condensed abelian group underlying $A$ is represented by a discrete abelian group, then \[c_{\Sp}\left(\condgrpcoh{(\kappa\text{-})}(G,A)\right)\cong \icondgrpcoh{(\kappa\text{-})}(G,A).\]
\end{cor}
We prove this below. Together with \cref{solidcohomologyhomotopyinvariant}, it implies the following more concrete statement:  
\begin{cor}\label{solidcohomologyhomotopyinvariantmoreconcrete} Suppose $G$ is a topological group. 
\begin{romanenum}
\item If $G$ is homotopy equivalent to a $\kappa$-light compact Hausdorff space, for all discrete continuous $G$-modules $A$, \[c_{\Sp}\left(\condgrpcoh{\kappa\text{-}}(\underline{G}_{\kappa},\underline{A}_{\kappa})\right)\cong \icondgrpcoh{\kappa\text{-}}(\underline{G}_{\kappa},\underline{A}_{\kappa}).\]
\item If $G$ is Hausdorff and homotopy equivalent to a compact Hausdorff space, for all discrete continuous $G$-modules $A$,  \[c_{\Sp}\left(\condgrpcoh{}(\underline{G},\underline{A})\right)\cong \icondgrpcoh{}(\underline{G},\underline{A}).\]
\end{romanenum}
\end{cor}
\begin{rem}For profinite groups, \cref{solidcohomologyhomotopyinvariantmoreconcrete} was established in \cite[Lemma 2.5]{Anschuetzsolidhomology}.  

\cref{solidcohomologyhomotopyinvariantmoreconcrete} in particular applies to ($\kappa$-light) locally compact Hausdorff abelian groups with finitely many connected components: 
By the structure theorem for locally compact abelian groups (see e.g.\ \cite[Theorem 4]{HofmannMorris+2023}), every locally compact Hausdorff abelian group $A$ is isomorphic to a product $\mathbb R^n\times B$, where $B$ admits a compact open subgroup $K\subseteq B$.
The space $B$ is homeomorphic to $\sqcup_{\pi_0B}K$ and hence compact of weight $\leq \wt(A)$ and $A$ is homotopy equivalent to $B$. 

Internal condensed group cohomology has a non-constant condensed structure in many cases. For example, for a continuous $G$-module $M$, \[\icondgrpcoh{(\kappa\text{-})}^0(G,M)=\underline{M^G}_{(\kappa)},\] is represented by the $G$-fixed-points of $M$ (with subspace topology), whereas \[ c_{\Sp}\left(\condgrpcoh{\kappa\text{-}}^0(\underline{G}_{\kappa},M)\right)=\underline{(M^{G})^{\delta}}_{\kappa}\] is represented by the discrete abelian group $(M^{G})^{\delta}$. 
\end{rem}   
\begin{proof}[Proof of \cref{internalgroupcohomologyconstantcondensed}]\label{proofinternalgroupcohomologyconstantcondensed}
    By assumption, there exists a ($\kappa$-light) compact Hausdorff space $K$ with $\mathbb Z[G]^{L\solid}\cong \mathbb Z[\underline{K}_{(\kappa)}]^{L\solid}$ in $\mathcal D(\Sol{(\kappa)})$. 
    This implies that for $p\in\mathbb N_0$, 
    \[\mathbb Z[{G}^p]^{L\solid}\cong \bigotimes^{L\solid}_{i=1, \ldots,p}(\mathbb Z[G]^{L\solid})\cong \mathbb Z[\underline{K^p}_{(\kappa)}]^{L\solid}.\] 
    Hence by \cref{internalcohomologydiscrete}, for all $p\in\mathbb N_0$ and all discrete abelian groups $A$, \[c_{\Sp}\left(\cckH{(\kappa\text{-})}(G^{p}, \underline{A}_{(\kappa)})\right)\cong \icckH{(\kappa\text{-})}(G^{p}, \underline{A}_{(\kappa)})\] via the adjunction counit for $c_{\Sp}\dashv \Gamma_{\Sp}$. 
    The statement now follows from \cref{internalgroupcohomologyconstant} 
    since $\Cond{(\kappa)}(\Sp)$ has small limits and colimits (\cite[Corollary 5.5.2.4]{highertopostheory}/\cref{condensedextremallydisconnected1}) and the symmetric monoidal structure on $\Cond{(\kappa)}(\Sp)$ is closed (\cref{closedmonoidalstructurecondensedcategories}).
\end{proof}
\subsection{Condensed and continuous group cohomology}
We now compare ($\kappa$-)condensed group cohomology of a (T1) topological group $G$ with its continuous group cohomology. 
    \begin{definition}    
    For a discrete group $G\in\Grp(\tau_{\leq 0}\topo{X})$ in a big topos $\topo{X}$ denote by \[S_*^{G}\to \mathbb Z \in \Ch(\LMod{\mathbb Z[G]}{\Ab(\tau_{\leq 0}\topo{X})})\] the simplicial resolution of $G\to *\in \tau_{\leq 0}\topo{X}$ (\cref{definitionsimplicialcomplex}).

        For a group $G\in\Grp(\Set)$ and a $G$-module $M\in \LMod{\mathbb Z[G]}{\Ab}$, the cochain complex 
        \[\Hom_{\LMod{\mathbb Z[G]}{\Ab}}(S_*^{G},M)\in \Ch(\LMod{\mathbb Z[G]}{\Ab})\] obtained by applying $\Hom_{\mathbb Z[G]}(-,M)$ to the simplicial resolution $S_*^{G}\in \Ch(\LMod{\mathbb Z[G]}{\Ab})$ componentwise is called  the \emph{complex of} $G$-\emph{cochains in} $M$.  
    \end{definition}

    Explicitly, for $i\in\mathbb N_0$, \[\Hom_{\LMod{\mathbb Z[G]}{\Ab}}(S_i^{G},M)=\Map_{\Set}(G^i,M)\] and for $i\in\mathbb N_1$ and $f\in \Map_{\Set}(G^{i-1},M)$, \begin{align*}\partial_{i-1}f\colon G^{i}&\to M, \\ (g_1, \ldots,g_{i})&\mapsto g_1f(g_2, \ldots,g_{i})+\sum_{j=1}^{i-1} (-1)^j f(g_1, \ldots, g_{j-1},g_jg_{j+1},g_{j+2}\ldots,g_i)+(-1)^{i}f(g_1,g_2, \ldots,g_{i-1}).\end{align*}

    \begin{definition}
    If $G$ is a topological group, a \emph{continuous} $G$-module $M$ is a topological abelian group $M$ with a continuous $G$-action $G\times M \to M$ by linear maps. 
    Denote by $\Cont{G}$ the category of continuous $G$-modules with $G$-equivariant, continuous linear maps.
\end{definition}

\begin{lemma}\label{inclusioncontinuouscochains}Suppose that $G$ is a topological group and $M$ is a continuous $G$-module.
The continuous/$k$-continuous/$\kappa$-continuous maps $G^*\to M$ define subcomplexes 
\[\mathcal C(G^*,M)\subseteq \mathcal C_k(G^*,Y)\subseteq \mathcal C_{\kappa}(G^*,M)\subseteq \Hom_{\LMod{\mathbb Z[G]}{\Ab}}(S_*^{G},M)\] of the complex of $G$-cochains in $M$. 
\end{lemma}
\begin{proof}This follows from continuity of the $G$-action on $M$ and the group multiplication of $G$. \end{proof}
\begin{definition}\label{definitioncontinuousgroupcohomology}
    If $G$ is a topological group and $M$ is a continuous $G$-module, we denote by \[\mathcal C(G^*,M)\subseteq \mathcal C_k(G^*,Y)\subseteq \mathcal C_{\kappa}(G^*,M)\] the complexes of \emph{continuous}/\emph{k-continuous}/$\kappa$-\emph{continuous cochains} of $G$ in $M$. 
    Their cohomology groups are called the \emph{continuous}/\emph{k-continuous}/$\kappa$-\emph{continuous group cohomology} of $G$ with coefficients in $M$, and are denoted by
    \begin{align*}
    \contgrpcoh{}^*(G,M)&\coloneqq H^*(\cont(G^*,M)), \\ \contgrpcoh{k\text{-}}^*(G,M)&\coloneqq H^*(\cont_k(G^*,M)), \\
    \contgrpcoh{\kappa\text{-}}^*(G^*,M)&\coloneqq H^*(\cont_{\kappa}(G^*,M)),
    \end{align*} respectively. 
\end{definition}
\cref{inclusioncontinuouscochains} yields natural transformations\[ \contgrpcoh{}^*(G,-)\to \contgrpcoh{k\text{-}}^*(G,-)\to \contgrpcoh{\kappa\text{-}}^*(G,-)\in \Fun(\Cont{G}, \grAb).\]
\cref{underlinefunctor} and \cref{kappacontinuousfullyfaithfullyintocondensed} imply: 
\begin{cor}\label{continuouscohomologysimplicialcomplex}
Suppose that $G$ is a topological group and $M$ is a continuous $G$-module. 

\begin{romanenum} \item For all uncountable cardinals $\kappa$,  
\[\cont_{\kappa}(G^*,M)\cong \Hom_{\Cond{\kappa}(\mathbb Z[\underline{G}_{\kappa}])}(S_*^{\underline{G}_{\kappa}}, \underline{M}_{\kappa}).\] 
\item If $G$ and $M$ are $\To$, then 
\[\cont_{k}(G^*,M)\cong \Hom_{\Cond{}(\mathbb Z[\underline{G}])}(S_*^{\underline{G}}, \underline{M}).\]
\end{romanenum}
\end{cor}
This implies that the \v{C}ech-to-cohomology spectral sequence for the cover $\underline{G}_{(\kappa)}\to *$ (\cref{chechtocohomologyspectralsequencegroupcohomology}) interpolates between $\kappa/k$-continuous and $(\kappa)$-condensed group cohomology, cf.\ \cref{edgehomomorphismgroupcohomology}: 
\begin{cor}\label{spectralsequenceidentifycontinuousgroupcohomology}

    \begin{romanenum}
    \item Suppose $G$ is a Hausdorff topological group. 
    There is an upper quadrant spectral sequence of functors $\TCont{G}\to \Ab$ with \[E_{1}^{p,q}=\ccH^q(\underline{G^p},-)\text{ and } E_2^{*,0}=\contgrpcoh{k\text{-}}^*(G,-),\] which converges to $\condgrpcoh{}^*(\underline G, \underline{-})$. The edge homomorphism is a natural transformation 
    \[ \contgrpcoh{k\text{-}}^*(G,-)\to \condgrpcoh{}^*(\underline G, \underline{-})\in \Fun(\TCont{G}, \grAb).\] 

    \item Suppose $G$ is a topological group and $\kappa$ is an uncountable cardinal. 
    There is an upper quadrant spectral sequence of functors $\Cont{G}\to \Ab$ with \[ E_{1}^{p,q}=\ckH^q(\underline{G^p}_{\kappa},-)\text{ and }E_2^{*,0}=\contgrpcoh{\kappa\text{-}}^*(G,-)\] which converges to $\condgrpcoh{\kappa\text{-}}^*(\underline G_{\kappa}, \underline{-}_{\kappa})$.
    The edge homomorphism is a natural transformation
    \[ \contgrpcoh{\kappa\text{-}}^*(G,-)\to \condgrpcoh{\kappa\text{-}}^*(\underline G_{\kappa}, \underline{-}_{\kappa})\in \Fun(\Cont{G}, \grAb).\] 
    \end{romanenum} 
\end{cor}
\begin{cor}[{\cite[Lemma 2.1]{Anschuetzsolidhomology}}]\label{condensedequalscontinuousprofinite}

    \begin{romanenum}
    \item Suppose $G$ is a Hausdorff topological group which is homotopy equivalent to a coproduct of profinite sets. 
    If $M$ is a $\To$ continuous $G$-module such that $\underline{M}$ is solid, then the spectral sequence from \cref{chechtocohomologyspectralsequencegroupcohomology} collapses at the $E_2$-page and yields \[\contgrpcoh{k\text{-}}(G,M)\cong \condgrpcoh{}^*(\underline{G}, \underline{M}) .\] 
    \item 
    Suppose that $G$ is a topological group which is homotopy equivalent to a coproduct of $\kappa$-light profinite sets. If $M$ is a continuous $G$-module such that $\underline M_{\kappa}\in \Sol{\kappa}$, the spectral sequence from \cref{chechtocohomologyspectralsequencegroupcohomology} collapses at the $E_2$-page and yields  
 \[\contgrpcoh{\kappa\text{-}}^*(G,M)\cong \condgrpcoh{\kappa\text{-}}^*(\underline{G}_{\kappa}, \underline{M}_{\kappa}).\] 
    \end{romanenum}
\end{cor}
\begin{ex}
By van Dantzig's theorem, the above corollary in particular applies to ($\kappa$-light) locally profinite groups. 
\end{ex}
\begin{proof}
If $G\to \sqcup_{i\in I}X(i)$ is a homotopy equivalence with $X(i)\in\Pro(\Fin)_{(\kappa)}$ for all $i\in I$, for all $p\in\mathbb N_1$, 
\[G^p\to (\sqcup_{i\in I}X(i))^p=\sqcup_{(i_1, \ldots,i_p)\in I^p}X({i_1})\times X({i_2})\times \ldots \times X({i_p})\] is a homotopy equivalence and $X({i_1})\times X({i_2})\times \ldots \times X({i_p})\in \Pro(\Fin)_{(\kappa)}$. 
By \cref{solidcohomologyhomotopyinvariant}, this induces an isomorphism \[\mathbb Z[\underline{G}_{(\kappa)}]^{L\solid}\cong \oplus_{i\in I^p}\mathbb Z[\prod_{k=1}^p\underline{X({i_k})}_{(\kappa)}]^{L\solid}\in\mathcal D(\Sol{(\kappa)}).\] 
\cref{solidificationonderivedmappingspacesforacyclicrings} now implies that for $M\in\Sol{(\kappa)}$, 
    \begin{align*}\ckH(\underline{G^p}_{(\kappa)},M)&\cong \map_{\mathcal D(\Cond{(\kappa)}(\Ab))}(\mathbb Z[\underline{G^p}_{(\kappa)}],M)\\ &\cong \map_{\mathcal D(\Cond{(\kappa)}(\Ab))}(\mathbb Z[\underline{G^p}_{(\kappa)}]^{L\solid},M)\\&\cong \map_{\mathcal D(\Cond{(\kappa)}(\Ab))}(\oplus_{i\in I^p}\mathbb Z[\prod_{k=1}^n\underline{X({i_k})}_{(\kappa)}]^{L\solid},M)\\ 
    &\cong \prod_{i\in I^n}\map_{\mathcal D(\Cond{(\kappa)}(\Ab))}(\mathbb Z[\prod_{k=1}^p\underline{{X({i_k})}}_{(\kappa)}]^{L\solid},M)\\
&\cong \prod_{i\in I^n}\map_{\mathcal D(\Cond{(\kappa)}(\Ab))}(\mathbb Z[\prod_{k=1}^p\underline{{X({i_k})}}_{(\kappa)}],M), \end{align*} where we used that coproducts in $\Cond{(\kappa)}(\Ab)$ are exact (\cref{condensedabeliangroupsgrothendieckaxioms}). 
The right-hand side is concentrated in degree $0$ by \cref{solidderived}. 
\end{proof}
\cref{condensedandsheafcohomology1,sheafcohomologywvspacecoeffvanonparcomplc} imply:
\begin{cor}\label{continuousequalscondensedonvectorspacecoefficients}
    Suppose that $G$ is a locally compact Hausdorff topological group and $V$ is a finite-dimensional normed $\mathbb R$-vector space with a linear, continuous $G$-action.  
    For $\kappa>\wt(G)$, the spectral sequence from \cref{spectralsequenceidentifycontinuousgroupcohomology} collapses at the $E_2$-page and yields 
    \[ \contgrpcoh{}^*(G,V)\cong \contgrpcoh{\kappa\text{-}}^*(G,V)\cong \condgrpcoh{\kappa\text{-}}^*(\underline{G}_{\kappa}, \underline{V}_{\kappa}).\] 

    Analogously, 
    \[ \contgrpcoh{}^*(G,V)\cong \contgrpcoh{k\text{-}}^*(G,V)\cong \condgrpcoh{}^*(\underline{G}, \underline{V}).\]
\end{cor}
\begin{proof}
This follows from 
\cref{condensedandsheafcohomology1,sheafcohomologywvspacecoeffvanonparcomplc}, since every locally compact Hausdorff topological group is paracompact, see e.g.\ \cite[page 1]{topologicalgroupslocalversusglobal}. 
\end{proof}

\begin{rem}
Continuous group cohomology is insensitive with respect to the higher terms \[(E_1^{p,q}=\cckH{(\kappa\text{-})}^q(\underline{G}^p, \res{G}{*, \Sp}-))_{p,q\geq 1}\] of the spectral sequence from \cref{chechtocohomologyspectralsequencegroupcohomology}, which are nontrivial in many cases. 
For example, \cite[Corollary II.16.21]{Bredon} and \cref{condensedandsheafcohomology1} imply that if $G$ is locally compact Hausdorff (and $\kappa$-light), then $\cckH{(\kappa\text{-})}^*(\underline{G}_{\kappa},A)$ is concentrated in degree $0$ for all discrete abelian groups $A$ if and only if $G$ is locally profinite. 
This suggests that continuous group cohomology is a coarser invariant than condensed group cohomology. 

Their difference can be made more explicit: 
Continuous group cohomology with coefficients in a discrete continuous $G$-module only depends on the group of connected components $\pi_0G$ of $G$: For a discrete continuous $G$-module $M$, \[\contgrpcoh{}^*(G,M)\cong \contgrpcoh{}^*(\pi_0G,M),\] where $\pi_0G$ is endowed with the quotient topology, and analogously for $k$-continuous/$\kappa$-continuous group cohomology provided that $G^i$ is compactly generated/$\kappa$-compactly generated for all $i\in\mathbb N_0$ or alternatively, that $G$ is locally path-connected. 

We will see in \cref{section:condensedcohomologyandcohomologyofclassifyingspaces}, specifically \cref{condensedequalscontinuousimpliescontractible} that $\condgrpcoh{}^*(\underline{G}, \underline{M})$ also encodes information on the homotopy type of the connected component of the identity $G^0\subseteq G$.   
\end{rem}
Condensed group cohomology has better exactness properties than continuous group cohomology: 
\begin{ex}
    The exponential sequence $0\to \mathbb Z \to \mathbb R \to \sphere{1}\to 0$ induces a short exact sequence of ($\kappa$-)condensed abelian groups \[0\to \underline{\mathbb Z}_{(\kappa)}\to \underline{\mathbb R}_{(\kappa)}\to \underline{\sphere{1}}_{(\kappa)}\to 0.\] 
    Suppose $G$ is a ($\To$) topological group and consider the above as condensed $\mathbb Z[\underline{G}_{(\kappa)}]$-modules with trivial $G$-action. This induces a long exact sequence of abelian groups
    \begin{align*}0\to\mathbb Z\to \mathbb R\to \sphere{1}\to \condgrpcoh{(\kappa\text{-})}^1(\underline{G}_{(\kappa)}, \underline{\mathbb Z}_{(\kappa)})& \to \condgrpcoh{(\kappa\text{-})}^1(\underline{G}_{(\kappa)}, \underline{\mathbb R}_{(\kappa)})\to \condgrpcoh{(\kappa\text{-})}^1(\underline{G}_{(\kappa)}, \underline{\sphere{1}}_{(\kappa)})\\& \to \condgrpcoh{(\kappa\text{-})}^2(\underline{G}_{(\kappa)}, \underline{\mathbb Z}_{(\kappa)})\to \ldots.\end{align*} 
    
    Continuous group cohomology does not fit into such a long exact sequence: 
    For a continuous $G$-module with trivial $G$-action $M$, $H^1_{\operatorname{cont}}(G,M)$ is the group of continuous group homomorphisms $G\to M$. In particular, if $G$ is compact, then $H^1_{\operatorname{cont}}(G, \mathbb R)=0$, whereas $H^1_{\operatorname{cont}}(G, \sphere{1})=\widehat{G}$ is the Pontrjagin-dual of $G$. 
    Since $\contgrpcoh{}^*(G, \mathbb Z)\cong \contgrpcoh{}^*(\pi_0G, \mathbb Z)$, this shows that if $0\neq G$ is a compact and connected topological group, there is no exact sequence 
    \[ \ldots \to \contgrpcoh{}^1(G, \mathbb Z)\to \contgrpcoh{}^1(G, \mathbb R)\to \contgrpcoh{}^1(G, \sphere{1})\to\contgrpcoh{}^2(G, \mathbb Z)\to \ldots .\] 
\end{ex}
 
\paragraph{Comparison with group cohomology in the gros topos}\label{comparisongrostoposgroupcohomology}
\cite{flach} studied group cohomology in the gros topos of topological spaces (\cref{grostoposdefinition}). Ibid.\ used universe enlargement to deal with size issues, which corresponds to taking $\lambda$ a regular strong limit cardinal.
The spectral sequence provided by \cref{chechtocohomologyspectralsequencegroupcohomology} in the gros topos, and the equivalence of sheaf and gros topos cohomology (\cref{grostoposequalssheaf}) imply that if  $G$ is a $\lambda$-small topological group and $M$ is a continuous $G$-module such that $H^q_{\sheaf}(G^p, \mathcal C(-,M))=0$ for all $p,q\in \mathbb N_1$, then 
\[\contgrpcoh{}^*(G,M)\cong \grpcoh{\grostop}^*(h_G,h_M)\] where $h_G,h_M$ denote the sheaves represented by $G$ and $M$ respectively. This was shown in \cite{flach}. 
It applies in the following two cases (both shown by \cite{flach}): 
\begin{romanenum}
\item If $G$ is locally profinite, in which case $\mathbb Z[h_{G^p}]$ is projective in $\Shv_{LS}(\Top^{\kappa}, \Ab)=\stab{\grostop}^{\heart}$: By van Dantzig's theorem, every locally profinite group is a coproduct of profinite sets and for $X\in\Pro(\Fin)_{\lambda}$,  
$\mathbb Z[h_X]$ is projective in $\Shv_{LS}(\Top^{\lambda}, \Ab)$ since for $i>0$, \begin{align*}\Ext^i_{\Shv_{LS}(\Top^{\lambda}, \Ab)}(\mathbb Z[h_X],-)|_{\Shv_{LS}(\Top^{\lambda},\Ab)}&\cong \cH{\widehat{\grostop}}^i(h_X,-)|_{\Shv_{LS}(\Top^{\lambda}, \Ab)}\\ &\cong \cH{\grostop}^i(h_X,-)|_{\Shv_{LS}(\Top^{\lambda}, \Ab)}\\ &\cong \cH{\sheaf}^i(X, \pi_0\stradg{t}-)=0.\end{align*} 
Here $\widehat{\grostop}\subseteq \grostop$ denotes the subcategory on hypercomplete objects. 
The vanishing of sheaf cohomology is \cite[\href{https://stacks.math.columbia.edu/tag/0A3F}{Tag 0A3F}]{stacks-project}, the identifications hold by \ref{cohomologyisext5}, \cref{geometricmorphismcohomology,grostoposequalssheaf}. 
\item If $G$ is locally compact Hausdorff and $M$ is a finite-dimensional, continuous real $G$-representation. This follows from \cref{sheafcohomologywvspacecoeffvanonparcomplc} since locally compact Hausdorff groups are paracompact (\cite[page 1]{topologicalgroupslocalversusglobal}). 
\end{romanenum}

Suppose $\kappa$ is an uncountable cardinal and $\lambda\geq \kappa$ is a strong limit cardinal. 
In \cref{grostoposgeometricmorphismcondensed} (page \pageref{grostoposgeometricmorphismcondensed}), we described a geometric morphism \begin{equation}\label{grostopcondensedgeometricmorphism}\radg{j}\colon \Cond{\kappa}(\an)\leftrightarrows \grostop\colon \ladg{j}, \ladg{j}\dashv \radg{j}.\end{equation}
For $G\in \Grp(\grostop)$, this yields a geometric morphism 
\begin{equation}\label{grostoposcondensedgobjects} \Gradg{j}\colon \Act{\ladg{j}G}{\Cond{\kappa}(\an)}\leftrightarrows \Act{G}{(\grostop)}\colon \Gladg{j}, \end{equation} cf.\ \cref{adjunctioninducesadjunctionongobjects}. In particular, 
\[\condgrpcoh{\kappa\text{-}}(\ladg{j}G,-)\cong \grpcoh{\grostop}(G, \Gstradg{j}-)\] by \cref{geometricmorphismcohomology}.  
If $G=h_{G}$ is represented by a $\lambda$-small topological group $G$, then $\ladg{j}h_G=\underline{G}_{\kappa}$ by \cref{grostoposcondensedzariskirepresentables}, and if $A$ is a continuous $G$-module $A$, then $\pi_0\Gstradg{j}\underline{A}_{\kappa}\cong \mathcal C_{\kappa}(-,A)$ with the obvious $h_G=\mathcal C(-,G)$-action. 

\begin{lemma}\label{comparisongrostoposgroupcohomologyl}
    Suppose $\kappa\leq \lambda$ are uncountable cardinals and $G$ is a $\lambda$-small topological abelian group.\begin{romanenum}
    \item If $A\in\Cond{\kappa}(\mathbb Z[\underline{G}_{\kappa}])$ is such that the comparison map \[ \cH{\sheaf}(G^p,\pi_0\stradg{i}A)\to \ckH(\underline{G}^p_{\kappa},A) \,\, (\ref{comparisonmapkappacondensedsheaf}) \]  is an isomorphism for all $p\in\mathbb N_0$, then \[ \grpcoh{\grostop}(h_G, \pi_0\Gstradg{j}A)\cong\condgrpcoh{\kappa\text{-}}(\underline{G}_{\kappa},A).\]
    \item If $A$ is a topological group such that the comparison map \[ \cH{\sheaf}(G^p,A)\to \ckH(\underline{G}^p_{\kappa},\underline{A}_{\kappa}) \, \, (\ref{comparisonmapkappacondensedsheaftopologicalwithoutkappa})\] is an isomorphism for all $p\in\mathbb N_0$, then \[ \grpcoh{\grostop}(h_G, h_A)\cong\condgrpcoh{\kappa\text{-}}(\underline{G}_{\kappa},\underline{A}_{\kappa}).\]
   \end{romanenum}
    \end{lemma}
This applies for example in the following two cases: 
\begin{romanenum}
\item $G$ is $\kappa$-light, locally compact Hausdorff and $A$ is a product of a discrete group and a finite-dimensional normed $\mathbb R$-vector space (\cref{condensedandsheafcohomology1}), 
\item $G$ is homotopy equivalent to a locally contractible or $\kappa$-light locally compact Haudorff space and $A$ is discrete (\cref{condensedandsheafcohomologylocallycontractible,condensedandsheafcohomologylocallycompacthausdorffhomotopy}).
\end{romanenum}
    \begin{proof}We only prove the first statement, the second statement can be shown completely analogously. 
By \cref{spectralenrichmentcocontinuous}, \[\condgrpcoh{\kappa\text{-}}(\underline{G}_{\kappa},-)\cong \clim{\Delta}\cH{\Act{\underline{G}_{\kappa}}{\Cond{\kappa}(\an)}}(\check{C}(\underline{G}_{\kappa}\to *),-)\] 
and \[\grpcoh{\grostop}(h_G,-)\cong \clim{\Delta}\cH{\Act{h_G}{(\grostop)}}(\check{C}(h_G\to *),-).\]  
By \cref{geometricmorphismcohomology}, the geometric morphism \ref{grostoposcondensedgobjects} induces a map of cosimplicial spectra 
\[\cH{\Act{h_G}{(\grostop)}}(\check{C}(h_G\to *), \pi_0\Gstradg{j}A)\to \cH{\Act{h_G}{(\grostop)}}(\check{C}(h_G\to *), \Gstradg{j}A)\cong \cH{\Act{\underline{G}_{\kappa}}{\Cond{\kappa}(\an)}}(\check{C}(\underline{G}_{\kappa}\to *),A).\] 
For $n\in\mathbb N_0$,
\[ \cH{\Act{\underline{G}_{\kappa}}{\Cond{\kappa}(\an)}}(\check{C}(\underline{G}_{\kappa}\to *),-)[n]\cong \ckH(\underline{G}_{\kappa}^n, \res{G}{*, \Sp}-),\] and 
\[ \cH{\Act{h_G}{\Cond{\grostop}}}(\check{C}(h_G\to *),-)[n]\cong \cH{\grostop}(h_{G^n}, \res{G}{*, \Sp}-)\]
by \cref{adjunctionsspectrallyenriched,indgspisstabilization}. 
As $\res{G}{*, \Sp}\circ \Gstradg{j}=\stradg{j}\circ \res{G}{*, \Sp}$ and $\res{G}{*, \Sp}$ is $t$-exact, it follows that 
\[\grpcoh{\grostop}(\check{C}(h_G\to *), \pi_0\Gstradg{j}-)([n])\cong  \cH{{\grostop}}(h_{G^n}, \res{G}{*, \Sp}\pi_0\Gstradg{j}-)=\cH{{\grostop}}(h_{G^n}, \pi_0\stradg{j}\res{G}{*, \Sp}-).\] 
Under these identifications, the above morphism of cosimplicial spectra agrees with the map 
\[ \cH{\grostop}(-,\pi_0\stradg{j}-)\to \cH{\grostop}(-,\stradg{j}-)\cong \ckH(\ladg{j}-,-)\] induced by $\pi_0\stradg{j}\to \stradg{j}$ and the equivalence from \cref{geometricmorphismcohomology}. By \cref{characterizationexactnessgrostopos}, this is an equivalence if and only if the comparison map \[ \cH{\sheaf}(G^p,\pi_0\stradg{i}A)\to \ckH(\underline{G^p}_{\kappa},A) \,\, (\ref{comparisonmapkappacondensedsheaf})\] is an isomorphism for all $p\in\mathbb N_0$. 
\end{proof}
  
\subsection{Condensed group cohomology and cohomology of classifying spaces}\label{section:condensedcohomologyandcohomologyofclassifyingspaces}
In this section, we compare condensed group cohomology of a topological group $G$ with the condensed cohomology of a classifying space $BG$ for numerable principal $G$-bundles. 
We refer to \cite[section 14.4]{tomdiecktopology} for basic material on classifying spaces and principal bundles. 
The main result is that condensed group cohomology with coefficients in a solid abelian group $M$ with trivial $G$-action is isomorphic to the condensed cohomology of $BG$ with coefficients in $M$. This will follow from our general discussion around \cref{groupcohomologyofbundles} and homotopy invariance of condensed cohomology with solid coefficients (\cref{solidcohomologyhomotopyinvariant}). 

We start with some preliminaries on principal $G$-bundles. Recall that for a topological group $G$, the Milnor join construction (\cite{de56d803-a9ff-36ca-abf7-383d35d44527}) yields a universal numerable principal $G$-bundle $EG\to BG$. 
\begin{lemma}\label{classifyingspacet1}If $G$ is a Hausdorff topological group, the classifying space $BG$ of numerable principal $G$-bundles obtained through the Milnor join construction is Hausdorff.
\end{lemma}
\begin{proof}We first show that $EG$ is Hausdorff.  
Recall that elements of $EG$ are represented by tuples 
$(g_i,t_i)_{i\in\mathbb N_0}\in (G\times [0,1])^{\mathbb N}$ such that $t_i=0$ for all but finitely many $i\in\mathbb N_0$, $\sum_{i\in \mathbb N_0}t_i=1$, and two such tuples $(g_i,t_i), (h_i,s_i)$ represent the same element if and only if $t_i=s_i$ for all $i\in\mathbb N_0$ and for all $i\in\mathbb N_0$ with $t_i\neq 0$, $g_i=h_i$.
The topology on $EG$ is the coarsest topology such that the projections 
\begin{align*}p_i\colon EG\to [0,1], (g_i,t_i)_i\mapsto t_i\end{align*} and \begin{align*}q_i\colon p_i^{-1}(\rbrack0,1\rbrack)\to G, (g_i,t_i)\mapsto g_i\end{align*} are continuous for all $i\in\mathbb N_0$. 

Suppose now that $(g_i,t_i), (h_i,s_i)$ represent two distinct elements in $EG$. 
If there exists $i\in \mathbb N_0$ with $t_i\neq s_i$, choose disjoint open neighborhoods $t_i\in U\subseteq [0,1], s_i\in V\subseteq [0,1]$. 
Then \[(g_i,t_i)_i\in p_i^{-1}(U)\] and \[(h_i,s_i)\in p_i^{-1}(V)\] are disjoint open neighborhoods in $EG$. 
Suppose now that $t_i=s_i$ for all $i\in \mathbb N_0$. As $(g_i,t_i), (h_i,s_i)$ represent distinct elements, there has to exist $n\in \mathbb N_0$ such that $t_n\neq 0$ and $g_n\neq h_n$. 
Since $G$ is Hausdorff, there exist disjoint open neighborhoods $g_n\in U\subseteq G, h_n\in V\subseteq G$. Now \[(g_i,t_i)_i\in q_n^{-1}(U), (h_i,s_i)_i\in q_n^{-1}(V)\] are disjoint open neighborhoods in $EG$. This shows that $EG$ is Hausdorff. 

Suppose now that $[(g_i,t_i)_i],[(h_i,s_i)_i]\in BG=EG/G$ represent two distinct elements. 
The projections $(p_i)_{i\in\mathbb N_0}$ descend to continuous maps $\overline{p}_i\colon BG\to [0,1]$. 
If there exists $i\in\mathbb N_0$ such that $t_i\neq s_i$, then choose disjoint open neighborhoods $t_i\in U\subseteq [0,1], s_i\in V\subseteq [0,1]$. 
Then $[(g_i,t_i)_i]\in \overline{p}_i^{-1}(U)$ and $[(h_i,s_i)]\in \overline{p}_i^{-1}(V)$ are disjoint open neighborhoods in $BG$.
Suppose now that $t_i=s_i$ for all $i\in \mathbb N_0$ and choose $i\in \mathbb N_0$ with $t_i>0$. 
The proof of \cite[Theorem 3.1]{de56d803-a9ff-36ca-abf7-383d35d44527} shows that \[V_i\coloneqq \overline{p}_i^{-1}(\rbrack 0,1\lbrack)\subseteq BG\] is an open subset of $BG$ over which the bundle $\pi\colon EG\to BG$ trivializes. 
Since \[\pi^{-1}(V_i)\cong G\times V_i\subseteq EG\] is a subspace, it is Hausdorff. This implies that $V_i$ is Hausdorff. As $[(t_i,g_i)], [(s_i,h_i)]\in V_i$, they can be separated by disjoint open neighborhoods in $V_i$ and hence in $BG$. 
\end{proof}

\begin{lemma}\label{bundlesarecondensedbundles}
    Suppose $G$ is a topological group and $E\to B$ is a principal $G$-bundle. 
    \begin{romanenum}
    \item For all uncountable cardinals $\kappa$, $\underline{E}_{\kappa}\to\underline{B}_{\kappa}$ is a principal $\infty$-$\underline{G}_{\kappa}$-bundle in $\Cond{\kappa}(\an)$ (\cref{definitionprincipalinfinitybundle}). 

    \item If $B$ and $G$ are $\To$, then so is $E$ and $\underline{E}\to\underline{B}$ is a principal $\infty$-$\underline{G}$-bundle in $\Cond{}(\an)$. 
    \end{romanenum} 
\end{lemma}
\begin{proof}
    Choose an open cover $B=\cup_{i\in I}U_i$ such that for all $i\in I$, $E\times_{B}U_i\cong G\times U_i$ as topological spaces with $G$-action over $U_i$.
    Then $E=\bigcup_{i\in I}E\times_{B} U_i$ is an open cover. In particular, if $B$ and $G$ are $\To$, then $E\times_B U_i\cong G\times U_i$ is $\To$ for all $i\in I$, which implies that $E$ is $\To$. 
    As the local section topology is coarser than the $\kappa$-condensed topology on $\Top^{\lambda}$ (\cref{grostopos}), $p\colon \sqcup_{i\in I}\underline{U_i}_{(\kappa)}\to \underline{B}_{(\kappa)}$ is an effective epimorphism in $\Cond{(\kappa)}(\an)$. 
    Since $\underline{(-)}_{(\kappa)}$ commutes with pullbacks and small colimits are universal in $\Cond{(\kappa)}(\an)$ (\cref{Girauds axioms}),
    \begin{align*}\underline{E}_{(\kappa)}\times_{\underline{B}_{(\kappa)}}\sqcup_{i\in I}\underline{U_i}_{(\kappa)}&\cong \sqcup_{i\in I}\underline{E}_{(\kappa)}\times_{\underline{B}_{(\kappa)}}\underline{U_i}_{(\kappa)}\\&\cong \sqcup_{i\in I}\underline{E\times_B U_i}_{(\kappa)}\\ &\cong \sqcup_{i\in I}\underline{G\times_B U_i}_{(\kappa)}\\&\cong \underline{G}_{(\kappa)}\times_{\underline{B}_{(\kappa)}}\sqcup_{i\in I}\underline{U_i}_{(\kappa)}.\qedhere\end{align*}
\end{proof}
By \cref{bundlesarecondensedbundles} and \cref{characterizationbundles}, every numerable principal $G$-bundle $E\to B$ determines a map $\underline{B}_{(\kappa)}\cong \underline{E}_{(\kappa)}//\underline{G}_{(\kappa)}\to \mathbb B\underline{G}_{(\kappa)}$ to the classifying space $\mathbb B\underline{G}_{(\kappa)}$ of principal $\infty$-$\underline{G}_{(\kappa)}$-bundles in $\Cond{(\kappa)}(\an)$. 
By \cref{grpcohomologyiscohomologyonclassifyingspaces}, this induces a map \begin{align}\label{condensedgroupcohomologycohomologybgcomparisonmap}\condgrpcoh{(\kappa\text{-})}(\underline{G}_{(\kappa)}, \triv_{\Sp}-)\cong \cckH{\kappa\text{-}}(\mathbb B\underline{G}_{(\kappa)},-)\to \cckH{\kappa\text{-}}(\underline{BG}_{\kappa},-).\end{align} (Without $\kappa$, both sides are only defined if $G$ and $BG$ are $\To$). 

The condensed anima represented by a classifying space of numerable principal $G$-bundles is very different from the classifying space of principal $\infty$-$\underline{G}_{(\kappa)}$-bundles in $\Cond{(\kappa)}(\an)$: 
If $EG\to BG$ is a universal numerable principal $G$-bundle, then $\underline{BG}_{(\kappa)}\in\tau_{\leq 0}\Cond{\kappa}(\an)$ is $0$-truncated. On the other hand, $\pi_1(\mathbb B\underline{G}_{\kappa})=\underline G_{(\kappa)}$ (\cref{deloopingclassifyingspace}), and in particular $\underline{BG}_{\kappa}\not\cong \mathbb B\underline{G}_{\kappa}$ if $G\neq 0$. 
However, since the total space $EG$ of a universal numerable principal $G$-bundle is contractible as space without group action (see e.g.\ \cite[Theorem 14.4.12]{tomdiecktopology}), homotopy invariance of cohomology with solid coefficients (\cref{groupcohomologyofbundles}) and \cref{solidcohomologyhomotopyinvariant} imply the following:
\begin{cor}\label{condensedcohomologyiscohomologyofclassifyingspaces}
Suppose $G$ is a topological group and $BG$ is a classifying space for numerable principal $G$-bundles. 
Denote by $f$ the forget functor \[\mathcal D(\Sol{(\kappa)})\to\mathcal D(\Cond{(\kappa)}(\Ab))\cong\LMod{H\mathbb Z}{\Cond{(\kappa)}(\Ab)}\to \stab{\Cond{(\kappa)}(\Ab)}.\] 
\begin{romanenum}
\item For $M\in \mathcal D(\Sol{\kappa})$, \[ \ckH(\underline{BG}_{\kappa},fM)\cong \condgrpcoh{\kappa\text{-}}(\underline{G}_{\kappa}, \triv_{\Sp} fM)\] via the map \ref{condensedgroupcohomologycohomologybgcomparisonmap}. 

\item If $G$ is Hausdorff and $BG$ is a classifying space for numerable principal $G$-bundles $BG$ which is $\To$ (e.g.\ the Milnor join construction), for all $M\in\mathcal D(\Sol{})$, \[ \ccH(\underline{BG}_{},fM)\cong \condgrpcoh{}(\underline{G}, \triv_{\Sp} fM)\] via the map \ref{condensedgroupcohomologycohomologybgcomparisonmap}. 
\end{romanenum}
\end{cor}
\begin{proof}Choose a universal numerable principal $G$-bundle $EG\to BG$. Since $EG$ is contractible (\cite[Theorem 14.4.12]{tomdiecktopology}), by \cref{solidcohomologyhomotopyinvariant}, for all $M\in\Sol{(\kappa)}$ and $p\in\mathbb N_0$, the projection $G^p\times EG\to G^p$ induces an equivalence \[\cckH{(\kappa\text{-})}(\underline{G}_{(\kappa)}^p,M)\cong \cckH{(\kappa\text{-})}(\underline{G}_{(\kappa)}^p\times\underline{EG}_{(\kappa)},M)\] (without $\kappa$, both sides are only defined if $G$ and $EG$ are $\To$).  
The statement now follows from \cref{groupcohomologyofbundles}. 
\end{proof}
\subsubsection{Cohomology of classifying spaces}\label{section:Cohomologyclassifyingspaces}
We now show that for many topological groups, the condensed cohomology of their classifying space agrees with its sheaf and singular cohomology. 
Together with \cref{condensedcohomologyiscohomologyofclassifyingspaces}, this yields many examples where condensed group cohomology differs from continuous group cohomology, see \cref{condensedequalscontinuousimpliescontractible}. 

Homotopy invariance of gros topos/sheaf cohomology with constant coefficients (\cref{grostoposcohomologyomotopyinvariant}) and homotopy invariance of condensed cohomology with solid coefficients (\cref{solidcohomologyhomotopyinvariant}) imply the following: 
\begin{lemma}\label{covercohomologybg}
    Suppose that $M$ is a discrete abelian group, $G$ is a topological group and $BG$ is a classifying space for numerable principal $G$-bundles. 
\begin{romanenum}
\item If for all $p\in\mathbb N_0$, the
map \[\cH{\sheaf}(G^p,M)\to \ckH(\underline{G^p}_{\kappa},\underline{M}_{\kappa})\,\,  (\ref{comparisonmapkappacondensedsheaftopological}) \] is an equivalence, then \[ \cH{\sheaf}(BG,M)\cong \ckH(\underline{BG}_{\kappa},\underline{M}_{\kappa})\] via the map (\ref{comparisonmapkappacondensedsheaftopological}).  
\item If $G$ is $\To$ and $BG$ is $\To$ (e.g.\ the classifying space obtained through the Milnor-join construction), and for all $p\in\mathbb N_0$, the map \[\cH{\sheaf}(G^p,M)\to\cH{\sheaf}(G^p,M_{\kappa})\to \ccH(\underline{G^p},\underline{M}) \, \,(\ref{comparisonmapkappacondensedsheaftopologicalwithoutkappa}) \] is an equivalence, then \[ \cH{\sheaf}(BG,M)\cong \ccH(\underline{BG},\underline{M})\] via the map (\ref{comparisonmapkappacondensedsheaftopologicalwithoutkappa}).
\end{romanenum}
\end{lemma}
\begin{proof}Denote by $EG\to BG$ a universal principal $G$-bundle. Since $EG$ is contractible (\cite[Proposition 14.4.6]{tomdiecktopology}), our assumption and \cref{deltakappaexactnesshomotopyinvariant} imply that $G^p\times EG$ is $\epsilon$-($\kappa$-)$M$-exact (\cref{deltakappaexactgoodcover}) for all $p\in\mathbb N_0$. As $EG\to BG$ is a local section cover, and for $p\in\mathbb N_0$, \[\check{C}(EG\to BG)([p])\cong G^p\times EG,\] the statement now follows from \cref{deltakappaexactgoodcover}. 
\end{proof}

\begin{cor}\label{condensedgroupcohomologyissheafcohomologyofbg}
    Suppose $G$ is a topological group with one of the following properties: 
    \begin{romanenum}
        \item $G$ is homotopy equivalent to a locally $\kappa$-light compact Hausdorff topological space. 
        \item $G$ is homotopy equivalent to a space which is locally contractible. 
    \end{romanenum}
If $BG$ is a classifying space for numerable principal $G$-bundles, for all discrete abelian groups $M$, 
\[\cH{\sheaf}(BG,M)\cong \ckH(\underline{BG}_{\kappa},\underline{M}_{\kappa})\] via the map (\ref{comparisonmapkappacondensedsheaftopological}).  
If in addition to the above, $G$ is Hausdorff and $BG$ is a $\To$ classifying space for numerable principal $G$-bundles (e.g.\ the Milnor join construction), then for all discrete abelian groups $M$, 
\[\cH{\sheaf}(BG,M)\cong \ccH(\underline{BG},\underline{M})\] via the map (\ref{comparisonmapkappacondensedsheaftopologicalwithoutkappa}).  
\end{cor}
\begin{proof}It suffices to prove the $\kappa$-condensed statement, then the condensed statement follows from \cref{condensedcohomologycanbecomputedonfinitestage}. 
By \cref{deltakappaexactnesshomotopyinvariant,condensedandsheafcohomology1,condensedandsheafcohomologylocallycontractible}, our conditions on $G$ imply that for all $p\in\mathbb N_0$ and all discrete abelian groups $M$, the map (\ref{comparisonmapkappacondensedsheaftopological})
\begin{align}\label{obstruction} \cH{\sheaf}(G^p,M)\to \ckH(\underline{G^p}_{\kappa},\underline{M}_{\kappa})\end{align} is an equivalence. 
The statement now follows from \cref{covercohomologybg}. 
\end{proof}
\begin{lemma}
Suppose $G$ is a topological group which is homotopy equivalent to a locally contractible topological space and $BG$ is a classifying space for numerable principal $G$-bundles.
\begin{romanenum} 
\item For all solid abelian groups $M\in\Sol{\kappa}$, 
\[ \ckH(\underline{BG}_{\kappa},\underline{M(*)^{\delta}}_{\kappa})\cong \ckH(\underline{BG}_{\kappa},M),\] where $M(*)^{\delta}$ denotes the discrete abelian group $M(*)$.  
\item If $G$ is Hausdorff and $BG$ is $\To$  (e.g.\ the Milnor join construction), then for all solid abelian groups $M\in\Sol{}$, 
\[\ccH(\underline{BG},\underline{M(*)^{\delta}})\cong \ccH(\underline{BG},M).\] 
\end{romanenum}
\end{lemma}
\begin{proof}Fix a universal numerable principal $G$-bundle $EG\to BG$. 
As $\underline{EG}_{(\kappa)}\to \underline{BG}_{(\kappa)}$ is a principal $\infty$-$\underline{G}_{\kappa}$-bundle (\cref{bundlesarecondensedbundles}), 
\[  \cckH{(\kappa-)}(\underline{BG}_{(\kappa)},-)=\clim{\Delta}\cckH{(\kappa-)}(\check{C}(\underline{EG}_{(\kappa)}\to \underline{BG}_{(\kappa)}),-)\] by \cref{homotopyquoptientbundle,characterizationbundles,spectralenrichmentcocontinuous}.
As $EG$ is contractible, our assumption on $G$ implies that for all $p\in\mathbb N_0$, \[\check{C}(EG\to BG)([p])=G^p\times EG\] is homotopy equivalent to a locally contractible space, whence the counit $\underline{M(*)^{\delta}}_{(\kappa)}\to M$ of the constant sheaf/global sections geometric morphism induces an equivalence of cosimplicial spectra 
\[ \cckH{(\kappa-)}(\check{C}(\underline{EG}_{(\kappa)}\to \underline{BG}_{(\kappa)}),\underline{M(*)^{\delta}}_{(\kappa)})\cong \cckH{(\kappa-)}(\check{C}(\underline{EG}_{(\kappa)}\to \underline{BG}_{(\kappa)}),M)\] by \cref{condensedandsheafcohomologylocallycontractible}. (Without $\kappa$, the terms are only defined if $G$ and $BG$ are $\To$, in which case $EG$ is $\To$ as well). 
\end{proof}
\begin{lemma}\label{singularissheafcohomologybg}
    Suppose $G$ is a topological group and $BG$ is a classifying space for numerable principal $G$-bundles.
    For all discrete abelian groups $M$, \[ \cH{\sheaf}^*(BG,M)\cong H^*_{\operatorname{sing}}(BG,M)\] and naturally in $M$.  
\end{lemma}
\begin{proof}
    We claim that $BG$ is semi-locally contractible (in the sense of \cite{sella2016comparison}), then the statement follows from \cite{sella2016comparison}. Fix an open subset $U\subseteq BG$. As $U\times_{BG}EG\to U$ is a principal $G$-bundle, there exists an open cover $U=\bigcup_{i\in I}U_i$ such that $EG\times_{BG}U_i\cong G\times U_i$ as spaces with $G$-action over $U_i$. 
    As $BG$ is a classifying space for principal $G$-bundles, this implies that the inclusion $U_i\to BG$ is homotopic to a constant map for all $i\in I$.
\end{proof}

The identification of condensed cohomology with cohomology of classifying spaces yields many examples where condensed group cohomology differs from continuous group cohomology: 
\begin{cor}\label{condensedequalscontinuousimpliescontractible}
Suppose $G$ is a path-connected topological group which satisfies one of the following conditions:
\begin{romanenum}
    \item $G$ is homotopy equivalent to a ($\kappa$-light) locally compact Hausdorff space. 
    \item $G$ is homotopy equivalent to a locally contractible topological space.
    \item There exists a classifying space for numerable principal $G$-bundles which is homotopy equivalent to a locally contractible space or to a locally ($\kappa$-light) compact Hausdorff space. 
\end{romanenum}
\begin{itemize}
\item If for all discrete abelian groups $M$, 
\[ \contgrpcoh{\kappa\text{-}}^*(G,M)\cong \condgrpcoh{\kappa\text{-}}(\underline{G}_{\kappa}, \underline{M}_{\kappa}), \] then $G$ is weakly contractible. 

\item If $G$ is Hausdorff and for all discrete abelian groups $M$, 
\[ \contgrpcoh{k\text{-}}^*(G,M)\cong \condgrpcoh{}(\underline{G}, \underline{M}), \] then $G$ is weakly contractible. 
\end{itemize}
\end{cor}
\begin{proof}
    Since $G$ is path-connected, so is $G^i$ for all $i\in\mathbb N_0$ and in particular every $\kappa$/$k$-continuous map $G^i\to M$ to a discrete abelian group is constant. 
    This implies that for all discrete continuous $G$-modules $M$, \[M\cong \contgrpcoh{}^*(G,M)\cong \contgrpcoh{k\text{-}}^*(G,M)\cong \contgrpcoh{\kappa\text{-}}^*(G,M)\] is concentrated in degree $0$. 

    We now show that under the assumptions on $G$, for all discrete abelian groups $M$, \[\condgrpcoh{(\kappa-)}^*(\underline{G}_{(\kappa)},\underline{M}_{(\kappa)})\cong H^*_{\operatorname{sing}}(BG,M).\footnote{ 
    If $G$ is Hausdorff, we choose $BG$ to be $\To$, so that condensed cohomology is well-defined. This is possible by \cref{classifyingspacet1}.}\] 
    By \cref{condensedcohomologyiscohomologyofclassifyingspaces,profinitesolid}, for all discrete abelian groups $M$ (with trivial $G$-action), \[\condgrpcoh{(\kappa\text{-})}(\underline{G}_{(\kappa)}, \underline{M}_{(\kappa)})\cong \cckH{(\kappa-)}(\underline{BG}_{(\kappa)}, \underline{M}_{(\kappa)}).\]  
    If $G$ is homotopy equivalent to a ($\kappa$-light) locally compact Hausdorff space or to a locally contractible space, for all discrete abelian groups $M$, \[\cH{\sheaf}(BG,M)\cong \cckH{(\kappa-)}(\underline{BG}_{(\kappa)},\underline{M}_{(\kappa)})\] by \cref{condensedgroupcohomologyissheafcohomologyofbg}, and if there exists a classifying space $BG$ for numerable principal $G$-bundles which is homotopy equivalent to a locally contractible or locally ($\kappa$-)light compact Hausdorff space, then \[\cH{\sheaf}(BG,M)\cong \cckH{(\kappa-)}(\underline{BG}_{(\kappa)},\underline{M}_{(\kappa)})\] by \cref{condensedandsheafcohomologylocallycontractible,condensedandsheafcohomologylocallycompacthausdorffhomotopy}, respectively. 
    By \cref{singularissheafcohomologybg},
     \[H^*_{\operatorname{sing}}(BG,M)\cong \cH{\sheaf}^*(BG,M).\]

    In particular, if \[ \contgrpcoh{\kappa/k\text{-}}^*(G,M)\cong \condgrpcoh{(\kappa\text{-})}^*(\underline{G}_{(\kappa)}, \underline{M}_{(\kappa)}), \] then $H^*_{\operatorname{sing}}(BG,M)$ is concentrated in degree $0$ for all discrete abelian groups $M$. 
    Since $BG$ is connected and $\pi_1(BG)=\pi_0(G)=0$ (we assumed $G$ to be path-connected), it follows from the Hurewicz theorem that $BG$ is weakly contractible. 
    As $G$ is weakly homotopy equivalent to $\Omega BG$, this implies that $G$ is weakly contractible. 
\end{proof}
\subsection{Solid group cohomology}
When working with solid coefficients, it is natural to consider the following variant of group cohomology:
\begin{definition}[Solid group cohomology]
 Suppose $G\in \Grp(\Cond{(\kappa)}(\Set))$ is a ($\kappa$)-condensed group and denote by $\Sol{(\kappa)}(\mathbb Z[G])$ the category of solid $\mathbb Z[G]$-modules (\cref{definitionsolid}).  
 By \cref{moregeneralenrichmentssolid}, $\mathcal D(\Sol{(\kappa)}(\mathbb Z[G]))$ is $\mathcal D(\Ab)$-enriched. Denote by \[\map_{\mathcal D(\Sol{(\kappa)}(\mathbb Z[G]))}(-,-)\colon\mathcal D(\Sol{(\kappa)}(\mathbb Z[G]))^{\operatorname{op}}\times \mathcal D(\Sol{(\kappa)}(\mathbb Z[G]))\to\mathcal D(\Ab)\] the enriched mapping functor and consider $\mathbb Z$ as $\mathbb Z[G]$-module with trivial $G$-action. 

 \emph{Solid group cohomology} of $G$ is defined as  
 \[ \solgrpcoh{(\kappa\text{-})}(G,-)\coloneqq \map_{\mathcal D(\Sol{(\kappa)}(\mathbb Z[G]))}(\mathbb Z,-), \] and the solid cohomology groups are  
 \[ \solgrpcoh{(\kappa\text{-})}^q\coloneqq \pi_{-q}\solgrpcoh{(\kappa\text{-})}(G,-)=\operatorname{Ext}^q_{\Sol{(\kappa)}(\mathbb Z[G])}(\mathbb Z,-), q\in\mathbb Z.\]
\end{definition}
The category $\Sol{(\kappa)}(\mathbb Z[G])$ admits plenty projectives (\cref{solidenoughprojectiveskappa,solidenoughprojectiveswithoutkappa}), which makes solid cohomology more computationally accessible than condensed group cohomology. 
In this section, we compare continuous, condensed and solid group cohomology and show that for a broad class of topological groups, solid group cohomology recovers continuous group cohomology. For locally profinite groups, solid, condensed and continuous group cohomology were identified in \cite{Anschuetzsolidhomology}.

\begin{notation}
Suppose that $G\in\Grp(\Cond{(\kappa)}(\Set))$ is a ($\kappa$-)condensed group. 
We denote by $(-)^{\solid G}\colon \Cond{(\kappa)}(\mathbb Z[G])\to \Sol{(\kappa)}(\mathbb Z[G])$ the left adjoint of the forget functor. 
If the left adjoint of the forget functor \[f\colon \mathcal D(\Sol{(\kappa)}{\mathbb Z[G]})\to \mathcal D(\Cond{(\kappa)}(\mathbb Z[G]))\] exists, we denote it by $(-)^{L\solid G}$ and call it derived solidification. 

We denote by $(-)^{\solid}\colon \Cond{(\kappa)}(\Ab)\to\Sol{(\kappa)}$ and $(-)^{L\solid}\colon\mathcal D(\Cond{(\kappa)}(\Ab))\to\mathcal D(\Sol{(\kappa)})$ the left adjoints of the forget functors. \end{notation} 
Recall from \cref{derivedsolificationexistsstronglimitcardinal}, \cref{derivedsolidificationwithoutkappa}, \cref{lightringssflat,derivedsolidificationexistssflatrings} that the derived solidification exists over $\Cond{}(\mathbb Z[G])$ as well as over $\Cond{\kappa}(\mathbb Z[G])$ for strong limit cardinals $\kappa$ and for $\kappa=\aleph_1$. It also exists if $G=\underline{H}_{(\kappa)}$ for a Hausdorff topological group $H$ (and arbitrary $\kappa$) by \cref{derivedsolidificationexistssflatrings,examplessflatrings}. 

\subsubsection{Solid and condensed group cohomology}
Suppose $G\in\Grp(\Cond{(\kappa)}(\Set))$. 
Recall from \cref{condensedgroupcohomologyisext} that ($\kappa$-)condensed group cohomology enhances to a functor \[ \condgrpcoh{(\kappa\text{-})}(G,-)\colon \mathcal D(\Cond{(\kappa)}(\mathbb Z[G]))\to \mathcal D(\Ab).\] 

We want to compare \[ \condgrpcoh{(\kappa\text{-})}(G,-)\circ f\colon\mathcal D(\Sol{(\kappa)}(\mathbb Z[G]))\to \mathcal D(\Ab)\] with solid group cohomology. 
If the derived solidification $(-)^{L\solid G}$ exists, the counit $ \mathbb Z^{L\solid G}\to \mathbb Z\in \mathcal D(\Sol{(\kappa)}{\mathbb Z[G]})$ induces a map 
\[ \solgrpcoh{\kappa\text{-}}(G,-)\to \map_{\mathcal D(\Sol{\kappa}(\mathbb Z[G]))}(\mathbb Z^{L\solid G},-)\cong \condgrpcoh{\kappa\text{-}}(G,-)\circ f.\] 
\begin{lemma}\label{solidequalscondensedifGacyclic}
    Suppose that $G\in\Grp(\Cond{(\kappa)}(\Set))$ is a ($\kappa$-)condensed group such that the derived solidification $(-)^{L\solid G}$ exists. 
    Then 
    \[ \solgrpcoh{(\kappa\text{-})}(G,-)\cong \condgrpcoh{(\kappa\text{-})}(G,-)\circ f \in \Fun\bigl(\mathcal D(\Sol{(\kappa)}(\mathbb Z[G])), \mathcal D(\Ab)\bigr)\] if and only if $\mathbb Z^{L\solid G}\cong \mathbb Z$. 
\end{lemma}
\begin{proof}
    Denote by $\pi_0\colon\mathcal D(\Sol{(\kappa)}(\mathbb Z[{G}]))_{\geq 0}\to\Sol{(\kappa)}(\mathbb Z[G])$ the left adjoint of the inclusion. Then $\pi_0\circ (-)^{L\solid G}\cong (-)^{\solid G}$ since their right adjoints are equivalent. As $ \mathbb Z^{\solid G}\cong \mathbb Z$, this implies that $\mathbb Z^{L\solid G}\cong \mathbb Z$ if and only if the counit $\mathbb Z^{L\solid G}\to\mathbb Z$ is an equivalence, which proves the if-statement. 

Denote by $g\colon\mathcal D(\Ab)\to \Sp$ the forget functor.  
Since \[\Omega^{\infty}g\solgrpcoh{(\kappa\text{-})}(G,-)\cong \Map_{\mathcal D(\Sol{(\kappa)}(\mathbb Z[G]))}(\mathbb Z,-)\] and 
\[\Omega^{\infty}g\condgrpcoh{(\kappa\text{-})}(G,-)\circ f\cong\Map_{\mathcal D(\Cond{\kappa}(\mathbb Z[G]))}(\mathbb Z,-)\circ f\cong \Map_{\mathcal D(\Sol{(\kappa)}(\mathbb Z[G]))}(\mathbb Z^{L\solid G},-), \] the only-if statement follows from the Yoneda lemma. 
\end{proof} 
\begin{cor}\label{groupacycliccondensedeuqalssolid}
    Suppose that $G\in\Grp(\Cond{(\kappa)}(\Set))$ is such that $\mathbb Z[G]$ is ($\kappa$-)-s-flat (\cref{definitionsflatrings}) and $\mathbb Z[G]^{L\solid}\cong \mathbb Z[G]^{\solid}$. 
    Then \[ \solgrpcoh{(\kappa\text{-})}(G,-)\cong \condgrpcoh{(\kappa\text{-})}(G,-)\circ f \in \Fun\bigl(\mathcal D(\Sol{(\kappa)}(\mathbb Z[G])), \mathcal D(\Ab)\bigr).\] 
\end{cor}
\begin{proof}
Denote by $g^{\solid}\colon \mathcal D(\Sol{(\kappa)}(\mathbb Z[G]))\to\mathcal D(\Sol{(\kappa)})$ and $g\colon \mathcal D(\Cond{(\kappa)}(\mathbb Z[G]))\to\mathcal D(\Cond{(\kappa)}(\Ab))$ the forget functors. 
By \cref{derivedsolidificationexistssflatrings}, the left adjoint $(-)^{L\solid G}$ exists and by \cref{derivedsolidificationunderlyingforacyclicrings}, $g^{\solid}(\mathbb Z^{L\solid G})\cong \mathbb Z^{L\solid}\cong\mathbb Z$. 
Since $g^{\solid}$ is conservative, the statement now follows from \cref{solidequalscondensedifGacyclic}.
\end{proof}
\begin{notation}
    For a topological space $X$ denote by $\pi_0X$ the set of connected components of $X$ and endow it with the quotient topology induced by the projection $X\to\pi_0X$. \end{notation}
\begin{cor}\label{explicitcondensedequalssolid}
Suppose that $G$ is a topological group such that the quotient map $G\to\pi_0G$ is a homotopy equivalence and $\pi_0G$ is locally compact. 

For cardinals $\kappa>\wt(G)$, \[ \solgrpcoh{\kappa\text{-}}(\underline{G}_{\kappa},-)\cong \condgrpcoh{\kappa\text{-}}(\underline{G}_{\kappa},-)\circ f.\]

If $G$ is $\To$, then also 
\[ \solgrpcoh{}(\underline{G},-)\cong \condgrpcoh{}(\underline{G},-)\circ f\] via the above natural transformation. 
\end{cor}
\begin{ex}\cref{explicitcondensedequalssolid} in particular applies to locally profinite groups. The statement for locally profinite groups was shown in \cite{Anschuetzsolidhomology}. 
\end{ex}
\begin{proof}By \cref{weightofconnectedcomponents}, $\pi_0G$ is a totally disconnected, locally compact Hausdorff topological group of $\wt(\pi_0G)\leq \wt(G)$. By van Dantzig's theorem, $\pi_0G$ is a coproduct of profinite sets. In particular, $\mathbb Z[\underline{\pi_0G}_{(\kappa)}]$ is $(-)^{\solid}$-acyclic for all $\kappa>\wt(G)$ by \cref{profiniteprojectiveinsolidkappa}/\cref{freeabeliangroupsacyclicforsolidificationwithoutkappa}.
    
By \cref{solidcohomologyhomotopyinvariant}, \[\mathbb Z[\underline{G}_{(\kappa)}]^{L\solid}\cong\mathbb Z[\underline{\pi_0G}_{(\kappa)}]^{L\solid}\cong \mathbb Z[\underline{\pi_0G}_{(\kappa)}]^{\solid},\] whence $\mathbb Z[\underline{G}_{(\kappa)}]$ is ($\kappa$-)$s$-flat (\cref{definitionsflatrings}) by \cref{examplessflatrings}. In particular, $(-)^{L\solid G}$ exists by \cref{derivedsolidificationexistssflatrings}. 
The statement now follows from \cref{groupacycliccondensedeuqalssolid}.
\end{proof}
\subsubsection{Solid and continuous group cohomology}
\label{section:continuousandsolidgroupcohomology}
\begin{notation}\label{solidcontinuousgmoduledef}
    For an uncountable cardinal $\kappa$ and a topological group $G$ denote by \[\Cont{G}_{\solid, \kappa}\subseteq \Cont{G}\] the full subcategory on continuous $G$-modules $M$ such that $\underline{M}_{\kappa}$ is solid. 

    For a Hausdorff topological group $G$ denote by \[\TCont{G}_{\solid}\subseteq \TCont{G}\] the full subcategory on $\To$ continuous $G$-modules $M$ such that $\underline{M}$ is solid. 
\end{notation}
The goal of this section is to prove the following result:
\begin{thm}\label{solidequalscontinuousgoodgroups}
    Suppose that $G$ is a topological group which is a finite product of topological groups $\{G_i\}_{i=1, \ldots,n}$ each of which satisfies one of the following properties, respectively: 
    \begin{romanenum}\label{goodgroups}
    \item $G_i$ is homotopy equivalent to a coproduct of ($\kappa$-light) compact Hausdorff spaces.
    \item $G_i$ is homotopy equivalent to a Hausdorff space and to a locally contractible space. 
    \item $G_i$ is locally connected and locally ($\kappa$-light) compact Hausdorff.
    \end{romanenum}

    Then there is a natural isomorphism \[ \contgrpcoh{\kappa\text{-}}^*(G,-)\cong \solgrpcoh{}^*(\underline{G}_{\kappa},-)\circ \underline{(-)}_{\kappa}\in \Fun(\Cont{G}_{\solid, \kappa}, \grAb).\]
    If in addition to the above, $G$ is Hausdorff, then \[\contgrpcoh{k\text{-}}^*(G,-)\cong \solgrpcoh{}^*(\underline{G},-)\circ \underline{(-)}\in \Fun(\TCont{G}_{\solid}, \grAb).\]
\end{thm}
\begin{rem}\begin{romanenum} The conditions of \cref{goodgroups} are satisfied by a broad class of topological groups: 
        \item If $G$ is locally compact abelian, then $G\cong \mathbb R^n\times H$ where $H$ admits a compact open subgroup (see e.g.\ \cite[Theorem 4]{HofmannMorris+2023}), whence $G$ satisfies \cref{goodgroups}. 
        \item Lie groups are locally contractible and hence satisfy \cref{goodgroups}. 
        \item By van Dantzig's theorem, locally profinite groups are coproducts of compacta, whence \cref{solidequalscontinuousgoodgroups} applies to them as well.
        The statement for locally profinite groups was established in \cite{Anschuetzsolidhomology}. \end{romanenum}
The author is unaware whether there exist Hausdorff topological groups whose solid group cohomology is not isomorphic to their $k$-continuous group cohomology. 
\end{rem}
 
\cref{solidequalscontinuousgoodgroups} is an immediate consequence of \cref{solidequalscontinuousifzgprojective} and \cref{goodgroups} below. It follows from the following two observations:  
\begin{romanenum}
\item For large classes of topological groups, including those satisfying \cref{goodgroups}, the degreewise solidification $(S_*^{\underline{G}_{(\kappa)}})^{\solid \underline{G}_{(\kappa)}}$ of the simplicial resolution is a resolution of $\mathbb Z$ (\cref{solidcohomologysimplicialresolution}).
\item Under the assumptions of \cref{solidequalscontinuousgoodgroups}, this solidified resolution is a projective resolution of $\mathbb Z$ in $\Sol{(\kappa)}(\mathbb Z[\underline{G}_{(\kappa)}])$ (\cref{goodgroups1}). 
\end{romanenum}   
\begin{lemma}\label{solidcohomologysimplicialresolution}
    Suppose $G$ is a topological group. 
    \begin{romanenum}    
    \item If $G$ is homotopy equivalent to a Hausdorff space, the degreewise (underived) solidification of the simplicial complex $(S_*^{\underline{G}_{\kappa}})^{\solid\underline{G}_{\kappa}}$ is a resolution of $\mathbb Z$ in $\Sol{\kappa}(\mathbb Z[\underline{G}_{\kappa}])$.  
    \item If $G$ is Hausdorff, the degreewise (underived) solidification of the simplicial complex $(S_*^{\underline{G}})^{\solid\underline G}$ is a resolution of $\mathbb Z$ in $\Sol{}(\mathbb Z[\underline{G}])$.
    \end{romanenum} 
\end{lemma} 

\begin{proof}
    Since $\Sol{(\kappa)}(\mathbb Z[\underline{G}_{(\kappa)}])\cong \LMod{\mathbb Z[\underline{G}_{(\kappa)}]^{\solid}}{\Sol{(\kappa)}}$ and the symmetric monoidal structure on $\Sol{(\kappa)}$ is cocontinuous in both variables (the symmetric monoidal structure is closed \cref{monoidalstructureonsolidabeliangroupsclosed}), the forget functor $\Sol{(\kappa)}(\mathbb Z[\underline{G}_{(\kappa)}])\to \Sol{(\kappa)}$ creates small limits and colimits by \cref{forgetfreeadjunctionmodules}. 
    It therefore suffices to show that the chain complex of solid abelian groups underlying $(S_*^{\underline G})^{\solid\underline G_{(\kappa)}}$ is a resolution of $\mathbb Z$. 
    By \cref{underivedsolidification}/\cref{underivedsolidification}, this chain complex is given by applying $(-)^{\solid}$, the left adjoint of $\Sol{(\kappa)}\subseteq \Cond{(\kappa)}(\Ab)$, degreewise to $S_*^{\underline{G}}$.
    For a $(\kappa$-)condensed set $X$ denote by $S_*^{X}\in \Ch(\Cond{(\kappa)}(\Ab))$ the simplicial resolution of $X\to *$ (\cref{definitionsimplicialcomplex}). 
    A homotopy equivalence $G\to X$ induces a chain map $S_*^{\underline{G}_{\kappa}}\to S_*^{\underline{X}_{\kappa}}$, and by \cref{solidcohomologyhomotopyinvariant}, the induced map $(S_*^{\underline{G}_{\kappa}})^{\solid}\to (S_*^{\underline{X}_{\kappa}})^{\solid}$ is an isomorphism of chain complexes. It therefore suffices to show that for a Hausdorff space $X$, $(S_*^{\underline{X}_{(\kappa)}})^{\solid}$ is a resolution of $\mathbb Z$. 
    Fix $n\in\mathbb N_1$ and denote by $\pi_i\colon X^n\to X, 1\leq i\leq n$ the projections.  
    For a compact Hausdorff space $K$ and a continuous map $f\colon K\to X^i$, \[f(K)\subseteq (\cup_{i=1}^n\pi_if(K))^n\subseteq X^n,\] and $\cup_{i=1}^n\pi_i(f(K))\subseteq X$ is Hausdorff (since $X$ is) and compact as quotient of $\sqcup_{i=1}^n K$. Hence by \cref{weightofquotients}, $\wt(\cup_{i=1}^n\pi_i(K))\leq \wt(\sqcup_{i=1}^n K)\leq \wt(K)$. 
    This shows that for $n\in\mathbb N_1$, \[\colim{\substack{K\subseteq X\\ K \in \CH_{(\kappa)}}} \underline{K^n}_{(\kappa)}\cong \colim{\substack{K\to X^n\\ K\in \operatorname{CH}_{\kappa}}}\underline{K}\cong \underline{X^n}_{(\kappa)}\] (\cref{approximationbycompacty}), and hence \[S_*^{\underline{X}_{(\kappa)}}\cong \colim{\substack{K\subseteq X\\ K\in \CH_{(\kappa)}}} S_*^{\underline{K}_{(\kappa)}}.\] 
    As solidification is a left adjoint, it follows from \cref{solidificationcompacthausdorffspace} that  
    \begin{equation}\label{simplicial as colimit}
        (S_*^{\underline X_{(\kappa)}})^{\solid} \cong \colim{\substack{K\subseteq X \\ K\in \CH_{(\kappa)}}} (S_*^{\underline{X}_{(\kappa)}})^{\solid}\cong \colim{\substack{K\subseteq X \\ K \in \operatorname{CH}_{(\kappa)}}} (S_*^{\underline{\pi_0K}_{(\kappa)}})^{\solid},
    \end{equation} where $\pi_0K$ denotes the topological space of connected components of $K$ with quotient topology from $K\to\pi_0K$. By \cref{connectedcomponentsofchisch}, $\pi_0K$ is profinite and by \cref{weightofquotients}, $\wt(\pi_0K)\leq \wt(K)$.
    By \cref{profiniteprojectiveinsolidkappa}/\cref{freeabeliangroupsacyclicforsolidificationwithoutkappa} for $P\in \Pro(\Fin)_{(\kappa)}$, 
    $S_*^{\underline P_{(\kappa)}}$ is a $(-)^{\solid}$-acyclic resolution of $\mathbb Z$ in $\Cond{(\kappa)}(\Ab)$, and hence \[(S_*^{\underline {P}_{(\kappa)}})^{\solid}\cong \mathbb Z^{L\solid}\cong \mathbb Z \in \mathcal D(\Sol{(\kappa)}),\] cf.\ \cref{derivedsolidificationisderivedfunctorofsolidification}. 
    This means that for $P\in\Pro(\Fin)_{(\kappa)}$, the degreewise solidification $(S_*^{\underline{P}_{(\kappa)}})^{\solid}$ of $S_*^{\underline{P}_{(\kappa)}}$ is a resolution of $\mathbb Z$. 
    Since filtered colimits in $\Cond{(\kappa)}(\Ab)$ are exact (\cref{condensedabeliangroupsgrothendieckaxioms}), and $\Sol{(\kappa)}\subseteq \Cond{(\kappa)}(\Ab)$ is closed under limits and colimits (\cref{solidclosedunderlimitscolimitskappa}, \cref{underivedsolidificationwithoutkappa}), filtered colimits in $\Sol{(\kappa)}$ are exact. 
    It now follows that for $X$ Hausdorff, $(S_*^{\underline{X}_{(\kappa)}})^{\solid}$ is a resolution of $\mathbb Z$.
\end{proof}
Together with \cref{continuouscohomologysimplicialcomplex}, this yields a comparison map from $\kappa$-/$k$-continuous group cohomology with solid coefficients to solid group cohomology for Hausdorff topological groups:  
\begin{recollection}
    Suppose $G$ is a Hausdorff topological group. 
    As $\Sol{(\kappa)}(\mathbb Z[\underline{G}_{(\kappa)}])$ has enough projectives (\cref{solidenoughprojectiveskappa,solidenoughprojectiveswithoutkappa}), there exists a projective resolution $P_*\to\mathbb Z$ in $\Sol{(\kappa)}(\mathbb Z[\underline{G}_{(\kappa)}])$. This computes solid group cohomology (see e.g.\ \cite[\href{https://stacks.math.columbia.edu/tag/06XR}{Tag 06XR}]{stacks-project}), i.e.\ \[H^*(\Hom_{\Sol{(\kappa)}(\mathbb Z[\underline{G}_{(\kappa)}])}(P_*,-))\cong \solgrpcoh{(\kappa\text{-})}^*(\underline{G}_{(\kappa)},-).\] 
    As $(S_*^{\underline{G}_{(\kappa)}})^{\solid \underline{G}_{(\kappa)}}$ is a resolution of $\mathbb Z$ (\cref{solidcohomologysimplicialresolution}), the identity $\id_{\mathbb Z}$ lifts to a chain map $P_*\to (S_*^{\underline{G}_{(\kappa)}})^{\solid \underline{G}_{(\kappa)}}$ (\cite[Proposition 2.2.6]{Weibel_1994}). 
    By \cref{continuouscohomologysimplicialcomplex}, \begin{align*} \contgrpcoh{\kappa\text{-}}^*(G,-)|_{\Cont{G}_{\solid, \kappa}}&\cong H^*(\Hom_{\Cond{\kappa}(\mathbb Z[\underline{G}_{\kappa}])}(S_*^{\underline{G}_{\kappa}}, \underline{-}_{\kappa}))|_{\Cont{G}_{\solid, \kappa}}\\ &\cong H^*(\Hom_{\Sol{\kappa}(\mathbb Z[\underline{G}_{\kappa}])}((S_*^{\underline{G}_{\kappa}})^{\solid\underline{G}_{\kappa}}, \underline{-}_{\kappa}))|_{\Cont{G}_{\solid, \kappa}}\end{align*} and \begin{align*} \contgrpcoh{k\text{-}}^*(G,-)|_{\TCont{G}_{\solid}} & \cong H^*(\Hom_{\Cond{}(\mathbb Z[\underline{G}_{}])}(S_*^{\underline{G}}, \underline{-}))|_{\TCont{G}_{\solid}}\\& \cong H^*(\Hom_{\Sol{}(\mathbb Z[\underline{G}_{}])}((S_*^{\underline{G}})^{\solid\underline{G}}, \underline{-}))|_{\TCont{G}_{\solid}}.\end{align*} Pullback along the chain map $P_*\to (S_*^{\underline{G}_{(\kappa)}})^{\solid \underline{G}_{(\kappa)}}$ therefore defines natural transformations \[ \contgrpcoh{\kappa\text{-}}(G,-)\longrightarrow H^*(\Hom_{\Cond{\kappa}(\mathbb Z[\underline{G}_{\kappa}])}(P_*,-))\cong \solgrpcoh{\kappa\text{-}}^*(\underline{G}_{\kappa}, \underline{(-)}_{\kappa})\in \Fun(\Cont{G}_{\solid, \kappa}, \grAb)\] and
\[ \contgrpcoh{k\text{-}}(G,-)\longrightarrow H^*(\Hom_{\Cond{}(\mathbb Z[\underline{G}_{}])}(P_*,-))\cong \solgrpcoh{}^*(\underline{G}_{}, \underline{(-)})\in \Fun(\TCont{G}_{\solid}, \grAb),\] respectively. 
\end{recollection}
\begin{rems}
The natural transformation \[ \contgrpcoh{\kappa\text{-}}(G,-)\longrightarrow H^*(\Hom_{\Cond{\kappa}(\mathbb Z[\underline{G}_{\kappa}])}(P_*,-))\cong \solgrpcoh{\kappa\text{-}}^*(\underline{G}_{\kappa}, \underline{(-)}_{\kappa})\in \Fun(\Cont{G}_{\solid, \kappa}, \grAb)\] exists more generally if $G$ is not Hausdorff, but merely homotopy equivalent to a Hausdorff topological space. 
Since the chain map $P_*\to (S_*^{\underline{G}_{(\kappa)}})^{\solid \underline{G}_{(\kappa)}}$ is unique up to chain homotopy equivalence (\cite[Proposition 2.2.6]{Weibel_1994}), the comparison map is independent of choices.

The composition \[ \contgrpcoh{\kappa/k\text{-}}^*(G,-)|_{\Sol{(\kappa)}(\mathbb Z[\underline{G}_{(\kappa)}])}\to\solgrpcoh{\kappa\text{-}}(\underline{G}_{(\kappa)}, \underline{(-)}_{(\kappa)})\to \condgrpcoh{}^*(\underline{G}_{\kappa}, \underline{(-)}_{(\kappa)})\] is the map induced by the edge homomorphism from \cref{spectralsequenceidentifycontinuousgroupcohomology}. 
\end{rems}
\begin{lemma}\label{criterionprojectiveresolution}For a ($\kappa$-)condensed group $G\in \Grp(\Cond{(\kappa)}(\Set))$, the degreewise solidification $(S_*^{G})^{\solid G}$ of the simplicial resolution is a chain complex of projectives in $\Sol{(\kappa)}(\mathbb Z[G])$ if and only if $\mathbb Z[G]^{\solid}$ is projective in $\Sol{(\kappa)}$. 
\end{lemma}
\begin{proof}
    Denote by \[F\colon \Sol{(\kappa)}\to \Sol{(\kappa)}(\mathbb Z[G])\cong \LMod{\mathbb Z[G]^{\solid}}{\Sol{(\kappa)}}\] the free $\mathbb Z[G]^{\solid}$-module functor (\cref{solidmodulesismodulesinsolid}). 
    Since the forget functor (the right adjoint of $F$) reflects exact sequences, a solid abelian group $P\in\Sol{(\kappa)}$ is projective in $\Sol{\kappa}$ if and only if $F(P)$ is projective in $\Sol{(\kappa)}(\mathbb Z[\underline{G}_{(\kappa)}])$.
    In particular, 
    \[(S_i^{G})^{\solid G}=\mathbb Z[{G^{i+1}}]^{\solid G}\cong F(\mathbb Z[{G^i}]^{\solid})\] is projective in $\Sol{(\kappa)}(\mathbb Z[G])$ if and only if 
    $\mathbb Z[G^i]^{\solid}\cong \otimes^i_{\Sol{(\kappa)}}\mathbb Z[G]^{\solid}$ is projective in $\Sol{(\kappa)}$. Since tensor products of projectives are projective in $\Sol{(\kappa)}$ (\cref{solidprojectivestensorproduct}), this holds for all $i\in\mathbb N_0$ if and only if $\mathbb Z[G]^{\solid}$ is a projective solid abelian group. 
\end{proof} 
\begin{cor}\label{solidequalscontinuousifzgprojective}
    Suppose $G$ is a topological group. 
    \begin{romanenum}
        \item If $G$ is homotopy equivalent to a Hausdorff space and $\kappa$ is an uncountable cardinal such that $\mathbb Z[\underline{G}_{\kappa}]^{\solid}$ is projective in $\Sol{\kappa}$, then \[\solgrpcoh{\kappa\text{-}}^*(\underline G_{\kappa},-)\cong H^*(\Hom_{\Cond{\kappa}(\mathbb Z[\underline{G}_{\kappa}])}(S_*^{\underline G_{\kappa}},-))|_{\Sol{\kappa}(\mathbb Z[\underline{G}_{\kappa}])}.\] 
         \item If $G$ is a Hausdorff topological group and  $\mathbb Z[\underline{G}]^{\solid}$ is projective in $\Sol{}$, then \[\solgrpcoh{}^*(\underline G,-)\cong H^*(\Hom_{\Cond{}(\mathbb Z[\underline{G}])}(S_*^{\underline G},-))|_{\Sol{}(\mathbb Z[\underline{G}])}.\] 
    \end{romanenum} 
    In particular, under the above assumptions on $G$, \[\contgrpcoh{\kappa\text{-}}^*(G,-)\cong \solgrpcoh{\kappa\text{-}}^*(\underline{G},-)\circ \underline{(-)}_{\kappa}\colon \Cont{G}_{\solid, \kappa}\to \grAb, \] and \[\contgrpcoh{k\text{-}}^*(G,-)\cong \solgrpcoh{}^*(\underline{G},-)\circ \underline{(-)}\colon \TCont{G}_{\solid}\to \grAb, \] respectively.  
\end{cor} 
\begin{proof}
    By \cref{solidcohomologysimplicialresolution} and \cref{criterionprojectiveresolution}, our assumptions imply that $(S_*^{\underline{G}_{(\kappa)}})^{\solid\underline{G}_{(\kappa)}}$ is a projective resolution of $\mathbb Z$, whence  
    \[ \solgrpcoh{(\kappa\text{-})}^*(\underline{G}_{(\kappa)},-)\cong H^*(\Hom_{\Sol{\kappa}(\mathbb Z[\underline{G}_{(\kappa)}])}((S_*^{\underline{G}_{(\kappa)}})^{\solid\underline{G}_{(\kappa)}},-)\] by e.g.\ \cite[\href{https://stacks.math.columbia.edu/tag/06XR}{Tag 06XR}]{stacks-project}. 
    The right-hand side is isomorphic to \[H^*(\Hom_{\Cond{(\kappa)}(\mathbb Z[\underline{G}_{(\kappa)}])}(S_*^{\underline{G}_{(\kappa)}},-)|_{\Sol{(\kappa)}(\mathbb Z[\underline{G}_{(\kappa)}])}\] by definition of $(-)^{\solid\underline{G}_{(\kappa)}}$. 
    The comparison with $k/\kappa$-continuous group cohomology is  \cref{continuouscohomologysimplicialcomplex}.
\end{proof}

\cref{solidequalscontinuousgoodgroups} is now a consequence of the following observation: 
 \begin{lemma}\label{goodgroups1}
    Suppose that $G$ is a topological group which is a finite product of topological groups $\{G_i\}_{i=1, \ldots,n}$ each of which satisfies one of the following properties, respectively.
    \begin{romanenum}
    \item \label{Gcoproductcompactagoodgroups}$G_i$ is homotopy equivalent to a coproduct of ($\kappa$-light) compact Hausdorff spaces.
    \item \label{Glocallycontractiblegoodgroups} $G_i$ is homotopy equivalent to a locally contractible topological space.
    \item \label{Glocallyconnectedlocallycompact}$G_i$ is locally connected and locally ($\kappa$-light) compact Hausdorff. 
    \end{romanenum}
    Endow $\pi_0G$ with the quotient topology from the projection $G\to\pi_0G$. 

    Then the following hold: 
    \begin{itemize}[label={--}]
    \item $\pi_0G$ is a coproduct of ($\kappa$-light) compact Hausdorff spaces. 
    \item The quotient map induces an isomorphism $\mathbb Z[\underline{G}_{\kappa}]^{\solid}\cong \mathbb Z[\underline{\pi_0G}_{\kappa}]^{\solid}$ and $\mathbb Z[\underline{G}_{\kappa}]^{\solid}$ is projective in $\Sol{\kappa}$. 
    \item If in addition to the above, $G$ is Hausdorff, then $\mathbb Z[\underline{G}]^{\solid} \cong \mathbb Z[\underline{\pi_0G}]^{\solid}$ is projective in $\Sol{}$.  
    \end{itemize}
 \end{lemma}
\begin{proof}
    We first reduce to the case that $G$ satisfies one of the above conditions. 
    Write \[G=\prod_{k=1}^3 \prod_{i=1}^{n_i}G_{k,i}\] as a product of topological groups such that $G_{k,i}$ satisfies property $k$ of \cref{goodgroups1}. 
    Then $\prod_{i=1}^{n_k}G_{i,k}$ satisfies property $(k)$ for $1\leq k\leq 3$, so we can assume that $G=\prod_{k=1}^3G_{k}$ for $G_k$ a group satisfying property $(k)$.
    We first show that \[\pi_0(G_1\times G_2\times G_3)\to\pi_0G_1\times \pi_0G_2\times \pi_0G_3\] is a homemorphism, i.e. that $G_0\times G_1\times G_3\to \pi_0G_1\times \pi_0G_2\times \pi_0G_3$ is a quotient map. If $h\colon G_2\to X$ is a homotopy equivalence to a locally contractible topological space, then $h$ descends to a continuous bijection $\pi_0G\to\pi_0X$. As $\pi_0X$ is discrete, this implies that $\pi_0G$ is discrete, i.e. $G_2$ is locally connected. 
    This implies that $G_4\coloneqq G_2\times G_3$ is locally connected. Denote by $C\subseteq G_4$ the connected component of the identity. Then $G_4\cong C\times \pi_0G_4$ as topological spaces (not necessarily as groups). 
    In particular, $G_1\times G_4\to G_1\times \pi_0G_4$ is a quotient map. As $\pi_0G_4$ is discrete, $G_1\times \pi_0G_4\to \pi_0G_1\times \pi_0G_4$ is a quotient map as well, which shows that $G_1\times G_2\times G_3\to \pi_0G_1\times \pi_0G_2\times \pi_0G_3$ is a quotient map. 

    As \[\mathbb Z[\underline{(-)}_{(\kappa)}]^{\solid}\colon (\To)\Top\to \Sol{(\kappa)}(\Ab)\] is symmetric monoidal and projective solid abelian groups are stable under solid tensor products, we are now reduced to proving the statement for groups satisfying one of the conditions $i)-iii)$ of the lemma, respectively. We do this case-by-case. 
    \begin{romanenum}
    \item If $h\colon G\to \sqcup_{i\in I}K(i)$ is a homotopy equivalence to a coproduct of $(\kappa\text{-})$light compact Hausdorff spaces, then \[\mathbb Z[\underline{G}_{(\kappa)}]^{\solid}\cong \oplus_{i\in I}\mathbb Z[\underline{K(i)}_{(\kappa)}]^{\solid}\] by \cref{solidcohomologyhomotopyinvariant}. By \cref{solidificationcompacthausdorffspace}, $\mathbb Z[\underline{K(i)}_{(\kappa)}]^{\solid}\cong \oplus_{i\in I}\mathbb Z[\underline{\pi_0K(i)}_{(\kappa)}]^{\solid}$, where $\pi_0K(i)$ is endowed with quotient topology.
    The space $\pi_0K(i)$ is profinite (\cref{connectedcomponentsofchisch}), and by \cref{weightofquotients}, $\wt(\pi_0K(i))\leq\wt(K(i))$. \cref{solidenoughprojectiveskappa,solidenoughprojectiveswithoutkappa} now imply that $\mathbb Z[\underline{\pi_0K(i)}_{(\kappa)}]^{\solid}$ is projective in $\Sol{(\kappa)}$.
    It remains to show that $h$ descends to a homeomorphism $\pi_0G\cong \sqcup_{i\in I}\pi_0K(i)$. 
    Choose a homotopy inverse $f$ of $h$. Then $G=\sqcup_{i\in I}f^{-1}(K(i))$ as topological space, which implies that $f$ and $h$ restrict to mutually inverse homotopy equivalences $K(i)\simeq f^{-1}(G_i)$ for all $i\in I$. The map $h$ descends to continuous bijections $h_i\colon \pi_0K(i)\to \pi_0f^{-1}(G_i), i\in I$, which shows that $\pi_0G=\sqcup_{i\in I}\pi_0f^{-1}(K(i))$ is a coproduct of compacta. As $\pi_0G$ is Hausdorff (\cref{weightofconnectedcomponents}), $\pi_0f^{-1}(K(i))$ is Hausdorff for all $i\in I$. In particular, the maps $h_i$ are homeomorphisms. 
    \item Suppose now that $G$ is homotopy equivalent to a locally contractible topological space. We explained above that $\pi_0G$ is discrete in this case.
    By \cref{condensedandsheafcohomologylocallycontractible}, for all solid abelian groups $M$, \begin{align*}\Hom_{\Sol{(\kappa)}}(\mathbb Z[\underline{G}_{(\kappa)}]^{\solid},M)&\cong \Hom_{\Cond{(\kappa)}(\Ab)}(\mathbb Z[\underline{G}_{(\kappa)}],M)\\ &\cong \cckH{(\kappa\text{-})}^0(\underline{G}_{(\kappa)},M)\\ &\cong \cH{\sheaf}^0(G,M(*)^{\delta})\\ &\cong \mathcal C(G,M(*)^{\delta})\\ &\cong \mathcal C(\pi_0G,M(*)^{\delta}),\end{align*} and 
    \begin{align*}\Hom_{\Sol{(\kappa)}}(\mathbb Z[\underline{\pi_0G}_{(\kappa)}]^{\solid},M)&\cong \Hom_{\Cond{(\kappa)}(\Ab)}(\mathbb Z[\underline{\pi_0G}_{(\kappa)}],M)\\ &\cong \cckH{(\kappa\text{-})}^0(\pi_0G,M)\\ & \cong\cH{\sheaf}^0(\pi_0G,M(*)^{\delta})\\ &\cong \mathcal C(\pi_0G,M(*)^{\delta}), \end{align*} whence $\mathbb Z[\underline{G}_{(\kappa)}]^{\solid}\cong\mathbb Z[\underline{\pi_0G}_{(\kappa)}]^{\solid}\cong \oplus_{\pi_0G}\mathbb Z$ is projective in $\Sol{(\kappa)}$.

  \item Suppose now that $G$ is locally ($\kappa$-light) compact Hausdorff and locally connected and denote by $G^0$ the connected component of the identity. 
  Then $\pi_0G$ is discrete and $G\cong \sqcup_{x\in \pi_0G}G^0$ as topological space. This implies that $\mathbb Z[\underline{G}_{(\kappa)}]^{\solid}\cong \oplus_{\pi_0G}\mathbb Z[\underline{G^0}_{(\kappa)}]^{\solid}$ and the quotient map \[\mathbb Z[\underline{G}_{(\kappa)}]^{\solid}\to \oplus_{\pi_0G}\mathbb Z=\mathbb Z[\underline{\pi_0G}_{(\kappa)}]^{\solid}\] is the coproduct of the projections $\mathbb Z[\underline{G^0}_{(\kappa)}]^{\solid}\to\mathbb Z^{\solid}$. It is therefore enough to show that $\mathbb Z[\underline{G^0}_{(\kappa)}]^{\solid}=\mathbb Z$. 
  Below, we will construct an ascending sequence $K_n\subseteq K_{n+1}$ of ($\kappa$-light) compact connected subspaces such that $G^0=\bigcup_{n\in\mathbb N_0}K_n$ and for all compact subspaces $C\subseteq G^0$, there exists $N\in\mathbb N_0$ with $C\subseteq K_N$. 
  Then it follows from \cref{approximationbycompacty} that \[\underline{G^0}_{(\kappa)}\cong \colim{\substack{C\subseteq G^0\\ C\in \CH_{(\kappa)}}}\underline{C}_{(\kappa)}\cong \colim{n\in\mathbb N_0}\underline{K_n}_{(\kappa)}\] and in particular $\mathbb Z[\underline{G^0}_{(\kappa)}]^{\solid}=\colim{n\in\mathbb N_0}\mathbb Z[\underline{K_n}_{(\kappa)}]^{\solid}\cong \colim{n\in\mathbb N_0}\mathbb Z=\mathbb Z$ by \cref{solidificationcompacthausdorffspace}.

We now construct such an exhaustion by compact, connected subspaces. Choose a ($\kappa$-light) compact neighborhood $e_G\in V\subseteq G$ of the identity and an open, connected neighborhood $e_G\in U\subseteq V^{\circ}$, where $V^{\circ}$ denotes the interior of $V$. Then $\overline{U}\subseteq V$, $\overline{U}$ is compact Hausdorff and connected (as closure of a connected set) and $\wt(\overline{U})\leq \wt(V)$. Since \[\tau\colon G\times G\to G, (g,h)\mapsto gh^{-1}\] is continuous, $K_0\coloneqq \tau(\overline{U}\times\overline{U})\subseteq G$ is a compact connected subspace. 
As continuous surjections between compact Hausdorff spaces are quotient maps, \[\wt(K_0)\leq \wt(\overline{U})\leq \wt(V)\] by \cref{weightofquotients}. Define recursively \[K_{n+1}\coloneqq \tau(K_n\times K_n)\subseteq G.\] 
It follows by induction that $K_{n}$ is compact and connected and that $\wt(K_n)\leq \wt(V)$ (\cref{weightofquotients}). Moreover, $K_n=\tau(K_n\times\{1\})\subseteq K_{n+1}$. 
    If $x,y\in K_{\infty}\coloneqq \bigcup_{n\in\mathbb N_0}K_n$ choose $n\in\mathbb N_0$ with $x,y\in K_n$. Then $xy^{-1}\in K_{n+1}$, which shows that $K_{\infty}\subseteq G^0$ is a subgroup. 
    Since $U\subseteq K_{\infty}$ and $K_{\infty}$ is a subgroup, $K_{\infty}=\bigcup_{k\in K_{\infty}}kU$ which shows that $K_{\infty}\subseteq G^0$ is an open subgroup. This implies that $G^0=\sqcup_{g\in G^0/K_{\infty}}gK_{\infty}$, whence $G^0=K_{\infty}$ by connectedness of $G^0$. Since $G^0=\bigcup_{k\in K_{\infty}}Uk^{-1}$ is an open cover, for $C\subseteq G$ compact, there exists a finite subset $F\subseteq K_{\infty}$ such that $C\subseteq \bigcup_{f\in F}^n Uf^{-1}$. 
    If $N\in\mathbb N_0$ is such that $F\subseteq K_N$, then $C\subseteq \tau(U\times K_N)\subseteq K_{N+1}$. \qedhere
\end{romanenum}
\end{proof}

For topological groups as in \cref{goodgroups}, 
solid/continuous group cohomology can be identified with the continuous group cohomology of the group $\pi_0G$ (endowed with the quotient topology from $G\to\pi_0G$):  
 \begin{lemma}
    Suppose $G$ is a topological group as in \cref{goodgroups}. 
    \begin{romanenum}
    \item Restriction along $p\colon G\to\pi_0G$ yields an equivalence \[p^*\colon \Cont{\pi_0G}_{\solid, \kappa}\cong \Cont{G}_{\solid, \kappa}, \] and 
\[ \contgrpcoh{}^*(\pi_0G,-)\cong \contgrpcoh{\kappa\text{-}}^*(G,-)\circ p^*\colon \Cont{\pi_0G}_{\solid, \kappa}\to\grAb.\] 
\item If $G$ is in addition Hausdorff, then restriction along $p\colon G\to\pi_0G$ is an equivalence \[ p^*\colon \TCont{\pi_0G}_{\solid}\cong \TCont{G}_{\solid}, \] and \[\contgrpcoh{}^*(\pi_0G,-)\cong \contgrpcoh{k\text{-}}^*(G,-)\circ p^* \colon \TCont{\pi_0G}_{\solid}\to \grAb.\] 
    \end{romanenum}
 \end{lemma}
 \begin{rem}
 In contrast to the above condensed group cohomology also depends on the connected component $G^0$, see the discussion around \cref{condensedequalscontinuousimpliescontractible}. This yields many examples where condensed group cohomology differs from solid group cohomology. 
 \end{rem}
\begin{proof}We only prove the $\kappa$-condensed statement, the condensed statement can be shown completely analogously.
Since $p\colon G\to\pi_0G$ is surjective, 
\[p^*\colon \Cont{\pi_0G}\to \Cont{G}\] is fully faithful. Suppose that $M\in \Cont{G}$ is such that $\underline{M}_{\kappa}\in\Sol{\kappa}$. 
The action map $G\times M\to M$ induces a map of condensed abelian groups $\mu\colon \mathbb Z[\underline{G}_{\kappa}]\otimes_{\mathbb Z}\underline{M}_{\kappa}\to \underline{M}_{\kappa}$. 
Since $\underline{M}_{\kappa}$ is solid, this extends to a map \[\mu^{\solid}\colon (\mathbb Z[\underline{G}_{\kappa}]\otimes_{\mathbb Z}\underline{M}_{\kappa})^{\solid}\cong\mathbb Z[\underline{G}_{\kappa}]^{\solid}\otimes_{\mathbb Z}^{\solid}\underline{M}_{\kappa}\to \underline{M}_{\kappa}.\] 
By \cref{goodgroups1}, $\mathbb Z[\underline{G}_{\kappa}]^{\solid}\cong\mathbb Z[\underline{\pi_0G}_{\kappa}]^{\solid}$ via the projection, this implies that $\mu$ factors over a map \[\mathbb Z[\underline{\pi_0G}_{\kappa}]\otimes_{\mathbb Z}\underline{M}_{\kappa}\to \underline{M}_{\kappa}, \] and hence the $G$-action $\underline{G}_{\kappa}\times \underline{M}_{\kappa}\to \underline{M}_{\kappa}$ factors as $\underline{G}_{\kappa}\times \underline{M}_{\kappa}\to \underline{\pi_0G}_{\kappa}\times \underline{M}_{\kappa}\to \underline{M}_{\kappa}$. 
The underlying map $\pi_0G\times M\to M$ is $\kappa$-continuous and hence continuous as $\pi_0G$ is $\kappa$-compactly generated by \cref{goodgroups1}. 
This shows that $p^*\colon \Cont{\pi_0G}_{\solid, \kappa}\to \Cont{G}_{\solid, \kappa}$ is essentially surjective.  

By \cref{kappacontinuousfullyfaithfullyintocondensed} and \cref{goodgroups1}, for a topological abelian group $M$ with $\underline{M}_{\kappa}\in\Sol{\kappa}$, 
\begin{align*}\mathcal C_{\kappa}(G^i,M)&\cong \Hom_{\Cond{\kappa}(\Ab)}(\mathbb Z[\underline{G^i}_{\kappa}],M)\\ &\cong \Hom_{\Sol{\kappa}}(\mathbb Z[\underline{G^i}_{\kappa}]^{\solid},M)\\ 
    &\cong \Hom_{\Sol{\kappa}}(\mathbb Z[\underline{(\pi_0G^i)}_{\kappa}]^{\solid}, \underline{M}_{\kappa}) \\&\cong \Hom_{\Cond{\kappa}(\Ab)}(\mathbb Z[\underline{\pi_0(G^i)}_{\kappa}], \mathbb Z)\\&\cong\mathcal C_{\kappa}(\pi_0(G^i),M)\end{align*} via the projection $G^i\to\pi_0G^i$.
    We have shown in the proof of \cref{goodgroups} that $\pi_0(G^i)\cong (\pi_0G)^i$ for all $i\in\mathbb N_0$, and hence 
    $\contgrpcoh{\kappa\text{-}}(G,-)\circ p^*\cong \contgrpcoh{\kappa\text{-}}(\pi_0G,-)\circ p^*$.
    As $\pi_0(G^i)$ is $\kappa$-compactly generated for all $i\in\mathbb N_1$, 
     \[ \contgrpcoh{}^*(\pi_0G,-)\cong \contgrpcoh{\kappa\text{-}}^{*}(\pi_0G,-).\qedhere\]
\end{proof}
\newpage
\appendix
\section{Appendix}
\subsection{Adjunctions and filtered colimits of categories}
In this section, we record some basic results on adjunctions and filtered colimits of categories which we repeatedly used to generalise results from presentable to big presentable categories. 
For categories $\mathcal C, \mathcal D$ denote by $\Fun^{L}(\mathcal C, \mathcal D)\subseteq\Fun(\mathcal C, \mathcal D)$ the full subcategory on left adjoint functors. 
\begin{lemma}\label{leftadjointsstableundercolimits}
    Suppose that $\mathcal C, \mathcal D$ are categories and $I$ is a category such that $\mathcal D$ admits $I$-indexed colimits and $\mathcal C$ admits $I^{\operatorname{op}}$-indexed limits. Then the following hold: 
    \begin{romanenum}
        \item $\Fun^{L}(\mathcal C, \mathcal D)\subseteq \Fun(\mathcal C, \mathcal D)$ is closed under $I$-indexed colimits and $I$-indexed colimits can be computed pointwise. 
        \item  $\Fun^R(\mathcal D, \mathcal  C)\subseteq \Fun(\mathcal D, \mathcal C)$ is closed under $I$-indexed limits and they can be computed pointwise. 
    \end{romanenum}
\end{lemma}
\begin{proof}It suffices to prove the statement on left adjoint functors, then the statement on right adjoint functors follows from replacing $\mathcal C$ and $\mathcal D$ by their opposites. 
    Denote by \[(-)^R\colon \Fun^L(\mathcal C, \mathcal D)\cong \Fun^R(\mathcal D, \mathcal C)^{\operatorname{op}}\] the equivalence provided by \cite[Proposition 5.2.6.2]{highertopostheory} and suppose $F\colon I\to \Fun^{L}(\mathcal C, \mathcal D)$ is some diagram.  
    As $\mathcal D$ has all $I$-indexed colimits, $\Fun(\mathcal C, \mathcal D)$ has all $I$-indexed colimits and they can be computed pointwise by \cite[Corollary 5.1.2.3]{highertopostheory}. 
    Analogously, $\Fun(\mathcal D, \mathcal C)$ has all $I^{\operatorname{op}}$-indexed limits and they can be computed pointwise. 
    This implies that for $d\in\mathcal D$, $c\in\mathcal C$, 
    \begin{align*}\Map_{\mathcal D}((\colim{i\in I}F_i)(c),d)& \cong \Map_{\mathcal D}(\colim{i\in I}F_i(c),d) \cong \clim{i\in I^{\operatorname{op}}}\Map_{\mathcal D}(F_i(c),d) \cong \clim{i\in I}\Map_{\mathcal C}(c,F_i^{R}(d)) \\&\cong \Map_{\mathcal C}(c, \clim{i\in I^{\operatorname{op}}}(F_i^{R}(d)))\cong \Map_{\mathcal C}(c,(\clim{i\in I^{\operatorname{op}}}F_i^{R})(d)).\end{align*} 
    This equivalence is clearly natural in $c,d$ and exhibits $\colim{i\in I}F_i$ as left adjoint to $\clim{i\in I^{\operatorname{op}}}F_i^{R}$. 
\end{proof}
We repeatedly used the following result to \textit{glue} left adjoints. 
\begin{lemma}\label{adjunctionsbigpresentable}Suppose that $\Lambda$ is a possibly large category, $l\colon \mathcal C\to\mathcal D\in \Fun(\Lambda, \vlCat)$ and for all $\lambda\in\Lambda$, $l({\lambda})$ admits a right adjoint such that the mate of the commutative diagram \begin{center}
        \begin{tikzcd}
            \mathcal{C}_{\lambda}\arrow[r,"C(\lambda\to\mu)"]\arrow[d,"l({\lambda})"] & \mathcal{C}_{\mu}\arrow[d,"l({\mu})"]\\ 
            \mt{D_{\lambda}}\arrow[r,"D(\lambda\to \mu)"]&\mt{D}_{\mu}
        \end{tikzcd}
    \end{center} commutes. 
\begin{romanenum}
    \item There exist $r\colon \mt{C}\to \mt{D}\in \Fun(\Lambda, \vlCat)$ with $r(\lambda)=r_{\lambda}$ for all $\lambda\in\Lambda$ and $\epsilon\colon l\circ r\to \id_{\mt{C}}$ such that for all $\lambda\in\Lambda$, $\epsilon(\lambda)$ exhibits $r(\lambda)$ as right adjoint to $l({\lambda})$.

\item By taking the colimit over $\Lambda$, we obtain functors \[l_{\infty}\colon \colim{\Lambda}\mt{D}_*\to \colim{\Lambda}C_*\text{ and }r_{\infty}\colon\colim{\lambda} C_*\to\colim{\Lambda}\mt{D}_*\] and a natural transformation $\epsilon_{\infty}\colon l_{\infty}r_{\infty}\to \id$ which exhibits $l_{\infty}$ as left adjoint to $r_{\infty}$.
\end{romanenum}
\end{lemma}
\begin{proof}
    This follows from a characterization of parametrized adjoint functors \cite[Corollary 3.2.3, 3.2.8/Corollary 3.2.3]{martini2024colimitscocompletionsinternalhigher}/ the discussion of relative adjunctions in \cite[section 7.3]{higheralgebra}.
    Work in the large universe $\mathcal U_1$ so that $\Lambda$ is small and 
    identify \[ \Fun(\Lambda, \vlCat)\cong \Fun^{\mathcal U_1-\operatorname{lim}}(\Psh(\Lambda)^{\operatorname{op}}, \vlCat), \] where presheaves are taken in the universe $\mathcal U_1$ and $\Fun^{\mathcal U_1-\operatorname{lim}}$ refers to the category of $\mathcal U_1$-small limits-preserving functors.
    By \cite[Corollary 3.2.3, 3.2.8]{martini2024colimitscocompletionsinternalhigher}, there exist \[r\colon \mt{D}\to \mt{C}\in  \Fun^{\mathcal U_1-\operatorname{lim}}(\Psh(\Lambda^{\operatorname{op}})^{\operatorname{op}}, \vlCat)\] with $r(\lambda)=r_{\lambda}$ for all $\lambda\in\Lambda$ and natural transformations \[\epsilon\colon l\circ r\to \id, \eta\colon \id\to r\circ l\] satisfying the triangle identities. 
    It follows that \[\colim{\lambda\in\Lambda}\epsilon(\lambda)\colon \colim{\lambda\in\Lambda}l(\lambda)\circ\colim{\lambda\in\Lambda}r(\lambda)\to \id_{\colim{\lambda\in\Lambda}\mathcal D(\lambda)}\] and \[\colim{\lambda\in\Lambda}\eta(\lambda)\colon \id_{\colim{\lambda\in\Lambda}\mathcal C(\lambda)}\to \colim{\lambda\in\Lambda}r(\lambda)\colim{\lambda\in\Lambda}l(\lambda)\] satisfy the triangle identities and hence exhibit $\colim{\lambda\in\Lambda}l(\lambda)$ as left adjoint to $\colim{\lambda\in\Lambda}r(\lambda)$. 
\end{proof}


\begin{lemma}\label{filteredcolimitsofcategoriesappendix}
Suppose that $\Lambda$ is a possibly large filtered category and $\mt{B}_*\colon \Lambda\to \vlCat$ is a diagram of large categories such that for all $\lambda\to\mu\in \Lambda$, $\mt{B}_{\lambda}\to\mt{B}_{\mu}$ is fully faithful.
Denote by $\mt{B}_{\infty}$ its colimit in the very large category of large categories $\vlCat$. 

\begin{romanenum}
\item For all $\lambda\in\Lambda$, $\mt{B}_{\lambda}\to \mt{B}_{\infty}$ is fully faithful. 
\item If for all $\lambda\to \mu\in\Lambda$, $\mt{B_{\lambda}}\to\mt{B}_{\mu}$ is a left adjoint, then $\mt{B_{\lambda}}\to {\mt{B}_{\infty}}$ is a left adjoint and in particular preserves all colimits. 
\item If $I$ is a category such that for all $\lambda\to\mu\in\Lambda$, $\mt{B_{\lambda}}\to\mt{B}_{\mu}$ preserves $I$-indexed (co)limits, then $\mt{B_{\lambda}}\to{\mt{B}_{\infty}}$ preserves $I$-indexed (co)limits, respectively. 
\item Suppose $I$ is a category such that for all $F\colon I\to \mt{B}_{\infty}$, there exists $\lambda\in\Lambda$ such that $F$ factors over $I\to \mt{B_{\lambda}}\hookrightarrow \mt{B}_{\infty}$. If $\mt{B_{\lambda}}$ has $I$-indexed (co)limits for all $\lambda\in\Lambda$ and for all $\lambda\to\kappa\in\Lambda$, $\mt{B_{\lambda}}\to \mt{B_{\kappa}}$ preserves $I$-indexed (co)limits, then ${\mt{B}_{\infty}}$ has all $I$-indexed (co)limits.
\end{romanenum} 
\end{lemma}
\begin{proof}\begin{romanenum}
    \item For $\kappa\to\lambda\in\Lambda$ denote by \[ \mt{B}_{\kappa}\xrightarrow{i^{\lambda}_{\kappa}}\mt{B}_{\lambda}\xrightarrow{i_{\lambda}}\mt{B}_{\infty}\] the canonical maps. 
    For $\lambda\in\Lambda$ and $b,c\in \mt{B}_{\lambda}$, \[\Map_{\mt{B_{\infty}}}(i_{\lambda}b,i_{\lambda}c)\cong \colim{\mu\in \phantom{}_{\lambda\backslash}\Lambda}\Map_{B_{\mu}}(i_{\lambda}^{\mu}b,i_{\lambda}^{\mu}c), \] see e.g \cite{Rozenblyumfilteredcolimitsofcategories} in the large universe.
    By assumption, this is a filtered colimit over equivalences, which shows fully faithfulness of $i_{\lambda}$. 

    \item Suppose now that $\lambda$ is such that for all $\lambda\to\mu\in\Lambda$, $i_{\lambda}^{\mu}\colon \mt{B_{\lambda}}\to\mt{B}_{\mu}$ is a left adjoint. We want to construct a left adjoint to $\mt{B_{\lambda}}\to\mt{B_{\infty}}$. 
    Since $\Lambda$ is filtered, $\colim{\kappa\in\Lambda}\mt{B}_{\kappa}\cong \colim{\mu\in \lambda\backslash\Lambda}\mt{B}_{\mu}$. 
    Denote by $c\mt{B}_{\lambda}\colon \phantom{}_{\lambda\backslash}\Lambda\to \vlCat$ the constant functor with value $\mt{B}_{\lambda}$, and denote by \[c\colon c\mt{B}_{\lambda}\to \mt{B}_*\in \Fun(_{\lambda\backslash}\Lambda, \vlCat)\] the canonical natural transformation. 
    For $\lambda\to\mu\in \lambda\backslash \Lambda$, choose a right adjoint \[r_{\mu}\colon\mt{B}_{\mu}=\mt{B}(\mu)\to \mt{B_{\lambda}}=c\mt{B}(\mu)\] of $c(\mu)=i_{\lambda}^{\mu}=\mt{B}(\lambda\to\mu)$.  
    We claim that for $(\lambda\to\kappa)\to(\lambda\to\mu)\in \phantom{}_{\lambda\backslash}\Lambda$, the mate of 
    \begin{center}
        \begin{tikzcd}
                c\mt{B}_{\lambda}(\kappa)=\mt{B}_{\lambda}\arrow[d,"c(\kappa)"]\arrow[r,"\id"]& c\mt{B_{\lambda}}(\mu)=\mt{B}_{\lambda}\arrow[d,"c(\mu)"]\\ 
                \mt{B_{\kappa}}\arrow[r,"\mt{B}(\kappa\to\mu)"]&\mt{B}_{\mu}
            \end{tikzcd}
\end{center} commutes, i.e.\ 
\[ \id_{\mt{B}_{\lambda}}r^{\kappa}_{\lambda}\to r^{\mu}_{\lambda}k(\mu) \id_{B_{\lambda}}r^{\kappa}_{\lambda}\cong r^{\mu}_{\lambda}c(\mu)c(\kappa)r^{\kappa}_{\lambda}\to r^{\mu}_{\lambda}\mt{B}(\lambda\to\mu)\] is an equivalence. The first map is induced by a unit whitnessing $k(\mu)=\mt{B}(\lambda\to \mu)\dashv r^{\mu}_{\lambda}$ and hence an equivalence since $\mt{B}(\lambda\to\mu)$ is fully faithful. 
Since both are right adjoint to $\mt{B}_{\lambda}\to \mt{B}_{\kappa}\to \mt{B}_{\mu}$, $r^{\mu}_{\lambda}\cong r^{\kappa}_{\lambda}\circ r^{\mu}_{\kappa}$. 
As $\mt{B}(\kappa\to\mu)$ is fully faithful, under this equivalence, the right map becomes 
$r^{\kappa}_{\lambda}(c(\kappa)r^{\kappa}_{\lambda}\to \id)$ which is an equivalence by fully faithfulness of $c(\kappa)=\mt{B}(\lambda\to\kappa)$. 
It now follows from \cref{adjunctionsbigpresentable}, that \[i_{\lambda}\colon \mt{B_{\lambda}}=\colim{\lambda\backslash \Lambda}c\mt{B}_{\lambda}\to \colim{\lambda\backslash \Lambda}\mt{B}_{\lambda}\] is left adjoint to $\colim{\mu\in\phantom{}_{\lambda\backslash}\Lambda}r^{\mu}_{\lambda}$. 
    
\item Suppose that $I$ is a diagram such that for all $\lambda\to \mu$ in $\Lambda$, $\mt{B_{\lambda}}\to\mt{B}_{\mu}$ preserves $I$-indexed (co)limits. Suppose $F\colon I\to \mt{B_{\lambda}}$ has a limit $\clim{i\in I}^{\lambda}F$ in $\mt{B_{\lambda}}$. 
The projections $\clim{i\in I}^{\lambda}F\to F_i$ induce a natural transformation \[\Map_{\mt{B}_{\infty}}(-,i_{\lambda}\clim{i\in I}^{\lambda}F)\to \clim{i\in I}\Map_{\mt{B_{\infty}}}(-, i_{\lambda}F_i)\in \Fun(\colim{\Lambda}\mt{B}_*^{\operatorname{op}}, \an).\] 
For all $\lambda\to\mu\in\Lambda$ and all $b\in\mt{B}_{\mu}^{\operatorname{op}}\subseteq \mt{B}_{\infty}^{\operatorname{op}}$, this is an equivalence since $\mt{B_{\lambda}}\subseteq\mt{B}_{\mu}$ preserves $I$-indexed limits and $\mt{B}_{\mu}\subseteq \mt{B}_{\infty}$ is fully faithful. This shows that \[\Map_{\mt{B_{\infty}}}(-,i_{\lambda}\clim{i\in I}^{\lambda}F)\cong \clim{i\in I}\Map_{\mt{B_{\infty}}}(-, i_{\lambda}F_i)\in \Fun(\mt{B}_{\infty}^{\operatorname{op}}, \an), \] i.e.\ $i_{\lambda}(\clim{i\in I}^{\lambda}F)$ is the limit of $i_{\lambda}F$ in $\mt{B_{\infty}}$. 
The dual argument shows the statement on colimits. 
\item The fourth statement is an immediate consequence of the third.
\end{romanenum} 
\end{proof}
\begin{lemma}\label{limitsincat}
 If $\mathcal C_*\colon I\to \vlCat$ is a possibly large diagram of large categories, for $i\in I$ denote by $p_i\colon\clim{i\in I}\mathcal C\to \mathcal C_i$ the canonical map. 
    For $x,y\in \clim{i\in I}\mathcal C_*$, \[\Map_{\clim{i\in I}\mathcal C_*}(x,y)\cong \clim{i\in I}\Map_{\mathcal C_i}(p_i(x),p_i(y)).\] 
\end{lemma}
\begin{proof}
    This holds as $\Map_{\clim{I}\mathcal C_*}(x,y)$ is the pullback of \[ *\xrightarrow{(x,y)}\clim{i\in I}\mathcal C_*\times \clim{i\in I\mathcal C_*}\xleftarrow{(\text{s,t})}\Fun(\Delta^1, \clim{I}\mathcal C_*), \] in $\Cat$ where $s,t$ denote the source and target maps, and 
    $\Fun(\Delta^1, \clim{I}\mathcal C_*)\cong \clim{i\in I}\Fun(\Delta^1, \mathcal C_i)$. 
\end{proof}
We used the following elementary observation in our discussion of accessible (hyper)\-sheaves. 
\begin{lemma}\label{everyobjectcompactinanaccessiblecategory}
Suppose that $\mathcal X$ is an accessible category. 
For every object $x\in\mathcal X$, there exists a regular cardinal $\kappa$ such that $x$ is $\kappa$-compact. 
\end{lemma}
\begin{proof}Choose a regular cardinal $\lambda$ such that $\mathcal X$ is $\lambda$-accessible. 
 By \cite[Proposition 5.4.2.2]{highertopostheory}, every object in $\mathcal X$ is a $\lambda$-filtered colimit of $\lambda$-compact objects. For $x\in\mathcal X$ choose a regular cardinal $\kappa\geq \lambda$ such that $x$ is a $\kappa$-small colimit of $\lambda$-compact objects. 
 As $\kappa$-filtered colimits commute with $\kappa$-small limits in $\an$ (\cite[Proposition 5.3.3.3]{highertopostheory}), this implies that $x$ is $\kappa$-compact.\end{proof}

 \subsubsection{Derived functors} 

 We now record two results about derived functors which we used in our discussion of solid modules. 
 \begin{lemma}\label{existenceunboundedderivedfunctors}
    Suppose that $\mathcal A$ is a Grothendieck abelian category with enough projectives and $\mathcal D$ is a stable, presentable category with an accessible, left and right-complete $t$-structure $(\mathcal D_{\geq 0},\mathcal D_{\leq 0})$ which is compatible with filtered colimits. 
    Denote by \[\Fun^{\operatorname{r-t-ex},\operatorname{p},\operatorname{colim}}(\mathcal D(\mathcal A),\mathcal D)\subseteq \Fun(\mathcal D(A),\mathcal D)\] the full subcategory on small colimits preserving, right $t$-exact functors $F\colon\mathcal D(A)\to \mathcal D$ such that for all projective objects $P\in\mathcal A$, $F(P)\in\mathcal D^{\heart}$. 

    Then \[ \Fun^{\operatorname{r-t-ex},\operatorname{p},\operatorname{colim}}(\mathcal D(A),\mathcal D)\to \Fun^{\operatorname{colim}}(\mathcal A,\mathcal D^{\heart}), F\mapsto H_0\circ F|_{\mathcal A}\] is an equivalence. 
\end{lemma}
\begin{proof}Let $\mathcal D(A)_{>-\infty}\coloneqq \colim{n\to -\infty}\mathcal D(A)_{\geq n}$. By \cite[Theorem 1.3.3.2, Theorem 1.3.5.24]{higheralgebra}, 
     \[ \Fun^{\operatorname{r-t-ex},\operatorname{p}}(\mathcal D(A)_{>-\infty},\mathcal D)\to \Fun^{\operatorname{r-ex}}(\mathcal A,\mathcal D^{\heart}), F\mapsto H_0\circ F|_{\mathcal A}\] is an equivalence, where the left-hand side denotes right $t$-exact functors $F\colon\mathcal D(A)_{>-\infty}\to \mathcal D$ such that for all projective objects $P\in\mathcal A$, $F(P)\in\mathcal D^{\heart}$. 

     By \cite[Proposition 4.4.2.7]{highertopostheory}, a right exact functor $\mathcal A\to \mathcal D^{\heart}$ preserves small colimits if and only if it preserves coproducts, and a right $t$-exact functor $\mathcal D(\mathcal A)_{>-\infty}\to\mathcal D$ preserves small colimits if and only if it preserves small coproducts. As small coproducts in $\mathcal A$ are exact, under the above equivalence, a functor $\mathcal D(\mathcal A)_{>-\infty}\to\mathcal D$ preserves small coproducts if and only if $H_0\circ F|_{\mathcal A}$ does. 
     Whence the above restricts to an equivalence \[ \Fun^{\operatorname{r-t-ex},\operatorname{p},\operatorname{colim}}(\mathcal D(A)_{>-\infty},\mathcal D)\cong \Fun^{\operatorname{r-ex},\operatorname{colim}}(\mathcal A,\mathcal D^{\heart}).\] 
Restriction defines an equivalence \[\Fun^{\operatorname{colim},\operatorname{r-t-ex},\operatorname{p}}(\mathcal D(A)_{>-\infty},\mathcal D)\cong \Fun^{\operatorname{colim},\operatorname{r-t-ex},\operatorname{p}}(\mathcal D(A)_{\geq 0},\mathcal D_{\geq 0}).\]  
As the $t$-structures on $\mathcal D(\mathcal A)$ and $\mathcal D$ are right-complete (\cite[Theorem 1.3.5.21]{higheralgebra}) and $\mathcal D(\mathcal A)_{\geq 0},\mathcal D_{\geq 0}$ are Grothendieck prestable (\cite[Definition C.1.4.2]{SAG}), by \cite[Proposition C.3.1.1 and Remark C.1.2.10]{SAG}, restriction is an equivalence 
    \[ \Fun^{\operatorname{colim},\operatorname{r-t-ex}}(\mathcal D(\mathcal A),\mathcal D)\cong \Fun^{\operatorname{colim}}(\mathcal D(\mathcal A)_{\geq 0},\mathcal D_{\geq 0}).\] 
    This restricts to an equivalence 
    \[ \Fun^{\operatorname{colim},\operatorname{r-t-ex},\operatorname{p}}(\mathcal D(\mathcal A),\mathcal D)\cong \Fun^{\operatorname{colim},\operatorname{p}}(\mathcal D(\mathcal A)_{\geq 0},\mathcal D_{\geq 0})\] between those functors which carry projectives in $\mathcal A$ to $\mathcal D^{\heart}$. 
\end{proof}
\begin{cor}\label{existencetotalderivedfunctorsonaccessiblesheaves1cat}
     Suppose that $\mathcal D$ is a stable, big presentable category which has small colimits, and $(\mathcal D_{\geq 0},\mathcal D_{\leq 0})$ is a left-complete $t$-structure such that the following holds: \begin{itemize}\item There is an exhaustion $\mathcal D_*\colon M\to \Pr^L$ for $\mathcal D$ by presentable categories such that for all $\mu\in M$, $(\mathcal D_{\mu}\cap \mathcal D_{\geq 0},\mathcal D_{\mu}\cap \mathcal D_{\leq 0})$ is an accessible, right-complete $t$-structure on $\mathcal D_{\mu}$ which is compatible with filtered colimits. \end{itemize}

     Suppose that $(\mathcal C,S)$ is a hyperaccessible explicit covering site, $\mathcal C$ is a 1-category and $R\in\Alg(\Shv_{S}(\mathcal C,\Ab))$ is such that $\mathcal A\coloneqq \LMod{\algebra{R}}{\Shv_S(\mathcal C,\Ab)}$ has enough projectives. 

     Then \[ \Fun^{\operatorname{p},\operatorname{colim},\operatorname{r-t-ex}}(\mathcal D(A),\mathcal D)\to \Fun^{\operatorname{colim}}(\mathcal A,\mathcal D^{\heart}), F\mapsto H_0\circ F|_{\mathcal A}\] is an equivalence. 
\end{cor}
\begin{proof}
    By \cite[Theorem 1.3.3.2, Theorem 1.3.3.7]{higheralgebra}, 
     \[ \Fun^{\operatorname{r-t-ex},\operatorname{p}}(\mathcal D(A)_{>-\infty},\mathcal D)\to \Fun^{\operatorname{r-ex}}(\mathcal A,\mathcal D^{\heart}), F\mapsto H_0\circ F|_{\mathcal A}.\] 

     As coproducts in $\mathcal A$ are exact (this follows from \cref{forgetfreeadjunctionmodules,limitscolimitscanbecomputedfinitestagesheaves}), the same argument as in the proof of \cref{existenceunboundedderivedfunctors} shows that this restricts to an equivalence \[\Fun^{\operatorname{colim},\operatorname{p},\operatorname{r-t-ex}}(\mathcal D(\mathcal A)_{>-\infty},\mathcal D)\cong \Fun^{\operatorname{colim}}(\mathcal A,\mathcal D^{\heart}).\] 
     It remains to show that restriction defines an equivalence \[\Fun^{\operatorname{colim},\operatorname{r-t-ex}}(\mathcal D(\mathcal A),\mathcal D)\cong \Fun^{\operatorname{colim}}(\mathcal D(\mathcal A)_{\geq 0},\mathcal D_{\geq 0}).\]
     Denote by $\Lambda_{\mathcal C}$ the large poset of regular cardinals $\lambda$ such that $(\mathcal C,S)$ is $\lambda$-hyperaccessible and choose $r\in\Lambda_{\mathcal C}$ with \[\algebra{R}\in\Alg(\Shv_S(\mathcal C_{r},\Ab))\subseteq \Alg(\Shv^{\operatorname{acc}}_S(\mathcal C,\Ab).\] 
     Recall from \cref{derivedcategorymodules} that there is a $t$-exact equivalence \[\mathcal D(\mathcal A)\cong \LMod{\algebra{R}}{\hypershvacc_S(\mathcal C,\Sp)}\cong \colim{\lambda\in\Lambda_{\mathcal C,\geq r}}{\LMod{\algebra{R}}{\hypershv_S(\mathcal C_{\lambda},\Sp)}}.\] 
     By construction (see the proof of \cref{derivedcategorymodules}) this is a filtered colimit of $t$-exact equivalences 
     \[ \mathcal D(\mathcal A_{\lambda})\cong \LMod{\algebra{R}}{\hypershv_S(\mathcal C_{\lambda},\Sp)}, \lambda\in\Lambda_{\mathcal C,\geq r},\] where $\mathcal A_{\lambda}\coloneqq \LMod{\algebra{R}}{\Shv_S(\mathcal C_{\lambda},\Ab)}$. 

     \cref{forgetfreeadjunctionmodules,limitscolimitscanbecomputedfinitestagesheaves,closednessmonoidalstructureaccessiblesheaves} 
     imply that for all $\lambda\in \Lambda_{\mathcal C,\geq r}$, \[\mathcal D(\mathcal A_{\lambda})\cong \LMod{\algebra{R}}{\hypershv_S(\mathcal C_{\lambda},\Sp)}\to \LMod{\algebra{R}}{\hypershvacc_S(\mathcal C,\Sp)}\cong \mathcal D(\mathcal A)\] preserves small colimits and that small colimits in $\mathcal D(\mathcal A)\cong \LMod{\algebra{R}}{\hypershvacc_S(\mathcal C,\Sp)}$ can always be computed on some stage \[\mathcal D(\mathcal A_{\lambda})\cong\LMod{\algebra{R}}{\hypershv_S(\mathcal C_{\lambda},\Sp)}\subseteq  \LMod{\algebra{R}}{\hypershvacc_S(\mathcal C,\Sp)}\cong\mathcal D(\mathcal A).\] 
     It now follows from \cref{tstructurespectrumobjects} that 
     \[ \Fun^{\operatorname{colim},\operatorname{r-t-ex}}(\mathcal D(\mathcal A),-)\cong \clim{\lambda\in\Lambda}\Fun^{\operatorname{colim},\operatorname{r-t-ex}}(\mathcal D(\mathcal A_{\lambda}),-)\] and that \[ \Fun^{\operatorname{colim}}(\mathcal D(\mathcal A)_{\geq 0},-)\cong \clim{\lambda\in\Lambda}\Fun^{\operatorname{colim}}(\mathcal D(\mathcal A_{\lambda})_{\geq 0},-).\]

     We are therefore reduced to showing that for all $\lambda\in\Lambda_{\mathcal C,\geq r}$, restriction 
     \[\Fun^{\operatorname{colim},\operatorname{r-t-ex}}(\mathcal D(\mathcal A_{\lambda}),\mathcal D)\to\Fun^{\operatorname{colim}}(\mathcal D(\mathcal A_{\lambda})_{\geq 0},\mathcal D_{\geq 0})\] is an equivalence.
     Fix $\lambda\in\Lambda$. By \cite[Expos{\'e} 2, 6.7]{SGA4}, $\mathcal A_{\lambda}$ 
     is Grothendieck abelian. Choose an exhaustion $\mathcal D_*\colon  M\to \Pr^L$ by presentable categories such that the induced $t$-structure on $\mathcal D_{\mu}$ is right-complete, accessible and compatible with filtered colimits for all $\mu\in M$. 
     As for all $\mu\to\kappa\in M$, $\mathcal D_{\mu}\to\mathcal D_{\kappa}$ preserves colimits, is fully faithful and $t$-exact, for all $\lambda\in\Lambda$,
\[ c_1\colon \colim{\mu\in M}\Fun^{\operatorname{colim},\operatorname{r-t-ex}}(\mathcal D(\mathcal A_{\lambda}),\mathcal D_{\mu})\to \Fun^{\operatorname{colim},\operatorname{r-t-ex}}(\mathcal D(\mathcal A_{\lambda}),\mathcal D)\] and \[c_2\colon \colim{\mu\in M}\Fun^{\operatorname{colim}}(\mathcal D(\mathcal A_{\lambda})_{\geq 0},\mathcal D_{\mu,\geq 0})\to \Fun^{\operatorname{colim}}(\mathcal D(\mathcal A_{\lambda})_{\geq 0},\mathcal D_{\geq 0})\] 
 are fully faithful by \cref{filteredcolimitsofcategoriesappendix}. By \cref{presentablestablesubcategorycontainedinfinitestage,derivedcategoryisstable} and \cite[Theorem 1.3.5.21]{higheralgebra}, $c_1$ is essentially surjective. 
The $t$-structures on $\mathcal D(\mathcal A_{\lambda})$ and $\mathcal D_{\mu}$ are right-complete (by \cite[Theorem 1.3.5.21]{higheralgebra} and assumption on $\mathcal D_{\mu}$), and $\mathcal D(\mathcal A)_{\geq 0},\mathcal D_{\geq 0}$ are Grothendieck prestable (\cite[Definition C.1.4.2]{SAG}). \cite[Proposition C.3.1.1 and Remark C.1.2.10]{SAG} now implies that \[\colim{\mu\in M}\Fun^{\operatorname{colim},\operatorname{r-t-ex}}(\mathcal D(\mathcal A_{\lambda}),\mathcal D_{\mu})\cong \colim{\mu\in M}\Fun^{\operatorname{colim}}(\mathcal D(\mathcal A_{\lambda})_{\geq 0},\mathcal D_{\mu,\geq 0})\] via restriction. 

It remains to show that $c_2$ is essentially surjective. 
As $\mathcal D(\mathcal A_{\lambda})$ is presentable, there exists an small collection of objects $\mathcal C$ which generates $\mathcal D(\mathcal A_{\lambda})$ under small colimits. 
For a functor $F\colon \mathcal D(\mathcal A_{\lambda})_{\geq 0}\to\mathcal D$ choose $\mu\in M$ such that $\oplus_{c\in \mathcal C}F(c)\in\mathcal D_{\mu,\geq 0}$. As $\mathcal D_{(\mu),\geq 0}\subseteq \mathcal D_{(\mu)}$ is closed under small colimits, this implies that for all $c\in\mathcal C$, \[F(c)=\Cofib(\oplus_{d\in\mathcal C}F(d)\to \oplus_{d\in\mathcal C, c\neq d}F(d)\to \oplus_{d\in \mathcal C}F(d))\in\mathcal D_{\mu}\subseteq \mathcal D.\] 
As $\mathcal D_{\mu}\subseteq \mathcal D$ is closed under small colimits and $F$ is cocontinuous, this implies that $F$ factors over $\mathcal D_{\mu,\geq 0}$, which shows that $c_2$ is essentially surjective. 
\end{proof}

\begin{lemma}\label{fullyfaithfulnessderivedfunctors}
    Suppose that $\mathcal A,\mathcal B$ are Grothendieck abelian categories such that countable products in $\mathcal A$ and $\mathcal B$ are exact, and $l\colon \mathcal A\to\mathcal B$ is an exact, fully faithful functor which preserves small colimits and countable products. 
    Suppose that $\mathcal A$ has enough projectives, and for all projectives $p\in \mathcal A$ and all $a\in\mathcal A$, \[\Ext^
    *_{\mathcal B}(lp,lb)\cong \Ext^0_{\mathcal B}(lp,lb)\] is concentrated in degree $0$. 

    Then the derived functor of $l$, $L\colon\mathcal D(\mathcal A)\to\mathcal D(\mathcal B)$ is fully faithful. 
\end{lemma}
\begin{proof}We proceed analogous to the proof of \cite[Proposition 1.3.3.7]{higheralgebra}. 
    
    Since $\mathcal A$ and $\mathcal B$ are Grothendieck abelian, the $t$-structures on their derived categories  $\mathcal D(\mathcal A),\mathcal D(\mathcal B)$ are right-complete by \cite[Proposition 1.3.5.21]{higheralgebra}. 
    As countable products in $\mathcal A$ and $\mathcal B$ exist and are exact, $\mathcal D(\mathcal A)$ and $\mathcal D(\mathcal B)$ admit countable products and they can be computed on the level of representing chain complexes. In particular, $\mathcal D(\mathcal A)_{\geq 0}\subseteq \mathcal D(\mathcal A)$ and $\mathcal D(\mathcal B)_{\geq 0}\subseteq \mathcal D(\mathcal B)$ are stable under countable products. It now follows from \cite[Proposition 1.2.1.19]{higheralgebra} that the $t$-structures on $\mathcal D(\mathcal A)$ and $\mathcal D(\mathcal B)$ are left-complete.
    
    This implies that \begin{align*}\Map_{\mathcal D(\mathcal A)}(-,-)&\cong \clim{m\to -\infty}\clim{n\to \infty}\Map_{\mathcal D(\mathcal A)}(\tau_{\geq m}-, \tau_{\leq n}-)\\ &\cong \clim{m\to -\infty}\clim{n\to \infty}\Map_{\mathcal A}(\tau_{\leq n}\tau_{\geq m}-, \tau_{\geq m}\tau_{\leq n}-)\end{align*} and 
    \begin{align*}\Map_{\mathcal D(\mathcal B)}(-,-)&\cong \clim{m\to -\infty}\clim{n\to \infty}\Map_{\mathcal D(\mathcal B)}(\tau_{\geq m}-, \tau_{\leq n}-)\\ &\cong \clim{m\to -\infty}\clim{n\to \infty}\Map_{\mathcal B}(\tau_{\leq n}\tau_{\geq m}-, \tau_{\geq m}\tau_{\leq n}-).\end{align*} 

    As coproducts and countable products in $\mathcal A$ and $\mathcal B$ are exact, and $l$ preserves small colimits, finite limits and countable products, $L$ preserves coproducts and countable products (which can be computed degreewise on representing chain complexes). $L$ also preserves mapping cones and hence arbitrary small colimits by \cite[Proposition 4.4.2.7]{highertopostheory}. 
    It follows from \cref{derivedcategoryisstable} that $L$ also preserves finite limits and hence countable limits by \cite[Proposition 4.4.2.7]{highertopostheory}. 
    Since $L$ is $t$-exact, it follows that under the above identifications, 
    \[ L_*\colon \Map_{\mathcal D(\mathcal A)}(-,-)\to \Map_{\mathcal D(\mathcal B)}(L-,L-)\] is the limit of the natural transformations 
    \begin{align*} \Map_{\mathcal D(\mathcal A)}(\tau_{\geq m}\tau_{\leq n}-,\tau_{\geq m}\tau_{\leq n}-)\to &\Map_{\mathcal D(\mathcal A)}(L\tau_{\geq m}\tau_{\leq n}-,L\tau_{\geq m}\tau_{\leq n}-)\\ &\cong \Map_{\mathcal D(\mathcal A)}(\tau_{\geq m}\tau_{\leq n}L-,\tau_{\geq m}\tau_{\leq n}L-)\end{align*} induced by $L_*$. 

    It therefore suffices to show that $L_*$ restricts to a fully faithful functor $\mathcal D^b(\mathcal A)\to\mathcal D(\mathcal B)$. 
    Denote by $R\subseteq \mathcal D^b(\mathcal A)$ the full subcategory on objects $r$ such that for all $a\in\mathcal D^b(\mathcal A)$, \[ L_*\colon \Map_{\mathcal D(\mathcal A)}(a,r)\to \Map_{\mathcal D(\mathcal B)}(La,Lr)\] is an equivalence. 
    As $R$ is closed under finite limits, shifts and extensions, it suffices to show that $\mathcal A\subseteq R$, then it follows that $R=\mathcal D^b(\mathcal A)$.  
    Fix $a\in\mathcal D^b(A)$ and represent it by a chain complex of projectives $p_*$ which is concentrated in degrees $\geq -n$. 
    Denote by $p_*^{\leq 0}$ the stupid truncation of $p_*$.
    As $\operatorname{Cofib}(p_*^{\leq 0}\to p_*)\in \mathcal D^b(\mathcal A)_{\geq 1}$, for $r\in\mathcal A$, \[ \Map_{\mathcal D
    (\mathcal A)}(p_*,r)\cong \Map_{\mathcal D(\mathcal A)}(p_*^{\leq 0},r)\] and \[ \Map_{\mathcal D(\mathcal B)}(Lp_*,Lr)\cong \Map_{\mathcal D(\mathcal B)}(L(p_*^{\leq 0}),Lr),\] where we used again that $L$ is $t$-exact. 

    We are therefore reduced to showing that for a bounded complex of projectives $p_*$ and $r\in\mathcal A$, \[ L_*\colon \Map_{\mathcal D(\mathcal A)}(p_*,r)\to \Map_{\mathcal D(\mathcal B)}(Lp_*,Lr)\] is an equivalence.
    We show the stronger statement that for a bounded complex of projectives $p_*$ and $r\in\mathcal A$, \[ L_*\colon \Ext^k_{\mathcal A}(p_*,r)\to \Ext^k_{\mathcal B}(Lp_*,Lr)\] is an equivalence for all $k\in \mathbb Z$. 

    Fix $r\in\mathcal A$ and denote by $G_r\subseteq \mathcal D(\mathcal A)$ the category on objects $a$ such that \[ L_*\colon \Ext^k_{\mathcal A}(a,r)\to \Ext^k_{\mathcal B}(La,Lr)\] is an equivalence for all $k\in \mathbb Z$. 
    As $L$ is cocontinuous, $G_r$ is closed under colimits, extensions and shifts. It therefore suffices to show that for all projective objects $p\in\mathcal A$, $p\in G_r$, then it follows that $G_r$ contains all bounded complexes of projectives.
    Fix a projective object $p\in\mathcal A$. 
    As \[L(p)=l(p),L(r)=l(r)\in\mathcal D(\mathcal B)^{\heart}\] and $p,r\in\mathcal D(\mathcal A)^{\heart}$, 
    for $k<0$, \[\Ext^k_{\mathcal A}(p,r)=0=\Ext^k_{\mathcal B}(Lp,Lr).\]  
    For $k>0$, $\Ext^k_{\mathcal A}(p,r)=0$ since $p$ is projective and $\Ext^k_{\mathcal B}(Lp,Lr)=0$ by our assumption on $l$. 
    The map $L_*\colon \Ext^0_{\mathcal A}(p,r)\to \Ext^0_{\mathcal A}(Lp,Lr)$ is an equivalence by fully faithfulness of $l$. 
\end{proof}

\subsection{Operadic adjunctions}
\label{sec:operadsandmodulecategories}
In this section, we recall how a monoidal adjunction/localization induces an adjunction/localization on module categories. Everything is a well-known and straightforward consequence of the results on relative adjunctions from \cite[section 7.3.2]{higheralgebra}. 
We freely use the notation established in \cite[Chapters 2,4]{higheralgebra}. 
We first recall the results for a general operad before specialising to module categories. 
    \begin{definition}
        Suppose $\mathcal O^{\otimes}\to N(\Fin_*)$ is an operad (\cite[Definition 2.1.1.10]{higheralgebra}). 
        
        \begin{romanenum}
        \item We refer to $\mathcal O\coloneqq \mathcal O^{\otimes}\times_{N(\Fin)_*}\{<0>\}$ as the category underlying $\mathcal O^{\otimes}$. 
        \item For an $\mathcal O^{\otimes}$-monoidal (\cite[Definition 2.1.2.13]{higheralgebra}) category $\mathcal C^{\otimes}\to \mathcal O^{\otimes}$, its underlying category is $\mathcal C\times_{N(\Fin_*)}\{[0]\}$. 

        \item If $\mathcal C^{\otimes}\to \mathcal O^{\otimes}, \mathcal D^{\otimes}\to \mathcal O^{\otimes}\in \oc{{\Cat}}{\mathcal O^{\otimes}}$ are $\mathcal O^{\otimes}$-monoidal categories, a lax $\mathcal O^{\otimes}$-monoidal functor is an operad map (\cite[Definition 2.1.2.7]{higheralgebra}) $\mathcal C^{\otimes}\to \mathcal D^{\otimes}$ over $\mathcal O^{\otimes}$. 
        We denote by $\Fun^{\mathcal O^{\otimes}-\operatorname{lax}}(\mathcal C^{\otimes}, \mathcal D^{\otimes})\subseteq \Fun_{{{\Cat}_{/\mathcal O^{\otimes}}}}(\mathcal C^{\otimes}, \mathcal D^{\otimes})$ the full subcategory on lax $\mathcal O^{\otimes}$-monoidal functors.
        \end{romanenum}
    \end{definition}

    \begin{lemma}\label{equivalencesfiberwise}
        Suppose $\mathcal C^{\otimes}\to \mathcal O^{\otimes}, \mathcal D^{\otimes}\to \mathcal O^{\otimes}\in \oc{{\Cat}}{\mathcal O^{\otimes}}$ are $\mathcal O^{\otimes}$-monoidal categories and $F,G\colon \mathcal C^{\otimes}\to \mathcal D^{\otimes}\in \Fun^{\mathcal O^{\otimes}-\operatorname{lax}}(\mathcal C, \mathcal D)$ are lax $\mathcal O^{\otimes}$-monoidal functors.     
        Denote by \[i\colon \mathcal O=\mathcal O^{\otimes}\times_{N(\Fin_*)}\{<1>\}\to\mathcal O^{\otimes}\] the canonical functor. 
    
        A natural transformation $\eta\colon F\to G\in\Fun_{\mathcal O^{\otimes}}(\mathcal C^{\otimes}, \mathcal D^{\otimes})$ is an equivalence if and only if for all $o\in\mathcal O$, the induced natural transformation $\eta_{io}\colon F_{io}\to G_{io} \in \Fun(\mathcal C^{\otimes}\times_{\mathcal O^{\otimes}}\{io\}, \mathcal D^{\otimes}\times_{\mathcal O^{\otimes}}\{io\})$ is an equivalence.
    \end{lemma}
    \begin{proof}
           Denote by $\mathcal O^{\circ}\subseteq \mathcal O^{\otimes}$ the wide subcategory on inert morphisms (\cite[Definition 2.1.2.3]{higheralgebra}), and let \[\mathcal C^{\circ}\coloneqq \mathcal C^{\otimes}\times_{\mathcal O^{\otimes}}\mathcal O^{\circ}, \, \, \mathcal D^{\circ}\coloneqq \mathcal D^{\otimes}\times_{\mathcal O^{\otimes}}\mathcal O^{\circ}.\] 
           As cocartesian fibrations are stable under pullback (\cite[Proposition 2.4.2.3]{highertopostheory}), $\mathcal C^{\circ}\to\mathcal O^{\circ}$ and $\mathcal D^{\circ}\to\mathcal O^{\circ}$ are cocartesian fibrations.
           A morphism $a\to b\in\mathcal C^{\circ}$ is $(\mathcal C^{\circ}\to\mathcal O^{\circ})$-cocartesian if and only if its image in $\mathcal C^{\otimes}$ is $(\mathcal C^{\otimes}\to \mathcal O^{\otimes})$-cocartesian (\cite[\href{https://kerodon.net/tag/01UF}{Tag 01UF}]{kerodon}). Hence by \cite[Proposition 2.4.1.3]{highertopostheory}, the $(\mathcal C^{\circ}\to\mathcal D^{\circ})$-cocartesian morphisms in $\mathcal C^{\circ}$ are precisely the inert morphisms, and analogously for $\mathcal D^{\otimes}$. 
           Pullback along $\mathcal O^{\circ}\to\mathcal O^{\otimes}$ therefore defines a functor \[(-)^{\circ}\colon \Fun^{\mathcal O^{\otimes}-\operatorname{lax}}(\mathcal C^{\otimes}, \mathcal D^{\otimes})\to \Fun_{\operatorname{CoCart}/\mathcal O^{\circ}}(\mathcal C^{\circ}, \mathcal D^{\circ})\] from the category of lax $\mathcal O^{\otimes}$-monoidal functors $\mathcal C^{\otimes}\to\mathcal D^{\otimes}$ to the category of functors of coartesian fibrations $\mathcal C^{\circ}\to\mathcal D^{\circ}$ over $\mathcal O^{\circ}$. 
        As $\mathcal C^{\circ}\to \mathcal C^{\otimes}$ is essentially surjective, this functor is conservative. 
        By 2-categorical straightening unstraightening (\cite[Proposition A.2.9]{blans2024chainrulegoodwilliecalculus}), the right-hand side is equivalent to the category $\operatorname{Nat}(\St_{\mathcal O^{\circ}}\mathcal C^{\circ}, \St_{\mathcal O^{\circ}}\mathcal D^{\circ})$ of natural transformations 
        $\St_{\mathcal O^{\circ}}\mathcal C^{\circ}\to\St_{\mathcal O^{\circ}}\mathcal D^{\circ}$ of functors $\mathcal O^{\circ}\to \operatorname{Cat}_{(\infty,2)}$.
        By \cite[Definition 2.1.2.13]{higheralgebra}, for $n\in\mathbb N_0$ and \[o=(o_1, \ldots,o_n)\in\mathcal O^{\otimes}_{<n>}\cong \mathcal O^n, \] the inert morphisms $o\to o_i\in\mathcal O^{\circ}$ induce an equivalence \[\St_{\mathcal O^{\circ}}(\mathcal C^{\circ})(o)=\mathcal C^{\circ}\times_{\mathcal O^{\circ}}\{o\}=\mathcal C^{\otimes}\times_{\mathcal O^{\otimes}}\{o\}\cong \prod_{i=1}^n\mathcal C^{\otimes}\times_{\mathcal O^{\otimes}}\{ o_i\}=\prod_{i=1}^n\St_{\mathcal O^{\circ}}(\mathcal C^{\circ})(o_i), \] and analogously for $\mathcal D^{\otimes}$. 
        Since $F^{\circ},G^{\circ}$ are maps of cocartesian fibrations, under these identifications, 
        \[\St_{\mathcal O^{\circ}}F^{\circ}(o)\cong \prod_{i=1}^n \St_{\mathcal O^{\circ}}F^{\circ}(o_i), \, \,  \St_{\mathcal O^{\circ}}G^{\circ}(o)\cong \prod_{i=1}^n \St_{\mathcal O^{\circ}}G^{\circ}(o_i), \] and \[\St_{\mathcal O^{\circ}}\eta^{\circ}(o)=\prod_{i=1}^n \St_{\mathcal O^{\circ}}\eta^{\circ}(o_i).\] 
        In particular, $\eta^{\circ}$ is an equivalence if and only if for all $o\in\mathcal O^{\circ}_{<1>}=\mathcal O^{\otimes}_{<1>}$, \[\St_{\mathcal O^{\circ}}\eta^{\circ}(o)\colon \St_{\mathcal O^{\circ}}F^{\circ}(o)\to \St_{\mathcal O^{\circ}}G^{\circ}(o)\in\Fun(\mathcal C^{\otimes}\times_{\mathcal O^{\otimes}}\{o\}, \mathcal D^{\otimes}\times_{\mathcal O^{\otimes}}\{o\})\] is an equivalence. 
    \end{proof}

\begin{lemma}[{\cite[Corollary 7.3.2.7]{higheralgebra}}]\label{adjunctionfiberwise}
    Suppose that $\mathcal O^{\otimes}$ is an $\infty$-operad and denote by \[i\colon\mathcal O=\mathcal O^{\otimes}\times_{N(\Fin_*)}\{< 1>\}\to\mathcal O^{\otimes}\] the canonical functor. 
    Suppose that $p_{\mathcal C}\colon \mathcal C^{\otimes}\to \mathcal O^{\otimes}, \, p_{\mathcal D}\colon \mathcal D^{\otimes}\to\mathcal O^{\otimes}$ are $\mathcal O^{\otimes}$-monoidal categories and $L\colon\mathcal C^{\otimes}\to\mathcal D^{\otimes}$ is a strong $\mathcal O^{\otimes}$-monoidal functor (\cite[Definition 2.1.3.7]{higheralgebra}) such that for all $o\in\mathcal O$, the induced functor 
    \[L_{io}\colon\mathcal C^{\otimes}\times_{\mathcal O^{\otimes}}\{io\}\to\mathcal D^{\otimes}\times_{\mathcal O^{\otimes}}\{io\}\] admits a right adjoint. 

    \begin{romanenum}
        \item There exists a lax $\mathcal O^{\otimes}$-monoidal functor $R\colon \mathcal D^{\otimes}\to\mathcal C^{\otimes}$ which is right adjoint to $L$ relative to $\mathcal O^{\otimes}$, i.e.\ there exists a unit transformation $\eta\colon 1\to RL$ such that $p_{\mathcal C}\circ \eta\colon p_{\mathcal C}\to p_{\mathcal C}RL$ is an equivalence. 
        \item For all $o\in\mathcal O^{\otimes}$, $\eta_o\colon \id_{\mathcal C^{\otimes}\times_{\mathcal O^{\otimes}}\{o\}}\to R_oL_o$ exhibits $R_o$ as right adjoint to $L_o$. 
    \end{romanenum}
\end{lemma}
\begin{rem}
In the situation of the above lemma, the counit $\epsilon LR\to \id$ is also relative to $\mathcal O^{\otimes}$, i.e. $p_{\mathcal D}(\epsilon)$ is an equivalence. 
Indeed: As $R$ is lax $\mathcal O^{\otimes}$-monoidal, $p_{\mathcal D}=p_{\mathcal C}\circ R$. 
By the triangle identities, for $d\in\mathcal D^{\otimes}$, $R(\epsilon_d)\circ \eta_{Rd}=id_{Rd}$, whence \[id_{p_{\mathcal C}(Rd)}=p_{\mathcal D}(\epsilon_{d})\circ p_{\mathcal C}(\eta_{Rd}).\] As $p_{\mathcal C}(\eta_{Rd})$ is an equivalence, this shows that $p_{\mathcal D}(\epsilon_d)$ is an equivalence. 
The triangle identities moreover imply that for all $o\in\mathcal O^{\otimes}$, $\epsilon_{o}$ is a counit for the adjunction $R_{o}\dashv L_{o}$.
\end{rem}
\begin{cor}\label{adjunctionfullyfaithfulonoperads}In the situation of \cref{adjunctionfiberwise}, denote by $R\colon\mathcal D^{\otimes}\to\mathcal C^{\otimes}$ the right adjoint relative to $\mathcal O^{\otimes}$. 
    \begin{romanenum}
\item A relative adjunction unit $\id_{\mathcal C^{\otimes}}\to RL$ is an equivalence if and only if $L_{io}$ is fully faithful for all $o\in\mathcal O$. 
\item A relative adjunction counit $L R\to \id_{\mathcal D}$ is an equivalence if and only if for all $o\in\mathcal O$, $L_{io}$ admits a fully faithful right adjoint. 
    \end{romanenum}
\end{cor} 
\begin{proof}
    By \cref{equivalencesfiberwise}, a relative adjunction unit $\eta$ is an equivalence if and only if for all $o\in\mathcal O$, \[\eta_{io}\colon \id_{\mathcal C_{io}}\to (RL)_{io}\cong R_{io}L_{io}\] is an equivalence. Since for all $o\in\mathcal O$, $\eta_{io}$ exhibits $R_{io}$ as right adjoint to $L_{io}$ (\cref{adjunctionfiberwise}), this holds if and only if $L_{io}$ is fully faithful for all $o\in\mathcal O$. 
    The same argument shows that the adjunction counit is an equivalence if and only if $R_{io}$ is fully faithful for all $o\in\mathcal O$. 
\end{proof}
    By definition of $\mathcal O^{\otimes}$-algebras (\cite[Definition 2.1.3.1]{higheralgebra}), pushforward along a lax $\mathcal O^{\otimes}$-monoidal functor $f\colon \mathcal C^{\otimes}\to\mathcal D^{\otimes}$ defines a map \begin{align*}f_{\Alg}\colon \Alg_{/\mathcal O}(\mathcal C^{\otimes})\to\Alg_{/\mathcal O}(\mathcal D^{\otimes}).\end{align*} 
    In particular, a lax monoidal functor $\mathcal C^{\otimes}\to\mathcal D^{\otimes}$ induces a map 
    $\Alg(\mathcal C^{\otimes})\to\Alg(\mathcal D^{\otimes})$, and a lax symmetric monoidal functor $\mathcal C^{\otimes}\to\mathcal D^{\otimes}$ induces a map $\CAlg(\mathcal C^{\otimes})\to\CAlg(\mathcal D^{\otimes})$. 
\begin{cor}\label{adjuncttioninducesadjunctiononalgebraobjects}
    Suppose that $\mathcal C^{\otimes}\to \mathcal O^{\otimes}, \mathcal D\to\mathcal O^{\otimes}$ are $\mathcal O^{\otimes}$-monoidal categories and $L\colon\mathcal C^{\otimes}\to\mathcal D^{\otimes}$ is a strong $\mathcal O^{\otimes}$-monoidal functor such that for all $o\in\mathcal O$, $L_o$ admits a right adjoint. 
    Denote by $R\colon\mathcal D^{\otimes}\to\mathcal C^{\otimes}$ its right adjoint relative to $\mathcal O^{\otimes}$ provided by \cref{adjunctionfiberwise} and by $\eta\colon 1\to RL$ an adjunction unit relative to $\mathcal O^{\otimes}$.  
    \begin{romanenum}
       \item $\eta$ exhibits $L_{\Alg}\colon\Alg_{/\mathcal O}(\mathcal C^{\otimes})\to \Alg_{/\mathcal O}(\mathcal D^{\otimes})$ as left adjoint to $R_{\Alg}$. 
       \item If for all $o\in\mathcal O^{\otimes}$, $L_{io}$ respectively $R_{io}$ is fully faithful, then so are $L_{\Alg}$ and $R_{\Alg}$, respectively.
    \end{romanenum} 
\end{cor}
\subsubsection{Module categories}\label{modulecategories}
We now specialize the above results to module categories, see \cite[section 4.2]{higheralgebra} for the relevant definitions. 
Suppose that $L\colon\mathcal C^{\otimes}\to\mathcal D^{\otimes}$ is a strong monoidal functor such that the underlying functor of categories $l\colon\mathcal C\to\mathcal D$ admits a right adjoint. By \cref{adjunctionfiberwise}, $L$ admits a lax monoidal right adjoint $R\colon \mathcal D^{\otimes}\to\mathcal C^{\otimes}$.   
The functors $L,R$ pull back to functors \[ L_{\LM}\colon\mathcal C^{\otimes}\times_{\Assoc^{\otimes}}\LM^{\otimes}\leftrightarrows \mathcal D^{\otimes}\times_{\Assoc^{\otimes}}\LM^{\otimes}\colon R_{\LM}\] over $\LM^{\otimes}$. The functor $L_{\LM}$ is strong $\LM^{\otimes}$-monoidal and $R_{\LM}$ is lax $\LM$-monoidal.
 The relative adjunction unit provided by \cref{adjunctionfiberwise} induces a natural transformation \[\eta_{\LM}\colon 1\to R_{\LM}L_{\LM}\] which exhibits $R_{\LM}$ as right adjoint to $L_{\LM}$ relative to $\LM^{\otimes}$.
Hence by \cref{adjuncttioninducesadjunctiononalgebraobjects}, 
pushforward along $L_{\LM},R_{\LM}$ defines an adjoint pair 
\[ L_{\LM}\colon \LMod{}{\mathcal C}\rightleftarrows \LMod{}{\mathcal D}\colon R_{\LM}\] and pushforward along $\eta_{\LM}$ is a unit for this adjunction. 
In particular, if $l$ or $r$ is fully faithful then so are $L_{\LM}, R_{\LM}$, respectively by \cref{adjunctionfullyfaithfulonoperads}, \cref{adjuncttioninducesadjunctiononalgebraobjects}. 

\begin{lemma}\label{localisationinduceslocalisationonmodulecategories}
 For $A\in\Alg(\mathcal C)$, the adjunction $L_{\LM}\dashv R_{\LM}$ lifts to an adjoint pair 
    \[ L_{LM}^A\colon \LMod{A}{\mathcal C}\to \LMod{L^{\Alg}A}{\mathcal D}\colon R_{LM}^A.\] 

    The right adjoint $R_{LM}^{A}$ factors as \[ \LMod{L^{\Alg}A}{\mathcal D}\to \LMod{R^{\Alg}L^{\Alg}A}{\mathcal C}\xrightarrow{\eta_{A}^*} \LMod{A}{\mathcal C}, \] where $\eta_{A}^*$ is restriction of scalars along the adjunction unit $\eta_{\Alg}(A)\colon A\to R^{\Alg}L^{\Alg}A\in \Alg(\mathcal C)$. 
\end{lemma}
\begin{proof}
By construction, pushforward along $L,R$, respectively $L_{\LM},R_{\LM}$ yields commutative diagrams 
\begin{center}
    \begin{tikzcd}
    \LMod{}{\mathcal C}\arrow[r,"L_{\LM}"]\arrow[d] & \LMod{}{\mathcal D}\arrow[d] & \LMod{}{\mathcal C}\arrow[d]& \arrow[l,"R_{\LM}"']\LMod{}{\mathcal D}\arrow[d]\\ 
    \Alg(\mathcal C)\arrow[r,"L^{\Alg}"] & \Alg(\mathcal D) & \Alg(\mathcal C)& \arrow[l,"R^{\Alg}"'] \Alg(\mathcal D)
    \end{tikzcd}
\end{center}
where the bottom horizontal maps are the adjoint pair provided by \cref{adjuncttioninducesadjunctiononalgebraobjects} for the operad $\Assoc^{\otimes}$. Denote by $\eta_{\Alg}\colon \id\to R^{\Alg}L^{\Alg}$ a relative adjunction unit whitnessing $L^{\Alg}\dashv R^{\Alg}$.  
The above diagrams imply that for $A\in\Alg(\mathcal C)$, $L_{\LM},R_{\LM}$ pull back to functors 
\[ L^A_{LM}\colon \LMod{A}{\mathcal C}\to \LMod{L^{\Alg}A}{\mathcal D}\] 
and \[ \tilde{R}^{A}_{LM}\colon \LMod{L^{\Alg}A}{\mathcal D}\to \LMod{R^{\Alg}L^{\Alg}A}{\mathcal D}.\]
Denote by $R^{A}_{\LM}$ the composition of $\tilde{R}^A_{\LM}$ with restriction
of scalars \[ \eta_{\Alg}^*\colon \LMod{R^{\Alg}L^{\Alg}A}{\mathcal D}\to \LMod{A}{\mathcal D}\] along the adjunction unit $\eta_{\Alg}(A)\colon A\to R^{\Alg}L^{\Alg}A \in\Alg(\mathcal C)$. 
By \cite[Corollary 4.2.3.2]{higheralgebra}, $(\LMod{}{\mathcal C}\to\Alg(\mathcal C))$ is a cartesian fibration. For a cartesian lift $M\to N\in \LMod{}{\mathcal C}$ of an algebra homomorphism $\phi \colon A\to B$, $\phi^*(N)=M\in\LMod{A}{\mathcal C}$. 
In particular, if $M\in\LMod{A}{\mathcal C}$ and $c\colon (A,N)\to R_{\LM}L_{\LM}(A,M)$ is a cartesian lift of $A\to R^{\Alg}L^{\Alg}A$, then \[N=R^{A}_{\LM}L^{A}_{\LM}M\in\LMod{A}{\mathcal C}.\] Since $c\colon N\to R_{\LM}L_{\LM}M$ is cartesian, there exists a unique filler $\eta^A_{M}\colon (A,N)\to R^{A}_{\LM}L^{A}_{\LM}M$ for \begin{center}\begin{tikzcd}
    & (A,M)\arrow[dl, dashed, "\eta^A_M"',"\exists!"]\arrow[d,"\eta_{\LM}"] \\ 
     (A,N)\arrow[r,"c"] & R_{\LM}L_{\LM}(A,M)
\end{tikzcd}
\end{center} lifting  \begin{center}\begin{tikzcd}  & A\arrow[d,"\eta_{\Alg}"] \arrow[dl,"\id"']\\     A \arrow[r,"\eta_{\Alg}"'] & R^{\Alg}L^{\Alg}A.\end{tikzcd}\end{center} As $\eta^A_{M}$ lifts $id_A$, $\eta^A_{M}$ defines a morphism $\eta^{A}(M)\colon M\to R^{A}_{LM}L^{A}_{LM}(M)\in\LMod{A}{\mathcal C}$. We now show that these morphisms exhibit $R_{LM}^A$ as right adjoint to $L_{LM}^A$ in the sense of \cite[Proposition 5.1.10]{landinfinity}.

For $T\in\LMod{L^{\Alg}A}{\mathcal C}$, choose a cartesian lift \[l\colon (A,S)\to R_{\LM}(L^{\Alg}A,T)\] of $\eta_{\Alg}(A)\colon A\to R^{\Alg}L^{\Alg}A$.
Since $\LMod{A}{\mathcal C}\to\LMod{}{\mathcal C}$ maps $R_{\LM}^{A}(T)$ to $(A,S)$, \[\Map_{\LMod{A}{\mathcal C}}(-,R_{\LM}^{A}(T))=\Map_{\LMod{}{\mathcal C}}(-,(A,S))\times_{\Map_{\Alg(\mathcal C)}(A,A)}\{\id_A\}.\] Pushforward along the cartesian morphism $l$ induces an equivalence \begin{align*}\Map_{\LMod{}{\mathcal C}}(-,(A,S))&\times_{\Map_{\Alg(\mathcal C)}(A,A)}\{\id_A\}\\&\cong \Map_{\LMod{}{\mathcal C}}(-,R_{\LM}(L^{\Alg}A,T))\times_{\Map_{\Alg(\mathcal C)}(A,R^{\Alg}L^{\Alg}A)}\{\eta_{\Alg}(A)\}.\end{align*} 
$R_{\LM}$ and the adjunction unit $\eta_{\LM}$ yield an equivalence 
\[  \Map_{\LMod{}{\mathcal C}}(L_{\LM}-,-)\cong \Map_{\LMod{}{\mathcal C}}(-,R_{\LM}-)\] 
lifting the equivalence 
\[ \Map_{\Alg(\mathcal D)}(L^{\Alg}-,-)\cong \Map_{\Alg(\mathcal C)}(-,R^{\Alg}-) \] provided by $R^{\Alg}$ and the adjunction unit $\eta_{\Alg}$. 
This induces an equivalence
\begin{align*} \Map_{\LMod{}{\mathcal C}}&(L_{\LM}-,-)\times_{\Map_{\Alg(\mathcal C)}(L^{\Alg}A,L^{\Alg}A)}\{\id_{L^{\Alg}(A)}\}\\ &\cong \Map_{\LMod{}{\mathcal C}}(-,R_{\LM}-)\times_{\Map_{\Alg(\mathcal C)}(A,R^{\Alg}L^{\Alg}A)}\{ \eta_{\Alg}(A)\}.\end{align*}
The left-hand side is $\Map_{\LMod{L^{\Alg}A}{\mathcal C}}(L_{\LM}^A-,-)$. 
For $M\in\LMod{A}{\mathcal C}$, the equivalence constructed above factors as \begin{align*} \Map_{\LMod{L^{\Alg}A}{\mathcal C}}(L_{\LM}^A(M),-) \xrightarrow{R^{A}_{\LM,*}} &\Map_{\LMod{A}{\mathcal C}}(R_{\LM}^AL_{\LM}^A(M),R_{\LM}^A-)\\ &\xrightarrow{(\eta^A_{M})^*}\Map_{\LMod{A}{\mathcal C}}(M,R_{\LM}^A-),\end{align*} which shows that $\eta_{\LM}^A$ exhibits $R_{\LM}^A$ as right adjoint to $L_{\LM}^A$. 
\end{proof}
\begin{rem}\label{naturalitylocalizationmodulecategories} 
The adjunction $L_{\LM}^A\dashv R^A_{\LM}$ has the following properties: 
\begin{romanenum}
    \item For an algebra map $A\to B\in\Alg(\mathcal C)$, 
    there are commutative diagrams 
    \begin{center}
            \begin{tikzcd}
                \LMod{B}{\mathcal C}\arrow[d,"\phi^*"]\arrow[r,"L_{LM}^B"] & \LMod{L^{\Alg}B}{\mathcal D}\arrow[d,"(L^{\Alg}\phi)^*"] \\ 
                \LMod{A}{\mathcal C}\arrow[r,"L_{LM}^A"]& \LMod{L^{\Alg}A}{\mathcal D}
            \end{tikzcd}
    and 
    \begin{tikzcd}
    \LMod{L^{\Alg}B}{\mathcal D}\arrow[d,"(L^{\Alg}\phi)^*"]\arrow[r,"\overline{R}_{LM}^B"] \arrow[rr,"R^{A}_{\LM}", bend left]& \LMod{R^{\Alg}L^{\Alg}B}{\mathcal D}\arrow[d,"R^{\Alg}L^{\Alg}\phi^*"] \arrow[r,"\eta_{\LM}(B)^*"]& \LMod{B}{\mathcal D}\arrow[d,"\phi^*"]\\ 
    \LMod{L^{\Alg}A}{\mathcal D}\arrow[rr,"R^{A}_{\LM}", bend right]\arrow[r,"\overline{R}_{LM}^A"] & \LMod{R^{\Alg}L^{\Alg}A}{\mathcal C} \arrow[r,"\eta_{\LM}(A)^*"] & \LMod{A}{\mathcal C}
    \end{tikzcd}.
    \end{center}

    \item If the functor $l\colon \mathcal C\to\mathcal D$ underlying $L$ is fully faithful, then so are 
    \begin{align*}L^{\Alg}\colon\Alg(\mathcal C^{\otimes})\to\Alg(\mathcal D^{\otimes}), \\ L_{\LM}\colon \LMod{}{\mathcal C}\to \LMod{}{\mathcal D} \end{align*} by \cref{adjuncttioninducesadjunctiononalgebraobjects}. This implies that  
    \[ L_{\LM}^A\colon \LMod{A}{\mathcal C}\to \LMod{L^{\Alg}A}{\mathcal D}\] is fully faithful for all $A\in \Alg(\mathcal C).$
    \item If $l$ admits a fully faithful right adjoint $r$, then the counit $1\to LR$ is an equivalence, i.e.\ $R$ is fully faithful by \cref{adjunctionfullyfaithfulonoperads}. 
    \item Suppose that $R$ is fully faithful. Then \[R^{\Alg}\colon\Alg(\mathcal D^{\otimes})\to\Alg(\mathcal C^{\otimes})\] and \[R_{\LM}\colon \LMod{}{\mathcal D}\to \LMod{}{\mathcal C}\] are fully faithful by \cref{adjuncttioninducesadjunctiononalgebraobjects}. This implies that \[ \overline{R}_{\LM}^A\colon \LMod{A}{\mathcal D}\to \LMod{R^{\Alg}A}{\mathcal C}\] is fully faithful for all $A\in\Alg(\mathcal D).$ 
    As for $A\in\Alg(\mathcal C)$, $L^{\Alg}R^{\Alg}A\cong A$ (by fully faithfulness of $R^{\Alg}$),    
    \[ R_{\LM}^{R^{\Alg}A}\colon \LMod{A}{\mathcal D}\cong \LMod{L^{\Alg}R^{\Alg}A}{\mathcal D}\to \LMod{R^{\Alg}A}{\mathcal C}\] is fully faithful as well. This is right adjoint to $L^{R^{\Alg}A}_{\LM}$. 
\end{romanenum}
\end{rem}\label{modulesarelefttensored}
   Suppose that $\mathcal C$ is a symmetric monoidal category and $A\in \Alg(\mathcal C)$. By \cite[Remark 4.3.3.7, Example 4.3.1.15]{higheralgebra}, for $A\in \Alg(\mathcal C)$, $\LMod{A}{\mathcal C}$ is right-tensored over $\mathcal C$. 
    As $\mathcal C^{\otimes}\cong \mathcal C^{\otimes}_{\operatorname{rev}}$, we can also consider $\LMod{A}{\mathcal C}$ as left-tensored over $\mathcal C$ by \cite[Proposition 4.6.3.15, Corollary 4.3.2.8]{higheralgebra}. 
\begin{lemma}\label{modulesinmodulecategories}
    Suppose that $\mathcal C$ is a symmetric monoidal category with geometric realizations ($\Delta^{\operatorname{op}}$-indexed colimits) whose tensor product preserves geometric realizations in both variables. 
    Endow $\Alg(\mathcal C)$ with the symmetric monoidal structure described in \cite[Example 3.2.4.4]{higheralgebra}. 
    For $A,B\in\Alg(\mathcal C)$, 
    \[\LMod{B}{\LMod{A}{\mathcal C}}\cong \LMod{A\otimes B}{\mathcal C}, \] where we consider $\LMod{A}{\mathcal C}$ as left-tensored over $\mathcal C$ as described above. 
\end{lemma}
\begin{proof}
    By \cite[Corollary 4.3.2.8, Proposition 4.6.3.11, Proposition 4.3.2.7]{higheralgebra}, \[\LMod{A\otimes B}{\mathcal C}\cong \phantom{}_{A\otimes B}\BiMod{\mathcal C}_{1}\cong \phantom{}_{A}\BiMod{\mathcal C}_{B^{\operatorname{rev}}}\cong \RMod{B^{\operatorname{rev}}}{\LMod{A}{\mathcal C}}.\]  
    The involution $(-)^{\operatorname{rev}}\colon \mathcal{BM}\to \mathcal{BM}$ \cite[Construction 4.6.3.1]{higheralgebra} restricts to an equivalence $r\colon \LM\cong\RM$. Postcomposition with $r$ yields an equivalence \[\RMod{}{\LMod{A}{\mathcal C}}\cong \LMod{}{\LMod{A}{\mathcal C}}\] which pulls back to an equivalence 
    \[\RMod{B^{\operatorname{rev}}}{\LMod{A}{\mathcal C}}\cong \LMod{B}{\LMod{A}{\mathcal C}}.\qedhere\] 
\end{proof}
\subsection{Tensor products of presentable categories and symmetric monoidal structures} 

By construction of the tensor product on $\Pr^L$, the forget functor $\Pr^L\to \vlCat$ admits a lax symmetric monoidal structure.
For presentable categories $\mathcal C$, $\mathcal D$, the lax monoidal structure evaluates to a functor \[\mu_{\mathcal C, \mathcal D}\colon\mathcal C\times \mathcal D\to \mathcal C\otimes_{\Pr^L}\mathcal D.\] We want to describe these functors explicitly.  
Recall from \cite[Proposition 4.8.1.17]{higheralgebra} and \cite[Proposition 5.2.6.2, Remark 5.5.2.10]{highertopostheory}/\cref{identifysheaveswithtensorproductinprl}, that for presentable categories $\mathcal C, \mathcal D$, \[\mathcal C\otimes_{\Pr^L}\mathcal D\cong \Fun^{\operatorname{lim}}(\mathcal C^{\operatorname{op}}, \mathcal D).\] 
\begin{lemma}\label{idenrtificationtensorproductprllimitpreservingfunctors}
    Suppose that $\mathcal C, \mathcal D$ are presentable categories and $L\colon\mathcal P(\mathcal C_0)\to\mathcal C$ is a localization.
    For $c\in\mathcal C$ and $d\in\mathcal D$, the fully faithful functor 
    \[\mathcal C\otimes_{\Pr^L} \mathcal D\cong \Fun^{\operatorname{lim}}(\mathcal C^{\operatorname{op}}, \mathcal D)\xhookrightarrow{L^*} \Fun^{\operatorname{lim}}(\mathcal P(\mathcal C_0)^{\operatorname{op}}, \mathcal D)\cong \Fun(\mathcal C_0^{\operatorname{op}}, \mathcal D)\] sends $\mu_{\mathcal C, \mathcal D}(c,d)$ to  
\[ \mathcal C_0^{\operatorname{op}}\xrightarrow{\Map_{\mathcal C}(Ly-,c)}\an\xrightarrow{\operatorname{const}_d}\mathcal D, \] where $y\colon \mathcal C_0\to \mathcal P(\mathcal C_0)$ is the Yoneda embedding and $\an\xrightarrow{\operatorname{const}_d}\mathcal D$ denotes the unique cocontinuous functor with $\operatorname{const}_d(*)=d$. 
\end{lemma}
\begin{proof}We will deduce this from the proof of \cite[Proposition 4.8.1.15]{higheralgebra}. 
    Suppose first that $\mathcal C=\mathcal P(\mathcal C_0)$ and $\mathcal D=\mathcal P(\mathcal D_0)$ are categories of presheaves on small categories. 
    By construction of the symmetric monoidal structure on $\Pr^L$, \[\mathcal P(\mathcal C_0)\otimes_{\Pr^L}\mathcal P(\mathcal D_0)\cong \mathcal P(\mathcal C_0\times\mathcal D_0), \] and under this identification, \[\mathcal P(\mathcal C_0)\times\mathcal P(\mathcal D_0)\to \mathcal P(\mathcal C_0\times\mathcal D_0)\] is the left Kan extension of the Yoneda embedding $\mathcal C_0\times\mathcal D_0\hookrightarrow \mathcal P(\mathcal C_0\times\mathcal D_0)$ along the Yoneda embedding \[\mathcal C_0\times\mathcal D_0\to \mathcal P(\mathcal C_0)\times \mathcal P(\mathcal D_0).\] Denote by $y\colon \mathcal C_0\to\mathcal P(\mathcal C_0)$ the Yoneda embedding and by $\operatorname{const}_{*}\colon\an\to\mathcal P(\mathcal D_0)$ the unique cocontinuous functor with $\operatorname{const}_{*}(*)=*$ and let \[\operatorname{const}_{-}\colon \mathcal P(\mathcal D_0)\to \Fun(\an,\mathcal P(\mathcal D_0)), d\mapsto d\times \operatorname{const}_{*}.\] For $d\in\mathcal P(\mathcal D_0)$, $\operatorname{const}_{d}$ is the unique cocontinuous functor $\an\to\mathcal P(\mathcal D_0)$ with $\operatorname{const}_{d}(*)=d$. 
    The functor \begin{align*}\nu\colon \mathcal P(\mathcal C_0)\times\mathcal P(\mathcal D_0)& \to\Fun(\mathcal C_0^{\operatorname{op}}, \mathcal P(\mathcal D_0))\\ (c,d)&\mapsto \operatorname{const}_d\circ \Map_{\mathcal P(\mathcal C_0)}(y-,c)\end{align*} is cocontinuous in both variables. As the composition \[\mathcal C_0\times\mathcal D_0\hookrightarrow \mathcal P(\mathcal C_0)\times\mathcal P(\mathcal D_0) \xrightarrow{\nu}\Fun(\mathcal C_0^{\operatorname{op}}, \mathcal P(\mathcal D_0))\cong \mathcal P(\mathcal C_0\times\mathcal D_0)\] is equivalent to the Yoneda embedding, $\nu=\mu_{\mathcal P(\mathcal C_0),\mathcal P(\mathcal D_0)}$ which proves the statement in case $\mathcal C$ and $\mathcal D$ are presheaf categories. 

    Suppose now that $\mathcal D$ is an arbitrary presentable category and $\mathcal C=\mathcal P(\mathcal C_0)$ is a presheaf category. Choose a localization $L\colon \mathcal P(\mathcal D_0)\rightleftarrows \mathcal D\colon f$.  
    The proof of \cite[Proposition 4.8.1.15]{higheralgebra} shows that 
    \begin{center}
    \begin{tikzcd}
    \mathcal P(\mathcal C_0)\times \mathcal D\arrow[d,"\id\times f"]\arrow[rr,"\mu_{\mathcal P(\mathcal C_0){,}\mathcal D}"] && \mathcal P(\mathcal C_0)\otimes_{\Pr^L} \mathcal D\\ 
    \mathcal P(\mathcal C_0)\times \mathcal P(\mathcal D_0)\arrow[rr,"\mu_{\mathcal P(\mathcal C_0){,}\mathcal P(\mathcal D_0)}"]&& \mathcal P(\mathcal C_0)\otimes_{\Pr^L}\mathcal P(\mathcal D_0)\arrow[u,"\id\otimes L"] 
    \end{tikzcd}
    \end{center} commutes, and the right vertical map exhibits $\mathcal P(\mathcal C_0)\otimes_{\Pr^L} \mathcal D$ as localization of $\mathcal P(\mathcal C_0)\otimes_{\Pr^L} \mathcal P(\mathcal D_0)$ at the maps \[\mu_{\mathcal P(\mathcal C_0), \mathcal P(\mathcal D_0)}(c,d\to fLd), \, c\in\mathcal P(\mathcal C_0), d\in\mathcal P(\mathcal D_0).\] 
    By the above, under the identification $\mathcal P(\mathcal C_0)\otimes_{\Pr^L}\mathcal P(\mathcal D_0)\cong \Fun(\mathcal C_0^{\operatorname{op}}, \mathcal P(\mathcal D_0))$ these are precisely the maps \[\operatorname{const}^{\mathcal P(D_0)}_{d}\circ \Map_{P(\mathcal C_0)}(y-,c)\to \operatorname{const}^{\mathcal P(\mathcal D_0)}_{fLc}\circ \Map_{\mathcal P(\mathcal C_0)}(y-,c), \, d\in\mathcal P(\mathcal D_0), c\in\mathcal P(\mathcal C_0)\] induced by the unit $\id\to fL$. 
    As for $d\in\mathcal P(\mathcal D_0)$, 
    pushforward along $\operatorname{const}_{d}$ is left adjoint to pushforward along $\Map_{\mathcal P(\mathcal D_0)}(d,-)$, $F\in\Fun(\mathcal C_0^{\operatorname{op}}, \mathcal P(\mathcal D_0))$ is local for the maps \[\operatorname{const}^{\mathcal P(\mathcal D_0)}_d\circ \Map_{\mathcal C}(-,c)\to \operatorname{const}^{\mathcal P(\mathcal D_0)}_{fLd}\circ \Map_{\mathcal C}(-,c), \, c\in\mathcal P(\mathcal C_0),d\in\mathcal D\] induced by the unit $\id\to fL$ if and only if $F\in\Fun(\mathcal C^{\operatorname{op}}, \mathcal D)\subseteq \Fun(\mathcal C^{\operatorname{op}}, \mathcal P(\mathcal D_0))$, i.e.\ $F\cong fL\circ F$. 
    We therefore obtain a commutative diagram 
    \begin{center}
    \begin{tikzcd}
    \mathcal P(\mathcal C_0)\times \mathcal D\arrow[d,"\id\times f"]\arrow[rr,"\mu_{\mathcal P(\mathcal C_0){,}\mathcal D}"] && \mathcal P(\mathcal C_0)\otimes_{\Pr^L} \mathcal D\arrow[r,"\cong"] &\Fun(\mathcal C_0^{\operatorname{op}}, \mathcal D) \\ 
    \mathcal P(\mathcal C_0)\times \mathcal P(\mathcal D_0)\arrow[rr,"\mu_{\mathcal P(\mathcal C_0){,}\mathcal P(\mathcal D_0)}"] && \mathcal P(\mathcal C_0)\otimes_{\Pr^L}\mathcal P(\mathcal D_0)\arrow[r,"\cong"]\arrow[u] & \Fun(\mathcal C_0^{\operatorname{op}}, \mathcal P(\mathcal D_0))\arrow[u,"L_*"].
    \end{tikzcd}
\end{center} 
    This implies that for $c\in\mathcal P(\mathcal C_0),d\in\mathcal D$, 
    \begin{align*}\mu_{\mathcal P(\mathcal C_0), \mathcal D}(c,d)=L\circ \mu_{\mathcal P(\mathcal C_0), \mathcal P(\mathcal D_0)}(c,fd)&\cong L\circ \operatorname{const}_{fd}^{\mathcal P(\mathcal D_0)}\circ \Map_{\mathcal P(\mathcal C_0)}(-,c)\\&\cong \operatorname{const}_{d}^{\mathcal D}\circ \Map_{\mathcal P(\mathcal C_0)}(-,c), \end{align*} which shows the statement in case $\mathcal C=\mathcal P(\mathcal C_0)$ a presheaf topos. 

    Suppose now that $\mathcal C$ is a general presentable category and choose a localization \[L\colon\mathcal P(\mathcal C_0)\rightleftarrows \mathcal C\colon f.\]
    Then $L^*\colon \Fun^{\operatorname{lim}}(\mathcal C^{\operatorname{op}}, \mathcal D)\subseteq \Fun^{\operatorname{lim}}(\mathcal P(\mathcal C_0)^{\operatorname{op}}, \mathcal D)$ is an equivalence onto the full subcategory on functors for which the unit $id\to fL$ induces an equivalence $F\cong F\circ fL$. 
    By the above and the proof of \cite[Proposition 4.8.1.15]{higheralgebra}, \[\mathcal C\otimes_{\Pr^L} \mathcal D\subseteq \mathcal P(\mathcal C_0)\otimes_{\Pr^L}\mathcal D\cong \Fun^{\operatorname{lim}}(\mathcal P(\mathcal C_0)^{\operatorname{op}}, \mathcal D)\] is the full subcategory on objects which are local for the morphisms \[\operatorname{const}_d\circ \Map_{\mathcal P(\mathcal C_0)}(-,c)\to \operatorname{const}_d\circ \Map_{\mathcal P(\mathcal C_0)}(-,fLc), \, c\in\mathcal P(\mathcal C_0), d\in\mathcal D\] induced by the unit $c\to fLc$. 
    As pushforward along $\operatorname{const}_d$ is left adjoint to pushforward along $\Map_{\mathcal D}(d,-)$, a functor \[F\colon \mathcal P(\mathcal C_0)^{\operatorname{op}}\to \mathcal D\] is local for these maps if and only if the unit $\id\to fL$ induces an equivalence $F\cong F\circ fL$. We therefore a commutative diagram 
    \begin{center}
        \begin{tikzcd}\mathcal C\otimes_{\Pr^L}\mathcal D \arrow[d,hookrightarrow,"i"]\arrow[r,"\cong"]& \Fun^{\operatorname{lim}}(\mathcal C^{\operatorname{op}}, \mathcal D)\arrow[d,hookrightarrow, "L^*"]\\ 
            \mathcal P(C_0)\otimes_{\Pr^L} \mathcal D\arrow[r,"\cong"]& \Fun^{\operatorname{lim}}(\mathcal P(\mathcal C_0)^{\operatorname{op}}, \mathcal D)\arrow[r,"\cong","y^{*}"'] &\Fun(\mathcal C_0^{\operatorname{op}}, \mathcal D). \end{tikzcd}\end{center}
            Denote by $l\colon  \mathcal P(\mathcal C_0)\otimes_{\Pr^L}\mathcal D\to\mathcal C\otimes_{\Pr^L}\mathcal D$ the left adjoint of $i$. 
            It was shown in the proof of \cite[Proposition 4.8.1.17]{higheralgebra} that $\mu_{\mathcal C, \mathcal D}(c,d)=l(\mu_{\mathcal P(\mathcal C_0), \mathcal D}(c,fd))$.
            We showed above that for $(c,d)\in\mathcal P(\mathcal C_0)\times \mathcal D$, \[\mu_{\mathcal P(\mathcal C_0), \mathcal D}(c,d)=\operatorname{const}_d\circ \Map_{\mathcal P(\mathcal C_0)}(y-,c)\in \Fun(\mathcal C_0^{\operatorname{op}}, \mathcal D).\] 
            In particular, for $(c,d)\in\mathcal C\times\mathcal D\subseteq \mathcal P(\mathcal C_0)\times\mathcal D$, \[\mu_{\mathcal P(\mathcal C_0), \mathcal D}(fc,d)\in \mathcal C\otimes_{\Pr^L}\mathcal D\subseteq \mathcal P(\mathcal C_0)\otimes_{\Pr^L}\mathcal D, \] and hence $i(\mu_{\mathcal C, \mathcal D}(c,d))=\mu_{\mathcal P(\mathcal C_0), \mathcal D}(fc,d)$ and 
            \[(yL)^*\mu_{\mathcal C, \mathcal D}(c,d)=y^*\mu_{\mathcal P(\mathcal C_0), \mathcal D}(fc,d)\cong  \operatorname{const}_d\circ \Map_{\mathcal P(C_0)}(y-,fc)\cong \operatorname{const}_d\circ \Map_{\mathcal C}(Ly-,c).\qedhere\] 
\end{proof}
\subsubsection{Monoidal categories and enrichment}
In this section, we collect examples for categories which are enriched in other monoidal categories. 
\begin{recollection}\label{symmetricmonoidalstructurepseudoenriched}
    By \cite[Corollary 4.2.3.2]{higheralgebra}, a monoidal functor \[\phi\colon \mathcal C^{\otimes}\to\mathcal D^{\otimes}\in \Alg(\Cat)\] induces a functor $\phi^*\colon \LMod{\mathcal D^{\otimes}}{\Cat}\to \LMod{\mathcal C^{\otimes}}{\Cat}$. 

    For a monoidal category $\mathcal D^{\otimes}\to\Assoc^{\otimes}$, the cocartesian fibration $\mathcal D^{\otimes}\times_{\Assoc^{\otimes}}\LM^{\otimes}\to \LM^{\otimes}$ of operads exhibits $\mathcal D$ as left $\mathcal D^{\otimes}$-module. 
    The cocartesian fibration $q\colon \mathcal M\to \LM^{\otimes}$ classifying \[\phi^*(\mathcal D^{\otimes}\times_{\Assoc^{\otimes}}\LM^{\otimes})\in \LMod{\mathcal C^{\otimes}}{\Cat}\] exhibits $\mathcal D$ as left-tensored (\cite[Definition 4.2.1.19]{higheralgebra}) over $\mathcal C$.  
    \end{recollection}
    \begin{ex}\label{internalhomrecoversenrichment}Suppose that $\mathcal C^{\otimes}$ is a closed monoidal category and denote by $\imap_{\mathcal C}(-,-)$ the internal Hom. If $L^{\otimes}\colon \mathcal D^{\otimes}\to\mathcal C^{\otimes}$ is a monoidal left adjoint (i.e.\ a monoidal functor such that the underlying functor $L\colon \mathcal C\to\mathcal D$ admits a right adjoint $R$), then the induced left tensoring exhibits $\mathcal C^{\otimes}$ as enriched (\cite[Definition 4.2.1.28]{higheralgebra}) over $\mathcal D^{\otimes}$ as for $c\in\mathcal C$, 
    \[R\circ \imap_{\mathcal C}(c,-)\vdash c\otimes L(-).\]  
    In particular, $R\circ \imap_{\mathcal C}(c,-)$ is the $\mathcal C$-enriched mapping space functor (\cite[Remark 4.2.1.30]{higheralgebra}).
    
    The equivalences \[(-_1\otimes_{\mathcal C}-)\circ (-_2\otimes_{\mathcal C}-)\cong (-_1\otimes_{\mathcal C}-_2)\otimes_{\mathcal C} -\] provided by the monoidal structure on $\mathcal C$ determine an equivalence 
    \[ \imap_{\mathcal C}(-_2, \imap_{\mathcal C}(-_1,-))\cong \imap_{\mathcal C}(-_1\otimes -_2,-).\] 
    In particular, the adjunction $x\otimes_{\mathcal C} -\dashv \imap_{\mathcal C}(x,-)$ is $\mathcal D$-enriched. \end{ex}

    \begin{ex}\label{enrichmenttensorproductpresentablecats}By \cite[Proposition 3.2.4.3, Example 3.2.4.4]{higheralgebra}, $\CAlg(\Pr^L)$ admits a symmetric monoidal structure such that $\CAlg(\Pr^L)\to \Pr^L$ is symmetric monoidal. In particular, if $\mathcal C^{\otimes}, \mathcal D^{\otimes}\in\CAlg(\Pr^L)$ are two presentably symmetric monoidal categories, $C^{\otimes}\otimes \mathcal D^{\otimes}$ inherits a closed monoidal structure and the unit maps $\an\to\mathcal C^{\otimes}, \mathcal D^{\otimes}$ yield symmetric monoidal left adjoints \[\mathcal C^{\otimes}\to (\mathcal C\otimes_{\Pr^L}\mathcal D)^{\otimes}, \mathcal D^{\otimes}\to (\mathcal C\otimes_{\Pr^L}\mathcal D)^{\otimes}.\]  
This exhibits $\mathcal C^{\otimes}\otimes_{\CAlg(\Pr^L)}\mathcal D^{\otimes}$ as enriched over $\mathcal C^{\otimes}$. 
    \end{ex}
    
Using \cref{idenrtificationtensorproductprllimitpreservingfunctors}, we can recover the internal Hom of $\mathcal C\otimes_{\Pr^L}\mathcal D$ from its $\mathcal D$-enrichment:
\begin{cor}\label{enrichmentrecoversinternalhom}
Suppose that $\mathcal C^{\otimes}, \mathcal D^{\otimes}$ are presentably symmetric monoidal categories. 
Then $\mathcal C^{\otimes}\otimes_{\CAlg(\Pr^L)}\mathcal D^{\otimes}$ is an enhancement of $\mathcal C^{\otimes}\otimes_{\Pr^L}\mathcal D^{\otimes}$ to a presentably symmetric monoidal category (\cite[Example 3.2.4.4]{higheralgebra}). Denote by $\iMap_{\mathcal C\otimes\mathcal D}(-,-)$ its internal Hom, and by \[i_{\mathcal C}\colon\mathcal C\to\mathcal C\otimes_{\Pr^L}\mathcal D\] the functor induced by the unit $\an\to\mathcal D, *\mapsto 1_{\mathcal D}$. 
Denote by $\map_{\mathcal C\otimes\mathcal D}(-,-)$ the $\mathcal D$-enrichment of $\mathcal C\otimes_{\Pr^L}\mathcal D$, by $\iMap_{\mathcal C\otimes\mathcal D}(-,-)$.
Under the equivalence $\mathcal C\otimes_{\Pr^L} \mathcal D\cong \Fun^{\operatorname{lim}}(\mathcal C^{\operatorname{op}}, \mathcal D)$, \[\iMap_{\mathcal C\otimes\mathcal D}(-_1,-_2)\cong \map_{\mathcal C\otimes\mathcal D}(-_1\otimes i_{\mathcal C}(-),-_2).\] 
\end{cor}

\begin{proof}
    Denote by $i_{\mathcal C}\colon\mathcal C\to \mathcal C\otimes\mathcal D, i_{\mathcal D}\colon\mathcal C\to \mathcal D$ the functors induced by the units $\an\to\mathcal C, \an\to\mathcal D$.
    The forget functor $\Pr^L\to \vlCat$ is lax symmetric monoidal and therefore induces a lax symmetric monoidal functor $\CAlg(\Pr^L)\to \CAlg(\vlCat)$ by \cite[Construction 3.2.4.1, Proposition 3.2.4.2]{higheralgebra}. 
    (By \cite[Construction 3.2.4.1]{higheralgebra}, a lax symmetric monoidal functor $\mathcal T\to\mathcal S$ induces a functor $\CAlg(\mathcal T)\to\CAlg(\mathcal S)$ of categories over $N(\Fin_*)$, and by \cite[Proposition 3.2.4.2.2)]{higheralgebra}, this is a lax symmetric monoidal functor.)
    By construction, this fits into a commutative diagram of lax symmetric monoidal functors 
    \begin{center}
        \begin{tikzcd}
        \CAlg(\Pr^L)\arrow[r]\arrow[d] & \CAlg(\vlCat)\arrow[d]\\ 
        \Pr^L\arrow[r] & \vlCat,
        \end{tikzcd}
    \end{center} where the vertical maps are the forget functors $\CAlg(\Pr^L)\to \Pr^L, \CAlg(\vlCat)\to\vlCat$ which are symmetric monoidal (\cite[Example 3.2.4.4]{higheralgebra}). 
    In particular, we obtain a symmetric monoidal functor \[\mu^{\otimes}_{\mathcal C, \mathcal D}\colon \mathcal C\otimes_{\CAlg(\vlCat)}\mathcal D\to \mathcal C\otimes_{\CAlg(\Pr^L)}\mathcal D\] enhancing the functor \[\mu_{\mathcal C, \mathcal D}\colon \mathcal C\times\mathcal D\to\mathcal C\otimes_{\Pr^L}\mathcal D\] from \cref{idenrtificationtensorproductprllimitpreservingfunctors}. 
    Since the functor $\mu^{\otimes}_{\mathcal C, \mathcal D}$ is  symmetric monoidal and natural in $\mathcal C, \mathcal D$ with respect to cocontinuous, symmetric monoidal functors \[\mu_{\mathcal C, \mathcal D}(c,d)\cong \mu_{\mathcal C, \mathcal D}(c,1_{\mathcal D})\otimes_{\mathcal C\otimes_{\Pr^L}\mathcal D} \mu_{\mathcal C, \mathcal D}(1_{\mathcal D},d)\cong i_{\mathcal C}(c)\otimes_{\mathcal C\otimes_{\Pr^L}\mathcal D} i_{\mathcal D}(d).\] 

    Choose a localization $L\colon \mathcal P(\mathcal C_0)\rightleftarrows \mathcal C\colon f$. By \cref{idenrtificationtensorproductprllimitpreservingfunctors}, for $c\in\mathcal C,d\in\mathcal D$, under the embedding \[\mathcal C\otimes_{\Pr^L}\mathcal D\cong \Fun^{\operatorname{lim}}(\mathcal C^{\operatorname{op}}, \mathcal D)\hookrightarrow \Fun^{\operatorname{lim}}(\mathcal C_0^{\operatorname{op}}, \mathcal D), \] $\mu_{\mathcal C, \mathcal D}(c,d)=i_{\mathcal C}(c)\otimes i_{\mathcal D}(d)$ becomes \[\operatorname{const}_d\circ \Map_{\mathcal C}(Ly-,c), \] where $y$ denotes the Yoneda embedding $\mathcal C_0\subseteq \mathcal P(\mathcal C_0)$. As pushforward along $\operatorname{const}_{d}\colon\an\to\mathcal D$ is left adjoint to \[\Map_{\mathcal D}(d,-)_*\colon \Fun(\mathcal C_0^{\operatorname{op}}, \mathcal D)\to \Fun(\mathcal C_0^{\operatorname{op}}, \an), \] it follows that for $c\in\mathcal C$ and $d\in\mathcal D$, 
    \begin{align*}&\Map_{\mathcal C\otimes_{\Pr^L}\mathcal D}\left(i_{\mathcal C}(c)\otimes i_{\mathcal D}(d), \map_{\mathcal C\otimes\mathcal D}(i_{\mathcal C}(-)\otimes x,a)\right)\\ & \cong \Map_{\Fun(\mathcal C_0^{\operatorname{op}}, \mathcal D)}\left(\operatorname{const}_d\circ \Map_{\mathcal C}(Ly-,c), \map_{\mathcal C\otimes\mathcal D}(i_{\mathcal C}(Ly-)\otimes x,a)\right)\\  
    & \cong \Map_{\Fun(\mathcal C_0^{\operatorname{op}}, \an)}\left(\Map_{\mathcal C}(Ly-,c), \Map_{\mathcal D}(d, \map_{\mathcal C\otimes\mathcal D}(i_{\mathcal C}(Ly-)\otimes x,a))\right)\\ 
     & \cong \Map_{\Fun^{\operatorname{lim}}(\mathcal P(\mathcal C_0)^{\operatorname{op}}, \an)}\left(\Map_{\mathcal C}(L-,c), \Map_{\mathcal D}(d, \map_{\mathcal C\otimes\mathcal D}(i_{\mathcal C}(L(-)\otimes x,a))\right)\\ 
     &\cong  \Map_{\Fun^{\operatorname{lim}}(\mathcal P(C_0)^{\operatorname{op}}, \an)}\left(\Map_{\mathcal P(\mathcal C_0)}(-,f(c)), \Map_{\mathcal D}(d, \map_{\mathcal C\otimes \mathcal D}(i_{\mathcal C}(L-)\otimes x,a))\right)\\ 
     &\cong \Map_{\mathcal D}\left(d, \map_{\mathcal C\otimes\mathcal D}i_{\mathcal C}(Lf(c))\otimes x,a\right)\\ 
     &\cong \Map_{\mathcal D}\left(d, \map_{\mathcal C\otimes\mathcal D}(i_{\mathcal C}(c)\otimes x,a)\right),
\end{align*} where we used that \[L^*\colon \Fun^{\operatorname{lim}}(\mathcal C^{\operatorname{op}}, \mathcal D)\hookrightarrow \Fun^{\operatorname{lim}}(\mathcal P(\mathcal C_0)^{\operatorname{op}}, \mathcal D)\cong \Fun(\mathcal C_0^{\operatorname{op}}, \mathcal D)\] is fully faithful, that $\Map_{\mathcal C}(L-,c)\cong \Map_{\mathcal P(\mathcal C_0)}(-,fc)$, fully faithfulness of the Yoneda embedding, and that $f\colon\mathcal C\subseteq \mathcal P(\mathcal C_0)$ is fully faithful. 

By definition of $\map_{\mathcal C\otimes\mathcal D}(-,-)$ and $\iMap_{\mathcal C\otimes\mathcal D}(-,-)$, the right-hand side is isomorphic to \[\Map_{\mathcal C\otimes_{\Pr^L}\mathcal D}(i_{\mathcal D}(d)\otimes i_{\mathcal C}(c)\otimes x,a) \cong \Map_{\mathcal C\otimes_{\Pr^L}\mathcal D}(i_{\mathcal D}(d)\otimes i_{\mathcal C}(c), \iMap_{\mathcal C\otimes\mathcal D}(x,a)).\]
This identification is natural in $c$, $d$, $x$ and $a$. As $\{i_{\mathcal C}(c)\otimes i_{\mathcal D}(d), \, c\in\mathcal C,d\in\mathcal D\}$ generates $\mathcal C\otimes\mathcal D$ under small colimits, this implies that \[\map_{\mathcal C\otimes\mathcal D}(i_{\mathcal C}(-)\otimes x,a)\text{ and } \iMap_{\mathcal C\otimes\mathcal D}(x,a)\] represent the same functor, whence are isomorphic. Their identification is clearly natural in $x$ and $a$. 
\end{proof}
\subsubsection{Monoidal structure on presheaves with values in a presentably symmetric monoidal category}
    We now recall an explicit description of the symmetric monoidal structure on presheaves with values in a presentably symmetric monoidal category which we used in the proof of \cref{closednessmonoidalstructureaccessiblesheaves}. 
        Suppose that $\mathcal C$ is a small symmetric monoidal category and $\mathcal D$ is a presentably symmetric monoidal category. Then $\mathcal P(\mathcal C)\otimes_{\Pr^L}\mathcal D\cong \Fun(\mathcal C^{\operatorname{op}}, \mathcal D)$, cf.\ \cref{identifysheaveswithtensorproductinprl}. Endow $\mathcal P(\mathcal C)$ with the Day convolution symmetric monoidal structure (\cite[Corollary 4.8.1.2]{higheralgebra}). 
        
        By \cite[Remark 2.4.2.7]{higheralgebra}, $\mathcal C$ induces a symmetric monoidal structure on $\mathcal C^{\operatorname{op}}$. By \cite[Proposition 2.2.6.16]{higheralgebra}, $\Fun(\mathcal C^{\operatorname{op}}, \mathcal D^{\otimes})$ inherits a symmetric monoidal structure from the symmetric monoidal structures on $\mathcal C^{\operatorname{op}}$ and $\mathcal D$.
        These two symmetric monoidal structures are equivalent: 
        \begin{lemma}\label{formulasymmetricmonoidalstructureontensorproductofpresentablecatsasdayconcvolution}
The equivalence $\mathcal P(\mathcal C)\otimes_{\Pr^L}\mathcal D\cong \Fun(\mathcal C^{\operatorname{op}}, \mathcal D)$ enhances to a symmetric monoidal equivalence \[\mathcal P(\mathcal C)^{\operatorname{Day}}\otimes_{\CAlg(\Pr^L)}\mathcal D^{\otimes}\cong \Fun(\mathcal C^{\operatorname{op}}, \mathcal D)^{\otimes}.\] 

In particular, for $F,G\in \Fun(\mathcal C^{\operatorname{op}}, \mathcal D)$ and  $c\in\mathcal C$, $(F\otimes G)(c)$ is the left Kan extension of \[ \mathcal C^{\operatorname{op}}\times \mathcal C^{\operatorname{op}}\xrightarrow{F\times G}\mathcal D\times \mathcal D\xrightarrow{-\otimes_{\mathcal D}-}\mathcal D\] along (the opposite of)\[ \mathcal C\times\mathcal C\xrightarrow{-\otimes_{\mathcal C}-}\mathcal C.\] 
    \end{lemma}
    \begin{proof}
        For $\mathcal D=\an$, the equivalence of the two symmetric monoidal structures is remarked in \cite[Remark 4.8.1.13]{higheralgebra}. It can be seen as follows: By \cite[section 3]{Glasman2013DayCF}, the Yoneda embedding enhances to a symmetric monoidal functor $y^{\otimes}\colon \mathcal C^{\otimes}\to \Fun(\mathcal C^{\operatorname{op}},\mathcal D)$. It is straightforward to check that the tensor product on $\Fun(\mathcal C^{\operatorname{op}},\mathcal D)$ preserves colimits in both variables. Whence by the universal property of $\mathcal P(\mathcal C)^{\otimes}$ (\cite[Corollary 4.8.1.12]{higheralgebra}), $y^{\otimes}$ extends uniquely to a symmetric monoidal functor $\mathcal P(\mathcal C)^{\otimes}\to \Fun(\mathcal C^{\operatorname{op}},\mathcal D)^{\otimes}$. The underlying functor $\mathcal P(\mathcal C)\to \Fun(\mathcal C^{\operatorname{op}},\mathcal D)$ is a colimits-preserving extension of the Yoneda embedding and hence an equivalences. 

    Suppose now that $\mathcal D_0$ is a small symmetric monoidal category and $\mathcal D^{\otimes}=\mathcal P(\mathcal D_0)^{\operatorname{Day}}$ is the symmetric monoidal category of presheaves provided by \cite[Corollary 4.8.1.12]{higheralgebra}.
    By construction of the symmetric monoidal structure on $\Pr^L$, \[\mathcal P(\mathcal C)^{\operatorname{Day}}\otimes_{\CAlg(\Pr^L)}\mathcal P(\mathcal D_0)^{\operatorname{Day}}\cong \mathcal P(\mathcal C\times \mathcal D_0)^{\operatorname{Day}}.\] 
    By \cite[Remark 4.8.1.13]{higheralgebra}, the right-hand side is equivalent to the symmetric monoidal category $\Fun(\mathcal C^{\operatorname{op}}\times\mathcal D^{\operatorname{op}}, \an)$ provided by \cite[Proposition 2.2.6.16]{higheralgebra}, where we view $\mathcal C^{\operatorname{op}}\times\mathcal D^{\operatorname{op}}=(\mathcal C\times\mathcal D)^{\operatorname{op}}$ as a symmetric monoidal category via the cartesian monoidal structure on $\Cat$ and \cite[Remark 2.4.2.7]{higheralgebra}.
    By \cite[Remark 2.2.6.8]{higheralgebra} and again \cite[Remark 4.8.1.13]{higheralgebra}, the symmetric monoidal category $\Fun(\mathcal C^{\operatorname{op}}\times\mathcal D^{\operatorname{op}}, \an)^{\otimes}$ is equivalent to $\Fun(\mathcal C^{\operatorname{op}}, \mathcal P(\mathcal D_0)^{\operatorname{Day}})$ with the symmetric monoidal structure provided by \cite[Proposition 2.2.6.16]{higheralgebra}. 
    We now deduce from this the general statement. By \cite[Theorem 2.2]{Nikolaus-Sagave}, for every presentably symmetric monoidal category $\mathcal D^{\otimes}$, there exists a small symmetric monoidal category $\mathcal D_0$ and a symmetric monoidal localization $L\colon \mathcal P(\mathcal D_0)^{\operatorname{Day}}\to \mathcal D^{\otimes}$. Denote by $W\subseteq \Fun(\Delta^1, \mathcal P(\mathcal D_0))$ the $L$-equivalences. 
    It is shown in the proof of \cite[Proposition 4.8.1.15]{highertopostheory} that $\mathcal P(\mathcal C)^{\operatorname{Day}}\otimes_{\CAlg(\Pr^L)}\mathcal D^{\otimes}$ is the symmetric monoidal localization of $\mathcal P(\mathcal C)^{\operatorname{Day}}\otimes_{\CAlg(\Pr^L)}\mathcal P(\mathcal D_0)$ at the morphisms $1\otimes w, w\in W$ and \[\mathcal P(\mathcal C)^{\operatorname{Day}}\otimes_{\CAlg(\Pr^L)}\mathcal P(\mathcal D_0)\cong \mathcal P(\mathcal C\times \mathcal D_0)^{\operatorname{Day}}.\] 
    
    The symmetric monoidal localization $\mathcal P(\mathcal D_0)^{\operatorname{Day}}\to \mathcal D^{\otimes}\in \CAlg(\Pr^L)$ and the equivalence established above yield a symmetric monoidal functor \[ c^{\otimes}\colon \mathcal P(\mathcal C)^{\operatorname{Day}}\otimes_{\CAlg(\Pr^L)}\mathcal P(\mathcal D_0)^{\otimes}\cong \Fun(\mathcal C^{\operatorname{op}}, \mathcal P(\mathcal D_0)^{\operatorname{Day}})^{\otimes}\to \Fun(\mathcal C^{\operatorname{op}}, \mathcal D)^{\otimes}.\] 
    Since the underlying functor $\mathcal P(\mathcal C)^{\operatorname{Day}}\otimes_{\Pr^L}\mathcal P(\mathcal D_0)^{\otimes}\to \Fun(\mathcal C^{\operatorname{op}}, \mathcal D)\cong \mathcal P(\mathcal C)\otimes\mathcal D$ equals $\id_{\mathcal P(\mathcal C)}\otimes L$, $c^{\otimes}$ inverts the morphisms $1\otimes w, w\in W$ and hence factors over a symmetric monoidal functor 
    \[ \mathcal P(\mathcal C)\otimes_{\CAlg(\Pr^L)}\mathcal D\to \Fun(\mathcal C^{\operatorname{op}}, \mathcal D).\] 
    As the forget functor $\CAlg(\Pr^L)\to \Pr^L$ is conservative and symmetric monoidal, this is an equivalence.
    The explicit description of $F\otimes G$ as left Kan extension is \cite[Remark 2.2.6.15]{higheralgebra}. 
    \end{proof}
\subsubsection{Modules in cartesian categories}
In this section, we recall from \cite[section 4.2.2]{higheralgebra} simplicial models for modules, algebras and groups in cartesian monoidal categories. 
\begin{definition}[{\cite[2.4.1.1]{higheralgebra}}]
    Suppose $\mathcal C$ is a category with finite products. 
        Denote by $\mathcal C^{\times}\to N(\Fin_*)$ the cartesian monoidal category associated to $\mathcal C$ and by $\pi\colon\mathcal C^{\times }\to\mathcal C$ the corresponding cartesian structure \cite[Proposition 2.4.1.5]{higheralgebra}. 

    For an operad $\mathcal O^{\otimes}\to N(\Fin_*)$, denote by $\Fun^{\operatorname{lax}}(\mathcal O^{\otimes}, \mathcal C)\subseteq \Fun(\mathcal O^{\otimes}, \mathcal C)$ the full subcategory on functors $F\colon \mathcal O^{\otimes}\to \mathcal C$ such that for all $o\coloneqq (o_1, \ldots,o_n)\in\mathcal O^{\otimes}_{<n>}\cong \mathcal O^n$, 
    the Segal maps $F(o)\to F(o_i)$ exhibit $F(o)\cong \prod_{i=1}^n F(o_i)$. 
    \end{definition}
    
    \begin{proposition}[{\cite[Proposition 2.4.1.7]{higheralgebra}}]\label{algebrasarelaxfunctors} The composition \[\Alg_{/\mathcal O}(\mathcal C^{\times}\times_{N(\Fin_*)}\mathcal O^{\otimes})\subseteq \Fun(\mathcal O^{\otimes}, \mathcal C^{\times}\times_{N(\Fin_*)}\mathcal O^{\otimes})\to \Fun(\mathcal O^{\otimes}, \mathcal C^{\times})\xrightarrow{\pi_*}\Fun(\mathcal O^{\otimes}, \mathcal C)\] induces an equivalence 
    \[ \Alg_{/\mathcal O}(\mathcal C^{\times}\times_{N(\Fin_*)}\mathcal O^{\otimes})\cong \Fun^{\operatorname{lax}}(\mathcal O^{\otimes}, \mathcal C).\] \end{proposition}
    \begin{proof}This is \cite[Proposition 2.4.1.7]{higheralgebra}, recall that a trivial Kan fibration of simplicial sets is a weak equivalence in the Joyal model structure on simplicial sets and hence an equivalence of associated categories. 
    \end{proof}
\begin{definition}[{\cite[Definition 4.1.2.5]{higheralgebra}}]
     Denote by $\Mon(\topo{X})\subseteq \Fun(\Delta^{\operatorname{op}}, \mathcal X)$ the full subcategory of monoid objects, that is simplicial objects $M\colon \Delta^{\operatorname{op}}\to\mathcal X$ such that for all $n\in\mathbb N_0$, the face maps 
    $\{ M([n])\to M(\{i-1,i\})\}_{1\leq i\leq n}$ exhibit $M([n])= \prod_{i=1}^n M(\{i-1,i\})$.
\end{definition}
\begin{proposition}[{\cite[Proposition 4.1.2.10]{higheralgebra}}]\label{characterizationgroupobjects}
   There exists a functor $i\colon \Delta^{\operatorname{op}}\to \Assoc^{\otimes}$ such that for all categories $\topo{X}$ with finite products, $i^*\colon \Fun^{\operatorname{lax}}(\Assoc^{\otimes}, \mathcal X)\to \Fun(\Delta^{\operatorname{op}}, \mathcal X)$ factors over an equivalence $\Fun^{\operatorname{lax}}(\Assoc^{\otimes}, \mathcal X)\cong \Mon(\topo{X})\subseteq \Fun(\Delta^{\operatorname{op}}, \topo{X})$. 
\end{proposition} 
The composite $\Alg(\topo{X}^{\times})=\Alg_{/\Assoc}(\topo{X}^{\times})\cong \Fun^{\operatorname{lax}}(\Assoc^{\times}, \topo{X})\cong \Mon(\topo{X})$ sends an algebra $A\in\Alg(\topo{X}^{\times})$ to the simplicial object  $[n]\mapsto A^n$ with 
\begin{align*}d_i\colon A^n&\to A^{n-1}, \\(a_1, \ldots,a_n)&\mapsto (a_1, \ldots, a_{i}a_{i+1},a_{i+2}, \ldots,a_n)\end{align*} and \begin{align*}s_i\colon A^{n}& \to A^{n+1}, \\ (a_1, \ldots,a_n)& \mapsto (a_1, \ldots, a_{i-1},1,a_i, \ldots,a_n).\end{align*} 
In particular, under the above equivalence, a group object is an monoid object $M\colon\Delta^{\operatorname{op}}\to\mathcal X$ such that $M(\{0\})=*$ and if $[n]=S_1\cup S_2$ with $|S_1\cap S_2|=1$, then $M([n])\to M(S_i), i=1,2$ exhibit $M([n])= M(S_1)\times M(S_2)$. This condition can be checked in the homotopy category of $\topo{X}$, whence \[ \Grp(\topo{X})\cong \Mon(\topo{X})\times_{\Mon(\operatorname{ho}(\topo{X}))}\Grp(\operatorname{ho}(\topo{X})).\]
\begin{lemma}
Suppose $\topo{X}$ is a big topos. 
A monoid $M\in\topo{X}$ is a group object if and only if $\tau_{\leq 0}M\in \Mon(\tau_{\leq 0}\topo{X})$ is a group. 
\end{lemma}
\begin{proof}
Choose an exhaustion $\topo{X}_*\colon\Lambda\to \Pr^L$ of $\topo{X}$ by topoi. 
As for all $\lambda\in\Lambda$, $\topo{X}_{\lambda}\subseteq \topo{X}$ is closed under finite limits and the truncation functor $\tau_{\leq 0}\colon \topo{X}\to\tau_{\leq 0}\topo{X}$ restricts to $\topo{X_{\lambda}}\to\tau_{\leq 0}\topo{X}_{\lambda}$ (\cref{truncationbigtopos}), it suffices to prove the statement for $\topo{X}$ a topos.
As $\Mon(\an)\subseteq \Fun(\Delta^{\operatorname{op}},\an)$ is closed under small limits (which are computed pointwise), the equivalence 
\[\Fun(\Delta^{\operatorname{op}},\topo{X})\cong \Fun(\Delta^{\operatorname{op}},\Fun^{\lim}(\topo{X}^{\operatorname{op}},\an))\cong \Fun^{\lim}(\topo{X}^{\operatorname{op}},\Fun(\Delta^{\operatorname{op}},\an))\] restricts to an equivalence $\Mon(\topo{X})\cong \Fun^{\lim}(\topo{X}^{\operatorname{op}},\Mon(\an))$. Here $\Fun^{\lim}$ denotes the categories of small limits preserving functors.  
The above equivalence restricts further to an equivalence $\Grp(\topo{X})\cong \Fun^{\lim}(\topo{X}^{\operatorname{op}},\Grp(\an))$. 
By \cite[Example 5.2.6.4 and Remark 5.2.6.5]{higheralgebra}, \[\Grp(\an)\cong \Mon(\an)\times_{\Mon(\Set)}\Grp(\Set),\] 
whence \[\Fun^{\lim}(\topo{X}^{\operatorname{op}},\Grp(\an))\cong \Fun^{\lim}(\topo{X}^{\operatorname{op}},\Mon(\an))\times_{\Fun^{\lim}(\topo{X}^{\operatorname{op}},\Mon(\Set))}\Fun^{\lim}(\topo{X}^{\operatorname{op}},\Grp(\Set)).\] 
The right-hand side is equivalent to $\Mon(\topo{X})\times_{\Mon(\tau_{\leq 0}\topo{X})}\Grp(\tau_{\leq 0}\topo{X})$. 
\end{proof}

\begin{definition}{\cite[Definition 4.2.2.2]{higheralgebra}}
For a category $\mathcal X$ denote by \[\LMon{}(\mathcal X)\subseteq \Fun(\Delta^1, \Fun(\Delta^{\operatorname{op}}, \mathcal X))\] the full subcategory of maps of simplicial objects $M\to A\in \Fun(\Delta^{\operatorname{op}}, \mathcal X)$ such that 
\begin{romanenum}\item $A\in\Mon(\mathcal X)$.
    \item For all $n\in\mathbb N_0$, the Segal map $M([n])\to M(\{n\})\cong M(\{0\})$ and the map $M([n])\to A([n])$ exhibit $M([n])$ as \[M([n])\cong A([n])\times M(\{0\}).\] 
\end{romanenum}
    For $A\in\Mon(\topo{X})$ let $\LMon_{A}(\mathcal X)\coloneqq \LMon{}(\mathcal X)\times_{\Mon(\topo{X})}\{A\}$. 
\end{definition} 
By \cite[Proposition 4.2.2.9]{higheralgebra}, there exists a functor $\Delta^{\operatorname{op}}\times\Delta^1\to \LM^{\otimes}$ such that for all categories $\mathcal C$ with finite products, pullback along $j$ induces an equivalence of categories \[\Fun^{\operatorname{lax}}(\mathcal LM^{\otimes}, \topo{C})\cong \LMon(\topo{C}).\] 
Composing with the equivalence from \cref{algebrasarelaxfunctors} we obtain an equivalence 
    \[\LMod{}{\mathcal C^{\times}}\cong \LMon(\topo{C}).\]
    Since $\Delta^{\operatorname{op}}\times \{0\}\hookrightarrow \Delta^{\operatorname{op}}\times\Delta^1\to \LM^{\otimes}$ factors as $\Delta^{\operatorname{op}}\xrightarrow{i}\Assoc^{\otimes}\to\LM$, the diagram 
    \begin{center}
        \begin{tikzcd}
        \LMod{}{\mathcal C^{\times}}\arrow[r,"\cong"]\arrow[d] & \Fun^{\operatorname{lax}}(\LM^{\otimes}, \mathcal C)\arrow[r]\arrow[d] & \LMon(\mathcal C)\arrow[d]\\ 
        \Alg(\mathcal C^{\times})\arrow[r,"\cong"]& \Fun^{\operatorname{lax}}(\Assoc^{\otimes}, \mathcal C)\arrow[r] & \Mon(\mathcal C)
        \end{tikzcd}
    \end{center} commutes. 
    This implies:
\begin{cor}\label{modulecategoriesincartesian}Suppose $\topo{X}$ is a category with finite products. 
    For $G\in\Alg(\topo{X}^{\times})\cong \Mon(\topo{X})$, there is an equivalence 
\[ \LMod{G}{\topo{X}}\cong \LMon_{G}(\topo{X}).\]
This is natural in the algebra $G$ and with respect to finite products-preserving functors. \end{cor}

\begin{ex}{\cite[page 656]{higheralgebra}}\label{compositionmonoidalstructure}
    The large category of small categories $\Cat$ is cartesian closed with internal Hom $\Fun(-,-)$ (\cite[Remark 4.8.1.6]{higheralgebra}). 
    \begin{romanenum}
    \item For categories $\mathcal C, \mathcal D$ and $\mathcal E$, this yields 
    \textit{evaluation} 
    \[ \Fun(\mathcal C, \mathcal D)\times\mathcal C\to\mathcal D\] and \textit{composition maps}
    \[\Fun(\mathcal C, \mathcal D)\times\Fun(\mathcal D, \mathcal E)\to\Fun(\mathcal C, \mathcal E).\] 
    \item For a category $\mathcal C$, the composition \[\Fun(\mathcal C, \mathcal C)\times \Fun(\mathcal C, \mathcal C)\to \Fun(\mathcal C, \mathcal C)\] allows to define a simplicial object \[\Fun(\mathcal C, \mathcal C)^*\in \Fun(\Delta^{\operatorname{op}}, \Cat), n\mapsto \Fun(\mathcal C, \mathcal C)^{n}\] which exhibits $\Fun(\mathcal C, \mathcal C)$ as a monoid in $\Cat$. 
    We refer to the corresponding monoidal structure $\Fun(\mathcal C, \mathcal C)^{\circ}$ on $\Fun(\mathcal C, \mathcal C)$ as the composition monoidal structure. 
    \item The evaluation $\Fun(\mathcal C, \mathcal C)\times\mathcal C\to\mathcal C$ allows to construct a simplicial object \begin{align*}\Fun(\mathcal C, \mathcal C)^*\times\mathcal C\in\oc{\Fun(\Delta^{\operatorname{op}}, \Cat)}{\Fun(\mathcal C, \mathcal C)^*}, \\ [n]\mapsto \Fun(\mathcal C, \mathcal C)^n\times \mathcal C\end{align*} which exhibits $\mathcal C$ as left $\Fun(\mathcal C, \mathcal C)^*$-monoid. Via the equivalence \[\LMod{\Fun(\mathcal C, \mathcal C)^{\circ}}{\Cat}\cong \LMon_{\Fun(\mathcal C, \mathcal C)^*}(\Cat)\] from \cref{modulecategoriesincartesian}, this exhibits $\mathcal C$ as left-tensored over $\Fun(\mathcal C, \mathcal C)^{\circ}$. 

    \item If $\mathcal C^{\otimes}$ is a monoidal category, denote by $\mathcal C^{\otimes}_{\operatorname{rev}}$ its reverse (\cite[Remark 4.1.1.7]{higheralgebra}) and consider it as object of $\Mon(\mathcal C)$. 
    The monoidal structure on $\mathcal C$ allows to enhance \[\mathcal C\to \Fun(\mathcal C, \mathcal C), \, c\mapsto c\otimes-\] to a map \[\mathcal C_{rev}^{\otimes}\to \Fun(\mathcal C, \mathcal C)^{*}\in \Mon(\vlCat).\] 
    We thereby obtain a symmetric monoidal functor $\mathcal C^{\otimes}_{\operatorname{rev}}\to \Fun(\mathcal C, \mathcal C)^{\circ}$. 
    \end{romanenum}

\end{ex}
\subsection{Sheaf topoi and morphisms of sites}\label{section:topologies}
In this section, we record basic terminology and results on sheaf topoi. 
All results are immediate consequences from their 1-categorical analogues and the discussion in \cite[section 6.2.2]{highertopostheory}. 

\begin{definition}[{\cite[Definition 6.2.2.1]{highertopostheory}}]
Suppose $\mathcal C$ is a small category.
\begin{itemize}
\item A \emph{sieve} on an object $c\in\mathcal C$ is a full subcategory $S\subseteq \oc{\mathcal C}{c}$ such that for all $a\to b\in \oc{\mathcal C}{c}$ with $b\in S$: $a\in S$. 

\item A Grothendieck topology $\tau$ on $\mathcal C$ consists of a collection of sieves $\Cov_{\tau}(c)$ on $c$ for all objects $c\in\mathcal C$ such that 
\begin{romanenum}
\item For all $c\in\mathcal C$, $\oc{\mathcal C}{c}\in \Cov_{\tau}(c)$
\item For $f\colon b\to c\in\mathcal C$ and $\mathcal C^0\subseteq \oc{\mathcal C}{c}\in \Cov_{\tau}(c)$, $f^*(\mathcal C^0)\coloneqq \oc{\mathcal C}{b}\times_{\oc{\mathcal C}{c}}\mathcal C^0\in\Cov_{\tau}(b)$. 
\item Suppose that $S\in\Cov_{\tau}(c)$ and $\mathcal C^0\subseteq \oc{\mathcal C}{c}$ is a sieve on $c$. If for all $f\colon s\to c\in S\subseteq \oc{\mathcal C}{c}$, $f^*(C^0)\coloneqq  \oc{\mathcal C}{b}\times_{\oc{\mathcal C}{c}}\mathcal C^0\in\Cov_{\tau}(s)$, then $C^0\in \Cov_{\tau}(\mathcal C)$. 
\end{romanenum}
\item A small category together with a Grothendieck topology is called a site.
\end{itemize} 
\end{definition}

A sieve $S$ on an object $c\in\mathcal C$ defines a monomorphism $S\subseteq \Map_{\mathcal C}(-,c)\in \mathcal P(\mathcal C)=\Fun(\mathcal C^{\operatorname{op}}, \an)$. This is pointwise an inclusion of connected components, and determines the sieve $S$. We will in the sequel freely identify sieves with the associated subfunctors of $\Map_{\mathcal C}(-,c)$. 
We denote by $h_c\coloneqq \Map_{\mathcal C}(-,c)\colon\mathcal C^{\operatorname{op}}\to \an$ the presheaf represented by $c$. 
\begin{definition}
    If $(\mathcal C, \tau)$ is a Grothendieck topology, denote by $\Shv_{\tau}(\mathcal C)\subseteq \mathcal P(\mathcal C)$ the full subcategory on functors such that for all $c\in\mathcal C$ and all sieves $S\in\Cov_{\tau}(c)$, the inclusion \[S\subseteq \Map_{\mathcal C}(-,c)=:h_c\] induces an equivalence 
    \[ F(c)=\Map_{\mathcal P(\mathcal C)}(h_c,F)\cong \Map_{\mathcal P(\mathcal C)}(S,F).\] 
    \end{definition}
    By \cite[Proposition 6.2.2.7]{highertopostheory}, $\Shv_{\tau}(\mathcal C)\subseteq \mathcal P(\mathcal C)$ admits a left-exact left adjoint, whence $\Shv_{\tau}(\mathcal C)$ is a topos. 
    We call this left adjoint $\tau$-sheafification. 

    It is often convenient to characterise Grothendieck topologies/sheaf topoi in terms of coverages.
\begin{definition}\label{definitioncoverage}
    Suppose that $\mathcal C$ is a small category. 
    \begin{romanenum}
        \item A cover for an object $c\in\mathcal C$ is a set of maps $\mathcal T=\{ c_i\to c\}_{i\in I}$. 
        \item For a cover $\mathcal T=\{c_i\to c\}_{i\in I}$ denote by $S(\mathcal T)\subseteq \oc{\mathcal C}{c}$ the full subcategory on objects $x\to c$ such that there exists $i\in I$ such that $q$ factors as $x\to c_i\to c$. This is a sieve on $c$. 
        \item A coverage $\tau$ on $\mathcal C$ consists of a set $\tau(c)$ of covers for every object $c\in\mathcal C$ such that the following is satisfied: 
For every cover $\{ c_i\to c\}_{i\in I}\in\tau(c)$ and every morphism $f\colon b\to c\in\mathcal C$, there exists a covering family $\{ b_j\to b\}_{j\in J}\in \tau(b)$ such that for all $j\in J$, $b_j\to b\to c\in S(\{ c_i\to c\})$, i.e.\ there exists $i(j)\in I$ such that 
$b_j\to b\to c$ factors as $b_j\to c_{i(j)}\to c$.  

\item We call a small category $\mathcal C$ together with a coverage $\tau$ a quasi-site. 
    \end{romanenum}
\end{definition} 
Every Grothendieck topology determines a coverage by allowing those covering families $\{ x_i\to x\}_{i\in I}$ for which $S(\{x_i\to x\}_{i\in I})$ is a sieve. Conversely, every coverage determines a Grothendieck topology: 
\begin{definition}\label{definitiongrothendiecktopologyviacovers}
 Suppose that $(\mathcal C, \tau)$ is a quasi-site. 
 \begin{romanenum}
\item If $\rho$ is a Grothendieck topology on $\mathcal C$, we say that $\tau\subseteq \rho$ if for all $c\in\mathcal C$ and all covering families $\mathcal X=\{ c_i\to c\}\in \tau(c)$, $S(\mathcal X)\in \Cov_{\rho}(c)$. 
For $c\in\mathcal C$, let \[\lbrack \tau \rbrack(c)\coloneqq \bigcap_{\tau\subseteq \rho} \Cov_{\rho}(c), \] where the intersection runs over all Grothendieck topologies $\rho$ on $\mathcal C$ with $\tau\subseteq \rho$. 
We call $\lbrack \tau \rbrack$ the Grothendieck topology generated by $\tau$. 
It is straightforward to check that $\lbrack \tau\rbrack$ is a Grothendieck topology. As the collection of all sieves forms a Grothendieck topology, the intersection on the right-hand side is over a non-empty set of Grothendieck topologies. 

\item If for all $c\in \mathcal C$, \[\lbrack \tau \rbrack(c)=\{ S\in \operatorname{Sieve}(c)\, | \, \exists \mathcal X\in \tau(c) \text{ s.t. } S(\mathcal X)\subseteq S\}, \] we say that $\tau$ is a \emph{Grothendieck coverage}.
\item If $\tau$ is a coverage on $\mathcal C$, a functor $F\colon\mathcal C^{\operatorname{op}}\to\an$ is a $\tau$-sheaf of for all $c\in\mathcal C$ and all covering families $\{c_i\to c\}\in\tau(c)$, the inclusion  $S(\{c_i\to c\}_{i\in I})\subseteq \Map_{\mathcal C}(-,c)$ induces an equivalence \[F(c)\cong \Map_{\mathcal P(\mathcal C)}(S(\{c_i\to c\}_{i\in I}),F).\]
 \end{romanenum} 
\end{definition}
\begin{rem}\label{remarkgrothendieckcoverage}
Suppose that $\tau$ is a Grothendieck topology, $S\in \Cov_{\tau}(c)$ and $\tilde S\in\operatorname{Sieve}(c)$ is a sieve with $S\subseteq \tilde{S}$. 
For $f\colon x\to c\in S, \oc{\mathcal C}{c}=f^*(S)\subseteq f^*(\tilde S)$, which implies that $f^*(\tilde S)=\oc{\mathcal C}{x}\in\Cov_{\tau}(x)$, and hence $\tilde S\in \Cov_{\tau}(c)$ by the third axiom in the definition of Grothendieck topologies. 
This implies that for all coverages $\tau$, \[\lbrack \tau \rbrack(c)\supseteq\{ S\in \operatorname{Sieve}(c)\, | \, \exists \mathcal X\in \tau(c) \text{ s.t. } S(\mathcal X)\subseteq S\}.\]

If $\tau$ is a Grothendieck topology and $\tilde{\tau}$ the associated coverage, then \[\tilde{\tau}(c)=\{ S\in \operatorname{Sieve}(c)\, | \, \exists \mathcal X\in \Cov_{\tau}(c) \text{ s.t. } S(\mathcal X)\subseteq S\}, \] i.e. $\tilde{\tau}$ is a Grothendieck coverage. 

If the category $\mathcal C$ has \textit{enough} pullbacks, there are more explicit characterizations of Grothendieck coverages, see e.g.\ \cite[Definition A.1]{Pstragowskisynthetic}.
\end{rem}

The assignment of a Grothendieck topology to a covering family is not one-to-one (but onto). 
We now show that it is one-to-one on associated categories of sheaves.  
\begin{lemma}\label{Grothendiecktopologyofcoverage}
\begin{romanenum}
\item If $(\mathcal C, \tau)$ is a quasi-site, then $(\mathcal C, \lbrack \tau\rbrack)$ is a site. 
\item A functor $F\colon \mathcal C^{\operatorname{op}}\to \an$ is a $\lbrack \tau\rbrack$-sheaf if and only if it is a $\tau$-sheaf. 
\end{romanenum}
\end{lemma}
\begin{proof}This is analogous to {\cite[Theorem 3.14]{landcondensed}}, who proved the statement for $\mathcal C$ a 1-category. 
    As the collection of all sieves forms a Grothendieck topology, there exist Grothendieck topologies which contain $\tau$.
    By definition, every $\lbrack \tau \rbrack$-sheaf is a $\tau$-sheaf, so it remains to prove the converse. 
    For $c\in\mathcal C$ denote by $\Cov_{\tau}(c)$ the collection of sieves $S$ on $c$ such that for all $f\colon b\to c\in\mathcal C$ and all $\tau$-sheaves $F$, \[F(b)\cong \Map_{\mathcal P(\mathcal C)}(h_b,F)\cong \Map_{\mathcal P(\mathcal C)}(f^*S,F).\] 
We claim that this defines a Grothendieck topology on $\mathcal C$. 
Clearly, $\oc{\mathcal C}{c}\in \Cov_{\tau}(c)$ for all $c\in\mathcal C$ and for $S\in \Cov_{\tau}(c)$ and $f\colon b\to c\in \mathcal C$, $f^*(S)\in \Cov_{\tau}(b)$. 
Suppose now that $R\subseteq \oc{\mathcal C}{c}$ is a sieve and $S\in \Cov_{\tau}(c)$ such that for all $f\colon s\to c\in S$, $f^*(R)\in \Cov_{\tau}(s)$.
We have to show that for all $q\colon b\to c\in\mathcal C$ and $\tau$-sheaves $F$, 
$q^*(R)\subseteq h_b$ induces an equivalence 
\[ F(b)\cong \Map_{\mathcal P(\mathcal C)}(q^*(R),F), \] then it follows that $R\in \Cov_{\tau}(c)$.  
Since $q^*(S)\in \Cov_{\tau}(b)$ and for all $t\colon a\to b\in q^*(S)$, 
\[ t^*q^*(R)=(q\circ t)^*(R)\in \Cov_{\tau}(a)\] by assumption on $R$, it suffices to check this for $q=id$. 
Let $\tilde R\coloneqq R\times_{h_c}S$. This is a subpresheaf of $R$ and $S$. 
For $f\colon b\to c\in R$, $f^*(\tilde R)\cong f^*(S)\in \Cov_{\tau}(b)$, so there is a commutative diagram of presheaves
\begin{center}
    \begin{tikzcd}
    f^*(\tilde R)=f^*(S)\arrow[d,"f"]\arrow[r] & h_b=f^*(R)\arrow[d]\\ 
    \tilde R\arrow[r] & R 
    \end{tikzcd}
\end{center}
This yields a commutative diagram 
\begin{center}
    \begin{tikzcd}
    \colim{f\colon b\to c\in R}f^*(\tilde R)\arrow[d,"f"]\arrow[r] & \colim{b\to c\in R} h_b\arrow[d,"\cong"]\\ 
    \tilde R\arrow[r] & R 
    \end{tikzcd}
\end{center}
The left vertical map is an equivalence as $\colim{f\colon b\to c\in R}f^*\tilde{R}=R\times_{h_c}\tilde R$ (by universality of colimits in $\mathcal P(\mathcal C)$) and $R\times_{h_c}\tilde R=\times_{h_C}R\times_{h_C}S\cong R\times_{h_C}S\cong \tilde R$ since $R\to h_c$ is a monomorphism. 
Applying a $\tau$-sheaf $F$ yields a commutative diagram 
\begin{center}
    \begin{tikzcd}
    \clim{f\colon b\to c\in R}\Map_{\mathcal P(\mathcal C)}(f^*\tilde R,F) & \arrow[l,"\cong"] \clim{f\colon b \to c\in R}F(b) \\ 
    \Map_{\mathcal P(\mathcal C)}(\tilde R,F)\arrow[u,"\cong"] & \arrow[l]\arrow[u,"\cong"]\Map_{\mathcal P(\mathcal C)}(R,F). 
    \end{tikzcd}
\end{center}
The top horizontal map is an equivalence as for all $f\colon b\to c\in R$, $f^*(\tilde R)=f^*(S)\in \Cov_{\tau}(b)$.
It therefore suffices to show that for all $\tau$-sheaves $F$, $\tilde R\subseteq R\subseteq h_c$ induces an equivalence \[F(c)\cong \Map_{\mathcal P(\mathcal C)}(\tilde R,F),\] then it follows that $F(c)\cong \Map_{\mathcal P(\mathcal C)}(R,F)$ for all $\tau$-sheaves $F$.  
By universality of colimits in $\mathcal P(\mathcal C)$, 
$\tilde R=\colim{f\colon s\to c\in S}R\times_{h_c}h_s=\colim{f\colon s\to c\in S}f^*(R)$. 
Since for all $f\colon s\to c\in S$, $f^*(R)\in \Cov_{\tau}(s)$ by assumption, this implies that for all $\tau$-sheaves $F$, \[\Map_{\mathcal P(\mathcal C)}(\tilde R,F)\cong \clim{f\colon s\to c\in S}\Map_{\mathcal P(\mathcal C)}(f^*R,F)\cong \clim{f\colon s\to c\in S}\Map_{\mathcal P(\mathcal C)}(s,F)\cong \Map_{\mathcal P(\mathcal C)}(S,F).\] 
As $S\in \Cov_{\tau}(c)$, the right-hand side is $F(c)$. 
This proves that $\Cov_{\tau}$ defines a Grothendieck topology on $\mathcal C$. 
By construction, every $\tau$-sheaf is a $\Cov_{\tau}$-sheaf and $\tau\subseteq \Cov_{\tau}$. 
This implies that $\lbrack \tau \rbrack \subseteq \Cov_{\tau}$ and in particular, every $\tau$-sheaf is a $[\tau]$-sheaf. 
\end{proof}

\begin{rem}\label{subsieveascolimitoverchechnerve}
    For a covering $\mathcal X=\{ c_i\to c\}_{i\in I}\in\mathcal C$, \[S(\mathcal X)\cong \colim{\Delta^{\operatorname{op}}}\check{C}(\sqcup_{i\in I}\Map_{\mathcal C}(-,c_i)\to \Map_{\mathcal C}(-,c))\in \mathcal P(\mathcal C).\] 
    Indeed: Recall that a map $A\to B\in\mathcal P(\mathcal C)$ is monomorphism if and only if $A\cong A\times_B A$.  Universality of colimits in $\mathcal P(\mathcal C)$ therefore implies that 
    \[ \colim{\Delta^{\operatorname{op}}}\check{C}(\sqcup_{i\in I}\Map_{\mathcal C}(-,c_i)\to \Map_{\mathcal C}(-,c))\to \Map_{\mathcal D}(-,c)\] is a monomorphism, i.e.\ pointwise an inclusion of connected components. 
    At $t\in\mathcal C$, its image consists precisely of the maps $t\to c$ which factor over $c_i$ for some $i\in I$. 
\end{rem}
\begin{definition}\label{sheavesquasitopology}
If $\tau$ is a coverage on an essentially small category $\mathcal C$, we denote by \[\hypershv_{\tau}(\mathcal C)\subseteq \Shv_{\tau}(\mathcal C)\subseteq \mathcal P(\mathcal C)\] the categories of (hypercomplete) $\lbrack \tau \rbrack$-sheaves. 

By \cite[Proposition 6.2.2.7, 6.2.1.1/p.666]{highertopostheory} and \cref{Grothendiecktopologyofcoverage}, the inclusions \[\hypershv_{\tau}(\mathcal C)\subseteq \Shv_{\tau}(\mathcal C)\subseteq \mathcal P(\mathcal C)\] admit left-exact left adjoints, called $\tau$-(hyper)-sheafification. 
\end{definition}
Next, we recall how functors between (quasi-)sites induce functors on the associated sheaf topoi. We used this to model $\kappa$-condensed animae on other categories than $\Pro(\Fin)_{\kappa}$ and to compare condensed and sheaf cohomology.  
\begin{definition}\label{definitionmorphismsites}

    \begin{itemize} \item Suppose $(\mathcal C, \tau_{\mathcal C}),(\mathcal D, \tau_{\mathcal D})$ are quasi-sites. 
    A functor $f\colon \mathcal C\to\mathcal D$ is 
    \begin{romanenum}\item \emph{continuous} if for all $c\in\mathcal C$ and $\{ c_i\to c\}_{i\in I}\in \tau_{\mathcal C}(c)$,
    $\{ f(c_i)\to f(c)\}_{i\in I}\in \tau_{\mathcal D}(f(c))$ is a cover. 
    \item \emph{covering flat} if  the composition \[\Shv_{\tau_{\mathcal C}}(\mathcal C)\subseteq \mathcal P(\mathcal C)\xrightarrow{f_{!}}\mathcal P(\mathcal D)\xrightarrow{L_{\mathcal D}}\Shv_{\tau_{\mathcal D}}(\mathcal D)\] of the left Kan extension along $f$ with $\tau_{\mathcal D}$-sheafification preserves finite limits.
    \item A \emph{morphism of quasi-sites} if it is continuous and covering flat.
    \end{romanenum}

    \item Suppose $(\mathcal C, \tau_{\mathcal C}),(\mathcal D, \tau_{\mathcal D})$ are sites. 
    A functor $f\colon \mathcal C\to\mathcal D$ is 
    \begin{romanenum}\item \emph{continuous} if for all $c\in\mathcal C$ and $S\in \Cov_{\tau_{\mathcal C}}(c)$, \[f_*(S)\coloneqq \{ x\in \oc{\mathcal C}{f(c)}\, |\, \exists s\in S \text{ and } x\to f(s)\in\oc{C}{f(c)}\}\] is a $\tau_{\mathcal D}$-sieve.  
    \item \emph{covering flat} if  the composition \[\Shv_{\tau_{\mathcal C}}(\mathcal C)\subseteq \mathcal P(\mathcal C)\xrightarrow{f_{!}}\mathcal P(\mathcal D)\xrightarrow{L_{\mathcal D}}\Shv_{\tau_{\mathcal D}}(\mathcal D)\] of the left Kan extension along $f$ with $\tau_{\mathcal D}$-sheafification preserves finite limits.
    \item A \emph{morphism of sites} if it is continuous and covering flat.
    \end{romanenum}
\end{itemize}
\end{definition}

\begin{ex}
If $\tau_{1}, \tau_{2}$ are two Grothendieck topologies (coverages) on a category $\mathcal C$ such that for all $c\in\mathcal C$, $\tau_{1}(c)\subseteq \tau_2(c)$, the \textit{identity} $(\mathcal C, \tau_1)\to(\mathcal C, \tau_2)$ is a morphism of (quasi-)sites. 
\end{ex}

\begin{lemma}\label{criterioncoveringflat}
Suppose that $(\mathcal C, \tau_{\mathcal C}),(\mathcal D, \tau_{\mathcal D})$ are (quasi-)sites, $\mathcal C$ has finite limits and $f\colon\mathcal C\to\mathcal D$ is a functor which preserves finite limits. 
\begin{romanenum}
\item Then $f$ is covering-flat. 
\item Moreover, if $f$ is a continuous functor of quasi-sites, then the induced functor of sites \[(\mathcal C,[\tau_{\mathcal C}])\to (\mathcal D,[\tau_{\mathcal D}])\] is continuous. 
\end{romanenum}
\end{lemma}
\begin{proof} 
    Since $\mathcal C\xrightarrow{f}\mathcal D\hookrightarrow\mathcal P(\mathcal D)$ preserves finite limits, the left Kan extension along $f$, 
\[ \mathcal P(\mathcal C)\to \mathcal P(\mathcal D)\] preserves finite limits by \cite[Proposition 6.2.3.20]{highertopostheory}, hence $f$ is covering flat. 

Suppose now that $f$ is a continuous functor of quasi-sites.
We claim that \[\Cov_{f^*[\tau_{\mathcal D}]}(c)\coloneqq \{ S\in \operatorname{Sieve}(c)\, |\, f_*(S)\in \Cov_{[\tau_{\mathcal D}]}(f(c))\}\] defines a  Grothendieck topology $f^*[\tau_{\mathcal D}]$ on $\mathcal C$. 
As $\tau_{\mathcal C}\subseteq f^*[\tau_{\mathcal D}]$ and $f(f^*[\tau_{\mathcal D}])\subseteq [\tau_{\mathcal D}]$, this then implies that $f\colon (\mathcal C,[\tau_{\mathcal C}])\to (\mathcal D,[\tau_{\mathcal D}])$ is continuous. 
For all $c\in\mathcal C$, $f_*(\oc{\mathcal C}{c})=\oc{\mathcal D}{f(c)}$, whence $\oc{\mathcal C}{c}\in f^*[\tau_{\mathcal D}](c)$. 
Suppose that $c\in\mathcal C$, $S\subseteq \oc{\mathcal C}{c}$ is a covering sieve and $\alpha\colon b\to c\in \mathcal C$. 
We claim that \[(f\alpha)^*(f_*(S))= f_*(\alpha^*(S)).\] 
For $x\to f(b)\in (f\alpha)^*(f_*(S))$, there exist $t\colon s\to c\in S$ and a morphism \[(x\to f(b)\xrightarrow{f\alpha} f(c))\to (f(s)\xrightarrow{ft} f(c))\in\oc{\mathcal C}{f(c)}.\] 
Since $\mathcal C$ has pullbacks and $f$ preserves them, this implies that 
$x\to f(b)$ factors over \[x\to f(b\times_c s)\to f(b).\] 
As $S$ is a sieve, $b\times_c s\to b\in \alpha^*(S)$, which shows that $x\to f(b)\in f_*(\alpha^*(S))$. 
Conversely, suppose that $u\colon x\to f(b)\in f_*(\alpha^*(S))$. 
Then there exist $i\colon a\to b\in\mathcal C$, a morphism \[(a\to b\xrightarrow{\alpha} c)\to (s\to c)\in\oc{\mathcal C}{c}\] with $s\to c\in S$, and a morphism $(x\xrightarrow{u} f(b))\to (f(a)\xrightarrow{f i} f(b))\in \oc{\mathcal D}{f(b)}$. 
This shows that $(u\colon x\to f(b))\in (f\alpha)^*(f_*(S))$, and hence \[(f\alpha)^*(f_*(S))= f_*(\alpha^*(S)).\] 
As $[\tau_{\mathcal D}]$ is a Grothendieck topology, it follows that for all $c\in\mathcal C, S\in \Cov_{f^*[\tau_{\mathcal D}]}(c)$ and $\alpha\colon b\to c\in\mathcal C$, \[\alpha^*(S)\in \Cov_{f^*[\tau_{\mathcal D}]}(b).\] 
Suppose that $S,T\subseteq \oc{\mathcal C}{c}$ are sieves on $c$ such that $S\in \Cov_{f^*[\tau_{\mathcal D}]}(c)$ and for all $\alpha\colon b\to c\in S$, 
\[\alpha^*(T)\in \Cov_{f^*[\tau_{\mathcal D}]}(b).\] 
We claim that $T\in \Cov_{f^*[\tau_{\mathcal D}]}(b)$, i.e. $f_*(T)\in [\tau_{\mathcal D}]$, this then proves that $f^{*}[\tau_{\mathcal D}]$ is a Grothendieck topology. 
As $[\tau_{\mathcal D}]$ is a Grothendieck topology, by the third condition on Grothendieck topologies it suffices to show that for $\beta\colon x\to f(c)\in f_*(S)$, \[\beta^*(f_*(T))\in \Cov_{[\tau_{\mathcal D}]}(b).\] 
Suppose that $\beta\colon x\to f(c)\in f_*(S)$.
By definition of $f_*(S)$, $\beta$ factors as \[x\xrightarrow{j} f(s)\xrightarrow{f(\alpha)}f(c)\] with $\alpha\colon s\to c\in S\subseteq \oc{C}{c}$, and hence \[\beta^*(f_*(T))=j^*(f\alpha)^*(f_*T).\] 
We showed above that $(f\alpha)^*(f_*T)=(f_*(\alpha^*(T))).$ By our assumption that $\alpha^*(T)\in \Cov_{f^*[\tau_{\mathcal D}]}(s)$, $f_*(\alpha^*(T))\in [\tau_{\mathcal D}]$, and hence \[\beta^*(f_*(T))=j^*(f_*(\alpha^*(T)))\in [\tau_{\mathcal D}]\] by the second condition for Grothendieck topologies. 
\end{proof}

\begin{cor}\label{morphismofquasisitesinducesfunctoronsheaves}
\begin{romanenum}\item If $f\colon (\mathcal C, \tau_{\mathcal C})\to (\mathcal D, \tau_{\mathcal D})$ is a continuous functor of (quasi-)sites, then  
\begin{align*} \mathcal P(\mathcal D)&\to \mathcal P(\mathcal C), \\ F& \mapsto F\circ f\end{align*} restricts to a functor  
\[ f_*\colon \Shv_{\tau_{\mathcal D}}(\mathcal D)\to \Shv_{\tau_{\mathcal C}}(\mathcal C).\]
\item If $f$ is in addition covering flat, i.e.\ a morphism of (quasi-)sites, then $f_*$ restricts further to a functor \[ \hypershv_{\tau_{\mathcal D}}(\mathcal D)\to \hypershv_{\tau_{\mathcal C}}(\mathcal C), \] and \[f_*\colon \mywidehatshv_{\tau_{\mathcal D}}(\mathcal D)\to \mywidehatshv_{\tau_{\mathcal C}}(\mathcal C)\] is the right adjoint of a geometric morphism.
\end{romanenum} 
\end{cor}
\begin{notation}\label{functoronsheavesinducedbymorphismonquasisites}
If $f\colon  (\mathcal C, \tau_{\mathcal C})\to (\mathcal D, \tau_{\mathcal D})$ is a morphism of (quasi-)sites, we denote by $f_*\colon \mywidehatshv_{D}(\mathcal D)\to \mywidehatshv_{\tau}(\mathcal C), F\mapsto F\circ f$ the functor denoted by $f^*$ in the above lemma and by $f^*$ its left adjoint.
\end{notation}
\begin{proof}
    The first statement is obvious from the definition of continuity. 
    Sheafification after left Kan extension \[ f_{!}^L\colon \Shv_{\tau_{\mathcal C}}(\mathcal C)\subseteq \mathcal P(\mathcal C)\xrightarrow{f_{!}}\mathcal P(\mathcal D)\xrightarrow{L_{\tau_{\mathcal D}}}\Shv_{\tau_{\mathcal D}}(\mathcal D)\] is a left adjoint of \[ f_*\colon \Shv_{\tau_{\mathcal D}}(\mathcal D)\to \Shv_{\tau_{\mathcal C}}(\mathcal C).\] If $f$ is covering flat, then $f_{!}^L$ is left-exact. By \cite[Proposition 6.5.1.16.(4)]{highertopostheory}, this implies that $f_*$ preserves hypercomplete objects. 
    Hypersheafification after left Kan extension \[ \hat{f}_{!}^L\colon \hypershv_{\tau_{\mathcal C}}(\mathcal C)\subseteq \mathcal P(\mathcal C)\xrightarrow{f_{!}}\mathcal P(\mathcal D)\xrightarrow{L_{\tau_{\mathcal D}}}\Shv_{\tau_{\mathcal D}}(\mathcal D)\to \hypershv_{\tau_{\mathcal D}}(\mathcal D)\] is a left-exact left adjoint of 
    \[ f_*\colon \hypershv_{\tau_{\mathcal D}}(\mathcal D)\to \hypershv_{\tau_{\mathcal C}}(\mathcal C).\] 
    (Hypercompletion is left-exact by \cite[Proposition 6.2.1.1, Proposition 6.5.1.16]{highertopostheory}.) 
\end{proof}
\begin{definition}\label{definitioncoveringliftingproperty}
    Suppose $(\mathcal C, \tau_{\mathcal C})$ and $(\mathcal D, \tau_{\mathcal D})$ are quasi-sites. 
    A functor $f\colon \mathcal C\to \mathcal D$ has the \emph{covering lifting property} if for all $\{p_i\colon d_i\to d\}_{i\in I}\subseteq \tau_{\mathcal D}$, and all $f(c)\to d\in \mathcal C/d$ there exists a cover $\{q_j\colon c_j\to c\}_{j\in J}\subseteq \tau_{\mathcal C}$ such that for all $j\in J$, there exists $i(j)\in I$ such that $h_{f(q_j)}$ factors as \[h_{f(c_j)}\to h_{d_{i(j)}}\times_{h(d)}h_{f(c)}\to h_{f(c)}.\] 
\end{definition}
\begin{rem}
    If $\mathcal D$ has all pullbacks and $\tau_{\mathcal D}$ is stable under pullbacks, i.e.\ for $\{ d_i\to d\}\in\tau_{\mathcal D}(d)$ and $e\to d\in\mathcal D$. $\{ d_i\times_d e\to e\}\in\tau_{\mathcal D}(e)$, a functor has the covering lifting property if and only if for every $c\in\mathcal C$ and every cover $\{p_i\colon d_i\to f(c)\}_{i\in I}\subseteq \tau_{\mathcal D}$, there exists a cover  $\{q_j\colon c_j\to c\}_{j\in J}\subseteq \tau_{\mathcal C}$ such that for all $j\in J$, there exists $i(j)\in I$ such that $f(i_j)$ factors as $f(c_j)\to d_{i(j)}\to f(c)$. 
\end{rem}
\begin{lemma}\label{rightkanextensionsheaves}
    Suppose $(\mathcal C, \tau_{\mathcal C})$ and $(\mathcal D, \tau_{\mathcal D})$ are quasi-sites and $f\colon \mathcal C\to\mathcal D$ is a functor which has the covering lifting property. 
    Right Kan extension $f_*\colon \mathcal P(\mathcal C)\to\mathcal P(D)$ along $f$ restricts to functors 
    \[ {\Shv}_{\tau_{\mathcal C}}(\mathcal C)\to {\Shv}_{\tau_{\mathcal D}}(\mathcal D)\text{ and } \hypershv_{\tau_{\mathcal C}}(\mathcal C)\to \hypershv_{\tau_{\mathcal D}}(\mathcal D).\]
    Both functors are right adjoints of a geometric morphism. 
\end{lemma} 
\begin{proof}
    Denote by $f^*\colon\mathcal P(\mathcal D)\to\mathcal P(\mathcal C), F\mapsto F\circ f$ the pullback along $f$. 
    If $\mathcal X\coloneqq \{d_j\to d\}_{j\in J}$ is a covering in $\mathcal D$, denote by $S(\mathcal X)\subseteq h_d\in\mathcal P(\mathcal D)$ the subpresheaf on maps $c\to d$ which factor over $c\to d_i\to d$ for some $i\in I$. We will show below that $L_{\tau_{\mathcal C}}(f^*S(\mathcal X)\to f^*h_d)\in \Shv_{\tau_{\mathcal C}}(\mathcal C)$ is an equivalence. This then implies that for $F\in \Shv_{\tau_{\mathcal C}}(\mathcal C)$, 
    \begin{align*}\Map_{\mathcal P(\mathcal D)}(h_d,f_*F)& \cong \Map_{\mathcal P(\mathcal C)}(f^*h_d,F)\\ & \cong \Map_{\Shv_{\tau_{\mathcal C}}}(L_{\tau_{\mathcal C}}f^*h_d,F)\\ & \cong \Map_{\Shv_{\tau_{\mathcal C}}}(L_{\tau_{\mathcal C}}f^*S(\mathcal X),F)\\ & \cong \Map_{\mathcal P(\mathcal C)}(f^*S(\mathcal X),F)\\ & \cong \Map_{\mathcal P(\mathcal D)}(S(\mathcal X),f_*F), \end{align*}
    which shows that $f_*F$ is a $\tau_{\mathcal D}$-sheaf. Since colimits in $\mathcal P(\mathcal C)$ are universal, $f^*h_d$ is a colimit of representables, and $L_{\tau_{\mathcal C}}\colon\mathcal P(\mathcal C)\to \Shv_{\tau_{\mathcal C}}(\mathcal C)$ preserves colimits, it suffices to show that for all $h_c\to f^*h_d, c\in\mathcal C$, \[f^*S(\mathcal X)\times_{f^*h_d}h_c\to h_c\] is an $L_{\tau_{\mathcal C}}$-equivalence. Fix \[t\in \Map_{\Psh(\mathcal C)}(h_c,f^*h_d)\cong f^*h_d(c)\cong h_d(f(c))\cong\Map_{\mathcal D}(f(c),d).\] By construction, $f^*S(\mathcal X)\times_{f^*h_d}h_c\subseteq h_c$ is the full subpresheaf on maps $r\colon x\to c$ such that there exists $i\in I$ such that \[f(x)\to f(c)\to d\] factors over $d_i\to d$. Equivalently, $f^*S(\mathcal X)\times_{f^*h_d}h_c(x)\subseteq h_c(x)$ consists of those maps $q\colon x\to c$ such that $h_{f(q)}\colon h_{f(x)}\to h_{f(c)}$ factors over $h_{d_i}\times_{h_{d}} h_{f(c)}\to h_{f(c)}$ for some $i\in I$. 
By assumption, there exists a cover $\mathcal Y\coloneqq \{c_j\to c\}_{j\in J}\in \tau_{\mathcal C}(c)$ such that for all $j\in J$, $h_{f(c_j)}\to h_{f(c)}$ factors over \[h_{f(c_j)}\to h_{d_i}\times_{h_d}h_{f(c)}\to h_{f(c)}\] for some $i=i(j)\in  I$.  
    Denote by $S(\mathcal Y)\subseteq h_c$ the associated full subpresheaf, i.e.\ $S(\mathcal Y)(x)\subseteq \Map_{\mathcal C}(x,c)$ is the subspace on maps $x\to c$ which admit a factorization $x\to c_j\to c$ for some $j\in J$.
    Then \[S(\mathcal Y)\subseteq f^*S(\mathcal X)\times_{f^*h_d}h_c\subseteq h_c.\] 
    As $\{c_i\to c\}_{i\in I}$ is a $\tau_{\mathcal C}$-cover, $S(\mathcal Y)\subseteq j_c$ induces an equivalence $L_{\tau_{\mathcal C}}S(\mathcal Y)\cong L_{\tau_{\mathcal C}}h_c$.
    This implies that \[L_{\tau_{\mathcal C}}(f^*S(\mathcal X)\times_{f^*h_d}h_c)\to L_{\tau_{\mathcal C}}h_c\] admits a section. Since sheafification preserves pullbacks, it preserves monomorphisms (recall that a map $A\to B$ in a category $\mathcal T$ is a monomorphism if and only if $A\cong A\times_B A$). 
    In particular, $L_{\tau_{\mathcal C}}(f^*S(\mathcal X)\times_{f^*h_d}h_c)\to L_{\tau_{\mathcal C}}h_c$ is a monomorphism, and hence an equivalence as it admits a section. This shows that $L_{\tau_{\mathcal C}}f^*S(\mathcal X)\cong L_{\tau_{\mathcal C}}f^*h_d$, and hence right Kan extension restricts to a functor 
    \[ f_*\colon \Shv_{\tau_{\mathcal D}}(\mathcal D)\to \Shv_{\tau_{\mathcal C}}(\mathcal C).\] 
    Since sheafification is left  exact and $f^*\colon\mathcal P(\mathcal D)\to\mathcal P(\mathcal C)$ is also a right adjoint, the composite 
    \[f^{*,L}\colon  \Shv_{\tau_{\mathcal D}}(\mathcal D)\subseteq \mathcal P(\mathcal D)\xrightarrow{f^*}\mathcal P(\mathcal C)\xrightarrow{L_{\tau_{\mathcal D}}}\Shv_{\tau_{\mathcal D}}(\mathcal D)\] is a left-exact left adjoint of $f_*$. 
    By \cite[Proposition 6.5.1.16]{highertopostheory}, $f^{*,L}$ preserves $\infty$-connective morphisms, whence $f_*$ preserves hypersheaves, i.e.\ restricts to a functor 
    \[ \hat{f}_*\colon \hypershv_{\tau_{\mathcal C}}(\mathcal C)\to \hypershv_{\tau_{\mathcal D}}(\mathcal D).\]
    Hypersheafification after restriction $\hypershv_{\tau_{\mathcal D}}(\mathcal D)\subseteq \Shv_{\tau_{\mathcal D}}(\mathcal D)\xrightarrow{f^*}\Shv_{\tau_{\mathcal C}}(\mathcal C)\xrightarrow{\widehat{(-)}_{\mathcal D}}\hypershv_{\tau_{\mathcal C}}(\mathcal C)$ is a left-exact left adjoint to $\hat{f}_*$. Hypercompletion is left-exact by \cite[Proposition 6.2.1.1, Proposition 6.5.1.16]{highertopostheory}. 
\end{proof} 
    If in the situation of \cref{rightkanextensionsheaves}, $f$ is fully faithful, then right Kan extension along $f$ is a fully faithful functor $f_*\colon \mathcal P(\mathcal C)\to \mathcal P(\mathcal D)$: For every $c\in\mathcal C$, $c$ is a terminal object of $\mathcal C\times_{\mathcal D}\{c\}$, hence the counit $f^*f_*\to \id$ is an equivalence.  
    This then implies that right Kan extension restricts to fully faithful functors \[ \Shv_{\tau_{\mathcal C}}(\mathcal C)\hookrightarrow{\Shv}_{\tau_{\mathcal D}}(\mathcal D), \text{ and } \hypershv_{\tau_{\mathcal C}}(\mathcal C)\hookrightarrow\hypershv_{\tau_{\mathcal D}}(\mathcal D).\]
    \begin{cor}\label{criterionrightkanextensionequivalence}Suppose that $i\colon (\mathcal C, \tau_{\mathcal C})\to (\mathcal D, \tau_{\mathcal D})$ is a continuous functor of quasi-sites which has the covering lifting property and $i\colon\mathcal C\to\mathcal D$ is fully faithful. 
\begin{romanenum}
    \item If the restriction $i^*\colon \Shv_{\tau_{\mathcal D}}(\mathcal D)\to \Shv_{\tau_{\mathcal C}}(\mathcal C), F\mapsto F\circ i$ (\cref{morphismofquasisitesinducesfunctoronsheaves}) is conservative, then \[\Shv_{\tau_{\mathcal D}}(\mathcal D)\cong \Shv_{\tau_{\mathcal C}}(\mathcal C)\text{ and } \hypershv_{\tau_{\mathcal C}}(\mathcal C)\cong \hypershv_{\tau_{\mathcal D}}(\mathcal D)\] via restriction and right Kan extension. 

    \item If the left adjoint of right Kan extension $\hat{i}^*\colon \hypershv_{\tau_{\mathcal D}}(\mathcal D)\to \hypershv_{\tau_{\mathcal C}}(\mathcal C)$ is conservative, then $\hypershv_{\tau_{\mathcal C}}(\mathcal C)\cong \hypershv_{\tau_{\mathcal D}}(\mathcal D)$ via right Kan extension.
\end{romanenum}
    \end{cor}
    \begin{proof}
        By \cref{morphismofquasisitesinducesfunctoronsheaves} and \cref{rightkanextensionsheaves}, pullback and right Kan extension along $i$ restrict to an adjoint pair \[i_*\colon \Shv_{\tau_{\mathcal C}}(\mathcal C)\leftrightarrows \Shv_{\tau_{\mathcal D}}(\mathcal D)\colon i^*.\] 
        As $i$ is fully faithful, the right Kan extension $i_*$ is fully faithful.
        Together with conservativity of $i^*$, this implies that $i^*$ and $i_*$ are mutually inverse equivalences. (It from the triangle identities that the unit and counit are equivalences).
        In particular, $i^*$ and $i_*$ restrict to mutually inverse equivalences $\hypershv_{\tau_{\mathcal C}}(\mathcal C)\cong \hypershv_{\tau_{\mathcal D}}(\mathcal D)$. 

        By \cref{rightkanextensionsheaves}, right Kan extension also restricts to a functor $\hypershv_{\tau_{\mathcal C}}(\mathcal C)\cong \hypershv_{\tau_{\mathcal D}}(\mathcal D)$. This is a fully faithful right adjoint. 
        If its left adjoint is conservative, then $i_*$ is an equivalence by the same reasoning as above. 
    \end{proof}
    \begin{cor}\label{criterionrightkanextensionisocovers}
    Suppose that $f\colon (\mathcal C, \tau_{\mathcal C})\to (\mathcal D, \tau_{\mathcal D})$ is a continuous functor of quasi-sites that has the covering lifting property. 
    Suppose that $f\colon \mathcal C\to \mathcal D$ is fully faithful, and for all $d\in\mathcal D$, there exists a cover $\{ f(c_i)\to d\}_{i\in I}\in \tau_{\mathcal D}(d)$ with $c_i\in\mathcal C$. 

    Then right Kan extension along $f$ restricts to equivalences 
    \[  \hypershv_{\tau_{\mathcal C}}(\mathcal C)\cong \hypershv_{\tau_{\mathcal D}}(\mathcal D).\]
    \end{cor}
    \begin{proof}
        By \cref{criterionrightkanextensionequivalence}, it suffices to show that the left adjoint $f^*$ of right Kan extension is conservative. 
        This left adjoint $f^*$ sends a hypersheaf $G\in \hypershv_{\tau_{\mathcal D}}(\mathcal D)$ to the hypersheafification of $G\circ f$. 
        Since $f$ is continuous, restriction along $f$ restricts to a functor \[f_0^*\colon \Shv_{\tau_{\mathcal D}}(\mathcal D)\to \Shv_{\tau_{\mathcal C}}(\mathcal C), A\mapsto A\circ f.\] It therefore suffices to show that if $\phi\colon A\to B\in\hypershv_{\tau_{\mathcal D}}(\mathcal D)$ is such that $f_0^{*}(\phi)\colon A\circ f\to B\circ f$ is $\infty$-connective, then $\phi$ is an equivalence (equivalently, $\infty$-connective).  
        Being a left-exact left adjoint, $f^*_0$ commutes with homotopy groups. This follows from \cite[Proposition 5.5.6.28]{highertopostheory} as for $n\in\mathbb N_0$, $\pi_n=\pi_0\Omega^n=\tau_{\leq 0}\Omega^n$.
        We are therefore reduced to showing that $f_0^*$ restricts to a conservative functor $\tau_{\leq 0}\Shv_{\tau_{\mathcal D}}(\mathcal D)\to \tau_{\leq 0}\Shv_{\tau_{\mathcal C}}(\mathcal C)$. 

        Fix $d\in\mathcal D$ choose a $\tau_{\mathcal D}$-cover $\{f(c_i)\to d\}$. Denote by $h\colon \mathcal D\to\mathcal P(\mathcal D)$ the Yoneda embedding. 
        For $A\in \Shv_{\tau_{\mathcal D}}(\mathcal D)$, \[A(d)\cong \clim{\Delta_{s}}\Map_{\mathcal P(\mathcal D)}(\check{C}(\sqcup_{i\in I}h(f(c_i))\to h(f(d))),A).\] As $\Delta_{s,\leq 1}\subseteq \Delta_{s}$ is left 1-cofinal (\cite[Example 6.14]{du2023reshapinglimitdiagramscofinality}), for 
        $A\in \tau_{\leq 0}\Shv_{\tau_{\mathcal D}}(\mathcal D)$, 
        \begin{align*} A(d)\cong \Eq\left(\prod_{i\in I}A(f(c_i))\rightrightarrows \prod_{i,j\in I}\Map_{\mathcal P(\mathcal D)}(h(f(c_i))\times_{h(d)}h(f(c_j)),A)\right).\end{align*} 
        For $i,j\in I$ choose a small set $(d_q)_{q\in  Q_{i,j}}\subseteq \mathcal D$ with an effective epimorphism \[\pi\colon \sqcup_{q\in Q_{i,j}}h(d_q)\to h(f(c_i))\times_{h(d)}h(f(c_j)).\] 
        For $q\in Q_{i,j}$ choose a cover $\{ f(c_k)\to d_q\}_{k\in K_q}$. 

        By left 1-cofinality of $\Delta_{s,\leq 1}\subseteq \Delta_s\subseteq \Delta$ (\cite[Lemma 6.5.3.7]{highertopostheory}, \cite[Example 6.14]{du2023reshapinglimitdiagramscofinality}), for $A\in \tau_{\leq 0}\Shv_{\tau_{\mathcal D}}(\mathcal D)$, 
        \begin{align*}\Map_{\mathcal P(\mathcal D)}&(h(f(c_i))\times_{h(d)}h(f(c_j)),A)\cong \clim{\Delta}\Map_{\mathcal P(\mathcal D)}(\check{C}(\pi),A)\\
            &\cong\Eq\left(\prod_{q\in Q_{i,j}}A(d_q)\rightrightarrows \prod_{q,r\in Q_{i,j}}\Map_{\mathcal P(\mathcal D)}(h(d_q)\times_{h(f(c_i))\times_{h(d)}h(f(c_j))}h(d_r),A)\right)\\& \hookrightarrow \prod_{q\in Q_{i,j}}A(d_q)\end{align*} is injective. 
        The above argument shows that for $q\in Q_{i,j}$, and $A\in \tau_{\leq 0}\Shv_{\tau_{\mathcal D}}(\mathcal D)$, 
        \[A(d_q)\hookrightarrow \prod_{k\in K_q}A(f(c_k))=\prod_{k\in K_q}f_0^{*}A(c_k)\] is injective.  
        This implies that for $A\in \tau_{\leq 0}\Shv_{\tau_{\mathcal D}}(\mathcal D)$, \[ A(d)\cong \Eq\left(\prod_{i\in I}f_0^*A(c_i)\rightrightarrows \prod_{i,j\in I}\prod_{q\in Q_{i,j}}\prod_{k\in K_q}f_0^*A(c_k)\right),\] which proves that $f_0^*\colon \tau_{\leq 0}\Shv_{\tau_{\mathcal D}}(\mathcal D)\to \tau_{\leq 0}\Shv_{\tau_{\mathcal C}}(\mathcal C)$ is injective.  

    \end{proof}

\subsubsection{The condensed and local section topology}

\begin{lemma}\label{condensedtopologyonprofinite}
    For a $\kappa$-light profinite set $X$ denote by $\Cov_{\condo}(X)$ the collection of sieves $S\subseteq \oc{\Pro(\Fin)_{\kappa}}{X}$ such that there exists a finite family $\{ X_i\to X\}_{i=1}^n\subseteq S\subseteq \oc{{\Pro(\Fin)_{\kappa}}}{X}$ such that $\sqcup_{i=1}^n X_i\to X$ is onto. 
    This defines a Grothendieck topology on $\Pro(\Fin)_{\kappa}$. 
\end{lemma}
\begin{proof}
    This follows from \cite[Proposition A.3.2.1]{SAG}.
    Since $\Pro(\Fin)_{\kappa}$ has all pullbacks and coproducts and the forget functor $\Pro(\Fin)_{\kappa}\subseteq \CH\to \Set$ preserves them and is conservative, it is straightforward to check that $\Pro(\Fin)_{\kappa}$ and the continuous surjections satisfy the conditions of \cite[Proposition A.3.2.1]{SAG}, which implies that there exists a Grothendieck topology $\tau$ on $\Pro(\Fin)_{\kappa}$ such that a sieve $S\subseteq \oc{{\Pro(\Fin)_{\kappa}}}{c}$ is a covering sieve if and only if there exists a finite collection of morphisms $\{ c_i\to c\}_{i=1}^n\subseteq S\subseteq \oc{{\Pro(\Fin)_{\kappa}}}{c}$ such that 
    $\sqcup_{i=1}^n c_i\to c$ is surjective. 
\end{proof}
\begin{definition}
Fix an uncountable strong limit cardinal $\kappa$ and suppose that $\mathcal T\subseteq \operatorname{Top}_{\kappa}$ is a full subcategory. 
We say that a family $\{ X^i\to X\}_{i\in I}$ of continuous maps in $\mathcal T$ is a 
\begin{itemize}
    \item local section cover if $\bigsqcup_{i\in I} X^i\to X$ admits local sections.
    \item a $\kappa$-condensed cover if $\bigsqcup_{i\in I} \underline{X^i}_{\kappa}\to \underline X_{\kappa}$ is an epimorphism of $\kappa$-condensed sets. 
\end{itemize}
\end{definition}
\begin{lemma}\label{grostopos}Suppose $\lambda\geq \kappa$ are uncountable cardinals and denote by $\Top^{\lambda}$ the category of $\lambda$-small topological spaces. 
\begin{romanenum}\item Open coverings (respectively local section coverings) constitute a coverage on $\operatorname{Top}_{\lambda}$, respectively. 
\item The local section covers form a Grothendieck coverage. 
\item Both coverages generate the same topology, and this topology is coarser than the $\kappa$-condensed topology on $\Top^{\lambda}$. 
\end{romanenum}
\end{lemma}
    \begin{proof}
        If $f\colon Y\to X$ is continuous and $X=\bigcup_{i\in I} X_i$ is an open cover, then $Y=\bigcup_{i\in I} f^{-1}(X_i)$ is an open cover and therefore open covers constitute a coverage.  
        As isomorphisms are local section maps, and compositions and pullbacks of local section maps are local section maps, local section maps constitute a coverage on $\Top^{\lambda}$.
        Every open cover is a local section cover. Conversely, if $\{ X_i\to X\}_{i\in I}$ is a local section cover, choose an open cover $X=\bigcup_{j\in J} U_j$ such that sections $s_j$ of the identity $\operatorname{id}_{U_j}$ exist, \[ U_j\xrightarrow{s_j} \bigsqcup_{i\in I} X_i\to U_j.\] Then $\{ s_j^{-1}(X_i)\}_{i\in I,j\in J}$ is an open cover and a refinement of $\{ X_i\to X\}_{i\in I}$. 
        This shows that local section covers and open covers generate the same topology. 
        The local section covers form a Grothendieck coverage by \cite[Definition A.1, Proposition A.5]{Pstragowskisynthetic}. 
            
        Suppose $K$ is a $\kappa$-light profinite set, $f\colon K\to X$ is a continuous map to a $\lambda$-small topological space and $X=\bigcup_{i\in I}U_i$ is an open cover. 
        For $i\in I$, $f^{-1}(U_i)\subseteq K$ is open. Since $K$ is locally compact, for $k\in K$ there exists a compact neighborhood $k\in N_k$ with $N_k\subseteq f^{-1}(U_i)$ for some $i\in I$. 
        By compactness of $K$, there exists a finite set of elements $k_1, \ldots,k_n\in K$ such that $K=\bigcup_{i=1}^n N_{k_i}$.
        Hence by \cref{characterizationcondensedcovers}, $\sqcup_{i\in I}\underline{U_i}_{\kappa}\to \underline{X}_{\kappa}$ is an epimorphism of $\kappa$-condensed sets.   
        \end{proof}

\subsubsection{Slices of sheaf topoi}
In this section, we recall that the slice category of a sheaf topos over a representable sheaf is equivalent to the category of sheaves on the respective slice. 
\begin{definition}\label{topologyonslice}
    Suppose that $(\mathcal C, \tau)$ is a site and $c\in\mathcal C$. 
    For $x\in \oc{\mathcal C}{c}$ denote by $\Cov_{\tau,c}(x)$ the collection of sieves $\{ x_i\to x\}\in(\oc{\mathcal C}{c})_{/x}$ such that the image in $\oc{\mathcal C}{x}$ is a $\tau$-sieve.

    Denote by $f\colon\oc{\mathcal C}{c}\to \mathcal C$ the forget functor. 
    As \[\oc{{\mathcal C_{/c}}}{(-)}\cong \oc{\mathcal C}{(-)}\circ f\colon \oc{\mathcal C}{c}\to \Cat, \] $\tau_{c}$ defines a Grothendieck topology on $\oc{\mathcal C}{c}$. 
    \end{definition}
    
    \begin{rem}\label{remarkslicegrothendieckcoverage} For general coverages, there is no handy description of the induced topology on the slice, but for Grothendieck coverages, the following holds: 
    Suppose that $\tau$ defines a Grothendieck coverage on a category $\mathcal C$ and $c\in\mathcal C$. 
    For $d\to c\in\oc{\mathcal C}{c}$, we say that a family \[\{x_i\to d\}_{i\in I}\subseteq (\oc{\mathcal C}{c})_{/d}\cong \oc{\mathcal C}{d}\] is a $\tau_c$-cover if the underling family $\{x_i\to d\}\subseteq\oc{\mathcal C}{d}$ is a $\tau$-cover. 
    The $\tau_c$-covers form a coverage on $\oc{\mathcal C}{c}$ and the associated Grothendieck topology is the topology pulled back from the topology $[\tau]$ on $\mathcal C$ generated by $\tau$. It follows that the $\tau_c$-covers form a Grothendieck coverage.  
    \end{rem}

    By \cite[Corollary 5.1.6.12]{highertopostheory}, the functor $\oc{\mathcal C}{c}\to\oc{\mathcal P(\mathcal C)}{h_c}$ induced by the Yoneda embedding \[h_{-}\colon \mathcal C\to\mathcal P(\mathcal C)\] extends to an equivalence 
    $\mathcal P(\oc{\mathcal C}{c})\cong \oc{\mathcal P(\mathcal C)}{c}$. 
    The sheaficiation ${\mathcal P(\mathcal C)}\xrightarrow{L_{\tau}} \Shv_{\tau}(\mathcal C)$ induces a functor 
    \[ \oc{\mathcal P(\mathcal C)}{h_c}\to \oc{\Shv_{\tau}(\mathcal C)}{L_{\tau}h_c}.\]  
    \begin{lemma}\label{sheavesonslice}
        Suppose $(\mathcal C, \tau)$ is a site. 
        \begin{romanenum}
        \item The functor \[\mathcal P(\oc{\mathcal C}{c})\cong \oc{\mathcal P(\mathcal C)}{h_c}\xrightarrow{L_{\tau/c}} \oc{\Shv_{\tau}(\mathcal C)}{L_{\tau}h_c}\] 
        factors over an equivalence \[{\Shv}_{\tau_{c}}(\mathcal C_{/c})\cong \oc{{\Shv}_{\tau}(\mathcal C)}{L_{\tau}h_c}.\] 
        \item If $L_{\tau}h_c$ is hypercomplete, then this restricts to an equivalence 
        \[\hypershv_{\tau_{c}}(\oc{\mathcal C}{c})\cong \oc{\hypershv_{\tau_{c}}(\mathcal C)}{L_{\tau}h_c}.\] 
        \end{romanenum}
    \end{lemma}
    \begin{rem}
    If $\mathcal C=\tau_{\leq n}\mathcal C$ is an $(n,1)$-category ($n\in\mathbb N$), then $L_{\tau}c\in\tau_{\leq n}\topo{X}$ is hypercomplete for all $c\in\mathcal C$ by \cite[Lemma 6.5.2.9]{highertopostheory}. 
    \end{rem} 
    \begin{proof}[Proof of \cref{sheavesonslice}]
    This is a straightforward consequence of \cite[Corollary 5.1.6.12]{highertopostheory} and the description of the sheafification functor provided by \cite[Proposition 6.2.2.7]{highertopostheory}, but we could not find a reference. 
    For $F\in\oc{\mathcal P(\mathcal C)}{h_c}$, the unit $F\to L_{\tau}h_c$ and the map $F\to c$ induce a map $F\to L_{\tau}F\times_{L_{\tau}h_c}h_c$. These maps exhibit \[\mathcal P(\oc{\mathcal C}{c})\cong \oc{\mathcal P(\mathcal C)}{h_c}\xrightarrow{L_{\tau/h_c}} \oc{\Shv_{\tau}(\mathcal C)}{L_{\tau}h_c}\] as left adjoint to \[r\colon \oc{\Shv_{\tau}(\mathcal C)}{L_{\tau}h_c}\xrightarrow{-\times_{L_{\tau}h_c}h_c}\oc{\Shv_{\tau}(\mathcal C)}{h_c}\hookrightarrow \oc{\mathcal P(\mathcal C)}{h_c}.\]
    Since $L_{\tau}$ is left-exact, $L_{\tau/h_c}\circ r=1$, i.e.\ $L_{\tau/h_c}$ exhibits $\oc{\Shv_{\tau}(\mathcal C)}{L_{\tau}h_c}$ as localization of $\oc{\mathcal P(C)}{h_c}$ at the morphisms $F\to L_{\tau}F\times_{L_{\tau}h_c}h_c$. We will show that for a presheaf $F\in\mathcal P(\oc{C}{c})\cong \oc{\mathcal P(C)}{h_c}$, the map $F\to L_{\tau}F\times_{L_{\tau}h_c}h_c$ exhibits $L_{\tau}F\times_{L_{\tau}h_c}h_c$ as $\tau_c$-sheafification of $F$, then the statement follows. 
    
    Recall from \cite[section 6.2.2.7]{highertopostheory} that if $(\mathcal D, \rho)$ is a site, $\rho$-sheafification can be realised as transfinite composition of the functor \[(-)^{\dagger, \rho}\colon\mathcal P(\mathcal D)\to \mathcal P(\mathcal D)\] which sends a presheaf $F$ to $x\mapsto\colim{S\in\Cov_{\rho}(X)^{\operatorname{op}}}F(S)$.
        
    For $x\to  c\in\oc{\mathcal C}{c}$, \begin{align*}(F^{\dagger, \tau}\times_{h_c^{\dagger, \tau}}h_c)(x\to c)&=\Map_{\oc{\mathcal P(C)}{h_c}}(h_x\to h_c, F^{\dagger, \tau}\times_{h_c^{\dagger, \tau}}h_c)\\& \cong \Map_{\oc{\mathcal P(C)}{h_c^{\dagger, \tau}}}(h_x\to h_c\to h_c^{\dagger, \tau}, F^{\dagger, \tau})\\ & \cong \Map_{\mathcal P(C)}(h_x,F^{\dagger, \tau})\times_{\Map_{\mathcal P(C)}(h_x,h_c^{\dagger, \tau})}\{ h_x\to h_c\to h_c^{\dagger, \tau}\}\\ 
    &\cong \colim{S\in\Cov_{\tau}(x)^{\operatorname{op}}}\Map_{\mathcal P(C)}(S,F)\times_{\colim{T\in\Cov_{\tau}(x)^{\operatorname{op}}}\Map_{\mathcal P(\mathcal C)}(T,h_c)}\{ h_x\to h_c\}.\end{align*}
    
    For two covering sieves $S,T\in \Cov_{\tau}(x)$, $S\cap T\in \Cov_{\tau}(x)$ as for all $f\colon t\to x\in T$, \[f^*(S\cap T)=f^*(S)\cap f^*(T)=f^*(S)\in\Cov_{\tau}(t).\] This shows that $\Cov_{\tau}(x)$ is cofiltered, whence \begin{align*}\colim{S\in\Cov_{\tau}(x)^{\operatorname{op}}}&\Map_{\mathcal P(C)}(S,F)\times_{\colim{T\in\Cov_{\tau}(x)^{\operatorname{op}}}\Map_{\mathcal P(\mathcal C)}(T,h_c)}\{ h_x\to h_c\}\\ &\cong \colim{S\in\Cov_{\tau}(x)^{\operatorname{op}}}(\Map_{\mathcal P(C)}(S,F)\times_{\Map_{\mathcal P(C)}(S,h_c)}\{S\subseteq h_x\to h_c\})\\& \cong \colim{S\in\Cov_{\tau}(x)^{\operatorname{op}}}\Map_{\oc{\mathcal P(C)}{h_c}}(S\subseteq h_x\to h_c,F\to h_c)\\ &\cong F^{\dagger, \tau_{c}}(x\to c).\end{align*} 
        We now deduce from this and the description of the sheafification functors (\cite[Proposition 6.2.2.7]{highertopostheory}) that  $F\to L_{\tau}F\times_{L_{\tau}h_c}h_c$ is the $\tau_{c}$-sheafification of $F$.  
        For ordinals $\beta$ and $F\in\mathcal P(C)$ let $T^{\mathcal C}_{\beta+1}F\coloneqq T^{\mathcal C}_{\beta}F^{\dagger, \tau}$ and if $\gamma$ is a limit ordinal let $T^{\mathcal C}_{\gamma}\coloneqq \colim{\alpha<\gamma}T^{\mathcal C}_{\alpha}F$. 
        For $F\in\oc{\mathcal P(C)}{h_c}$, define $T_{\alpha}^{\oc{\mathcal C}{c}}F$ analogously. 
        We show by transfinite induction that for all ordinals $\alpha$, the map $F\to T^{\mathcal C}_{\alpha}F\times_{T^{\mathcal C}_{\alpha}h_c}h_c$ factors over an equivalence \[F\to T^{\oc{\mathcal C}{c}}_{\alpha}F\cong T^{\mathcal C}_{\alpha}F\times_{T^{\mathcal C}_{\alpha}h_c}h_c.\] 
    The zero case was shown above, the successor case follows from the above as 
        \begin{align*}T_{\alpha+1}^{\oc{\mathcal C}{c}}(F)\cong (T_{\alpha}^{\mathcal C}(F)\times_{T_{\alpha}^{\mathcal C}h_c}h_c)^{\dagger, \tau_c}& \cong (T_{\alpha}^{\mathcal C}(F)\times_{T_{\alpha}^{\mathcal C}(h_c)} h_c)^{\dagger, \tau}\times_{h_c^{\dagger, \tau}}h_c\\ & \cong (T_{\alpha+1}^{\mathcal C}(F)\times_{T_{\alpha+1}^{\mathcal C}(h_c)}h_c^{\dagger, \tau})\times_{h_c^{\dagger, \tau}}h_c\cong T_{\alpha+1}^{\mathcal C}(F)\times_{T_{\alpha+1}^{\mathcal C}(h_c)}h_c, \end{align*} where we used the induction hypothesis, the zero case and that $(-)^{\dagger, \tau}$ is left-exact (cf.\ the proof of \cite[Proposition 6.2.2.7]{highertopostheory}). The limit case follows from this since the identifications are compatible with the maps $T_{\alpha}\to T_{\beta}$ for $\alpha\leq \beta$. 
        Choose a regular cardinal $\kappa$ such that the following hold: 
        For every $x\in\mathcal C$ and every covering $S\in\Cov_{\tau}(x)$, 
        the functor \[\mathcal P(\mathcal C)\to \an, F\mapsto \Map_{\mathcal P(\mathcal C)}(S,F)\] commutes with $\kappa$-filtered colimits, and for every $x\to c\in\oc{\mathcal C}{c}$ and every covering $S\in\Cov_{\tau_c}(x\to c)$, \[\mathcal P(\mathcal C/c)\to \an, F\mapsto \Map_{\mathcal P(C/c)}(S,F)\] commutes with 
        $\kappa$-filtered colimits. This is possible since $\mathcal C$ is an essentially small category. 
        Then for all $F\in \mathcal P(\oc{C}{c})\cong \oc{\mathcal P(C)}{h_c}$, \[L_{\tau_{c}}(F)\cong T_{\kappa}F\cong T_{\kappa}F\times_{T_{\kappa}C}\{c\}\cong L_{\tau}F\times_{L_{\tau}h_c}\{h_c\}\in\mathcal P(\oc{C}{c})\cong \oc{\mathcal P(C)}{h_c}\] by \cite[Proposition 6.2.2.7]{highertopostheory} and the above. 
        The statement on hypercomplete topoi follows from \cref{sliceshypercomplete} below. 
    \end{proof} 
    
    \begin{lemma}\label{sliceshypercomplete}
    Suppose that $\topo{X}$ is a topos denote by $\widehat{\topo{X}}\subseteq \topo{X}$ the category of hypercomplete objects. 
    For $x\in\widehat{\topo{X}}$, 
    $\oc{\widehat{\topo{X}}}{x}\subseteq \oc{\topo{X}}{x}$ is the full subcategory of hypercomplete objects in $\oc{\topo{X}}{x}$. 
    \end{lemma}
    \begin{proof}
        Denote by $L_{\operatorname{hyp}}\colon\topo{X}\to \widehat{\topo{X}}$ and $L_{\operatorname{hyp},x}\colon \oc{\topo{X}}{x}\to \widehat{(\oc{\topo{X}}{x})}$ the hypersheafification functors. 
        Since $x$ is hypercomplete, the hypercompletion $L_{\operatorname{hyp}}\colon\topo{X}\to \widehat{\topo{X}}$ induces a functor 
        \[L_{\operatorname{hyp}/x}\colon \oc{\topo{X}}{x}\to \oc{\widehat{\topo{X}}}{x}\] left adjoint to the inclusion. 
        By definition (\cite[p. 657]{highertopostheory}), the forget functor $\oc{\topo{X}}{x}\to \topo{X}$ reflects $\infty$-connective morphisms. 
        In particular, for $y\to x\in \oc{\topo{X}}{x}$, the unit $y\to L_{\operatorname{hyp}/x}(y)\in \oc{\widehat{\topo{X}}}{x}$ is $\infty$-connective. 
        As $L_{\operatorname{hyp}/x}$ is left adjoint to the inclusion, for $q\colon t\to x\in\oc{\widehat{\topo{X}}}{x}\subseteq \oc{\topo{X}}{x}$, 
        \[\Map_{\oc{\topo{X}}{x}}(-,q)\cong \Map_{\oc{\widehat{\topo{X}}}{x}}(-,q)\circ L_{\operatorname{hyp}/x}\] inverts $\infty$-connective morphisms, whence $\oc{\widehat{\topo{X}}}{x}\subseteq \oc{\topo{X}}{x}$ consists of hypercomplete objects. 
        This shows that for all $y\in\oc{\topo{X}}{x}$, the unit $y\to L_{\operatorname{hyp}/x}(y)$ induces an equivalence \[L_{\operatorname{hyp}, \oc{\topo{X}}{x}}(y)\cong L_{\operatorname{hyp},x}(L_{\operatorname{hyp}/x}(y))\cong L_{\operatorname{hyp}/x}(y), \] i.e. $L_{\operatorname{hyp},x}\cong L_{\operatorname{hyp}/x}$.  
    \end{proof}

    \subsection{Homotopy invariance of gros topos cohomology}
    In this section, we show that sheaf/gros topos cohomology with constant coefficients is homotopy invariant. 
    We used this to prove \cref{deltakappaexactnesshomotopyinvariant}. 

    \begin{notation}
    Suppose that $\topo{X}$ is a topos. For a morphism $f\colon Y\to X\in \topo{X}$ let \[f^*\coloneqq Y\times_X-\colon \oc{\topo{X}}{X}\to \oc{\topo{X}}{Y}\] and denote by $f_*$ its right adjoint. 
        As $f^{*}$ is left-exact, $f^*\dashv f_{*}$ stabilizes to an adjoint pair $f^{*}_{\Sp}\dashv f_{*,\Sp}$ (\cref{geometricmorphismstabilization}). 

    For $Z\in\topo{X}$ let $\pi_{Z}^*\coloneqq Z\times -\colon \topo{X}\to \oc{\topo{X}}{Z}$. This functor is right adjoint to the forget functor $\pi_{Z,!}\colon \oc{\topo{X}}{Z}\to\topo{X}$. 
Denote by \[\pi_{Z,\Sp}^*\colon \stab{\topo{X}}\to \stab{\oc{\topo{X}}{Z}}\] the stabilization (\cref{definitionstabilization}) of $\pi_Z^*$. 
By \cite[Corollary 1.4.4.4]{higheralgebra}, $\pi_{Z,\Sp}^*$ admits a left adjoint $\pi_{Z,\Sp,!}$ and $\pi_{Z,\Sp,!}\circ \Sigma^{\infty}_{+}\cong \Sigma^{\infty}_{+}\circ \pi_{Z,!}$ (\cref{existenceleftadjointmoregeneral}). 
    \end{notation}
    \begin{lemma}\label{cohomologyisocheckedonslice}
        Suppose $f\colon X\to Y\in\topo{X}$ is a morphism in a topos $\topo{X}$. 

        For $A\in\stab{\topo{X}}$, the map $f^*\colon \cH{\topo{X}}(Y,A)\to \cH{\topo{X}}(X,A)$ is an equivalence if and only if the unit $A\to f_{*,\Sp}f^*_{\Sp}$ induces an equivalence 
        \[ \cH{\oc{\topo{X}}{X}}(X,\pi_{X,\Sp}^*A)\to \cH{\oc{\topo{X}}{X}}(X,f_{*,\Sp}f^*_{\Sp}\pi_{X,\Sp}^*A).\] 
    \end{lemma}
    \begin{proof} 
        We claim that \begin{equation}\label{diagrambeckchevalleygrostoposcohomology}
        \begin{tikzcd}[cramped,sep=small]
            \map_{\stab{\oc{\topo{X}}{X}}}(-,\pi_{X,\Sp}^*-)\arrow[r,"\eta_{f,*}"] \arrow[d,"\cong"]& \map_{\stab{\oc{\topo{X}}{X}}}(-,f_{*,\Sp}f^{*}_{\Sp}\pi_{X,\Sp}^{*}-)\arrow[r,"\cong"] & \map_{\stab{\oc{\topo{X}}{Y}}}(f^*_{\Sp}-,\underbrace{f^{*}_{\Sp}\pi_{X,\Sp}^{*}}_{\cong \pi_{Y,\Sp}^{*}}-)\arrow[d]\\ 
            \map_{\stab{\topo{X}}}(\pi_{X,\Sp,!}-,-)\arrow[rr,"\beta^{*}-"]&& \map_{\stab{\topo{X}}}(\pi_{Y,\Sp,!}\circ f^{*}_{\Sp},-) 
        \end{tikzcd}
        \end{equation} commutes, where \[\beta\colon \pi_{Y,\Sp,!}\circ f^{*}_{\Sp}\to \pi_{X,\Sp,!}\] denotes the Beck-Chevalley transformation of $f^{*}_{\Sp}\circ \pi^{X,*}_{\Sp}\cong\pi^{Y,*}_{\Sp}$ and $\eta_f$ is an adjunction unit for $f^{*}_{\Sp}\dashv f_{*,\Sp}$. 

        All functors $\stab{\oc{\topo{X}}{X}}^{\operatorname{op}}\times \stab{\topo{X}}\to \Sp$ in the above diagram preserve small limits in both variables: 
        $\pi_{X,\Sp}^*,f^*_{\Sp},f_{*,\Sp}$ and $\pi_{Y,\Sp}^*$ are stabilizations of limit-preserving functors and in particular preserve limits. 
        $f^*_{\Sp},\pi_{X,\Sp,!}$ and $\pi_{Y,\Sp,!}$ are left adjoints and therefore preserve colimits. 
        The cohomology functors preserve small limits in both variables by \cref{spectralenrichmentcocontinuous}. Hence by  \cite[Proposition 1.4.2.21]{higheralgebra}/\cref{spectraclosetoidempotentbigpresentable}, it suffices to show that the diagram obtained by applying $\Omega^{\infty}$ to the above diagram is commutative, i.e. that 
        \begin{equation}\label{diagrambeckchevalleygrostoposcohomology2}
            \begin{tikzcd}[cramped,sep=small]
            \Map_{\stab{\oc{\topo{X}}{X}}}(-,\pi_{X,\Sp}^*-)\arrow[r,"\eta_{f,*}"] \arrow[d,"\cong"]& \Map_{\stab{\oc{\topo{X}}{X}}}(-,f_{*,\Sp}f^{*}_{\Sp}\pi_{X,\Sp}^{*}-)\arrow[r,"\cong"] & \Map_{\stab{\oc{\topo{X}}{Y}}}(f^*_{\Sp}-,\underbrace{f^{*}_{\Sp}\pi_{X,\Sp}^{*}}_{\cong \pi_{Y,\Sp}^{*}}-)\arrow[d]\\ 
            \Map_{\stab{\topo{X}}}(\pi_{X,\Sp,!}-,-)\arrow[rr,"\beta^{*}-"]&& \Map_{\stab{\topo{X}}}(\pi_{Y,\Sp,!}\circ f^{*}_{\Sp}-,-) 
        \end{tikzcd}\end{equation} commutes.
        By the triangle identities for $f^{*}_{\Sp}\dashv f_{*,\Sp}$, the top horizontal arrow is the map induced by $f^{*}_{\Sp}$.
        For all $T\in \stab{\topo{X}}$, $ \Map_{\stab{\oc{\topo{X}}{X}}}(-,\pi_{X,\Sp}^*T)$ and \[\Map_{\stab{\topo{X}}}(\pi_{Y,\Sp,!}\circ f^{*}_{\Sp}-,T)\cong \Map_{\stab{\topo{X}}}(-,\pi^*_{Y,\Sp}f_{*,\Sp}-)\] are representable.  
        Hence by definition of the vertical identifications and fully faithfulness of the Yoneda embedding \[\Fun(\stab{\topo{X}}, \stab{\oc{\topo{X}}{X}})\to \Fun(\stab{{\topo{X}}},\Fun(\stab{\oc{\topo{X}}{X}}^{\operatorname{op}},\an)),\]
        it suffices to describe a homotopy between \[\pi_{Y,\Sp,!}\circ f^*_{\Sp}\circ \pi_{X,\Sp}^*\cong \pi_{Y,\Sp,!}\circ \pi_{Y,\Sp}^*\xrightarrow{\epsilon_{\pi_Y}}\id\] and
        \[\pi_{Y,\Sp,!}f^{*}_{\Sp}\pi_{X,\Sp}^*\xrightarrow{\beta}\pi_{X,\Sp,!}\pi_{X,\Sp}^*\xrightarrow{\epsilon_{\pi_X}}\id.\] 
        This exists by \cite[Lemma 2.3.3.(4)]{carmeli2022ambidexterity}. 
        This proves that \ref{diagrambeckchevalleygrostoposcohomology2} commutes, and hence so does \ref{diagrambeckchevalleygrostoposcohomology}.
        By evaluating \ref{diagrambeckchevalleygrostoposcohomology} at $\Sigma^{\infty}_{+}X$ in the first variable, we obtain a commutative diagram 
        \begin{center}
         \begin{tikzcd}
            \cH{{\oc{\topo{X}}{X}}}(X,\pi_{X,\Sp}^*-)\arrow[r,"\eta_{*}"] \arrow[d,"\cong"]& \cH{{\oc{\topo{X}}{X}}}(X,f_{*,\Sp}f^{*}_{\Sp}\pi_{X,\Sp}^{*}-)\arrow[r,"\cong"] & \cH{{\oc{\topo{X}}{Y}}}(Y,\underbrace{f^{*}_{\Sp}\pi_{X,\Sp}^{*}}_{\cong \pi_{Y,\Sp}^{*}}-)\arrow[d]\\ 
            \cH{{\topo{X}}}(X,-)\arrow[rr,"f^{*}"]&& \cH{{\topo{X}}}(Y,-). 
        \end{tikzcd}
        \end{center}
        In particular, for $M\in \stab{\topo{X}}$, $f^*\colon\cH{\topo{X}}(Y,M)\to \cH{\topo{X}}(X,M)$ is an equivalence if and only if 
        $\eta_f$ induces an equivalence \[ \cH{\oc{\topo{X}}{X}}(X,\pi_{X,\Sp}^{*}M)\cong \cH{\oc{\topo{X}}{X}}(X,f_{*,\Sp}f^{*}_{\Sp}\pi_{X,\Sp}^{*}M).\qedhere\] 
    \end{proof}

    \begin{recollection}A continous map $f\colon X\to Y$ defines a morphism of sites \begin{align*}f^{-1}\colon \Op(Y)&\to\Op(X),\\ V&\mapsto f^{-1}(V).\end{align*} 
      Denote by \[ f_*\colon \Shv_{\Zar}(Y)\leftrightarrows \Shv_{\Zar}(X)\colon f^*, f^*\dashv f_*\] the induced geometric morphism (\cref{morphismofquasisitesinducesfunctoronsheaves}) and by $f_{*,\Sp}$ the stabilization of $f_*$ (\cref{geometricmorphismstabilization}). By \cref{geometricmorphismcohomology},
        \[ \cH{\sheaf}(Y,A)\cong \cH{\sheaf}(X,f_{*,\Sp}A).\]         
        For $A\in\Shv_{\Zar}(Y)$, $f_*A=A\circ f^{-1}$. In particular, for $A\in \Ab(\tau_{\leq 0}\Shv_{\Zar}(Y))$, $\pi_0f_{*,\Sp}A=A\circ f^{-1}$. 

        The unit $\eta_f\colon \id\to f_{*,\Sp}f^*_{\Sp}$ induces a natural transformation 
        \[ \cH{\sheaf}(X,-)\to \cH{\sheaf}(X,f_{*,\Sp}f^*_{\Sp}-)\cong \cH{\sheaf}(Y,f^*_{\Sp}-).\] 
    \end{recollection}

    For a topological space $X$ denote by $c_{X}\colon \Sp\to \stab{\Shv_{\Zar}(X)}$ the stabilization of constant sheaf functor. 
    \begin{lemma}\label{homotopyinvariancesheafcohomologyconstantcoefficients}Suppose that $X$ is a topological space and denote by $\pi\colon X\times [0,1]\to X$ the projection.  
        \begin{romanenum}
    \item The unit $\id_{\Shv(X\times [0,1],\Sp)}\to \pi_{*,\Sp}\pi^*_{\Sp}$ is an equivalence. 

    \item In particular, $\pi_{*,\Sp}\circ c_{X\times [0,1]}\cong c_{X}$. 
        \end{romanenum}
    \end{lemma}
    \begin{proof}
        By \cite[Variant A.2.10]{higheralgebra}, the unit $\id\to \pi_{*}\pi^*$ is an equivalence. 
        Hence by construction of the unit for $\pi_{*,\Sp}\vdash \pi^*_{\Sp}$ (\cref{geometricmorphismstabilization}), $\id_{\Shv(X\times [0,1],\Sp)}\to \pi_{*,\Sp}\pi^*_{\Sp}$ is an equivalence. 

        As $\pi^*$ is a left-exact left adjoint, the constant sheaf functor for $X\times [0,1]$ factors as \[\an\xrightarrow{\text{const}_{X}}\Shv(X)\xrightarrow{\pi^*}\Shv(X\times [0,1]),\] and in particular, $\pi^*_{\Sp}\circ c_{X}=c_{X\times [0,1]}$. 
        It follows from the above that \[c_{X}\cong \pi_{*,\Sp}\pi^{*}_{\Sp}c_{X}\cong \pi_{*,\Sp}c_{X\times [0,1]}.\qedhere\]
    \end{proof}

    \begin{lemma}\label{smallsitesbigsitesfullyfaithful}
    Suppose $\lambda$ is an uncountable cardinal. For a $\lambda$-small topological space $X$, denote by 
    \[ \xladg{t}{X}\colon \Shv(X)\to \oc{\grostop}{h_X}\] the left adjoint of the geometric morphism induced by the morphism of sites $\Op(X)\subseteq \oc{\Top^{\lambda}}{X}$ (page \pageref{grostoposgeometricmorphismcondensed}). 
    This is fully faithful. 
    \end{lemma}
    \begin{proof}
        The functor $\xladg{t}{X}$ factors as \[\Shv(X)\subseteq \mathcal P(\Op(X))\xrightarrow{t_{!}}\mathcal P(\oc{\Top^{\lambda}}{X})\xrightarrow{L_{LS}}\Shv_{LS,x}(\oc{\Top^{\lambda}}{X})\cong \oc{\grostop}{h_X}\] where $t_{!}$ is left Kan extension along $t\colon\Op(X)\subseteq \oc{\Top^{\lambda}}{X}$, $L_{LS}$ is sheafification with respect to the Grothendieck topology induced by the local section topology on $\Top^{\lambda}$ and the right equivalence is \cref{sheavesonslice}. 
        As the local section covers form a Grothendieck coverage on $\Top^{\lambda}$ (\cref{grostopos}), by \cref{remarkslicegrothendieckcoverage}, the local section covers define a Grothendieck coverage on $\oc{\Top^{\lambda}}{X}$, and the associated topology agrees with the one pulled back from the local section topology. 
        The open covers define a Grothendieck pretopology (\cite[Definition A.1]{Pstragowskisynthetic}) on $\Op(X)$, and the local section covers define a Grothendieck pretopology on $\oc{\Top^{\lambda}}{X}$. The inclusion $t\colon \Op(X)\subseteq \oc{\Top^{\lambda}}{X}$ defines a morphism of $\infty$-sites in the sense of \cite[Definition A.3]{Pstragowskisynthetic} which has the covering lifting property. 
        Hence by \cite[Proposition A.13]{Pstragowskisynthetic}, the precomposition \[t_*\colon \mathcal P(\oc{\Top^{\lambda}}{{X}})\to \mathcal P(\Op(X)), F\mapsto F\circ t\] commutes with sheafifications. 
        In particular, for $Y\in \Shv(X)$, \[\xradg{t}{X}\xladg{t}{X}(Y)=\xradg{t}{X}(L_{\text{LS}}t_{!}Y)=t_*(L_{\text{LS}}t_{!}Y)\cong L_{\Zar}(t_*t_{!}Y).\] As $t$ is fully faithful, the unit $\id\to t_*t_{!}$ is an equivalence. 
        This implies that the unit $\id\to \xradg{t}{X}\xladg{t}{X}$ is an equivalence as well.  
    \end{proof}
    \begin{cor}\label{grostoposcohomologyomotopyinvariant}
        Suppose $\lambda$ is an uncountable cardinal. 
        Denote by $\text{const}_{\grostop,\Sp}\colon\Sp\to\stab{\grostop}$ the stabilization of the constant sheaf functor and by $h_{-}\colon \Top^{\lambda}\to \grostop$ the Yoneda embedding. 
   For $A\in\Sp$,  
    \[ \cH{\grostop}(h_{-},\text{const}_{\grostop,\Sp}A)\colon (\grostop)^{\operatorname{op}}\to \Sp\] inverts homotopy equivalences. 
    \end{cor}
    \begin{ex}
    The stabilized constant sheaf functor \[\text{const}_{\grostop,\Sp}\colon \Sp\to \stab{\grostop}\subseteq \Fun((\Top^{\lambda})^{\operatorname{op}},\Sp)\] sends an abelian group $A$ to the sheaf \[\mathcal C(-,A^{\delta})\colon (\Top^{\lambda})^{\operatorname{op}}\to \Ab\] of continuous functions into the discrete abelian group $A^{\delta}$.\end{ex}  
    \begin{proof}
        Note that a functor $F\colon (\Top^{\lambda})^{\operatorname{op}}\to \Sp$ inverts homotopy equivalences if and only if for $X\in\Top^{\lambda}$, the projection $p\colon X\times [0,1]\to X$ induces an equivalence $F(X\times [0,1])\cong F(X)$. 
        Fix $X\in\operatorname{\Top}^{\lambda}$ and denote by \[p^*\coloneqq (X\times [0,1])\times_{X}-\colon \oc{\grostop}{h_X}\rightleftarrows \oc{\grostop}{h_{X\times [0,1]}}\colon p_*, p^*\dashv p_*\] the associated basechange geometric morphism.  
    By \cref{cohomologyisocheckedonslice}, \[p^*\colon \cH{\grostop}(h_{X},A)\to\cH{\grostop}(h_{X\times [0,1]},A)\] is an equivalence if and only if the unit $\id\to p_{*,\Sp}p^*_{\Sp}$ induces an equivalence 
    \[ \cH{\oc{\grostop}{h_X}}(h_X,\pi_{X}^*A)\cong\cH{\oc{\grostop}{h_X}}(h_X,p_{*,\Sp}p^*_{\Sp}\pi_{X}^*A).\] 

    Denote by $\xladg{t}{X}\colon \Shv(X)\rightleftarrows \oc{\grostop}{h_X}\colon \xradg{t}{X}$ the geometric morphism induced by the morphism of sites $\Op(X)\subseteq \oc{\Top^{\lambda}}{X}$ (see page \pageref{grostoposgeometricmorphismcondensed}), and let \[\pi_{X}^*\coloneqq h_X\times -\colon \grostop\to\oc{\grostop}{h_X},\] \[\pi_Y^*\coloneqq h_Y\times -\colon \grostop\to \oc{\grostop}{h_Y}.\] 
    By \cref{geometricmorphismcohomology}, 
    \[ \cH{\oc{\grostop}{h_X}}(h_X,\pi_{X,\Sp}^*A)\to\cH{\oc{\grostop}{h_X}}(h_X,p_{*,\Sp}p^*_{\Sp}\pi_{X,\Sp}^*A)\] is an equivalence if and only if 
    \[ \xstradg{t}{X}(\eta)\colon \cH{\sheaf}(X,\xstradg{t}{X}\pi_{X,\Sp}^*A)\to\cH{\sheaf}(X,\xstradg{t}{X}p_{*,\Sp}p^*_{\Sp}\pi_{X,\Sp}^*A)\] is. The comutative diagram of sites 
        \begin{center}
        \begin{tikzcd}
            \Op(X)\arrow[d,"p^{-1}"] \arrow[r]& \oc{\Top^{\lambda}}{X}\arrow[d,"{(X\times [0,1])\times_X-}"]\\ 
            \Op(X\times [0,1])\arrow[r]&\oc{\Top^{\lambda}}{X\times [0,1]}
        \end{tikzcd}\end{center} yields a commutative diagram 
        \begin{center}
        \begin{tikzcd}
        \Shv(X)\arrow[r,"\xladg{t}{X}"]\arrow[d,"p^*_{\Shv}"]& \oc{\grostop}{h_X}\arrow[d,"p^*"]\\ 
        \Shv(X\times [0,1])\arrow[r,"\xladg{t}{X\times [0,1]}"] & \oc{\grostop}{h_{X\times [0,1]}}.
        \end{tikzcd}
    \end{center}
    The horizontal maps are fully faithful by \cref{smallsitesbigsitesfullyfaithful}. Denote by $\xradg{t}{X},\xradg{t}{X\times [0,1]}$ their right adjoints. 
    Fully faithfulness of $\xladg{t}{X}, \xladg{t}{X\times [0,1]}$ implies that \[\xradg{t}{X}\circ p_{*}\circ \xladg{t}{X}\] is right adjoint to $p^{*}_{\Shv}$ and \[ \id\cong\xradg{t}{X}\xladg{t}{X} \xrightarrow{\xradg{t}{X}({\eta})} p_{*,\Sp}p^{*}_{\Sp}\xstladg{t}{X}\] is a unit for this adjunction, where $\eta$ denotes adjunction unit for $p^*\dashv p_{*}$ and the left equivalence is the unit for $\xladg{t}{X}\dashv \xradg{t}{X}$.      
     Passing to stabilizations, we obtain that  \[\tilde{\eta}\colon \id\cong\xstradg{t}{X}\xstladg{t}{X} \xrightarrow{\xstradg{t}{X}({\eta})} p_{*,\Sp}p^{*}_{\Sp}\xladg{t}{X}\] exhibits \[\xstradg{t}{X}\circ p_{*,\Sp}\circ \xstladg{t}{X}\] as right adjoint to $p^{*}_{\Shv,\Sp}$. (Note that all functors are left-exact.)
     By \cref{homotopyinvariancesheafcohomologyconstantcoefficients}, $\tilde{\eta}$ is an equivalence, whence  
\[ \xstradg{t}{X}(\eta)\colon \cH{\sheaf}(X,\xstradg{t}{X}\pi_{X,\Sp}^*A)\to\cH{\sheaf}(X,\xstradg{t}{X}p_{*,\Sp}p^*_{\Sp}\pi_{X,\Sp}^*A)\] is an equivalence for all $A\in \stab{\grostop}$ with $\pi_{X,\Sp}^*A\in \operatorname{\im}(\xstladg{t}{X}\colon \Shv(X,\Sp)\hookrightarrow \stab{\oc{\grostop}{h_X}})$. 

We have shown above that this implies that \[p^*\colon \cH{\grostop}(h_{X},A)\to\cH{\grostop}(h_{X\times [0,1]},A)\] is an equivalence for all $A\in\stab{\grostop}$ with $\pi_{X,\Sp}^*A\in \operatorname{\im}(\stladg{t}{X}\colon \Shv(X,\Sp)\hookrightarrow\stab{\oc{\grostop}{h_X}})$.
As \[\pi_{X}^*\circ \text{const}_{\grostop}\cong \text{const}_{\oc{\grostop}{h_X}}\cong \xladg{t}{X}\circ \text{const}_{\Shv(X)}\] (all functors are left-exact left adjoints), it follows that \[p^*\colon \cH{\grostop}(h_{X},\text{const}_{\grostop,\Sp}A)\to\cH{\grostop}(h_{X\times [0,1]},\text{const}_{\grostop,\Sp}A)\] is an equivalence for all $A\in\Sp$. 
    \end{proof}
\subsection{\texorpdfstring{Proof of \cref{cohomologyclausenscholze}}{Proof of Theorem 2.5.8}}\label{section:proofofcohomologyclausenscholze}

In this section, we recall the proof of \cref{cohomologyclausenscholze} from \cite[Lecture 3]{Scholzecondensed}. 
We start with some preparation. 
\begin{lemma}\label{colimitcomparisonmapinjectivewithdenseimage}
    Suppose $A$ is the product of a discrete abelian group and a finite-dimensional normed $\mathbb R$–vector space. For a filtered system of finite sets $(X_k)_{k\in K}$, the canonical map 
\[\colim{k\in K}\mathcal C(X_k,A)\to \mathcal C(\clim{k\in K}X_k,A)\] is injective and has dense image where both sides are endowed with the compact open topology. 
\end{lemma}
\begin{proof}For $k\in K$ denote by $\pi_k\colon \clim{k}X_k\to X_k$ the projection. 
Suppose that \[(f\colon X_k\to A), (g\colon X_l\to A)\] with $g\circ \pi_l=f\circ \pi_k$. 
As $K$ is filtered, we can assume that $k=l$. 
Let $C\coloneqq X_k\setminus \im(\pi_k)$. This is finite. As $\emptyset=\clim{l\in K_{/k}}(X_l\times_{X_k}C)$ and $K_{/l}$ is filtered, by \cref{inverselimitnonempty}, there exists $l\to K\in K$ with $X_l\times_{X_k}C=\emptyset$, i.e. $\operatorname{\im}(X_l\to X_k)\subseteq \operatorname{\im}(\pi_k\colon \clim{k}X_k\to X_k)$. 
This implies that the two maps \[f\circ X(l\to k), g\circ X(l\to k)\colon X_l\to X_k\to A\] agree, which proves that the map is injective. 
We now show that it has dense image. Choose a metric on $A$ (defining the topology). 
Supppose that $f\colon \clim{k\in K}X_k\to A$ and fix $\epsilon>0$. 
For $x\in \clim{k\in K}X_k\to A$, there exists an open neighbourhood $x\in U$ such that $f(U)\subseteq D_{\epsilon}(f(x))$. By definition of the limit topology and since $K$ is filtered, there exist $k_x\in K$, $U_x\subseteq X_{k_x}$ open with $\pi_{k_x}^{-1}(U_{x})\subseteq U$. 
By compactness of $\clim{k\in K}X_k$, there exists a finite collection $\{x_1,\ldots,x_m\}\subseteq \clim{k\in K}X_k$ with $\cup_{j=1}^m \pi_{k_{x_j}}^{-1}(U_{{x_j}})$. Using filteredness of $K$ again, we can reduce to the case where $k_{x_j}=k_{x_1}=:k$ for all $1\leq j\leq m$. 
That is, we have found $k\in K$, a finite set $\{x_1,\ldots,x_m\}\subseteq \clim{k\in K}X_k$, and a collection of opens $(U_j)_{1\leq j\leq m}\subseteq X_k$ such that $\bigcup_{j=1}^m\pi_k^{-1}(U_j)=\clim{k\in K}X_k$ and for all $1\leq j\leq m$, $f(\pi_{k}^{-1}(U_j))\subseteq D_{\epsilon}(f(x_j))$.  
Let \begin{align*}f_k\colon F_k&\to A\\ x&\mapsto \begin{cases} f(x_j) & x\in U_j\setminus \cup_{l=1}^{j-1}U_l \\ 
0 & x\notin \cup_{k=1}^m U_k.
\end{cases}\end{align*} 
As $\clim{k\in K}X_k=\cup_{j=1}^m \pi_k^{-1}(U_j\setminus \cup_{l=1}^{j-1}U_l)$, and for $x\in \pi_{k}^{-1}(U_j\setminus \cup_{l=1}^{j-1}U_l)$, \[f_k(\pi_k(x))=f(x_j)\in D_{\epsilon}(f(x))\] (by choice of $U_j$), it follows that 
$\max_{x\in \clim{k\in K}X_k}d(f(x),f_k\circ \pi_k(x))<\epsilon$, which proves that the image of \[\colim{k\in K}\mathcal C(X_k,A)\hookrightarrow \mathcal C(\clim{k\in K}X_k,A)\] is dense. 
\end{proof}
\begin{notation}
For a continuous map $Y\to X$ and a topological abelian group $A$, denote by $\mathcal C(\check{C}(Y\to X),A)$ the cosimplicial abelian group given by continous maps from the \v{C}ech nerve of $Y\to X$ to $A$, and by $[\check{\mathcal C}(Y\to X),A]$ the unnormalized Moore complex (\cite[Definition 1.2.3.8]{higheralgebra}) of this cosimplicial abelian group.\end{notation} 
\begin{lemma}[{\cite[Theorem 3.2, 3.3]{Scholzecondensed}}]\label{exactnesscheckcomplex}
    \vspace{0pt}\noindent
\begin{romanenum}
    \item Suppose $A$ is a discrete abelian group and $S\to X$ is a continuous surjection between profinite sets. 
    Then $[\check{\mathcal C}(S\to X),A]$ augments to a resolution of $\mathcal C(X,A)$.
\item If $X$ is a compact Hausdorff space and $S\to X$ is a continuous surjection from a profinite set $S$, then 
$[\check{\mathcal C}(S\to X),\mathbb R]$ augments to a resolution of $\mathcal C(X,\mathbb R)$.  
\end{romanenum} 
\end{lemma}
\begin{proof}
We proceed as in the proofs of {\cite[Theorem 3.2, 3.3]{Scholzecondensed}}.
Suppose first that $S\to X$ is a continuous surjection of profinite sets and $A$ is a discrete abelian group.  
We can write \[S\to X=\clim{j\in J}(S_j\to X_j)\] as filtered limit over surjections of finite sets $f_j\colon S_j\to X_j$ (see e.g.\ \cite[Lemma 2.1.6]{LucasMannthesis}).

For all $j\in J$, $[\check{C}(S_j\to X_j),A]$ is a resolution of $\mathcal C(X_j,A)$:
\cref{globalsectionstexactcocontinuous,cohomologyglobalsections} imply that for all finite sets $F$, \[\ckH(\underline{F}_{\kappa},\underline{A}_{\kappa})\cong \prod_{f\in F}\ckH(*,\underline{A}_{\kappa})\cong \mathcal C(F,A)\] is concentrated in degree $0$. Together with \cref{Bousfieldkanspectralsequencehomology}, this implies that \[H^q([\check{C}(S_j\to X_j),A])\cong \ckH^q(\underline{X_j}_{\kappa},\underline{A}_{\kappa})\] for all $q\in\mathbb Z$, i.e. 
$[\check{C}(S_j\to X_j),A]$ augments to a resolution of $\mathcal C(X_j,A)$. 
It follows that $\colim{j\in J}[\check{C}(S_j\to X_j),A]$ augments to a resolution of $\colim{j\in J}\mathcal C(S_j,A)$.
Since $A$ has the discrete topology, the compact open topology on \[\mathcal C(\check{C}(S\to X)([p]),A)=\mathcal C(\underbrace{S\times_{X}\ldots \times_{X}S}_{p+1 \text{ times }},A)\] is discrete for all $p\in\mathbb N_0$. \cref{colimitcomparisonmapinjectivewithdenseimage} now implies the canonical map \[\colim{j\in J}[\check{C}(S_j\to X_j),A]\to [\check{C}(S\to X),A]\] is an isomorphism of chain complexes, whence $[\check{C}(S\to X),A]$ augments to a resolution of \[\colim{j\in J}\mathcal C(X_j,A)=\mathcal C(X,A).\]  

Suppose now that $S\to X$ is a surjection from a profinite set $S$ to a compact Hausdorff space $X$. 
As \[H^0([\check{C}(S\to X),\mathbb R])=\operatorname{Ker}\left(\pi_2^*-\pi_1^*\colon\mathcal C(S,\mathbb R)\to C(S\times_X S,\mathbb R)\right)=\mathcal C(X,\mathbb R),\] it suffices to show that \begin{align}\label{exactnesschaincomplex} H^n([\check{C}(S\to X),\mathbb R])=0\text{ for } n>0.\end{align}
For $n\in\mathbb N_0$, endow \[[\check{C}(S\to X),\mathbb R]_n=\mathcal C(\underbrace{S\times_X \ldots \times_X S}_{n+1 \text{ times }},\mathbb R)\] with the sup-norm. 
Fix $n\geq 1$. 
We claim that for $0\neq f\in \ker(d_n)$ and $1>\epsilon>0$, there exists \[b\in \mathcal C(\check{C}(S\to X)([n-1]),\mathbb R)\] with $||b||<(1+\epsilon)||f||$ and $||f-d_{n-1}b||<\epsilon||f||$. 
Let us first explain how this implies \ref{exactnesschaincomplex}: Fix $0<\epsilon<1$, let $f^0\coloneqq f$ and define a sequence of elements \[(b_k)_{k\in\mathbb N_1}\in \mathcal C(\check{C}(S\to X)([n-1]),\mathbb R),\,  (f^k)_{k\in\mathbb N_0}\in \ker(d_n)\] recursively as follows: 
If $(b_k)_{k\leq m-1}$ and $(f^k)_{k\leq m-1}$ have been defined, choose $b_{m}$ such that \[||b_m||<(1+\epsilon)||f^{m-1}||\]and \[||f^{m-1}-d_{n-1}b_m||<\epsilon||f^{m-1}||,\] and set $f^m\coloneqq f^{m-1}-d_{n-1}b_m$. 
It follows inductively that for all $k\in\mathbb N_0$, $||f^k||<\epsilon||f^{k-1}||<\epsilon^k||f||$, and hence \[||b_k||<(1+\epsilon)\epsilon^{k-1}||f||.\] In particular, $b\coloneqq \sum_{k\in\mathbb N_0}b_k$ exists. 
As \begin{align*}||f-d_{n-1}b||&<||\underbrace{f-\sum_{j=1}^md_{n-1}b_j}_{=f^{m}}||+\sum_{j> m}||d_{n-1}b_j||\\ &<\epsilon^m||f||+\sum_{j>m}(n+1)||b_j||\\&<\epsilon^m||f||+(n+1)\sum_{j> m}(1+\epsilon)\epsilon^{j-1}||f||\\&=\epsilon^{m}||f||(1+(n+1)\frac{1+\epsilon}{1-\epsilon})\xrightarrow{m\to\infty}0,\end{align*}  $f=d_{n-1}b$.
Note also that \[||b||\leq \sum_{k\in\mathbb N_0}||b_k ||\leq (1+\epsilon)\sum_{k\in\mathbb N_0}\epsilon^k||f||= \frac{1+\epsilon}{1-\epsilon}||f||.\] 

We now prove the claim. Suppose first that $X$ is profinite and write $(S\to X)=\clim{j\in J}(S_j\to X_j)$ as inverse limit of surjections. 
For $j\in J$, choose a section $s_j\colon X_j\to S_j$. 
This defines a splitting of the augmented simplicial object $\check{C}(S_j\to X_j)\to X_j$. This splitting induces a contracting chain homotopy $h_j$ of the associated augmented cochain complex $\mathcal C(X_j,\mathbb R)\to [\check{C}(S_j\to X_j),\mathbb R]$, see e.g.\ \cite[\href{https://kerodon.net/tag/04SE}{Tag 04SE}]{kerodon}. 
In degree $m$, \[h_{j,m}\colon \mathcal C(\check{C}(S_j\to X_j)([m]),\mathbb R)\to \mathcal C(\check{C}(S_j\to X_j)([m-1]),\mathbb R)\] is pullback along a map $\check{C}(S_j\to X_j)[m-1]\to \check{C}(S_j\to X_j)[m]$ and in particular of norm $\leq 1$. 
For $m\in\mathbb N_0$, the boundary map \[d_m\colon  [\check{C}(S\to X),\mathbb R]_{m}\to [\check{C}({S\to X}),\mathbb R]_{m+1}\] is the alternating sum of the restrictions along the $m+2$ projections \[ \check{C}({S\to X})([{m+1}])=\underbrace{S\times_X \times \ldots\times_X S}_{m+2 \text{ times }}\to \underbrace{S\times_X \times \ldots\times_X S}_{m+1 \text{ times }}=\check{C}(S\to X)([m])\] and in particular has norm $\leq m+2$. 
Suppose $f\in \ker(d_n)$ and choose $j\in J$ and \[g_0\in[\check{C}(S_j\to X_j),\mathbb R]_{n}\] with $||g_0\circ \pi_j-f||<\frac{\epsilon}{4^{n}}||f||$. This is possible since \[\colim{j}[\check{C}(S_j\to X_j),\mathbb R]_{n}\subseteq [\check{C}(S\to X),\mathbb R]_{n}\] is dense (\cref{colimitcomparisonmapinjectivewithdenseimage}). 
As $d_n$ has norm $\leq n+2$, \[||d_n(g_0\circ \pi_j)||=||d_n(g_0\circ \pi_j-f)||\leq (n+2)||g_0\circ \pi_j-f ||< \frac{n+2}{4^n}\epsilon||f||.\] 
Choose $k\to j\in J$ such that $\im(\pi_j\colon S\to S_j)=\im(S_k\to S_j)$. This exists by \cref{inverselimitnonempty}.
Then \[||g_0\circ \pi_{j}||_S=||g_0\circ (S_k\to S_j)||_{S_k}\] and \[||d_n(g_0\circ \pi_j)||_{S}=||d_ng_0\circ (S_k\to S_j)||_{S_k}.\]  
As $h_{k,n+1}\colon [\check{C}(S_k\to X_k),\mathbb R]_{n+1}\to [\check{C}(S_k\to X_k),\mathbb R]_{n}$ has norm $\leq 1$, it follows that 
\[||h_{k,n+1}(d_n(g_0\circ (S_k\to S_j)))||_{S_k}\leq ||d_n(g_0\circ (S_k\to S_j))||_{S_k}=||d_n(g_0\circ \pi_j)||_{S}<\frac{n+2}{4^n}\epsilon||f||_S.\] 
Let \[b\coloneqq h_{k,n}(g_0\circ (S_k\to S_j))\circ \pi_k.\] Since $h_{k,n}$ has norm $\leq 1$, \[||b||_S\leq ||g_0\circ (S_k\to S_j)\circ \pi_k||_S=||g_0\circ \pi_j||_S\leq ||f||_S+||g_0\circ \pi_j-f||_S< ||f||_S+\frac{\epsilon}{4^n}||f||_S<(1+\epsilon)||f||_S.\] 
As $\id_{[\check{C}(S_k\to X_k),\mathbb R]_n}=d_{n-1}h_{k,n}+h_{k,n+1}d_{n}$, \[d_{n-1}b=g_0\circ \pi_j-(h_{k,n+1}d_n(g_0\circ (S_k\to S_j)))\circ \pi_k,\] whence \begin{align*}||f-d_{n-1}b||_S&\leq ||f-g_0\circ \pi_j||_S+ ||(h_{k,n+1}d_n(g_0\circ (S_k\to S_j)))\circ \pi_k||_S\\ & \leq ||f-g_0\circ \pi_j||_S+ ||(h_{k,n+1}d_n(g_0\circ (S_k\to S_j)))||_{S_k}\\ &< (\frac{1}{4^n}+\frac{n+2}{4^n})\epsilon||f||_S\leq\epsilon||f||_S\end{align*} by the above. 
This proves the claim for $X$ profinite. 

Suppose now that $X$ is a general compact Hausdorff space. 
For $n\in\mathbb N_0$ let \[S(n)\coloneqq \check{\mathcal C}(S\to X)([n])=\underbrace{S\times_X \ldots \times_X S}_{n+1 \text{ times }}.\] Fix $n\in\mathbb N_1$. 
For $x\in X$, $\check{C}(S\to X)\times_X \{x\}\cong \check{C}(S\times_{X}\{x\}\to \{x\})$. As $S\times_{X}\{x\}$ is profinite, the above applied to $\tilde{\epsilon}\coloneqq \frac{\epsilon}{4}$ shows that for all $x\in X$, there exists $b_x\colon S(n-1)\times_{X}\{x\}\to \mathbb R$ with $(d_{n-1}b_x)(x)=f_n(x)$ and \[||b_x||\leq \frac{1+\tilde{\epsilon}}{1-\tilde{\epsilon}}||f||_{S({n-1})\times_{X}\{x\}}< (1+\epsilon)||f||_{S({n-1})}.\]
By Tietzes extension theorem we can extend $b_x$ to a map $\tilde{b}_x\colon S({n-1})\to \mathbb R$ with \[||\tilde{b}_x||_{S({n-1})}\leq ||b_{x}||_{S(n-1)\times_{X}\{x\}}< (1+\epsilon)||f||.\]  
Choose an open neighbourhood $x\in U_x\subseteq X$ such that $||d_{n-1}\tilde b_x-f||_{S({n-1})\times_{X}U_x}<\epsilon||f||_X$. 
By compactness of $X$ we find $x_1,\ldots, x_k\in S(n)$ such that $X=\bigcup_{j=1}^k U_{x_j}$. Choose a partition of unity $1=\sum_{j=1}^k\rho_j$ subordinate to this cover, i.e. for all $1\leq j\leq k$, $\rho_j\colon X\to [0,1]$ has closed support contained in $U_j$.  
For $1\leq j\leq k$ let \[\rho^{n-1}_j\colon S({n-1})\to X\xrightarrow{\rho_j}[0,1]\] and define $\rho^n_j$ analogously. 
Let $b\coloneqq \sum_{j=1}^k\rho^{n-1}_j\tilde{b}_{x_{j}}$.
Then \[||b||\leq \max_{j=1,\ldots,k}||\tilde{b}_{x_j}||< (1+\epsilon)||f||\] and $d_{n-1}(b)=\sum_{j=1}^k\rho^n_j d_{n-1}(\tilde{b}_{x_j})$, whence 
\[||f-d_{n-1}(b)||=||\sum_{j=1}^k\rho_j(f-d_{n-1}\tilde{b}_{x_{j}})||< \epsilon||f||.\]
This proves the claim for $X$ compact Hausdorff. 
\end{proof}
\begin{cor}[{\cite[Theorem 3.2, 3.3]{Scholzecondensed}}]\label{cohomologyclausenscholzeprofinite}
    \vspace{0pt}\noindent
\begin{romanenum}
    \item If $X$ is a ($\kappa$-light) profinite space and $A$ is a product of a finite-dimensional normed $\mathbb R$-vector space with a discrete abelian group, then $\cckH{(\kappa-)}(\underline{X}_{\kappa},\underline{A}_{\kappa})\cong \mathcal C(X,A)$ is concentrated in degree $0$.   
\item If $X$ is a compact Hausdorff space and $A$ is a finite-dimensional normed $\mathbb R$-vector space, then $\cckH{(\kappa-)}(\underline{X}_{\kappa},\underline{A}_{\kappa})\cong \mathcal C(X,A)$ is concentrated in degree $0$.   
\end{romanenum}  
\end{cor}
\begin{proof}We proceed as in the proofs of {\cite[Theorem 3.2, 3.3]{Scholzecondensed}}. By \cref{condensedcohomologycanbecomputedonfinitestage}, it suffices to prove the $\kappa$-condensed statement. 
Recall from \cref{condensedoncompactextremallydisconnected} that $\Cond{\kappa}(\an)\cong \widehat{\Shv}_{\condo}(\CH_{\kappa})$ 
is equivalent to the category of hypersheaves on $\CH_{\kappa}$ with respect to the topology generated by finite, jointly surjective maps. 
\cref{compacthausdorffquotientoftdch} implies that the covers $\{X_i\to X\}_{i\in I}$ with $I$ finite and $\sqcup_{i\in I}X_i\in\Pro(\Fin)_{\kappa}$ are cofinal among the covers of this topology on $\CH_{\kappa}$. 
This implies the following, see e.g.\ \cite[\href{https://stacks.math.columbia.edu/tag/03F9}{Tag 03F9}]{stacks-project}:
Suppose $F\in\Cond{\kappa}(\Ab)$ and $X\in \CH_{\kappa}$. If for all covers $\{X_i\to X\}_{i\in I}$ with $I$ finite and $\sqcup_{i\in I}X_i\in\Pro(\Fin)_{\kappa}$, the \v{C}ech complex $[\check{C}(\{X_i\to X\}),F]$, that is the unnormalized Moore complex (\cite[Definition 1.2.3.8]{higheralgebra}) of the cosimplicial abelian group \[\Map_{\Cond{\kappa}(\Set)}(\check{C}(\sqcup_{i\in I}\underline{X_i}_{\kappa}\to \underline{X}_{\kappa}),F)=F(\check{C}(\sqcup_{i\in I}X_i\to X)),\] is a resolution of $F(X)$, then $\ckH^*(\underline{X}_{\kappa},F)=F(X)$ is concentrated in degree $0$. 

It now follows from \cref{exactnesscheckcomplex,kappacontinuousfullyfaithfullyintocondensed} that for $A$ discrete and $X$ $\kappa$-light profinite, $\ckH(X,A)=A(X)=\mathcal C(X,A)$ is concentrated in degree $0$. Analogously, for $X$ $\kappa$-light compact Hausdorff and $A=\mathbb R$, $\ckH(X,A)=A(X)=\mathcal C(X,A)$ is concentrated in degree $0$. 
As condensed cohomology commutes with products in the coefficient groups, this implies the statement. 
\end{proof}
    \begin{recollection}[Stalk functors]\label{stalkfunctors}
Suppose $x\in X$ is a point. The inclusion $i_x\colon \{x\}\subseteq X$ defines a morphism of sites $\Op(X)\to \Op(\{x\})$.
Denote by $x^*\colon\Shv(X)\to \Shv(\{x\})\cong \an\colon x_*, x^*\dashv x_*$ the induced geometric morphism (\cref{morphismofquasisitesinducesfunctoronsheaves}). 

As $x^*$ is the left Kan extension of \begin{align*}\Op(X)\to \an, U\mapsto \begin{cases}* & x\in U\\ \empty & x\notin U \end{cases}\end{align*} along the Yoneda embedding $h_U\colon \Op(X)\to \Shv(X)$, \[x^*(-)=\colim{x\in U\in \Op(X)^{\operatorname{op}}}\Map_{\Shv(X)}(h_U,-).\] 
By \cite[Lemma A.3.9]{higheralgebra}, the stalk functors $x^*\colon\Shv(X)\to \an$ jointly detect $\infty$-connective morphisms. 
In particular, the induced functors $x^*\colon \widehat{\Shv(X)}\to \an, x\in X$ are jointly conservative. 

For $x\in X$ denote by $(-)_x\colon \stab{\Shv(X)}\to \Sp$ the stabilization of $x^*$.
As for all $U\in \Op(X)$, $\cH{\Shv(X)}(h_U,-)\colon \stab{\Shv(X)}\to \Sp$ is the stabilization of $\Map_{\Shv(X)}(h_U,-)\colon \Shv(X)\to \an$ and $x\in U\in\Op(X)^{\operatorname{op}}$ is filtered, 
 \[(-)_x=\colim{x\in U\in\Op(X)^{\operatorname{op}}}\map_{\stab{\Shv(X)}}(\Sigma^{\infty}_{+}h_U,-)= \colim{\substack{x\in U\\ U\in\Op(X)^{\operatorname{op}}}}\cH{\Shv(X)}(h_U,-).\]
The inclusion $\widehat{\Shv}(X)\subseteq \Shv(X)$ stabilizes to a fully faithful functor $\stab{\widehat{\Shv}(X)}\subseteq \stab{\Shv(X)}$. 
The stalk functors $(-)_{x,\Sp}\colon \stab{\widehat{\Shv}(X)}\to \Sp$ are jointly conservative: For $i\in\mathbb N_0$, \[\Omega^{\infty-i}\circ (-)_x=x^*\circ \Omega^{\infty-i}\text{ by \cref{geometricmorphismstabilization}},\] and the functors $x^*\circ \Omega^{\infty-i}\colon \widehat{\Shv}(X)\to \an, x\in X, i\in\mathbb N_0$ are jointly conservative by \cite[Lemma A.3.9]{higheralgebra} and \cite[Proposition 1.3.2.27]{SAG}.
    \end{recollection}
    
\cref{cohomologyclausenscholze} is now a straightforward consequence: 
\begin{proof}[Proof of \cref{cohomologyclausenscholze}]
By \cref{condensedcohomologycanbecomputedonfinitestage}, it suffices to prove the $\kappa$-condensed statement. 
Note that if $X$ is ($\kappa$-light) compact Hausdorff, for all topological abelian groups $A$, $\mathcal C_{\kappa/k}(-,A)=\mathcal C(-,A)$, respectively. 
By construction of the comparison map, it suffices to show that for a $\kappa$-light compact Hausdorff space $X$, the map $\pi_0\stradg{j}\underline{A}_{\kappa}\to \stradg{j}A\in\Shv(X,\Sp)$ is an equivalence. As $\radg{j}$ is a right adjoint of a geometric morphism and $\Cond{\kappa}(\an)$ is hypercomplete,  \[\stradg{j}A\in \hypershv(X,\Sp)\coloneqq \stab{\hypershv(X)}\subseteq \stab{\Shv(X)}\] is hypercomplete by \cite[Proposition 6.5.1.16 (4)]{highertopostheory}, so it suffices to show that $\pi_k\stradg{j}(A)=0$ for $k\neq 0$. This can be checked on stalks (\cref{stalkfunctors}). 
By \cref{geometricmorphismcohomology}, the stalk at $x\in X$ is \[(\stradg{i}\underline{A}_{\kappa})_x=\colim{x\in U\in \Op(X)^{\operatorname{op}}}\cH{\sheaf}(U,\stradg{i}\underline{A}_{\kappa})\cong \colim{x\in U\in\Op(X)^{\operatorname{op}}}\ckH(\underline{U}_{\kappa},\underline{A}_{\kappa}).\] 
As every open neighbourhood of $x$ contains a closed neighbourhood of $x$ and vice versa, \[\colim{x\in U\in\Op(X)^{\operatorname{op}}}\ckH(\underline{U}_{\kappa},\underline{A}_{\kappa})\cong \colim{\substack{x\in C\\ C\text{ clsd nbhd}}}\ckH(\underline{C}_{\kappa},\underline{A}_{\kappa})\] where the colimit on the right is over the poset of closed neighbourhoods $C$ of $x$ in $X$. 
We claim that pullback along the inclusions $x\to C$ defines an isomorphism \[ \colim{\substack{x\in C\\ C\text{ clsd nbhd}}}\ckH(\underline{C}_{\kappa},\underline{A}_{\kappa})\cong \ckH(\{x\},\underline{A}_{\kappa})=\underline{A}_{\kappa}[0],\] this then implies that $\pi_k\stradg{j}(A)=0$ for $k\neq 0$. 

As both sides commute with finite products in $A$, it suffices to prove this for $A=\mathbb R$ and $A$ discrete. 
If $A=\mathbb R$, then \[\ckH(\underline{C}_{\kappa},\underline{\mathbb R}_{\kappa})\cong \mathcal{C}(C,\mathbb R)\] is concentrated in degree $0$ by \cref{cohomologyclausenscholzeprofinite}, whence \[\colim{\substack{x\in C\\ C\text{ clsd nbhd}}}\ckH(\underline{C}_{\kappa},\underline{\mathbb R}_{\kappa})\cong \colim{\substack{x\in C\\ C\text{ clsd nbhd}}}\mathcal C(C,\mathbb R)\cong \colim{{x\in U\in \Op(X)^{\operatorname{op}}}}\mathcal C(U,\mathbb R)\cong \mathcal C(\{x\},\mathbb R).\]  

Suppose now that $A$ is discrete and choose a continuous surjection $Y\to X$ from a $\kappa$-light profinite set. This exists by \cref{compacthausdorffquotientoftdch}. Then $\check{C}(q)\to X$ is a condensed hypercover of $X$ by $\kappa$-light profinite sets, and for a  closed subspace $Z\subseteq X$, \[q_Z\colon \underline{Z\times_Y X}_{\kappa}\to\underline{Z}_{\kappa}\] is an effective epimorphism in $\Cond{\kappa}(\an)$. 
As for all $m\in\mathbb N_0$, $\check{C}(q_Z)([m])=\underline{\check{C}(q)([m])\times_X Z}_{\kappa}$ and $\check{C}(q)([m])\times_X Z$ is a $\kappa$-light profinite space (as closed subspace of the $\kappa$-light profinite space $Y^{m+1}$), \cref{cohomologyclausenscholzeprofinite,Bousfieldkanspectralsequencehomology} imply that for a closed subspace $Z\subseteq X$,
\[\ckH^q(\underline{Z}_{\kappa},\underline{A}_{\kappa})=H^q([\check{C}(Y\times_X Z\to  Z),A])\] is the cohomology of the unnormalized Moore complex (\cite[Definition 1.2.3.8]{higheralgebra}) of the cosimplicial abelian group $\mathcal C(\check{C}(Y\times_X Z\to Z),A)$.
Since $A$ is discrete, the map \[ \colim{\substack{x\in C\\ C\text{ clsd nbhd}}}[\check{C}(Y\times_X C\to C),A]\to [\check{C}(Y\times_X \{x\}\to \{x\}),A]\] is an isomorphism of chain complexes.  
\end{proof}

\newpage

\printbibliography
\end{document}